\documentclass[11pt]{report}
\usepackage[english]{babel}
\usepackage[utf8]{inputenc}

\usepackage{amsmath,amsfonts,amssymb,amsthm,mathrsfs}
\usepackage{amsbsy,mathtools,stmaryrd,upgreek,mathabx,eufrak}

\renewcommand{\leq}{\leqslant}
\renewcommand{\geq}{\geqslant}

\usepackage[mono=false]{libertine}

\newcommand*{\match}[4]{%
\ifx#1(%
	#3(#2#4)%
\else\ifx#1[%
	#3[#2#4]%
\else\ifx#1.%
	#3\{#2#4\}%
\else\ifx#1<%
	#3\langle#2#4\rangle%
\else\ifx#1/%
	#3\lfloor#2#4\rfloor%
\else\ifx#1!%
	#3\lceil#2#4\rceil%
\else%
	#3#1#2#4#1
\fi\fi\fi\fi\fi\fi}

\newcommand{\indset}[1]{\mathbf{1}_{#1}}

\usepackage{wasysym}

\usepackage[dvipsnames]{xcolor}
\usepackage[a4paper,vmargin={3.5cm,3.5cm},hmargin={2.5cm,2.5cm}]{geometry}
\usepackage[font=sf, labelfont={sf,bf}, margin=1cm]{caption}
\usepackage{graphicx,graphics}
\usepackage{epsfig}
\usepackage{latexsym}
\usepackage{stmaryrd}
\usepackage{ae,aecompl}

\usepackage[english]{babel}
 \usepackage[colorlinks=true]{hyperref}
\usepackage{pstricks}
\usepackage{enumerate}
\usepackage{cleveref}

\newcommand{\Z}{\mathbf{Z}}
\renewcommand{\P}{\mathbb{P}}
\newcommand{\E}{\mathbb{E}}
\DeclareFixedFont{\beaupetit}{T1}{ftp}{b}{n}{2cm}

\newtheorem{theorem}{Theorem}[chapter]
\newtheorem*{result}{Informal Result}
\newtheorem{definition}[theorem]{Definition}
\newtheorem{proposition}[theorem]{Proposition}
\newtheorem{assumption}[theorem]{Assumption}

\newtheorem{lemma}[theorem]{Lemma}
\newtheorem{corollary}[theorem]{Corollary}

\theoremstyle{definition}
\newtheorem{remark}[theorem]{Remark}
\newtheorem{notation}[theorem]{Notation}
\newtheorem{example}[theorem]{Example}

\def\llbracket{[\hspace{-.10em} [ }
\def\rrbracket{ ] \hspace{-.10em}]}

\newcommand{\eps}{\varepsilon}

\usepackage{todonotes}

\newcommand{\TT}{\mathbb{T}}

\newcommand{\N}{\mathbb{N}}
\newcommand{\R}{\mathbb{R}}
\newcommand{\U}{\mathbb{U}}
\newcommand{\Q}{\mathbb{Q}}

\newcommand{\esp}{\varepsilon}
\def\d{{\mathrm d}}

\newcommand{\e}{{\mathrm{e}}}

\newcommand{\GW}{ {\normalfont\text{\tiny{GW}}}}

\newcommand{\dd}{\d }

\usepackage{titling} 

\thanksmarkseries{fnsymbol} 

\makeatletter
\renewcommand*{\@fnsymbol}[1]{\ensuremath{\ifcase#1\or *\or \dagger\or \ddagger\or
   \mathsection\or \mathparagraph\or \|\or **\or \dagger\dagger
   \or \ddagger\ddagger \else\@ctrerr\fi}}
\makeatother

\title{\Huge{\bf Self-similar Markov trees\\ and scaling limits}}
\author{Jean Bertoin\thanks{Universit\"at Z\"urich, Institute of Mathematics, Switzerland.\hfill  \texttt{jean.bertoin@math.uzh.ch}} ~ and ~ Nicolas Curien\thanks{Universit\'e Paris-Saclay, LMO, France.\hfill  \texttt{nicolas.curien@gmail.com}} ~ and ~ Armand Riera\thanks{Sorbonne Universit\'e, LPSM, France.\hfill  \texttt{riera@lpsm.paris}}.} 
\date{}
             
\begin{document}
  
     \maketitle
     \thispagestyle{empty}
   \centerline{To the memory of John W. Lamperti (1932-2024)}

             \begin{abstract}

Self-similar Markov trees constitute a remarkable family of random compact real trees carrying a decoration function that is positive on the skeleton. As the terminology suggests, they are self-similar objects that further satisfy a Markov branching property. They are built from the combination of the recursive construction of real trees by gluing line segments with the seminal observation of Lamperti, which relates positive self-similar Markov processes and Lévy processes via a time change. They carry natural length and harmonic measures, which can be used to perform explicit spinal decompositions. Self-similar Markov trees encompass a large variety of random real trees that have been studied over the last decades, such as the Brownian CRT, stable Lévy trees, fragmentation trees, and growth-fragmentation trees. We establish general invariance principles for Galton--Watson trees with integer types and illustrate them with many combinatorial classes of random trees that have been studied in the literature.
\end{abstract}
\tableofcontents

\chapter*{Index of some notation and terminology}
\addcontentsline{toc}{chapter}{Index of some notation and terminology}

We first recall a sample of a few standard notation which are implicitly used throughout the text.
\medskip

 \paragraph{Basic notation}\ \\ 

\begin{tabular}{cl}
$ \mathbb{N}$& set of positive natural numbers $\{1,2,3, ...\}$ \\
$ \mathbb{Z}_{+}$ &set of nonnegative natural numbers $\{0,1,2, \ldots\}$\\
$ \mathbb{R}_{+}$ &set of nonnegative real numbers $[0, \infty)$\\
$ \# S$ & cardinality of some set $S$\\
$(x)_{+}$ & positive part of a real number $x$, $\max(x,0)$\\
$x \vee y$ & maximum of two real numbers $x$ and $y$, $\max(x,y)$\\
$x \wedge y$ & minimum of two real numbers $x$ and $y$, $\min(x,y)$\\
 $ x_{n} = O(y_{n})$ & the real sequence $x_{n}/y_{n}$ is bounded as $n \to \infty$\\
  $ x_{n} = o(y_{n})$ & the real sequence $x_{n}/y_{n}$ tends to $0$ as $n \to \infty$\\
    $ x_{n} \sim y_{n}$ & the real sequence $x_{n}/y_{n}$ tends to $1$ as $n \to \infty$\\
  $ x_{n} \asymp y_{n}$ & $x_{n}= O(y_{n})$ and $y_{n} = O(x_{n})$\\
  $ \P, \Q, P$& probability measures on some space \\
   $ \E$& mathematical expectation under a probability measure $\P$ or $\Q$ \\
   $ E$& mathematical expectation under a probability measure $P$ \\
$:= \mbox{ or } =:$ & definition of a mathematical object (set, function, number, ...)\\

\end{tabular}

\medskip 
 We also use
the same Landau's notation $O,o,\sim$ and $\asymp$  to compare the asymptotic behavior  functions as for of real sequences. In addition, we follow the standard convention $\int_a^b=\int_{(a,b]}$.
\\

 To facilitate the navigation through the text, we next  list  some more specific notation that are used across several sections.
\eject

 \paragraph{Chapter \ref{chap:topology}}
 
 \ 
 \medskip 
 
 \begin{tabular}{cl}
 $(T,d_T,\rho)$ & real tree with root $\rho$ and distance $d_T$, often simply denoted by $T$ \\
 $\mathtt{T}=(T, d_T,\rho,g)$ & decorated tree where $g : T \to \mathbb{R}_+$ is an usc decoration\\
   $\mathbf{T}$ & measured decorated tree, also denoted by $(\mathtt{T},\nu)$ or $(T, d_T,\rho,g, \nu)$\\
    $\uplambda_T$ & length (or Lebesgue) measure on the skeleton of $T$\\
 $ \llbracket x,y \rrbracket $ & segment between $x,y \in T$\\
 $T_x$& fringe subtree above the point $x \in T$\\
 $ \partial T$ & leaves (points of degree $1$) of $T$\\
 $\U=\bigcup_{n\geq 0} \N^n$ & Ulam tree; its vertices are often called individuals and denoted by $u$ \\
 $\varnothing$ &empty sequence,  referred to as the ancestor or the root of $\U$\\
 $|u|$ & generation of $u\in \U$\\
 $u\preceq v$ or $v\succeq u$ & the vertex $u\in \U$ is a predecessor (or a prefix) of $v\in \U$\\
 $u\prec v$ or $v \succ u$ & $u\prec v$ and $u\neq v$ \\
 $(f_u)_{u \in \mathbb{U}}, (t_u)_{u \in \mathbb{U}^*}$ &building blocks for the construction of $ \mathtt{T}$ by gluing hypographs, where \\ 
  & $f_u$ gives the decoration along the segment labelled by $u$ and \\
  &$t_u$ specifies the location of a gluing point on that segment\\
  $z_u$ & length of the segment labelled by $u$, i.e. lifetime of the individual $u$ \\
  $ \mathrm{Gluing}$ & gluing operator \\
  $\dagger$ &designates a fictitious element that can be discarted \\
  &(needed for some formalism only)\\
  $\uprho(u)$, $\uprho(u,t)$ & location on $T$ corresponding to  individual $u$ at birth, resp. at age $t$\\
$T^{n}$ & tree obtained by gluing the first $n$ generations only in Theorem \ref{T:recolinfty}\\
$\mathtt{T}_u$ & decorated subtree stemming from an individual $u\in \U$\\
$  m_u$ &  $\upmu$-mass of the subtree $T_u$\\
 $\mathrm{Hyp}(g)$ &hypograph of an usc decoration $g:K\to \R_+$\\
  $\d_{\mathrm {Haus}}$ & Hausdorff distance between compact subspaces of a metric space\\
  $ \d_{\mathrm {Prok}}$ & Prokhorov distance between finite measures on a metric space\\
  $ \d_{\mathrm {Hyp}}$ & hypograph distance between two decoration functions \\
   $ \mathbb{T}_m$ & space of equivalence classes of measured decorated compact metric spaces\\
    $ \mathbb{T}$ & space of equivalence classes of (non-measured) decorated compact metric spaces\\
        $ \mathbb{T}^\bullet, \mathbb{T}^{\bullet I}$ & space equivalence classes of pointed  decorated compact metric spaces
 \end{tabular}
 
\eject
 
  \paragraph{Chapter \ref{chap:generalBP}}
 \
 \medskip 
 
  \begin{tabular}{cl}
$ \eta_u(\d t, \d x)$ & reproduction point process of an individual $u\in \U$\\
$ \chi(u)$ & type of the individual $u\in \U$, most often given by $\chi(u)=f_u(0)$\\
$(f_u, \eta_u)$ & decoration-reproduction processes (d.r.p.) of an individual $u$\\
$ P_x$ & law of the d.r.p.~$(f_u, \eta_u)$ for an individual $u$ of type $\chi(u)=x>0$\\
$ \mathbb{P}_x$ & law of the family of d.r.p.'s $(f_u, \eta_u)_{u \in \mathbb{U}}$ when $\chi( \varnothing)=x$\\
$( \mathcal{P})$ & property of the family of decoration-reproduction processes requested\\  &for the construction of a decorated tree $\texttt{T}$ in Section \ref{sec:2.1}\\
 $  \mathbb{Q}_{x}$ & law on  $\mathbb{T}$ of the equivalence class of the decorated tree $ \mathtt{T}$ under $ \mathbb{P}_{x}$\\
$ \mathcal{S}_1$ & space of non-increasing sequences in $[-\infty, \infty)$\\
$ \mathcal{S}$ & space $[-\infty, \infty) \times \mathcal{S}_1$ \\
$ \boldsymbol{ \Lambda}$ &  generalized L\'evy measure on $ \mathcal{S}$\\
$ { \Lambda_0}$ & first marginal of $ \boldsymbol{ \Lambda}$, which satisifes $\int(1\wedge y^2) \Lambda_0(\d y)< \infty$\\
$ { \boldsymbol{ \Lambda}_1}$ & second marginal of $ \boldsymbol{ \Lambda}$ on $\mathcal{S}_1$\\
$ \mathrm{k}$ &  killing rate, i.e. $ \mathrm{k}= \boldsymbol{ \Lambda}(\{-\infty\} \times \mathcal{S}_1)$\\
$( \sigma^2, \mathrm{a}, \boldsymbol{ \Lambda} ; \alpha)$ & characteristic quadruplet, where $\alpha>0$ is the scaling exponent\\
$ \psi$ & L\'evy-Khintchine exponent associated to  $( \sigma^2, \mathrm{a},{ \Lambda_0})$ via \eqref{E:LKfor}\\
$\xi$ & L\'evy process with L\'evy--Khintchine exponent $\psi$\\
$X$ & pssMp decoration constructed from $\xi$ by Lamperti transformation\\
$ \kappa$ & cumulant function defined via \eqref{E:cumulant}\\
$\uplambda^\gamma$ & weighted length measure $g^{\gamma-\alpha} \cdot \lambda_T$ defined in Proposition \ref{prop:lengthmeasures}\\
$\upmu$ & harmonic measure defined in Lemma \ref{L:MUI} under Assumption \ref{A:omega-}\\
$ (\omega_-,\omega_+)$ & interval on which $\kappa$ is negative\\
  \end{tabular}

   \paragraph{Chapter \ref{chap:example}}
 
 \ 
 \medskip 
 
  \begin{tabular}{cl}
  $ \mathrm{a_{can}}= {\mathrm a} - \int \Lambda_0(\dd y) y {\mathbf 1}_{|y|\leq 1}$ &  canonical drift when $\int (1\wedge |y|)\Lambda_0( \mathrm{d}y)<\infty$\\
finite branching activity & $ \boldsymbol{ \Lambda}_1\big( \mathcal{S}_1 \backslash \{(-\infty, -\infty, \ldots)\}\big) < \infty$\\
  non-increasing & $ \mathrm{a_{can}}  \leq 0$ and $\Lambda ( \{ (y_0, (y_i)_{i \geq 1}) : \exists y_j >0\}) =0$\\
 conservative & $\boldsymbol{\Lambda} \left( \left\{(y_{0},(y_{1}, ... )) \in \mathcal S:  \sum_{j=0}^{\infty} \e^{y_j}  \ne 1 \right\} \right) =0$\\
  fragmentation&  non-increasing conservative, $ \mathrm{a_{can}} =0$ and $ \mathrm{k} =0$\\
  binary & $\boldsymbol{\Lambda} \left( \left\{(y_{0},(y_{1}, ... ))  \in \mathcal S :  {y_3} \ne -\infty  \right\} \right) =0$\\
  growth-fragmentation &  conservative, binary and $ \mathrm{k}=0$
\end{tabular}

\eject

   \paragraph{Chapter \ref{chap:markov}}
 
 \ 
 \medskip 
 
  \begin{tabular}{cl}
  $(\tau)_{i \in I}$ & family of subtrees dangling from a base subtree\\
  $F(\varepsilon)$ &line of individuals $u\in \U$ with type $\chi(u)<\varepsilon$ such that $\chi(v)\geq \varepsilon$ for all $v\prec u$\\
$  T^{[ \varepsilon]}$ & subtree pruned at the threshold $\varepsilon$, i.e. keeping only individuals below $F(\varepsilon)$\\
$B_a(T)$& ball of radius $a$ centered at $\rho \in T$\\
  \end{tabular}

       \paragraph{Chapter \ref{chap:spinal:deco}}
 
 \ 
 \medskip 
  
  \begin{tabular}{cl}
 $ ( \sigma^{2}, \mathrm{a}_{\gamma}, \boldsymbol{ \Lambda}_{\gamma} ; \alpha)$  & $\gamma$-tilted characteristic quadruplet defined in Lemma \ref{L:LKtilde}\\
 $ ( X_{\gamma}, \eta_{\gamma})$ & decoration-reproduction process along the tagged branch\\
 $ \psi_{\gamma}(z) = \kappa(\gamma+z)$ & L\'evy--Khintchine exponent of $X_{\gamma}$\\
 $ \mathtt{T}^{\bullet} = ( \mathtt{T}, \rho^{\bullet})$ & a pointed decorated tree\\
 $ (\widetilde{\mathbb{Q}}_{x}^{\gamma})_{x >0}$& laws of pointed decorated tree defined in \eqref{eq:defspinebiaised}\\
 $\approx$ &equivalence relation for characteristic quadruplets (bifucators)\\
 $ \mathrm{ord}$ & the map $ (y, (y_{i})_{i \geq 1}) \in \mathcal{S} \mapsto (y_{j}: j \geq 0)^{\downarrow} \in \mathcal{S}_{1}$\\
 $ ( \sigma^{2}, \mathrm{a}_{\ast}, \boldsymbol{ \Lambda}_{\ast} ; \alpha)$ & the locally largest bifurcator\\ 

  \end{tabular}

   \paragraph{Chapter \ref{chap:6}}
 
 \ 
 \medskip 
 
  \begin{tabular}{cl}
  $\Z(n)=(Z_j(n))_{j\geq 1}$ &Galton--Watson process with types $j\in \N$ at generation $n\geq 0$.\\
  $T_{\GW}=(T_{\GW}, d_{T_{\GW}},\rho)$ & genealogical tree of a Galton--Watson process\\
  $\texttt{T}_{\GW}=(T_{\GW},g_{\GW})$ & decorated genealogical tree of a Galton--Watson process\\
  $\varsigma(\nu)=(v_0,v)$ & selection rule applied to the measure $\nu$ encoding a progeny\\
  $\Q^{\GW}_j$ &distribution of (the equivalence class of) the decorated tree $\texttt{T}_{\GW}$\\ 
   $P^{\GW}_j$ & distribution of the decoration-reproduction process \\  
  $\mathbb{P}^\GW_j$ & law of the family of d.r.p.~$(f_{u}, \eta_{u})_{u \in \mathbb{U}}$ \\
  & with ancestral individual of type $j \geq 1$\\
  $\boldsymbol{m}_{\GW}$ & mean matrix of the Galton--Watson process\\

  $(\boldsymbol{\pi}^{\GW}_j)_{j\geq 1}$ & reproduction kernel of the Galton--Watson process\\
   $\boldsymbol{G}_{\GW}$ &generator operator for the Galton--Watson process\\
   $\varpi: \N\to \R_+$ & weight function\\
    $l^{ \varpi}_\GW$  &mean mass of $T_{\GW}$ for the weighted length measure $\uplambda^{\varpi}_{\GW}$\\
   $\uplambda^{\varpi}_{\GW}= \varpi\circ g_{\GW} \cdot \uplambda_{\GW}$ &weighted length measure for the decorated tree $\texttt{T}_{\GW}$\\
    $\mathbf{T}_{\GW}= (\texttt{T}_{\GW}, \uplambda^{\varpi}_{\GW})$ & decorated Galton--Watson tree with a weighted length measure
    \end{tabular}

    \medskip 
\eject

     Rescalings under \eqref{eq:linmap}\\
\noindent      \begin{tabular}{cl}
    ${T}^{(n)}_{\GW}=(T_{\GW}, n^{-\alpha}d_{T_{\GW}},\rho)$ & rescaled version of ${T}_{\GW}$\\
 $\texttt{T}^{(n)}_{\GW}=(T^{(n)}_{\GW}, g_{\GW}^{(n)})$ & rescaled version of $\texttt{T}_{\GW}$ with  rescaled decoration   $g_{\GW}^{(n)}= n^{-1}g_{\GW}$  \\
 $( f^{(n)}_u, \eta^{(n)}_u)$ & rescaled version of the d.r.p.~$(f_{u}, \eta_{u})$\\
 $P_{\GW}^{(n)}$ &distribution of the rescaled decoration-reproduction process\\
 $ \P^{(n)}_{\GW}$ & law of the family of decoration-reproduction processes $(f^{(n)}_u, \eta^{(n)}_u)_{u \in \mathbb{U}}$\\
    $\Q_{\GW}^{(n)}$ &distribution of $\texttt{T}^{(n)}_{\GW}$ \\
    &  with ancestral individual of type $n \geq 1$\\

    $\varpi^{(n)}=\varpi(\lfloor n\cdot \rfloor)/\varpi(n)$ & rescaled weight function\\
\end{tabular}

\
\medskip
  
   \paragraph{Chapter \ref{C:scaling}}
 \ 
 \medskip 
 
  \begin{tabular}{cl}
  $\mathbb{D}^{\dagger}$ & space of rcll functions $\omega: [0,z]\to \R_+$ with arbitrary lifetime $z$, \\
  &equipped with Skorokhod's topology\\
   $\boldsymbol{\Lambda}^{(n)}$ & discrete generalized L\'evy measure\\
 $\kappa^{(n)}$ & discrete cumulant\\
\end{tabular}

\chapter*{Introduction}\label{chap:intro}
\addcontentsline{toc}{chapter}{Introduction}

Since the early 1990s and the introduction of the ubiquitous Brownian Continuum Random Tree (CRT) by Aldous \cite{Ald91a}, random real trees have become central objects in probability theory. Apart from their obvious applications to scaling limits of population models, they emerge in a variety of areas ranging from superprocesses, combinatorial optimization (minimal spanning trees), and analysis of algorithms to Liouville Quantum Gravity and the construction of the Brownian sphere \cite{LeGallICM}. We refer to \cite{Eva08,DLG02,LG05} for some general literature on the subject. In particular, two important classes of random real trees have been known for a long time: the Lévy trees of Duquesne--Le Gall--Le Jan \cite{Du03,LGLJ98}, which are the scaling limits of discrete Bienaymé--Galton--Watson trees, and the fragmentation trees of Haas--Miermont \cite{HM12}, which are the genealogical trees underlying self-similar fragmentation processes \cite{Ber06}. The intersection of these two classes consists of the one-parameter family of stable trees, which generalize the Brownian CRT and have proved to be crucial objects in random geometry \cite{CKlooptrees}. 

In this work, we considerably broaden the former and introduce the \textit{self-similar Markov trees}. We prove that they appear as scaling limits of many natural discrete models of trees that already popped-up in the literature, and pave the way for their systematic study. \bigskip

\begin{figure}[!h]
 \begin{center}
  \includegraphics[width=11.1cm]{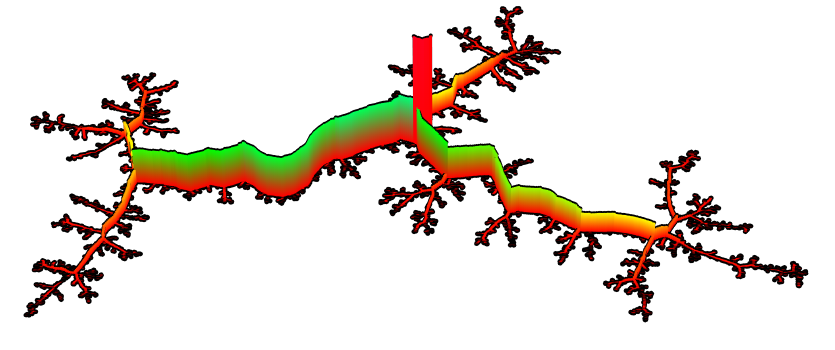}
 \caption{Illustration of a self-similar Markov tree  (embedded in the plane $ \mathbb{R}^{2}$) where its decoration function is represented in the third (vertical) dimension.} \label{fig:ssMt}
 \end{center}
 \end{figure}

\subsection*{I. Self-Similar Markov Trees} 
\begin{figure}[!h]
 \begin{center}
 \includegraphics[width=12cm]{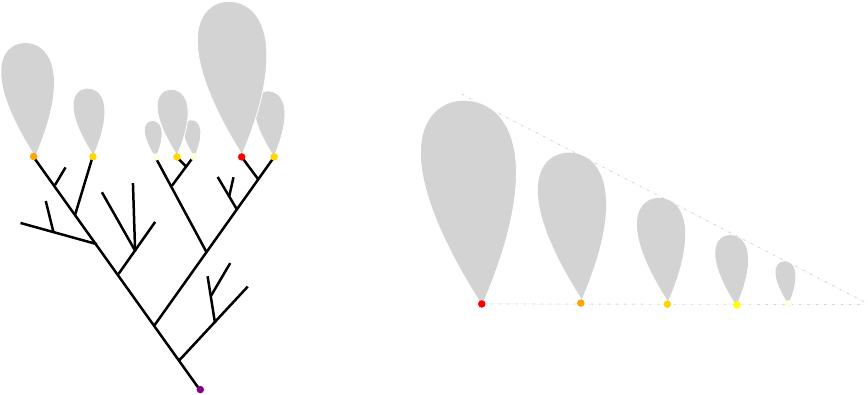}
 \caption{Illustration of the Markov property (left) and the self-similar property (right). The color of the root of a tree is meant for its decoration. Left: Conditionally on the black structure of the tree up to height $h$ and on the decoration (colors) of the root vertices, the dangling subtrees (in gray) are independent.} \label{fig:intro-prop}
 \end{center}
 \end{figure}

The first goal of this monography is to define and study the random rooted real trees $(T,d_T,\rho)$, where $\rho \in T$ is a distinguished point called the \textbf{root}, which  further support  a function $ g: T \to \mathbb{R}_{+}$, positive on its skeleton and referred to as a \textbf{decoration}. See Figure \ref{fig:ssMt} for an illustration. Informally,  a self-similar Markov tree (ssMt in short) is a family of laws $ (\mathbb{Q}_x)_{x >0}$ on the Polish space of  decorated compact random trees  (see Chapter \ref{chap:topology} for a presentation of the topology) such that  the following properties hold:
\begin{itemize}
\item \textbf{Initial decoration.} Under $ \mathbb{Q}_x$, the decoration at the root is $ g(\rho) =x$.
\item \textbf{Self-similarity}. There is a (unique) real number $\alpha >0$, called the \textbf{self-similarity index}, such that
for every $x>0$, the law  under $ \mathbb{Q}_1$ of the rescaled version $( T, x^\alpha \cdot d_T, \rho, x \cdot g)$ is  $ \mathbb{Q}_x$.
\item \textbf{Markov property}. For any height $h \geq 0$,  conditionally on the  subtree $\{ u \in T: d_T(\rho,u) \leq h\}$ up to height $h$ and on its decoration, the decorated subtrees above height $h$  are independent and each has the law $ \mathbb{Q}_{y}$, where $y$ is the decoration at the root of the subtree (see Chapter~\ref{chap:markov} for proper statements and Figure \ref{fig:intro-prop} for an illustration).
\end{itemize}

A related Markov property (without decoration) was employed in  \cite{weill2007regenerative} to characterize L\'evy trees.  Of course, random real trees carrying functions have been considered before in the literature. One may think  e.g. of the genealogical trees underlying superprocesses \cite{li2011measure} as constructed in \cite{DLG02}. One important difference though, is that usually for superprocesses, the spatial displacements is merely superposed on the branching structure, whereas in our case, the branching and the decoration are intimately tighten. Still, we borrow a lot from this theory, in particular in our rigorous treatment of the Markov property. When the tree is merely a segment, the decoration is given by a \textit{positive self-similar Markov process} (pssMp). The latter have been studied in the pioneer work of Lamperti \cite{Lam72}, who notably identified positive self-similar Markov processes as being  time changed of  exponential of  L\'evy processes, see \cite[Chapter 5]{kyprianou2022stable} and our Section \ref{sec:2.2}. It follows from the well-known L\'evy--Khintchine--It\^{o} decomposition of L\'evy processes that the distribution of  pssMp is hence determined by a self-similarity exponent $\alpha  \in \mathbb{R}$  together with the characteristic triplet $(\sigma^2 \geq 0, \mathrm a \in \mathbb{R},  \Lambda)$ of the underlying L\'evy process $\xi$. More precisely,  $ \Lambda$ is the so-called L\'evy measure on $ \mathbb{R} \cup \{-\infty\}$, it integrates $1 \wedge x^2$ and its mass at $-\infty$ denoted by  $ \mathrm{k}:=\Lambda(\{-\infty\})$, serves as killing rate. The drift depends on the choice of cutoff in the L\'evy--Khintchine formula, and in this work we use the standard cutoff $y\mapsto y \mathbf{1}_{|y|\leq 1} $. We have 
$$ \mathbb{E}\big(\exp( \gamma \xi_t)\big) = \exp\big(t \psi(\gamma)\big),$$ where 
 \begin{equation*}
 \psi(\gamma) \coloneqq- {\mathrm k} +\frac{1}{2}\sigma^2 \gamma^2+ {\mathrm a} \gamma + \int_{\R^*}   \left( \mathrm{e}^{\gamma y} -1-  \gamma y \mathbf{1}_{|y|\leq 1} \right)\Lambda ( \d y).
 \end{equation*}
 Furthermore, to ensure that the pssMp gets absorbed at $0$ in finite time a.s., $\alpha$ must be positive and the L\'evy process $\xi$ must either drift to $-\infty$ or be killed.

The first part of this work can be seen as the branching analog of the seminal contribution of Lamperti that we just sketched. We provide a rigorous setting (in particular a topology on the space  decorated real trees, see  \Cref{chap:topology})  and construct what we believe are essentially the most general (positive) self-similar Markov trees. In a nutshell, we shall see a decorated tree $(T,g)$ as a closed subset $ \mathrm{Hyp}(g)$ of the space $T \times \mathbb{R}_{+}$ where  the base space is the tree $T$, via its hypograph
$$  \mathrm{Hyp}(g) \coloneqq \big\{(u,x) \in T \times \mathbb{R}_{+} : x \leq g(u)\big\}.$$
To ensure compactness of the latter, we shall impose that $T$ is compact and $g$ is upper-semi continuous. See the numerous illustrations below for a visualization of  this concept.  As for pssMp,  ssMt are characterized by their self-similarity index $\alpha >0$  together with a Gaussian coefficient $\sigma^2\geq 0$, a drift coefficent $\mathrm a \in \R$,  and what we call a \textbf{generalized L\'evy measure} $\boldsymbol{\Lambda}$ which is a measure  on the product space $\mathcal{S}=[-\infty,\infty)\times \mathcal{S}_1$, where $\mathcal{S}_1$ is the space 
of non-increasing sequences $$ \mathcal{S}_1\coloneqq \big\{  \mathbf{y} = ( y_i)_{ i \geq 1} :  y_1 \geq y_2 \geq ...  \in \mathbb{R}\cup \{ -\infty\}\big\}.$$  We also require that the image measure $ \Lambda_0$ of $ \boldsymbol{ \Lambda}$ by the first projection $ (y, \mathbf{y})\to y$ is  a usual L\'evy measure, that is, integrates $1 \wedge y^2$. We then refer to $(\sigma^2, \mathrm a, \boldsymbol{\Lambda}; \alpha)$ as a \textbf{characteristic quadruplet}. 
Roughly speaking, we extend the idea of Lamperti and construct ssMt through a time change of the exponential of \textbf{branching Lévy processes}, which have been introduced and studied recently by Bertoin--Mallein \cite{bertoin2019infinitely}. Heuristically, the evolution of the decoration along distinguished branches of the tree are pssMp given as the time-change exponential of the L\'evy process $\xi$ with characteristics $(\sigma^{2}, \mathrm{a}, \Lambda_{0})$. Let us give an informal interpretation of the generalized L\'evy measure. Recall first that in the pssMp case, if the process is at state $x>0$, then it jumps to state $ x\cdot \mathrm{e}^{y}$ with a rate $x^{\alpha} \cdot \Lambda_{0}( \mathrm{d}y)$. If one then imagines a ssMt as the genealogical tree of a cloud of independent particles, then each particle of mass $x>0$ becomes a particle of  $ x\cdot \mathrm{e}^{y} $ and in the same time gives rise to a cloud of new particles of mass $x\cdot (\mathrm{e}^{y_{1}}, \mathrm{e}^{{y_{2}}}, ...) $ at a rate $x^{\alpha} \cdot \boldsymbol{ \Lambda}( \mathrm{d}y,  \mathrm{d}(y_{i})_{i \geq 1})$. 

Many properties of ssMt are encapsulated  by the so-called \textbf{cumulant function} defined by
$$ \kappa(\gamma) \coloneqq   \psi(\gamma) + \int_{\mathcal{S}} \boldsymbol{\Lambda}( \d y, \d  \mathbf y ) \left(\sum_{i=1}^{\infty} \e^{\gamma y_i} \right)\nonumber,$$
which can also be seen as the Biggins transform or moment generating function of the underlying branching L\'evy process. To ensure non-explosion, we assume that $\kappa$ takes negative values, which will sometimes be referred to as sub-criticality, and we leave open the construction of ssMt in the critical case $\min \kappa =0$, see Section \ref{sec:commentsGBP}.

\begin{result}[Construction of self-similar Markov trees] For every  characteristic quadruplet $(\sigma^2, \mathrm a, \boldsymbol{\Lambda}; \alpha)$ such that $\kappa(\gamma) < 0$ for some $\gamma >0$, there exists a family of laws $( \mathbb{Q}_{x})_{x >0}$ of decorated random trees $ \mathtt{T}=(T, d_T,\rho, g)$ which are self-similar with exponent $\alpha>0$ and fulfill the Markov property. 
\end{result}

As already mentioned, the decoration $g$ along branches of $T$ evolves according to a pssMp with characteristics $(\sigma^2,  \mathrm a, \Lambda_{0}; \alpha)$, and the generalized L\'evy measure $ \boldsymbol{ \Lambda}$ induces a way to ``explore" branches within $T$. However, contrary to the case of pssMp, different  characteristic quadruplets can produce the same ssMt, and we identify precisely when this happens in  \Cref{sec:bifurcators}, using a concept of bifurcators that is adapted from \cite{pitman2015regenerative,shi2017growth}.  

\subsection*{Properties of ssMt and their random measures}  Regarding the Markov property, we will in fact describe various Markov decompositions  in  \Cref{chap:markov}. We then establish several basic properties of ssMt, including the computation of their Hausdorff dimension, by studying  natural \textbf{finite measures} on ssMt.  Real trees are naturally equipped with the length measure, i.e. the $1$-dimensional Hausdorff measure, which may be thought of as the Lebesgue measure on $T$ and is therefore denoted by $\uplambda$. Although  $\uplambda$ is not even locally finite  in most cases of interest, the decoration function $g$ enables us to circumvent this issue. We consider the measures $  \mathrm{d}\uplambda^{\gamma} := g^{\gamma-\alpha} \mathrm{d} \lambda$ supported by the skeleton of $T$ and which we call \textbf{weighted length measures}. We  show that the latter has a finite total mass provided $\kappa$ takes negative values before $\gamma$, see \Cref{prop:lengthmeasures}. In particular, by self-similarity, we have 
$$ \uplambda^{\gamma}(T) = \int_{T} \mathrm{d}\uplambda \ g^{\gamma-\alpha} \quad \mbox{ under } \mathbb{Q}_{x} \qquad  \overset{(d)}{=} \qquad x^{\gamma}  \int_{T} \mathrm{d}\uplambda \  g^{\gamma-\alpha}  = x^{\gamma} \uplambda^{\gamma}(T)\quad \mbox{ under } \mathbb{Q}_{1},$$ so that the weighted length measures $\uplambda^\gamma$ are self-similar with exponent $\gamma$. In most situations, there is yet another natural  finite measure on $T$, denoted by $\upmu$, which is supported by the leaves of $T$. Its construction requires a more stringent hypothesis on $\kappa$, that we call the first Cram\'er hypothesis. This requires the existence of  $\omega_{-} \in (0,\infty)$ so that $$\kappa(\omega_{-})=0 \quad \text{and}\quad  \kappa(q)<0\text{ for some }q>\omega_-.$$ The formal slightly more restrictive condition is given by Assumption \ref{A:omega-}.

\begin{result}[Harmonic measure] Suppose $(\sigma^2, \mathrm a, \boldsymbol{\Lambda}; \alpha)$  satisfies the first Cram\'er hypothesis. Then, as $\gamma \to \omega_-$ the renormalized length measure $- \kappa(\gamma) \cdot \uplambda^\gamma$ converge in probability towards a measure $\upmu$ supported by the leaves of $T$.
\end{result}

The measure $\upmu$ is called the \textbf{harmonic measure} since it is connected to the so-called harmonic or additive martingale in the branching random walk underlying our construction of $ \mathtt{T}=(T, d_T, \rho, g)$. Alike the weighted  length measures, the harmonic measure is self-similar with exponent $\omega_{-}$. The harmonic measure is natural in many respects and it is for example used as a Frostman measure to compute the Hausdorff dimension of $T$:
\begin{result}[Hausdorff dimension] Suppose $(\sigma^2, \mathrm a, \boldsymbol{\Lambda}; \alpha)$  satisfies the first Cram\'er hypothesis. Then almost surely, $T$ has Hausdorff dimension 
$$ \mathrm{dim}_{H}( T) = \left(1 \vee \frac{\omega_{-}}{\alpha}\right), \qquad a.s.$$
\end{result}
This result considerably extends the fragmentation case \cite{HM04} or the growth-fragmentation case \cite{rembart2018recursive}.

 We then use those finite measures to provide \textbf{spine decompositions} of our decorated random trees. Roughly speaking, we use $-\kappa(\gamma)\cdot \uplambda^{\gamma}$, for $\gamma$ such that $\kappa(\gamma)<0$, or $\upmu$, which can  thought as  the extremal case $\gamma = \omega_{-}$, 
  to distinguish a point $\rho^{\bullet}$ at random in $T$. The branch $\llbracket \rho, \rho^{\bullet} \rrbracket$ connecting $\rho$ and $\rho^{\bullet}$ is then called the spine. Using the Markov property, we shall see that the decorated trees dangling from the spine are, conditionally on their initial decoration, independent ssMt with characteristics $(\sigma^{2}, \mathrm{a}, \boldsymbol{ \Lambda} ; \alpha)$. However, the evolution of the decoration along the spine is now governed by another set of characteristics $(\sigma^{2}, \mathrm{a}_{\gamma}, \boldsymbol{ \Lambda}_{\gamma} ; \alpha)$ explicitly given in terms of  $(\sigma^{2}, \mathrm{a}, \boldsymbol{ \Lambda} ; \alpha)$, and in particular, the L\'evy Khintchine exponent $\psi_{\gamma}$ of the  L\'evy process underlying the pssMp evolution along the tagged branch is simply given by 
  $$ \psi_{\gamma}(z) \coloneqq \kappa( \gamma + z).$$
  As the reader may know, spinal decompositions are essential tools in branching process theory and are also instrumental in many of our proofs. The spinal decomposition can also be seen as a more intrinsic and geometric description of the law of the ssMt. Indeed, as we alluded to above, several characteristic quadruplets $( \sigma^2, \mathrm{a}, \boldsymbol{ \Lambda} ; \alpha)$ can yield to the same  ssMt. However we shall see in Corollary \ref{C:tiltedunique} that the quadruplet  $(\sigma^{2}, \mathrm{a}_{\gamma}, \boldsymbol{ \Lambda}_{\gamma} ; \alpha)$ generically uniquely specifies the law $( \mathbb{Q}_x)_{x>0}$ of the ssMt. Furthermore, although the pssMp associated to $( \sigma^2, \mathrm{a}, \boldsymbol{ \Lambda} ; \alpha)$ is quite arbitrary, those appearing as the decoration along the tagged branch of a ssMt have special properties, see Section \ref{sec:spinalex} for details.
  
  \subsection*{II. Examples}
  
\begin{figure}[!h]
   \begin{center}
      \includegraphics[height=5cm]{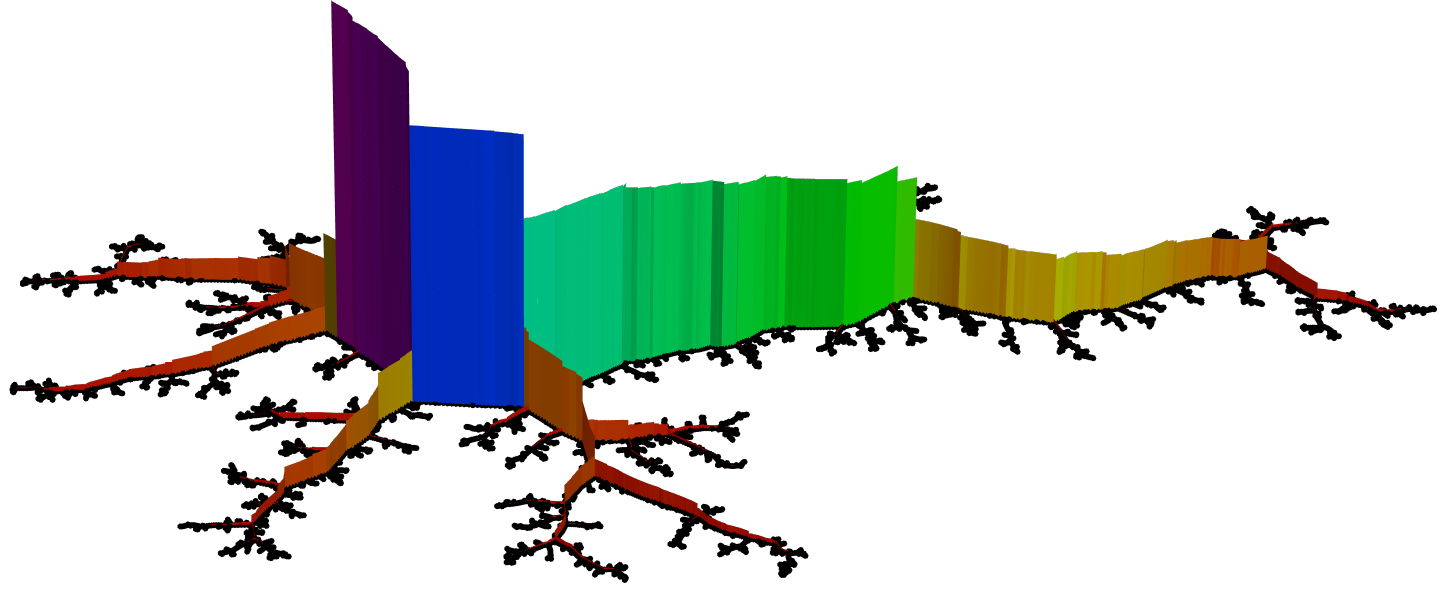}
   \caption{A simulation of a Brownian CRT. The tree is embedded (non-isometrically in $ \mathbb{R}^{2}$) and the decoration function representing the $\upmu$-mass above each point is depicted in the vertical coordinate. \label{fig:massbrointro}}
   \end{center}
   \end{figure}

  Self-similar Markov trees include many examples of important random (undecorated) real trees that have appeared in the literature. They encompass in particular the fragmentations trees constructed by Haas and Miermont \cite{HM04} as the genealogical trees underlying self-similar conservative fragmentations processes \cite{Ber02}, and the generalized fragmentation trees introduced recently by Stephenson  \cite{stephenson2013general} to cover the dissipative case. If one considers ssMt as branching analogs of L\'evy processes, then the (generalized) fragmentation trees would correspond to the subordinator case. More precisely, they consist of ssMt for which the decoration along branches is decreasing. The fragmentation trees of Haas and Miermont had ``no erosion'' and were conservative: the mass of particles is conserved in splitting events which in our case means that 
  $$\boldsymbol{\Lambda} \Big( \Big\{\big(y_{0},(y_{1}, ... )\big) \in \mathcal{S}:  \sum_{j=0}^{\infty} \e^{y_j}   \ne 1 \Big\} \Big) =0.$$
  In such situations, it is plain that the first Cram\'er hypothesis is always satisfied with $\omega_{-}=1$, and in fact the decoration $g(u)$ corresponds to the harmonic mass of the fringe subtree above point $u \in T$. In particular, we have $\upmu(T) =x$ under $ \mathbb{Q}_{x}$, so that the total harmonic mass is deterministic. Fragmentation trees already include many interesting examples such as the Brownian Continuum Random Tree (CRT) or more generally  the stable trees. Specifically, it is well known that the \textbf{ Brownian CRT} $ \mathcal{T}_{1}$ can be seen as the real tree constructed from a standard Brownian excursion of length $1$, say  $(\e_{1}(s))_{0\leq s \leq 1}$ as in \cite{DLG05,LG06}. We  denote its contour measure  by $\gamma_{ \e_{1}}$ and endow $ \mathcal{T}_{1}$ with the decoration which assigns to each vertex $ v \in \mathcal{T}_{1}$ the contour-mass $\gamma_{\e_{1}}( \mathcal{T}_{1,v})$ of the fringe-subtree above point $v$, see Figure \ref{fig:massbrointro}. 

  It follows then from \cite{Ber02,uribe2009falling} that the decorated Brownian CRT is indeed a self-similar Markov tree with self-similar with index $\alpha=1/2$, no erosion, no Gaussian part and generalized L\'evy measure  given by 
  \begin{eqnarray} \label{eq:GLMbrownianfrag} \int_{ \mathcal{S}} F\big(   \mathrm{e}^{y_{0}}, (\mathrm{e}^{y_{1}}, \ldots )\big)  \  \boldsymbol{ \Lambda}_{ \mathrm{Bro}}( \mathrm{d}y_{0},\dd  ( y_{i})_{i \geq 1}) \coloneqq    \sqrt{\frac{2}{\pi}} \int_{1/2}^{1} F(x,1-x,0,0, ...) \frac{ \mathrm{d}x}{(x(1-x))^{3/2}},  \end{eqnarray}
where $F$ stands for a generic nonnegative functional on $\mathcal S$. See Example \ref{ex:brownian}. This well-known and useful interpretation has already been used many times in the literature see e.g. \cite{borga2023power}. 

\begin{figure}[!h]
   \begin{center}
      \includegraphics[height=6cm]{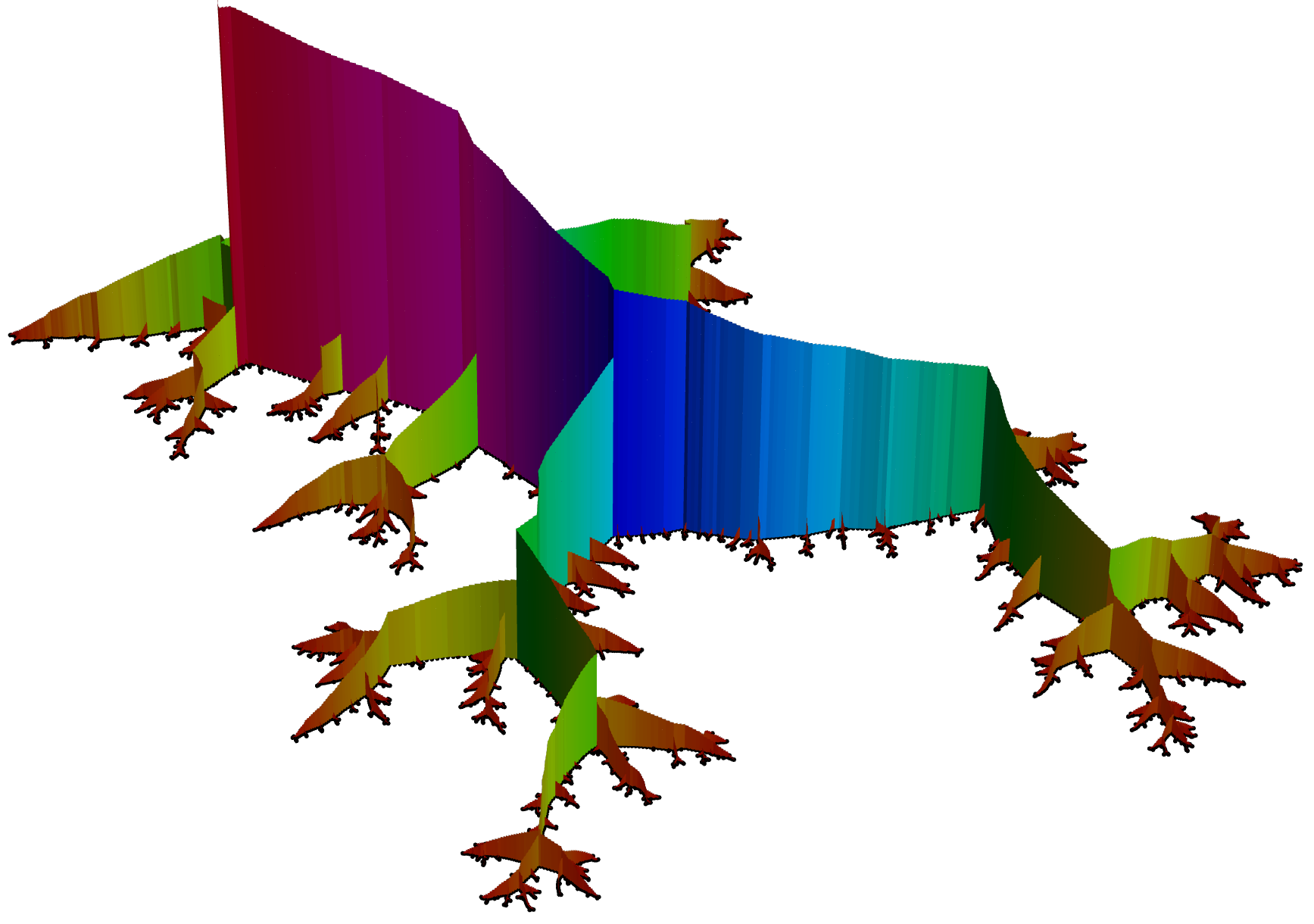}
   \caption{A simulation of a Brownian CRT normalized by the height. The tree is embedded non-isometrically in $ \mathbb{R}^{2}$; the decoration function represents the height of fringe subtrees and is depicted in the vertical coordinate. \label{fig:Dinointro}}  
   \end{center}
   \end{figure}

\begin{figure}[!h]
 \begin{center}
 \includegraphics[width=12cm]{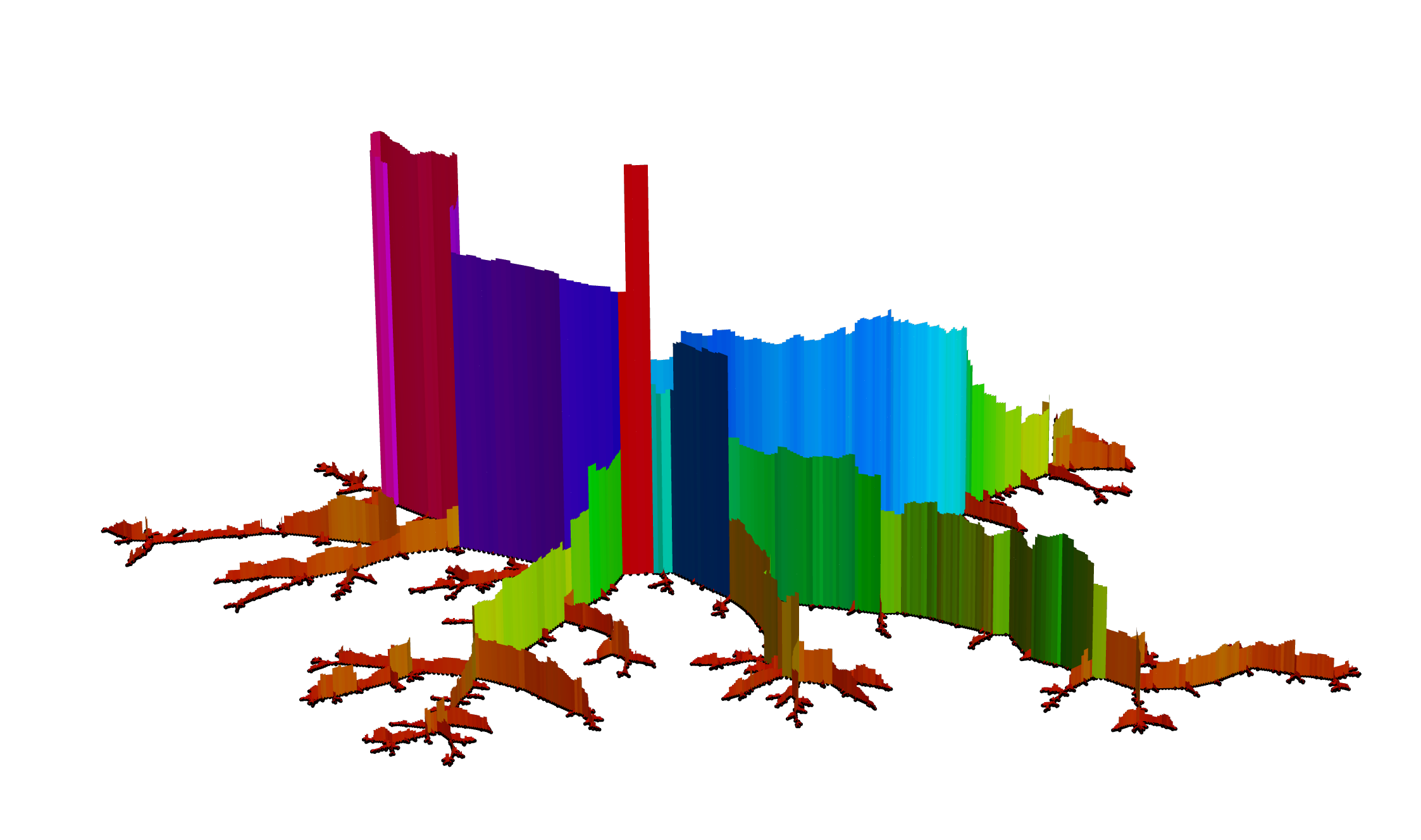}
 \caption{The decorated random tree $ \mathcal{T}_{\e}$ associated with a half-planar Brownian excursion. The tree is embedded non-isometrically in $ \mathbb{R}^{2}$; the decoration function represents the horizontal $X$ displacement in fringe subtrees and is depicted in the vertical coordinate. \label{fig:adsintro}}
 \end{center} 
 \end{figure}

Perhaps more surprisingly, there are other representations of (a small variant of the) Brownian CRT as a self-similar Markov tree which is not anymore a fragmentation tree.
 Consider this time the Brownian CRT $ \mathcal{T}^{(1)}$ \textit{of height $1$}, i.e. the tree $ \mathcal{T}_{\e^{(1)}}$ coded by an Brownian excursion $\e^{(1)}$ of height $1$.
We can endow $\mathcal {T}^{(1)}$ with the deterministic decoration which assigns to each vertex $v\in \mathcal {T}^{(1)}$  the \textit{height} of the fringe-subtree ${\mathcal T}^{(1)}_v$ rooted at $v$.  See Figure \ref{fig:Dinointro} for an illustration.
  It follows then by a classical decomposition of $\e^{(1)}$ due to  David  Williams, that the decorated tree $ \mathcal{T}^{(1)}$ is indeed a self-similar Markov tree with index $\alpha=1$,  no Brownian part, constant erosion and generalized L\'evy measure given by 
 $$\int_{\mathcal S} F\big(   \mathrm{e}^{y_{0}}, (\mathrm{e}^{y_{1}}, \ldots )\big) \boldsymbol{\Lambda}_{\mathrm{Height}}(\dd y, \dd\mathbf y)= 2\int_{0}^1 F(1,(x,0,0, \ldots))  \frac{ \dd x}{x^{2}},$$
 where $F$ denotes  a generic nonnegative functional on $\mathcal S$.  The cumulant function is easy to compute and equals $\kappa_{\mathrm{Height}}(\gamma) = -\gamma + 2/(\gamma-1)$ for  $\gamma>1$. In particular, Cram\'er  assumption holds with $\omega_{-}=2$. As the reader may expect, the harmonic measure $\upmu$ then coincides with 
 (a multiple of) the contour measure $\gamma_{\e^{(1)}}$ on $\mathcal T^{(1)}$, and in particular its total mass is now random. See Example \ref{ex:brownianheight}.  These are not the only representations of -- variations of --  the Brownian CRT 
that can be obtained using ssMt. In particular, in Example \ref{ex:ADS}, we present another representation based on the recent work of Aidekon and Da Silva \cite{aidekon2022growth}, which connects growth-fragmentations to planar Brownian excursions. In this example, the Brownian CRT is, roughly speaking, seen as a ssMt  where the decoration processes along branches are variants of the symmetric Cauchy process, see Figure \ref{fig:adsintro}.  We refer to Example \ref{ex:ADS} for  details.

\begin{figure}[!h]
  \begin{center}
  \includegraphics[width=12cm]{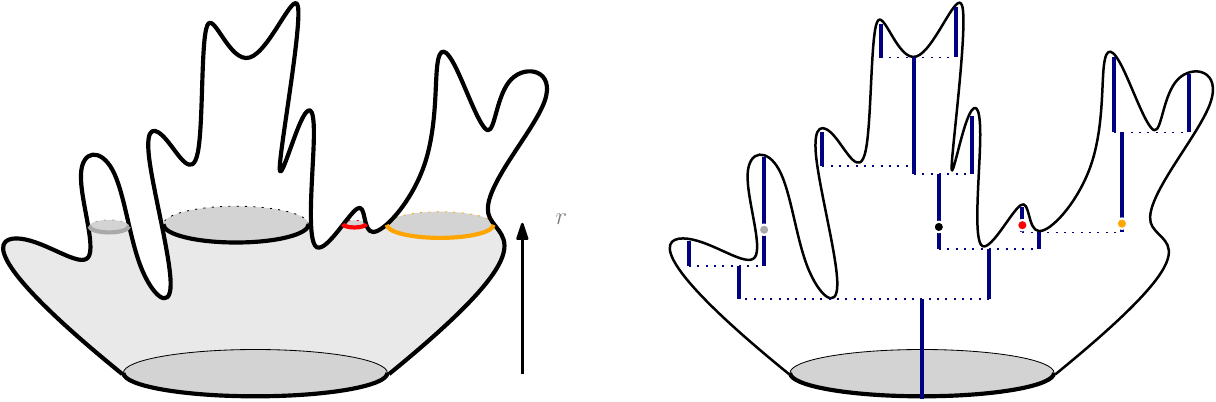}
  \caption{Illustration of the Cactus of a surface $S$ with a boundary $ \partial S$. The ball of radius $r$ (measured from $ \partial S$) is depicted in light gray and it has several boundary components. When each of these components has a ``size'', this enables us to decorated the cactus tree (on the right). \label{fig:cactusintro}}
  \end{center}
  \end{figure}

Luckily, there is more to life than variations on the Brownian CRT. The primary incentive for developing the general theory of ssMt is its connection with random planar geometry. Motivated in part by $2$-dimensional quantum gravity, the last few decades have witness spectacular developments in this field, notably around the Brownian sphere. This is a random compact metric space almost surely homeomorphic to the $2$-sphere \cite{LGP08} but with fractal dimension $4$ \cite{LG07}. It has first been constructed as the scaling limit of random planar quadrangulations by Le Gall \cite{LG11} and Miermont \cite{Mie11} and since then appeared as a universal scaling limit model for many planar graphs models see e.g. \cite{LG11,CLGmodif,BJM13,albenque2020scaling}. The Brownian sphere has also been shown to be the random metric induced by exponentiating a planar Gaussian Free Field (GFF) to the proper power, see the works of Miller \& Sheffield \cite{miller2015liouville,MS15} yielding to a so-called Liouville Quantum Gravity metric \cite{ding2020tightness,gwynne2020existence}. The Brownian sphere has a cousin, the Brownian disk \cite{BM15}, which, as the name suggests, has the topology of a disk and is better suited to make the connection with ssMt. Informally, the Brownian disk is a random compact metric space  $(S,d)$ homeomorphic to the closed unit disk of the complex plane, and can be obtained as the scaling limit of generic random planar maps with a large boundary \cite{BM15}. In particular, we can define its boundary  $\partial S$ as the set of all points that have no neighborhood homeomorphic to the open unit disk. We can then consider, for every $r\geq 0$, the ball $B_r = \{ x \in S : d(x, \partial S) \leq r\}$.
 The topological boundary $\partial B_r$ of these balls generically has infinitely many connected components $(C_r^i: i \geq 1)$, which all have the topology of a circle and have fractal dimension $2$. However it is possible to associate  a natural notion of ``boundary length" $(|  C_r^i|: i \geq 1)$ to them extending the classical notion of perimeter, see \cite{le2020growth}. Furthermore, as $r$ grows, the connected components $(C_r^i: i \geq 1)$ describe a tree structure , see Figure \ref{fig:cactusintro} for an illustration and Example \ref{ex:3/2stable} for details. This tree is called the \textbf{Brownian cactus} and was first  studied in \cite{CLGMcactus}. We can then decorate each point of the Brownian cactus with the associated boundary length.

  It follows  from the recent work \cite{le2020growth} (see e.g. \cite{BCK18,BBCK18}) that the resulting decorated tree is\footnote{In the case when the base surface is a free Brownian disk, see \cite{BM15}} a (multiple of the) self-similar Markov tree called \textbf{Brownian growth-fragmentation tree}. This is the ssMt with exponent $ \frac{1}{2}$, with no Brownian part and generalized L\'evy measure given by  
 \begin{eqnarray} \label{eq:GLMgrowthBrownian} \int_{\mathcal S} F\big(  \mathrm{e}^{{y_{0}}}, (\mathrm{e}^{y_{1}}, ...) \big)  \  \boldsymbol{ \Lambda}_{\mathrm{BroGF}}(\mathrm{d} y_{0}, \mathrm{d} \mathbf{y})  \quad  := \quad    \frac{3}{ 4 \sqrt{\pi}} \int_{1/2}^{1} F(x,1-x,0,0, ...) \frac{ \mathrm{d}x}{(x(1-x))^{5/2}}.  \end{eqnarray}
 Notice that the L\'evy measure  ${\Lambda}_{\mathrm{BroGF}, 0}( \mathrm{d}y)$ does not integrate $y$ in the vicinity of $0$ and so compensation is involved in the definition. In particular, the drift term is non-trivial (and explicit) and the cumulant function is equal to 
$ \kappa_{\mathrm{BroGF}}(\gamma) = \frac{\Gamma(\gamma- 3/2)}{\Gamma(\gamma-3)},$ for   $\gamma > 3/2.$ So that $\omega_-=2$ and the first Cram\'er hypothesis holds. We refer again to Example \ref{ex:3/2stable}. The same ssMt but with self-similarity exponent $\alpha = 3/2$ also appears as the scaling limits of many labeled trees arising in the combinatorial literature associated with polynomial equations with one catalytic variable \cite{bousquet2006polynomial}, we refer to the forthcoming Part II for examples. We also refer to Chapter \ref{chap:example} for other examples of ssMt related to $\alpha$-stable processes for $\alpha \in (1, 3/2]$ and which appear in random planar geometry as well.

\begin{figure}[!h]
 \begin{center}
 \includegraphics[width=8cm,angle=-3]{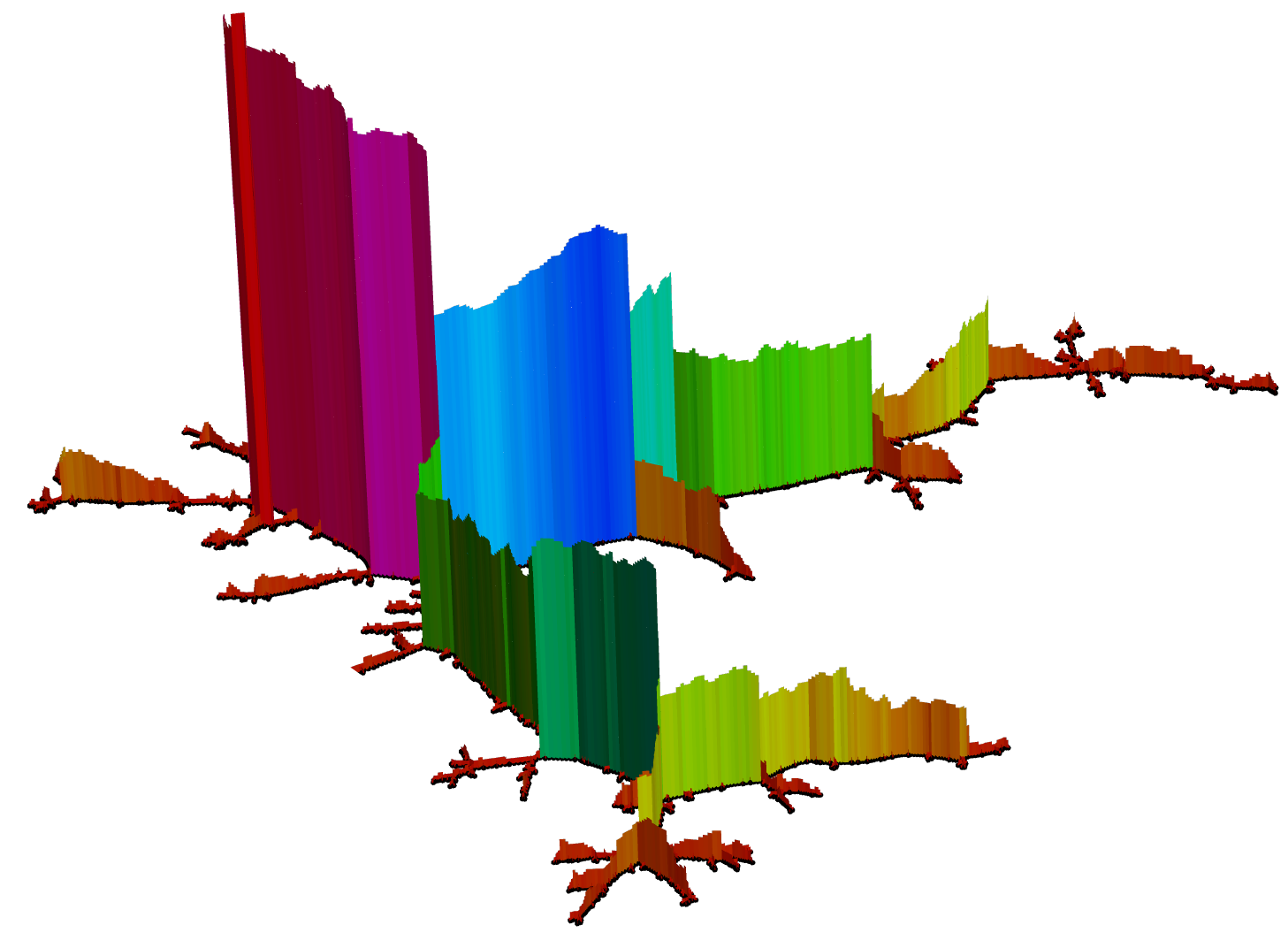}
 \caption{A simulation of the Brownian growth-fragmentation tree. The self-similar Markov tree is binary and conservative: at each splitting event, the total mass is conserved and split between two children.}
 \end{center}
 \end{figure}

\subsection*{III. Invariance principles for multi-type Galton--Watson processes} 
Any scaling limit inherently exhibits self-similarity. This fundamental insight, due to Lamperti, 
 leads to the central problem of identifying classes of processes 
 that, under appropriate rescaling, converge to a given self-similar process. For example, the Brownian CRT
 is the universal scaling limit as the size goes to infinity of Galton--Watson trees conditioned on their sizes, for arbitrary critical reproduction laws with a  finite variance.
 Our main objective in the second part of this  work is to develop robust invariance principles under which multi-type Galton--Watson trees converge towards a self-similar Markov tree. 
 More formally, suppose that we have a collection of particles which evolve as  a multi-type Galton--Watson process with types  in $ \mathbb{N} = \{1,2, ... \}$. We denote  the law of the process, starting from a single particle of type $j\geq 1$, by $  \mathbb{P}^{\GW}_{j}$ and suppose aperiodicity for simplicity. Here also we stress that the type is intimately tighten to the branching mechanism and should not be thought as a superposed spatial displacement as for many superprocesses. It is easy to interpret such a branching system  as a random decorated tree $(T_{\GW},d_{T_{\GW}}, \rho,g)$, see Figure \ref{fig:GWtoBPlight} below.

In this direction, it will be convenient to systematically  distinguish one child particle in each non-empty progeny, 
for instance the child with the largest type, and then gather the remaining children as a non-increasing sequence. Then,  a reproduction event can be represented in the form 
$j \mapsto (k_0, (k_1, \ldots, k_\ell))$, meaning that a particle of type $j$ gives birth to $\ell+1$ particles with types $k_0 \geq k_1 \geq ... \geq k_\ell$, and the first one with type $k_0$ has been distinguished. 
The reproduction law of the Galton--Watson process induces a family of (sub)probability measures $(\pi_j)_{j \geq 1}$, where 
$$\pi_j\left( (k_0, (k_1, \ldots, k_\ell)\right)$$
is the probability of the reproduction event $j \mapsto (k_0, (k_1, \ldots, k_\ell))$ for any given non-increasing finite sequence $ k_0 \geq k_1 \geq ... \geq k_\ell$ in $\N$. 
\begin{figure}[!h]
 \begin{center}
 \includegraphics[width=12cm]{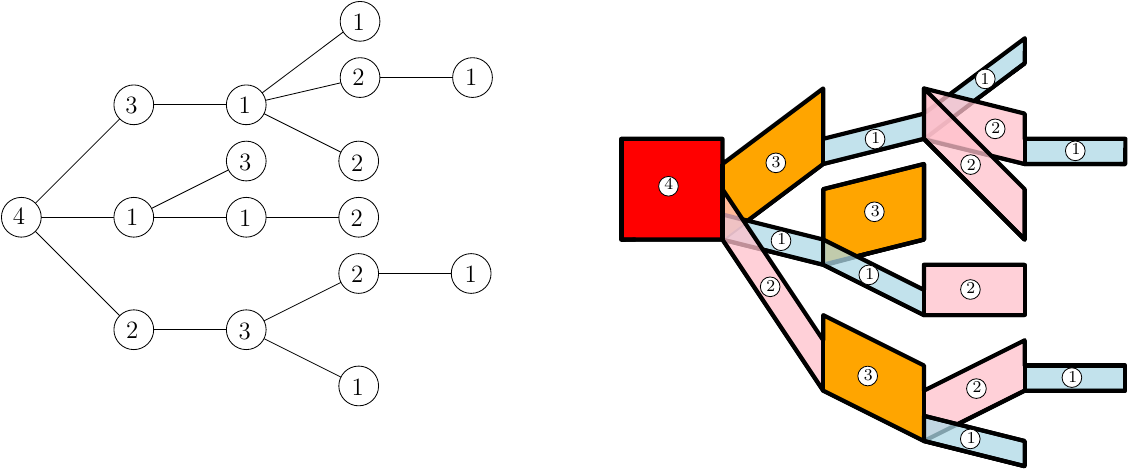}
 \caption{A representation of a Galton--Watson process with integer types as a random decorated tree, see the forthcoming Part II for more details. \label{fig:GWtoBPlight}}
 \end{center}
 \end{figure} 
 We now present the conditions on the reproduction kernel $(\pi_{j})_{ j \geq 1}$ ensuring the convergence of the rescaled discrete decorated trees towards a $(\sigma^2,  \mathrm{a}, \boldsymbol{\Lambda} ; \alpha)$-ssMt. Recall that $ \mathbb{P}_{1}$ is then the law of $  \mathtt{T}=(T,d_T, \rho,g)$ when the decoration at the root is $1$. The most obvious necessary condition is the vague convergence of the renormalized kernel towards the generalized L\'evy measure $ \boldsymbol{\Lambda}$, namely
\begin{align*}
  &  \lim_{n \to \infty}  n^{\alpha} \cdot \sum_{k_0\geq  \cdots \geq k_{\ell}}\pi_n\big(k_0, (k_1, \ldots, k_\ell)\big) f\left( \log \frac{k_0}{n}, \left( \log \frac{k_1}{n}, ... , \log \frac{k_\ell}{n} \right) \right) \\ &= \int_{ \mathcal{S}}  \boldsymbol{\Lambda}( \dd y_0, \dd \mathbf y) f( y_0,\mathbf{y}),
  \tag{${\color{red} \boldsymbol{\heartsuit}}$}
 \end{align*}
  for continuous functionals $f : \mathcal{S} \to \mathbb{R}_{+}$ with compact support avoiding $(0,(-\infty, \ldots))$. In particular, the above assumption gives a clear meaning to the generalized L\'evy measure: typically a particle of large type  $n$ gives birth to an essentially single particle of type close to $n$, but with probability of order $n^{-\alpha}\cdot  \boldsymbol{\Lambda} ( \mathrm{d} \mathbf{s}) $ it gives rise to several particles $( n \mathrm{e}^{y_0}, n \mathrm{e}^{y_1}, \ldots)$ of type comparable to $n$.  We further need a control on the drift and variance of the first coordinate, see Part II for the proper definition and denote  those assumptions by $ ({\color{red} \boldsymbol{\spadesuit}})$ in this introduction. Extending a result of Bertoin \& Kortchemski \cite{BK14}, the previous two conditions implies the convergence of rescaled decoration-reproduction processes over long branches of $T_{\GW}$ towards the  Markov decoration-reproduction process $X$ associated with the characteristic quadruplet $(\sigma^2, \mathrm{a}, \boldsymbol{\Lambda} ; \alpha)$ defined in Section \ref{sec:2.2}.  However, these two conditions are only asymptotic in the particle-type and cannot account alone for the convergence of rescaled trees and does not even guarantee that the tree $T_{\GW}$ is finite under $ \mathbb{P}_{j}^{\GW}$. The latter is ensured  by requiring the existence of a super-harmonic function for the multi-type Galton--Watson process and the integrability condition 
\begin{align*}\limsup_{n \to \infty} ~n^{\alpha} \cdot  \sum_{(k_1, ... , k_{\ell})}\pi_n(k_1, ..., k_\ell) \left(\sum_{i=1}^\ell \left(\frac{k_i}{n}\right)^q -1\right)<0, 
  \tag{${\color{red} \boldsymbol{\diamondsuit}}$}
  \end{align*}
  in some  open interval of $(0,\infty)$ for the parameter $q$.  The latter assumption is the discrete counterpart to the sub-criticality assumption $ \min \kappa <0$ in the continuous. We prove that the previous set of three assumptions is sufficient to imply the convergence of the rescaled decorated tree in the Gromov--Hausdorff decorated sense, see Part II for the proper statement.
  \begin{result}[Scaling limits without  mass measure]  Suppose (${\color{red} \boldsymbol{\heartsuit}}$),(${\color{red} \boldsymbol{\spadesuit}}$) and (${\color{red} \boldsymbol{\diamondsuit}}$). Then as $n\to \infty$, we have the following convergence for the decorated Gromov--Hausdorff topology  
$$\left(T_{\GW}, \frac{d_{T_\GW}}{n^\alpha}, \rho,   \frac{g}{n}\right) \quad \mbox{ under } \mathbb{P}_{n}^{\GW} \quad \xrightarrow[n\to\infty]{(d)} \quad (T,d_T,\rho, g) \quad \mbox{ under }\P_1.$$
\end{result}

The topology for this convergence is that developed in  \Cref{chap:topology}; it is adapted from the classical Gromov--Hausdorff topology to involve the decoration as well. Let us now consider more delicate versions of this result that incorporate measures. We fix a non-zero \textbf{weight function} $\varpi:\mathbb{N}\to \mathbb{R}_+$  regular varying with exponent $\gamma-\alpha$  and write  $\upmu^\varpi$ for the measure on $T_\GW$ that assigns a mass $\varpi(k)$ to each particle with type $k$. We  let $l^{\varpi}_{\GW}(n)$ denote the  expected mass of a multi-type tree starting from a single particle of type $n$, i.e.
$$ l^{\varpi}_{\GW}(n) := \mathbb{E}^{\GW}_{n}\big(\upmu^\varpi \big( T^{\GW}\big) \big).$$
We naturally aim for a scaling limit result of the measure decorated tree $\left(\frac{T_\GW}{n^{\alpha}}, d_{T_{{\GW}}}, \rho,  \frac{g}{n},  \frac{  \mu^{ \varpi}}{l^{\varpi}_\GW(n)}\right)$. As in the continuous case, there is a big difference between the case $\gamma > \omega_{-}$ and $\gamma \leq \omega_{-}$. More precisely, if $\gamma > \omega_{-}$, the  measure $\mu^{ \varpi}$ is essentially carried by the branches of $T_{\GW}$ and converge after renormalization towards the weighted length measure $ \uplambda^\gamma$ (jointly with the converge in the above result). However, in the case when $\gamma \leq \omega_{-}$ the measure $\mu^{ \varpi}$ is now carried by the ``leaves'' of $T_{\GW}$ and  converges towards the harmonic measure $\upmu$. Actually, proving this convergence requires the discrete counterpart of the first Cram\'er hypothesis, see the forthcoming Part II for details.  Our discrete invariance principles recover those of Haas--Miermont \cite{HM12} in the fragmentation case. This was actually one of the main source of inspiration for our convergence results. 
  In particular, the requirement  $({\color{red} \boldsymbol{\heartsuit}})$ is the analog of the fundamental hypothesis (H) of Haas--Miermont \cite{HM12}. 
But  the most interesting, by far, applications of our invariance principles concern the discrete multi-type Galton--Watson trees converging to ssMt where the decoration can exhibit growth. An important example in this direction is given by the peeling trees underlying Boltzmann stable maps already considered in \cite[Section 6]{BBCK18} (see \cite{BCK18} in the triangulation case). Let us describe that setting more precisely. Given a non-zero sequence $ \mathbf{q}:= (q_{k})_{ k \geq 1}$ of non-negative numbers we define a measure $w$ on the set of all bipartite planar maps $ \mathbf{m}$ (finite graph embedded in the sphere up to homeomorphisms, given with a distinguished oriented edge) by the formula
\begin{equation} \label{eq:defface}
w( \mathbf{m}) := \prod_{f\in \mathrm{Faces}( \mathbf{m})} q_{\deg(f)/2}.
\end{equation}
We shall suppose that the weight sequence is admissible (the above measure is finite) and critical (roughly speaking, the weights cannot be increased while staying admissible), see \cite[Chapter 5 ]{CurStFlour}. We shall furthermore suppose that $ \mathbf{q}$ is non-generic, in the sense that it satisfies
 \begin{eqnarray} \label{eq:asymptoticqk} q_{k} \sim c\, \kappa^{k-1} \, k^{-\beta} \quad \mbox{ as } k \to \infty, \qquad \mbox{ for }  \beta\in \left( \frac{3}{2}, \frac{5}{2}\right),  \end{eqnarray} for some $c, \kappa>0$.\footnote{Actually, we will establish the result for a slightly more general form of non-generic maps.} For this choice of $ \mathbf{q}$, a random $ \mathbf{q}$-Boltzmann  $ \mathbf{M}_{n}$ conditioned on having $n$ faces (provided that the conditioning make sense) possesses ``large faces'' and converge in the scaling limit 
 $$ (\mathbf{M}_{n} , n^{ -\frac{1}{2 \beta-1}} \cdot d_ {\mathrm{gr}}) \xrightarrow[n\to\infty]{(d)} \mathcal{S}_\beta,$$ in law for the Gromov--Hausdorff distance, where $\mathcal{S}_\beta$ is the $\beta$-stable carpet/gasket introduced implicitly  by Le Gall and Miermont  in \cite{LGM09} and whose uniqueness has been very recently proved in \cite{CRM22}. In particular, the limiting case $\beta = 2$ corresponds to  the Brownian sphere, see \cite{marzouk2018scaling}. On the other hand, the dual maps $ \mathbf{M}_{n}^{\dagger}$ (with large degree vertices) are much less understood\footnote{  For readers familiar with the theory of random maps, this is due to the lack of suitable bijective encodings  \`a la Schaeffer.}. Although the typical distances in $ \mathbf{M}_{n}^{\dagger}$ are known to be of order $n^{\frac{\beta-2}{\beta-1/2}}$ when $\beta >2$, see \cite{BC16}, it is not known whether the metrics spaces $( \mathbf{M}_{n}^{\dagger}, n^{-\frac{\beta-2}{\beta-1/2}}  \cdot d_{ \mathrm{gr}})$ are tight. One important tool in the theory of random planar maps is the so-called \textbf{peeling process}, which is a Markovian exploration procedure that discovers a map step-by-step, see \cite{CurStFlour}. This exploration actually encodes a planar map into a binary labeled plane  tree, which under the $ \mathbf{q}$-Boltzmann measure is a multi-type Galton--Watson tree, see the forthcoming Part II  for details. In the critical non-generic case above, we prove that those trees satisfy our standing assumptions (${\color{red} \boldsymbol{\heartsuit},\boldsymbol{\spadesuit}, \boldsymbol{\clubsuit}})$ with $\alpha = (\beta -1)$ for an explicit subcritical characteristic quadruplet $( \sigma^{2}, \mathrm{a}, \boldsymbol{ \Lambda} ; \alpha)$. As an immediate corollary of this convergence, we deduce that the diameter of $( \mathbf{M}_{n}^{\dagger}, n^{-\frac{\beta-2}{\beta-1/2}}  \cdot d_{ \mathrm{gr}})$ is tight. The corresponding ssMt should play the role of the cactus tree in the potential scaling limits of $ \mathbf{M}_n^\dagger$ and this convergence thus lays the foundations for the definition of their scaling limits when $\beta >2$.

\subsection*{Relation to previous works} As we said above, the main source of inspiration for this monography is the work of Haas--Miermont \cite{HM04,HM12}, where they constructed the self-similar fragmentation trees and established invariance principles, see \cite{haas2018scaling} for a beautiful set of lecture notes. To be more specific, the fragmentation trees of Haas--Miermont \cite{HM04}, later generalized by Stephenson \cite{stephenson2013general} correspond to the family of self-similar Markov trees where the decoration $g$ is decreasing along branches. A powerful invariance principle for discrete fragmentation was established in \cite{HM12} for the Gromov-Hausdorff topology (using quite different tools as ours) and proved to be very useful in the study of random trees \cite{BM13,Riz11}. 

The definition of growth-fragmentation processes was given by Bertoin \cite{bertoin2016compensated,Ber15} first in the binary conservative case. The properties of  growth-fragmentation processes were then studied in much details, see e.g. \cite{shi2017growth,watson2023growth,dadoun2017asymptotics,bertoin2020strong,bertoin2016local} (the list is by no mean exhaustive). It was soon realized that these processes encode the genealogy of real trees, and they were constructed by Rembart--Winkel using the ``string of beads'' construction \cite{rembart2018recursive}, which also inspired the construction of \Cref{chap:generalBP}. In the self-similar case, those (decorated) trees are indeed examples of ssMt. Growth-fragmentation processes were shown to appear in peeling exploration of random planar maps \cite{BCK18,BBCK18},  directly within Brownian geometry \cite{le2020growth} and in random plane Brownian excursions \cite{aidekon2022growth}. We revisit these results in light of the ssMt theory. Let us in particular point to the paper of Dadoun \cite{dadoun2019self} which builds upon the techniques of \cite{BCK18} to establish an invariance principles in the Gromov--Hausdorff topology for discrete (binary, conservative) self-similar Markov trees, appropriately truncated. As stated above, much of these works focus on the binary conservative case until Bertoin--Mallein layed the foundations of general branching L\'evy processes  \cite{bertoin2018biggins,bertoin2019infinitely} which  are closely related to our construction of self-similar Markov trees. 

An important caveat, is that in many of the previous works the  random mass measures and the decorations on random continuous and discrete trees was neglected (it is ``trivial'' in the fragmentation case of Haas--Miermont). Establishing invariance principle for measured decorated random trees pushed us to introduce an appropriate topology, find the appropriate assumptions and develop new and wider proofs ideas. Most of our approach relies on construction of trees via ``stick-breaking'' construction, i.e.~by the recursive gluing of (decorated) branches. This is old technique pioneered by Aldous \cite{Ald91a} and  which was revived recently \cite{CH17,rembart2018recursive,blanc2023compactness} as an alternative to the ``contour function approach'' \cite{DLG05}.\\

\noindent \textbf{Acknowledgments.} We thank most warmly Alejandro Rosales-Ortiz for his careful reading of this text and his comments. We are also grateful
 to Timothy Budd, Igor Kortchemski, and Thomas Lehéricy for stimulating discussions at various stages of this project, and  to Elie Aidekon, Yueyun Hu, Bastien Mallein, and Zhan Shi for insightful exchanges on branching random walks, particularly in the critical case. The second author is supported by SuPerGRandMa, the ERC CoG 101087572. Finally, we  acknowledge  the events \textit{Décima Escuela de Probabilidad y Procesos Estocásticos de México: 2022} ,  \textit{LMS Research School on Probability at the University of Liverpool: 2023} and \textit{ CIRM Self-similar Markov trees in random geometry
:  2024} where preliminary versions of this work were presented.

\part{\sc{Self-similar Markov trees}}
\chapter{Decorated real trees and their topologies}
\label{chap:topology}
Random decorated real trees are the central objects in this work. A decorated real tree consists of a compact real tree $T$  equipped with a root $\rho$  and nonnegative function $g$ on $T$, which we refer to as a decoration and is typically discontinuous at branching points. Often, decorated trees will further be measured, i.e. endowed with a finite Borel measure.

After recalling some basic background on real trees, we present in the first section a natural gluing operation, using marks on the base tree to specify locations where the other trees are glued. This enables us in the second section to present a recursive construction of decorated real trees from so-called building blocks, indexed by the Ulam tree, by gluing iteratively line segments on which some functions have been defined (Theorem \ref{T:recolinfty}).
\begin{figure}[!h]
 \begin{center}
  \includegraphics[width=15cm]{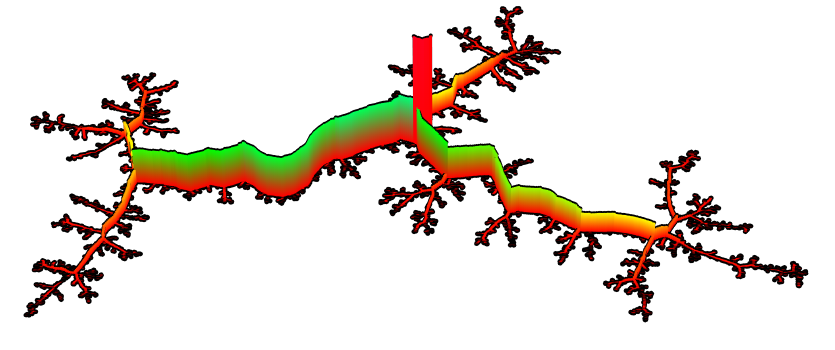}
 \caption{Hypograph of a  tree  (embedded in the plane) decorated by an usc function in the third (vertical) dimension.}\label{fig:hypogr}
 \end{center}
 \end{figure}
 
 We  will be interested in the convergence of sequences of such decorated real trees, which incites us to define a notion of closeness for two decorations  on two different domains. One of the purposes of this chapter is thus to introduce a  convenient formalism and,
to stay on safe ground as soon as randomness and measurability will be  involved, we represent these as elements of some Polish space.   
The topological framework that is needed to define a notion of convergence for sequences of 
 decorated real trees is developed in the third section. In this direction, we shall always impose that the decoration $g$ is upper semi-continuous, i.e.,  the superlevel sets $\{x:~g(x)\geq r \}$ are closed for every $r\geq 0$. The main idea there is to represent usc functions by their hypographs; see Figure~\ref{fig:hypogr}.
In the situation where the domains can be viewed as subsets of a same metric space,
the Hausdorff metric yields a natural distance in terms of  hypographs.  We then adapt an idea of Gromov and consider isometric embeddings  to measure the distance between the structures induced by two usc decorations with unrelated domains (Theorem \ref{theo:Polish}).

\section{Decorated real trees and a gluing operator}\label{sec:1.1}

We start by recalling a few basic features on real compact trees. The reader is referred to the lecture notes by Evans \cite{Eva08}, by Le Gall \cite{LG05}, or to the manuscript  by Duquesne \cite{duquesne2006coding} for detailed accounts and further properties. 
 
A \textbf{real tree} is a metric space  $(T,d_{T})$ such that for  any $x,y \in T$, there exists a unique isometry
 $$\phi_{x,y}: [0, d_{T}(x,y)] \hookrightarrow T \quad \text{ with }
  \phi_{x,y}(0)=x \text{ and }\phi_{x,y}\big(d_{T}(x,y)\big)=y;$$ 
  we shall refer to $\phi_{x,y}$ as the \textbf{path} from $x$ to $y$ in $T$.
 Furthermore,  the image of any continuous injective map $\psi:[0,1] \to T$ with $\psi(0)=x$ and $\psi(1)=y$ coincides with $\phi_{x,y}([0, {d}_{T}(x,y)])$. The   image of $\phi_{x,y}$ is called the \textbf{segment} from $x$ to $y$ in $T$ and is denoted by $\llbracket x,y \rrbracket$ in the sequel.  
 
  A useful characterization is that a connected metric space  $(T,d_{T})$ is a real tree if and only if  the so-called four point condition holds, that is:
  \begin{equation}\label{E:4points}
  d_T(x_1,x_2)+d_T(x_3,x_4)\leq (d_T(x_1,x_3)+d_T(x_2,x_4))\vee (d_T(x_1,x_4)+d_T(x_2,x_3)),
  \end{equation}
   for every $x_1, x_2, x_3, x_4$ in $T$. The real trees appearing in this work are typically \textbf{compact} and  \textbf{rooted} at some distinguished point $\rho\in T$, even though for the sake of simplicity we shall often omit the distance $d_{T}$ and the root $\rho$ from the notation. The distance $d_T(\rho,x)$  of a point $x\in T$ to the root  is called the \textbf{height} of $x$, and then the height of $T$, denoted by  $\mathrm{Height}(T)$, is the maximal height of points in $T$.
 The root enables us to endow $T$ with a partial order: for any two points $x,y\in T$ we write $x\preceq y$ and we say that $x$ is an \textbf{predecessor} of $y$, or  that $y$ is a \textbf{descendant} of $x$, if $x$ belongs to the segment $\llbracket \rho,y \rrbracket$. The \textbf{fringe subtree} $T_x=\{y\in T: x\preceq y\}$ induced by a point $x\in T$ is the subset of $T$ consisting of  all the descendants of $x$ (including $x$ itself). This fringe subtree is naturally equipped with the metric induced by the restriction of $d_T$, and rooted at $x$. 
 We also use the following standard nomenclature and notation for points in a  real tree. The \textbf{degree} of a point $x \in T$ is the number (possibly infinite) of connected components of $T \backslash \{x\}$, and then:
\begin{itemize}
\item A point $x\in T$ is a \textbf{leaf} if it has degree $1$. We denote the set of leaves of $T$ by $\partial T$.
\item The \textbf{skeleton} of $T$ is the subset $T\setminus\partial T$.
\item A point $x\in T$ is a \textbf{branching point} if it has degree at least $3$. 
\end{itemize}
  We will further equip real trees with some \textbf{upper semi-continuous} (usc in the sequel) functions\footnote{We enforce this requirement because usc-functions are well-behaved when gluing trees (see below), and also for their compatibility with Gromov-type topologies (see Section \ref{sec:1.3}).}.  
  
  \begin{definition}\label{Def:dectree}  A \textbf{decorated real tree} is a
quadruplet 
  $$ \normalfont{\texttt{T}}= ( T, {d}_{T}, \rho, g)$$
where  $(T, {d}_{T},\rho)$ is a rooted compact  real tree and   $g: T\to \R_+$ an usc function referred to as the \textbf{decoration}.
 \end{definition}
 Let us give a couple of simple examples. Consider a rooted compact real tree $( T, {d}_{T}, \rho)$. The function $g: T\to \R_+$ that assigns to every $x\in T$ 
 its height, that is $g(x)=d_T(\rho,x)$, is  continuous and can thus be taken as a decoration. If we  further equip $( T, {d}_{T}, \rho)$ with a finite measure $\upnu$, then 
 the function $g': T\to \R_+$ that assigns to every $x\in T$ the $\nu$-mass of the fringe subtree induced by this vertex, $g'(x)=\upnu(T_x)$, is clearly usc, and can also be used as a decoration of $T$.

  In the sequel, it will sometimes be  convenient to consider decorated compact trees not just with a single distinguished point (the root),  but more generally endowed  with a further finite or countable family of distinguished points in $T$,  referred to as  \textbf{marks}. Specifically, marks consist of a family $(r_i)_{i\in I}$  of points in $T$, where $I$ is a finite or countable set of indices. We stress that we do not request $r_i\neq r_j$ for $i\neq j$, and the same mark may arise for different indices. In particular, the notion of marks allow us to introduce  the \textbf{gluing operation} which will lie at the heart of the construction of decorated compact real trees from building blocks in the next section. 
  
\begin{figure}[!h]
 \begin{center}
  \includegraphics[width=15cm]{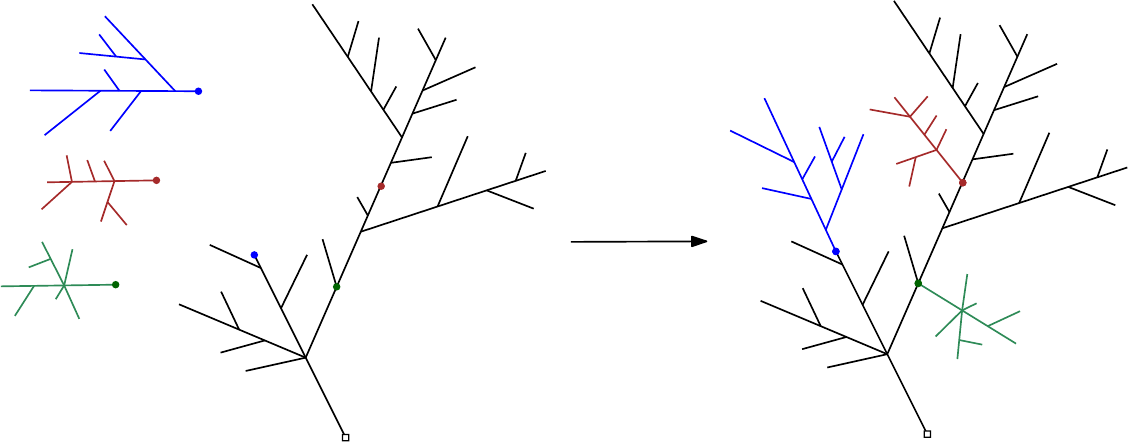}
 \caption{Gluing three subtrees on a base tree with three marked points.}\label{fig:gluingtree}
 \end{center}
 \end{figure}

Let us clarify what we mean by gluing decorated trees. In this direction, consider  a rooted compact real tree $(T', {d}_{T'},\rho')$ with marks $(r_i)_{i\in I}$. The tree $T'$ serves as a base, marks specify the locations on $T'$ where gluing will be performed.
Let further $(T_i, {d}_{T_i},\rho_i)_{i\in I}$  be a family of rooted compact real trees which are pairwise disjoint, and also disjoint from $T'$.  Roughly speaking, we will  glue each $T_i$ on $T'$ by identifying  the mark $r_i$ and the root $\rho_i$ of $T_i$, and equip the resulting space with the distance induced by $d_{T'}$ and the $d_{T_i}$ on each component. See Figure \ref{fig:gluingtree} for an illustration.

In order to describe rigorously this operation, we first introduce the (disjoint) union of those  trees
$$T^{\sqcup}\coloneq T' \sqcup \left(\bigsqcup_{i\in I}T_i\right).$$ 
We next define $ d^{\circ}:~T^{\sqcup}\times T^{\sqcup} \to \R_+$ by
$$d^{\circ}(x,y)\coloneq \left\{ \begin{matrix}
d_{T'}(x,y) &\text{if }x,y\in T',\\
d_{T_i}(x,y) &\text{if }x,y\in T_i \ \text{for some }i\in I,\\
d_{T'}(x, r_i) + d_{T_i}(\rho_i,y)&\text{if }x \in T' \text{ and }y\in T_i \ \text{for some }i\in I,\\
d_{T_i}(x,\rho_i) +d_{T'}(r_i,y) &\text{if }x \in T_i \ \text{for some }i\in I \text{ and }y\in T' ,\\
d_{T_i}(x,\rho_i) +d_{T'}(r_i,r_j)+ d_{T_j}(\rho_j,y) &\text{if }x\in T_i \text{ and }y\in T_j\ \text{for some }i\neq j\in I.
\end{matrix} \right.
$$
It is immediately seen that $d^{\circ}$ is a pseudo-distance on $T^{\sqcup}$ such that for any distinct points  $x,y$ in $T^{\sqcup}$, we have 
$d^{\circ}(x,y)=0$  if and only if, either $x=r_i$ and $y=\rho_i$ for some $i\in I$, or vice versa, or
$x=\rho_i$ and $y=\rho_j$ for some $i\neq j\in I$ such that $r_i=r_j$.

 We write $T$ for  the quotient space of $T^{\sqcup}$ for the equivalence relation 
 $$x \sim y\ \Longleftrightarrow \ d^{\circ}(x, y) = 0,$$ and  equip this quotient space with the metric $d_T$ induced by $d^{\circ}$, that is
$$d_T(\tilde x, \tilde y)=d^{\circ}(x,y), \qquad \tilde x, \tilde y\in T,$$
 for  any representatives $x\in \tilde x$ and $y\in \tilde y$ of these equivalence classes.
 In other words, $(T, d_T)$ is obtained from $(T^{\sqcup}, d^{\circ})$ by identifying $r_i$ and $\rho_i$ for every $i\in I$. 

\begin{lemma} \label{L:newL1} Suppose that
\begin{equation}\label{E:nullheight}
\left({\mathrm{Height}}(T_i)\right)_{i\in I} \
\text{ is a \textbf{null family}},
\end{equation}
meaning that for any $h>0$, the set of indices $i\in I$ with $\mathrm{Height}(T_i)\geq h$ is finite. 
Then $(T, d_T)$ is a compact real tree. We further use the equivalence class of $\rho'$ as the root $\rho$ of $T$.
\end{lemma} 

\begin{proof} We argued above that $(T,d_T)$ is a metric space; let us now check that it is a real tree, which should be intuitively clear. We use the obvious notation $\tilde z$ for the equivalence class of a point $z$  in $T^{\sqcup}$. Pick any $x,y$ in $T^{\sqcup}$. If both points belong to the same tree before the gluing,
say $x,y\in T_i$ for some $i\in I$, then $d^{\circ}(x,y)=d_{T_i}(x,y)= d_T(\tilde x, \tilde y)$, and  the path $\phi_i: [0,d_{T_i}(x,y)] \hookrightarrow T_i$ from $\phi_i(0)=x$ to $\phi_i(d_{T_i}(x,y))=y$ in $T_i$ yields in the obvious notation
a path $\tilde \phi_i: [0,d_{T}(\tilde x,\tilde y)] \hookrightarrow T$ from $\tilde \phi_i(0)=\tilde x$ to $\tilde \phi_i(d_{T}(\tilde x, \tilde y))=\tilde y$ in $T$. The case when, say $x\in T'$ and $y\in T_i$ (or \textit{vice versa}), is treated by  concatenating two paths, the first from $x$ to $r_i$  in $T'$ and the second from $\rho_i$ to $y$ in $T_i$. Finally the case when $x\in T_i$ and $y\in T_j$ with $i\neq j$ involves the concatenation of three paths,  the first from $x$ to $\rho_i$ in $T_i$,  the second (possibly degenerate)  from $r_i$ to $r_j$ in $T'$,  and the third from $\rho_j$ to $y$  in $T_j$. 

The initial trees being pairwise disjoint, uniqueness of the path $\tilde \phi$ from $\tilde x$ to $\tilde y$ in $T$ should be plain. Specifically, we may assume without loss of generality that $\tilde x \neq \tilde y$, pick an element in each equivalence class, say $x\in \tilde x$ and $y \in \tilde y$.
Suppose first that both $x$ and $y$ can be chosen in the same tree, say for simplicity $x,y\in T'$. Write $\widetilde{T'}\subset T$
for the set of equivalence classes of points in $T'$. If $\tilde \phi$ entered $T \backslash \widetilde{T'}$, then by continuity we could find two times $0\leq s < t \leq d_{T}(\tilde x,\tilde y)$ such that  $\tilde \phi(s)=\tilde \phi(t)=\tilde \rho_i$,
which is absurd. Thus $\tilde \phi$ stays in $\widetilde{T'}$ and hence defines unambiguously a path $\phi$ from $x$ to $y$ in $T'$. By uniqueness of the latter, we conclude that 
$\tilde \phi$ is also unique. The case when $x$ and $y$ belong to different initial  trees can be treated similarly. Essentially the same argument also shows that if $\psi: [0,1]\to T$ is a continuous injective map, then its image must coincides with the segment from $\psi(0)$ to $\psi(1)$ in $T$.

Finally, we check (sequential) compactness. Let $(\tilde y_n)_{n\geq 1}$ be a sequence in $T$; pick for each $n$ an element $y_n\in \tilde y_n$.
In the case when there exists some tree, say $T_i$, such that $y_n\in T_i$ for infinitely many $n$'s, then, since $T_i$ is compact,  we can extract from $(y_n)_{n\geq 1}$ a subsequence which converges in $T_i$, say towards $y$,
and it follows that there is a subsequence extract from $(\tilde y_n)_{n\geq 1}$ that converges towards the equivalence class $\tilde y$ of $y$. Next consider the complementary case when for all trees $T_i$, there are only finitely many $n$'s with $y_n\in T_i$,
and the same also holds for $T'$. We can then extract from $(y_n)_{n\geq 1}$ a subsequence such that each $y_n$ (along this subsequence) belongs to a different tree, say $y_n\in T_{i(n)}$. 
The assumption \eqref{E:nullheight} entails that $d_{T_{i(n)}}(y_n, \rho_{i(n)})$ converges to $0$ as $n\to \infty$ along this subsequence. By compactness of $T'$,
we can extract a further subsequence from $(r_{i(n)})$  that converges to a point $x\in T'$.  We conclude that $(\tilde y_n)_{n\geq 1}$ contains a subsequence that converges towards the equivalence class $\tilde x$ of $x$, and hence $T$ is compact.
\end{proof}

We next extend the gluing operation to decorations; recall that these must be given by usc functions on trees. Specifically, we consider  a decorated real tree $\texttt{T}'=(T', {d}_{T'},\rho',g')$ with marks $(r_i)_{i\in I}$ and 
 a family   of decorated real  trees, $\texttt{T}_i=(T_i, {d}_{T_i},\rho_i, g_i)$ for $i\in I$, such that the $T_i$, $i\in I$, are pairwise disjoint and also disjoint from $T'$. 
We define a function $ g$ on the glued tree $T$ which coincides of course with $g'$ on $T'\backslash \{r_i: i\in I\}$,   with  $g_i$ on $T_i\backslash \{\rho_i\}$, and is adjusted at the glued points in order to ensure upper semi-continuity. More precisely, define first the map $g^{\sqcup}: T^{\sqcup}\to \R_+$ by $g^{\sqcup}(y)=g'(y)$ if $y\in T'$ and 
$g^{\sqcup}(y)=g_i(y)$ if $y\in T_i$ for $i\in I$, and then 
$$ g(\tilde y)\coloneqq \sup_{y\in \tilde y} g^{\sqcup}(y), \qquad \tilde y \in T.$$

\begin{lemma} \label{L:newL2} Suppose \eqref{E:nullheight} and further that
\begin{equation}\label{E:nulldeco}
\left(\max_{y\in T_i} g_i(y)\right)_{i\in I} \
\text{ is a null family}.
\end{equation}
Then $\normalfont{\texttt{T}}\coloneqq (T, d_T, \rho, g)$ is a decorated real tree.
\end{lemma} 

\begin{proof} We have to verify that $g$ is usc. Let $\tilde y \in T$ and $a > g(\tilde y)$. We distinguish three possible situations. 

The simplest is when $\tilde y=\{y_i\}$ for some $y_i\in T_i$, with $i\in I$. Then $g(\tilde y)= g_i(y_i)$, and since necessarily $y_i\neq \rho_i$, we can find a small 
neighborhood $V_i$ of $y_i$ in $T_i$ such that $\rho_i\not\in V_i$ and  $g_i(z)\leq a$ for all $z\in V_i$. In particular the equivalence class $\tilde z$ of any $z\in V_i$ is reduced to $\{z\}$,
so $g(\tilde z)=g_i(z)\leq a$. In the obvious notation, $\widetilde V_i$ is a neighborhood of $\tilde y$ in $T$, and we conclude that $g$ is usc at $\tilde y$. 
The same argument applies when $\tilde y=\{y\}$ for some $y\in T'$ which does not belong to the closure of $\{r_i: i\in I\}$ in $T'$. 

Next suppose that $\tilde y$ is the equivalence class of a marked point in $T'$, that is 
$$\tilde y=\{x\}\sqcup\{\rho_j: j\in J\},$$ where $x$ is one of the marked points in $T'$ and $J=\{j\in I: r_j= x\}\neq \emptyset$. 
Let $V$ be a neighborhood of $x$ in $T'$ such that $g'(z)\leq a$ for all $z\in V$. Similarly, for every $j\in J$, there is a neighborhood $V_j$ of $\rho_j$ in $T_j$ such that $g_j(z)\leq a$ for all $z\in V_j$. 
The assumption \eqref{E:nulldeco} ensures that the set $J_a\coloneqq \{j\in J: \max_{z\in T_j} g_j(z)> a\}$ is finite. 
We deduce that 
$$G\coloneqq V\sqcup\left( \bigsqcup_{j\in J_a} V_j \right) \sqcup\left( \bigsqcup_{j\in J\backslash J_a} T_j \right)$$
is a neighborhood of $x$ in $T^{\sqcup}$ for the pseudo-distance $d^{\circ}$ and $g^{\sqcup}(z)\leq a$ for all $z\in G$. If we write $\widetilde G$ for the set of equivalence classes that $G$ induces in $T$, then $\widetilde G$ is a neighborhood of $\tilde y$ and $g(\tilde z) \leq  a$ for all $\tilde z\in \widetilde G$. 

The last case is when $\tilde y=\{x\}$ for some $x\in T'$ which is not marked, but adherent in $T'$ to the family of marked points. Recall from \eqref{E:nulldeco} that $I_a\coloneqq \{i\in I: \max_{z\in T_i} g_i(z)> a\}$ is finite. 
Since $x$ is not marked, we can find a neighborhood $V'$ of $x$ in $T'$ that avoids the mark $r_i$ for every $i\in I_a$, and then a possibly smaller neighborhood $V''\subset V'$  such that $g'(z)\leq a$ for all $z\in V''$. 
The set 
$$U\coloneqq V'' \sqcup \left( \bigsqcup_{i\not \in I_a} T_i\right)$$
is then a neighborhood of $x$ in $T^{\sqcup}$ for the pseudo-distance $d^{\circ}$, and we have $g^{\sqcup}(z)\leq a$ for all $z\in U$. If we write $\widetilde U$ for the set of equivalence classes $\tilde z$ for $z\in U$,   then $\widetilde U$ is a neighborhood of $\tilde y$  in $T$, and $g(\tilde z) \leq  a$ for all $\tilde z\in \widetilde U$. 

We have thus verified that $g$ is indeed usc on $T$. An appeal to Lemma \ref{L:newL1} completes the proof. \end{proof} 

When we want to record specifically all the elements involved in the gluing operator, we shall use the notation
$$\texttt{T}= \mathrm{Gluing}\Big( \texttt{T}', (r_i)_{i\in I}, \big(\texttt{T}_i\big)_{i\in I}\Big).$$

It will be convenient  to call \textbf{degenerate} a decorated real tree $\texttt{T}_i=(T_i,d_{T_i},\rho_i, g_i)$ such that $T_i=\{\rho_i\}$ is merely a singleton decorated with  $g_i(\rho_i)=0$; in that case we shall use the notation $\texttt{T}_i\sim \dagger$. 
Plainly, gluing a degenerate tree $\texttt{T}_i$ on a decorated real tree $\texttt{T}'$ is a neutral operation without incidence on the outcome (apart from removing the associated mark), and we may as well discard degenerate elements in gluing operations. We also stress that if the decorated real trees $\texttt{T}_i$, which are glued onto  on $\texttt{T}'$, are themselves marked, say $(r_{i,j})_{j\in J_i}$ on  $T_i$ for every $i\in I$, 
where we use pairwise disjoint set of indices $J_i$ for $i\in I$, 
then slightly more generally the gluing process yields a decorated real tree with marks indexed by $J=\bigsqcup_{i\in I}J_i$. The latter is naturally denoted by
$$ \mathrm{Gluing}\Big(\texttt{T}', (r_i)_{i\in I}, \big(\texttt{T}_i, (r_{i,j})_{j\in J_i}\big)_{i\in I}\Big).$$
In particular, this enables us to iterate the gluing operation, notably as it will be done in the next section.\par

\section{Construction of decorated real trees by gluing  building blocks} \label{sec:1.2}

We first introduce some standard  notation related to the \textbf{Ulam tree}, whose set of vertices\footnote{$\U$ is not a \textit{real} tree, but rather a tree in the sense of combinatorics.} consists of 
 finite sequence of positive integers,
$$\U\coloneq \bigcup_{k\geq 0}\N^k,$$ 
with the convention $\N^0\coloneq\{\varnothing\}$. The vertices are often referred to as individuals, since one may think of $\U$ as a population structured by its genealogy, where each individual has an infinite offspring. We use related vocabulary; for instance the edges 
of the Ulam tree connect parents to their children. 
  We write  $|u|$ for the generation  of $u\in \U$ (the length of the sequence $u$),  $u(j)$ for the predecessor at generation $j$ of $u$  (that is the prefix of $u$ with length $j$) whenever $|u|\geq j$, and  $uv$ for the sequence of length $|u|+|v|$ resulting from appending $v$ to $u$ (in particular, for $i\in \N$, $ui$ is viewed as the $i$-th child of $u$). We also write $\U^*=\U\backslash \{\varnothing\}$ for the set of nonempty finite sequences of positive integers
and  $u-$ for the parent of $u\in \U^*$, i.e. $u-=u(|u|-1)$. 

The Ulam tree will be used to index the building blocks appearing in the construction of decorated real trees through successive gluing. In order to unify the presentation, it is convenient to allow certain of these building blocks to be degenerate, where, as in Section \ref{sec:1.1},  degenerate elements can always be discarded in the construction. For simplicity, we also call a vertex $u$ \textbf{fictitious} if it indexes a degenerate block
and then write $u\sim \dagger$ for a fictitious vertex. We agree that only fictitious vertices $uv$ appear in 
the descent of a fictitious vertex $u$; as a consequence the subspace of non-fictitious vertices forms a subtree of $\U$. 

For the sake of simplicity, we shall present first the construction of decorated real trees by recursive gluing without measures, and then later on discuss natural Borel measures on such spaces.
The \textbf{building blocks} for the construction by gluing consist of a pair of families  
$$(f_u)_{u\in\U}\quad , \quad(t_{u})_{u\in \U^*}
 $$
that satisfy the following two properties for any non-fictitious $u\in \U$:
\begin{itemize}
\item  $f_u: [0,{z}_u]\to \R_+$ is a rcll (right-continuous with left limits) function with $z_u>0$; 
\item if a child $ui$ of $u$ is not fictitious, then $t_{ui}\in[0,z_u]$.
\end{itemize}
We further use the notation 
$$\sup f_u\coloneq \sup_{0\leq x \leq {z}_u} f_u(x) \quad \text{and}\quad \|f_u\|\coloneq{z}_u + \sup f_u,$$ 
with the agreement  for definitiveness that if $u\sim \dagger$ is fictitious, then $f_u$ is degenerate in the sense that  $z_u=0$, $f_u(0)=0$, and $\|f_u\|=0$.

We always assume that 
\begin{equation} \label{E:fam0}
\big(\| f_{u}\| : u\in  \mathbb{U}\big) \text{ is a null family}.
\end{equation}

 Before we formally define the gluing construction, let us provide an informal explanation of the roles of the  building blocks. A (non-fictitious) interval $[0,z_u]$ will become a segment of the resulting real tree $T$, and  the real number $t_{ui}\neq \dagger$ will be used in the gluing construction to specify the location on  the parent segment where the $i$-th child segment  labeled by $ui$ is glued. In particular, $t_{ui}$ will yield a branching point of the real tree
 whenever $0<t_{ui}<z_u$ and $z_{ui}>0$. 
  When $u$ has only fictitious children, the construction regarding the subtree rooted at $u$ stops at that point, in the sense that no segments are glued on the segment labeled by $u$. Finally, the family $(f_u)_{u\in\U}$ will correspond to the usc-decoration. In this direction,  we must address a technical  issue regarding the regularity of trajectories. When dealing with most random processes in continuous time (e.g. a Feller process), one usually works with rcll versions, which is the reason why we requested the functions $f_u$ above to be rcll. Nonetheless,   we have to consider rather their ucs versions in order to conform to the framework introduced in Section~\ref{sec:1.1}. Formally, if $f:[0,z]\to \mathbb{R}_+$ is rcll function, then its  usc version $\check{f}$ is defined by the relation:
  \begin{eqnarray} \label{eq:cadlagtousc}\check{f}(t)\coloneq\max\big\{f(t-),f(t)\big\},\quad \text{for } t\in[0,z],  \end{eqnarray}
 with the convention $f(0-)=f(0)$. This only affects the values of the function at times when it has a negative jump (including possibly at lifetime).

We also agree to write simply $d$ for the usual distance on any interval. For every non-fictitious $u\in \mathbb{U}$,  we view $f_u$ as a decorated compact interval $\big([0,z_u],d,0, \check{f}_u\big)$, and the family of non-fictitious times $t_{ui}$ for $i\geq 1$ as marks on $[0,z_u]$. The gluing operation described in the preceding section uses pairwise disjoint real trees, and we shall therefore  introduce disjoint isomorphic copies of the preceding for different $u\in \U$. Namely, we use the Ulam tree to differentiate these segments, and, in the notation of Definition \ref{Def:dectree}, we consider for every non-fictitious  $u\in \mathbb{U}$ the decorated  segment 
$$\texttt{S}_u\coloneqq(S_u, d_u, \rho_u, \hat f_u),$$
where $S_u\coloneqq\big\{(u,t):~ t\in [0,z_u]\big\}$,  $ \rho_u\coloneqq(u,0)$, and  the metric $d_u$ and usc decoration $\hat f_u$ are defined by
$$d_u((u,s), (u,t))\coloneqq d(s,t) \quad \text{ and } \quad \hat f_u(u,t)\coloneqq \check{f}_u(t), \quad \text{for } s,t\in [0,z_u]. $$
In particular,  the segments $S_u$ are pairwise disjoint sets and $\texttt{S}_u$ is isomorphic to $\big([0,z_u],d,0, \check{f}_u\big)$, in the sense that  there is a bijective  isometry $\varphi_u : [0,z_u] \to S_u$ with $\varphi_u(0)= \rho_u$, and $\hat f$ is the usc version of $f\circ \varphi_u^{-1}$. We stress that  $\varphi_u$ is actually unique, since we are dealing with oriented segments.
 Last but not least, we further mark each segment 
using the family of (non-fictitious) $t_{ui}$ with $i\geq 1$. More precisely, if we let $I_u=\{uj: uj\neq \dagger\}$ 
for the non-fictitious offspring of the individual $u$ (introducing this notation is needed as we want to use different set of indices for the marks arising from different blocks), then the family of marks on $S_u$  is given by  $(\hat{t}_{v}\coloneqq (u,t_v))_{v\in I_u}$.

We can now proceed with the formal construction of a decorated real tree from building blocks, where at first, measures are discarded. We introduce first the disjoint union of segments
$$T^{\sqcup}\coloneqq \bigsqcup_{u\not \sim \dagger} S_u,$$
where the union is taken over the set of non-fictitious vertices $u\in \U$. We use $\rho_{\varnothing}$ as root, and write $d^0$ for the natural distance on $T^{\sqcup}$ which is given by $d_u$ on each segment $S_u$ and such that
$d^0(x,y)=\infty$ when $x$ and $y$ belong to different segments (i.e. $x$ and $y$ are in different connected components). We also define unambiguously  the map $g^{\sqcup}: T^{\sqcup} \to \R_+$ by
$$g^{\sqcup}(x)\coloneqq \hat{f}_u(x),\qquad\text{for every }x\in S_u \text{ and }u\not \sim \dagger.$$

We then construct recursively\footnote{Instead of developing a recursive construction, one could have given directly  an explicit formula for $d^n$ in terms of distances $d_w$ for all the vertices $w$ on the segment from $u$ to $v$ in $\U$. However the formula would be a bit cumbersome to state precisely and not quite transparent, at least for the first reading.}  a sequence $(d^n)_{n\geq 1}$ of pseudo-distances on $T^{\sqcup}$, using the gluing operator of Section \ref{sec:1.1}; note that \eqref{E:fam0} ensures the requirements  \eqref{E:nullheight} and \eqref{E:nulldeco}  in this setting.
Specifically, $d^1$ is the pseudo-distance on $T^{\sqcup}$ obtained by identifying  $\rho_{i}$ and the 
marked point $\hat{t}_{i}$ for each non-fictitious vertex $i$ at the first generation, that is
$i\in I_{\varnothing}$.  Segments at generations $2$ and more are not affected, in particular $d^1(x,y)=\infty$ for any 
$x\in S_u$ and $y\in S_v$ with $u\neq v$ and $|u|\vee |v|\geq 2$. 
Next $d^2$ is the pseudo-distance on $T^{\sqcup}$ obtained by further identifying $\rho_{ij}$ and the 
marked point $\hat{t}_{ij}$ for each vertex $ij\in I_i$ and each $i\in I_{\varnothing}$. And so on, and so forth, generation after generation. 

Then consider any $x,y\in T^{\sqcup}$, say $x\in S_u$ and $y\in S_v$. Plainly, if $u=v$, then
$d^n(x,y)=d^0(x,y)$ for all $n\geq 0$. If  $u\neq v$, then $d^n(x,y)=\infty$ when $n<|u|\vee |v|$, whereas
$d^n(x,y)=d^{n'}(x,y)< \infty$ for all $n,n'\geq |u|\vee |v|$. We can then set
$$d^{\circ}(x,y)\coloneqq \lim_{n\to \infty} d^n(x,y),$$
and $d^{\circ}$ defines a pseudo-distance on $T^{\sqcup}$ which is now everywhere finite (this can be interpreted as connectivity). 
Alternatively, we could also define directly $d^{\circ}$ as the largest pseudo-distance on $T^{\sqcup}$ which coincides with $d_u$ on each segment $S_u$, and such that $d^{\circ}(\hat t_v,\rho_v)=0$ for every non-fictitious 
$v\in \U^*$. 

We next define $T^{\circ}$ as the quotient space for the equivalence relation 
$$x\sim y\ \Longleftrightarrow \ d^{\circ}(x,y)=0, \qquad  x,y\in T^{\sqcup}.$$
We write $d_{T^{\circ}}$ for the distance on $T^{\circ}$ induced by $d^{\circ}$ and
 also define the map $g^{\circ}: T^{\circ}\to \R_+$ given by
$$g^{\circ}(\tilde x) \coloneqq \sup_{x\in \tilde x}g^{\sqcup}(x), \qquad \tilde x\in T^{\circ}.$$
Note that $g^{\circ}$ is actually bounded, thanks to \eqref{E:fam0}.

The following claim can be viewed as the analog of Lemma \ref{L:newL2} for  infinite iterations of the gluing operation, and its proof relies heavily on this lemma. Recall that a metric space is called pre-compact if, for every $\varepsilon>0$, it can be covered by finitely many balls with radius $\varepsilon$.
\begin{lemma}\label{L:newL3} Assume  \eqref{E:fam0} and 
\begin{equation}\label{E:seriehypbis}
\lim_{k \to \infty} \sup\Big\{ \sum_{n=k}^{\infty} z_{ \bar{u}(n)}: \bar{u} \in \mathbb{N}^{\mathbb{N}}\Big\} =0,
\end{equation}
where  $\bar{u}(n)\in \U$ stands for the prefix at generation $n$ of an infinite sequence $ \bar{u} \in \mathbb{N}^{\mathbb{N}}$ of positive integers. 
Then $(T^{\circ}, d_{T^{\circ}})$ is a pre-compact real tree and $g^{\circ}$ an usc function. 
\end{lemma}
We briefly postpone the proof of Lemma \ref{L:newL3}  and mention that since plainly
$$ \sup\Big\{ \sum_{n=k}^{\infty} z_{ \bar{u}(n)}: \bar{u} \in \mathbb{N}^{\mathbb{N}}\Big\}  \leq \sum_{n=k}^{\infty} \sup \big\{{z}_u: u\in \N^n\big\} , $$
we shall often use 
 in the sequel  the simpler but stronger requirement 
\begin{equation}
 \label{E:seriehyp}\sum_{n=0}^{\infty} \sup \big\{{z}_u: u\in \N^n\big\} <\infty,
 \end{equation}
which implies  \eqref{E:seriehypbis}.

Taking Lemma \ref{L:newL3} for granted, we can easily finalize the construction of decorated real trees from building blocks as follows.   We let $(T,d_T)$ be the completion of the metric space $(T^{\circ}, d_{T^{\circ}})$ and write simply  $\rho$ for the equivalence class of $\rho_{\varnothing}$ in $T^{\circ}$.
We also extend $g^{\circ}$ 
to the boundary $T\backslash T^{\circ}$ and
define the map $g: T\to \R_+$ by $g(y)=g^{\circ}(\tilde y)$ if $y=\tilde y\in T^{\circ}$,  and $g(y)=0$ if $y\in T\backslash T^{\circ}$.

\begin{theorem}\label{T:recolinfty} Assume that the building blocks fulfill  \eqref{E:fam0} and \eqref{E:seriehypbis}. Then $ \mathtt{T} = (T,d_T,\rho, g)$ is a decorated real tree. \end{theorem}
\begin{proof} The completed space $(T,d_T)$ plainly inherits connectivity and the  four point condition \eqref{E:4points}
from $(T^{\circ}, d_{T^{\circ}})$, and is therefore a real tree. We also know from Lemma \ref{L:newL3} that $(T^{\circ}, d_{T^{\circ}})$ is a pre-compact space, so its completion $(T,d)$ is compact. We turn our attention to the extension $g$ of $g^{\circ}$. Take first $x\in T^{\circ}$ and a real number $a>g^{\circ}(x)$. Since $g^{\circ}$ is usc at $x$, we can choose $ \varepsilon >0$ small enough such that $g(y)=g^{\circ}(y)< a$
for all $y\in T^{\circ}$ such that $d_T(x,y)<\varepsilon$. Since $g(y)=0$, for every $y\in T\backslash T^{\circ}$, it follows that $g$ is usc at $x$.  Suppose finally that $x\in T\backslash T^{\circ}$, and take any $a>0$. 
The distance from $x$ to any segment $S_u$ with $\|f_u\|\geq a$ is bounded away from $0$, since  by \eqref{E:fam0}, there is only finitely many such segments.
We can again choose $ \varepsilon >0$ small enough such that $g(y)=g^{\circ}(y)< a$ for all $y\in T^{\circ}$ such that $d_T(x,y)<\varepsilon$, 
which shows  that $g$ is usc at $x$. So $g$ is an usc function on $T$, and the proof is complete. 
\end{proof}
\begin{remark}\label{R:boundaryvsbranch} Beware that, in spite of what the notation  might suggest,  $T^{\circ}$ should not be thought of as the interior of $T$. It should be plain that every boundary point of $T^{\circ}$ is always a leaf of $T$, but the converse usually fails. In other words, we have
$$T\backslash T^{\circ} \subseteq \partial T  ,$$
and the inclusion can be strict in general. 
\end{remark}

Let us now establish Lemma \ref{L:newL3}.

\begin{proof}[Proof of Lemma \ref{L:newL3}] By construction, $(T^{\circ}, d_{T^{\circ}})$ is a connected metric space; let us show that it is a real tree by checking the four point condition \eqref{E:4points}. Pick any $x_1, x_2, x_3, x_4$ in $T^{\sqcup}$, and let $k$ denote the maximal generation of the indices $u\in \U$ such that the segment $S_u$ contains at least one of those points. 
Since we know by iteration from Lemma \ref{L:newL1}  that gluing the segments up to the $k$-th generation produces a tree and that trees satisfy the four point condition, we have
$$d^k(x_1,x_2)+d^k(x_3,x_4)\leq (d^k(x_1,x_3)+d^k(x_2,x_4))\vee (d^k(x_1,x_4)+d^k(x_2,x_3)).$$
Recall that  $d^{\circ}(x,y)=d^k(x,y)$ for any points $x,y$ in segments indexed by vertices at generations at most $k$, 
so we have as well
$$d^{\circ}(x_1,x_2)+d^{\circ}(x_3,x_4)\leq (d^{\circ}(x_1,x_3)+d^{\circ}(x_2,x_4))\vee (d^{\circ}(x_1,x_4)+d^{\circ}(x_2,x_3)).$$
Denote by $\tilde x_i$ the equivalence class  of $x_i$ for the pseudo-distance $d^{\circ}$, with $i=1, 2,3, 4$. We can rewrite
the above inequality as 
$$d_{T^{\circ}}(\tilde x_1,\tilde x_2)+d_{T^{\circ}}(\tilde x_3,\tilde x_4)\leq \big(d_{T^{\circ}}(\tilde x_1,\tilde x_3)+d_{T^{\circ}}(\tilde x_2,\tilde x_4)\big)\vee \big(d_{T^{\circ}}(\tilde x_1,\tilde x_4)+d_{T^{\circ}}(\tilde x_2,\tilde x_3)\big),$$
and we conclude that $(T^{\circ}, d_{T^{\circ}})$ is a real tree.

We next turn our attention to the pre-compactness assertion. Fix $\varepsilon>0$ arbitrarily small; the Assumption \ref{E:seriehypbis} allows us to pick $k\geq 1$ sufficiently large so that
$$\sum_{n=k+1}^{\infty} z_{ \bar{u}(n)}< \varepsilon/2,$$
for any infinite sequence $ \bar{u} \in \mathbb{N}^{\mathbb{N}}$. 
It follows from the construction by gluing that for  every $y \in T^{\sqcup}$, we can find some $x(y)\in T^{\sqcup k}\coloneqq \bigsqcup_{|u|\leq k} S_u$ such that 
$$d^{\circ} (y,x(y)) \leq \varepsilon / 2.$$ 
Indeed, this claim is trivial if $y\in S_v$ for some vertex $v$ at generation $|v|\leq k$, and otherwise we can take for
$x(y)$  the marked point $\hat t_{v(k+1)}$ on $S_{v(k)}$ (recall that $v(\ell)$ denotes the predecessor of the vertex $v\in \U$ at generation $\ell\leq |v|$) so that then
\begin{align*}d^{\circ} (x(y),y) &\leq \sum_{n=k+1}^{|v|-1} d_{v(n)}(\rho_{v(n)}, \hat t_{v(n+1)}) + d_v(\rho_v,y)\\
&\leq  \sum_{n=k+1}^{|v|-1}  t_{v(n+1)} + z_v\\
&\leq \sum_{n=k+1}^{|v|} z_{ v(n)}.
\end{align*}
Denote by $T^{\circ k}\subset T^{\circ}$ the subset of the equivalence classes of points in $T^{\sqcup k}$ and infer
from above that for any $\tilde y\in T^{\circ}$, there is some $\tilde x(\tilde y)\in T^{\circ k}$ such that 
$$d_{T^{\circ}}(\tilde y, \tilde x(\tilde y)) \leq \varepsilon/2.$$
On the other hand, we deduce by induction  from Lemma \ref{L:newL1} that $T^{\circ k}$ endowed with the distance $d_{T^{\circ}}$ is compact. Thus there exists a finite sequence in  $T^{\circ k}$, say $\tilde x_1, \ldots, \tilde x_n$ such that the sequence of balls with radius $\varepsilon/2$ centered at those points cover $T^{\circ k}$.
We conclude from the triangle inequality that the sequence of balls with radius $\varepsilon$ and centered at the $\tilde x_i$ now cover the whole $T^{\circ}$, and $T^{\circ}$ is hence pre-compact. 

We finally check that $g^{\circ}$ is an usc function on $T^{\circ}$. Fix $\varepsilon>0$; from \eqref{E:fam0} we can choose $k\geq 1$ sufficiently large so that
$\sup f_u \leq \varepsilon$ for any $u\in \U$ at generation $|u| >k$. Define the function $g^{\circ}_k$ on $T^{\circ k}$ by
$$g^{\circ}_k(\tilde x)=  \sup\left \{g^{\sqcup}(y):  y\in \tilde x \cap T^{\sqcup k} \right\}.$$
We know from Lemma \ref{L:newL2} and an iteration argument that $g^{\circ}_k$ is an  usc function on $T^{\circ k}$. Therefore, for any $\tilde x\in T^{\circ k}$, there is $\varepsilon'>0$ such that
$$g^{\circ}_k(\tilde x') \leq g^{\circ}_k(\tilde x)+ \varepsilon, \quad \text{for every $\tilde x'\in T^{\circ k}$ with }d_{T^{\circ}}(\tilde x, \tilde x')< \varepsilon'.$$
It then follows from the choice of $k$ that
$$g^{\circ}(\tilde x') \leq g^{\circ}(\tilde x)+ \varepsilon, \quad \text{for every $\tilde x'\in T^{\circ}$ with }d_{T^{\circ}}(\tilde x, \tilde x')< \varepsilon',$$
proving that $g^{\circ}$ is indeed an usc function at any $\tilde x\in T^{\circ}$.
\end{proof}

It will be convenient to introduce a notation for points in $T^\circ$ given by equivalence classes of marks (possibly fictitious) on the initial segments. For every   vertex $v\in \U^*$, consider the mark $\hat t_v$ on the segment $S_{v-}$ (recall that $v-$ stands for the parent of the vertex $v$). We then write
$\uprho(v)\in T$ for  the equivalence class of  $\hat t_v$, that is also the equivalence class of the root $\rho_v$ of the segment $S_v$ as those two points are identified in $T$, and by convention set $\uprho(\varnothing):=\rho$.
The reader will easily check that any branching point of $T$, say $b$, is of the form $b= \uprho(v)$ for some non-fictitious vertex $v\in \U$. We also stress that in the converse direction, $\uprho(v)$ is not necessarily a branching point of $T$ (counter-examples arise in the situation where $t_v=z_{v-}$ and there are no other non-fictitious vertices aside $v$ in the offspring of the parent $v-$). 
A bit more generally, we will sometimes use the following notation to identify points in $T^{\circ}$ even when they have not been marked. Consider any vertex $v\in \mathbb{U}$ and any $t\in[0,z_v]$.
We write $\rho_v(t)$ for the unique point on the segment $S_v$ at distance $t$ from the root $\rho_v$, and then $\uprho(v,t)$ for the corresponding point in $T^{\circ}$ (strictly speaking, $\uprho(v,t)$ is the equivalence class of $\rho_v(t)$, which is actually reduced to the singleton $\{\rho_v(t)\}$ except if $\rho_v(t)$ has been marked). 
 
We also point out that the proof of Lemma \ref{L:newL3} also yields some useful information, notably  about the height of $T$ and the maximum of the decoration $g$, which we record in the next statement.
In this direction, recall that we wrote there $T^{\circ k}$ for  the subset of $ T^{\circ}$ given by the equivalence classes of points in $T^{\sqcup k}=\bigsqcup_{|u|\leq k} S_u$. On the other hand, we infer from Lemma \ref{L:newL1} that $T^{\sqcup k}$ is sequentially compact for the pseudo-distance $d^k$, that is we can always extract from any sequence in $T^{\sqcup k}$  a subsequence which converges in $T^{\sqcup k}$  for the pseudo-distance $d^k$. \textit{A fortiori}, this holds as well for the  pseudo-distance $d^{\circ}$, and we deduce that $T^{\circ k}$ is actually closed for the distance $ d_{T^{\circ}}$. We shall henceforth use the simpler notation $T^k=T^{\circ k}$, and view $T^k$ as a closed subset of $T$. 

\begin{proposition} \label{T:recolheightsup} Assuming  \eqref{E:fam0} and \eqref{E:seriehypbis}, we have for every $k\geq 0$ that
$$d_T(y, T^k) \leq  \sup\Big\{ \sum_{n=k+1}^{\infty} z_{\bar{u}(n)}: \bar{u} \in \mathbb{N}^{\mathbb{N}}\Big\} , \qquad \text{ for all }y\in T$$
  and 
$$
\mathrm{Height}(T) \leq \sup\Big\{ \sum_{n=0}^{\infty} z_{\bar{u}(n)}: \bar{u} \in \mathbb{N}^{\mathbb{N}}\Big\}.
$$
Moreover, we have also
$$\sup_{y\in T \backslash T^k}g(y) \leq  \max \big\{\sup f_u:  u\in \U, |u|\geq k+1 \big\}.$$
 \end{proposition}
 
We also note the following direct consequence concerning Hausdorff dimensions.

\begin{lemma} \label{lem:dimensionabstrait}  Assume  \eqref{E:fam0} and \eqref{E:seriehypbis}, and write
\begin{equation}\label{eq:partial:0:T}
\partial_0 T\coloneq\{ x \in \partial T : g(x) =0\},
\end{equation}
for the set of leaves with decoration $0$. Then, we have 
$$\mathrm{dim}_{H}(T) \leq \mathrm{dim}_{H}(\partial_0 T) \vee 1.$$
\end{lemma}

\begin{proof} Recall that for any $n\geq 0$, the subset $T^n\subset T$ induced by the countable collection of segments $S_u$ with generation $|u|\leq n$ is a closed subtree. Therefore  its Hausdorff dimension cannot exceed $1$. Moreover, $T^\circ= \bigcup_{n\geq 0} T^n$ obviously contains the skeleton $T\backslash \partial T$ of $T$, and by definition we have  $g(x)=0$   for every $x\in T\setminus  T^\circ$. 
\end{proof}

We now conclude this section and record notation for two important families of (decorated) subtrees. The second and third have already appeared in this chapter, whereas the first will be used later.

\begin{notation}\label{N:subtrees} Let $\texttt{T}=(T, d_T, \rho,g)$ denote the decorated real tree constructed from the building blocks $(f_u)_{u\in \U}$ and $(t_u)_{u\in \U^*}$ 
 in Theorem \ref{T:recolinfty}.
\begin{enumerate}

\item[(i)] For every non-fictitious vertex $u\in \U$,  the building blocks indexed by descendants
of $u$, namely $(f_{uv})_{v\in \U}$, $(t_{uv})_{v\in \U^*}$ and $(m_{uv})_{v\in \U}$ also fulfill the requirements of Theorem \ref{T:recolinfty}.
We write $\texttt{T}_u=(T_u, d_{T_u}, \uprho(u), g_u)$ for the decorated real tree constructed from the latter by gluing\footnote{We make here a slight abuse in order to view $T_u$ as a subtree of $T$: we have previously defined $\uprho(u)\in T$ as  the equivalence class of the root $\rho_u$ of the segment $S_u$, which is actually larger than the equivalence class  obtained when the gluing construction  restricted to descendants of $u$. A similar minor abuse is made for the same reason in (ii). }.

\item[(ii)] For every generation $n\geq 0$, we write $\texttt T^n=(T^n, d_{T^n}, \rho, g^n)$ for the decorated real tree constructed by gluing the building blocks up to generation $n$ only,
that is from the building blocks $(f_u^n)_{u\in \mathbb{U}}$ and  $(t^n_u)_{u\in \mathbb{U}^*}$, where $f_u^n=f_u$ and  $t_u^n=t_u$  when $|u|\leq n$, whereas these quantities are fictitious when $|u|>n$.

\item[(iii)] We write $T^\circ = \bigcup_{n\geq 1} T^n$. The set $T\backslash T^\circ$ of adherence points of $T^\circ$ is included into the set of leaves $\partial T$. 
\end{enumerate}
\end{notation}
We stress that $T_u$ is always a subtree of the fringe subtree $T_{\uprho(u)}$ rooted at $\uprho(u)\in T$, and that the inclusion is strict in general. Last, note also that although  $g^n$ is always dominated by the restriction of $g$ to $T^n$, these two functions may be different only at marked points.

\section{Measured decorated trees} \label{sec:1:measures}

We will now equip the compact real tree $T$ constructed in Theorem \ref{T:recolinfty} with some finite measures. First, as any compact real tree, our tree $T$ is naturally equipped with the length measure $\uplambda_T$, which is given by the one-dimensional Hausdorff measure on the skeleton, see \cite[Section 2.4]{EPW06}.  In words,   $\uplambda_T$ mirrors on each segment  of $T$ the Lebesgue measure on intervals, where we recall that segments of $T$ are isometric to intervals of $\mathbb{R}_+$. We stress that $\uplambda_T$ is usually only a sigma-finite measure. However, 
we can consider a density function $f : T\to \R_+$ in $L^1(\uplambda_T)$ and then take $\upnu(\d x)=f(x)\uplambda_T(\d x)$. Densities that will appear naturally in this work depend on the decoration, that is,  they are of the type $f=\varpi\circ g$, where $\varpi:\mathbb{R}_+\to \mathbb{R}_+$ is  a measurable function. The resulting measure on $T$ is then denoted by $ \varpi\circ g \cdot \uplambda_T$ and often called the \textbf{weighted length measure}. We stress that the length measure, and therefore \textit{a fortiori} $ \varpi\circ g \cdot \uplambda_T$ as well, gives no mass to the set of leaves $\partial T$.  We also stress that, since the set of  points  $\{\uprho(v):~v\in \mathbb{U}\}$ is countable and thus receives no mass from $\uplambda_T$, we have
\begin{equation}\label{E:newEmassum} 
\int_T \varpi\circ g(x) \uplambda_T(\dd x) = \sum_u \int_0^{z_u} \varpi\circ f_u(t) \dd t,
\end{equation}
where  the sum in the right-hand side is implicitly taken over non-fictitious vertices. In particular, the requirement $\varpi\circ g \in L^1(\uplambda_T) $
is equivalent to the finiteness of the previous display.

We shall now  construct another class of measures on $T$, which, at the opposite, are typically  carried by  $\partial T$. 
This requires to introduce, on top of the building blocks,  a further family of nonnegative real numbers
$$(m_u)_{u\in \mathbb{U}}$$
with $m_u=0$ when $u\sim \dagger$ is a fictitious vertex, and most importantly,
\begin{equation}\label{E:massspreads}
\sum_{i=1}^{\infty} m_{ui}=m_u \qquad \text{for all }u\in \U.
\end{equation} 
Note that  if $m_u>0$, then \eqref{E:massspreads} implies that $u$ has at least one infinite line of descent with only non-fictitious vertices.
Next recall that  $\uprho(u)$ denotes the point in $T$ associated to a mark $t_u$.
We then define for every $n\geq 0$ the purely atomic measure $\upmu^n$ on $T$ given by
\begin{equation}\label{eq:def:mu:n}
\upmu^n\coloneqq \sum_{|u|=n} m_u \cdot \delta_{\uprho(u)},
\end{equation}
where the notation $\delta_x$ is used for the Dirac point mass at $x$, and we implicitly agree to ignore fictitious vertices in the sum.
It is immediately seen by induction from the requirement \eqref{E:massspreads}  that $\upmu^n(T)=m_{\varnothing}$ for any $n\geq 0$.

\begin{proposition} \label{P:newPmass} Assume that the building blocks fulfill \eqref{E:fam0}, \eqref{E:seriehypbis} and \eqref{E:massspreads}.
The sequence $(\upmu^n)_{n\geq 1}$ converges in the sense of Prokhorov towards a  Borel measure on $T$ denoted by $\upmu$. One has $\upmu(T)=m_{\varnothing}$.
\end{proposition}
\begin{proof} Recall that for any $n\geq 0$, the set $T^n$ denotes the subtree of $T$ obtained by gluing the collection of segments $S_u$ with generation $|u|\leq n$.
We write $p^{n}: T\to T^n$ for the projection on $T^n$, that is for any $y\in T$, $p^n(y)$ is the right endpoint of the segment $ \llbracket \rho,y \rrbracket \cap T^n$. 
On the one hand, the first claim in Proposition \ref{T:recolheightsup} entails
$$d_T\big(y,p^n(y)\big) \leq \sup\Big\{ \sum_{k=n+1}^{\infty} z_{\bar{u}(k)}: \bar{u} \in \mathbb{N}^{\mathbb{N}}\Big\}.$$
On the other hand,  we infer by iteration from \eqref{E:massspreads} that for any $n'\geq n$, the push-forward of $\upmu^{n'}$ by $p^n$ coincides with $\upmu^n$.

We deduce from these two observations that the Prokhorov distance between $\upmu^n$ and $\upmu^{n'}$ is at most 
$$\d_{\mathrm{Prok}}(\upmu^n, \upmu^{n'}) \leq \sup\Big\{ \sum_{k=n+1}^{\infty} z_{\bar{u}(k)}: \bar{u} \in \mathbb{N}^{\mathbb{N}}\Big\}.$$
Since \eqref{E:seriehypbis} requests the right-hand side to converge to $0$ as $n\to \infty$, the sequence $(\upmu^n)_{n\geq 1}$ is Cauchy on the space of measures on $T$ with total mass $m_{\varnothing}$,
and our claim follows. 
\end{proof} 

We also point out that, under a minor additional hypothesis, the measure $\upmu$ constructed above
is carried by the boundary points of the pre-compact tree $T^\circ$. Recall that the latter denotes the equivalence class of $T^{\sqcup}$; see
also Remark \ref{R:boundaryvsbranch}.

\begin{proposition}\label{P:newPleaves} Assume  \eqref{E:fam0}, \eqref{E:seriehypbis} and \eqref{E:massspreads}. Suppose further that for  any non-fictitious vertex $v\in \U^*$, the  mark $t_v$  and the length $z_v$ are strictly positive. Then we have $\upmu(T^\circ) =0$, and as a consequence, $\upmu$ is carried by the subset of leaves $T\setminus T^\circ$.
\end{proposition}
\begin{proof} We shall check first that  the root $\rho$ is not an  atom of $\upmu$, which should be intuitively clear.
 In this direction, write $\tilde S_{\varnothing}:=\{\upvarrho(\varnothing, t):~t\in[0,z_\varnothing]\}$ for the segment in $T$ induced by the equivalence class of the ancestral segment $S_{\varnothing}$. Then for any $a>0$, consider the subset $T_{<a}$  of points $x\in T$ 
 such that the length (i.e. the $\uplambda_T$ measure) of the segment $\tilde S_{\varnothing}\cap \llbracket \rho,x \rrbracket $ is less than $a$.
  By \eqref{E:massspreads} and the construction by gluing, we have
 $$\upmu(T_{<a}) = \sum_{j\geq 1} \indset{t_j<a}m_j,$$
 and since we assumed that the (non-fictitious) marks $t_j$ are strictly positive, we have
 $$\lim_{a\to 0+} \upmu(T_{<a})=0.$$
Also by construction, the open ball in $T$ centered at $\rho$ and with radius $a$ is contained in  $T_{<a}$. Letting $a\to 0+$, we conclude that $\upmu(\{\rho\})=0$.

Next,  recall that for every non-fictitious vertex $j\geq 1$ at the first generation, $\uprho(j)\in \tilde S_{\varnothing}$ denotes the equivalence class of the mark $\hat t_j$ on the ancestral segment $S_{\varnothing}$. 
On the one hand,  the argument above shows that we have also $\upmu(\{\uprho(j)\})=0$.
On the other hand, the requirement \eqref{E:massspreads} entails that the entire mass of $\upmu$ is carried by the union  for $j\geq 1$ of the fringe subtrees rooted at the $\uprho(j)$.
This shows that $\upmu$ assigns zero mass to  the equivalence class of ancestral segment, $\upmu(\tilde S_{\varnothing})=0$. We conclude by iteration on generations that $\upmu$ 
assigns zero mass to  $T^{\circ}$,  and \textit{a fortiori} to the skeleton $T\backslash \partial T$. 
\end{proof}

\section{Hypographs, topologies, and isomorphic identifications} \label{sec:1.3}

We have developed so far general material on (measured) decorated real trees and their construction, and roughly speaking, we now would like to compare two different decorated real trees one  with the other.
More precisely, our main motivation is to give a rigorous definition of a notion of convergence for sequences of these objects. Actually, although dealing with usc functions on compact metric spaces is fundamental to our approach,  tree structures are essentially irrelevant for  this question, and we shall develop first a more general framework  that could also be used for other purposes.
It will be only at the end of this section that the special case of real trees will be addressed  more specifically.

Let  $(Y, d_Y)$ be a Polish space and ${\mathcal K}(Y)$ denote the set of non-empty compact subspaces in~$Y$. 
Consider  an usc function $g:K\to \R_+$ for some $K\in {\mathcal K}(Y)$ which we refer to as the domain of $g$. 
Upper semi-continuity enables us to view $g$ as a compact subspace of the larger Polish space $Y\times \R_+$ by introducing  the \textbf{hypograph} 
$$\mathrm{Hyp}(g) \coloneq \big\{(x,r): x\in K \text{ and } 0\leq r \leq  g(x)\big\} \subset Y\times \R_+.$$ 
Specifically, the product space   $Y\times \R_+$ is naturally equipped with the distance
 \begin{eqnarray} \label{eq:disthypo}
d_{Y\times \R_+}\big((y,r),(y',r')\big)\coloneqq d_{Y}(y,y') \vee |r-r'|, 
  \end{eqnarray}
and  $\mathrm{Hyp}(g)$ is then a compact subspace\footnote{This is one of the main reasons why we considered usc functions.} of the Polish space $(Y\times \R_+, d_{Y\times \R_+})$. 
Plainly, the hypograph $\mathrm{Hyp}(g)$ determines $g:K\to \R_+$. Namely,  the domain  $K$ is the base of the  hypograph, that is the 
image of $\mathrm{Hyp}(g)$ by the first projection $p_1: Y\times \R_+\to Y$, and  
$$g(x)=\max\big\{r\geq 0: (x,r)\in \mathrm{Hyp}(g)\big\},\quad \text{ for every }x\in K. $$

In this work, we will interpret a decorated tree through its hypograph. Before continuing with the general study of hypographs and defining a suitable metric for comparing them, let us highlight the following consequence of Lemma~\ref{lem:dimensionabstrait} and establish a bound for the Hausdorff dimension of the hypograph of the decoration $g$.

\begin{corollary} \label{C:dimensionhyp} Assume  \eqref{E:fam0} and \eqref{E:seriehypbis}, and recall from \eqref{eq:partial:0:T} that $\partial_0 T$ denotes the subset  of leaves of $T$ on which the decoration $g$ vanishes.
 Then, the Hausdorff dimension of the hypograph of $g$ can be bounded by
$$ \mathrm{dim}_{H}( \mathrm{Hyp}(g)) \leq \mathrm{dim}_{H}(\partial_0 T) \vee 2,$$
 where $\mathrm{dim}_H(\partial_0 T)$ stands for the Hausdorff dimension of $\partial_0 T$ equipped with the restriction of $d_T$.
\end{corollary}
\begin{proof} Recall from Notation \ref{N:subtrees}, that for every $n\geq 1$, the notation $T^n$ stands for the subtree of $T$ built by gluing segments up to generation $n$ only. Both $T$ and $T^n$ are rooted at $\rho$ and
$T$ is decorated with the usc function $g$, while $T^n$ is decorated with  $g^n$ which verifies $g^n\leq g$ on $T^n$.  The hypograph $ \mathrm{Hyp}(g^n)$ 
can be constructed by gluing a countable family of hypographs of decorated segments, each having Hausdorff dimension smaller than $2$, and therefore 
$$ \mathrm{dim}_{H}\Big( \bigcup_{n\geq 0} \mathrm{Hyp}(g^n)\Big)\leq 2.$$
Since the family  $(\|f_u\ )_{u\in \mathbb{U}}$ is null, we infer that $ \bigcup_{n\geq 0} \mathrm{Hyp}(g^n)=  \bigcup_{n\geq 0} \mathrm{Hyp}(g_{T^n})$, where $g_{T^n}$ stands for the restriction of $g$ to $T^n$. Finally, recall also from the proof of Lemma \ref{lem:dimensionabstrait} that $T\setminus \bigcup_{n\geq 0} T^n$ is a subset of $\partial_0 T$, so we can write
$$ \mathrm{Hyp}(g) = \left( \bigcup_{n\geq 0} \mathrm{Hyp}(g_{\mid T^n})\right) \cup \Big(\partial_0 T \times\{0\}\Big),$$
and the desired bound follows.
\end{proof}

Our goal now is to define a metric to compare general hypographs. In this direction,    consider  two usc functions $g:K\to \R_+$ and $g': K'\to \R_+$ with respective domains $K$ and $K'$ in $ {\mathcal K}(Y)$. 
Their hypographs $\mathrm{Hyp}(g)$ and $\mathrm{Hyp}(g')$ are thus two elements of ${\mathcal K}(Y\times \R_+)$, and we can define 
$$\d_{\mathrm {Hyp}}(g,g')\coloneq \d_{\mathrm {Haus}}\left( \mathrm{Hyp}(g), \mathrm{Hyp}(g')\right),$$
where in the right-hand side, $ \d_{\mathrm {Haus}}$ denotes the Hausdorff distance between two compact subsets in $Y\times \R_+$.
In words, the hypograph distance $\d_{\mathrm {Hyp}}$ between two usc functions $g:K\to \R_+$ and $g': K'\to \R_+$ is at most $\varepsilon >0$ if and only if for every $x\in K$ we can find $x'\in K'$  such that $d_Y(x,x')\leq \varepsilon$ and $g(x)\leq g'(x')+\varepsilon$, and \textit{vice versa} when the roles of $g$ and $g'$ are permuted. The hypograph convergence for sequences of usc functions can be characterized as a kind of pointwise convergence, see \cite{beer1981natural} for details. This convergence, and the related convergence of epigraphs are also called $\Gamma$-convergence which  is widely used in optimization \cite{braides2002gamma}.

\begin{remark} \label{R:HypSko}
 In the case when $Y=\R$, it is natural to compare the distance between two functions in the sense of hypographs with other notions {\it\`a la} Skorokhod; see \cite[Chapter VI]{JS03} for background. Specifically, take $a<z$ and $a'<z'$, and let 
$g:[a,z]\to \R_+$ and $g': [a',z']\to \R_+$
be two usc functions. It is then readily checked that  for any  increasing bijection $\beta: [a,z]\to [a',z']$, one has
$$\d_{\mathrm {Hyp}}(g,g') \leq \max\big\{\| g-g'\circ \beta\| , \| \beta-\mathrm{Id}_{[a,z]}\|\big\},$$
where we use the notation $\|\cdot \|$ for the supremum norm on the space of bounded functions on $[a,z]$ and $\mathrm{Id}_{[a,z]}$ for the identity function on $[a,z]$. As a consequence, convergence in the sense of Skorokhod for sequences of rcll (right-continuous with left limits) functions defined on compact intervals entails convergence of the sequence of usc versions\footnote{ Recall that if $g:[a,z]\to \R_+$ is a rcll function, its usc version is the function defined by  $\check{g}(t)\coloneq\max\{g(t-),g(t)\}$ for $t\in[a,z]$, with the convention $g(a-)=g(a)$. } in the sense of hypographs. Of course, the converse fails. For instance, consider $g_n:[0,1]\to \R_+$ given by $g_n(x)=1+\cos(nx)$ for $n\geq 1$. Then the sequence $(g_n)_{n\geq 1}$ converges as $n\to \infty$ to the constant function $2$ on $[0,1]$ for the  distance $\d_{\mathrm {Hyp}}$, but does not converge in the sense of Skorokhod.  \end{remark}

Let us now extend the main nomenclature and notations introduced for real trees to the setting of general compact metric spaces. In this direction, consider $K \in {\mathcal K}(Y)$, and note that this set is naturally endowed with the distance $d_K$  induced by the restriction of $d_Y$ to $K$. We say that $(K, d_K)$ is rooted if it is equipped with a distinguished point $\rho$ in $K$, also referred to as the root.  Similarly, the compact space  $(K,d_K)$ is said  decorated (resp. measured)  if it is equipped with an usc function $g: K\to \R_+$ (resp. with a finite Borel measure $\upnu$ on $K$).We use short hand notations $\texttt{K} \coloneq (K, d_K, \rho, g)$ for a decorated compact space and $\mathbf{K}=(K, d_K, \rho,g, \upnu)$ for a measured decorated compact space. Obviously, a decorated compact space can always be seen as a measured decorated compact space once equipped with the null measure, so that we directly develop the formalism in the more general case.
We write $\mathbb{H}_{m}(Y)$ for the space of measured decorated compact spaces on $Y$, that is formally
$$\mathbb{H}_{m}(Y)\coloneq\bigsqcup_{K\in \mathcal K(Y)}
\left\{(K, d_K, \rho,g, \upnu): \rho \in K, g: K\to \R_+ \text{ usc},  \text{ and }  \upnu \in \mathcal{M}^{f}(K)\right\},$$
where  $\mathcal{M}^{f}(K)$ denotes the space of finite Borel measures on $K$.
We endow ${\mathbb{H}_{m}}(Y)$ with a natural metric defined  for any $\mathbf{K}, \mathbf{K}' \in \mathbb{H}_{m}(Y)$ in the obvious notation  by
$$\d_{{\mathbb{H}_{m}}(Y)} (\mathbf{K},\mathbf{K}')  \coloneq   d_Y\big(\rho, \rho'\big)\vee \d_{\mathrm {Hyp}}\left(g,g'\right) \vee \d_{\mathrm{Prok}}\big(\upnu,\upnu^\prime\big) ,$$
where $\d_{\mathrm{Prok}}$ stands for the Prokhorov distance on $\mathcal{M}^{f}(Y)$.

\begin{proposition} \label{p:hphypo} The space  ${\mathbb{H}_{m}}(Y)$ of measured decorated compact spaces on a  Polish space $(Y,d_Y)$ and equipped with the distance $\d_{{\mathbb{H}_{m}}(Y)}$ is also Polish. 
\end{proposition}
\begin{proof} 
The heart of the argument is the observation, which is doubtless well-known, that the space of hypographs of usc functions with compact domains in $Y$ is closed in ${\mathcal K}(Y\times \R_+)$. More precisely, we claim that if  $(g_n: K_n\to \R_+)_{n\geq 1}$ is a sequence of usc functions  with $K_n\in{\mathcal K}(Y)$, such that  the sequence of hypographs $\left(\mathrm{Hyp}(g_n)\right)_{n\geq 1}$ converges  as $n\to \infty$ to some $H$ for the Hausdorff distance  in ${\mathcal K}(Y\times \R_+)$, then the sequence $(K_n)_{n\geq 1}$ converges to some $K$  for the Hausdorff distance in ${\mathcal K}(Y)$. Furthermore, there exists  an usc function  $g:K\to \R_+$ such that $H=\mathrm{Hyp}(g)$, and hence the sequence of usc functions $(g_n)_{n\geq 1}$ converges to $g$ as $n\to \infty$ for the hypograph distance 
$\d_{\mathrm{Hyp}}$. 
 
  Indeed,  the projection on the first component $p_1: Y\times \R_+ \to Y$  is a $1$-Lipschitz map. 
 Set $K=p_1(H)$ for the image of $H$ by $p_1$, so $K$ is a nonempty compact set. We now also regard $p_1$ as a map from ${\mathcal K}(Y\times \R_+)$ to ${\mathcal K}(Y)$; plainly this still is a $1$-Lipschitz map, and therefore 
$K_n=p_1(\mathrm{Hyp}(g_n))$ converges to $K$ for the Hausdorff distance in ${\mathcal K}(Y)$.
 Then take any $x\in K$ and $r\geq 0$ such that  $(x,r)\in H$. There exists some sequence $(x_n, r_n)_{n\geq 1}$ with $(x_n,r_n)\in \mathrm{Hyp}(g_n)$ which converges to $(x,r)$ in $Y\times \R_+$.   Since $g_n(x_n)\geq r_n$, the segment $\{x_n\}\times [0,r_n]$ belongs to $\mathrm{Hyp}(g_n)$, and therefore $\{x\}\times [0,r]\subset H$. 
 We set $g(x)=\sup\{r\geq 0: \{x\}\times [0,r]\subset H\}$. Since $H$ is compact, $(x,g(x))\in H$, and we conclude that $H=\{(x,r): x\in K\text{ and } 0\leq r \leq g(x)\}$.
Now using that $H$ is closed, we infer that the function $g: K\to \R_+$ is usc, and $H=\mathrm{Hyp}(g)$.
 
 Next, let $E$ denote the space of rooted measured compact space on $(Y\times \R_+, d_{Y\times \R_+})$  endowed with the Hausdorff--Prokhorov  distance, say $\d_{\mathrm{HP}(E)}$. Identifying the usc function $g$ with its hypograph $\mathrm{Hyp}(g)$, the root $\rho$ with $(\rho,0)\in \mathrm{Hyp}(g)$ and similarly $\upnu$ with a measure $\bar \upnu$ on $Y\times \R_+$ supported by the  base $K\times \{0\}$ of  $\mathrm{Hyp}(g)$, yields a natural isometric embedding, say 
$$\Phi: \big({\mathbb{H}_{m}}(Y),\d_{{\mathbb{H}_{m}}(Y)}\big) \hookrightarrow \big(E, \d_{\mathrm{HP}(E)}\big).$$
Since it is well-known that $(E, \d_{\mathrm{HP}(E)})$ is a Polish space, all that is needed to establish Proposition~\ref{p:hphypo} is to verify that the image $\Phi({\mathbb{H}_{m}}(Y))$ of ${\mathbb{H}_{m}}(Y)$ under this embedding is closed in $(E, \d_{\mathrm{HP}(E)})$. So, consider a sequence $\big(\mathbf{K}_n\big)_{n\geq 1}$ in 
 ${\mathbb{H}_{m}}(Y)$ such that the embedded sequence $\left(\Phi (\mathbf{K}_n)\right)_{n\geq 1}$
converges in $(E, \d_{\mathrm{HP}(E)})$, say to $(H,d_H,\bar{\rho},\bar\upnu)$, where $H\in \mathcal{K}(Y\times \R_+)$, 
$d_H$ denotes the restriction of  $d_{Y\times \R_+}$ to $H$,  $\bar{\rho}\in Y\times \R_+$ and $\bar\upnu\in \mathcal{M}^f(Y\times \R_+)$. We have just seen above that $K_n$ converges to $K$
for the Hausdorff distance on $\mathcal{K}(Y)$ and that $H=\mathrm{Hyp}(g)$
for some usc function $g:K\to \R_+$. Plainly, since $\bar{\rho}$ is the limit of $(\rho_n,0)$ in $Y\times \R_+$ and
$\rho_n\in K_n$ for all $n\geq 1$, we have $\bar{\rho}=(\rho,0)$ for some $\rho\in K$. Similarly, by the Portemanteau theorem, the measure $\bar\upnu$ must be supported by $K\times \{0\}$, and can thus be identified as a finite measure $\upnu$ on $K$. Last, $d_H$ coincides with the restriction of $d_{K\times \R_+}$ to $H=\mathrm{Hyp}(g)$. 
Putting the pieces together, we get $(H,d_H,\bar{\rho},  \bar\upnu)=
\Phi(\mathbf{K})$ for some $\mathbf{K}=(K, d_K, \rho,g,\upnu)\in \mathbb{H}_{m}(Y)$, as we wanted to check.
\end{proof}

 Roughly speaking, we are only interested in the general  structure induced by a  decorated compact space rather than by a specific realization. More precisely, two
 decorated compact spaces, say
  $\mathbf{K}=(K, d_K, \rho, g,\upnu) \in {\mathbb{H}_{m}}(Y)$ and $\mathbf{K}'=(K', d_{K'},\rho', g^\prime,\upnu' )\in {\mathbb{H}_{m}}(Y')$ have the same structure and then are viewed as equivalent (or isomorphic) if there exists  a bijective isometry $\phi : (K, d_K) \to (K',d_{K'})$ with inverse denoted by $\phi^{-1}$,  such that $\rho'=\phi(\rho)$, $g^\prime=g\circ\phi^{-1}$   and $\upnu'=\upnu\circ\phi^{-1}$ is the pushed forward of $\upnu$ by $\phi$. 
  We then simply write $\mathbf{K}\approx \mathbf{K}'$.   
  Note that the underlying Polish spaces $(Y,d_Y)$ and $(Y', d_{Y'})$ play no role in this definition and we may simply take 
  $(Y,d_Y)=(K,d_K)$ and $(Y',d_{Y'})=(K', d_{K'})$ for definitiveness.  Let us also make some comments concerning decorated trees and the general construction presented in Section~\ref{sec:1.2}. First, note that due to the 4-point criterion, it is clear that a decorated real tree can only be isomorphic to another decorated real tree. Moreover, we emphasize that the gluing construction described in Section~\ref{sec:1.2} is not unique in the sense that  different building blocks, $((f_u)_{u\in \mathbb{U}}, (t_u)_{u\in \mathbb{U}^*})$  can produce isomorphic trees. For instance, one may apply a bijective isometry (for the graph distance) $\varsigma: \U\to \U$ that fixes the root. More complex rearrangements are also possible; for example, during a birth event, a mother particle might switch identities with one of its daughters. Such modifications, called \textbf{bifurcations}, are studied in depth in the self-similar Markov setting in Chapter \ref{chap:spinal:deco}. Informally, each bifurcation provides a different way to decompose the associated decorated tree in building blocks.
  
  \begin{figure}[!h]
   \begin{center}
   \includegraphics[width=13cm]{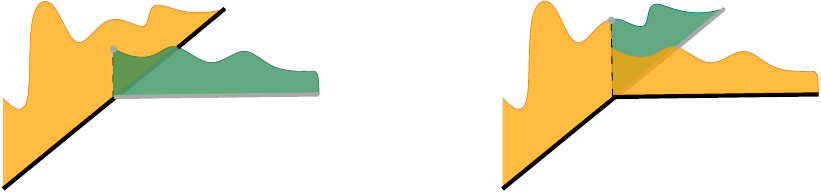}
   \caption{Illustration of a bifurcation event where the identity of two particles is exchanged during a birth event. Obviously the underlying decorated tree is unchanged. \label{fig:bifurcation}}
   \end{center}
   \end{figure}

Our goal now is to define a notion of distance between equivalence classes. In this direction, we denote the set of all equivalence classes of measured decorated   compact spaces by $\mathbb{H}_{m}$. 
 Our goal now is to endow $\mathbb{H}_{m}$ with a natural distance and make it a Polish space. In this direction, recall that
   the Gromov--Hausdorff--Prokhorov  distance between two rooted measured  compact spaces  $(K, d_K,\rho, \upnu)$
 and  $(K', d_{K'},\rho', \upnu')$ is defined as 
 \begin{align*}
 \d_{\mathrm{GHP}} \big(&(K, d_K,\rho, \upnu),(K', d_{K'},\rho', \upnu')\big) \\
&\coloneq \inf_{Y,\phi,\psi} \d_{\mathrm{HP}(Y)} \Big( \big(\phi(K), d_{\phi(K)}, \phi(\rho),\upnu\circ \phi^{-1} \big),  \big(\psi(K^\prime), d_{\psi(K^\prime)}, \psi(\rho^\prime),\upnu^\prime\circ \psi^{-1} \big)\Big),
 \end{align*}
 where  in the right-hand side, $\d_{\mathrm{HP}(Y)}$ stands for the Hausdorff--Prokhorov distance between  rooted measured  compact spaces in $Y$,  the infimum is over all the metric spaces $(Y,d_Y)$ and all the isometric embeddings $\phi: (K, d_K) \hookrightarrow (Y,d_Y)$ and $\psi: (K', d_{K'})\hookrightarrow (Y,d_Y)$, and  $\upnu\circ \phi^{-1}$ is the pushforward measure of $\upnu$ by $\phi$. The quantity  $\d_{\mathrm{GHP}} \big((K, d_K,\rho, \upnu),(K', d_{K'},\rho', \upnu')\big)$ is null if and only if there exists an isomorphism between $(K, d_K,\rho, \upnu)$ and $(K', d_{K'},\rho', \upnu')$, and $\d_{\mathrm{GHP}}$ defines a distance in the space of equivalence classes of rooted measured  compact spaces, see \cite{abraham2013note}.
 
This incites us to make the following definition. First,  for any measured decorated compact space $\mathbf{K}=(\texttt{K},\upnu)=(K,d_K,\rho,g,\upnu)$ and any isometric embedding $\phi: (K, d_K)\hookrightarrow (Y,d_Y)$, we define  
 $$g\circ \phi^{-1} : \phi(K) \to  \R_+\quad, \quad g\circ \phi^{-1}(\phi(x))\coloneq g(x) \text{ for all }x\in K. $$
 Plainly, $\phi(K)$ is a compact subset of $Y$ and 
 $g\circ \phi^{-1}$ is usc. We further  write
 $$ \phi(\mathbf{K})\coloneq \left( \phi(\texttt{K}), \upnu\circ \phi^{-1}\right)=\big(\phi(K), d_{\phi(K)}, \phi(\rho),g\circ \phi^{-1}, \upnu\circ \phi^{-1} \big).$$
For every $\mathbf{K}$ and $\mathbf{K}^{\prime}$ measured decorated compact spaces,  we then  set
 \begin{equation}\label{def:d:H}
\d_{\mathbb{H}_{m}}\left( \mathbf{K},\mathbf{K}'\right) \coloneq 
  \inf_{Y,\phi,\psi}   \d_{\mathbb{H}_{m}(Y)}\big( \phi(\mathbf{K}), \psi(\mathbf{K}')\big),
  \end{equation}
where again the infimum is taken  over all the Polish spaces $(Y,d_Y)$ and all the isometric embeddings $\phi: (K, d_K) \hookrightarrow (Y,d_Y)$ and $\psi: (K', d_{K'})\hookrightarrow  (Y,d_Y)$. 
We  point out that in the special case where $g$ and $g^\prime$ are identically zero on their respective domains, the distances  $\d_{\mathbb{H}_{m}}\big((K,d_{K},\rho,0,\upnu), (K^\prime,d_{K^\prime},\rho^\prime,0,\upnu^\prime)\big)$ and $\d_{\mathrm{GHP}} ((K,d_{K},\rho,\upnu), (K^\prime,d_{K^\prime},\rho^\prime,\upnu^\prime))$ coincide. Of course $\d_{\mathbb{H}_{m}}\left( \mathbf{K},\mathbf{K}^\prime\right)$ is invariant by isomorphisms, and therefore 
 can be viewed as a function $\mathbb{H}_{m}\times \mathbb{H}_{m}\to \R_+$ which we still denote by  $\d_{\mathbb{H}_{m}}$ for simplicity. Our goal now is to establish the following result.

\begin{theorem}\label{theo:Polish}
The map $\d_{\mathbb{H}_{m}}:\mathbb{H}_{m}\times \mathbb{H}_{m}\to \R_+$ defines a distance on $\mathbb{H}_{m}$ and   $(\mathbb{H}_{m},\d_{\mathbb{H}_{m}})$ is Polish.
\end{theorem}

Results of this nature for various classes of decorated compact metric spaces have appeared in numerous contexts in probability theory, and the techniques used to prove them have become standard. For a general framework establishing these types of results, see \cite{khezeli2018unified}. For completeness, we provide the full proof; and to prepare it we begin with a technical lemma.  Roughly speaking, it states that  given an arbitrary sequence $(\mathbf{K}_n)_{n\geq 1}$ of equivalence classes in $\mathbb{H}_{m}$, we can always find some Polish space $(Z,d_Z)$ and representatives of the equivalence classes  in $\mathbb{H}_{m}(Z)$,
such that the distance in $\mathbb{H}_{m}(Z)$ between the representatives of two consecutive equivalence classes  of the sequence is never significantly  larger than the distance between the equivalence classes $(\mathbf{K}_n)_{n\geq 1}$ in $\mathbb{H}_{m}$. Here is the formal statement. 

\begin{lemma}\label{lem:topo}
Let $(\eps_{n})_{n\geq 1}$ be a sequence in $\mathbb{R}_+$ and $(\mathbf{K}_n)_{n\geq 1}\coloneq(K_n,d_{K_n}, \rho_n,g_n, \upnu^n)_{n\geq 1}$ a sequence of measured decorated compact spaces with
$$\d_{\mathbb{H}_{m}}(\mathbf{K}_n,\mathbf{K}_{n+1})< \eps_n, \qquad \text{for all }n\geq 1.$$ Then there exist a Polish space $(Z,d_Z)$ and isometric embeddings $\varphi_1,\phi_2, \ldots$ respectively from $K_1, K_2, \ldots$ into $Z$ such that
$$ \d_{\mathbb{H}_{m}(Z)}\big(\phi_n(\mathbf{K}_n),\phi_{n+1}(\mathbf{K}_{n+1})\big)<\eps_n, \qquad \text{for all }n\geq 1.$$
\end{lemma}
Let us mention that similar results have already appeared in the literature for different variants of the Gromov--Hausdorff--Prokhorov distance. In this work, we  adapt the proof of Lemma 5.7 in \cite{greven2009convergence}, which establishes the analog of our Lemma \ref{lem:topo} for the Gromov--Prokhorov distance.
\begin{proof}
Without loss of generality, we may suppose that the compact spaces $K_n$ 
are pairwise disjoint. First,  by definition,  we can find for every $n\geq 1$ a Polish space $(Y_n,d_{Y_n })$ and two isometric embeddings  $\varphi_n:K_n \hookrightarrow Y_n$ and  $\psi_n:K_{n+1} \hookrightarrow Y_n$ such that:
$$\eta_n\coloneq \d_{\mathbb{H}_{m}(Y_n)}\big(\varphi_n(\mathbf{K}_n),\psi_n(\mathbf{K}_{n+1})\big)<\eps_n.$$
Then, we introduce the disjoint union 
$Z\coloneq\bigsqcup_{n\geq 1} K_n,$ 
and we endow $Z$ with the metric $d_Z$ defined as the largest distance that coincides with $d_{K_n}$ on $K_n$ for each $n\geq 1$ and 
such that 
$$d_Z(x,y)=d_{Y_n}\big(\varphi_n(x),\psi_n(y)\big) + (\esp_n-\eta_n)/2, \qquad \text{ for }x\in K_n \text{ and }y\in K_{n+1},$$ where the term $(\esp_n-\eta_n)/2$ ensures that $d_Z(x,y)>0$ when $x\in K_n$ and $y\in K_{n+1}$. Specifically, for $x_j\in K_j$ and $x_{j+k}\in K_{j+k}$ with $1\leq j,k$, we have 
$$d_Z(x_j,x_{j+k}) \coloneq \inf\Big\{\sum_{\ell=j+1}^{j+k} d_Z(x_{\ell-1},x_{\ell}): x_{\ell} \in K_{\ell} \text{ for every } \ell = j+1, \ldots, k-1\Big\}.$$

 Since the spaces $K_n$ are separable,  $(Z,d_Z)$ is also separable. By a slight abuse of notation, we still write $(Z,d_Z)$ for its completion, which is a Polish space.
We  claim that $(Z,d_Z)$ and the sequence $(\phi_n)_{n\geq 1}$, where  $\phi_n: K_n \hookrightarrow Z$ is a canonical embedding, satisfy the conclusion of the statement.

To begin with, note from the very definition of $d_Z$ that, for every $n\geq 1$, the distance in $Z$ between the roots of $K_n$ and $K_{n+1}$ satisfies
$$d_Z(\rho_n, \rho_{n+1}) = d_{Y_n}\big(\varphi_n(\rho_n),\psi_n(\rho_{n+1})\big) + (\esp_n-\eta_n)/2 \leq \eta_n+ (\esp_n-\eta_n)/2 < \eps_n. $$
Then, recall that the Prokhorov distance in $Y_n$ between $\upnu^n\circ \varphi_n^{-1}$ and $\upnu^{n+1}\circ \psi_n^{-1}$ is at most $\eta_n$. 
Take any $\eta > \eta_n$ and let $A$ be an arbitrary Borel subset of $Y_n$. Writing $A^{\eta}$ for the $\eta$-neighborhood of $A$ in $Y_n$, we have
$$\upnu^n\circ \varphi_n^{-1}(A) \leq \upnu^{n+1}\circ \psi_n^{-1}(A^{\eta})+\eta.$$
Let $B$ be an arbitrary Borel subset of $K_n$, which we also view as a subset of $Z$. Take $A=\varphi_n(B)$, so $A\subset Y_n$, and set $\eta'=\eta+ (\esp_n-\eta_n)/2$. One checks readily from the definition that 
$\psi_n^{-1}(A^{\eta}) =\phi_{n+1}^{-1}(B^{\eta'})$, where $B^{\eta'}$ is the $\eta'$-neighborhood of $B$ in $Z$. Thus the last displayed inequality yields 
$$ \upnu^n(B) = \upnu^n\circ \phi_n^{-1}(A) \leq \upnu^{n+1}\circ \psi_n^{-1}(A^{\eta})+\eta =  \upnu^{n+1}\circ \phi_{n+1}^{-1}(B^{\eta'})+\eta.$$
The very same argument also shows that for an arbitrary Borel subset $C$ of $K_{n+1}$,
$$ \upnu^{n+1}(C) \leq  \upnu^{n}\circ \phi_{n}^{-1}(C^{\eta'})+\eta.$$
Since $\eta$ can be chosen arbitrarily close to $\eta_n$, this entails that $$\d_{\mathrm{Prok}}(\upnu^n\circ \phi_n^{-1}, \upnu^{n+1}\circ \phi_{n+1}^{-1}) \leq \eta_n + (\esp_n-\eta_n)/2 < \eps_n.$$
\par 
Finally, we deal with the hypographs  and check that
\begin{equation}\label{eq:hyp:Z}
\d_{\mathrm{Hyp}}\Big(g_n\circ \phi_n^{-1},g_{n+1}\circ \phi_{n+1}^{-1}\Big)<\eps_n. 
\end{equation}
Indeed, for every $(x,r)\in\mathrm{Hyp}(g_n\circ \phi_n^{-1}) $, we have
\begin{align*}
&\inf \Big\{ d_Z(x,y) \vee |r-s|:~(y,s)\in \mathrm{Hyp}(g_{n+1}\circ \phi_{n+1}^{-1})\Big\}\\
&\leq (\eps_n-\eta_n)/2 + \inf \Big\{ d_{Y_n}\big(\phi_n(x),\psi_n(y)\big)\vee |r-s|:~(y,s)\in \mathrm{Hyp}(g_{n+1}) \Big\}.
\end{align*}
The second term in the sum   above is bounded by
$$ \d_{\mathrm{Hyp}}\big(g_n\circ \phi_n^{-1},g_{n+1}\circ \psi_n^{-1}\big) \leq  \d_{\mathbb{H}_{m}(Z)}\Big(\phi_n(\texttt{K}_n, \upnu^n),\psi_n(\texttt{K}_{n+1}, \upnu^{n+1})\Big)
=\eta_n,$$ 
and therefore
$$\sup \left\{ d_{Z\times \R_+} \left((x,r) , \mathrm{Hyp}(g_{n+1}\circ \phi_{n+1}^{-1})\right) : (x,r)\in\mathrm{Hyp}(g_n\circ \phi_n^{-1}) \right\}< \eps_n. $$
The same argument shows that
$$\sup\left\{ d_{Z\times \R_+} \left((y,s) , \mathrm{Hyp}(g_n\circ \phi_n^{-1})\right): (y,s)\in\mathrm{Hyp}(g_{n+1}\circ \phi_{n+1}^{-1}) \right\}< \eps_n. $$
This establishes \eqref{eq:hyp:Z} and completes the proof of the lemma.
\end{proof}
\noindent We can now proceed  with the proof of Theorem \ref{theo:Polish}.
\begin{proof}[Proof Theorem \ref{theo:Polish}] 
The proof  relies heavily on Proposition \ref{p:hphypo} and Lemma \ref{lem:topo}.
We shall use Lemma \ref{lem:topo} to represent measured decorated compact spaces  in the same well-chosen
Polish space $(Z,d^{Z})$. It is convenient for this purpose to let  $Z$ systematically appear as an exponent in the notation,  writing e.g. $\d_{\mathrm{Hyp}}^Z$ for the hypograph distance on $\mathbb{H}_{m}(Z)$, rather than using $Z$ as an index or omitting it like in Lemma \ref{lem:topo} and its proof.

We first establish that $\d_{\mathbb{H}_{m}}$ is a distance. 
Symmetry is clear; we now check the triangle inequality. Let  $\mathbf{K}_1, \mathbf{K}_2, \mathbf{K}_3$ three measured decorated compact spaces and  $\eps_1,\eps_2>0$ such that
$$\d_{\mathbb{H}_{m}}\big(\mathbf{K}_1,\mathbf{K}_2\big)<\eps_1~~~ \text{ and }~~~ \d_{\mathbb{H}_{m}}\big(\mathbf{K}_2,\mathbf{K}_3\big)<\esp_2.$$
By Lemma \ref{lem:topo}, there exist a Polish space $Z$ and isometric embeddings $\phi_1,\phi_2$ and  $\phi_3$, respectively from $K_1, K_2$ and $K_3$ into $Z$, such that
$$ \d_{\mathbb{H}_{m}(Z)}\big(\phi_n(\mathbf{K}_n),\phi_{n+1}(\mathbf{K}_{n+1})\big)<\eps_n, \qquad \text{for }n=1,2.$$
By Proposition \ref{p:hphypo}, $ \d_{\mathbb{H}_{m}(Z)}$ is a distance on $\mathbb{H}_{m}(Z)$ and then
 the triangle inequality gives
$$ \d_{\mathbb{H}_{m}}\big(\mathbf{K}_1,\mathbf{K}_3\big)\leq   \d_{\mathbb{H}_{m}(Z)}\big(\phi_1(\mathbf{K}_1),\phi_{3}(\mathbf{K}_3)\big)<\eps_1+\esp_2.$$
Passing to equivalent classes, we infer that $\d_{\mathbb{H}_{m}}:\mathbb{H}_{m}\times \mathbb{H}_{m}\to \R_+$  satisfies the triangle inequality.

We then check positivity. To this end we need to show that if  $\mathbf{K}$ and  $\mathbf{K}'$ are two measured decorated compact spaces such that $\d_{\mathbb{H}_{m}}(\mathbf{K},\mathbf{K}')=0$, then they have to be equivalent. Let us proceed. 
For every $n\geq 0$,  set  $\mathbf{K}_n=\mathbf{K}$ if $n$ is odd and  $\mathbf{K}_n=\mathbf{K}'$ if $n$ is even. Again by Lemma \ref{lem:topo}, we can find a Polish space $(Z,d^Z)$ and isometric embeddings $\phi_{2n-1}: K \hookrightarrow Z$ and 
$\phi_{2n} : K' \hookrightarrow Z$,  for all $n\geq 1$, such that
$$ \d_{\mathbb{H}_{m}(Z)}\big(\phi_n(\mathbf{K}_n),\phi_{n+1}(\mathbf{K}_{n+1})\big)<2^{-n}.$$
It follows that the sequence $\left(\phi_n(\mathbf{K}_n)\right)_{n\geq 1}$ is Cauchy, and 
by Proposition \ref{p:hphypo},   converges in $\mathbb{H}_{m}(Z)$ to, say, $\mathbf{K}^Z$. Specifying this for odd integers, we get in the obvious notation,
$$ \lim_{n\to \infty} d^Z( \phi_{2n-1}(\rho), \rho^Z) = 0\quad , \quad\lim_{n\to \infty} \d^Z_{\mathrm{Hyp}}(g\circ \phi_{2n-1}^{-1}, g^Z) = 0 \quad , \quad \lim_{n\to \infty} \d^Z_{\mathrm{Prok}}(\upnu\circ \phi_{2n-1}^{-1}, \upnu^Z) = 0. $$

Recall now from the proof of Proposition \ref{p:hphypo} that the convergence of usc functions
 for the  hypograph distance $ \d^Z_{\mathrm{Hyp}}$ entails the convergence of the domains for the Hausdorff distance, so
$$\lim_{n\to \infty} \d^Z_{\mathrm{Haus}}(\phi_{2n-1}(K), K^Z) = 0.$$
This implies that for any $x\in K$, the sequence $(\phi_{2n-1}(x))_{n\geq 0}$ is relatively compact in $Z$. On the other hand, the sequence of isometries $(\phi_{2n-1})_{n\geq 0}$ is of course equicontinuous.
These observations enable us to apply the Arzel\`a-Ascoli theorem, and we infer that
there is a strictly increasing sequence of odd integers $(n_k)_{k\geq 1}$ such that 
 $(\phi_{n_k})_{k\geq 1}$  converges uniformly to an isometric embedding $\phi: K \hookrightarrow Z$. 
 It is now immediate to check that $K^Z= \phi(K)$, $\rho^Z=\phi(\rho)$, $g^Z=g\circ \phi^{-1}$, and $\upnu^Z=\upnu\circ \phi^{-1}$, so $\mathbf{K}$
and $\mathbf{K}^Z$ belong to the same equivalence class in $\mathbb{H}_{m}$. 
The same argument shows that  $\mathbf{K}'$
and $\mathbf{K}^Z$ also belong to the same equivalence class in $\mathbb{H}_{m}$, and establish positivity. 

Finally, we check that the space $(\mathbb{H}_{m}, \d_{\mathbb{H}_{m}})$ is Polish.
Separability should be plain since the set -- of isometry classes -- of  measured decorated    rooted compact spaces  with a finite cardinality and  associated distances, measures and usc functions taking only rational values is dense in $(\mathbb{H}_{m},\d_{\mathbb{H}_{m}})$. 
Completeness is also immediate from  Lemma~\ref{lem:topo} and Proposition~\ref{p:hphypo}. Indeed, if   $(\mathbf{K}_n)_{n\geq 1}$ is a sequence of measured decorated compact spaces such that:
$$\lim_{n\to \infty} \sup_{\ell \geq 1} \d_{\mathbb{H}_{m}}(\mathbf{K}_n, \mathbf{K}_{n+\ell})=0,$$
then Lemma~\ref{lem:topo} enables us to embed $(\mathbf{K}_n)_{n\geq 1}$ into a Cauchy sequence $(\mathbf{K}^Z_n)_{n\geq 1}$  in $(\mathbb{H}_{m}(Z),\d_{\mathbb{H}_{m}(Z)})$ for some Polish space $(Z,d^Z)$. 
We know from  Proposition \ref{p:hphypo} that the latter converges, say to some
$\mathbf{K}^Z\in \mathbb{H}_{m}(Z)$. We now see from the definition \eqref{def:d:H} that the equivalent class of $\mathbf{K}_n$ converges to the equivalent class of  $\mathbf{K}$ in $\mathbb{H}_{m}$ as $n\to \infty$.
 \end{proof}
 
 Finally, we turn back our attention to measured decorated real trees; recall Definition \ref{Def:dectree}.

 \begin{corollary} \label{C:dmrtclosed} The set  $\mathbb{T}_{m}$ of equivalence (up to isomorphisms) classes of measured decorated  compact real trees  is closed in  $(\mathbb{H}_{m}, \d_{\mathbb{H}_{m}})$. \end{corollary}

\begin{proof}  Let $(\mathbf{T}_n)_{n\geq 1}$ be a sequence of measured decorated  compact real trees converging to some $\mathbf{T}$ in $\mathbb{H}_{m}$. We already know that $\mathbf{T}$ must be a decorated  compact space.
We just need to check that the metric space associated to the latter is a real tree, which is immediate, using e.g. the four point condition. See also \cite[Theorem~1]{EPW06}. \end{proof}

In the remainder of this work, we write $\d_{\mathbb{T}_{m}}$ for the restriction of $\d_{\mathbb{H}_{m}}$ to $\mathbb{T}_{m}$   and  equip $\mathbb{T}_{m}$ with the Borel sigma-field.  Let us mention that the class of decorated trees, with a trivial measure, can be viewed
as a closed subclass of $(\mathbb{T}_{m}, \d_{\mathbb{T}_{m}})$. This allow us to identify the set $ \mathbb{T}$ of equivalence of (non-measured) decorated trees  with the closed subset of $ \mathbb{T}_{m}$ of equivalence classes equipped with the null measure. We denote by $ \mathrm{d}_{ \mathbb{T}}$ the induced distance so that $( \mathbb{T}, \mathrm{d}_{ \mathbb{T}})$ is also a Polish space which we similarly equip with the Borel sigma-field. 
We will often abuse  terminology and refer to a decorated real tree $\mathtt{T}=(T,d_T,\rho, g)$ or a measured decorated tree $\mathbf{T}=(T,d_T,\rho, g, \upnu)$ instead of their equivalence class, implicitly identifying a decorated real tree and its equivalence class. 
When doing so, one must consider only notions that are invariant under isomorphisms, such as weighted length measures (if $\varpi\circ g \in L^1(\uplambda_T)$). Other examples of notions well-defined in $\mathbb{T}_{m}$ include the functions that associate to a decorated real tree its diameter, the maximal value of the upper semi-continuous decoration, or its total mass; which are all continuous with respect to $\d_{\mathbb{T}_{m}}$. We also stress that the  map $(T,d_T,\rho, g,\nu)\mapsto (T,d_T,\rho, g)$ that replaces $\nu$ by  the null measure,   is continuous by definition from $(\mathbb{T}_{m}, \d_{\mathbb{T}_{m}})$ to $(\mathbb{T}, \d_{\mathbb{T}})$.  It will also be convenient to write $\mathtt{0}$ (resp.  $\mathbf{0}$) for the element of $\mathbb{T}$ (resp. $ \mathbb{T}_{m}$) corresponding to a degenerate real tree reduced to a singleton with zero decoration (and the null measure).

  In the sequel, we will sometimes have to consider equivalence classes of measured decorated real trees  with marks. A minor difficulty however is that marking a decorated real tree is not  unambiguously  defined for equivalence classes,  and  as a remedy, we need to work with  an extension   of $\mathbb{T}_{m}$ and $ \mathbb{T}$ for marked  (measured) decorated real trees. This extension is straightforward and let us present the measured decorated case. 
 
  Let $\mathbf{T}=( T, d_T,\rho, g, \upnu)$ be a fixed measured decorated tree and  $I$ is a finite or countable set of indices.  For each $i\in I$, consider a point $x_i$ which is either an element of $T$ or is fictitious (i.e. absent), so we think of $(x_i)_{i\in I}$ as a family of points in $T$, possibly indexed by a strict subset of $I$.
  We then say that two pairs $(\mathbf{T},(x_i)_{i\in I})$ and $(\mathbf{T}^\prime,(x_i^\prime)_{i\in I})$ are equivalent and then write
  $\big( \mathbf T, (x_i)_{i\in I}\big) \approx \big(\mathbf T,   (x_i^\prime)_{i\in I}\big ) $,
    if there exists an isometric bijection $\varphi: T\to T'$ which induces an isomorphism between $\mathbf{T}$ and $\mathbf{T}^\prime$, such that furthermore $\varphi(x_i)=x_i^\prime$ for all $i\in I$ (we agree that this identity always holds when $x_i$ and $x'_i$ are both fictitious, and fails if only one of the two is fictitious). We denote the set of such equivalence classes  by $\mathbb{T}_{m}^{I\bullet}$, and simply $\mathbb{T}_{m}^{\bullet}$ when $I=\{1\}$.    In order to extend   the distance $\d_{\mathbb{T}_{m}}$ to the set $\mathbb{T}_{m}^{I\bullet}$ of equivalence classes of decorated real trees with marks, we fix some null family $(a_i)_{i\in I}$ of positive real numbers. 
   First, when  $(\mathbf{T},(x_i)_{i\in I})$  and $(\mathbf{T}^\prime,(x_i^\prime)_{i\in I})$ are two marked decorated real spaces in some Polish space $(Y,\d_Y)$, 
  we set first
  $$\d_{\mathbb{H}_{m}^{I\bullet}(Y)}\big((\mathbf{T}, (x_i)_{i\in I}),(\mathbf{T}^\prime,(x_i^\prime)_{i\in I})\big)\coloneqq \d_{\mathbb{H}_{m}(Y)}\big(\mathbf{T},\mathbf{T}^\prime)\vee  \Big(\sup \limits_{i\in I}\Big(\mathrm{d}_{Y}\big( x_i,  x'_i\big)\wedge a_i\Big)\Big), $$
  where  we agree that the distance $\mathrm{d}_{Y}\big( x_i,  x'_i\big)$ is $0$ when $x_i$ and $ x'_i$ are both fictitious, and infinite when only one is fictitious. 
  Then, in the general case, one defines
     \begin{align}\label{eq:def:d:bullet}
&\d_{\mathbb{T}_{m}^{I\bullet}}\big((\mathbf{T}, (x_i)_{i\in I}),(\mathbf{T}^\prime,(x_i^\prime)_{i\in I})\big)\nonumber\\ &\coloneq\mathop{\inf\limits_{\varphi:T \hookrightarrow Y}}_
{\varphi^\prime:T^\prime \hookrightarrow Y}\d_{\mathbb{H}_{m}(Y)}\big(\varphi(\mathbf{T}),\varphi^\prime(\mathbf{T}^\prime)\big)\vee  \Big(\sup \limits_{i\in I}\Big(\mathrm{d}_{Y}\big( \varphi(x_i),  \varphi^{\prime}(x_i)\big)\wedge a_i\Big)\Big), 
\end{align}
where  the infimum is over all the Polish spaces $(Y, d_Y)$ and all the isometric embeddings $\varphi:T \hookrightarrow Y$ and $\varphi^\prime:T' \hookrightarrow\ Y$.

By definition, $\d_{\mathbb{T}_{m}^{I\bullet}}$ induces a well defined function from $\mathbb{T}_{m}^{I\bullet}\times \mathbb{T}_{m}^{I\bullet}$ to $\mathbb{R}_+$, which we still denote by $\d_{\mathbb{T}_{m}^{I \bullet}}$ by a slight abuse of notation. It is straightforward to check that the space $(\mathbb{T}_{m}^{I\bullet}, \d_{\mathbb{T}_{m}^{I\bullet}})$ is  Polish, and we stress that the resulting topology on $\mathbb{T}_{m}^{I\bullet}$ does not depend on the specific choice of  $(a_i)_{i\in I}$. Specifically, Lemma~\ref{lem:topo} can be extended to this context using the exact same proof, and then the proof of Theorem  \ref{theo:Polish} can easily be adapted; we leave the details to the reader. Just as in the unmarked case, we write $ \mathbb{T}^{I\bullet}$ for the set of equivalence classes of marked, non-measured, decorated compact trees, and when no confusion is possible,  we identify a marked (measured) decorated compact tree  with its equivalence class in $\mathbb{T}_{m}^{I\bullet}$ or $\mathbb{T}^{I\bullet}$  and we use  the notation $(\mathbf{T},(x_i)_{i\in I})$  and $(\mathtt{T},(x_i)_{i\in I})$. Later on, when we will deal with random decorated trees, we shall always implicitly work on the canonical space $ \mathbb{T}, \mathbb{T}_{m}, \mathbb{T}^{I \bullet}$ or $ \mathbb{T}_{m}^{I\bullet}$ equipped with their Borel sigma-fields, which will be endowed with different laws.

 As an important example, we point out that the gluing operator defined in Section~\ref{sec:1.1}, which uses marks to specify the locations where gluing takes place on the base tree, can be made compatible with isomorphisms. Here is a formal statement. 

\begin{lemma}\label{L:gluingcompatible}  Let $\mathtt{T}'$ be a decorated real tree with marks $(x_i)_{i\in I}$ and $(\mathtt{T}_i)_{i\in I}$ 
 a family of decorated real  trees such that  the domains of   $(\mathtt{T}_i)_{i\in I}$ are pairwise disjoint and also disjoint from the one of $\mathtt{T}'$. Assume that \eqref{E:nullheight} and \eqref{E:nulldeco}  hold.
 
 Let also  $\breve{\mathtt{T}}'$ be another decorated real tree with marks $(\breve x_i)_{i\in I}$ and  $(\breve{\mathtt{T}}_i)_{i\in I}$
  another family of decorated real  trees such that the  domains of   $(\breve{\mathtt{T}}_i)_{i\in I}$ are pairwise disjoint and also disjoint from the one of $\breve{\mathtt{T}}'$.  Suppose that 
 $$\big( \mathtt{T}', (x_i)_{i\in I}\big) \approx \big(\breve{\mathtt{T}}', (\breve x_i)_{i\in I}\big) \quad \text{and} \quad \mathtt{T}_i\approx \breve{\mathtt{T}}_i\quad \text{for all }i\in I.$$
Then \eqref{E:nullheight} and \eqref{E:nulldeco}  also hold for the second family of decorated real trees, and we have
$$ \mathrm{Gluing}\Big(\big( \mathtt{T}', (x_i)_{i\in I}\big), \big(\mathtt{T}_i\big)_{i\in I}\Big) \approx 
 \mathrm{Gluing}\Big( \big(\breve{\mathtt{T}}^\prime, (\breve x_i)_{i\in I}\big), \big(\breve{\mathtt{T}}_i\big)_{i\in I}\Big).
 $$
\end{lemma} 
\begin{proof} With the obvious notation,  let $\varphi': T'\to \breve T'$ and $\varphi_i: T_i\to \breve T_i$, for $i\in I$, denote bijective isomorphisms underlying  the assumptions of the Lemma. By gluing these  isomorphisms
at the marks $(x_i)_{i\in I}$ in an obvious way, we can construct a function $\varphi$ from the domain of $\mathrm{Gluing}\big(\big( \mathtt{T}', (x_i)_{i\in I}\big), \big(\mathtt{T}_i\big)_{i\in I}\big)$ to the one of $ \mathrm{Gluing}\big( \big(\breve{\mathtt{T}}^\prime, (\breve x_i)_{i\in I}\big), \big(\breve{\mathtt{T}}_i\big)_{i\in I}\big)$. One readily checks that $\varphi$ is in turn a  bijective isomorphism for the glued decorated real trees.
\end{proof}
As a consequence, we  can henceforth view the gluing operator as a map from a sub-domain of $\mathbb{T}^{I\bullet} \times(\mathbb{T})^I$ to $\mathbb{T}$; and as usual we keep the same notation for the latter as for the former. By convention, we can extend the definition  when conditions \eqref{E:nullheight} and \eqref{E:nulldeco}  are not fulfilled by simply setting $ \mathrm{Gluing}\big( \big(\mathtt{T}^\prime, ( x_i)_{i\in I}\big), \big(\mathtt{T}_i\big)_{i\in I}\big)=\mathtt{0}.$

 \section{Comments and bibliographical notes}
 
The formalism presented here is chiefly inspired by the recursive construction of continuum random trees performed  by  Rembart and Winkel  \cite{rembart2018recursive} in terms of so-called strings of beads and by the general gluing of metric spaces along points performed by Senizergues, see \cite[Section 2]{senizergues2022growing}. Note in particular that gluing of real trees along points is a folklore operation in the literature on random real trees, see e.g. \cite[Section 2.4]{abraham2014exit}, which can at least  be traced back to Aldous \cite{Ald91} and the famous stick breaking construction of the Brownian Continuum Random Tree. 

Let us conclude with a discussion, for readers familiar with ``Gromov-type" topologies, on possible alternative topologies for decorated trees, and explain our choice. 
It is possible to define a complete separable topology on the set of compact metric spaces endowed with a continuous function taking values in a fixed Polish space and with a controlled regularity (e.g. a Lipschitz condition), see \cite[Section 3]{barlow2017subsequential} and especially Remark 3.2 there. However,  since we aim to consider functions that are typically discontinuous on trees, finding an adaptation ``\`a la Skorokhod''  for trees decorated with rcll functions in the spirit of \cite{khezeli2018unified} or \cite{piotrowiak2011dynamics} seemed complicated.  Another  approach inspired by the Brownian and Lévy snake constructions of Duquesne \& Le Gall \cite{DLG02}, involves viewing the label of a point $x$ of a ``decorated'' tree as a rcll path: $$\zeta_{x} :  [0,d_{T}(\rho,x)) \to \mathbb{R},$$ representing the entire history of the ``decoration'' from the root of $T$ to $x$. By doing so, the labeling becomes  Lipschitz over the real tree (for the appropriate topology of space of paths), allowing us to define a Gromov-type topology for such structures. The  resulting ``snake" topology\footnote{Snake trajectories \cite{ALG15,DLG02,riera2024structure} encode more structure than just the tree and the labels. Roughly speaking, they also encode a canonical contour function of the underlying tree. The analogy with snake trajectories here is that we consider the entire ancestral path as a label. This is also reminiscent of the notion of historical processes in the theory of superprocesses \cite{dynkin1991path}.}  differs from the topology discussed in this chapter: snake convergence roughly corresponds to the convergence of the (measured) rooted tree in the Gromov-Hausdorff-Prokhorov sense, combined with a Skorokhod convergence of the decorations along the macroscopic  branches, but it does not cover the decorations near the leaves. Conversely, our hypograph convergence uniformly controls the supremum of our function over any  non trivial interval of the tree, but  the topology is weaker ``in the interior of branches'', see Remark \ref{R:HypSko}. We chose the hypograph-type topology because, besides the fact that functions are defined on leaves, it also works well with the gluing operation and is particularly suited for studying Markov properties and spinal decompositions in the context of decorated trees, see Chapters  \ref{chap:markov} and \ref{chap:spinal:deco}.

\chapter{Branching processes with real types} \label{chap:generalBP}

Motivated by the evolution in continuous time of a population of individuals, Jagers \cite{jagers1989general} introduced 
general branching processes as Markov random fields indexed by the Ulam tree.  Roughly speaking, each individual is labeled according to its ancestral lineage, and  receives at birth a type  in some abstract space which determines the statistics of a so-called  life career.  The lifetime, the reproduction process and further traits of an individual (which typically may evolve with the age of the individual) are all viewed as measurable functions of its life career. 
The branching property requests that conditionally on their types, the life careers of individuals at generation $n$ are independent, and also independent of the life careers of individuals from the previous generations.

For the applications we have in mind, here types are positive real numbers\footnote{In Chapter \ref{chap:spinal:deco}, we shall also consider a bit more generally models with distinguished individuals, so that formally the type of an individual has then two components: a positive real number together with a distinction or an absence of distinction.}. Any  individual may beget infinitely many children\footnote{This situation is usually excluded in the setting of general branching processes, as it would be awkward from a biological point of view. It is nonetheless relevant when populations and types have rather a geometric interpretation as we shall see later on.}, but only finitely many with types greater than $\varepsilon$ for every $\varepsilon >0$. 
We apply results of the preceding chapter to construct, under fairly general assumptions,  a random real tree that encodes the evolution of such a branching process, and such that lengths of branches correspond to time durations. We decorate the latter with a nonnegative usc function
to represent some trait of every individual which may vary with the age. We further endow the resulting tree with different weighted length measures, and,  assuming further the existence of a harmonic function on the space of types, 
we also define a natural Borel measure carried by the set of leaves of the genealogical tree.  

The framework is made simpler when one further requests self-similarity and the Markov property in time.  More precisely, it is well-known feature that general branching processes are Markovian when viewed as processes indexed by generations; however the Markov property in the time variable fails, except in the very special case when lifetimes have an exponential distribution and reproduction processes are homogeneous Poisson processes. Nevertheless, the fact that we do not only consider the evolution of a branching population with types, but also endow individuals with a random decoration, will enable us to retrieve the Markov property in the time variable for a large family of self-similar models. Self-similar Markov branching processes  are introduced in Section \ref{sec:2.3} and  will be shown  to arise as scaling limits of a variety of discrete Markov branching models  in Part II. 
The construction relies on the so-called Lamperti transformation, which is classically applied to connect real L\'evy processes to positive self-similar Markov processes, and we shall also provide  some necessary background in this setting.

\section{General branching processes as random decorated real trees} \label{sec:2.1}

In this section, we shall apply general results from the preceding chapter and define the decorated real tree that depicts  a general branching process endowed   with some random decoration for individuals.  In this direction, we will have to specify the distribution of the building blocks, and to start with, we recall more precisely how a general branching process with types in $(0,\infty)$ can be  constructed.

We call  \textbf{decoration-reproduction process} a pair $(f,\eta)$ with $f: [0,z]\to \R_+$ a random rcll function  on a random interval $[0,z]$, and $\eta=\eta(\d t, \d y)$ a point process on $(0,z]\times (0,\infty)$. Strictly speaking, we mean by this that $\eta$ is a point process on $\R_+\times (0,\infty)$ such that 
\begin{equation} \label{E:pasnaissances}
 \eta(\{0\}\times (0,\infty))=0 \quad \text{and} \quad  \eta((z,\infty)\times (0,\infty))=0, \qquad \text{almost surely.}
 \end{equation}
 See Figure \ref{F:feta} below. We refer to  $f$, respectively to $\eta$,  as the decoration process, respectively the reproduction process. 
  We should think of $z$ as the lifetime of an individual,  of $f(t)$ as some trait of this individual at age $t$, and of $\eta$ as a reproduction process, in the sense that each atom of $\eta$, say $(t,y)$,  is interpreted as a birth event of a child of type $y$ occurring when the individual has reached age $t$ (notice that an individual can produce multiple offspring at the same time). The requirement \eqref{E:pasnaissances} means that the individual cannot beget at birth nor after its death which is a natural restriction (however it may happen that an individual produces offsprings at the exact time when it dies). Needless to say,  decoration and reproduction are generally not independent.

 \begin{figure}[!h]
 \begin{center}
 \includegraphics[width=7cm]{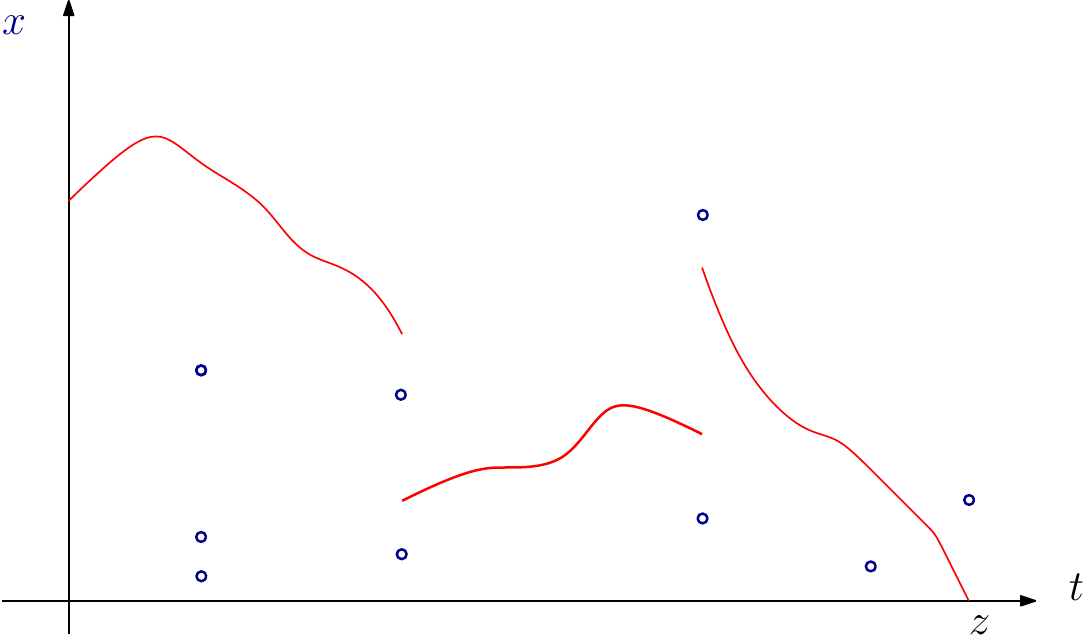}
 \caption{  Decoration and reproduction point process $(f,\eta)$. The centers of the small circles represent the locations of the atoms of $\eta$; observe that several atoms may share the same time-coordinate, and also that the time-coordinates of these atoms may or not be the time of a jump of the decoration. In the above case, there is a reproduction at the exact time of death of the individual.}
 \label{F:feta}
 \end{center}
 \end{figure}

 We then  call a  \textbf{decoration-reproduction kernel} a family of probability laws $(P_x)_{x>0}$, where for every  $x>0$, the distribution $P_x$ is the law  of a random decoration-reproduction process  $(f, \eta)$ for an individual with type $x$. 
 We shall implicitly assume that the map $x\mapsto P_x$ is a  measurable\footnote{This requirement will be automatically satisfied in the self-similar case.}
  function of the type $x$ and we use the notation $E_x$ for the mathematical expectation under $P_x$. 
 We always assume that for every $\varepsilon>0$,  any individual has finitely many children with type greater than $\varepsilon$, almost surely for $P_x$ for all $x>0$. In other words,  the types of the progeny of an individual always form a null family, even though the total progeny may be infinite. 
 
We now describe formally the construction of a general branching process with a given decoration-reproduction kernel, using 
 the Ulam  tree $ \mathbb{U}$ to encode the genealogy of individuals. Let us  first explain how to assign  to every individual a type, a reproduction process, and a random decoration.
By assigning types to individuals, we mean a random process $(\chi(u))_{ u\in \U}$ on the Ulam tree, so that $\chi(u)$ is the type of the individual labeled by $u$. 
 In this direction, it is convenient to add $0$ to the space of types;  the type $0$ will be assigned to fictitious individuals,
that  do not appear in the branching process, but nonetheless  have to be represented by some vertex $u\in \U$ for definitiveness.

Given the type of the ancestor, say $ \chi(\varnothing) = x>0$, we pick a pair $(f_{\varnothing}, \eta_{\varnothing})$ with law $P_x$ as above. We enumerate the atoms $(t_1,y_1), (t_2,y_2), \ldots$  of $\eta_{\varnothing}$ using some deterministic rule, for instance the co-lexicographical order\footnote{That is in non-increasing order of the types, and ties are broken using  the decreasing order of times.}. In the case when $\eta_{\varnothing}$ has only finitely many atoms, we complete with fictitious individuals of the form $(\dagger, 0)$ to get an infinite sequence. We then  set $\chi(j)=y_j$ for every individual $j\in\N$ at the first generation. 
Once the types of the individuals (fictitious or not) at the first generation have been assigned, we iterate the construction using the branching property. That is, conditionally on $(f_{\varnothing}, \eta_{\varnothing})$, we consider a sequence $(f_1, \eta_1), (f_2, \eta_2), \ldots $ of independent pairs distributed according to $P_{y_1}, P_{y_2},  \ldots$, where 
$y_j=\chi(j)$, and, for definitiveness, $P_0$ denotes the law of the trivial pair $(f,\eta)$ corresponding to  $z=0$, $f(0)=0$ and $\eta=0$. The construction by iteration for the next generations should now be obvious, using independent decoration-reproduction processes for different individuals. We stress that the type $\chi(u)$ of an individual at any generation $|u|\geq 1$ is determined by the reproduction process $\eta_{u-}$ of its parent $u-$. 

We write $\P_x$ for the probability law of the family of decoration-reproduction processes $\left( f_u, \eta_u\right)_{u\in \U}$ which results when the ancestor $\varnothing$ has the type $  x>0$; then we naturally write also $\E_x$ for the mathematical expectation under $\P_x$.
Let us spell out in this setting what we shall refer to as the \textbf{branching property} in the sequel. For any initial  type $x>0$ and for every  $n\geq 1$, the families $\left( f_u, \eta_u\right)_{|u|< n}$ and $\left( f_v, \eta_v\right)_{|v|\geq n}$ are conditionally independent under the law $\P_x$ given the family $\left(\chi(w)\right)_{|w|=n}$ of the types at  generation $n$. 
Specifically, the conditional law 
of each subfamily $\left( f_{wv}, \eta_{wv}\right)_{v\in \U}$ for a vertex $w$ at generation $n$ is $\P_{\chi(w)}$, and to different vertices at generation $n$  correspond  conditionally independent subfamilies. The proof of this branching property is immediate from the construction by a recursive argument.

Every vertex $u\in \U$ has now been assigned not only a type $\chi(u)$, but also a decoration ${f}_u: [0,z_u]\to \R_+$, and a reproduction process $\eta_u$.  
The latter encodes the age $t_{uj}$ of the individual labeled by $u$  at which its $j$-th child is born, and the type $y_{uj}=\chi(uj)$ of this child, for every $j\in \N$. We have therefore the  building blocks needed for the gluing construction of Section \ref{sec:1.2}, namely the families $(f_u)_{u\in \U}$ and $(t_{u})_{u\in \U^*}$:

  We stress that even though the decoration and reproduction processes above have been defined as random variables  under $\P_x$, the setting still makes  sense more generally for arbitrary (deterministic)  
  family $(f_u,\eta_u)_{u\in \mathbb{U}}$ of decoration-reproduction pairs and we can always consider their associated building blocks $(f_u)_{u\in \U}$ and $(t_{ui})_{ui\in \U^*}$ . 
  
 \begin{definition}[Property  \hypertarget{prop:P}{$(\mathcal{P})$}] \label{property:def:P} 
 We say that such a family
 $$(f_u,\eta_u)_{u\in \U}$$
verifies Property~{\hypersetup{linkcolor=black}\hyperlink{prop:P}{$(\mathcal{P})$}}  if the associated building blocks $(f_u)_{u\in \U}$ and  $(t_{ui})_{ui\in \U^*}$ satisfy the requirements  \eqref{E:fam0} and \eqref{E:seriehypbis}. 
  \end{definition}  
When Property~{\hypersetup{linkcolor=black}\hyperlink{prop:P}{$(\mathcal{P})$}}  is satisfied,  we write generically ${\texttt T}=(T,d_T,\rho,g)$ for the decorated real tree which then stems from an application of Theorem \ref{T:recolinfty}.
For the sake of notational simplicity, we still write $\texttt{T}$ for the respective  equivalence class up to isomorphisms (see Section \ref{sec:1.3}), which enables us to view the latter as  random variables with values in  the Polish space
 $ \mathbb{T}$ of equivalence classes of (non-measured) decorated compact trees; see Corollary \ref{C:dmrtclosed} and the discussion below it.  In this direction, we point out that  the function which maps a family  $\left(f_u,\eta_u\right)_{u\in \U}$ satisfying Property~{\hypersetup{linkcolor=black}\hyperlink{prop:P}{$(\mathcal{P})$}} into a decorated real tree ${\texttt T}=(T,d_T,\rho,g)$  is measurable when the set of such families indexed by $\U$  is equipped with the cylindrical sigma-algebra.\footnote{ In short, measurability is straightforward when we restrict our attention to families  having only finitely many non-fictitious elements, and the general setting follows by approximation, we leave details  to scrupulous readers familiar with Gromov-Hausdorff-Prokhorov topologies.} It may be tempting to  think of $\texttt{T}$ as  the genealogical tree of the general branching process, and then the decoration $g$ encodes some trait of individuals which may evolve with age. We stress however  that if $\texttt{T}^\prime$ is another random decorated real tree which is a.s. isomorphic to $\texttt{T}$,  then one cannot fully recover the general branching process from $\texttt{T}'$ as the precise genealogy of individuals may have been lost\footnote{ Recall the discussion around Figure \ref{fig:bifurcation}.}. Nonetheless $\texttt{T}'$ does keep track of  evolution as time passes  of the point process that records the values of traits of individuals in the population at any given time, which is sufficient for many applications. In general, we use the notation $\P$ to refer to distributions on families of decoration-reproduction processes, and we use $\mathbb{Q}$  instead  for distributions on the space of equivalence relations of decorated trees $\mathbb{T}$.

Our goal now is to introduce simple assumptions on the decoration-reproduction kernel $({P}_x)_{x>0}$ to ensure that Property~{\hypersetup{linkcolor=black}\hyperlink{prop:P}{$(\mathcal{P})$}}  is verified $\P_x$-a.s., for every $x>0$. These assumptions will become even more transparent in the case when the kernel $(P_x)_{x>0}$ is self-similar, as we
will see at the end of the section.

In this direction,  we consider first 
 the total intensity of children of given types which an individual of type $x$  begets, i.e. the measure $\imath_x$ on $(0,\infty)$ defined by
$$\imath_x(B)\coloneq E_x\left( \eta([0,z]\times B)\right), \qquad B\in \mathcal{B}((0,\infty)).$$
 We shall henceforth suppose  the existence of a 
function $\phi: [0,\infty)\to \R_+$ with $\phi(0)=0$ 
and a positive constant $c_{\imath}<1$, such that for all $x>0$:
\begin{equation} \label{E: AssumpH1}
  \int_{(0,\infty)} \phi(y) \imath_x(\d y)  \leq c_{\imath}  \phi(x).
\end{equation}
The function $\phi$ will be referred to as a \textbf{strictly excessive function}.  We shall also suppose 
that there exist $\gamma_0,\alpha>0$ and finite constants $c_z, c_f$ such that:
\begin{equation} \label{E: AssumpH2}
\left .\begin{matrix}
E_x(z^{\gamma_0/\alpha}) &\leq c_z  \phi(x), \\
  E_x(\sup f^{\gamma_0}) & \leq c_f  \phi(x),  \\
    \end{matrix} \quad \right\}
\end{equation}
for all $x>0$. The role of the first assumption \eqref{E: AssumpH1} is enlightened by the following elementary result.

\begin{lemma} \label{L:meanx} Assuming  \eqref{E: AssumpH1},  we have for every $x>0$
$$\E_x\left(  \sum_{u\in \U}\phi( \chi(u))\right) \leq \phi(x)/(1-c_{\imath}).$$
As a consequence,  the family
$\left( \phi( \chi(u))\right)_{ u\in \U}$ is null, $\P_x$-a.s.
\end{lemma}

\begin{proof} By the definition of the intensity measure $\imath_x$, there is the identity
$$ \E_x\left( \sum_{j=1}^{\infty} \phi( \chi(j))\right) = \int_{(0,\infty)} \phi(y) \imath_x(\d y).$$
We deduce from the branching property and \eqref{E: AssumpH1} that for $n \geq 1$ we have
\begin{equation} \label{eq:meanx} 
\E_x\left(  \sum_{|u|=n} \phi( \chi(u))\right) \leq c_{\imath}^{n}  \phi(x),
\end{equation}
and since $c_{\imath}<1$, the first  claim follows. In particular, the family $\left( \phi( \chi(u))\right)_{ u\in \U}$ is summable and \textit{a fortiori} null, $\P_x$-a.s.
\end{proof}
 We are now able to state the main result of this section:

\begin{theorem}\label{T:CMJ}  
Suppose that  \eqref{E: AssumpH2} is fulfilled
 for some strictly excessive function $\phi$ (i.e.  $\phi$ verifies \eqref{E: AssumpH1}) and some exponents $\alpha, \gamma_0>0$.    Then Property~{\hypersetup{linkcolor=black}\hyperlink{prop:P}{$(\mathcal{P})$}} is satisfied  by the family  $(f_u,\eta_u)_{u\in \U}$,   $\P_x$-a.s. 
 for every $x>0$. 
\end{theorem}
\begin{proof}
Fix $x>0$ and recall the notation from Section~\ref{sec:1.2}, and notably that $z_u\geq 0$ is the length of the interval on which $f_u$ is defined.  As a warm-up, we observe from the branching property, 
 \eqref{E: AssumpH1},  the first bound in \eqref{E: AssumpH2}, and Lemma \ref{L:meanx}, that
\begin{equation}\label{E:CMJz_unull} \E_x\left(  \sum_{u\in \U} z_u^{\gamma_0/\alpha} \right) \leq c_z  \E_x\left(  \sum_{u\in \U} \phi(\chi(u))\right) \leq \frac{c_z}{1-c_i} \phi(x)< \infty. 
\end{equation}
In particular,
$$
 \sum_{u\in \U} z_u^{\gamma_0/\alpha}<\infty, \qquad \P_x\text{-a.s.,}
$$
and \textit{a fortiori} $(z_u)_{u\in \U}$ is a null family, $\P_x$-a.s.

We  check similarly that  \eqref{E:seriehyp} (recall that this is a stronger requirement  than  \eqref{E:seriehypbis}) holds $\P_x$-a.s.  We distinguish two cases, depending on whether $\gamma_0< \alpha$ or $\gamma_0\geq \alpha$. 
In the first case where $\gamma_0<\alpha$, we simply write 
$$\left( \sum_{n=0}^{\infty} \sup_{|u|=n}z_u \right)^{\gamma_0/\alpha} \leq  \sum_{n=0}^{\infty} \sup_{|u|=n}z_u^{\gamma_0/\alpha} \leq  \sum_{u\in \U} z_u^{\gamma_0/\alpha}.
$$
We then take the expectation and invoke \eqref{E:CMJz_unull} to see that  \eqref{E:seriehyp} holds $\P_x$-a.s.

In the second case where $\gamma_0\geq \alpha$, we write from Minkovski's inequality
$$\E_x\left( \left( \sum_{n=0}^{\infty} \sup_{|u|=n}z_u \right)^{\gamma_0/\alpha} \right)^{\alpha/\gamma_0} \leq  \sum_{n=0}^{\infty} \E_x\left(  \sup_{|u|=n}z_u^{\gamma_0/\alpha}  \right)^{\alpha/\gamma_0},$$
and then, from  \eqref{E: AssumpH2} and the branching property
$$ \E_x\left(  \sup_{|u|=n}z_u^{\gamma_0/\alpha}  \right) \leq \sum_{|u|=n} \E_x(z_u^{\gamma_0/\alpha}) \leq c_z  \sum_{|u|=n} \E_x\left(   \phi(\chi(u))\right).$$
Thanks to \eqref{eq:meanx}, we can bound  the right-hand side by $c_z  c_{\imath}^{n} \phi(x)$ and infer 
$$\E_x\left( \left( \sum_{n=0}^{\infty} \sup_{|u|=n}z_u \right)^{\gamma_0/\alpha} \right) \leq \left(\sum_{n=0}^{\infty} c_z^{\alpha/\gamma_0} c_{\imath}^{n\alpha/\gamma_0} \phi(x)^{\alpha/\gamma_0}\right)^{\gamma_0/\alpha}
\leq c_z \phi(x) (1-c_{\imath}^{\alpha/\gamma_0})^{-\gamma_0/\alpha},$$
and again  \eqref{E:seriehyp} holds $\P_x$-a.s.

We check likewise that the requirement \eqref{E:fam0}  holds $\P_x$-a.s. 
Namely, we deduce from the branching property, the second bound in \eqref{E: AssumpH2}, and  Lemma \ref{L:meanx},  that 
\begin{equation}\label{theo:cmj;sup:f}
 \sum_{u\in \U} \E_x( \sup f_u^{\gamma_0}) \leq c_f \sum_{u\in \U} \E_x\left( \phi(\chi(u))\right)\leq \frac{c_z}{1-c_i} \phi(x)<\infty,
\end{equation}
 This ensures that $( \sup {f}_u)_{u\in  \mathbb{U}}$ is a null family, $\P_x$-a.s. Since we have already observed that $(z_u)_{u\in \U}$ is a null family, $\P_x$-a.s., the same holds for 
 $( \| f_u\|)_{u\in \mathbb{U}}$ as well. 
\end{proof}

When the conditions of Theorem \ref{T:CMJ} are fulfilled, we can consider the random decorated real tree ${\texttt T}=(T,d_T,\rho,g)$ which then stems from an application of Theorem \ref{T:recolinfty}. The above proof also shows:

 \begin{corollary} \label{C:CMJ}
 Under the assumptions of Theorem \ref{T:CMJ}, the random variables
$$\mathrm{Height}(T)^{\gamma_0/\alpha}\  \mbox{ and } \max_T g^{\gamma_0}\ $$
belong to $L^1(\P_x)$.
\end{corollary}

 The assumptions of Theorem \ref{T:CMJ} are more transparent when the decoration-reproduction kernel $(P_x)_{x>0}$ satisfies a  scaling property. Fix some $\alpha>0$. We start defining, for  any $c>0$ and a (rcll) function 
$f:[0,z]\to \R_+$, the rescaled function
$$f^{(c)}: [0,c^{\alpha}z] \to \R_+, \quad f^{(c)}(t)\coloneq c f(c^{-\alpha}t),$$
and also, for any measure $\eta$ on $[0,z]\times(0,\infty)$,
the rescaled  measure $\eta^{(c)}$ given by the push forward image of $\eta$ by the map
$$[0,z]\times (0,\infty) \to [0,c^{\alpha} z]\times (0,\infty) \ ,\  (t,y) \mapsto (c^{\alpha}t, cy).$$

\begin{definition} \label{D:SSGBP}
 We say that a decoration-reproduction kernel $(P_x)_{x>0}$  is \textbf{self-similar} with exponent $\alpha>0$ if for every $x>0$, the law $P_x$ coincides with the distribution of the rescaled pair 
 $(f^{(x)}, \eta^{(x)})$ under $P_1$. Then we also say that the general branching process is self-similar. 
\end{definition} 

The self-similarity of the decoration-reproduction kernel $(P_x)_{x>0}$ enables us to focus on the type $x=1$. Indeed it is immediate from the self-similarity assumption and the construction of general branching processes that the distribution of the  family of rescaled pairs 
$$\left(f_u^{(x)}, \eta_u^{(x)}\right)_{u\in \U}$$ under $\P_1$ is the same as that of $\left(f_u, \eta_u\right)_{u\in \U}$ under $\P_x$.  
We write for simplicity
$P\coloneq P_1$ and  $\imath\coloneq \imath_1$ for the total intensity of children of given types beget by an individual of type $1$.

We first point out that the self-similarity assumption has a simple
and important consequence for the process $\left(\chi(u)\right)_{u\in \U}$ that assigns types to individuals. 
For every $x>0$, the distribution  under $P_x$ of the family  of the logarithms of types $\left( \log \chi(j) \right)_{j\geq 1}$ of the individuals at the first generation, is identical to the law  under $P$ of the same family shifted by $\log x$. Therefore, if  we 
  consider  the point process on $\R$ induced by the logarithms of types at each generation,
 \begin{equation} \label{E:BRWlog}
 \sum_{|u|=n}\delta_{\log \chi(u)} , \qquad n\geq 0,
 \end{equation}
(implicitly the fictitious individuals with type $0$ are discarded in the sum), then we obtain 
 a \textbf{branching random walk}; see e.g. \cite{shi2015branching}.  This observation has simple consequences regarding the existence of  the strictly excessive function $\phi$ in \eqref{E: AssumpH1} (and later on that of a harmonic function) that we now explain. 
 
 We introduce then the function
\begin{equation}\label{E:Mellin}
{\mathcal M}(\gamma)\coloneq   \int_{(0,\infty)} y^\gamma  \imath(\d y),\quad \gamma\geq 0.
\end{equation}
In words, $\gamma \mapsto {\mathcal M}_{}(\gamma-1)$ is the Mellin transform of the intensity measure $\imath$. 
Note also that the push forward of  $\imath$ by the logarithm function
yields the  reproduction intensity of the branching random walk \eqref{E:BRWlog}, that is the intensity measure of the point process at the first generation,  $\sum_{j=1}^{\infty}\delta_{\log \chi(j)}$. Therefore 
${\mathcal M}_{}$ can also be viewed as the moment generating function of the latter. Except in the degenerate case where the reproduction process $\eta$ is merely a Dirac point mass with a fixed type $P$-a.s., this function is strictly $\log$-convex and takes its values in $  \mathbb{R}\cup \{\infty\}$. By convexity $ \log  {\mathcal M}_{}$ has at most two zeros. Notice that \eqref{E: AssumpH1} is satisfied for $\phi(y) = y^{\gamma_0}$  as soon as $ \mathcal{M}(\gamma_0) <1$.

In the self-similar case, one can  also easily bound from above the Hausdorff dimension of the set of leaves $\partial T$.
\begin{lemma}[Upper-bound on the Hausdorff dimension of leaves]\label{Lem:haus:upper} Suppose that the decoration-reproduction kernel $ (P_{x})_{x >0}$ is self-similar with exponent $\alpha$,
and that the assumptions of Theorem \ref{T:CMJ}  are satisfied with $ \phi(y) = y^{\gamma_0}$, so  $ \mathcal{M}(\gamma_0) < 1$. Assume  furthermore  that there exists $ \omega < \gamma_0$ with $ \mathcal{M}(\omega)=1$.  If we write $\mathtt{T}=(T,d_T, \rho, g)$ for the resulting random decorated tree, then   we have for every $x>0$ that 
$$ \mathrm{dim}_{H}(\partial T)\leq  \omega/\alpha, \quad \P_x\text{-a.s.}, $$
 where $\mathrm{dim}_H(\partial T)$ stands for the Hausdorff dimension of $\partial T$ equipped with the restriction of $d_T$.
\end{lemma}
In particular, combining Lemma \ref{Lem:haus:upper}  with Lemma \ref{lem:dimensionabstrait} and Corollary \ref{C:dimensionhyp},  we infer that $ \mathrm{dim}_{H}(T) \leq \max(1, \omega/\alpha)$
 and  $\mathrm{dim}_{H}(\mathrm{Hyp}(g)) \leq \max(2,\omega/\alpha),$
 $\P_x$-a.s., for every $x>0$.
\begin{proof} By Lemma \ref{lem:dimensionabstrait} and self-similarity, it suffices to prove the bound on the Hausdorff dimension of $\partial T$ under $\P_1$. In this direction, recall the notation $T^{n}$ for the tree obtained by performing the gluing of the first $n$ generation along Ulam's tree. We see $T^{n}$ as a subset of $T$. In the proof of Lemma \ref{lem:dimensionabstrait} we showed that $T^{n}$ is a countable union of segments and as a consequence of the gluing construction $T^{n}\cap \partial T$ is included in the union of the extremities of  those segments. We deduce that $T^{n}\cap \partial T$ is a countable set of points and thus  $\cup_{n\geq 0} T^{n}\cap \partial T$ has Hausdorff dimension~$0$. It remains to show that the Hausdorff dimension of $ \partial^{*} T = \partial T \backslash \cup_{n \geq 0} T^{n}$ is bounded above by $ \omega/\alpha$. To this end, for every $u\in \U$, consider the subtree $T_{u}$  obtained by performing the gluing in the subtree above $u$ in Ulam's tree, seen as a subset of $T$, and write $\mathrm{Diam}(T_u)$ for its diameter in $(T,d_T)$. Then, for every $n\geq 0$,  the collection $\big\{T_{u}:~u\in  \mathbb{N}^{n}\big\} $ is a covering of $\partial^{*} T$ and note that, for every $u\in \mathbb{N}^n$, we have $\mathrm{Diam}(T_u)\leq 2\cdot  \mathrm{Height}(T_u)$. Moreover, for every $u\in  \mathbb{N}^{n}$, the self-similarity of the decoration-reproduction process entails that $$ \mathbb{E}_1\left( \mathrm{Diam}(T_u)^{\gamma/\alpha}\right) \leq 2^{\gamma/\alpha}\cdot  \mathbb{E}_1\left( \chi(u)^{\gamma}\right)\cdot  \mathbb{E}_{1}\left( \mathrm{Height}(T)^{\gamma/\alpha}\right),$$ for every $\gamma \in ( \omega,\gamma_0)$. Now remark that by Corollary \ref{C:CMJ}, the quantity  $ \mathbb{E}_{1}\left( \mathrm{Height}(T)^{\gamma/\alpha}\right)\leq 1+  \mathbb{E}_{1}\left( \mathrm{Height}(T)^{\gamma_0/\alpha}\right)$ is finite  and by log-convexity of $\mathcal{M}$ we also have $ {\mathcal M}(\gamma) <1$,  for every  $\gamma \in ( \omega,\gamma_0)$. Hence,  by the branching property,  we infer that 
$$\mathbb{E}_1\left(\sum\limits_{u\in  \mathbb{N}^{n}} \mathrm{Diam}(T_u)^{\gamma/\alpha}\right)\leq 2^{\gamma/\alpha}\cdot  \mathcal{M}(\gamma)^{n}\cdot \mathbb{E}_{1}\left( \mathrm{Height}(T)^{\gamma/\alpha}\right)\xrightarrow[n\to \infty]{} 0,
 \quad \text{ for } \gamma \in ( \omega,\gamma_0).$$ 
 This proves that the Hausdorff dimension of $\partial^{*} T$ is bounded above by $\gamma/\alpha$, for every $\gamma \in ( \omega,\gamma_0)$,  and so by $ \omega/\alpha$.
\end{proof}

\section{L\'evy, It\^{o},  Lamperti,  and self-similar Markov trees} \label{sec:2.2}
We now proceed with the self-similar Markov case. Our objective is to construct a  general branching process endowed with a decoration that satisfies both self-similarity  and temporal Markov properties.  For this purpose, we first discuss pairs $(X, \eta)$, where $X$ is a self-similar Markov process started from $1$, which we interpret  as a decoration-reproduction process in the sense of Section~\ref{sec:2.1}.
  We then define the kernel $(P_x)_{x>0}$ through a scaling transformation of the law
   $P=P_1$ of $(X,\eta)$; in particular the decoration under $P_x$ has the law of the self-similar Markov process $X$ started from $x$. Note that the type of an individual now coincides with the initial value of its decoration; this was not necessarily so in the more general setting of Section \ref{sec:2.1}. We stress that the decoration-reproduction kernel $(P_x)_{x>0}$ is \textit{de facto} self-similar in the sense of Definition \ref{D:SSGBP}, and that  Markovian aspects will be analyzed in greater details in   Chapter \ref{chap:markov}.
 The framework that we develop here  lies at the heart of the construction of self-similar Markov trees in the next Section \ref{sec:2.3}, which constitute one of the primary objects of this work.
 
Introduce first the space $ \mathcal{S}_1$ of non-increasing sequences ${\mathbf y}=( y_1, y_2,  ... )$ in $[-\infty, \infty)$ with $\lim_{n\to \infty} y_n=-\infty$, and then set
 $ \mathcal{S}\coloneqq [-\infty,\infty)\times \mathcal{S}_1$. Agreeing that  $\log 0=-\infty$ and $\e^{-\infty}=0$,  we view  $\log: [0,\infty) \to [-\infty, \infty)$ as a bijection with reciprocal given by  the exponential function. Transforming an element $(y,{\mathbf y})$ of $\mathcal{S}$ by the exponential function
yields a sequence $(\e^{y}, \e^{y_1}, \e^{y_2}, \ldots)$ in $[0,\infty)$,  whose first term $\e^{y}$ is distinguished and the next ones form a non-increasing sequence that converges to $0$.  Hence, applying the exponential function to each term of the sequence enables us to  endow $ \mathcal{S}$ with the distance induced by the supremum norm on  the space  of real sequences converging to $0$. Then $\mathcal S$ equipped with this distance  is a Polish space. 

We can now introduce the following important terminology.

\begin{definition}\label{D:charquadr} Consider a measure   $\boldsymbol{\Lambda}= \boldsymbol{\Lambda}(\d y, \d \mathbf{y})$  on $\mathcal S$.
Write  $\Lambda_0=\Lambda_0(\d y)$ and $\boldsymbol{\Lambda}_1= \boldsymbol{\Lambda}_1(\d \mathbf{y})$ 
  for its push-forward images by the first projection $(y,{\mathbf y})\mapsto y$ from $\mathcal S$ to $[-\infty,\infty)$ and by the second projection $(y,{\mathbf y})\mapsto {\mathbf y}$ from $\mathcal S$ to $\mathcal S_1$, respectively.   We call $\boldsymbol{\Lambda}$ a \textbf{generalized L\'evy measure} provided that
$$  \int_{[-\infty, \infty)}(1 \wedge y^2)  \Lambda_0( \d  y)< \infty \ \text{ and }\ 
\boldsymbol{\Lambda}_1\left(\big\{\mathbf y \in \mathcal{S}_1 :  \mathrm{e}^{y_{1}}> \varepsilon\big\}\right)<\infty\quad\text{ for all } \varepsilon>0.$$

We then further call $(\sigma^2, \mathrm{a}, \boldsymbol{\Lambda} ; \alpha)$ a \textbf{characteristic quadruplet}, where
 $\alpha >0$ is a \textbf{self-similarity exponent}, $\sigma^2\geq 0$ a \textbf{Gaussian coefficient}, ${\mathrm a}\in\R$ a \textbf{drift coefficient}.
We also refer to $ \mathrm{k}\coloneqq\boldsymbol{\Lambda}(  \{-\infty\} \times  \mathcal{S}_{1}) <\infty$  as the  \textbf{killing rate}. 
\end{definition}

We will next associate a decoration process to a characteristic quadruplet $(\sigma^2, \mathrm{a}, \boldsymbol{\Lambda}; \alpha)$. We rely for this on the classical Lamperti construction of positive self-similar Markov processes 
 from (possibly killed) real Lévy processes and refer the reader to \cite{Lam72} and Chapter 13 in \cite{kyprianou2014fluctuations} for a complete account. 
We introduce a standard Brownian motion $B$ and a Poisson random measure $\mathbf{N}=\mathbf{N}(\d t, \d y, \d \mathbf{y})$ on $[0, \infty) \times \mathcal{S}$ with intensity measure $\d t \boldsymbol{\Lambda}(\d y, \d \mathbf{y})$. If the killing rate $\mathrm{k}$ is strictly positive, then we denote the time coordinate of the first atom of $\mathbf{N}$ belonging to $\{-\infty\} \times \mathcal{S}_{1}$  by $\zeta < \infty$, so that $\zeta$ is then an exponential variable with parameter the killing rate $\mathrm{k}$. If $\mathrm{k}=0$, then we agree that $\zeta=\infty$. We assume that $B$ and $\mathbf{N}$ are independent.

Let us quickly recall the construction of real  Lévy processes by the L\'evy-It\^{o} decomposition. Consider the first projection of ${\mathbf N}$ on $[0,\infty)\times \R$ and  write $N_0=N_0(\d t, \d y)$ for the resulting point process. By the mapping theorem for Poisson random measures, $N_0$ is a Poisson point process with intensity $\mathbf{1}_{y\in \mathbb{R}}\d  t \Lambda_0(\d y)$. Additionally, we introduce the compensated Poisson measure
$$N_0^{(c)}( \d  s,  \d y)\coloneq  N_0(\d  s,  \d y)- \d  s \Lambda_0(\d y).$$
As a consequence of the  conditions  fulfilled  by $\boldsymbol{\Lambda}$, the process $\xi$ defined for $0\leq t < \zeta$ by 
\begin{equation}\label{E:LevyIto}
\xi(t) \coloneqq \sigma B(t) + {\mathrm a}t + \int_{[0,t]\times \R}N_0( \d  s,  \d y)  ~y {\mathbf 1}_{|y|>1} + \int_{[0,t]\times \R} N_0^{(c)}( \d  s,  \d y)~ y {\mathbf 1}_{|y|\leq 1}
\end{equation}
 is a \textbf{L\'evy process}.  Furthermore, in the case $\mathrm{k}>0$,  it will be  convenient to declare that $\xi(t)=-\infty$ for $t\geq \zeta$, so that $\exp(\gamma\xi(t))=0$ whenever $\gamma>0$ and $t\geq \zeta$. In other words,   $\mathrm{k}$ is  the killing rate of the  L\'evy process $\xi$.  For every $t\geq 0$, we have
 \begin{eqnarray}E\big(\exp(\gamma\xi(t)) \big) = E\big(\exp(\gamma\xi(t)), t<\zeta\big) = \exp\big(t\psi(\gamma)\big) \in(0,\infty],   \label{eq:levykhintchine}\end{eqnarray}
 where $\psi$ is known as the \textbf{Laplace exponent} of $\xi$ and given by the \textbf{L\'evy-Khintchine formula}
 \begin{equation} \label{E:LKfor}
 \psi(\gamma) \coloneqq- {\mathrm k} +\frac{1}{2}\sigma^2 \gamma^2+ {\mathrm a} \gamma + \int_{\R^*}   \left( \mathrm{e}^{\gamma y} -1-  \gamma y \mathbf{1}_{|y|\leq 1} \right)\Lambda_0 ( \d y).
 \end{equation}
Slightly  more generally, one can start the Lévy process from any arbitrary $y \in \mathbb{R}$ by translating the entire process, i.e., by considering $y + \xi$.  

We now present the \textbf{Lamperti transformation} which allows  to construct a (positive) \textbf{self-similar Markov process}, for which the acronym pssMp is often used,  from a real L\'evy process $\xi$ and a positive\footnote{The case of  a negative exponent can then be obtained by applying the simple inverse $x\mapsto 1/x$ transformation to a pssMp.}  exponent of self-similarity, $\alpha>0$. We introduce first the exponential functional\footnote{ We shall soon view the random interval $[0,z]$ as the domain of a rcll function; then the notation $z$ has the same interpretation as in the preceding chapters. The reader should therefore not be worried about a possible confusion of notation.}
\begin{equation} \label{E:espilonlamperti}
\upepsilon(t)\coloneqq \int_0^ t \exp(\alpha \xi(s)) \d s \quad \text{for }0\leq t < \zeta, \quad z\coloneqq \upepsilon(\zeta-) =  \int_0^{\zeta} \exp(\alpha \xi(s)) \d s ,
\end{equation}  and note that $\upepsilon: [0,\zeta)\to [0,z)$ is an increasing bijection a.s.
We  then  define $\tau$ as the reciprocal bijection, so that
\begin{equation}\label{E:Lampertime} 
 \int_0^ {\tau(t)} \exp(\alpha \xi(s)) \d s = t, \quad \text{for any }0\leq t <z. 
\end{equation}
Lamperti \cite{Lam72} pointed out that the process  
\begin{equation} \label{E:Lamperproc}
X(t)= \exp\big(\xi(\tau(t))\big),\qquad 0\leq t < z,
\end{equation}
 obtained from the exponential of the L\'evy process by time-substitution based on $\tau$, is both Markovian and self-similar with scaling exponent $\alpha$ (in the literature, 
 one often calls $1/\alpha$ the  Hurst exponent).   More precisely, $X$ starts from $1$, and for every $x>0$, the rescaled process 
\begin{equation} \label{E:rescaling}
xX(x^{-\alpha}t),\qquad \text{for }0\leq t <  x^{\alpha}z,
\end{equation}
is a version of $X$ started from $x$ (that is, the underlying L\'evy process $\xi$ starts from $\log x$). 
Conversely, any pssMp with a positive exponent $\alpha$ can be realized from some real L\'evy process by this transformation. We stress that the lifetime $z$ of $X$ is finite a.s. if and only if either $\zeta<\infty$ a.s. or the L\'evy process drifts to $-\infty$ (i.e. $ \zeta=\infty$ a.s. and $\lim_{t\to \infty}  \xi(t)=-\infty$ a.s.). More precisely, the boundary point $0$ serves as a cemetery state for the self-similar Markov process $X$, it is reached by a jump when $\zeta<\infty$ (i.e. then $X(z-)>0$ a.s.) and continuously when $\xi$ drifts to $-\infty$ (i.e. then $X(z-)=0$ a.s.). By convention we take $X(z)\coloneqq0$, when $z<\infty$, to see  $X$ as a rcll process on the segment $[0,z]$.

For later use, we also observe that for any $c>0$, the L\'evy process, say $\tilde \xi$,
 constructed from the scaled characteristics $( c^2\sigma^{2}, c \mathrm{a}, c \boldsymbol{\Lambda})$ has the same law as $(\xi_{ct})_{t\geq 0}$. Therefore, if we write $\tilde X$ for the pssMp induced by the Lamperti transformation applied to $\tilde \xi$, then there is the identity in distribution 
  \begin{eqnarray} \label{eq:stretched} \big(\tilde{X}_{t} : t \geq 0\big) {\overset{(d)}{=}} \big(X_{c t} : t \geq 0\big).  \end{eqnarray}

 We now turn our attention to the reproduction process $\eta$. For convenience, we use the notation $[0,\zeta]=[0,\infty)$, if $\zeta=\infty$. The reproduction process $\eta$  is going to be defined using  the second projection of ${\mathbf N}$ on $[0,\zeta]\times {\mathcal{S}}_1$ that we denote by $ \mathbf{N}_{1}= \mathbf{N}_{1}(\d t, \d \mathbf{y})$. To this end, we expand each atom of the latter, say  $(s,{\mathbf y})$,  as a sequence $(s,y_{\ell})_{\ell\geq 1}$ in $[0,\zeta]\times [-\infty, \infty)$ and, as a first step, introduce a point process 
on $[0,\zeta]\times \R_+$ by
\begin{equation}\label{Eq:tildeeta:Markov}
\tilde \eta \coloneqq \sum \indset{\{s\leq \zeta\} }\delta_{(s, \exp(\xi(s-)+y_\ell))},
\end{equation}
where the  sum is taken over all the pairs $(s,y_{\ell})$ obtained by developing the atoms $(s,{\mathbf y})$ of ${\mathbf N}_1$, possibly repeated according to their multiplicities, with $\mathbf{y}\neq\{-\infty,-\infty, \dots\}$. In words, $\tilde \eta$ has an atom at $(s,x)$ for some $s \leq \zeta$ and $x>0$ if and only if the Poisson random measure ${\mathbf N}_1$ has an atom at $(s, \mathbf y)$, with $\mathbf{y}\neq\{-\infty,-\infty, \dots\}$,  such that $\log x -\xi(s-)$ is a component of $\mathbf{y}$. As a second step, we perform the Lamperti transformation and consider the push-forward of the measure $\tilde \eta$ by the Lamperti time-change. Namely, recall that $\tau$ is defined by \eqref{E:Lampertime}  as the inverse of the exponential functional $\upepsilon$ (in particular $\upepsilon(\zeta)=z$ is the lifetime of $X$), and set
\begin{equation}\label{Eq:eta:Markov}
\eta \coloneqq \sum \indset{\{\upepsilon(s)\leq z\} }\delta_{(\upepsilon(s), \exp(\xi(s-)+y_\ell))},
\end{equation}
where the same convention for the summation as above applies. 
In particular several atoms may occur at the same time.

We now write $P=P_1$ for the law of the pair  $(X,\eta)$ constructed above and recall that $\alpha>0$ denotes the self--similarity exponent. For every $x>0$, we denote by 
$P_x$ the image of $P$ by the scaling transformation $(s,y) \mapsto (x^{\alpha}s , xy)$  on $\R_+\times \R_+$, 
which, by \eqref{E:rescaling}, transforms the graph of $X$ into that of the version of the self-similar Markov process started from $x$, and the atoms $\delta_{(s,y)}$ of $\eta$ into $\delta_{(x^{\alpha} s, x y )}$. 
We call $(P_x)_{x>0}$  the \textbf{self-similar Markov decoration-reproduction kernel} with characteristic quadruplet $(\sigma^2, \mathrm{a}, \boldsymbol{\Lambda} ; \alpha)$;
the qualifier self-similar is taken in the sense of Definition~\ref{D:SSGBP}.
Note passing by that $\eta(\{0\}\times (0,\infty))=0$ and that,  in absence of killing, i.e. when $ \mathrm{k}=0$, the reproduction process also satisfies $ \eta(\{z\} \times (0, \infty)) = 0$, i.e. that no birth event can happen at the death of an individual.

Our goal now is to  associate a decorated compact real tree to every characteristic quadruplet $(\sigma^2, \mathrm{a}, \boldsymbol{\Lambda} ; \alpha)$ satisfying some additional requirements.  In this direction, we start by computing the moment generating function of the branching random walk associated with the self-similar Markov decoration-reproduction kernel $(P_x)_{x>0}$ with  characteristic quadruplet $(\sigma^2, \mathrm{a}, \boldsymbol{\Lambda} ; \alpha)$. Recall that for simplicity, the index for the starting point $x$ is omitted from the notation when $x=1$, and following \eqref{E:Mellin}, we write
 $$\mathcal M(\gamma)= E\left(\int_{[0,z]\times (0,\infty)} y^\gamma \eta(\dd t, \dd y) \right)=E\left( \sum_{j=1}^\infty \chi(j)^\gamma\right) =\int_{(0,\infty)} \imath(\d y ) y^\gamma, \qquad \gamma\geq 0,$$
for the Mellin transform of the total intensity measure $\imath$ of children of given types. Then, recalling the L\'evy-Khintchine formula \eqref{E:LKfor} for the Laplace exponent $\psi$ of $\xi$, we introduce the quantity
\begin{align}\label{E:cumulant}
\kappa(\gamma)&\coloneqq  \psi(\gamma) + \int_{\mathcal{S}} \boldsymbol{\Lambda}( \d y, \d  \mathbf y ) \left(\sum_{i=1}^{\infty} \e^{\gamma y_i} \right)\nonumber \\
&= \frac{1}{2}\sigma^2 \gamma^2+ {\mathrm a}\gamma + \int_{\mathcal{S}} \boldsymbol{\Lambda}( \d y,\d  \mathbf y )~ \left(\e^{\gamma y}-1- \gamma y \mathbf{1}_{|y|\leq 1} +\sum_{i=1}^{\infty} \e^{\gamma y_i} \right),
\end{align}
where we use again the convention $\mathrm{e}^{-\infty}=0$, so that the possible killing rate $\mathrm{k}$ is incorporated in the integral with respect to $\boldsymbol{\Lambda}$. 
The function $\kappa$ is called the  \textbf{cumulant function} of $(\sigma^2, \mathrm{a}, \boldsymbol{\Lambda} ; \alpha)$. We stress that the cumulant does not depend on the exponent of self-similarity $\alpha$, and that the functions $\psi$ and $\kappa$ are convex with $\psi \leq \kappa$. We will always assume that the cumulant is finite at least at some point $\gamma>0$ (we will actually soon impose more, see the forthcoming Assumption~\ref{A:omega-}). The cumulant function $\kappa$ enables us to compute the Mellin transform  $\mathcal M$. 
\begin{lemma}
\label{L:verCMJ1} The Mellin transform $ \mathcal{M}$ in \eqref{E:Mellin}  is given by 
\begin{equation}\label{eq:comp:M}
\mathcal M(\gamma)= 1-\kappa(\gamma)/\psi(\gamma), \quad \text{whenever } \psi(\gamma)<0.
\end{equation}
\end{lemma}
\begin{proof}
 Recall that $\imath$ is the measure on $(0,\infty)$ that describes the total intensity of the children with given types that an individual with type $1$ begets throughout its life. Note first that the connexion between $\eta$ and $\tilde \eta$ via the Lamperti transformation yields, using the same convention for the summations as in \eqref{Eq:eta:Markov},
$$\sum \indset{\{\upepsilon(s)\leq z\} } \left(X(\upepsilon(s)-)\exp(y_\ell)\right)^\gamma
= \sum \indset{\{s\leq \zeta\} } \exp\left(\gamma (\xi(s-)+y_\ell)\right).
$$
Next, the construction of $\tilde \eta$  in terms of the  random measure $\mathbf N_1$  and the  L\'evy process $\xi$ shows that:
$$
\sum \indset{\{s\leq \zeta\} } \exp(\gamma(\xi(s-)+y_\ell))=\int_{[0,\zeta] \times \mathcal{S}_1} \exp(\gamma\xi(s-)) \left( \sum_{i=1}^{\infty} \exp(\gamma y_i)\right) {\mathbf N_1}( \d  s,  \d  {\mathbf y}).
$$
Taking expectations  using \eqref{eq:levykhintchine}, we get by compensation  whenever $\psi(\gamma)<0$ that
\begin{align*}
\int_{(0,\infty)} y^\gamma  \imath (\d y) &=  E\left(\int_0^\zeta  \exp(\gamma \xi(s-))\d s\right) \left( \int_{\mathcal{S}_1}\sum_{i=1}^{\infty} \e^{\gamma y_i} \boldsymbol{\Lambda}_1(\d {\mathbf y}) \right)  \nonumber \\ 
& =  -\frac{1}{\psi(\gamma)} \left( \int_{\mathcal{S}_1}\sum_{i=1}^{\infty} \e^{\gamma  y_i} \boldsymbol{\Lambda}_1( \d \mathbf y)\right)= 1-\kappa(\gamma)/\psi(\gamma). 
\end{align*} 
\end{proof}
We now have all the ingredients to properly define the self-similar Markov trees (for short, ssMt). We fix a characteristic  quadruplet $(\sigma^2, \mathrm{a}, \boldsymbol{\Lambda} ; \alpha)$  and  we write $(P_x)_{x>0}$ for the associated self-similar Markov decoration-reproduction kernel. We also write $\P_x$ for the distribution of  the family  $(f_u, \eta_u)_{u\in \U}$  of decoration-reproduction processes of individuals in a general branching process governed by this kernel  when the ancestral individual has type $x>0$, that is such that $(f_\varnothing, \eta_\varnothing)$ has the law $P_x$. Let $\kappa$ denote the cumulant function associated by \eqref{E:cumulant} to the characteristic quadruplet  $(\sigma^{2}, \mathrm{a}, \boldsymbol{\Lambda} ; \alpha)$.  We make the following crucial assumption:
\begin{assumption}\label{A:gamma0} We say that a cumulant $\kappa$ is  \textbf{subcritical} if
there exists $\gamma_0>0$ such that
$$ \kappa(\gamma_0)<0.$$
\end{assumption}
Subcriticality of the cumulant is the only assumption on the  characteristic quadruplet needed to define  self-similar Markov trees as random variables with values
in the space $\TT$ of equivalence classes up to isomorphisms of decorated real trees (see Section \ref{sec:1.3}, and more specifically Corollary~\ref{C:dmrtclosed}). 

\label{sec:defssmb}

 \begin{proposition}[Construction of self-similar Markov trees] \label{P:constructionomega-}  Let Assumption~\ref{A:gamma0} be satisfied for some $\gamma_0>0$.   Then the following assertions hold: 
 \begin{itemize}
\item [(i)] The function $\phi(x)=x^{\gamma_0}$ is strictly excessive in the sense of \eqref{E: AssumpH1}. Assumption \ref{E: AssumpH2} is verified, and as a consequence, so does
 Property~{\hypersetup{linkcolor=black}\hyperlink{prop:P}{$(\mathcal{P})$}}, $\P_x$-a.s. for all $x>0$.
 
 \item [(ii)]  The equivalence class (up to isomorphism) of the  random decorated tree $\normalfont{\texttt{T}}=(T,d_T, \rho, g)$ without measures constructed in Theorem \ref{T:recolinfty}  is called a \textbf{self-similar Markov tree} with characteristic quadruplet $(\sigma^2, \mathrm{a}, \boldsymbol{\Lambda} ; \alpha)$. One has 
 $$g(\rho)= \limsup_{\begin{subarray}{c} {r\to \rho} \\  r\neq \rho \end{subarray} } g(r)=x, \qquad \P_x\text{-a.s. for all }x>0.$$
 \end{itemize}
  In this framework,   self-similarity means that for every $x>0$, the law of $\normalfont{\texttt{T}}$ under $\mathbb{P}_{x}$  is identical to that of 
 the rescaled version  $(T, x^{\alpha} \cdot d_T, \rho, x \cdot g)$  under $\P_1$, where the notation $x^{\alpha} \cdot d_T$ is for the distance on $T$ given by 
 $x^{\alpha} \cdot d_T(y,y')= x^{\alpha} d_T(y,y')$, and $x \cdot g$ denotes similarly the decoration with $x \cdot g(y)=xg(y)$. 
\end{proposition}

\begin{remark}[$\mathbb{P}$ and $\mathbb{Q}$] In the following we shall distinguish between the actual random decorated tree $ \texttt{T}$ constructed in the above proposition under the law $ \mathbb{P}_{x}$ and its equivalence class in $ \mathbb{T}$, still denoted by $ \texttt{T}$, whose law will be denoted by $ \Q_{x}$ (and similarly for the forthcoming measured versions). Under $ \mathbb{P}_{x}$ one has access to the construction of $ \texttt{T}$ from the labeled Ulam's tree as presented in the beginning of this section, whereas under $ \Q_{x}$ one can only deal with properties which are measurable in terms of decorated tree. Although this is a subtlety that may be overlooked for most applications, this will be important for the statement of the forthcoming Markov properties in Chapter \ref{chap:markov}. \label{rek:PQ} \end{remark}

\begin{proof} Let us check the assumptions of Theorem \ref{T:CMJ}, for the exponents $\alpha,\gamma_0$ and the function $\phi(y)=y^{\gamma_0}$. In this direction,  recall that $\kappa$ is convex, and since $\psi \leq \kappa$, we have $\psi(\gamma_0)<0$. Hence,  $\kappa(\gamma_0)/\psi(\gamma_0)\in(0,1)$ and we infer  from Lemma \ref{L:verCMJ1} that 
$$ \int_{(0,\infty)}  y^{\gamma_0} \imath(\d y)= \mathcal{M}(\gamma_0)< 1 ,$$ 
so  \eqref{E: AssumpH1} holds  with $\phi(y)=y^{\gamma_0}$. By self-similarity, it suffices now to verify that  $E_1\left(\sup X^{\gamma_0}\right)<\infty$ and  $E_1(z^{\gamma_0/\alpha})<\infty$, and by the Lamperti transformation,  this boils down to establishing that 
$$E_1\left(\sup_{t\geq 0} \exp\left(\gamma_0\xi(t)\right)\right)<\infty \quad \text{ and } \quad E_1\left(\left(\int_{0}^{\zeta}\exp(\alpha \xi(t)) \d t\right)^{\gamma_0/\alpha}\right)<\infty.$$

These kinds of results are part of the folklore of the theory of Lévy processes. For example, see Lemma~3 of Rivero \cite{rivero2012tail} for the second assertion when the killing rate is zero, and Patie and Savov \cite[Theorem 2.18]{PatieSavov} when killing is allowed. Since the display above follows from standard techniques, and because we believe it may be of independent interest, we provide details of the proof in the Appendix; see Lemma \ref{gene:levy:rivero} there.

The assertion that the value of the decoration at the root is $x$, $\P_x$-a.s. follows readily from \eqref{E:fam0} and the fact that $f_\varnothing$ is rcll with $f_{\varnothing}(0)=x$, $\P_x$-a.s.
Last, by construction,  the random decorated tree $\texttt{T}$ plainly inherits the self-similarity property from the kernel $(P_x)_{x>0}$; see Definition \ref{D:SSGBP} and the discussion thereafter.
 \end{proof}
 
   Let us point  right now at an important feature, which is closely related to the discussion around Figure \ref{fig:bifurcation}. We have argued here that one can associate 
a  self-similar Markov tree to any subcritical  characteristic quadruplet, and the characteristic quadruplet then determines the distribution of this self-similar Markov tree. Nonetheless we stress that different characteristic quadruplets may yield self-similar Markov trees with the same distribution, just as in Section \ref{sec:1.2} where different families of building blocks could produce isomorphic decorated trees. In that case, one says that the characteristic quadruplets belong to the same equivalence class of bifurcators. Thus, in short, a characteristic quadruplet determines the law of a self-similar Markov tree, and in the converse direction, a self-similar Markov tree determines an equivalence class of bifurcators, where each such bifurcator corresponds to a unique characteristic quadruplet. 
This matter will be discussed in details in Chapter~\ref{chap:spinal:deco}.

We also point out  from \eqref{eq:stretched} that for any $c>0$, the self-similar Markov tree $( \tilde{T}, d_{\tilde{T}}, \tilde{\rho} , \tilde{g})$  with dilated characteristic quadruplet $( c^{2} \sigma^{2}, c \mathrm{a} , c \boldsymbol{ \Lambda} ; \alpha)$ has the same distribution as 
 $ \left( T, c^{-1}\cdot d_{T}, \rho,g \right)$ where $\texttt{T}= \left( T, d_{T}, \rho,g \right)$ denotes the ssMt with characteristics $( \sigma^{2},\mathrm{a} ,\boldsymbol{ \Lambda} ; \alpha)$.

\section{Weighted length and harmonic measures}\label{sec:2.3}
Now that we have defined  self-similar Markov trees, our next purpose is to endow the latter with certain natural measures. We are  interested in measures compatible with  self-similarity, and which are also consistent with the Markov property of the decoration. This leads us to consider power functions with adequate exponents as weight functions for weighted length measures. To define a natural measure on the leaves, we will require a stronger  assumption, which ensures the existence of a simple harmonic function for the decoration-reproduction kernel. By analogy with the literature on branching random walk, we refer to this  assumption as the first Cramér's condition\footnote{The second Cramér's condition will appear when defining conditional versions of ssMt.}, and its statement and implications will be addressed in the second part of this section. Finally, we point out that the harmonic measure on leaves can also be obtained as a limit of weighted length measures (Proposition \ref{thm:cov:wei:leng}). Throughout this section, we use the notation of Sections~\ref{sec:2.1} and \ref{sec:2.2}. In particular,  we fix a characteristic  quadruplet $(\sigma^2, \mathrm{a}, \boldsymbol{\Lambda} ; \alpha)$  and  we write $(P_x)_{x>0}$ for the associated self-similar Markov decoration-reproduction kernel which has been defined in the preceding section. We also write $\P_x$ for the distribution of  the family  $(f_u, \eta_u)_{u\in \U}$  of decoration-reproduction processes of individuals in a general branching process governed by this kernel  when the ancestral individual has type $x>0$, that is such that $(f_\varnothing, \eta_\varnothing)$ has the law~$P_x$.

\subsection{Weighted length measures}\label{sec:2.3.1}
We start with weighted length measures which are simpler to construct. Recall Section \ref{sec:1:measures} and  focus on the case where $\varpi$ is a power function.

\begin{proposition}[Weighted length measures] \label{prop:lengthmeasures}   Let  Assumption \ref{A:gamma0} hold.  Then for every $x>0$, there is the identity
\begin{equation}\label{eq:mean:length:gamma0}
\mathbb{E}_x\left(\sum_{u\in\U} \int_{0}^{\infty} f_u(t)^{\gamma_0-\alpha} \d t\right)=-\frac{x^{\gamma_0}}{\kappa(\gamma_0)}.
\end{equation}
As a consequence,  for every $\gamma\geq \gamma_0$, the measure 
$\varpi\circ g\cdot \uplambda_T$ induced by
 the weight function\footnote{Our choice of  using $\gamma - \alpha$ instead of $\gamma$ as exponent is related to the Lamperti transformation. As we will see shortly, this will simplify  notation thereafter.} $\varpi(x)=x^{\gamma-\alpha}$ is  finite, $\P_x$-a.s. We denote it by $\uplambda^\gamma$ and  then $\left( T, d_{T}, \rho,g, \uplambda^\gamma \right)$ is a  random measured decorated compact tree.
\end{proposition}

\begin{proof} Let us start  establishing \eqref{eq:mean:length:gamma0}. Thanks to self-similarity, it is enough to treat the case $x=1$. We have 
$$\sum \limits_{u\in \mathbb{U}}\mathbb{E}_1\left( \int_{0}^{z_u}   f_{u}(t)^{\gamma_0-\alpha}\mathrm{d}t\right) =\mathbb{E}_1\left(\sum \limits_{u\in \mathbb{U}} \chi(u)^{\gamma_0}\right) E_1\left(\int_{0}^{z}   f(t)^{\gamma_0-\alpha}\mathrm{d}t\right).$$
By Assumption \ref{A:gamma0}, $\psi(\gamma_0)< \kappa(\gamma_0)<0$, and   Lemma \ref{L:verCMJ1} gives $ \mathcal{M}(\gamma_{0})=1-\kappa(\gamma_{0})/\psi(\gamma_{0})$. Then by  the very construction of the family $(f_u,\eta_u)_{u\in  \mathbb{U}}$, we get $$\E_1\left( \sum_{u\in \mathbb{N}^n} \chi(u)^{\gamma_0}\right) = \Big( 1 -  \kappa(\gamma_0)/\psi(\gamma_0)\Big)^n, \quad \text{ for } n\geq 0.$$ Since $1 -  \kappa(\gamma_0)/\psi(\gamma_0)\in(0,1)$, we deduce that 
$$\mathbb{E}_1\left(\sum \limits_{u\in \mathbb{U}} \chi(u)^{\gamma_0}\right)=\psi(\gamma_0)/\kappa(\gamma_0).$$

Next,  we have by the Lamperti transformation,
\begin{align*}  E_1\left(\int_{0}^{z}  f(t)^{\gamma_0-\alpha}  \mathrm{d}t\right) &=
E_1\left(\int_{0}^{z} \exp\left((\gamma_0-\alpha)\xi(\tau(t))\right)\mathrm{d}t \right)\\
&=E_1\left(\int_{0}^{\infty} \exp\left(\gamma_0\xi(s)\right)\mathrm{d}s \right)\\
&= \int_{0}^{\infty}   \exp\left(\psi(\gamma_0)s\right)\mathrm{d}s\\
&=-\frac{1}{\psi(\gamma_0)},
\end{align*}
where in the second equality we used that $\exp(\gamma\xi(s))=0$, for every $s\geq \zeta$, and  in the ultimate equality, that  $\psi(\gamma_0)\leq \kappa(\gamma_0)<0$. This completes the proof of \eqref{eq:mean:length:gamma0}. 

Let us now deduce the remaining claim. Recalling Section \ref{sec:1:measures} and in particular \eqref{E:newEmassum}, by Proposition~\ref{P:constructionomega-} and self-similarity, it suffices to verify that for every $\gamma>\gamma_0$, we have
$$\sum \limits_{u\in \mathbb{U}} \int_{0}^{z_u}   f_{u}(t)^{\gamma-\alpha}\mathrm{d}t<\infty, \qquad \P_1\text{-a.s.}$$
In this direction,   we notice that for every $\gamma_1>\gamma_0$, 
$$ \sup \limits_{\gamma\in[\gamma_0,\gamma_1]}\sum \limits_{u\in \mathbb{U}} \int_{0}^{z_u}   f_{u}(t)^{\gamma-\alpha}\mathrm{d}t \leq  \sup\Big\{1+f_u(t)^{\gamma_1-\gamma_0}:~u\in \mathbb{U}, t\in [0,z_u]\Big\}\cdot \left(\sum \limits_{u\in \mathbb{U}} \int_{0}^{z_u}   f_{u}(t)^{\gamma_0-\alpha}\d t\right). $$
Recalling  from Proposition \ref{P:constructionomega-} that $\left(\| f_{u}\|\right)_{ u\in  \U}$ is a null family; we infer from \eqref{eq:mean:length:gamma0} that:
$$  \sup \limits_{\gamma\in[\gamma_0,\gamma_1]}\sum \limits_{u\in \mathbb{U}} \int_{0}^{z_u}   f_{u}(t)^{\gamma-\alpha}\mathrm{d}t <\infty,\quad \P_1\text{-a.s.}$$
This completes the proof of the proposition.
\end{proof}
  In the setup of the above proposition, we have  $\uplambda^\gamma(T)<\infty$, $\P_x$-a.s. for every $x>0$,  and we even have 
\begin{equation}\label{mean:uplambda:gamma:T}
\mathbb{E}_x\big(\uplambda^\gamma(T)\big)= -\frac{x^{\gamma}}{\kappa(\gamma)}, \quad 	\mbox{ as soon as $\kappa(\gamma)<0$.}
\end{equation}
 We further stress that the self-similarity asserted in  Proposition~\ref{P:constructionomega-} then extends to these weighted length measures. Specifically,  for every $x>0$,  the distribution of the equivalence class in $ \mathbb{T}_{m}$ of $\left(T, d_T, \rho, g,  \uplambda^{\gamma}\right)$ under $\mathbb{P}_{x}$  is the same as that of $ \left(T, x^{\alpha} \cdot d_T, \rho, x \cdot g, x^{\gamma} \cdot \uplambda^{\gamma}\right)$ under $\P_1$. For this reason,  we say that $\uplambda^\gamma$ is self-similar with exponent $\gamma$, which motivates \textit {a posteriori} our choice for the 
parameter for the weight functions $\varpi$ in Proposition \ref{prop:lengthmeasures}.

\subsection{First Cram\'er's condition and  the harmonic measure on leaves}\label{sec:cramer:cond:1}
In this section, we define the natural measure on leaves for self-similar Markov trees. To do so, we shall rely on the construction presented in Section \ref{sec:1:measures} using an additional family $(m_{u})_{u \in \mathbb{U}}$ verifying \eqref{E:massspreads}.  Those $m_{u}$ will be constructed using a harmonic function for the underlying branching random walk. Specifically, recall from Section \ref{sec:2.1} the notation $(  f_{u}, \eta_{u})_{u \in \mathbb{U}}$ for the general branching process index by Ulam's tree, where the decoration-reproduction kernel $(P_{x})_{x >0}$ is self-similar and associated with the characteristic quadruplet $(\sigma^{2}, \mathrm{a}, \boldsymbol{\Lambda} ; \alpha)$. In this particular case, we have $f_{u}(0) = \chi(u)$, and, as noticed in \eqref{E:BRWlog},  the process 
$$ \sum_{|u|=n}\delta_{\log \chi(u)}, \qquad n \geq 0$$
 is a branching random-walk with moment generating function $ \mathcal{M}$ given by Lemma \ref{L:verCMJ1}. We deduce that if $\omega>0$ is such that $\kappa(\omega)=0$, then the function $h(x) = x ^{\omega}$ is a \textbf{harmonic function} of types, namely that $h(0)=0$ and 
\begin{equation}\label{E:mean-harmonic}
h(x)= \E_x\left(  \sum_{j=1}^{\infty}h(\chi(j)) \right) =\int_{(0,\infty)} h(y) \imath_x(\d y),\qquad \text{for all }x>0. 
\end{equation}

\begin{figure}[!h] 
 \begin{center}
 \includegraphics[height=4cm]{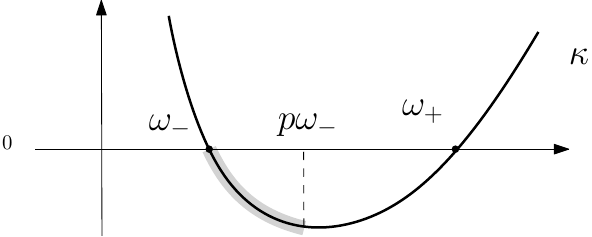}
 \caption {Illustration of the typical shape of the cumulant function $\kappa$. Assumption \ref{A:omega-} requires that $\kappa$ becomes strictly negative after its first zero (like in the gray region).}
 \label{F:kappa}
 \end{center}
 \end{figure}

It follows readily from  the branching property that the process
\begin{equation} \label{E:addmart}
M_n\coloneq \sum_{|u|=n} h(\chi(u)), \qquad n\geq 0,
\end{equation}
is then a  positive martingale, which is usually referred to as an additive martingale; it converges to its terminal value $\P_x$ almost surely. We have more generally for any vertex $u\in \U$
\begin{equation} \label{E:BiggKypu}
m_u := \lim_{n\to \infty} \sum_{|v|=n} h(\chi(uv)), \qquad \text{in } \P_x-a.s.
\end{equation}
Verifying the spread of mass condition \eqref{E:massspreads} is intimately connected to ensuring that the previous convergence also holds in mean. Whether or not the additive martingale \eqref{E:addmart} converges in $L^1(\P)$ is  well understood in the literature on branching random walks, see e.g. \cite{Big77a}. In our setup, this will be implied by the following assumption:

\begin{assumption}[First Cram\'er's condition]\label{A:omega-}
Suppose that there exist $\omega_->0$ and $p\in(1,2]$ such that 
$$ \mathrm{(i)} \quad \kappa(\omega_-)=0, \ \kappa(p\omega_-)< 0, \quad \mbox{ and } \quad \mathrm{(ii)} \quad    \int_{\mathcal{S}_1}\boldsymbol{\Lambda}_1(\d \mathbf y ) ~\left(\sum_{i=1}^{\infty} \e^{ y_i\omega_-}\right)^p< \infty.$$
\end{assumption}
\noindent From now on, under Assumption \ref{A:omega-}, we set:
\begin{equation}\label{eq:omega_+}
\omega_+:=\inf\big\{t >\omega_-:~\kappa(t)\geq 0\big\},
\end{equation}
with the usual convention $\inf\varnothing=\infty$. Observe that, by convexity of the cumulant, we have $p\omega_-<\omega_+$. Furthermore, Assumption \ref{A:omega-} ensures that both $\kappa$ and $\psi$ are negative over the interval $(\omega_-, \omega_+)$.  In particular, Assumption \ref{A:gamma0} holds for every $\gamma_0\in(\omega_-,\omega_+)$ and therefore we can define the associated ssMt as well as the measures $\uplambda^\gamma$,  $\gamma>\omega_-$.  Also, because of convexity,  $\kappa$ possesses a negative right-derivative at $\omega_-$ that we denote for simplicity by $\kappa^{\prime}(\omega_-)$, see  Figure~\ref{F:kappa} for an illustration. 
We are going to use the function $h(x)= x^{\omega_-}$ to define our measure on leaves. For this purpose, we need the following result:

 \begin{lemma}\label{L:MUI} Let Assumption \ref{A:omega-} hold. Then the function $h(x)=x^{\omega_-}$ is harmonic in the sense \eqref{E:mean-harmonic}, and the process
$$M_n(\omega_-)\coloneqq \sum_{u\in \N^n} \ \chi(u)^{\omega_-}, \qquad n\geq 0,$$
is a martingale bounded in $L^p(\P_x)$ for every $x>0$. As a consequence, the spread of mass condition \eqref{E:massspreads} is fulfilled for the family the $(m_{u})_{u \in \mathbb{U}}$ defined by \eqref{E:BiggKypu}.
\end{lemma}

\begin{proof}  Lemma \ref{L:verCMJ1} and Assumption \ref{A:omega-} entail that $\mathcal{M}(\omega_-)=1$, so 
$h(x)=x^{\omega_-}$ is a harmonic function, and we have already seen  that $M_n(\omega_-)$, $n\geq 0$, is then an additive  martingale of a branching random walk with moment generating function $\mathcal M$. Since  $\mathcal{M}(p\omega_-)=1-\kappa(p\omega_-)/\psi(p\omega_-) <1$, it is known that boundedness in $L^p$ will follow provided that $\E(M_1(\omega_-)^p)<\infty$; see  e.g. Theorem 3.1 of Alsmeyer and Kuhlbusch \cite{alsmeyerdouble}. 

In this direction, we write
$$M_1(\omega_-)=\sum_{j=1}^\infty \chi(j)^{\omega_-} = \sum \indset{\{s\leq \zeta\} } \exp\left((\xi(s-)+y_\ell)\omega_-\right),
$$
where the last sum is taken over all the pairs $(s,y_{\ell})$ obtained by expanding the atoms $(s,{\mathbf y})$ of ${\mathbf N}_1$. This shows that $M_1(\omega_-)$ is a so-called L\'evy-type perpetuity in the sense of Iksanov and Mallein \cite{iksanov2019result}, see Section 2 there. Specifying Theorem 3.1 of \cite{iksanov2019result} in our setting, all that we need to verify are,  first
$$E\left(\exp(\xi(1)p\omega_-) \right)=\exp(\psi(p\omega_-))<1,$$
and second, 
$$\int_{\mathcal{S}_1} \boldsymbol{\Lambda}_1(\d \mathbf y )\left( \sum_{i=1}^{\infty} \e^{ y_i\omega_-}\right)^p \indset{\left\{\sum_{i=1}^{\infty} \e^{ y_i\omega_-}>1\right\}} <\infty.\footnote{Theorem 3.1 in \cite{iksanov2019result} assumes that the killing rate is $0$. However, since the proof can be extended  for a positive killing  without modifications, we leave the details to the reader. }
$$
The first requirement holds since  $\psi(p\omega_-)<\kappa(p \omega_-)<0$, and the second is part of Assumption~\ref{A:omega-}.  
 Let us now prove the spread of mass condition \eqref{E:massspreads}. Introduce for any $u\in \U$, 
  \begin{eqnarray} \label{eq:martingaleu}M_n(u)\coloneq \sum_{|v|=n} h(\chi(uv)).
  \end{eqnarray}
Note from the Markov property that conditionally on $\chi(u)=y$, the martingale 
$(M_n(u))_{n\geq 0}$ has the same law as $(M_n)_{n\geq 0}$ under $\P_y$.
Plainly, we have on the one hand, 
 \begin{eqnarray}\label{eq:decomp1}M_{n+1}(u)= \sum_{j=1}^{\infty} M_n(uj), \end{eqnarray}
and on the other hand, there are the convergences a.s. under $\P_x$,
 $$m_u=\lim_{n\to \infty} M_n(u) \quad \text{ and }\quad m_{uj}=\lim_{n\to \infty} M_n(uj) \quad \text{for all }j\geq 1.$$
Hence Fatou's lemma yields 
  \begin{eqnarray} \label{eq:fatou1}m_u \geq \sum_{j=1}^{\infty} m_{uj}, \qquad  \P_x\text{-a.s.}  \end{eqnarray}
Moreover, we also know from the first point of the lemma  that these martingales converge in $L^{1}(\P_x)$, which implies 
$$\E_x(m_u)= \mathbb{E}_{x}(M_0(u)) = \E_x\left(h(\chi(u))\right)  \quad \text{ and }\quad \E_x(m_{uj})= \E_x(M_{0}(uj))=\E_x\left(h(\chi(uj))\right), $$ 
\noindent for all $j\geq 1$. Taking the expectation in \eqref{eq:decomp1} for $n=0$, we deduce that both sides in \eqref{eq:fatou1} have the same mean, and are thus almost surely equal.
 \end{proof}

In the sequel, we shall refer to $(M_n(\omega_-))_{n\geq 0}$ as the \textbf{intrinsic martingale}.   Under Assumption \ref{A:omega-} and using Proposition \ref{P:newPmass},  we can  endow the self-similar Markov tree $\texttt{T}=(T,d_T, \rho,g)$ with the measure induced by the $(m_{u})_{u \in \mathbb{U}}$;  recall  that the latter  is supported  by the set $\partial T$ of leaves of $T$. From now on, we will use the notation $\upmu$ for this measure and  recall that it is referred to as the \textbf{harmonic measure}. Observe that the total mass $\upmu(T)$ of the harmonic measure coincides with the terminal value $m_{\varnothing}$ of the intrinsic martingale. Plainly,  the self-similarity stated in Proposition \ref{P:constructionomega-} extends to $\upmu$ as follows. For any $x>0$, 
 the distribution,  the equivalence class in $ \mathbb{T}_{m}$, of $(T,x^\alpha\cdot d_T, \rho,x\cdot  g,x^{\omega_-}\cdot\upmu)$, under $\P_1$, coincides with that of $(T,d_T, \rho,g,\upmu)$, under $\P_x$. We say that $\upmu$ is self-similar with exponent $\omega_-$;
 note also that 
 \begin{equation}\label{eq:expec:mu:T:x}
 \E_x(\upmu(T))=x^{\omega_-},\quad \text{ for } x>0.
 \end{equation}

We mention that when the second root  $\omega_+$ in \eqref{eq:omega_+} is finite, $M_n(\omega_+)$, $n\geq 0$ is also a martingale. Nonetheless it follows by general results on branching random walks, see \cite{shi2015branching}, that its terminal value is always  $0$  a.s.  Even though the martingale $M_n(\omega_+)$ has also interesting applications (to be developed in a forthcoming work),   it cannot be used to construct a non-degenerate measure on self-similar Markov trees.

Proposition \ref{prop:lengthmeasures}  also allows us to endow  $T$ with the weighted length measures $\uplambda^\gamma$ for any  $\gamma>\omega_-$, as it was discussed at the end of Section \ref{sec:defssmb}. In this direction, one can infer from Lemma \ref{L:MUI}, that for every $\gamma'\leq \omega_-$, one has
$$\int_T g^{\gamma'-\alpha}(v)~ \uplambda_T(\d v)=\sum \limits_{u\in \U} \int_{0}^{z_u}f_{u}(t)^{\gamma'-\alpha} ~\d t=\infty,\quad \P_1\text{-a.s.}$$
In other words, the measures $c\cdot \uplambda^\gamma$, for $\gamma>\omega_-$ and some $c>0$, are the sole self-similar weighted length measures with a finite total mass. Additionally, since all the measures $\upmu$ and $\uplambda^\gamma$,  $\gamma>\omega_-$,  have distinct self-similarity exponents, a nontrivial linear combination cannot be self-similar. Consequently, up to an unimportant factor,  $\upmu$ and $\uplambda^\gamma$, $\gamma>\omega_-$,   are the only self-similar measures consistent with  our framework. They constitute the family of measures  that we will work with in the sequel, and will also naturally appear as scaling limits of discrete models in Part II (see Theorems \ref{T:mainunconds-length} and \ref{T:mainunconds-mass}).

\subsection{The harmonic measure as limit of weighted length measures}\label{sec:2.3.3}

The purpose of this section is to point out   that  the harmonic measure $\upmu$ on leaves can also be obtained as a limit of re-scaled versions of  weighted length measures $\uplambda^{\gamma}$, as $\gamma$ decreases towards $\omega_-$. One important motivation for establishing such a result  stems from the following observation.
Weighted length measures are given intrinsically in terms of the decorated real tree $\texttt T=(T, d_T, \rho, g)$ rather than in terms of a specific construction of $\texttt T$, such as by gluing building blocks in Section \ref{sec:1.2}. 
As a consequence, if $(T', d_{T'}, \rho', g')=\texttt T'$ is another decorated real tree isomorphic to $\texttt T$, then for any weight function $\varpi$ with $\varpi\circ g\in L^{1}(\uplambda_T)$, we have also in the notation of Definition \ref{Def:dectree} that 
$$(T, d_T, \rho, g, \varpi\circ g \cdot \uplambda_T ) \approx  (T', d_{T'}, \rho', g', \varpi\circ g' \cdot \uplambda_{T'}).
$$
At the opposite, the harmonic measure $\upmu$ on a self-similar Markov tree has been defined in terms the family $(m_u)_{u\in \U}$ in \eqref{E:BiggKypu}, and it is not clear \textit{a priori} whether
$\upmu$ can be given directly in terms of $\texttt T$ only. Therefore, it is also unclear whether  the compatibility with isomorphisms that we just stressed for weighted length measures remains valid for $\upmu$. Proposition \ref{thm:cov:wei:leng} below implies that actually, this is indeed the case, and the harmonic measure on leaves is also an intrinsic quantity for self-similar Markov trees.
 
\begin{proposition} \label{thm:cov:wei:leng} Let Assumption \ref{A:omega-} hold. Then there exists a sequence $(\gamma_n)_{n\geq 1}$  with $\gamma_n>\omega_-$ and 
$\lim_{n\to \infty} \gamma_n=\omega_-$, such that $\P_1$-a.s., 
$$\lim\limits_{n\to \infty}  -\kappa(\gamma_n)\cdot \uplambda^{\gamma_n}=\upmu,$$
in the sense of weak convergence for finite measures on $T$. 
\end{proposition}
Recall that $\kappa(\gamma)<0$  for $\gamma\in(\omega_-,\omega_+)$ so that  $-\kappa(\gamma)\cdot \uplambda^\gamma$  is a positive finite measure, and by \eqref{mean:uplambda:gamma:T}, the  factor $-\kappa(\gamma)>0$ is simply chosen so that 
$$\E_1\big(-\kappa(\gamma)\cdot \uplambda^\gamma(T)\big)=\E_1(\upmu(T))= 1.$$
 For technical reasons and to avoid extending the section unnecessarily, we restricted ourselves to convergence along a sequence $(\gamma_n)_{n\geq 1}$  in Proposition \ref{thm:cov:wei:leng}, which suffices for our purpose. We are nonetheless confident that the convergence should hold as $\gamma\downarrow \omega_-$.  
 
The rest of this section is devoted to the proof of Proposition \ref{thm:cov:wei:leng}; we implicitly take Assumption \ref{A:omega-} for granted. As a first step, we establish the convergence of the total mass.

 \begin{lemma}\label{lem:conv:mass} There is the convergence
 $$\lim_{\gamma \downarrow \omega_-}\kappa(\gamma) \uplambda^{\gamma}(T)
 =-\upmu(T),  \qquad \text{in } L^1(\P_1).  $$
 \end{lemma}

 \begin{proof}
 Recall from Lemma \ref{L:verCMJ1} that the Mellin transform of the reproduction intensity is given by  $\mathcal{M}(\gamma)=1-\kappa(\gamma)/\psi(\gamma)\in(0,1]$
  for every $\gamma\in[\omega_-, \omega_+]$. From the branching property,
  $$M_n(\gamma):=\mathcal{M}(\gamma)^{-n}\sum_{u\in \N^n} \ \chi(u)^{\gamma}=\mathcal{M}(\gamma)^{-n}\sum \limits_{u\in \mathbb{N}^n} f_u(0)^{\gamma},\quad n\geq 0,$$
 is then a nonnegative martingale, and we denote its terminal value by $M_\infty(\gamma)$. 
 We are going to establish:
 \begin{itemize}
 \item [$\mathrm{(i)}$]  $M_{\infty}(\gamma)$ converges in $L^{1}(\P_1)$ to $M_{\infty}(\omega_-)$ as $\gamma \downarrow \omega_-$,
\item [$\mathrm{(ii)}$]  $M_{\infty}(\gamma)+\kappa(\gamma)\cdot\uplambda^{\gamma}(T)$ converges in $L^{1}(\P_1)$ to $0$ as $\gamma \downarrow \omega_-$.
\end{itemize}
Combining these two claims implies the lemma, since, by definition $\upmu(T)=M_\infty(\omega_-)$. 
We start by proving the first item which follows by standard branching random walk technics\footnote{Similar results appeared in the theory of branching random walks, see notably \cite{biggins1992uniform}. Unfortunately, we have not been able to find a reference covering our specific situation.}. 

(i) Recall from Lemma \ref{L:MUI} that $M_n(\omega_-)\in L^p(\P_1)$ for some $p>1$ appearing in Assumption \ref{A:omega-} and  every fixed $n\geq 0$. We deduce readily from dominated convergence that 
\begin{equation}\label{eq:conv:L1:M:n}
\lim \limits_{\gamma \downarrow \omega_-}M_n(\gamma)= M_{n}(\omega_-), \quad \text{ in } L^{1}(\P_1).
\end{equation}
The idea now is to control the difference $M_\infty(\gamma)-M_n(\gamma)$ for $\gamma$ close enough of $\omega_-$. In this direction, fix $q\in(1,p\wedge 2)$ and take $\gamma_1\in(\omega_-,p \omega_-/q)$ such that 
 $$c:=\sup_{\gamma\in[\omega_-,\gamma_1]} \mathcal{M}(q \gamma)/\mathcal{M}(\gamma)^q<1.$$
This is indeed feasible since $\kappa(\omega_-)=0$ and $\kappa(p\omega_-)<0$. 

For every $\gamma\in [\omega_-,\gamma_1]$, a direct computation gives:
$$M_1(\gamma)^q \leq  \mathcal{M}(\gamma)^{-q} \left(\sup_{u\in \mathbb{N}} \chi(u)^{(\gamma-\omega_-) q}\right) M_1(\omega_-)^q \leq  \mathcal{M}(\gamma)^{-q}\left(1+M_1(\omega_-)^{p-q} \right) M_1(\omega_-)^q.$$ Since by Lemma \ref{L:MUI}, the variable $M_1(\omega_-)$ is in $L^{p}(\P_1)$ and $\inf_{[\omega_-,\gamma_1]}\mathcal{M}(\gamma)>0$, we deduce from H\"older's inequality that 
$$\sup_{\gamma\in [\omega_-,\gamma_1]}\E_1(M_1(\gamma)^q)\coloneqq K < \infty.$$ 
We can now apply \cite[Lemma 2(i)]{biggins1992uniform}, combined with Jensen inequality, to infer that for every $n\geq 0$ and $\gamma\in [\omega_-,\gamma_1]$, we have
\begin{equation}\label{eq:L1:control:mass:gamma:omega}
\E_1\left(|M_{\infty}(\gamma)-M_n(\gamma)|\right)\leq 2^{3} K^{1/q} c^{n/q}.
\end{equation}
Using \eqref{eq:conv:L1:M:n} and \eqref{eq:L1:control:mass:gamma:omega}, we arrive at
$$\limsup_{\gamma\to \omega_-}\E_1\left(|M_{\infty}(\gamma)-M_{\infty}(\omega_-)|\right) \leq 2^{4} K^{1/q} c^{n/q}.
$$
Finally, taking the  limit as $n\to \infty$ in the right-hand side, we obtain (i). 

(ii) From the very definition of the weighted length measure $\uplambda^{\gamma}$ and the canonical decomposition \eqref{E:newEmassum} of $T$ into line segments indexed by $\U$, we get
for every $\gamma\in(\omega_-,\gamma_1]$,
 $$\kappa(\gamma) \uplambda^{\gamma}(T)=\frac{\kappa(\gamma)}{\psi(\gamma)} \sum_{u\in \U} \chi(u)^{\gamma} A_u(\gamma),$$
 where  for non-fictitious vertex $u$, we write
 $$A_u(\gamma)=\psi(\gamma)\chi(u)^{-\gamma} \int_{0}^{z_u} f_u(t)^{\gamma-\alpha} \d t,$$
 and by convention we take $A_u=0$ if $u$ is fictitious. 
 Next, we infer  from  \eqref{eq:L1:control:mass:gamma:omega} that
\begin{align*}
\E_1\left( \left |M_{\infty}(\gamma) - \frac{\kappa(\gamma)}{\psi(\gamma)}\sum_{n\geq 0}\sum_{u\in\mathbb{N}^n}\chi(u)^{\gamma}\right |\right)&\leq \frac{\kappa(\gamma)}{\psi(\gamma)} \sum_{n\geq 0}\mathcal{M}(\gamma)^{n}\E_1\left(|M_n(\gamma)-M_\infty(\gamma)| \right) \\
&\leq \frac{ 2^{4}K^{1/q}\kappa(\gamma)}{\psi(\gamma)(1-\mathcal{M}(\gamma) c^{1/q})},\end{align*}
and the right-hand side converges to $0$ as $\gamma \downarrow \omega_-$. 
Therefore, it suffices to show that
\begin{equation} \label{Eq:keyconv}
\lim \limits_{\gamma\downarrow \omega_-} \E_1\left( \left |\frac{\kappa(\gamma)}{\psi(\gamma)}\sum_{n\geq 0}\sum_{u\in\mathbb{N}^n}\chi(u)^{\gamma}(A_u(\gamma)+1)\right|\right)=0. 
\end{equation}

In this direction, recalling that $q\in(1,2)$ and using the triangle and then the Jensen inequalities, we bound from above the expectation in \eqref{Eq:keyconv}  by
 \begin{equation}\label{eq:kappa:gamma:sum:chi:u:gamma:q}
\frac{\kappa(\gamma)}{\psi(\gamma)}\sum\limits_{n\geq 0}\E_1\left(\left|\sum_{u\in\mathbb{N}^n}\chi(u)^{\gamma }\left(A_u(\gamma)+1\right) \right|^q\right)^{1/q}.
 \end{equation}
 Now an application of the branching and the self-similarity properties,  combined with the Lamperti transformation, shows that, for $n\geq 0$, under $\P_1$ and conditionally on the types at generation $n$, $(\chi(u))_{u\in \mathbb{N}^n}$, the variables $A_u(\gamma)$, for $u\in \mathbb{N}^n$ with $\chi(u)\neq 0$, are i.i.d with the same law as  
 $$\psi(\gamma)\int_{0}^{\zeta}\exp(\gamma \xi(t)) \d t ,\text{ under }P_1.$$

We know from Tonelli's theorem that
  $$E_1\left(\int_{0}^{\zeta}\exp(\gamma \xi(t)) \d t \right) = \int_{0}^{\infty}E_1\left(\exp(\gamma \xi(t), t<\zeta \right) \d t=-1/\psi(\gamma),$$
  and we can apply the Marcinkiewicz-Zygmund inequality to infer that \eqref{eq:kappa:gamma:sum:chi:u:gamma:q} is bounded above by 
\begin{align*}
&c(q) \frac{\kappa(\gamma)}{\psi(\gamma)}\sum\limits_{n\geq 0}\E_1\left(\sum_{u\in\mathbb{N}^n}\chi(u)^{\gamma q}\right)^{1/q}\cdot E_1\left(\Big|\psi(\gamma)\cdot\int_{0}^{\infty}\exp(\gamma \xi(t)) \d t+1\Big|^q\right)^{1/q}\\
 &= c(q)  \frac{\kappa(\gamma)}{\psi(\gamma)\big(1-\mathcal{M}(q\gamma)^{1/q}\big)} E_1\left(\Big|\psi(\gamma)\int_{0}^{\infty}\exp(\gamma \xi(t)) \d t+1\Big|^q\right)^{1/q},
\end{align*}
where $c(q)<\infty$ is some constant depending on $q$  only.

On the one hand, $c(q)\kappa(\gamma)/\big(\psi(\gamma)\big(1-\mathcal{M}(q\gamma)^{1/q}\big)\big)$ converges to $0$ as $\gamma\downarrow \omega_-$.
On the other hand,  the bound $\e^{y\gamma} \leq \e^{y\omega_-}+ \e^{y \gamma_1}$  for all $y\in \mathbb{R}$ and $\gamma\in[\omega_-,\gamma_1]$ yields
\begin{align*}
&\sup_{\gamma\in[\omega_-,\gamma_1]}E_1\left(\left(\int_{0}^{\infty}\exp(\gamma \xi(t)) \d t\right)^q\right) \\
&\leq 2^{q} E_1\left(\left(\int \exp(\omega_- \xi(t)) \d t\right)^q \right)+2^{q} E_1\left(\left(\int \exp(\gamma_1 \xi(t)) \d t\right)^q \right),
\end{align*}
and  since $\psi(q\omega_-)<0$ and $\psi(q\gamma_1)<0$, the finiteness of the right-hand side above follows from standard properties of Lévy processes, see Lemma \ref{gene:levy:rivero} in the Appendix.
Putting the pieces together, we have checked \eqref{Eq:keyconv}, and the proof of (ii) is complete. 
 \end{proof}

We can now establish Proposition \ref{thm:cov:wei:leng}.
\begin{proof}[Proof of Proposition \ref{thm:cov:wei:leng}]
We work under $\P_1$ and equip the ssMt $\texttt{T}=(T,d_T, \rho, g)$ with the measures $\upmu$  and 
$\tilde{\uplambda}^\gamma:=-\kappa(\gamma)\cdot \uplambda^{\gamma}$ for  $\gamma\in(\omega_-,\omega_+)$. 
We infer from  Lemma \ref{lem:conv:mass} that there exists a sequence $(\gamma_n)_{n\geq 1}$ in $(\omega_-,p \omega_-)$ which converges to $\omega_-$ and  such that $\tilde{\uplambda}^{\gamma_n}( T) \to \upmu(T)$ a.s.  Our goal is to show that $\lim_{n\to \infty}\mathrm{d_{Prok}}( \tilde{\uplambda}^{ \gamma_n},\upmu)= 0$ a.s. 
Since the convergence for the total masses is  already know from Lemma \ref{lem:conv:mass},   it suffices to establish that, for every $\delta>0$ and every Borel set $A\subset T$,
we have 
\begin{equation}\label{Eq:ineqPro}
\upmu(A)\leq  \tilde{\uplambda}^{\gamma_n}(A^{\delta})+\delta, \qquad \text{for all $n$ sufficiently large},
\end{equation}
where we use the standard notation $A^\delta$ to denote the $\delta$-neighborhood of $A$ in $T$.  

In this direction, recall Notation \ref{N:subtrees}, and in particular that 
 for every $v\in \mathbb{U}$, $(T_v,  d_{T_v}, \uprho(v), g_v, \upmu_v)$ stands for the subtree encoded by the sub-family  $(f_{vu},\eta_{vu})_{u\in \U}$. 
  By the branching property,  the conditional  law of $(T_v,  d_{T_v}, \uprho(v), g_v, \upmu_v)$, in $\mathbb{T}$, given the type $\chi(u)=y$ is that of the ssMt under $\P_y$. 
 Therefore, by self-similarity, we get 
\begin{equation} \label{Eq:cvmasssub}
\lim \limits_{n\to \infty} \tilde{\uplambda}^{\gamma_n}(T_v)= \upmu_v(T_v)=\upmu(T_v),\quad \text{for every } v\in \mathbb{U}, \qquad \text{a.s.}
\end{equation}

Using again Notation \ref{N:subtrees}, we can decompose the tree $T$ at a generation  $k\geq 1$  as
$$T=T^k \cup\left(\bigcup_{|v|=k}T_v\right).$$
 We stress that  for any two distinct vertices at the same generation, say $v\neq w$ with $|v|=|w|=k$, 
the subtrees $T_v$ and $T_w$ are either disjoint, or they share the same root and  their intersection is then reduced to the latter.  In both cases, we have  $\upmu(T_v\cap T_w)=0$.
Combining this observation with the fact from Corollary \ref{C:CMJ}, that for any $k\geq 1$,  the harmonic measure assigns no mass to $T^k$, we conclude that
  there is the identity
 $$\upmu(A)= \sum_{v\in \N^k} \upmu(A\cap T_v).$$

  We then recall (from Propositions \ref{T:recolheightsup} and \ref{P:constructionomega-}, and Property ($\mathcal P$))   that  
 $$\lim_{k\to \infty} \sup\{\mathrm{Height}(T_v): v\in \N^k\} =0;$$
  we can therefore  almost surely   find a (random) integer $k$ sufficiently large such that
\begin{equation}\label{eq:diam:varkappa:theta}\mathrm{Height}(T_v) < \delta/2, \qquad \text{for all }v\in \N^k.
\end{equation}
  Next, since  $\upmu(T)=\sum_{|v|=k}m_v$ and $m_v=\upmu(T_v)$, we can select finitely many distinct vertices $v_1, \dots, v_M$ at generation $k$, such that
$$\sum_{1\leq i\leq M}\upmu\left(T_{v_i}\right)\geq \upmu(T)-\delta/2.$$
As a consequence of \eqref{Eq:cvmasssub}, we can now find an integer $N\geq 1$, such that
$$\upmu\left(T_{v_i}\right)\leq \tilde{\uplambda}^{\gamma_n}\left(T_{v_i}\right)+\delta/(2M), \qquad\text{ for every $1\leq i\leq M$ and $n\geq N$.}$$

Take any Borel subset $A$ of  $T$ and combine the preceding observations. We get
$$\upmu(A)\leq  \mathop{\sum \limits_{1\leq i\leq M}} \upmu\left(A\cap T_{v_i}\right)+\delta/2\leq  \mathop{\sum \limits_{1\leq i\leq M}}\limits_{A\cap  T_{v_i}\neq\emptyset} \tilde{\uplambda}^{\gamma_n}\left(T_{v_i}\right)+\delta.$$ 
Finally, \eqref{eq:diam:varkappa:theta} entails that 
for any vertex $v\in \N^k$  with $A\cap  T_{v}\neq\emptyset$, $T_v$ is included into the $\delta$-neighborhood of $A$, and 
the previous display is therefore bounded from above by $\tilde{\uplambda}^{\gamma_n}(A^\delta)+ \delta$;  here we used the facts that 
 $\tilde{\uplambda}^{\gamma_n}$ has no atoms and that the intersection of two different subtrees at the same generation is either empty or a singleton.  This completes the proof of our claim.
\end{proof}

\section{Comments and bibliographical notes}

\label{sec:commentsGBP}

 \paragraph{Construction of ssMt.} As mentioned in the previous chapter, the inspiration for the recursive random construction of self-similar Markov trees is the work of Rembart \& Winkel. In \cite{rembart2018recursive} those authors already constructed the underlying tree structure (but without the decoration) of binary growth-fragmentation processes (a subclass of our ssMt)  and gave a similar upper bound on the Hausdorff dimension as in Lemma \ref{Lem:haus:upper}. After the initial works on self-similar fragmentations \cite{Ber06,HM04}, the introduction of branching self-similar Markov processes can be traced back to  \cite{Ber15} (in particular Lemma \ref{L:verCMJ1} is adapted from  \cite[Lemma 4]{Ber15}) in the context of (binary) Growth-Fragmentation  and \cite{bertoin2019infinitely} in the context of branching L\'evy processes. Most of the framework of Section \ref{sec:2.2} is adapted from the literature on branching L\'evy processes \cite{bertoin2019infinitely}. A branching L\'evy process can be seen as the continuous-time version of a branching random walk which describes a particle system on the real line in which particles move and reproduce independently in a Poissonian manner. Just as for L\'evy processes, the law of a branching L\'evy process is determined by its characteristic triplet $(\sigma^{2}, \mathrm{a},  \boldsymbol{\Lambda})$ where the decorated L\'evy measure $ \boldsymbol{\Lambda}$ describes the intensity of the Poisson point process of births and jumps. In a nutshell, the self-similar Markov branching trees can be interpreted as the random decorated trees obtained after performing a Lamperti transformation with exponent $\alpha >0$ to the decorated trees coding for the genealogy of branching L\'evy processes. 

\paragraph{Critical case.} Our construction of ssMt in Proposition \ref{P:constructionomega-} assumes sub-criticality of the cumulant, i.e. that $\kappa(\gamma_{0})<0$ for some $\gamma_{0} >0$. Indeed, when $\kappa$ is strictly positive, the underlying branching random walk should witness local explosion, see \cite{bertoin2016local}, and it is hopeless to define a random compact tree from it. However, we left aside of this work the \textbf{critical case} when there exists $\omega >0$ for which $ \kappa(\omega)=0$ and otherwise $\kappa \geq 0$. In this case, we do not believe that our Theorem \ref{T:CMJ} can apply, since there are cases when 
$$ \sum_{n \geq 0} \mathbb{E}\big(\sup_{|u| =n}\chi(u)\big) =  \infty.$$
However, the application of Theorem \ref{T:recolinfty} requires the more flexible assumption \eqref{E:seriehypbis}, namely $ \lim_{k \to \infty} \sup\Big\{ \sum_{n=k}^{\infty} z_{ \bar{u}(n)}: \bar{u} \in \mathbb{N}^{\mathbb{N}}\Big\} =  0$ and its associated branching random walk analog 
$$ \lim_{k \to \infty} \sup\Big\{ \sum_{n=k}^{\infty}\chi(\bar{u}(n)): \bar{u} \in \mathbb{N}^{\mathbb{N}}\Big\} =  0.$$ 
In terms of \textit{growth-fragmentation process} (see Section \ref{sec:5.3} for details), although the particle system evolving in continuous time makes sense as well in the critical case, we do not know in general whether  its extinction time is finite a.s. or not (see \cite[Corollary 3]{Ber15}) and in turn whether the construction of Proposition \ref{P:constructionomega-} yields to a compact decorated random real tree.  A private communication of Elie Aidekon, Yueyun Hu and Zhan Shi reported progress in that direction and we believe that ssMt could  even be defined in the critical case\footnote{Notice then that Assumption \ref{A:omega-} cannot hold and natural harmonic measure $\upmu$ must be constructed from the \textit{derivative martingale}, \cite[Section 3.4]{shi2015branching}.}. See Remark \ref{rek:criticalbessel} and Example~\ref{ex:ADS}  for instances of critical cases.

\paragraph{On the intrinsic martingale.} 
Our construction of the harmonic measure on ssMt is based on a natural additive martingale, which is often referred to as the intrinsic martingale in the branching random walk literature. Biggins \cite{Big77a} and  Biggins \& Kyprianou \cite{biggins2005fixed} have given explicit  general criteria (extending the celebrated Kesten-Stigum theorem) which entail that the convergence also holds in mean of those martingales. In our case, Assumption \ref{A:omega-} and Lemma \ref{L:MUI} which ensures the convergence in mean builds on the recent works Alsmeyer \& Kuhlbusch \cite{alsmeyerdouble} and Iksanov \& Mallein \cite{iksanov2019result}. Although Lemma \ref{L:MUI} proves a convergence in $L^{p}(\P_1)$ for some $p>1$, Proposition \ref{P:constructionomega-} merely requires the uniform integrability. Even though we employ the condition of boundedness in $L^{p}(\P_1)$ to streamline some technical arguments, it might be possible to bypass it. Our results are likely to remain valid under less stringent assumptions ensuring the uniform integrability of the intrinsic martingale $(M_n(\omega_-))_{n\geq 1}$.

The law of $M_{\infty}^{-}$, the limit of the intrinsic martingale is subject of intense study in the literature: Passing \eqref{eq:decomp1} to the limit using the convergence in mean, we deduce the following recursive distributional equation
  \begin{eqnarray} M_{\infty}(\omega_-) \overset{(d)}{=} \sum_{u \in \mathbb{N}}  \chi(u)^{\omega_{-}} \cdot M_{\infty}^u(\omega_-),   \label{eq:smoothing}\end{eqnarray}
where on the right hand side the iid copies $M_{\infty}^u(\omega_-), u \in \mathbb{N}$ are independent of the vector $(\chi(u) : u \in \mathbb{N})$. The above fixed point equation, related to the equation “$X=AX+B$" is sometimes called a \textbf{smoothing  transform} and has been studied in depth in recent years, see \cite{dariusz2016stochastic}. Under our assumptions, we infer in particular from results due to  Biggins (see Section~2 in \cite{Big77a}), that this equation has a unique solution with given mean.  
In our Markovian setup, we can also write an ``infinitesimal version'' of the recursive distributional equation \eqref{eq:smoothing}, see \cite[bottom page 4]{bertoin2018biggins}, and deduce that the Laplace transform $w (\lambda) = \mathbb{E}\big(\mathrm{e}^{-\lambda M_{\infty}(\omega_-)}\big)$, $\lambda\geq 0$ , satisfies the following integro-differential equation, 
\begin{align} \label{eq:transfolaplacevolume}
0 = & - \mathrm{k}w(\lambda)+ \frac{1}{2} \sigma^{2} \lambda^{2}w''(\lambda) + a \, \lambda w'(\lambda) \nonumber \\
&+ \int \nu ( \mathrm{d} \mathbf{y}) \left(\prod_{i \geq 0} w(\lambda y_{i}^{\omega_{-}}) -w(\lambda) - y_{0}\lambda w'(\lambda) \mathbf{1}_{ \mathrm{e}^{-1} \leq y_{0} \leq  \mathrm{e}} \right), \end{align}
 where $\nu$ is the image measure of $ \boldsymbol{\Lambda}$ by $x \mapsto \mathrm{e}^{x}$.  Alas, despite the fact that the law of $M_{\infty}(\omega_-)$ satisfies \eqref{eq:smoothing} and \eqref{eq:transfolaplacevolume}, it seems to be difficult to identify its distribution in general. We shall however be able to identify the law of $M_{\infty}(\omega_-)$ in some special cases in particular via a surprising connection with random planar maps, see Chapter \ref{chap:example}.
\medskip

 \chapter{Examples} \label{chap:example}

The purpose of this chapter is to illustrate the construction of self-similar Markov trees and  to discuss some distinguished families of examples. Roughly speaking, many of those examples have the property that the decorations along branches are given by some versions of a stable L\'evy process. We stress here again -see also the end of Section \ref{sec:2.2}- that the law of a
self-similar Markov tree is determined by a characteristic quadruplet, and that different characteristic quadruplets in the same equivalence class of bifurcators induce self-similar Markov trees with the same distribution. Choosing one characteristic quadruplet rather than another one within an equivalence class of bifurcators is often only a matter of preferences. The reader may wish to have first a glance at the forthcoming Section \ref{sec:bifurcators}
that will provide a detailed account on bifurcators; notably Theorem \ref{sec:bifurcators} there describes these equivalence classes explicitly. 
Nonetheless, the examples treated in this chapter do no require any result from Section \ref{sec:bifurcators}. 

Before starting our list of examples and in order to give a purely analytic definition in terms of the characteristic quadruplet, we also need to make some comments on the concept of drift. The notion of drift  coefficient for a L\'evy process $\xi$  can be defined canonically when  the L\'evy measure $\Lambda_0(\dd y)$  integrates $1\wedge |y|$, but   depends  otherwise of the arbitrary choice of the cutoff function  in the L\'evy-Khintchin formula  \eqref{E:LKfor}.
More precisely, if  $\int (1\wedge |y|) \Lambda_0(\dd y) < \infty$, then no compensation is needed to make sense of the Poisson integral there and 
we can re-write the L\'evy-It\^{o} decomposition 
\eqref{E:LevyIto} in the simpler form
$$
\xi(t) \coloneqq \sigma B(t) + \mathrm a_{\mathrm{can}}t + \int_{[0,t]\times \R}N_0( \d  s,  \d y)  ~y , \qquad \mathrm a_{\mathrm{can}}\coloneqq {\mathrm a} - \int \Lambda_0(\dd y) y {\mathbf 1}_{|y|\leq 1}.
$$
Note that $\mathrm a_{\mathrm{can}}=\mathrm a$ when all the jumps of $\xi$ have size greater than $1$   (i.e. when $\Lambda_0([-1,1])=0$), but otherwise these two coefficients are usually different. 
It is well-known that $\mathrm a_{\mathrm{can}}$ is then a much more relevant quantity for the L\'evy process $\xi$ than the rather artificial  ${\mathrm a}$, and 
we shall therefore call $\mathrm{a_{can}}$ the \textbf{canonic drift coefficient}. Notice that in this case, the cumulant function takes the following simpler form
  \begin{eqnarray} \label{eq:kappasimple}\kappa(\gamma)
= \frac{1}{2}\sigma^2 \gamma^2+ {\mathrm{a_{can}}}\gamma + \int_{\mathcal{S}} \boldsymbol{\Lambda}( \d y_0,\d  \mathbf y )~ \left(\big(\sum_{i=0}^\infty\e^{\gamma y_i}\big)-1\right). \end{eqnarray}

\section{Finite branching activity} \label{sec:5.1}

One often says that a L\'evy process has a finite activity when almost surely, its sample paths have only a finite number of jumps along any time interval of finite duration; this is equivalent to the finiteness of the L\'evy measure.
In the setting of Section \ref{sec:2.2}, recall  that $ \boldsymbol{ \Lambda}_1$ denotes the image of the generalized L\'evy measure $ \boldsymbol{ \Lambda}$ by the second projection  from $\mathcal S=[-\infty,\infty)\times \mathcal{S}_1$ to $\mathcal{S}_1$, and that $ \boldsymbol{ \Lambda}_1$ bears
a close relation to  the reproduction process $\eta$; see  \eqref{Eq:tildeeta:Markov} and  \eqref{Eq:eta:Markov}. Recall also that  in a reproduction event, the degenerate sequence
$(-\infty, -\infty, \ldots)\in \mathcal{S}_1$ should be interpreted as empty. So strictly speaking, 
an atom of the Poisson random measure $\mathbf N$  of the form $(t,y,(-\infty, -\infty, \ldots))$
 is not associated to any reproduction event, but only to a jump of the decoration at time $t$.

\begin{figure}[!h]
 \begin{center}
 \includegraphics[height=7.5cm]{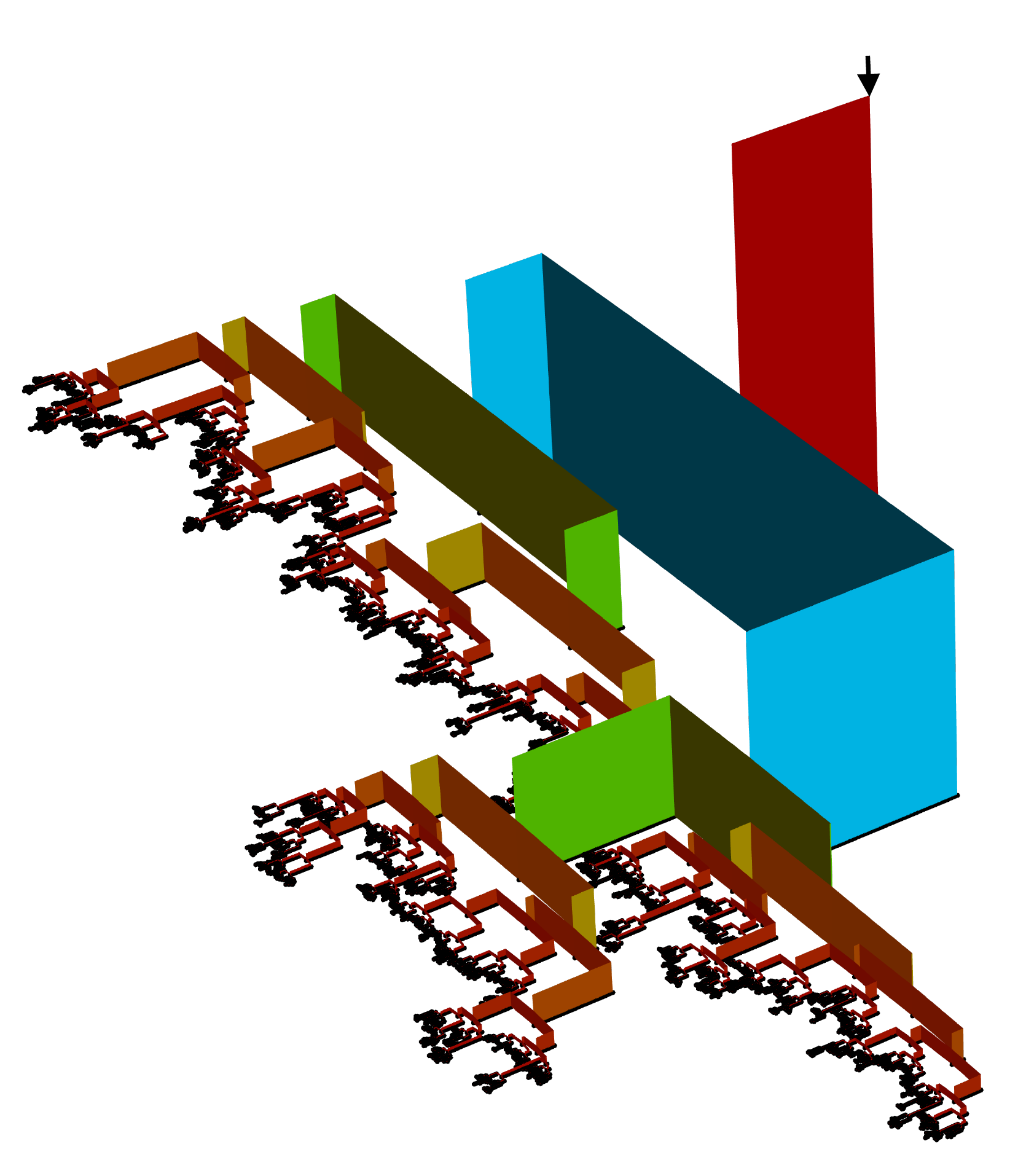}
 \caption{A simulation of the self-similar Markov tree  in Example \ref{Ex:stat} for $\alpha=0.4$. The tree $T$ is embedded (non-isometrically) in the plane, and the decoration is represented in the vertical dimension. The root is at the bottom of the right hand side, marked by an arrow.}
\label{Fig:static}
 \end{center}
 \end{figure}

\begin{figure}[!h]
 \begin{center}
   \includegraphics[width=10cm]{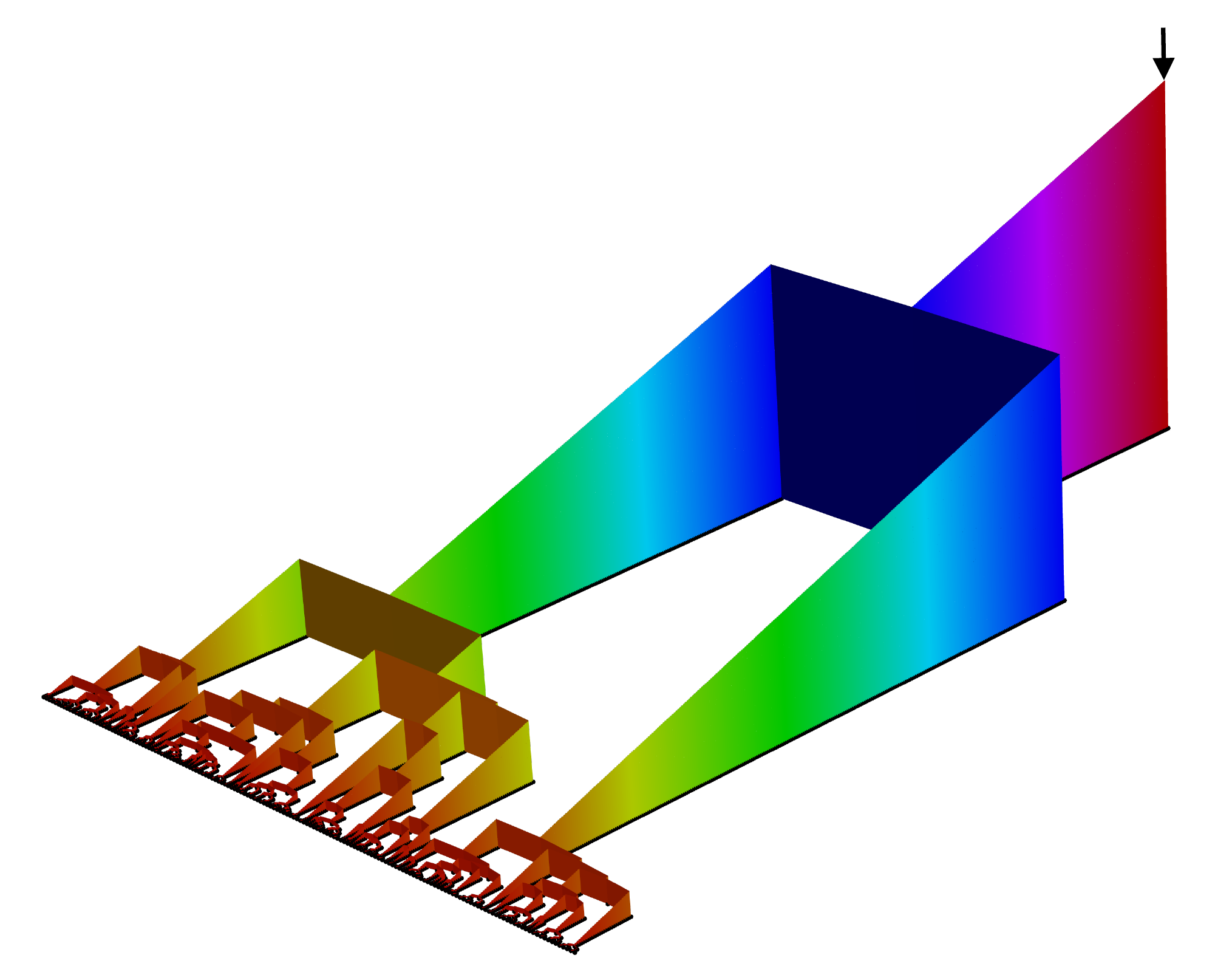}
 \caption{A simulation of reduced tree at height $1$ in the Brownian CRT. The reduced tree is embedded in the plane and the decoration (corresponding to the height of subtrees) corresponds to the vertical dimension. The root is placed on the right and marked with an arrow.}
  \label{Fig:reduced}
 \end{center}
 \end{figure}

We say that a self-similar Markov tree with characteristic quadruplet $( \sigma^2, \mathrm{a}, \boldsymbol{ \Lambda} ; \alpha)$ has a \textbf{finite branching activity} when
$$ \boldsymbol{ \Lambda}_1\big( \mathcal{S}_1 \backslash \{(-\infty, -\infty, \ldots)\}\big) < \infty.$$
This is the case precisely when the first atom of the reproduction process $\eta$ occurs at a (strictly) positive time a.s., and as a consequence, reproduction events for an individual can only (possibly) accumulate at the end of its lifetime. The tree has then essentially a discrete structure and
its branching points do not pile up, except possibly at  leaves; see Figures \ref{Fig:static} and \ref{Fig:reduced}.
Note also that since the killing rate $\mathrm k=\boldsymbol{ \Lambda}(\{-\infty\}\times \mathcal{S}_1)$
is always finite, a self-similar Markov tree has a finite branching activity if and only if 
its  generalized L\'evy measure  satisfies
\begin{equation} \label{Eq:finiteact}
\boldsymbol{ \Lambda}\big(\{(y,\mathbf y)\in \mathcal S:
 y=-\infty \text{ or } \mathbf y\neq(-\infty, -\infty, \ldots)\}\big) < \infty.
 \end{equation}

Having or not a finite branching activity is an intrinsic property of a self-similar Markov tree, hence  it does not depend on the choice of the characteristic quadruplet within a family of bifurcators.  Moreover, if a characteristic quadruplet $( \sigma^2, \mathrm{a}, \boldsymbol{ \Lambda} ; \alpha)$ has a finite branching activity, then we can always choose an equivalent bifurcator, which we denote again by 
$( \sigma^2, \mathrm{a}, \boldsymbol{ \Lambda} ; \alpha)$ for the sake of simplicity, such that the generalized L\'evy measure  now fulfills 
\begin{equation} \label{Eq:mortreprod}
\boldsymbol{ \Lambda} \big(\{(y,\mathbf y)\in \mathcal S:
 y\neq -\infty \text{ and } \mathbf y\neq(-\infty, -\infty, \ldots)\}\big)=0.
 \end{equation}
  This requirement means that reproduction events can only occur at
 the death of the parent,  which makes the reproduction process especially simple to depict.  
 In terms of general branching processes, considering such an equivalent bifurcator  amounts to the following somewhat artificial transformation. We decide to kill each individual at the first instant when it begets, and declare that at its death, it gives birth to an additional child then viewed as its reincarnation and whose decoration is thus defined by shifting the decoration of the parent. Although this transformation affects some genealogical aspects of the branching process, it has no impact on the decorated real tree which is induced.

We assume throughout the rest of this section  that the generalized L\'evy measure $\boldsymbol{ \Lambda} $
satisfies\footnote{Of course, \eqref{Eq:finiteact} and  \eqref{Eq:mortreprod} alone do not grant  Assumption \eqref{A:gamma0}
 and the latter is also needed to ensure the existence of the self-similar Markov tree.}
\eqref{Eq:finiteact} and  \eqref{Eq:mortreprod}. Then the self-similar Markov tree can be seen as the result of Lamperti's transformation applied to the decorated genealogical tree that records the evolution of a continuous time branching random walk in the sense of Biggins \cite[Section 5]{biggins1992uniform}. More precisely, the latter depicts a population of individuals living in $\R$, 
starting at time $0$  from a single ancestor located at the origin, and such that each individual lives for an exponentially distributed duration with a fixed parameter.
 During their lifetimes,  individuals move according to independent copies of a L\'evy process with characteristic triplet
 $( \sigma^2, \mathrm{a}, \Lambda_0)$ and  are thus killed with rate 
 $$\mathrm k=\Lambda_0(\{-\infty\})= \boldsymbol{ \Lambda}\big(\{-\infty\}\times \mathcal S_1\big)=\boldsymbol{ \Lambda}\big(\{(y,\mathbf y)\in \mathcal S:
 y=-\infty \text{ or } \mathbf y\neq(-\infty, -\infty, \ldots)\}\big),$$
 where the second identity stems from \eqref{Eq:mortreprod}. In particular, when a death event occurs, the individual which dies is chosen uniformly at random in the current population, independently of its location.

At the time of its death, the parent is replaced by its offspring. The distribution of the children positions relative to the parent is given by the normalized sub-probability measure  $\mathrm k^{-1}\boldsymbol{\Lambda}_1$ on $\mathcal{S}_1\backslash \{(-\infty, -\infty, \ldots)\}$, where the default of mass is the probability that an individual dies without begetting any child.
 The self-similar Markov tree is then obtained by interpreting the location of an individual in the continuous time branching random walk as a value of the (real) decoration on the genealogical tree, and then performing the Lamperti transformation on each line of descent. Last, we need to take the completion in order to deal with a compact structure.
 
 We will now describe in more detail three examples in which the motions of individuals and their reproductions for the continuous time branching random walk are particularly simple. In the first example, individuals are static and the sole motion occurs at birth. In the second, the motion of individuals is merely a linear drift and children are born at the same location as their parents. In the third, the displacements of individuals are governed by independent Brownian motions with drift, and again children are born at the same location as their parents. For the sake of simplicity, we mostly focus on binary branching, meaning that an individual begets exactly two children when it dies.

\begin{example}[Static, after Haas \cite{haas2018scaling}] \label{Ex:stat} Consider the characteristic quadruplet with $\sigma^2=0$, ${\mathrm{a_{can}}}=0$, $ \boldsymbol{\Lambda}_{ \mathrm{half}} = \delta_{(-\infty, (-\log 2, -\log 2,-\infty, \ldots))}$ and an arbitrary $\alpha>0$. 
In the continuous time branching random walk, each individual lives for a standard exponential duration and does not move until it dies. At death, each parent given birth to two children, both located at distance $\log 2$ at the left of the parent.

 The self-similar Markov tree is then obtained by performing the Lamperti transformation, which is elementary. Its structure  is very simple to described iteratively: it consists of a branch having a standard exponential length and decorated with the constant function $1$,  at the extremity of which two branches of independent  exponential lengths with mean $2^{-\alpha}$, each decorated with the constant function $1/2$ are grafted, and so on and so forth. Finally we take the closure; see Figure \ref{Fig:static} for an illustration.

The cumulant is simply
$$ \kappa(\gamma) = 2^{1-\gamma}-1, \quad \mbox{ for } \gamma >0.$$
 Assumption \ref{A:omega-} holds with $\omega_-=1$ and $\omega_+=\infty$, and since the sum of the decoration of the two children always equals the decoration of the parent, the intrinsic martingale is constant and the  total mass of the tree for the harmonic measure is merely $1$.
Turning our attention to weighted length measures, we observe from Proposition \ref{prop:lengthmeasures} that $\uplambda^\gamma$ is finite
for any $\gamma > 1$.  An application of the branching property shows that the total mass $\uplambda^\gamma(T)$
 satisfies the fixed-point equation in distribution $$ \uplambda^\gamma(T) \overset{(d)}{=}  \mathcal{E}(1) + 2^{-\gamma}( \lambda_1 +  \lambda_2)$$ where $ \lambda_1$ and $ \lambda_2$ are two independent copies of $\uplambda^\gamma(T)$ also independent of the exponential random variable $ \mathcal{E}(1)$.  \hfill $\diamond$
\end{example}

 The second example is a simple variation of the Yule process, which has also a natural interpretation in terms of  reduced stable trees, as it will be discussed afterwards.

 \begin{example}[Binary branching and drift] \label{ex:reduced} Consider here the characteristics $ \sigma^{2}=0$, ${\mathrm{a_{can}}} =-1$,   $\boldsymbol{ \Lambda}_{ \mathrm{two}} =  \delta_{(-\infty,(0,0,-\infty, \ldots))}$ and $ \alpha =1$. In the continuous time branching random walk,  individuals live in $\R_-$ and drift continuously with velocity $-1$; they die at rate $\mathrm k=1$ and then are replaced by two children at the same location. This is an elementary spatial version of the Yule process (see \cite[Chapter III.5]{AN72}); in particular, the number of individuals at any given time $t\geq 0$ has the geometric distribution with parameter $\e^{-t}$ and all these individuals are located at $-t$.  The sequence of times at which birth events occur is given by the partial sums $\epsilon_1+ \cdots + \epsilon_n$ for $n\geq 1$, where the variables $\epsilon_i$'s are independent and 
each $\epsilon_i$ has the exponential distribution with parameter $i$. As a consequence, the variables $\beta_i=\exp(-\epsilon_i)$ are independent beta$(i,1)$ variables with distribution functions $\P(\beta_i\leq r)=r^i$ for $r\in[0,1]$. In the present case, the  Lamperti transformation amounts to combining the time-change $\tau(r)=-\log(1-r)$ for $r\in[0,\beta_1)$ and the exponential map $y\to \e^y$ in space. In particular, the decoration process is simply given by $X(r)=\exp(-\tau(r))=1-r$, for $r\in [0,\beta_1)$ and $X(\beta_1)=0$. Furthermore, the ranked sequence of the heights of the branching points in the self-similar Markov tree  is $1-\beta_1\cdots \beta_n$, for $n\geq 1$.

Equivalently, the self-similar Markov tree can then be constructed recursively as follows; see Figure \ref{Fig:reduced} for an illustration.  We start with the line segment with unit length.  As a first step, we  glue at  height $1-\beta_1$ a segment with length $\beta_1$. Next we choose one of the two points at height $1-\beta_1\beta_2$ uniformly at random and glue there a segment with length $\beta_1\beta_2$, and so on, and so forth. That it, at the $n$-th step, we pick on the tree constructed so far one of its $n$ points at height $1-\beta_1\cdots \beta_n$ uniformly at random
 and glue there a segment with length $\beta_1\cdots \beta_n$. We end the construction by completing the tree and define the decoration as $1$ minus the height function. The decoration $g(x)$ at  vertex $x$  is then the height of the fringe subtree that stems from $x$.
 
We next turn our attention to the intrinsic martingale and the harmonic measure $\upmu$. By \eqref{E:cumulant}, we have 
$$ \kappa(\gamma) = 1-\gamma, \quad \mbox{ for } \gamma >0,$$ so that  Assumption \ref{A:omega-} holds with $\omega_{-} =1$ and $\omega_+=+\infty$. On the other hand, we see from  \eqref{Eq:eta:Markov} that reproduction process simply given by
$$\eta= 2\delta_{(U,U)},$$
where $U=1-\exp(-\epsilon_1)$ has the uniform distribution on $[0,1]$. Recall that the harmonic mass $\upmu(T)$ is given by the terminal value of the intrinsic martingale,
and we see from the branching property that the latter satisfies the distributional identity
$$\upmu(T) {\overset{(d)}{=}} \ U(\mu_1+ \mu_2),$$
where in the right-hand side, the variables $\mu_1$ and $\mu_2$ are independent copies of $\upmu(T) $, also independent of $U$.
It is easy to deduce (for instance, by considering Laplace transforms) that $\upmu(T) $ has the standard exponential distribution,
as we may expect from a well-known result on Yule processes (see, e.g. \cite[Problem 2 on page 136]{AN72}). 

 We consider as well the weighted length $\uplambda^{\gamma}(T)$ for $\gamma>1$, and get similarly the identity in distribution
 $$\uplambda^{\gamma}(T) {\overset{(d)}{=}} \gamma^{-1}(1-U^{\gamma}) + U^{\gamma}(\lambda_1+\lambda_2),$$
 where in the right-hand side, the variables $\lambda_1$ and $\lambda_2$ are independent copies of $\uplambda^{\gamma}(T) $, also independent of the uniform variable $U$.
 \hfill $\diamond$

\end{example}

 The self-similar Markov tree in Example \ref{ex:reduced} is well-known in the literature and notably appeared in connection with the Brownian motion and the Brownian CRT. 
 In this direction, let us first briefly recall basics about It\^{o}'s measure of Brownian excursions and the coding of real trees by excursions since this will serve in several examples below.

We consider the It\^{o} measure $ \mathbb{N}_{ \mathrm{Ito}}$ of positive Brownian excursions, which can be constructed using It\^{o} or William's decomposition, see \cite[Chapter XII, Theorem 4.2  and 4.5]{RY99}. Specifically, there exist a family of probability distributions $ (\mathbb{N}_\ell)_{\ell >0}$ on excursions of fixed length $\ell$, and a family $( \mathbb{N}^{(h)})_{h> 0}$ on excursions of fixed height $h$, which both inherit self-similarity from Brownian motion, such that 
  \begin{eqnarray} \label{def:itobro}  \mathbb{N}_{ \mathrm{Ito}} = \int_{0}^{\infty} \frac{ \mathrm{d}h}{2 h^{2}} \mathbb{N}^{(h)} = \int_{0}^{\infty} \frac{ \mathrm{d}\ell}{2 \sqrt{2 \pi \ell^{3}} } \mathbb{N}_{\ell}.  \end{eqnarray}

\label{sec:codagearbre}
 Next, recall from \cite{DLG05} that any continuous excursion, say $\e: [0,\ell]\to \R_+$ with $\e(0)=\e(\ell)=0$, encodes a rooted planar continuous tree, say $\mathcal T_{\e}$, via the contour function.
 The latter is a continuous surjection  $c_{\e}: [0,\ell]\to \mathcal T_{\e}$, and it is then natural to endow $\mathcal T_{\e}$
 with the push-forward image of the Lebesgue measure on $[0,\ell]$ by 
  $c_{\e}$, which we call here the contour measure\footnote{The same tree can be encoded by different excursions, and each such excursion induces a  contour measure.} on $\mathcal T_{\e}$ and denote by $\gamma_{\e}$.  
    We shall see in Examples \ref{ex:brownianheight} and \ref{ex:brownian} that the random tree $ \mathcal{T}_{\e}$ under $ \mathbb{N}_1$ or $ \mathbb{N}^{(1)}$ can be seen as self-similar Markov tree (see also Example \ref{ex:ADS} for a related example).

 At the heart of the connexions between Example \ref{ex:reduced} and Brownian excursions or CRT lies the observation that the sequence $(\beta_1\cdots \beta_n)_{n\geq 1}$ of the random lengths used in the construction by gluing
can be obtained as the sequence ranked in the decreasing order of the atoms of a Poisson point measure on $(0,1)$ with intensity $\ell^{-2}\dd \ell$. 
By Brownian excursion theory, the latter has the same statistics as the ranked sequence of the depths of the family of the excursions below $1$ for a Brownian motion started from $1$ and killed when hitting $0$. 
 Now imagine that instead of ranking these lengths in decreasing order, we keep the natural order induced by time) of occurrences of the excursions of the Brownian motion (informally speaking, this corresponds to a uniform random shuffling of the ranked sequence), and add a unit length 
 at the right-end to take into account the final excursion that brings the Brownian motion to $0$. By gluing each length $\ell <1$ at its bottom point to the first larger length at its right, 
 we recover precisely the recursive construction devised in Example  \ref{ex:reduced}.

 This construction can also be interpreted as follows. We work under the It\^{o} measure of positive Brownian excursion conditioned to have height at least $1$, which thanks to \eqref{def:itobro}  
 can be expressed in the form
$$  \int_{1}^{\infty} \frac{ \mathrm{d}h}{h^{2}} \mathbb{N}^{(h)} .$$ 
  Under this conditional probability measure, a sample path can be viewed as the contour process of a planar continuous tree with height at least $1$. Imagine that we reduce this tree by removing first all the vertices at height greater than $1$, and further by pruning all the remaining branches that do not reach height $1$. Then this reduced tree 
has the same distribution as the self similar Markov tree in  Example \ref{ex:reduced}; see e.g.~\cite[Section 2.1]{CLGharmo}.

This reduction process can be applied to the stable trees of the forthcoming Example \ref{ex:stable}. In this case we still have $\alpha=1$, the drift is ${\mathrm{a}}=-1$, the Gaussian coefficient $\sigma^{2}=0$ and the generalized L\'evy measure is given by 
$$ \boldsymbol{ \Lambda} = \sum_{k \geq 1} \delta_{(-\infty,(\underbrace{0,...,0}_{k \mathrm{\ times}}, -\infty, \ldots))} \cdot \frac{\beta \Gamma(k-\beta)}{k! \Gamma(2-\beta)},$$ so that the associated cumulant function is $$ \kappa(\gamma) = \frac{1}{\beta-1} -\gamma.$$ In particular  Assumption \ref{A:omega-} still holds with $\omega_{-}= \frac{1}{\beta-1} \in [1, \infty)$ and $\omega_+= +\infty$,   and in this case, the total $\upmu$-mass has Laplace transform 
 $$ \mathbb{E}\big( \mathrm{e}^{-\gamma M_{\infty}(\omega_-)}\big) = 1 - \left(1 +  \frac{1}{\gamma^{\beta-1}}\right)^{-1/(\beta-1)},$$
 see \cite{lin2014harmonic} and the references therein.

We end this subsection with a last example where now the motions of individuals are diffusion processes. 
\begin{example}[Branching Bessel processes] \label{Ex:bbess} Consider the characteristics $\sigma^{2}=1$, $ {\mathrm{a_{can}}} \in \mathbb{R}$, $ \boldsymbol{\Lambda}_{ \mathrm{two}}= \delta_{(-\infty,(0,0, -\infty, \ldots))}$,   and self-similarity parameter $\alpha = 2$. The continuous time branching process is then the well-known binary branching Brownian motion with drift $ {\mathrm{a_{can}}}$, and the cumulant function is 
$$ \kappa(\gamma) = \frac{\gamma^{2}}{2} + {\mathrm{a_{can}}}\gamma+1.$$
It satisfies Assumption \ref{A:omega-} as soon as $ {\mathrm{a_{can}}} < - \sqrt{2}$ and then $\omega_{-}= - {\mathrm{a_{can}}} - \sqrt{ {\mathrm{a_{can}}}^{2}-2}$  and $\omega_{+}= - {\mathrm{a_{can}}} + \sqrt{ {\mathrm{a_{can}}}^{2}-2}$. By the Lamperti transformation, the self-similar Markov process $X$ associated to the Brownian motion with drift $B_{t} + {\mathrm{a_{can}}}t$ is a Bessel process with dimension $d = 2 {\mathrm{a_{can}}} +2$.
See Figure \ref{Fig:Bessel} for an illustration.

Similarly to Example \ref{ex:reduced}, the spatial marginal of the reproduction measure is 
$$\eta(\R_+\times \dd x) = 2\delta_{\exp(B_\epsilon + {\mathrm{a_{can}}}\epsilon)}(\dd x),$$
where $\epsilon$ is an exponentially distributed random time independent of the Brownian motion $B$.
The density of the variable $B_\epsilon + {\mathrm{a_{can}}}\epsilon$ is given by
$$\P(B_\epsilon + {\mathrm{a_{can}}}\epsilon \in \dd y) = (2+{\mathrm{a_{can}}}^2)^{-1/2} \exp({\mathrm{a_{can}}}y-|y|\sqrt{2+{\mathrm{a_{can}}}^2})\dd y,$$
see \cite[(1.0.5) on page 256]{borodin2015handbook}. The terminal value of the intrinsic martingale, that is the mass of the harmonic measure, then satisfies the distributional identity
$$M_{\infty}(\omega_-) {\overset{(d)}{=}} \ \exp(\omega_-(B_\epsilon + {\mathrm{a_{can}}}\epsilon))(M_{\infty}(\omega_-)+M'_{\infty}(\omega_-)),$$
where in the right-hand side, the variables $B_\epsilon + {\mathrm{a_{can}}}\epsilon$, $M_{\infty}(\omega_-)$ and $M'_{\infty}$ are independent, and $M'_{\infty}$ is a copy of $M_{\infty}(\omega_-)$. 
It is not known to us whether an explicit expression for the distribution of $M_{\infty}(\omega_-)$ can be derived from this. 
\hfill $\diamond$
\begin{figure}[!h]
 \begin{center}
\includegraphics[width=10cm]{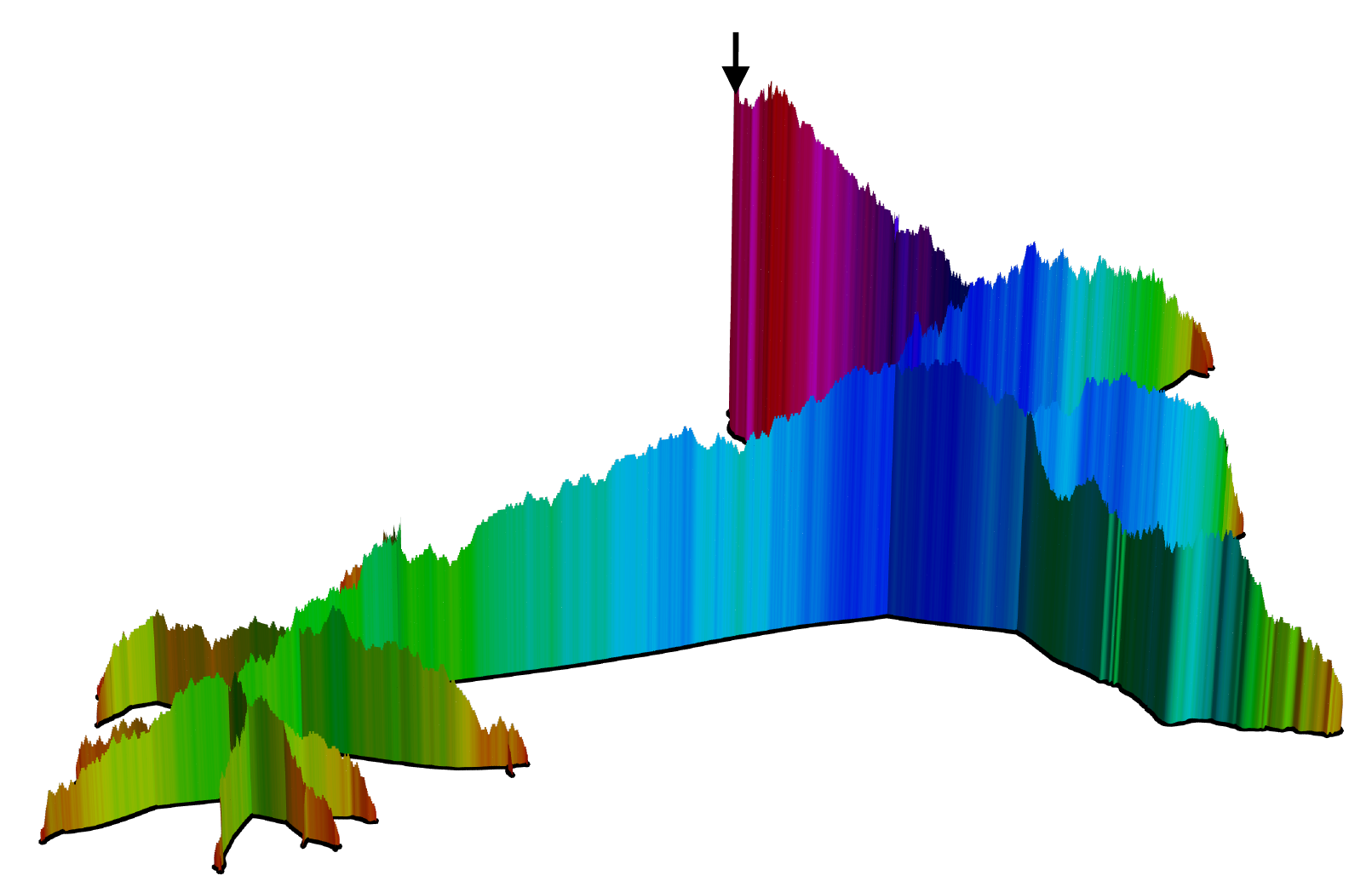}
 \caption{A simulation of the Branching Bessel process with $ {\mathrm{a_{can}}} =-3$. As usual, the decoration (the values of the Bessel processes) is displayed in the vertical coordinate. The root of the decorated tree is marked  by an arrow. \label{Fig:Bessel} }  
 \end{center}
 \end{figure}
\end{example}

\begin{remark}[Critical case]  \label{rek:criticalbessel}If we take ${\mathrm{a_{can}}} = - \sqrt{2}$ in the above example, then the overall minimum of the cumulant function $\kappa$ equals $0$, which we called the  \textit{critical case} in the comments section of Chapter \ref{chap:generalBP}. In this example, we leave open the fact that the procedure of Chapter \ref{chap:generalBP} actually produces a compact decorated random tree.
\end{remark}

\section{Non-increasing decorations and fragmentations} \label{sec:5.2}

We say that a decoration on rooted real tree  is non-increasing when its restriction  to any segment from the root defines a non-increasing function of the height. 
We then call a self-similar Markov tree \textbf{non-increasing} if its  decoration is non-increasing, almost-surely.
In that case, the self-similar Markov process $X$ that describes the decoration for a typical individual must plainly have non-increasing sample paths, that is, equivalently, the underlying L\'evy process $\xi$ that gives $X$ after the Lamperti transformation must be  the negative of a subordinator (possibly with killing). However, this requirement is clearly not sufficient; one needs  to impose further that the reproduction process $\eta$ has no atom, say at $(t,x)$, with $X(t-)<x$. The latter translates in terms of the generalized L\'evy measure as
\begin{equation}\label{Eq:GLMnoninc}
\boldsymbol{\Lambda}_1\big(\left\{\mathbf y=(y_1, \ldots)\in \mathcal{S}_1: y_1>0\right\}\big)=0.
\end{equation}
In the converse direction, it is immediately seen that if \eqref{Eq:GLMnoninc} holds and if the underlying L\'evy process $\xi$ is the negative of a subordinator, 
then the self-similar Markov tree is non-increasing. 
For instance, Examples \ref{Ex:stat} and \ref{ex:reduced}  are non-increasing self-similar Markov trees (with a finite branching activity), but not Example \ref{Ex:bbess}.  In particular, a (possibly killed) real L\'evy process is a subordinator -- that is, it has non-decreasing sample paths until it eventually dies almost surely -- if and only if its Gaussian coefficient is $\sigma^2=0$, its L\'evy measure $\Lambda_0(\dd y)$ gives zero mass to $(-\infty, 0)$ and integrates the function $1\wedge |y|$,
and finally, its canonic drift coefficient is nonnegative, $\mathrm a_{\mathrm{can}}\geq 0$. 

Putting the pieces together, we arrive at the following analytic definition in terms of the characteristic quadruplet\footnote{Whether or not a self-similar Markov tree is  non-increasing is obviously an intrinsic property, which  does not depend of the choice of the characteristic quadruplet within a family of bifurcators.}. A self-similar Markov tree is non-increasing if and only if its Gaussian coefficient is zero, $\sigma^2=0$,
its generalized L\'evy measure verifies
$$\boldsymbol{\Lambda}\big(\left\{(y,\mathbf y)\in \mathcal{S}: y>0 \text{ or } y_1>0\right\}\big)=0,$$
and 
\begin{equation} \label{Eq:lambda0int}\int_{(-\infty,0)} (1\wedge |y|) \Lambda_0(\dd y) <\infty, 
\end{equation}
and finally, its canonic drift coefficient is non-positive, 
$$\mathrm a_{\mathrm{can}}= {\mathrm a}  -  \int_{[-1,0)} y \Lambda_0(\d y) \leq 0.$$
 We stress however that at this stage,  the cumulant function $\kappa$ in  \eqref{E:cumulant}  could potentially be  infinite everywhere, and \textit{a fortiori} Assumption \eqref{A:gamma0} may fail in general.  \bigskip

We now present an example of a non-increasing self-similar Markov tree with an infinite branching activity which is related to Example \ref{ex:reduced}.  Recall the notation regarding It\^{o} excursion measure and contours of planar trees that was introduced after Example \ref{ex:reduced}.

    \begin{example}[Heights of Brownian sub-trees] \label{ex:brownianheight} Let us denote  the Brownian CRT of height $1$ by $ \mathcal{T}^{(1)}$, i.e. the tree $ \mathcal{T}_{\e^{(1)}}$ where $\e^{(1)}$ follows the law $ \mathbb{N}^{(1)}$ in \eqref{def:itobro}.
We further endow $\mathcal {T}^{(1)}$ with the deterministic decoration which assigns to each vertex $v\in \mathcal {T}^{(1)}$
  the height of the fringe-subtree ${\mathcal T}^{(1)}_v$ rooted at $v$ as defined in Section \ref{sec:1.1};  see Figure \ref{fig:Dino} for an illustration. We shall now argue that this decorated random real tree is a self-similar Markov tree and determine its characteristics.

      \begin{figure}[!h]
   \begin{center}
      \includegraphics[height=8cm]{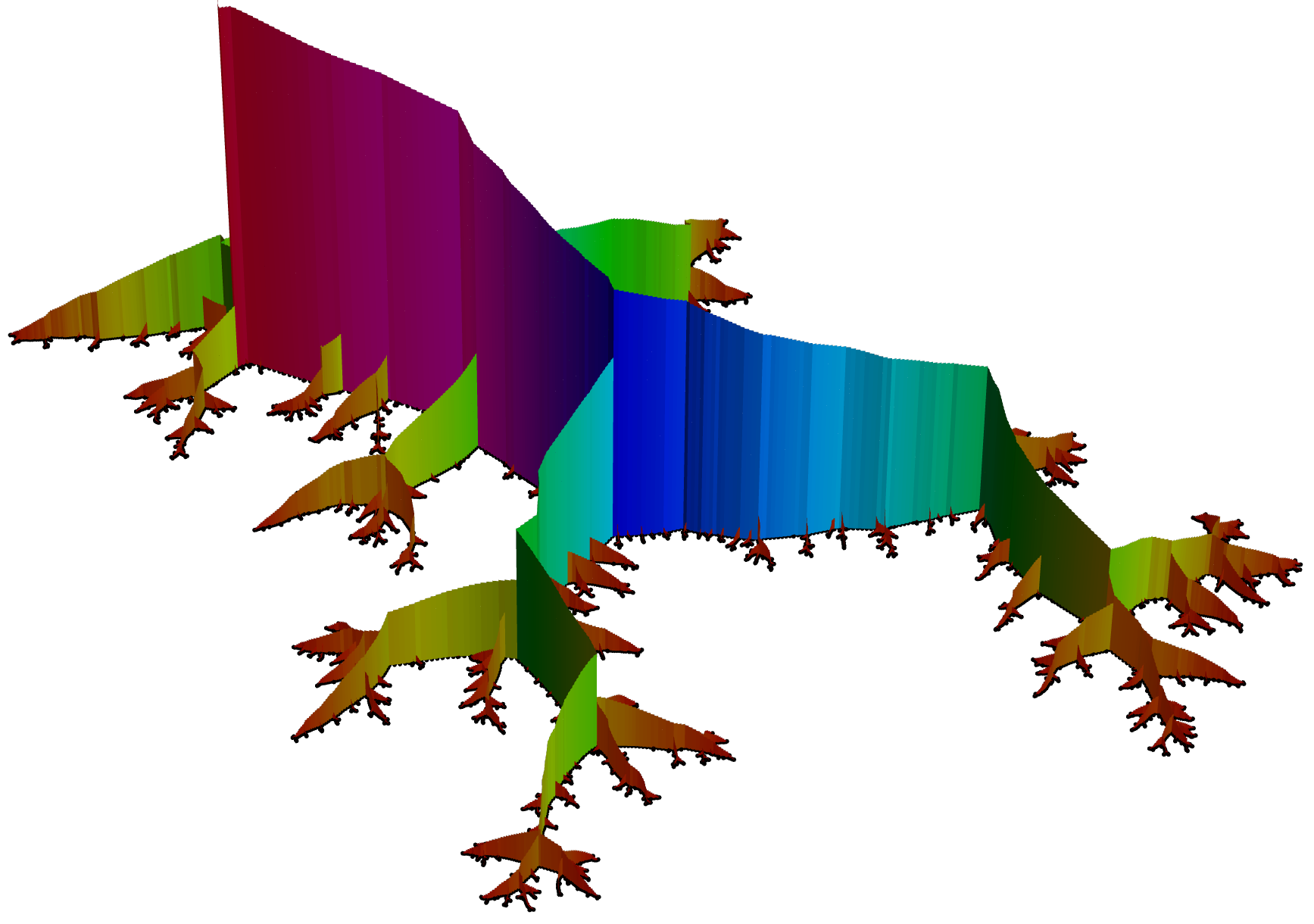}
   \caption{A simulation of a Brownian CRT normalized by the height. The tree is embedded non-isometrically in $ \mathbb{R}^{2}$; the decoration function represents the height of fringe subtrees and is depicted in the vertical coordinate.}  \label{fig:Dino}
   \end{center}
   \end{figure}  
  
  As a first step, recall from a well-known result of David  Williams, that the  Brownian excursion with height $1$ can be constructed by gluing back to back  the trajectories of two independent copies of 
 a   $3$-dimensional Bessel process started from $0$ and killed when hitting $1$, say   $R=(R(s))_{0\leq s \leq z}$. See \cite[Theorem 4.5 on page 499]{RY99}. We observe from Brownian excursion theory and classical relations between the Brownian motion and the $3$-dimensional Bessel process also due to Williams, that if we set $\sigma_u\coloneqq \sup\{s\geq 0: R(s)\leq u\}$ for the last passage time of $R$ below some level $u\in(0,1)$, then the following assertion holds. The point process of sub-excursions of $R$ above its future infimum, whose atoms are given by
  $$\left( u, (R_{\sigma_{u-}+s}-u: 0\leq s \leq \sigma_{u}-\sigma_{u-})\right)$$
  for those $u$ such that $\sigma_{u}>\sigma_{u-}$, 
  is a Poisson random measure with intensity 
 $$\indset{h(\e) < 1-u} \dd u\,  2 \mathbb{N}_{ \mathrm{Ito}}(\dd \e).$$
 Here   $\e$ stands for a generic excursion, $h(\e)=\max \e$ for the height of $\e$, and the restriction $h(\e) < 1-u$ stems from  the fact that the maximum of $R$ during the excursion interval $[\sigma_{u-}, \sigma_u]$ must be less than $1$. See for instance \cite{Pitman-Winkel}, and also \cite{Abraham-Delmas} for a more general version in the setting of  L\'evy's CRT. 
 If we think of $R$ as the left-contour of $\mathcal {T}^{(1)}$ until the highest vertex is reached, then an atom of this point measure, say at $(u,\e)$, should be interpreted as follows. The vertex at height $u$ on the branch   from the root to the highest vertex is a branchpoint of $\mathcal {T}^{(1)}$, and a contour function of the left subtree that stems from this vertex is induced by the excursion $\e$. 
 
 Now from \eqref{def:itobro} we have $2 \mathrm{N}_{\mathrm{Ito}}(h(\e)\in \dd x)= x^{-2} \dd x$. 
 Since by Williams' decomposition, the full contour process of $\mathcal {T}^{(1)}$ is obtained from 
  two independent copies of $R$,
  one now readily sees from the Brownian scaling property that the process that records the heights of the collection of the subtrees of $\mathcal {T}^{(1)}$ above height $t\in[0,1]$ is a general branching process in the sense of 
  Section \ref{sec:2.1}, whose decoration-reproduction kernel can be described as follows.
  For every $x>0$, under $P_x$, 
the decoration process is simply given by $X(t)=x-t$ for $0\leq t < x$, and the reproduction process $\eta(\dd u, \dd y)$ is a Poisson random measure on the triangle $\{(u,y): 0< y < x- u < x\}$ with intensity $2y^{-2}\dd u \dd y$
(the factor $2$ stems from the fact that two independent copies of $R$ are needed to
construct the contour function of $\mathcal {T}^{(1)}$). One immediately checks that the family of laws
$(P_x)_{x>0}$ is self-similar with exponent $\alpha=1$ in the sense of  Definition \ref{D:SSGBP}, and we shall now argue that it is actually one of the self-similar Markov decoration-reproduction kernels devised in Section   \ref{sec:2.2}.

In this direction, set $\sigma^2=0$, $\mathrm{a_{can}}= -1$ and $\alpha=1$, and 
introduce the generalized L\'evy measure $\boldsymbol{\Lambda}_{\mathrm{Height}}$ on $\mathcal S$ 
 given by
 $$\int_{\mathcal S} F(y_0, (y_1, y_2, \ldots)) \boldsymbol{\Lambda}_{\mathrm{Height}}(\dd y, \dd\mathbf y)= 2\int_{-\infty}^0 F(0,(y, -\infty, \ldots))  \e^{-y} \dd y,$$
 where $F$ denotes  a generic nonnegative functional on $\mathcal S$.  In particular, notice that we have $\mathrm a_{\mathrm{can}}=\mathrm a=-1$. Just as in Section \ref{sec:2.2},
 consider a Poisson random measure $\mathbf{N}(\dd t, \dd y, \dd \mathbf y)$ on $[0,\infty)\times \mathcal S$ with intensity $\dd t \boldsymbol{\Lambda}_{\mathrm{Height}}(\dd y, \dd\mathbf y)$.
 In this setting, the L\'evy process is merely linear, $\xi(t)=-t$ for $t\geq 0$, and  the exponential functional 
 $\epsilon(t)=1-\e^{-t}$ has inverse $\tau(s)=-\log(1-s)$ for $0\leq s <1$.
 After the Lamperti transformation, the positive self-similar Markov process is hence $X(t)=1-t$ for $0\leq t < 1$
 as we wished. 
Next observe that the push-forward image of the measure $\dd t \e^{-y} \dd y$ on $\R_+\times (-\infty,0)$ by the map
$(t,y)\mapsto (u,x)=(1-\e^{-t},\e^{-t+y})$ is the measure $x^{-2}\dd u \dd x$ on the triangle
 $\{(u,x): 0< x < 1- u < 1\}$.  It follows from the mapping theorem for Poisson point processes that the reproduction process defined by \eqref{Eq:eta:Markov} has the same distribution as the reproduction process $\eta$ above under $P_1$. Putting the pieces together, we can now conclude that the decorated tree $(\mathcal T^{(1)}, v\to \mathrm{Height}(\mathcal T^{(1)}_{v}), 0)$ is a self-similar Markov tree with characteristic quadruplet $(0, -1, \boldsymbol{\Lambda}_{\mathrm{Height}}; 1)$. 
 
As a consequence, we compute the cumulant by \eqref{eq:kappasimple}, which is simply 
  \begin{eqnarray} \label{eq:kappaheightbrow}\kappa_{\mathrm{Height}}(\gamma) = -\gamma +2\int_{-\infty}^0 \e^{(\gamma-1)y}\dd y = -\gamma + 2/(\gamma-1),\qquad \gamma>1.  \end{eqnarray}
 It follows that $\omega_-=2$ and $\omega_+=\infty$ and Assumption \ref{A:omega-} holds.
 We will now identify the harmonic  measure $\upmu$ in terms of the contour measure $\gamma_{\e^{(1)}}$ on $\mathcal T^{(1)}$, which has been defined above as  the push-forward image of the Lebesgue measure  by the contour function induced by $\e^{(1)}$.  
We claim that there is the identity
 \begin{equation}\label{eq:identmeasures}
 \upmu= \frac{3}{2} \gamma_{\e^{(1)}}.
 \end{equation}

  To establish this assertion, recall that the expectation of the lifetime $z$ of the killed Bessel process $R$ equals $1/3$. Therefore, the total mass of contour measure has expectation  $\E(\gamma_{\e^{(1)}}(\mathcal T^{(1)}))=2/3$; moreover we have also  $\mathrm{Var}(\gamma_{\e^{(1)}}(\mathcal T^{(1)}))<\infty$. 
  Now recall that $\chi(v)$ denotes the type of the individual labelled by $v\in \U$ in the general branching process,
  so that in the present setting, $\chi(v)$ is the height of some sub-excursion of $\e^{(1)}$.  We then deduce from the branching property and self-similarity that
 \begin{align*}
 \E\left( \left|  \gamma_{\e^{(1)}}(\mathcal T^{(1)})- \frac{2}{3}M_n(\omega_-) \right|^2\right) &=  \E\left(  \sum_{|v|=n}  \left|  \gamma_{\e^{(1)}}(\mathcal T_v^{(1)})- \frac{2}{3}  \chi(v)^2 \right|^2\right)\\
 &\leq c\cdot \E\left(  \sum_{|v|=n}  \chi(v)^4 \right)\\
 &\leq c\cdot \left(1+\kappa(4)/4\right)^n,
 \end{align*}
where the last line stems from Lemma \ref{L:verCMJ1}. Since $\kappa(4)<0$, we deduce that
$$\upmu(\mathcal T^{(1)}) = \lim_{n\to \infty} M_n(\omega_-) =  \frac{3}{2} \gamma_{\e^{(1)}}(\mathcal T^{(1)}).$$
Again by the branching property, we have more generally that
$$\upmu(\mathcal T_v^{(1)})= \frac{3}{2} \gamma_{\e^{(1)}}(\mathcal T_v^{(1)}), \qquad v\in \U,$$
and \eqref{eq:identmeasures} follows.  The law of $\gamma_{\e^{(1)}}(\mathcal T^{(1)})$ is itself well-known through its Laplace transform, 
$$ \mathbb{E}\left( \exp(-\lambda \gamma_{\e^{(1)}}(\mathcal T^{(1)})) \right) = \left(\frac{ \sqrt{2 \lambda}}{ \sinh( \sqrt{2\lambda})}\right)^{2},$$
see \cite[last equation]{LG10}
\hfill $\diamond$
 \end{example}

We next discuss an important family of non-increasing self-similar Markov trees that have received much attention in the literature; see for instance the survey \cite{Ber06} and references therein. A \textbf{self-similar fragmentation} process can be thought of as a model for an inert length that splits as time passes into smaller and smaller pieces called fragments. One assumes the branching property and self-similarity, in the sense that different fragments evolve independently the ones from the others and according to the same dynamics up to a proper rescaling of time and size. 
Intuitively speaking, the inertia of the length falling apart yields a natural genealogical structure of the family of fragments:  we view  a fragment present at time $t$ as a forebear of another fragment present at time $t'>t$ if the latter was part of the former at time $t$. Hence a fragmentation can be depicted by a continuous genealogical tree where branching points represent the sudden dislocations of a fragment, endowed with a decoration which records sizes.

Roughly speaking, self-similarity entails that the rate at which a fragment with size $x>0$ breaks into a (possibly infinite)
sequence of sizes $xs_1, xs_2, \ldots$ is given by $x^{-\alpha} \boldsymbol{\Xi}(\dd \mathbf s)$, where $\mathbf{s}=(s_1, s_2, \ldots)$ is a non-increasing sequence in $[0,1)$ and  $\boldsymbol{\Xi}$ is known as the \textbf{dislocation measure}. Beware that the convention  for the self-similarity exponent $\alpha$ that we follow in this text has the opposite sign of the one used in the literature on self-similar fragmentations, see e.g. \cite{haas2018scaling,HM04}. So for $\alpha>0$, small fragments break up at a higher pace than larger ones, hence faster and faster as time passes. An important consequence is that  the entire length becomes fully shattered after a finite time; in other words, there are no more fragments with positive length present in the system after a while.

The assumption of inertia of the length falling apart translates into the requirement that the sum of the sizes of the smaller fragments after a dislocation event never exceeds that of the fragment before the dislocation, and can even be strictly smaller in the case where some dust (fragments of infinitesimal sizes) is produced. So more precisely, $\boldsymbol{\Xi}$  is a  measure on the space of non-increasing sequences  $\mathbf{s}=(s_i)_{i\geq 1}$ in $[0,1)$ with $\sum_{i\geq 1} s_i \leq 1$, which furthermore fulfills the integral condition
\begin{equation} \label{Eq:xiint}\int (1-s_1) \boldsymbol{\Xi}(\dd \mathbf s) <\infty.
\end{equation}
Note that this requirement allows the dislocation measure  $\boldsymbol{\Xi}$ to be infinite, that is, the 
 fragmentation to have an infinite activity. 
If  \eqref{Eq:xiint}  failed, then the intensity of dislocations would then be too strong and any length would be instantaneously reduced to dust. 

Self-similar fragmentations naturally yield non-increasing self-similar Markov trees. More precisely, 
the dislocation measure $\boldsymbol{\Xi}$  of a fragmentation is related to the generalized L\'evy measure $\boldsymbol{\Lambda}$ as follows. The former is the push-forward image of the latter by the function 
which maps $(y_0,\mathbf y)\in \mathcal S$ to $ (s_i)_{i\geq 1}$, the version  ranked in  the non-increasing order of sequence $(\e^{y_0}, \e^{y_1}, \e^{y_2}, \ldots)$.  Note that \eqref{Eq:xiint} follows from \eqref{Eq:lambda0int} and that 
$$\boldsymbol{\Lambda} \left( \left\{(y_{0},(y_{1}, ... )) \in \mathcal{S}:  \sum_{j=0}^{\infty} \e^{y_j}  > 1 \right\} \right)= \boldsymbol{\Xi}\left( \left\{\mathbf s : \sum_{i=1}^{\infty} s_i>1\right\} \right) =0;$$
one then says that $\boldsymbol{\Lambda}$ is  \textbf{dissipative}. 
The last parameter of the model is the so-called \textbf{erosion coefficient}, a non-negative real number which accounts for rate at which fragments  continuously shrink and that is simply identified as the negative of the canonic drift coefficient  $\mathrm a_{\mathrm{can}}$. 
Putting the pieces together, a self-similar Markov tree encodes a  self-similar fragmentation if and only if it is  non-increasing and its generalized L\'evy measure is dissipative. Note from \eqref{eq:kappasimple} that the cumulant can also be expressed in terms of the dislocation measure as
$$ \kappa(\gamma) = \mathrm a_{\mathrm{can}}\gamma + \int \ \boldsymbol{\Xi}( \dd \mathbf s) \left( \sum_{i\geq 1}  s_{i}^{\gamma} -1\right).
$$
By dissipativity,   $\kappa$ is a non-increasing function with $\kappa(1)\leq 0$, and in particular, Assumption \eqref{A:gamma0} holds for $\gamma>1$ as soon as $\boldsymbol{\Xi}$ is not degenerated, i.e. $\boldsymbol{\Xi}\neq \delta_{(1,0,\dots)}$.

 The cumulant $\kappa$ has a simple probabilistic interpretation in terms of the so-called tagged fragment. 
 The latter is the process  that, as time passes,  records the size of the fragment that currently contains a point that has been initially tagged uniformly at random  in the length. It is easily seen from the branching property that this tagged fragment has the Markov property, and it is also naturally self-similar. The L\'evy process that underlies the latter via the Lamperti transformation is non-increasing, hence the negative of a subordinator. In this framework,  the Laplace exponent of this subordinator is 
 given by the function $\gamma \mapsto -\kappa(\gamma+1)$; see e.g. \cite[Corollary 3.1]{Ber06}. This observation is in close relation to the forthcoming Lemma \ref{L:LKtilde}.

A fragmentation is called pure when there is no erosion, that is when the canonic drift coefficient is  $\mathrm a_{\mathrm{can}}=0$, 
and then \textbf{conservative} when further total sizes are preserved at dislocation events, that is when the dislocation measure satisfies
$$ \boldsymbol{\Xi}\left( \left\{ \mathbf s=(s_i)_{i\geq 1} :  \sum_{i=1}^{\infty} s_{i} <1 \right\}\right)=0,$$
or equivalently in terms of the generalized L\'evy measure
  \begin{eqnarray} \label{def:conservativeLM} \boldsymbol{\Lambda} \left( \left\{(y_{0},(y_{1}, ... )) :  \sum_{j=0}^{\infty} \e^{y_j}  \ne 1 \right\} \right) =0.  \end{eqnarray}
 Example \ref{Ex:stat} is an elementary case of a conservative fragmentation.

For any conservative self-similar fragmentation, we have $\kappa(1)=0$, so Assumption \ref{A:omega-} holds with $\omega_-=1$ and $\omega_+=+\infty$. Much more precisely,
conservativeness readily entails in terms of the reproduction process $\eta$ that there is the identity
$$\int_{[0,\infty)\times (0,\infty)} x \eta(\dd t, \dd x)=1 \qquad\text{a.s.,}$$
and therefore the intrinsic martingale is trivial, $M_n(\omega_-)\equiv 1$. As a consequence, the decoration of a self-similar conservative fragmentation simply assigns to any vertex of the tree which is not a branching point, the $\upmu$-mass of the subtree that stems from this vertex.

Genealogical trees of self-similar fragmentations have been constructed first by Haas \& Miermont \cite{HM04} in the conservative setting. They were notably able, under mild hypotheses, to compute Hausdorff dimensions as well as the maximal H\"older exponents of the height functions. Specifically, they showed  that the Hausdorff dimension of the set of  leaves is $1/\alpha$. We will extend this result to our general framework of ssMt in Section \ref{sec:Hausdim}, and specifically in Proposition \ref{Prop:spine:Haus:bis}  (actually, \cite{HM04} has a mild assumption that is not needed Proposition  \ref{Prop:spine:Haus:bis}). 
Then this was extended to more general self-similar fragmentations where the conservativeness condition is dropped by Stephenson \cite{stephenson2013general}. We also mention that in the case ${\mathrm a}_{\mathrm{can}}=0$ without erosion, the total length measure $\uplambda^\gamma(T)$ has been introduced and studied in  \cite{bertoin2011area}.

 We already dealt in Example \ref{ex:brownianheight} (see also the discussion after Example \ref{ex:reduced}) with a variation  of the most distinguished member of the family of self-similar fragmentation trees, namely the ubiquitous Brownian Continuum Random Tree \cite{Ald91a}.  Here is a precise discussion. 

\begin{example}[Brownian CRT and its mass fragmentation] \label{ex:brownian} The real tree $ \mathcal{T}_{1}$ constructed from a standard Brownian excursion of length $1$, say  $(\e_{1}(s))_{0\leq s \leq 1}$ with the law $ \mathbb{N}_1$ defined in \eqref{def:itobro}, is known as the\textbf{ Brownian CRT} \footnote{Beware that, as in \cite{LG93} and many works in that field, there is a difference by a factor 2 between our definition of the Brownian CRT and that in \cite{Ald91a}.}. We  denote by $\gamma_{ \e_{1}}$ its contour measure and endow $ \mathcal{T}_{1}$ with the decoration which assigns to each vertex $ v \in \mathcal{T}_{1}$ the contour-mass $\gamma_{\e_{1}}( \mathcal{T}_{1,v})$ of the fringe-subtree above point $v$ (this is indeed a usc decoration), see Figure \ref{fig:massbro}. 

    \begin{figure}[!h]
   \begin{center}
      \includegraphics[height=6cm]{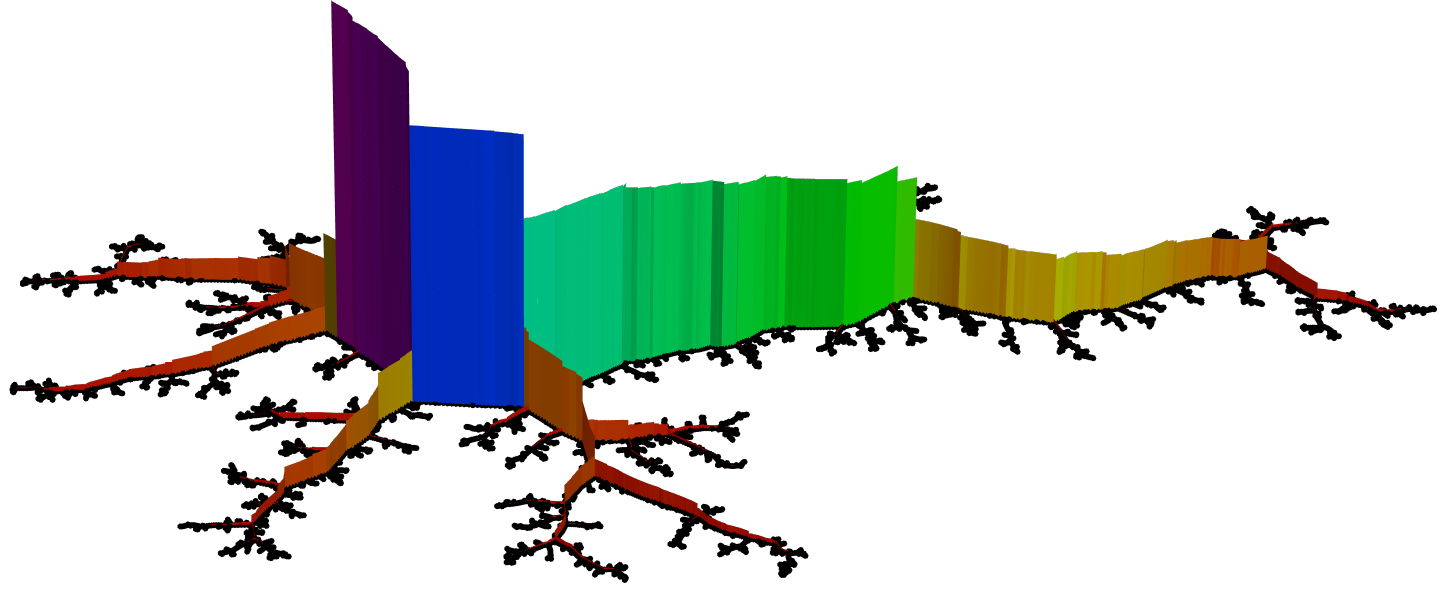}
   \caption{A simulation of a Brownian CRT. The tree is embedded (non-isometrically in $ \mathbb{R}^{2}$) and the decoration function representing the $\upmu$-mass above each point is depicted in the vertical coordinate. \label{fig:massbro}}
   \end{center}
   \end{figure}
   
   Just as  in Example \ref{ex:brownianheight}, this decorated random real tree is a self-similar Markov tree. By construction, it encodes a self-similar fragmentation whose state at time $t\geq 0$ 
   is given by the ranked sequence of decorations on the sphere $\{v\in  \mathcal{T}_{1}: d(\rho, v)=t\}$. 
Equivalently in terms of the excursion $\e_1$, for every $t\geq 0$, the random open set $\{s\in [0,1]: \e_1(s)>t\}$ can be decomposed into a (possibly empty) sequence of open intervals; the  process in the variable $t$ that records the sequence of the lengths of these intervals ranked in the decreasing order, is
the Brownian fragmentation. The characteristics are found using \cite[pages 338-340]{Ber02} and \cite{uribe2009falling}. 

 The Brownian fragmentation tree is self-similar with index $\alpha=1/2$, no erosion,  and binary  conservative  dislocation measure $ \boldsymbol{ \Xi}_{ \mathrm{Bro}}$ given by 
  \begin{eqnarray} \label{eq:GLMbrownianfrag} \int F(  s_1, s_2, \ldots )  \  \boldsymbol{ \Xi}_{ \mathrm{Bro}}(\dd \mathbf{s}) \coloneqq    \sqrt{\frac{2}{\pi}} \int_{1/2}^{1} F(x,1-x,0,0, ...) \frac{ \mathrm{d}x}{(x(1-x))^{3/2}},  \end{eqnarray}
where $F$ stands for a generic nonnegative functional of the sequence of fragments. 
The cumulant function can be evaluated using \eqref{eq:kappasimple}:
  \begin{eqnarray} \label{eq:kappabrownian}  \kappa_{ \mathrm{Bro}}(\gamma) =  \sqrt{ \frac{2}{\pi}}\int_{0}^{1}  (x^{\gamma}+(1-x)^{\gamma}-1) \frac{ \mathrm{d}x}{(x(1-x))^{{3/2}}} =   -2 \sqrt{2} \frac{\Gamma(\gamma- 1/2)}{\Gamma(\gamma-1)}, \quad \gamma>1.\end{eqnarray}
As for all conservative fragmentations, we have $\omega_-=1$ and $\omega_+=\infty$,  Assumption \ref{A:omega-} holds and the total $\upmu$-mass is trivially $x$ under $ \mathbb{P}_x$. In particular, we have the almost sure equality $\upmu = \gamma_{ \e_{1}}$ so that the decoration in $ \mathcal{T}_{1}$ can be recovered from its harmonic measure and vice-versa.

A possible choice for the generalized L\'evy measure in the family of bifurcators is then obtained by following the largest fragment at each dislocation. This amounts to defining first the measure
 $ \Lambda_{ \mathrm{Bro}, \max}$ on $\R_-$ as the push-forward image of the dislocation measure $ \boldsymbol{ \Xi}_{ \mathrm{Bro}}$ by the map $\mathbf s \mapsto \log x$ with $\mathbf s=(x,1-x, 0, \ldots)$ for $x\in[1/2,1)$. Concretely, this gives
 \begin{align*}
 \int_{\R_-} f(y)  \Lambda_{ \mathrm{Bro}, \max}(\dd y) &=  \sqrt{\frac{2}{\pi}} \int_{1/2}^{1} f(\log x) \frac{ \mathrm{d}x}{(x(1-x))^{3/2}}\\
 &=  \sqrt{\frac{2}{\pi}} \int_{-\log 2}^{0} f(y) \frac{ \dd y }{\sqrt{\e^y(1-\e^y)^3}}.
 \end{align*}

Then the generalized L\'evy measure $\boldsymbol{\Lambda}_{ \mathrm{Bro}, \max}$ is simply obtained by
 $$\int_{\R\times \mathcal{S}_1}
 F(y_0, \mathbf{y}) \boldsymbol{ \Lambda}_{ \mathrm{Bro}, \max}
 (\dd y, \dd \mathbf y)=  
\int_{\R_-} F(y, (\log(1-\e^{y}), -\infty, \ldots)) { \Lambda}_{ \mathrm{Bro}, \max} (\dd y),
$$
where now $F$ stands for a generic nonnegative functional on $\mathcal S$. 

   Notice also that using the fact that the Brownian CRT is coded by the Brownian excursion \cite{LG93,bertoin2000fragmentation}, the quantity $\uplambda^{3/2}(T) = \int_{T}  \mathrm{d} \lambda\ g$ can be interpreted as the area under a standard Brownian excursion, known as the Airy law and whose moments are explicit, see \cite{janson2007brownian} for details. Other length measures appeared recently in the literature \cite{fill2022sum,delmas2018cost,abraham2022global}. This Brownian fragmentation has been studied in depth in the literature and we did not try here to survey all its known properties, see e.g. \cite{borga2023power} for a recent application to the study of increasing subsequences in the Brownian separable permuton. Finally, it is interesting to recall that Aldous \cite{Ald91} has a construction of the Brownian CRT by recursive gluing of line segments that is somewhat similar to the construction in Example \ref{Ex:stat}. 
\hfill $\diamond$
   \end{example}

 One says that a conservative self-similar fragmentation is \textbf{binary} if, just like as in Example~\ref{ex:brownian},  exactly two fragments are produced at each dislocation event. Specifically, the dislocation measure  must satisfy
$$ \boldsymbol{\Xi}\left( \left\{ \mathbf s=(s_i)_{i\geq 1} :  s_{3}>0 \right\}\right)=0,$$
that is equivalently in terms of the generalized L\'evy measure to asking
\begin{equation}\label{Eq:binaryLM}
\boldsymbol{\Lambda} \left( \left\{(y_{0},(y_{1}, ... )) :  {y_3} \ne -\infty  \right\} \right) =0.
\end{equation}
We already observed in Example \ref{ex:brownian} that in the binary conservative case, the generalized L\'evy measure $\boldsymbol{ \Lambda}$ is entirely determined by its first marginal $\Lambda_0$ and  the identity
\begin{equation}\label{Eq:binaryLM2}
\int_{\R\times \mathcal{S}_1}
 F(y_0, \mathbf{y}) \boldsymbol{ \Lambda}
 (\dd y, \dd \mathbf y)=  
\int_{\R_-} F(y, (\log(1-\e^{y}), -\infty, \ldots)) { \Lambda}_0 (\dd y).
\end{equation}
Moreover, the reproduction process $\eta$ is  determined by the decoration process $X$ since we have 
\begin{equation}\label{Eq:reprodbinary}
\eta= \sum_{t>0} \delta_{(t,-\Delta X(t))}.
\end{equation}
Pitman and Winkel \cite{pitman2009regenerative} and \cite[Section 3]{pitman2015regenerative} devised first a recursive construction of binary conservative self-similar fragmentation trees using a so-called bead splitting processes that generalizes Aldous' line-breaking construction of the Brownian CRT.

   The stable CRT, introduced by Le Gall -- Le Jan \& Duquesne \cite{LGLJ98,DLG05}, form a one-parameter family indexed by $\beta\in(1,2]$ of random continuous trees, where the boundary case $\beta=2$ is a distance-rescaled version of  the Brownian CRT.  They appear as scaling limits of critical Galton--Watson trees whose offspring distribution belong to the domain of attraction of a spectrally positive stable law with index $\beta$. They also belong to the family of self-similar fragmentation trees of Haas and Miermont \cite{HM04}, and 
   satisfy the striking property that their distribution is invariant under  re-rooting at a random $\upmu$-point \cite{HMPW08}.
  A notable difference in the case $\beta<2$ is that branchpoints have infinite degrees almost surely, whereas branching are always binary for the Brownian CRT.

Just as in Examples \ref{ex:brownianheight} and \ref{ex:brownian}, the stable trees can be interpreted as self-similar Markov trees. Formally, for any $\beta \in (1,2]$, consider the excursion measure $ \mathbb{N}_{ \beta-\mathrm{stable}}$ of the $\beta$-stable L\'evy process with no negative jumps above its running infimum. For normalization, we suppose that the underlying L\'evy process $\xi$ has L\'evy-Khintchine exponent given by $ \mathbb{E}[ \exp(-\lambda \xi_{t})] = \exp(t \lambda ^{\beta})$ for $\lambda>0$ or equivalently that the L\'evy measure is given by 
   $$ \pi_{\beta}( \mathrm{d}r) = \frac{\beta(\beta-1)}{\Gamma(2- \beta)} \frac{ \mathrm{d}r}{r^{1 + \beta}}.$$
Following Le Gall and Le Jan, from the excursion process $\mathrm{e}$, one can construct another process (with the same duration), called the height process and denoted by $\mathrm{h} : \mathbb{R}_{+} \to \mathbb{R}_{+}$; see \cite{DLG02}. This process then encodes, as above, a random real tree $\mathcal{T}_{\mathrm{h}}$. As in \eqref{def:itobro}, the infinite excursion measure $\mathbb{N}_{\beta\text{-}\mathrm{stable}}$ can be disintegrated according to the length $\ell$ of the height process $\mathrm{h}$ as follows 
  \begin{eqnarray} \label{def:itostable}  \mathbb{N}_{ { \beta-\mathrm{stable}}} = \int_{0}^{\infty} \frac{ \mathrm{d}h}{(  (\beta-1) h)^{ 1 + \frac{1}{\beta-1}}} \mathbb{N}^{(h)}_{ \beta- \mathrm{stable}} = \int_{0}^{\infty} \frac{ \mathrm{d}\ell}{ \beta \Gamma(1- \beta^{-1}) \ell^{1+ \frac{1}{\beta}} } \mathbb{N}^{\beta- \mathrm{stable}}_{\ell}, \end{eqnarray}
 see \cite[Theorem 3.3]{Abraham-Delmas} or \cite[Section 3.2]{nassif2022zooming}. We stress that in the case  $\beta=2$, under $\mathbb{N}^{\beta- \mathrm{stable}}_{\ell}$, the process $2^{-1/2}\cdot \mathrm{h}$ is distributed according to $\mathbb{N}_{\ell}$.  We then have the analogs of Examples \ref{ex:brownianheight} and \ref{ex:brownian} in the stable case:

\begin{example}[Stable trees and their mass fragmentations] \label{ex:stable} The real tree $ \mathcal{T}_{ \mathrm{h}}$ when the height-process $ \mathrm{h}$ follows the law $\mathbb{N}^{\beta- \mathrm{stable}}_{1}$ is known as the  \textbf{ $\beta$-stable CRT} (of mass $1$). As in Example~\ref{ex:brownian}, we endow it  with the decoration given by the mass measure of its fringe trees, where the measure is the push-forward of the Lebesgue measure by the coding by $ \mathrm{h}$.
It then follows from the work of Miermont \cite{Mie03} that this is a self-similar Markov tree with  self-similar exponent 
$\alpha=1-1/\beta$, canonical drift $ \mathrm{a_{can}}=0$ and  no Brownian component. The generalized L\'evy measure  $\boldsymbol{ \Lambda}^{\beta \mathrm{-stable}}_{ \mathrm{mass}}$ is conservative but non-binary, and can be related to the Poisson-Dirichlet measure with parameter $(1/\beta, -1)$. More precisely, it can be given by 
$$\int  F\big(  \mathrm{e}^{{y_{0}}}, (\mathrm{e}^{y_{1}}, ...)\big)  \  \boldsymbol{ \Lambda}^{\beta \mathrm{-stable}}_{ \mathrm{mass}}(\mathrm{d} y_{0}, \mathrm{d} \mathbf{y})  \underset{ \text{\cite[Thm \  1]{Mie03}}}{:=}   \quad \frac{\beta^{2} \Gamma(2-1/\beta)}{\Gamma(2 -\beta)}\cdot   \mathbb{E}\left[ S_1 F \left( \frac{\Delta S_{t_i}}{S_1} : i \geq 1\right) \right],$$  where $\Delta S_{t_i}$ are the jumps of a  stable $1/\beta$ subordinator $(S_t)_{0\leq t \leq 1} $ started from $0$ on the unit time interval, ranked in decreasing order. 

\begin{figure}[!h]
 \begin{center}
 \includegraphics[width=13cm,angle=-1]{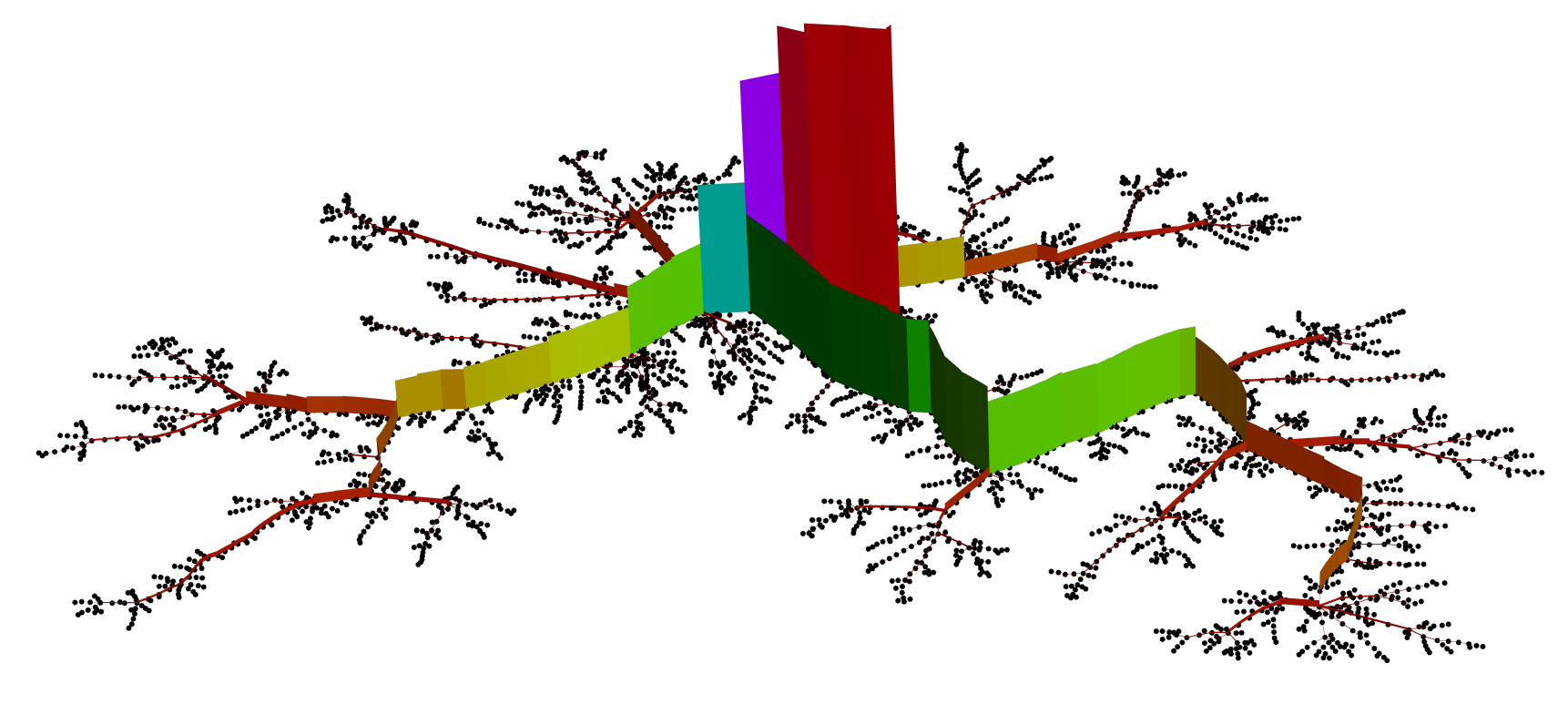}
 \caption{A simulation of a $3/2$-stable fragmentation tree embedded in the plane with the decoration in the third coordinate.}
 \end{center}
 \label{fig:3/2stabletree}
 \end{figure}

As for all conservative fragmentations we have $\omega_{-}=1$ and $\omega_+= + \infty$, Assumption \ref{A:omega-} holds with trivial total $\upmu$ mass, and as in Example \ref{ex:brownian}, $\upmu = \gamma_{ \mathrm{h}}$ and the decoration $g(x)$ is the $\upmu$-measure of the subtree above point $x \in T$. The length measures on $T$ have also been considered in \cite{fill2022sum,delmas2018cost,abraham2022global}. 

 The exact form of the cumulant function was computed in Section 3.4 in \cite{Mie03}, which established that
$$ \kappa_{\beta- \mathrm{stable}}(\gamma) \quad = \quad  - \beta\cdot \frac{\Gamma(\gamma- 1/\beta)}{\Gamma(\gamma-1)}.$$
In this normalization, the $\beta$-stable tree encodes the genealogy of the CSBP with branching mechanism $z \mapsto z^{\beta}$, and when $\beta = 2$ the associated random tree is  the Brownian CRT with distances rescaled by $\sqrt{2}$ (hence the disappearance of a factor of $ \sqrt{2}$ in \eqref{eq:kappabrownian}, see \cite{CH17} for careful normalizations). 
Last, we also refer to \cite{GH14} and \cite{rembart2018recursive} for constructions of stable CRTs by recursive random gluing of segments, in the same vein as in Aldous \cite{Ald91} and Example \ref{Ex:stat}.
\hfill $\diamond$
\end{example}

 \begin{example}[Heights of stable sub-trees] \label{ex:stableheight} As in Example \ref{ex:brownianheight} we now consider $ \mathcal{T}_{ \mathrm{h}}$ the  continuum random tree of height $1$ coded by $ \mathrm{h}$ under law $\mathbb{N}^{(1)}_{ \beta- \mathrm{stable}}$. The decoration is similarly given by the height of its fringe subtrees. It follows now from William's decomposition of the measure $ \mathbb{N}_{ \beta- \mathrm{stable}}$ established by Abraham \& Delmas \cite{Abraham-Delmas} that this yields a self-similar Markov tree. As in Example \ref{ex:brownianheight},  the self-similarity parameter is obviously $\alpha=1$, we have no Brownian part $\sigma^2=0$ and the canonical drift is $ \mathrm{a_{can}}=-1$. By \cite[Theorem 3.3]{Abraham-Delmas}\footnote{Notice the typo there where $\mathbb{N}( H_{\mathrm{max}} \leq m)$ should be replaced by $\mathbb{N}(H_{\mathrm{max}} \geq  m)$.} the generalized L\'evy measure can be described as follows. For $r >0$ fixed, we write $ \Pi_{r}$ for the law of the decreasing rearrangement of the single point $(1)$ completed with the vertical coordinates of the atoms of a Poisson measure $ \sum_{i \geq 1} \delta_{(x_{i},y_{i})}$ with intensity 
 $$ \mathrm{d}x \mathbf{1}_{[0,r]} \otimes \frac{\mathrm{d}y}{(y(\beta-1))^{1+ \frac{1}{\beta-1}}} \mathbf{1}_{[0,1]}(y).$$
 Then the generalized L\'evy measure $ \boldsymbol{ \Lambda}_{ \beta- \mathrm{stable}}^{ \mathrm{height}}$ is described as 
 $$ \int_{ \mathcal{S}}  \boldsymbol{ \Lambda}_{ \beta- \mathrm{stable}}^{ \mathrm{height}}( \mathrm{d}y_{0}, ( \mathrm{d}y_{1}, ...))  F( \mathrm{e}^{y_{0}}, ( \mathrm{e}^{y_{1}}, ...)) =  \int_{0}^{\infty} \frac{\beta (\beta-1)}{\Gamma(2-\beta)} \frac{\mathrm{d}r}{r^{\beta}} \exp\left(-r (\beta-1)^{\frac{-1}{\beta-1}}\right)  \Pi_{r}( F).$$
 As in \eqref{eq:kappaheightbrow}, the cumulant function is then easily computed from \eqref{eq:kappasimple} and equals
  \begin{eqnarray*} \kappa(\gamma) &=& -\gamma + \int_{0}^{\infty} \frac{\beta (\beta-1)}{\Gamma(2-\beta)} \frac{\mathrm{d}r}{r^{\beta}} \exp\left(-r (\beta-1)^{\frac{-1}{\beta-1}}\right) \cdot r \int_{0}^{1} \frac{\mathrm{d}y}{(y(\beta-1))^{1+ \frac{1}{\beta-1}}}  y^{\gamma}\\ &=& - \gamma + \frac{\beta}{\gamma(\beta-1)-1}, \quad \mbox{for } \gamma > \frac{1}{\beta-1}  \end{eqnarray*}
In particular, we have $\omega_-= \frac{\beta}{\beta-1}$ and $\omega_+=\infty$ and Assumption \ref{A:omega-} holds.  As in \eqref{eq:identmeasures}, the harmonic measure $\upmu$ is a constant multiple of the measure coming from the coding by $\mathrm{h}$. Let us now give a quick argument to justify it. By \cite[Sections 1.3 and 1.4]{DLG02}, for every $r> 0$, we can consider   $L_r$  the local time of $\mathrm{h}$ at height $r$, which in particular satisfies:
$$\gamma_{\mathrm{h}}(\mathcal{T}_{\mathrm{h}})=\int_0^{\infty} L_r ~\mathrm{d} r\quad \text{ and }\quad  \max \mathrm{h}= \sup\{r\geq 0:~L_r\neq 0\},  \quad \text{ under }  \mathbb{N}_{ { \beta-\mathrm{stable}}}. $$
Moreover, the process $(L_r)_{r\geq 0}$ can be characterized as follows. We have $L_0=0$, and
\begin{equation}\label{laplace:transform:N:beta:stable}
\mathbb{N}_{ { \beta-\mathrm{stable}}}\big(1-\exp(-\lambda L_r)\big)= \big(\lambda^{-(\beta-1)}+(\beta-1)r\big)^{-\frac{1}{\beta-1}}, \quad r>0 .
\end{equation}
Furthermore, conditionally on $L_r$ for some fixed $r>0$, the process $(L_{r+t})_{t\geq 0}$ is a branching process with branching mechanism  $\lambda\mapsto \lambda^\beta$. Hence,  we have:
\begin{align*}
\mathbb{N}_{ { \beta-\mathrm{stable}}}(\gamma_{\mathrm{h}}(\mathcal{T}_{\mathrm{h}}) \mathbf{1}_{\sup \mathrm{h}\leq 1})&=\int_0^1 \mathrm{d} r~\mathbb{N}_{ { \beta-\mathrm{stable}}}\big(L_r \mathbf{1}_{\sup \mathrm{h}\leq 1}\big)\\
&=\int_0^1 \mathrm{d} r~\mathbb{N}_{ { \beta-\mathrm{stable}}}\Big(L_r \exp\big(-\big((\beta-1)(1-r)\big)^{-\frac{1}{\beta-1}}\cdot L_r\big)\Big),
\end{align*}
where the last equality follows from the fact that the probability for a branching process with mechanism $\lambda\mapsto \lambda^\beta$, starting from $z>0$, to become extinct before time $t$ is given by $\exp(-\big((\beta-1)t\big)^{-\frac{1}{\beta-1}} z )$. We get from \eqref{laplace:transform:N:beta:stable} that $\mathbb{N}_{ { \beta-\mathrm{stable}}}\big(\gamma_{\mathrm{h}}(\mathcal{T}_{\mathrm{h}}) \mathbf{1}_{\sup \mathrm{h}\leq 1}\big)=\frac{\beta-1}{2\beta-1}.$  By disintegration and scaling, it follows that
$$\mathbb{N}_{ { \beta-\mathrm{stable}}}^{(1)}\big(\gamma_{\mathrm{h}}(\mathcal{T}_{\mathrm{h}})\big)= \frac{(\beta-1)^{\frac{2\beta-1}{\beta-1}} }{2\beta-1}.$$ 
This strongly suggests that $\upmu= ((2\beta-1)/(\beta-1)^{\frac{2\beta-1}{\beta-1}})\cdot  \gamma_{\mathrm{h}}$. Since  $\kappa(2\cdot \beta/(\beta-1))<0$,to complete the proof of the identity it suffices to adapt the argument for the Brownian tree (see Example~\ref{ex:brownianheight}) by establishing that $\mathrm{Var}( \gamma_{\mathrm{h}}(\mathcal{T}_{\mathrm{h}}))<\infty$, under $\mathbb{N}_{ { \beta-\mathrm{stable}}}^{(1)}$. In fact, one may notice that $\gamma_{\mathrm{h}}(\mathcal{T}_{\mathrm{h}})^2\mathbf{1}_{\sup \mathrm{h}\leq 1}=2\int_0^1\d r \int_r^1 \d s L_r L_s$ and an adaptation of the previous computations shows that $\mathbb{N}_{ { \beta-\mathrm{stable}}}(\gamma_{\mathrm{h}}(\mathcal{T}_{\mathrm{h}})^2 \mathbf{1}_{\sup \mathrm{h}\leq 1})<\infty$. The desired result then follows by disintegration and scaling. We leave the remaining details to the reader. Let us conclude this example by noting that \cite{duquesne2017decomposition} specifically study the diameter and maximum height of the $\beta$-stable CRT. In particular, equation (61) in \cite{duquesne2017decomposition} provides the exact distribution of the total height in terms of a series. By applying scaling properties and disintegration under the measure $\mathbb{N}_{\beta\text{-}\mathrm{stable}}$, it should be possible to recover the discussed results concerning the contour measure, and in particular, to derive a closed-form expression for the density law of the total mass.
\end{example}

We refer to \cite{HM12,haas2018scaling} for other examples of conservative fragmentation trees that we do not describe in these pages. Let us now present a few natural  generalized fragmentation trees that are not covered by the Haas--Miermont framework. One natural way to obtain them is to consider dissipative fragmentations obtained by trimming a (conservative) fragmentation using a local rule to keep only certain particles. Specifically, in the continuous world, a \textbf{trimming} rule is a function 
$$ \mathrm{Trim} : \left\{ \begin{array}{l} \mathcal{S} \times [0,1] \to \mathcal{S}\\
\Big( \big(y_{0}, (y_{1}, y_{2}, ... )\big) , \omega \Big) \mapsto \big( \tilde{y}_{0}, ( \tilde{y}_{1},  \tilde{y}_{2}, ...)\big), \end{array} \right.$$
which associates to a point $(y_{0}, (y_{1}, y_{2}, ... ))\in \mathcal{S}$ and an additional source of randomness $ \omega \in [0,1]$ a random variable $( \tilde{y}_{0}, ( \tilde{y}_{1},  \tilde{y}_{2}, ...))$ in $\mathcal{S}$ where $\tilde{y}_{0}$ is either $y_{0}$ or $-\infty$ (we interpret the latter case as a deletion) together with a -possibly finite\footnote{In that case, the finite sequence is completed by infinitely many's $-\infty$ for the sake of definitiveness.}- subsequence $( \tilde{y}_{1},  \tilde{y}_{2}, ...)$ excerpted from $(y_{1}, y_{2}, ... )$. Given a  trimming rule, we write
$$ \mathrm{Trim}( \boldsymbol{\Lambda}) :=  \mathrm{Trim}_{\#}  ( \boldsymbol{\Lambda} \otimes \mathrm{Leb}[0,1]),$$
for the image of a generalized L\'evy measure by this rule. 
We shall furthermore suppose that the killing rate does not explode after trimming, i.e.~that $$ \mathrm{k}_{ \mathrm{trim}} =  \mathrm{Trim}( \boldsymbol{\Lambda}) \left(  \{-\infty \} \times \mathcal{S}_{1} \right) < \infty;$$
plainly $ \mathrm{Trim}( \boldsymbol{\Lambda})$ is in turn  a generalized L\'evy measure. 
Remark that the cumulant function associated to $( \sigma^{2},  \mathrm{a} ,\boldsymbol{\Lambda} ; \alpha)$ is clearly upper bounded by that of $(  \sigma^{2},   \mathrm{a} ,\mathrm{Trim}(\boldsymbol{\Lambda}) ; \alpha)$ so that by Proposition \ref{P:constructionomega-}, we can construct  simultaneously on the same probability space $( T,d_T,\rho_T,g)$ and $( T_{ \mathrm{trim}},d_{T_{ \mathrm{trim}}}, \rho_{T_{ \mathrm{trim}}}, g_{ \mathrm{trim}})$ with characteristics $( \sigma^{2}, \mathrm{a}, \boldsymbol{\Lambda}; \alpha)$ and $( \sigma^{2}, \mathrm{a},   \mathrm{Trim}(\boldsymbol{\Lambda}) ; \alpha)$ so that  $T_{ \mathrm{trim}}$ is   a subtree of $T$ with the restriction of the decoration (the trees do not carry any measure).  Let us give a concrete example where explicit computations are doable:

\begin{example}[$k$-sampling in the Brownian case]  \label{ex:sampling} Consider the Brownian fragmentation of Example \ref{ex:brownian} which we trim as follows: Fix $k \geq 2$ and  for $(s_0,s_1) \in \mathcal{S}$, write $x_0 = \mathrm{e}^{s_0}$ and $x_1 = \mathrm{e}^{s_1}$. Next, using an independent source of randomness, sample $k$ i.i.d. Bernoulli variables $\epsilon_1, ... , \epsilon_k$ with  $ \mathbb{P}(\xi = i) = x_i$ for $i=0,1$. Then delete  $x_i$ if and only if none of the Bernoulli variables $\epsilon_j$ for $j=1, \ldots, k$ takes the value $i$. This yields a non-conservative binary self-similar fragmentation with characteristics $( \sigma^2=0,\mathrm a_{\mathrm{can}}=0, \boldsymbol{\Lambda}_{ \mathrm{Bro},k} ;  \alpha=1/2)$, where the generalized L\'evy measure is given by

\begin{align*} & \int_{\mathcal{S}}  \boldsymbol{ \Lambda}_{ \mathrm{Bro},k}( \mathrm{d}y_{0}, \mathrm{d}( y_{i})_{i \geq 1}) F\big(  \mathrm{e}^{y_0}, (\mathrm{e}^{y_1}, ...)\big) \\
&=  \sqrt{ \frac{2}{\pi}} \int_{1/2}^1 \frac{\dd x}{(x(1-x))^{3/2}} \left(\begin{array}{cl} &F\big(x,(1-x,0,...)\big) (1-x^k -(1-x)^k)\\ 
+& F\big(x,(0,0...)\big) x^k\\
+&F\big(0,(1-x,0, ...)\big) (1-x)^k\end{array} \right).
\end{align*}
In this example, we have:
\begin{itemize} 
\item For $k=2$, $$ \kappa_{ \mathrm{Bro},2}(\gamma) =   \sqrt{2}\cdot \frac{(1-2\gamma)\Gamma(\gamma+ 1/2)}{\Gamma(\gamma+1)},$$ so that Assumption \ref{A:omega-} holds with $\omega_-=1/2$ and $\omega_+ = +\infty$.  We believe that the total harmonic  measure $\upmu(T)$  should be distributed as a multiple of a Rayleigh law with density $x \mathrm{e}^{-x^2/2}$ on $ \mathbb{R}_+$. 
\item for $k=3$, $$ \kappa_{\mathrm{Bro},3}(\gamma) = \frac{(3-2\gamma(1+2\gamma)) \Gamma(\gamma+1/2)}{ \sqrt{2}\Gamma(2+\gamma)},$$
so that  Assumption \ref{A:omega-} holds with $\omega_- = \frac{ \sqrt{13}-1}{4}$ and $\omega_+ = +\infty$. The underlying ssMt is thus a rather natural subtree of the Brownian CRT with Hausdorff dimension $\frac{ \sqrt{13}-1}{2}$.
\end{itemize}
\hfill $\diamond$
\end{example}

\begin{figure}[!h]
 \begin{center}
 \includegraphics[height=7cm]{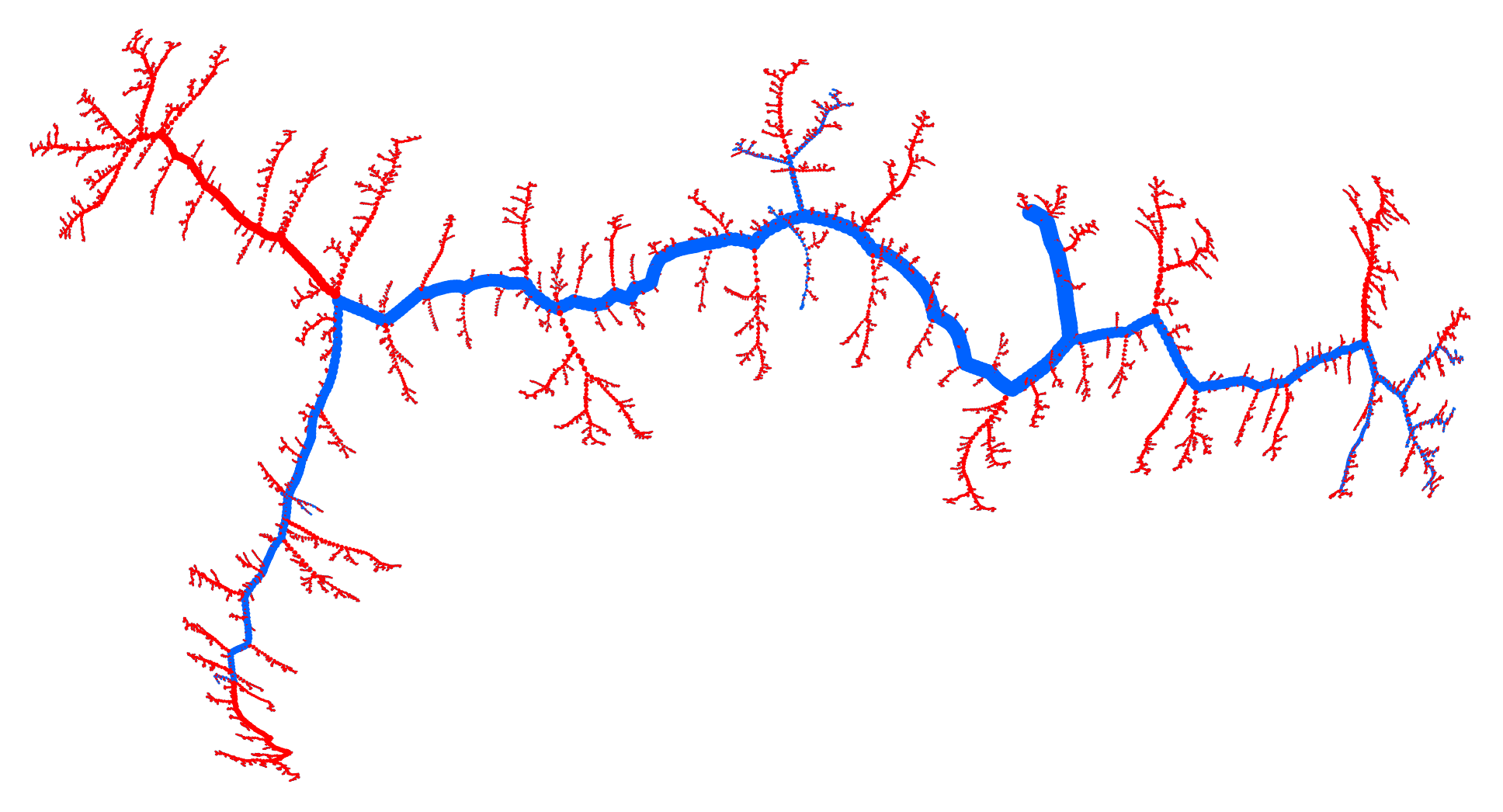}
 \caption{An illustration of the Trimmed Brownian tree  for $k=2$  (in blue) inside the standard Brownian tree of  Example \ref{ex:brownian}. The thickness indicates the labels.}
 \end{center}
 \end{figure}

\section{Conservative binary growth-fragmentations} \label{sec:5.3}
Markovian growth-fragmentations processes have been introduced in \cite{Ber15} as stochastic models describing the evolution of a system of living cells, where  at any given time cells are simply determined  by their sizes. Imagine that as time passes, cells may grow or shrink continuously and are further involved in birth events  at which a daughter cell split from the mother cell (this is called mitosis in biology). More precisely, let  $X(t)$ be the size of a typical cell at age $t$; we suppose that the process $X=(X(t))_{t\geq 0}$ is Markovian, has only negative jumps, and reaches $0$ continuously after some finite time which we view as the death of the cell. Suppose further that  each jump of $X$ is a birth event, such that if $\Delta {X}(s)=-x<0$, then $s$ is the birthtime of a daughter cell with size $x$ which then evolves independently and according to the same dynamics, i.e.  giving birth  in turn to great-daughters, and so on and so forth. 
This description fits the setting of general branching processes of Section \ref{sec:2.1}, where the reproduction process $\eta$ is simply the point process of the negative of the jumps of the decoration, see 	\eqref{Eq:reprodbinary}. We point out that the case when cells are actually inert (i.e. may split but otherwise do not grow nor shrink) corresponds to bead splitting processes introduced by Pitman and Winkel \cite{pitman2009regenerative,pitman2015regenerative}. 

Suppose further now that the Markov process $X$ is self-similar with exponent $\alpha>0$, and note from  \eqref{Eq:reprodbinary} that the reproduction process $\eta$ inherits self-similarity from $X$. Recall also that $X$ should have no positive jumps (to match the assumption that cells grow continuously and hence have only negative jumps at mitosis events) and die when reaching continuously a size $0$.   We are then in the conservative \eqref{def:conservativeLM} and binary \eqref{Eq:binaryLM} case. The L\'evy process  $\xi=(\xi(t))_{t\geq 0}$ that underlies $X$ in the Lamperti transformation is thus spectrally negative (i.e. its L\'evy measure $\Lambda_0$ is carried on $(-\infty,0)$), drifts to $-\infty$, and has also killing coefficient $\mathrm k=0$. Therefore, if we write $(\sigma^2, \mathrm a, \Lambda_0)$ for its characteristic L\'evy triplet,
then $\Lambda_0$ is a measure on $(-\infty, 0)$ with
$$\int_{(-\infty,0)}(1\wedge y^2) \Lambda_0(\dd y)< \infty \quad \text{and} \quad \mathrm a + \int_{(-\infty,-1)} y \Lambda_0(\dd y) <0.$$
We stress that the integrability condition \eqref{Eq:lambda0int} may fail, and that the second requirement above is the necessary and sufficient condition for $\xi$ to drift to $-\infty$. 
The generalized L\'evy measure is then given by \eqref{Eq:binaryLM2} and hence the cumulant by
$$\kappa(\gamma)= \frac{1}{2} \sigma^2 \gamma^2 + \mathrm a \gamma + \int_{(-\infty,0)}(\e^{\gamma y} + (1-\e^{y})^\gamma - 1 -\gamma y \indset{y\geq -1}) \Lambda_0(\dd y).$$
It can be shown that if $\kappa(\gamma)\leq 0$ for some $\gamma>0$, then almost-surely,   the family of the sizes of cells is null at any time $t\geq 0$, whereas if $\kappa(\gamma)\in (0,\infty]$ for all $\gamma>0$, then almost surely, there exist times $t>0$ at which there are infinitely many cells with size, say, larger than $1$. See \cite{Ber15} and \cite{bertoin2016local}. We stress that on top of the restriction above (binary, conservative, no killing), the main conceptual difference between the growth-fragmentation point of view and our approach in these pages is that a growth-fragmentation process is a process $ (\mathbf{X}(t) : t \geq 0)$ with values in the space of point measures on $ \mathbb{R}_+$ (or sequences of non-negative reals), whereas we focus here on the construction of random decorated trees. Specifically, when $ \kappa(\gamma)<0$ for some $\gamma>0$, that is in the subcritical case, the self-similar Markov tree with characteristics $( \sigma^{2}, \mathrm{a},  \boldsymbol{\Lambda}; \alpha)$ where $\boldsymbol{\Lambda}$ is given by \eqref{Eq:binaryLM2} can be thought of as the decorated genealogical tree underlying the self-similar growth-fragmentation process. This correspondence is similar to the one made between superprocesses and their encoding by random trees and snake trajectories \cite{DLG02}. It is still open to decide whether this correspondence also holds in the critical case, see the end of Chapter \ref{chap:generalBP}. Rembart and Winkel \cite[Section 4.2]{rembart2018recursive} were the first to describe a recursive construction of binary conservative self-similar growth-fragmentation trees that inspired our general construction by gluing in Section~\ref{sec:1.3}.

We now describe an interesting family of self-similar Markov trees for which  cumulant functions will always be expressed as ratio of Gamma functions (and simple trigonometric functions). We start with the most important example:

\begin{example}[Brownian growth-fragmentation tree] \label{ex:3/2stable}
The Brownian growth-fragmentation tree is a 
self-similar Markov tree related to a remarkable growth-fragmentation that
has appeared as the scaling limit of cactus trees inside random triangulations \cite{BCK18} or directly within the free Brownian disk and Brownian motion indexed by the Brownian tree \cite{le2020growth}.  It bears some obvious similarities with the Brownian fragmentation tree of Example~\ref{ex:brownian}.

\begin{figure}[!h]
 \begin{center}
 \includegraphics[width=12cm,angle=-3]{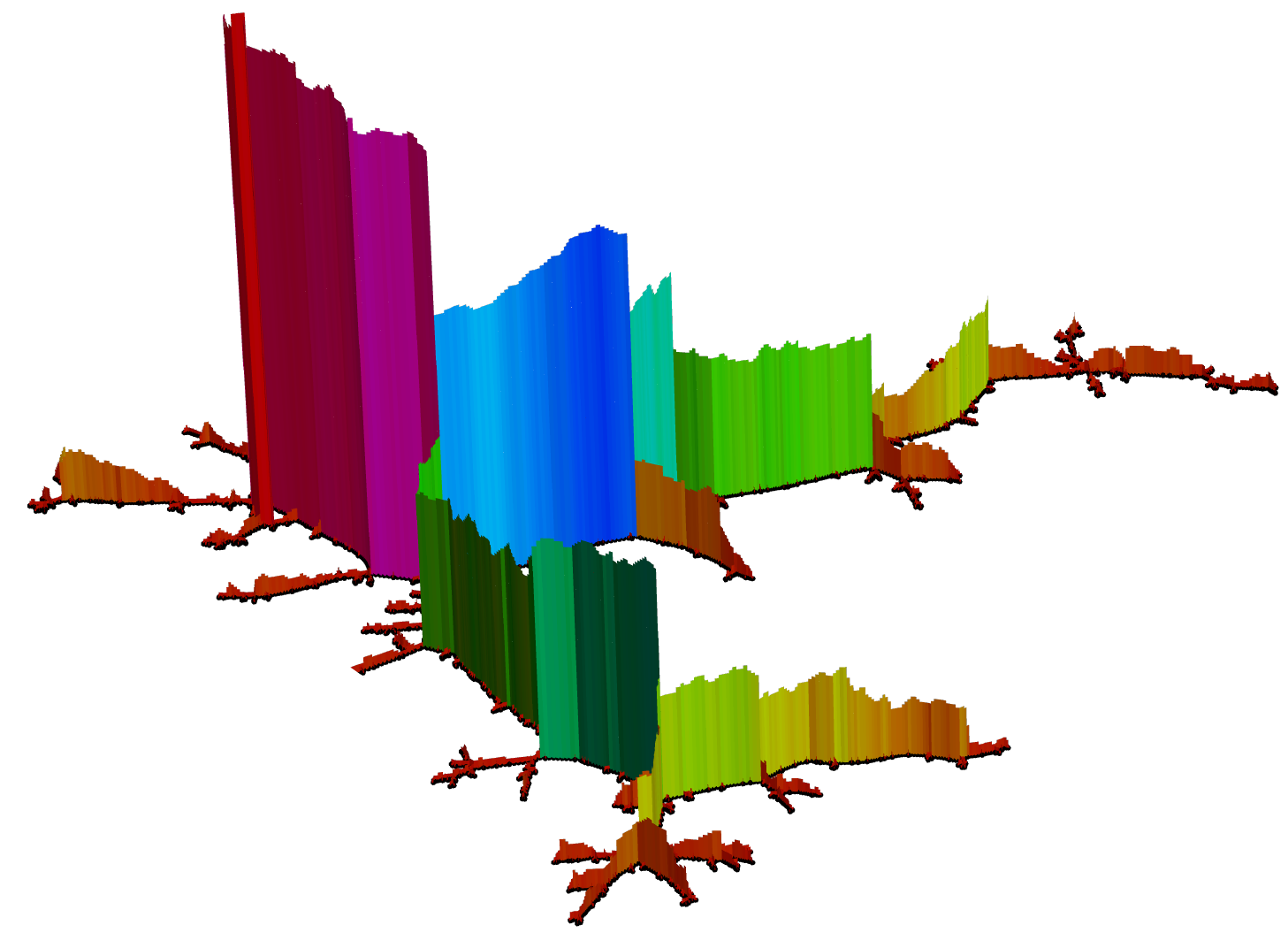}
 \caption{A simulation of the Brownian growth-fragmentation tree. The process is binary and conservative: at each splitting event, the total mass is conserved and split between two children.}
 \end{center}
 \end{figure}

The Brownian growth-fragmentation is self-similar with exponent $\alpha = 1/2$, and its 
 cumulant function is given by 
$$ \kappa_{\mathrm{BroGF}}(\gamma) = \frac{\Gamma(\gamma- 3/2)}{\Gamma(\gamma-3)}, \quad  \mbox{for }  \gamma > 3/2.$$
Its Gaussian coefficient is $\sigma^{2}=0$, its generalized L\'evy measure $ \boldsymbol{\Lambda}_{\mathrm{BroGF}}$ is binary conservative and given by 
  \begin{eqnarray} \label{eq:GLMgrowthBrownian} \int_{\mathcal S} F\big(  \mathrm{e}^{{y_{0}}}, (\mathrm{e}^{y_{1}}, ...)\big)  \  \boldsymbol{ \Lambda}_{\mathrm{BroGF}}(\mathrm{d} y_{0}, \mathrm{d} \mathbf{y})  \quad  := \quad    \frac{3}{ 4 \sqrt{\pi}} \int_{1/2}^{1} F\big(x,(1-x,0,0, ...)\big) \frac{ \mathrm{d}x}{(x(1-x))^{5/2}}.  \end{eqnarray}
Since $ \boldsymbol{\Lambda}_{\mathrm{BroGF}}(\dd y_0, \dd \mathbf{y}) $ does not integrate $1 \wedge |y_0|$, the drift coefficient is not canonic; for the cut-off function as in \eqref{E:LKfor} we have 
$$ \mathrm a_{\mathrm{BroGF}} = -\frac{2}{\sqrt{\pi}} + \frac{3}{4 \sqrt{\pi}} \int_{1/2}^{1}   \frac{ \log x + 1 -x}{(x(1-x))^{5/2}} \, \dd x= \frac{4(7-3\pi)}{3\sqrt{\pi}} ,$$
see \cite[Eq (32)]{BCK18} or \cite[Proof of Proposition 5.2]{BBCK18} for a different drift using a different cut-off function. 
 
 In particular, we have $\omega_{-} = 2$ and $\omega_+=3$, and Assumption \ref{A:omega-} holds. The distribution of  the total mass $\upmu(T)$ is known from \cite[Corollary 6.7]{BBCK18}, it satisfies the striking property that its size-biased transform is a $1/2$-stable law. Namely, under $\P_1$, the total mass $\upmu(T)$ has density:
 \begin{equation}\label{density:bro:gro}
 \frac{1}{\sqrt{2\pi} x^{\frac{5}{2}}}\exp\big(-\frac{1}{2x}\big)\mathbf{1}_{(0,\infty)}(x).
 \end{equation}
 Let us develop in more details the connection with Brownian geometry since this is one of central motivation for this whole work. The Brownian sphere is a random compact metric space almost surely homeomorphic to the $2$-sphere, and which has famously proved to be the scaling limit of various classes of random planar maps equipped with the graph distance \cite{Mie11,LG11}, or shown to be the random metric space obtained by exponentiating a planar Gaussian Free Field with the proper parameter $ \sqrt{8/3}$, see \cite{miller2015liouville,miller2021liouville,gwynne2020existence}. The Brownian sphere has a variant, having the topology of the disk, called the \textbf{Brownian disk}, see  \cite{BM15}. In particular,  the boundary of the Brownian disk may be defined as the set of all points that have no neighborhood homeomorphic to the open unit disk. Let us denote a \textbf{free Brownian disk} with boundary size $x$ (it has a random volume) by $ \mathbb{D}_x$  and  its boundary by  $\partial \mathbb{D}_x$. For every point $u \in \mathbb{D}_x$ we write $H(u)$ for the height of $u$, defined as its distance to $\partial \mathbb{D}_x$, and 
 \begin{figure}[!h]
  \begin{center}
  \includegraphics[width=14cm]{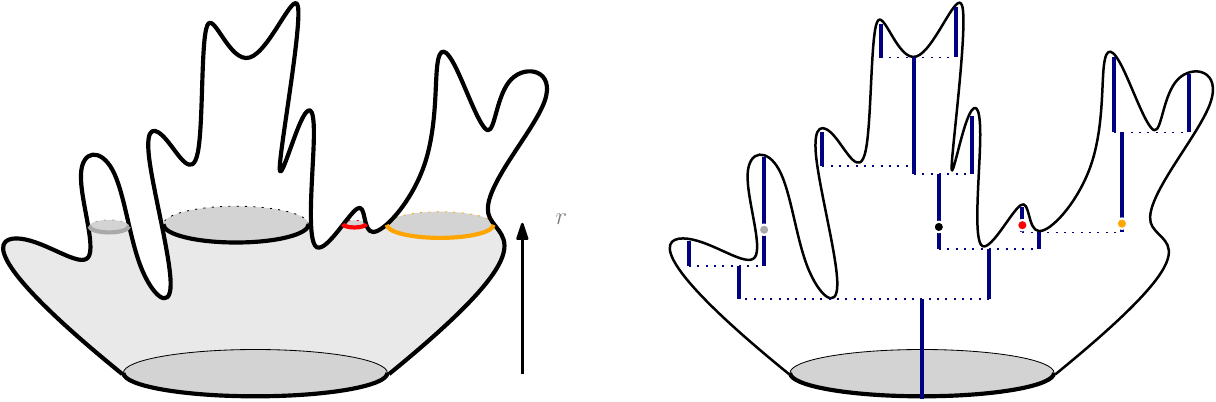}
  \caption{Illustration of the Cactus tree of a Brownian disk. The ball of radius $r$ (measured from $ \partial \mathbb{D}_{x}$) is depicted in light gray and it has several boundary components. Each of these components has a ``size'' which enables us to decorated the cactus tree (on the right).}
  \end{center}
  \end{figure}
for every $r>0$,  we consider the ``ball'' $ \mathrm{B}_{r} = \{ u \in \mathbb{D}_{x} : H(x) \leq r\}$. Its boundary $ \partial \mathrm{B}_{r}$ is  made of several components homeomorphic to circles. Each of these boundary components $ \mathcal{C}$ is a fractal curve (of dimension $2$) but it is possible to give a meaning to its size  $| \mathcal{C}|$ by approximation (a.s. simultaneously for every $r \geq 0$ and for each boundary component), see \cite[Theorem 3]{le2020growth}. In particular, the boundary $\partial \mathbb{D}_{x}$ has size $x$ a.s. Furthermore, those boundary components have a natural tree genealogy as $r$ varies. Specifically, consider the pseudo-metric $ d_{ \mathrm{Cac}}$ defined by 
  $$ d_{ \mathrm{Cac}}(u,v) = H(u) + H(v)- 2 \sup_{ \gamma : u \to v} \left( \min_{0 \leq t \leq 1} H( \gamma(t))\right),$$
  where the supremum is over all continuous curves $ \gamma :  [0,1] \to \mathbb{D}_{x}$ such that $\gamma(0) =u$ and $ \gamma(1)=v$. Then, we introduce  $ \mathrm{Cac}( \mathbb{D}_{x})$ the quotient space of $\mathbb{D}_{x}$ for the equivalence relation defined by setting $a\sim_{\mathrm{Cac}} b$ if and only if $d_{ \mathrm{Cac}}(a,b) =0$. The quotient  $ \mathrm{Cac}( \mathbb{D}_{x})$, equipped with the metric induced by $d_{ \mathrm{Cac}}$,  is a compact real tree called \textbf{the cactus} of $ \mathbb{D}_{x}$ seen from $ \partial \mathbb{D}_{x}$. We root  $ \mathrm{Cac}( \mathbb{D}_{x})$ at $\rho_{ \partial \mathbb{D}_{x}}$ the equivalence class of $ \partial \mathbb{D}_{x}$, see \cite[Section 2.2]{CLGMcactus}. It is easy to see that the equivalence classes for $\sim_{\mathrm{Cac}}$ are precisely the boundary cycles of $ \partial \mathrm{B}_{r}$ for all $r \geq 0$, so that the size function $ | \mathcal{C}|$ of boundary components is a well defined decoration  rcll on branches on the cactus tree and we denote  its usc modification by $g_{|\cdot|}$. Then, it follows from \cite{le2020growth} (see in particular Theorem 3 there) that the law of the equivalence class in $ \mathbb{T}$ of the decorated tree  
  $$ \big( \mathrm{Cac}( \mathbb{D}_{x}), d_{ \mathrm{Cac}}, \rho_{ \partial  \mathbb{D}_{x}}, g_{|\cdot|} \big) $$
  is that of $ \mathtt{T} =  (T, d_T, \rho, g)$ under $ \mathbb{P}_{x}$ for the  characteristic quadruplet $( 0,   \frac{2 \sqrt{2}}{ \sqrt{3}} \cdot \mathrm{a}_{\mathrm{BroGF}}, \frac{2 \sqrt{2}}{ \sqrt{3}} \cdot  \boldsymbol{ \Lambda}_{\mathrm{BroGF}} ;  \frac{1}{2})$. Actually, \cite{le2020growth} only deals with the growth-fragmentation point of view, but the results there are established using cell processes so that the previous display is a consequence of the arguments therein. Let us sketch the reasoning: It is proved in \cite[Proposition 16]{le2020growth} that the self-similar process obtained by following the locally largest exploration is precisely the pssMp $X$ constructed from $( 0,   \frac{2 \sqrt{2}}{ \sqrt{3}} \cdot \mathrm{a}_{\mathrm{BroGF}}, \frac{2 \sqrt{2}}{ \sqrt{3}} \cdot  \boldsymbol{ \Lambda}_{\mathrm{BroGF}} ;  \frac{1}{2})$ in Section \ref{sec:2.2}. In particular, in the growth-fragmentation case (binary conservative case), the decoration-reproduction $\eta$ is recovered from $X$. Finally, conditionally on this exploration, by \cite[Proposition 18]{le2020growth}, the decorated subtrees branching are conditionally independent given their initial  decorations. Using these two ingredients, one can couple the construction of Section \ref{sec:2.1} with the iterative locally largest exploration of a Brownian disk so that they coincide. We stress that the study of \cite{le2020growth} is possible thanks to the encoding of the Brownian disk through the Brownian motion indexed by the Brownian tree and the excursion theory developed  in \cite{ALG15}. 
  \hfill $\diamond$
  \end{example}
  
 Related to the above example, one can consider the self-similar Markov tree with the same first three characteristics $( \sigma^{2},  \mathrm a_{\mathrm{BroGF}}, \boldsymbol{ \Lambda}_{\mathrm{BroGF}})$ but with self-similarity parameter $3/2$ (instead of $1/2$). This modification has the effect of performing a length change along branches \textit{\`a la} Lamperti.
This tree actually appears as the scaling of the so-called \textbf{peeling} trees associated to random planar maps with small faces, see Theorem \ref{prop:GF} in Part II, and \cite[Section 6]{BBCK18} for previous work. We expect it to be further the scaling limit of several other discrete models such as critical fully parked trees \cite{ConCurParking} (this is proved in the bounded degree case in \cite{ConCurUni}, see Theorem \ref{thm:ConCur25} below), fighting-fish \cite{duchi2017fightingbis} or even the peeling trees of plane Weil-Petersson surfaces with $n$ punctures. See  Chapter \ref{chap:applications} for more details. 
The ssMt with characteristics $( \sigma^{2},  \mathrm a_{\mathrm{BroGF}}, \boldsymbol{ \Lambda}_{\mathrm{BroGF}}; 3/2)$ also appeared very recently as the genealogical tree of $2\pi/3$-cone excursions inside Brownian motion in a cone, see \cite{da2025growth}. Actually, \cite[Theorem 1.1]{da2025growth} is stated in terms of growth-fragmentation process but as in the preceding example, the ssMt can be extracted from their proof. We refer to Example \ref{ex:brownian} where another ssMt is constructed from excursions of a Brownian motion in a half-plane. \bigskip

At bit more generally, one can allow the trajectories of cells to have positive jumps, keeping up with the convention that only negative jumps are interpreted as mitosis events. Then \cite[Theorem 5.1]{BBCK18} presents the following generalization of the preceding example:

\begin{example}[The $(\beta, \varrho)$-stable family] \label{ex:stablefamily}There exists a family parametrized by $( \texttt{a}\in (0,1], \texttt{b} \in (0,1/2], \alpha\in(0,\infty))$ of  self-similar Markov trees with  cumulant functions given by 
  \begin{eqnarray}\kappa_{\texttt{a,b} }(\gamma)= - \frac{\Gamma(1 + 2 \texttt{a} + 2 \texttt{b} - \gamma) \Gamma(\gamma-\texttt{a} - \texttt{b} )}{ \Gamma(1 + \texttt{a} + 2 \texttt{b} - \gamma) \Gamma(\gamma-\texttt{a} - 2 \texttt{b})}, \quad \gamma \in (\texttt{a}+\texttt{b},2\texttt{a}+2\texttt{b}+1),   \label{eq:kappastable}\end{eqnarray}
  and where $\alpha$ is simply the self-similarity exponent.
   Those expressions may seem ad-hoc for the moment, but we will see in Section \ref{sec:examplespinal} that they appear in relation to so-called Lamperti-stable processes,  \cite[Section 4.3]{kyprianou2022stable},
and  are naturally associated with $\beta$-stable L\'evy process with positivity parameter $\varrho$ for  
 \begin{eqnarray} \label{eq:famy} \beta = \texttt{a}+ \texttt{b}  \quad \mbox{ and } \quad \varrho = \frac{\texttt{a}}{\texttt{a}+\texttt{b}}, \end{eqnarray} see Section \ref{sec:spinalex} for details. Specifically, for  $( \texttt{a}\in (0,1], \texttt{b} \in (0,1/2], \alpha\in(0,\infty))$, the associated  self-similar Markov tree can be obtained from a characteristic quadruplet $(\sigma^2_{ \texttt{a},\texttt{b}}, a_{ \texttt{a},\texttt{b}}, \boldsymbol{\Lambda}_{ \texttt{a},\texttt{b}}; \alpha)$ defined as follows. The Gaussian part is degenerate, $\sigma^{2}_{ \texttt{a},\texttt{b}}=0$, and
the generalized L\'evy measure $ \boldsymbol{\Lambda}_{ \texttt{a},\texttt{b}}$ is prescribed by three parts
$$ \int_{\mathcal S} F\big( \mathrm{e}^{y_0}, (\mathrm{e}^{y_{1}}, \mathrm{e}^{y_{2}}, \cdots)\big) \ \boldsymbol{\Lambda}_{ \texttt{a},\texttt{b}}( \dd y_0, \dd \mathbf{y}) \quad =$$
  \begin{align}&\frac{\Gamma(\beta+1) \sin(\pi \texttt{b})}{\pi} \int_{1/2}^{1}  \frac{ \mathrm{d}x}{(x(1-x))^{\beta+1}} F\big(x,(1-x, 0, \ldots)\big) & ( \mbox{conservative binary splitting})\nonumber \\
   + \quad &\frac{\Gamma(\beta+1) \sin(\pi \texttt{a})}{\pi}  \int_0^\infty \frac{ \mathrm{d}x}{(x(1+x))^{\beta+1}} F\big(1+x,(0, 0, \ldots)\big) & (\mbox{unique positive jump, growth})\nonumber \\
   + \quad  &\cos( \pi \texttt{b}) \frac{2 \Gamma(2\beta)}{\Gamma(\beta)} \cdot F\big(0,(0, 0, \ldots)\big) & (  \mbox{killing}),   \label{eq:stablefam}
   \end{align}
 for a generic positive function $F : \mathcal{S} \to \mathbb{R}_{+}$. The drift coefficient $\mathrm{a}_{ \texttt{a},\texttt{b}}$ depends on the cut-off $y\mathbf{1}_{|y|<1}$ and is  given by:
 \begin{equation}\label{drift:a:b}
  \mathrm{a}_{ \texttt{a},\texttt{b}}=\kappa^\prime_{\texttt{a,b} }(0)- \frac{\Gamma(\beta+1)}{\pi}\Big(  \sin(\pi \texttt{a})\cdot  \frac{\partial}{\partial r}\widetilde{\mathrm{B}}_{\frac{1}{2}}(-\beta, r)\Big|_{r=-\beta}+\sin(\pi \texttt{b})\cdot \frac{\partial}{\partial r}\mathrm{B}_{e^{-1}}(r,-\beta)\Big|_{r=2\beta+1} \Big),
  \end{equation}
 where $\mathrm{B}_{z}$ stands for the incomplete Beta function and $\widetilde{\mathrm{B}}_{z}=\mathrm{B}_{1}-\mathrm{B}_{z}$. A straightforward computation gives that the cumulant function of  $(\sigma^2_{ \texttt{a},\texttt{b}}, a_{ \texttt{a},\texttt{b}}, \boldsymbol{\Lambda}_{ \texttt{a},\texttt{b}}; \alpha_{ \texttt{a},\texttt{b}})$ is given by \eqref{eq:kappastable}. In particular, they satisfy $\omega_{-} = \texttt{a}+2\texttt{b}, \omega_{+} = \texttt{a}+2\texttt{b}+1$ and Assumption \ref{A:omega-} holds. To lighten notation, we write $\mathtt{k}_{ \texttt{a},\texttt{b}}=\cos( \pi \texttt{b}) \frac{2 \Gamma(2\beta)}{\Gamma(\beta)} $ for the killing rate.  We also stress that, in most of applications,  the self-similarity parameter is simply $\alpha= \beta$, see Chapter \ref{chap:applications}.

 Let us give a more explicit description of the characteristics in two cases (in {\color{red}red} and {\color{blue}blue} on Figure \ref{fig:famystable}):
 \paragraph{{\color{red}No killing}.} 
 We take $\texttt{b}=1/2$ in \eqref{eq:famy} so that  $\mathtt{k}_{ \texttt{a},\texttt{b}} = 0$, $ 1/2 < \beta = \texttt{a}+\texttt{b}\leq 3/2$ and $\beta(1-\varrho)= 1/2$. For definiteness, we take the exponent  $\alpha = \beta$. This gives the self-similar Markov tree with exponent,  no Gaussian component,  generalized L\'evy measure given by
 \begin{align} \label{eq:nokiling}  & \frac{\pi}{\Gamma(\beta+1)}  \int_{\mathcal S} F\big( \mathrm{e}^{y_0}, (\mathrm{e}^{y_{1}}, \mathrm{e}^{y_{2}}, \cdots)\big) \ \boldsymbol{\Lambda}_{\texttt{a},\frac{1}{2} }( \dd y_0, \dd \mathbf{y})\nonumber \\
 =& \int_{1/2}^{1}  \frac{ \mathrm{d}x}{(x(1-x))^{\beta+1}} F\big(x,(1-x, 0, \ldots)\big) \nonumber \\ 
 &+  \cos( (\beta+1) \pi) \cdot \int_0^\infty \frac{ \mathrm{d}x}{(x(1+x))^{\beta+1}} F\big(1+x,(0, 0, \ldots)\big),  \end{align} 
 and where the drift coefficient $\mathrm a$ is prescribed so that the cumulant function equals 
  \begin{eqnarray} \label{eq:kappanokilling} \kappa_{\texttt{a},\frac{1}{2} }(\gamma) = \frac{\cos(\pi ( \gamma - \beta))}{\sin(\pi (\gamma - 2 \beta))} \frac{\Gamma(\gamma - \beta)}{\Gamma(\gamma - 2 \beta)}, \quad \mbox{ for } \beta < \gamma < 2 \beta +1,  \end{eqnarray} in particular we have $\omega_{-} = \beta + \frac{1}{2}$ and $\omega_{+} = \beta + \frac{3}{2}$, see \cite[Section 5]{BBCK18}. When $\beta \in (1/2,1)$ the L\'evy measure $\Lambda_{0}$ integrates $x$ and the canonic drift coefficient is $ \mathrm{a_{can}}=0$. Assumption \ref{A:omega-} holds and the limiting mass measure has the law of a positive $1/(\beta+1/2)$-stable random variable biased by $x \mapsto 1/x$, see \cite[Proposition 4]{BC16}. Specifically, we have:

 $$\E_1\Big[\mu(T)\cdot \exp(-\lambda \mu(T))\Big]=\exp\Big(- \Big(\Gamma\big(\beta+\frac{3}{2}\big)  \lambda\Big)^{\frac{2}{2+\beta}}\Big)$$
and equivalently the density of $\mu(T)$, under $\P_1$, can be written in terms of special functions as:
 \begin{equation}\label{density:can:no:killing}
\frac{2 x^{-1}}{(2\beta+1)\Gamma(\beta+\frac{3}{2})}\cdot \Big(\frac{\Gamma\bigl(\beta+\frac{3}{2}\bigr)}{x}\Big)^{\frac{2\beta+3}{2\beta+1}}\cdot \, M_{\frac{2}{2\beta+1}}\Big(\Big(\frac{\Gamma\bigl(\beta+\frac{3}{2}\bigr)}{x}\Big)^{\frac{2}{2\beta+1}}\Big) \mathbf{1}_{x\in(0,\infty)},
\end{equation}
 where $M_{\nu}(z):=\sum_{n\geq 0}\frac{(-z)^n}{n! \Gamma(-\nu n +1-\nu)} $ is the $M$-Wright function. We can add to this family the limiting point $\texttt{b}=1/2, \texttt{a}=0$ corresponding to Example \ref{ex:brownian}. Those decorated trees appear as the scaling limit of the peeling trees in critical discrete stable planar maps, see Theorem \ref{prop:GF} in Chapter \ref{chap:applications} and \cite[Section 6]{BBCK18} for a previous work. They also appear  in Liouville Quantum Gravity: the results of  \cite{miller2022simple,miller2021non} says heuristically that the branching structure underlying the Markovian exploration of a $ \mathrm{CLE}_{\kappa}$, and where the labels are given by the $\sqrt{\kappa}$-LQG boundary length (with an independent GFF)  is the ssMt above with the relation $\beta = \frac{4}{ \kappa}.$

\begin{figure}[!h]
 \begin{center}
 \includegraphics[width=10cm]{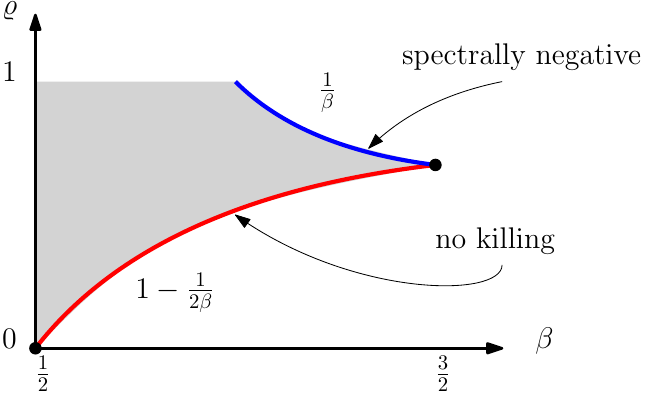}
 \caption{The two-parameters family of "stable" binary conservative self-similar Markov trees. See Section \ref{sec:examplespinal} for the explanation of the relation with $\beta$-stable processes with positivity parameter $\varrho$. The upper blue  boundary case corresponds to stable L\'evy process with no positive jumps (but the generalized L\'evy measure $ \boldsymbol{ \Lambda}$ has a positive killing). This case is encountered in the scaling limit of fully parked trees, see Theorem \ref{thm:parkingstable} in Chapter \ref{chap:applications}. The lower red boundary corresponds to the case with no killing. This case is encountered in critical stable random planar maps, see Theorem \ref{prop:GF} in Chapter \ref{chap:applications}. The case $\beta=3/2, \varrho= \frac{2}{3}$ where those two lines merge is the spectrally negative $3/2$-stable case discussed in Example \ref{ex:3/2stable}, but with self-similarity parameter $\alpha = \frac{3}{2}$. The case $\beta=1/2$ and $\varrho=0$ corresponds to the pure fragmentation of the Brownian CRT, see Example \ref{ex:brownian}.}
 \label{fig:famystable}
 \end{center}
 \end{figure}
 
\paragraph{{\color{blue} Conservative.}} We take $\texttt{a}=1$  in  \eqref{eq:famy} so that  $1< \beta = \texttt{a}+\texttt{b} \leq 3/2$ and $\varrho = \frac{1}{\beta}$. This gives the self-similar Markov tree with index $\alpha = \beta$, Gaussian coefficient  $\sigma^2=0$, and binary generalized L\'evy measure given by 
\begin{align} \label{eq:GLMstablespecneg}
& \int_{\mathcal S} F\big( \mathrm{e}^{y_0}, (\mathrm{e}^{y_{1}}, \mathrm{e}^{y_{2}}, \cdots)\big)  \boldsymbol{\Lambda}_{1,\texttt{b}}( \dd y_0, \dd \mathbf{y})\nonumber \\
&=  \frac{\Gamma(\beta+1) \sin (\pi (\beta-1))}{\pi}\left(\int_{1/2}^{1}  \frac{ \mathrm{d}x}{(x(1-x))^{\beta+1}} F\big(x,(1-x,0, \ldots)\big)\right)\nonumber\\
&+  \frac{2 \Gamma(2\beta) \cos( \pi (\beta-1))}{\Gamma(\beta)} \cdot F\big(0,(0, 0, \ldots)\big),
\end{align} and where the drift $ a \in \mathbb{R}$ is prescribed so that the cumulant function  is equal to 
$$\kappa_{1,\texttt{b}}(\gamma) = -\frac{\Gamma(1 + 2 \beta - \gamma) \Gamma(\gamma-\beta) \sin(\pi (2 \beta - \gamma))}{\pi}.  $$
Assumption  \ref{A:omega-} holds with $\omega_{-} = 2 \beta-1$. We show in Theorem \ref{thm:parkingstable} that these self-similar Markov trees (ssMt) arise as scaling limits of fully parked components in critical parking processes with stable car arrival distributions, as considered in \cite{chen2021enumeration}. In particular, we expect that the law of the total harmonic mass under $\P_1$ is given by the density
 \begin{equation}\label{density:can:no:positive:jump}
\sum \limits_{n\geq 0} \frac{c_{n}}{\Gamma(b_n-\beta)\Gamma(-\frac{b_n}{2\beta-1})} x^{-\frac{b_n}{2\beta-1}-1},\quad x>0, 
\end{equation}
where $\sum_{n\geq 0} c_n x^{b_n}=\sqrt{1-2\beta x^{2\beta-1}+(2\beta-1) x^{2\beta}}.$ 
The latter is strongly supported by the enumerative results of \cite{chen2021enumeration}, and in particular by the asymptotic expansion announced in \cite[Eq. (12–13)]{chen2021enumeration}, which gives exactly \eqref{density:can:no:positive:jump}. 

This family of self-similar Markov trees is also expected to arise naturally in continuous canonical models. In particular, they should appear in the study of excursions above the minimum of the stable Brownian snake (or Brownian motion indexed by the stable tree), as well as in related models of random geometry \cite{archer2024stable}. This connection is rigorously established in the ongoing work \cite{RRO2025}, where the density \eqref{density:can:no:positive:jump} arises naturally through the analysis of the stable Brownian snake, certain differential equations, and the excursion theory developed in \cite{RRO2024}. In particular, it is shown there that \eqref{density:can:no:positive:jump} indeed gives the density of the total harmonic mass under~$\P_1$.

We close this section by referring to  Section \ref{sec:examplespinal} for a related family of Examples with generalized L\'evy measures similar to \eqref{eq:GLMstablespecneg} but where $\beta \in (0,1/2]$.
\hfill $\diamond$
\end{example}

\section{An overlay on the stable family and a critical example}
Our final example is a family of binary non-conservative self-similar Markov trees, so that a trimmed version of which gives the family with no killing in Example \ref{ex:stablefamily}.  Those ssMt should appear in connection with $O(n)$-decorated random planar maps, see Part II for details. There is a critical case in this family  which is naturally associated with the Brownian CRT (Example \ref{ex:brownian}) in a rather surprising way.

\begin{example}[An overlay on the stable family] \label{ex:overlay} Recall the case $ \texttt{b}= 1/2$ in Example \ref{ex:stablefamily}. In particular, with the notation used there we have $ \alpha = \beta = \texttt{a}+\texttt{b}$. We consider now the ``augmented'' self-similar Markov tree obtained by adding; in case of positive jumps $x \mapsto x+y$, a new particle of mass $y$. Formally this is done by replacing the generalized L\'evy measure defined in \eqref{eq:nokiling} by $\boldsymbol{\Lambda}_{ \texttt{a}, \frac{1}{2}}^+$, where

\begin{align} \label{eq:nokilingbis}   
 & \frac{\pi}{\Gamma(\beta+1)}  \int_{\mathcal S} F\big( \mathrm{e}^{y_0}, (\mathrm{e}^{y_{1}}, \mathrm{e}^{y_{2}}, \cdots)\big) \ \boldsymbol{\Lambda}_{ \texttt{a}, \frac{1}{2}}^+( \dd y_0, \dd \mathbf{y})\nonumber \\
 &= \int_{1/2}^{1}  \frac{ \mathrm{d}x}{(x(1-x))^{\beta+1}} F\big(x,(1-x, 0, \ldots)\big) \nonumber \\ 
 &+  \cos( (\beta+1) \pi) \cdot \int_0^\infty \frac{ \mathrm{d}x}{(x(1+x))^{\beta+1}} F\big(1+x,(x, 0, \ldots)\big),  \end{align} 
  (notice the small difference with \eqref{eq:nokiling}: we replaced $F(1+x,(0, ...))$ there  by $ F(1+x,(x, 0, ...))$). The cumulant function is easily updated and becomes
 \begin{eqnarray*}&&\underbrace{\kappa(\gamma)}_{ \mathrm{in \ } \eqref{eq:kappanokilling}} +  \frac{\Gamma(\beta+1)}{\pi} \cos ((\beta-1)\pi) \int_{0}^{\infty}(x-1)^{\gamma} \frac{\mathrm{d}x }{(x(x-1))^{{\beta+1}}}\\  &=& - \frac{	\Gamma(\gamma - \beta) \mathrm{Sec}( \frac{\pi}{2} (\gamma - 2 \beta))\sin(\pi \gamma /2)}{\Gamma(\gamma- 2 \beta)}, \quad \mbox{ for } \beta < \gamma < 2 \beta+1,  \end{eqnarray*}
for which now the two roots are $\{ 2 \beta, 2 \}$ so that $\omega_{-} = \min \{ 2 \beta, 2 \}$  and $\omega_{+} = \max \{ 2 \beta, 2 \}$. Assumption \ref{A:gamma0} holds and so the self-similar Markov tree $ \mathtt{T}=(T,d_T,\rho, g)$ exists as soon as $\beta \ne 1$. Under this condition, Assumption \ref{A:omega-} also holds and the law of harmonic mass $\upmu(T)$ under $ \mathbb{P}_{1}$ has been identified in \cite[Theorem 4.3]{CCMOn} and related to an inverse-Gamma distribution.
Watson \cite{watson2023growth} extended the model where the exploration of positive jumps appears with probability $r \in [0,1]$ and found explicit cumulant functions as well.
\hfill $\diamond$
\end{example}

In the example above, when $\beta=1$ the cumulant function has a double root and stays non negative so that we cannot directly apply Proposition \ref{P:constructionomega-} to construct a self-similar Markov tree. However, the work of Aidekon \& Da Silva \cite{aidekon2022growth} identified directly the underlying ssMt as a variant of the Brownian CRT of Example \ref{ex:brownian}.

\begin{example}[Half-plane excursions after Aidekon \& Da Silva \cite{aidekon2022growth}] \label{ex:ADS} We consider a two-dimensional excursion $(X,\e)$  where $\e$ is a Brownian excursion of length $\ell$ and $X$ is a Brownian bridge of length $\ell$ going from $0$ to $1$ under the normalized excursion measure $ \mathbb{N}^{[1]}_{ \mathrm{ads}}$ for plane Brownian motion. This law can be expressed using It\^{o}'s positive Brownian excursion measure \eqref{def:itobro} as 
$$ \mathbb{N}^{[1]}_{\mathrm{ads}}\big( \mathrm{d} (X,\e)\big)= \int_{ \mathbb{R}_+} \mathrm{d}\ell ~\frac{ \mathrm{e}^{- \frac{1}{2\ell}}}{ 2 \ell^2} ~ P^{0\to 1}_\ell( \mathrm{d}X)\otimes \mathbb{N}_\ell ( \mathrm{d}\e),$$ where $ {P}_\ell^{0\to1}$ is the law of the Brownian bridge of length $\ell$ going from $0$ to $1$; see \cite[Proposition 2.9]{aidekon2022growth}. Given the pair $(X,\e)$ distributed according to  $\mathbb{N}^{[1]}_{\mathrm{ads}}$, we can first construct the rooted tree $ (\mathcal{T}_\e, d_{ \mathcal{T}_\e}, \rho)$ coded by the excursion $ \e : [0, \ell] \to \mathbb{R}_+$ as presented in Section \ref{sec:codagearbre} using \cite{DLG05}. The random real tree $ \mathcal{T}_{\e}$  is nothing but a mixture of Brownian CRT (Example \ref{ex:brownian}) whose size is distributed according to the $1/2$-stable law $\ell$. We then endow it with the following decoration using the process $X$: for any point $u \in \mathcal{T}_\e$, in the coding of $ \mathcal{T}_{\e}$ from $\e$, let us denote by $s_u,t_u$ respectively the minimal and maximal pre-images\footnote{Recall that almost surely, any point $u \in \mathcal{T}_{\e}$ has at most three pre-images by the projection $\pi : [0,\ell] \to \mathcal{T}_e$. The points having three pre-images correspond to the branch points of $ \mathcal{T}_{\e}$ and those having at least two pre-images correspond to the skeleton of $ \mathcal{T}_{\e}$.} of $u$ in $[0, \ell]$. The interval $[s_u,t_u]$ is then a \textit{subexcursion} of $ \e$ in the sense that  $ \e(r) \geq  \e(s_u) = \e(t_u)$ for all $r \in (s_u,t_u)$ and $t_u-s_u$ correspond to the size of the fringe subtree above $u$. The point $u \in \mathcal{T}_{\e}$ is  labeled by the $X$-displacement over the time interval $[s_u,t_u]$, i.e.
 $$\tilde{g}(u) := |X(t_u)-X(s_u)|.$$ \begin{figure}[!h]
  \begin{center}
  \includegraphics[width=13cm]{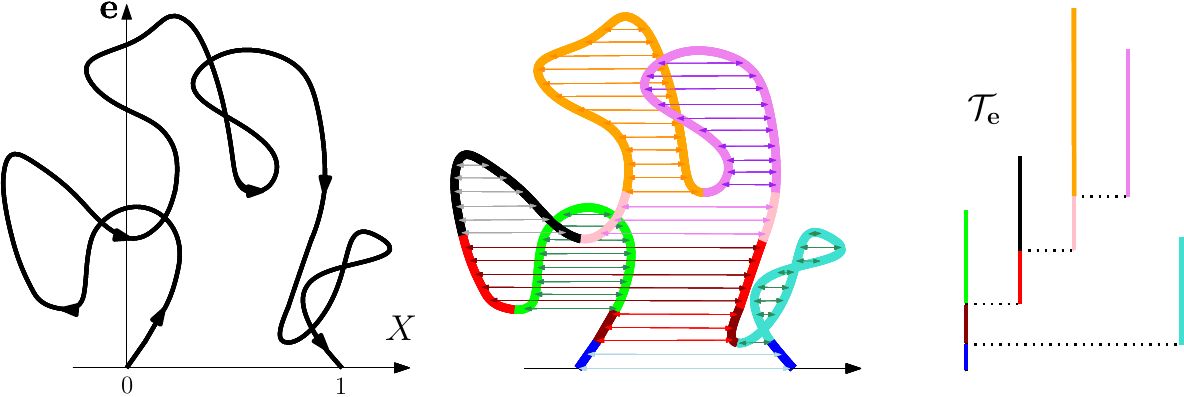}
  \caption{Illustration of the construction of a decorated tree from an excursion in the half-plane. The vertical coordinate encodes the tree structure, whereas the horizontal displacement encodes the decoration.}
  \end{center}
  \end{figure}
  
 It is not hard to see that $\tilde{g}$ is then a rcll function over the branches of $ \mathcal{T}_{\e}$ and  we thus consider the usc modification $g$ of $\tilde{g}$ as the decoration on $ \mathcal{T}_{\e}$ to fit the framework of Chapter \ref{chap:topology}. Remark in particular that we have $g(\rho) =1$ under the law $\mathbb{N}^{[1]}_{\mathrm{ads}}$. It follows from the arguments of \cite{aidekon2022growth}  that the law of the random decorated tree 
 $$ \big(\mathcal{T}_\e, d_{ \mathcal{T}_\e}, \rho, g \big) \quad \mbox{ under }\mathbb{N}^{[1]}_{\mathrm{ads}},$$
is that of the ssMt $ \mathtt{T}$ under $ \mathbb{P}_1$ with the characteristic quadruplet $(0, \mathrm{a}_{ \mathrm{ads}}, \boldsymbol{\Lambda}_{ \mathrm{ads}} ; 1)$ where 
  \begin{eqnarray*}&&\int F\big( \mathrm{e}^{y_{0}}, (\mathrm{e}^{y_{1}}, ... )\big) \boldsymbol{\Lambda }_{ \mathrm{ads}}( \mathrm{d}y_{0}, \mathrm{d}(y_{i})_{i\geq1})\\ &=&   \frac{2}{\pi} \left(\int_{1/2}^{1} \frac{\mathrm{d}x}{\big(x(1-x)\big)^{2}} F\big(x,(1-x,0,...)\big) + \int_{1}^{\infty} \frac{\mathrm{d}x}{\big(x(1+x)\big)^{2}} F\big(1+x,(x,0, ...)\big)\right),  \end{eqnarray*}
 $$ \mbox{ and } \quad  \mathrm{a}_{\mathrm{ads}} = - \frac{4}{\pi} + \frac{2}{\pi} \int_{- \log 2}^{\infty}  \mathrm{d}y \, ( \mathbf{1}_{|y| \leq 1} - ( \mathrm{e}^{y}-1)) \frac{ \mathrm{e}^{-y}}{( \mathrm{e}^{y}-1 )^{2}}. $$ 
\hfill $\diamond$
\end{example}
 \begin{figure}[!h]
 \begin{center}
 \includegraphics[width=14cm]{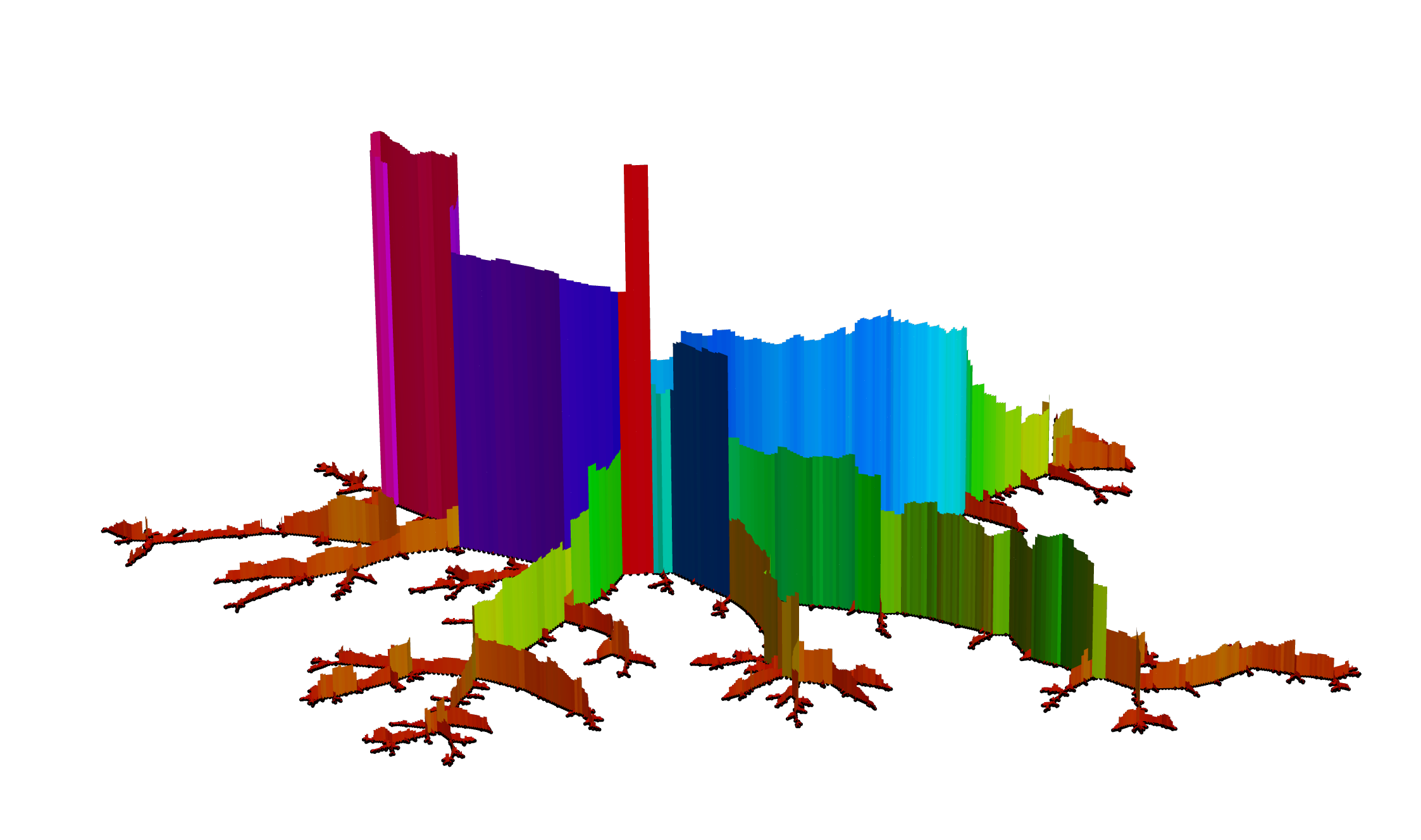}
 \caption{The decorated random tree associated with half-planar Brownian excursion. }
 \end{center} 
 \end{figure}
 That is, the decorated L\'evy measure is twice that of \eqref{eq:nokilingbis} and the drift coefficient is set so that we have 
$$ \kappa(\gamma) = 2(\gamma-2)\tan\left( \frac{\pi \gamma}{2}\right).$$
Some explanations are in order. First, in the work \cite{aidekon2022growth}, the authors only consider  ``positive excursions", i.e. they trimmed the tree at points $u \in \mathcal{T}_\e$ where $X(t_u)-X(s_u)$ becomes negative. If we do so, we instead recover the ssMt corresponding to \eqref{eq:nokiling} with $\beta =2$. However, we prefer to keep positive and negative horizontal excursions (their proof adapts easily) so that the underlying tree is exactly $ \mathcal{T}_{\e}$. Second, as in \cite{le2020growth}, the authors state their results in terms of growth-fragmentation processes, but we can argue similarly as in Example \ref{ex:3/2stable} using \cite[Theorem 3.3 and Theorem 3.6]{aidekon2022growth} and prove that we can couple the construction of Chapter \ref{chap:generalBP} with that of \cite{aidekon2022growth} so that the underlying branching process with types and decoration-reproduction is actually the same. We \textit{deduce} that in this case the gluing can be performed and yields a compact decorated tree without appealing to  Proposition \ref{P:constructionomega-} nor to Property~{\hypersetup{linkcolor=black}\hyperlink{prop:P}{$(\mathcal{P})$}}.  Indeed, Assumption \ref{A:gamma0} is \textbf{not} fulfilled and so one cannot use Proposition \ref{P:constructionomega-} to justify \textit{a priori} that $ \mathtt{T}$ is well-defined. Also, Assumption \ref{A:omega-}  is not fulfilled either and actually the ``natural measure on the leaves of $ \mathtt{T}$"  which coincides with  the contour measure $ \frac{1}{2}\cdot\gamma_{\e}$ should be constructed via the \textit{derivative martingale}, see \cite[Theorem 5.3]{aidekon2022growth}.
 
This example shows that the construction of self-similar Markov trees can, at least in some examples, still be performed in the critical case $\min \kappa = 0$ and gives credits to the discussion in the comments section of Chapter \ref{chap:generalBP}. Furthermore, it sheds yet another point of view on the usual Brownian CRT by showing that a version of it can be constructed as critical  self-similar Markov tree. 

\chapter{Markov properties}
\label{chap:markov}

In this chapter, we discuss several Markov properties of a self-similar Markov tree  with a given characteristic quadruplet, hence justifying \textit{a posteriori} the terminology. We state them as properties satisfied by the laws of the equivalence class in $ \mathbb{T}$ of the decorated trees $ \normalfont{\texttt{T}}=(T,d_T, \rho, g)$.  Heuristically, these properties claim that certain natural families of disjoint decorated subtrees  are conditionally independent given the values of the decoration at their respective roots, and that after a proper rescaling, they are  distributed as the original self-similar Markov tree. The chapter is divided as follows. First in Section \ref{Sec:Local:Decom}, we introduce the notion of local decomposition, which provides a general and  rigorous formulation of the Markov property. Then in Sections \ref{sub:section:gene} and \ref{sub:section:temp}, we apply this framework to the case of subtrees dangling from  spines, balls or hulls. The Markov properties will play a pivotal role in our  study of self-similar Markov trees, see Chapter~\ref{chap:spinal:deco}, and are also crucial for establishing invariance principles in Part II.

\section{Local decompositions}\label{Sec:Local:Decom}
Our goal here is to introduce a general framework that enables us to state rigorously Markov properties of self-similar Markov  trees. Throughout this section, we fix a probability space on which various random variables will be defined.  Heuristically, a local decomposition of a random decorated tree $\mathtt{T}$ is a way to reconstruct it from an initial random decorated  real tree $\mathtt{T}'$ with marks  and then by gluing on the latter a family  of random decorated real trees satisfying some branching property.  This notion is similar in a random framework to the deterministic one we used in Section \ref{sec:1.1}. See also the end of Section \ref{sec:1.3}, and notably Lemma \ref{L:gluingcompatible} therein, for the version of the gluing operator in the setting of equivalence classes modulo isomorphisms of decorated real trees.

Let us now provide a formal definition. Recall from Section \ref{sec:1.3} that $ \mathbb{T}$ is the Polish space  of all isomorphism classes of rooted decorated compact trees equipped with the distance $ \mathrm{d}_{\mathbb{T}}$.     Fix $(\mathbb{Q}_x)_{x\geq 0}$ a (measurable) kernel of probability measures on $\mathbb{T}$ and $I$ some countable  set of indices. In the same probability space, we let $(\mathtt{T}', ({r}_i)_{i\in {I}})$ be  a random decorated real tree endowed with a family of marks, that is a random variable in $\mathbb{T}^{I \bullet}$, $\big( \ell_i\big)_{i\in I}$  a family of random variables in $\mathbb{R}_+$,  and finally $({\uptau}_i)_{i\in {I}}$ a family of random decorated  real trees in $\mathbb{T}$. Recall the definition of the gluing operation as a map from $\TT^{I\bullet} \times (\TT)^I$ to $\TT$, which was given after Lemma \ref{L:gluingcompatible} in  Section~\ref{sec:1.3}.

\begin{definition} \label{Def:locdec}
We say that 
   $\Big((\mathtt{T}', ({r}_i)_{i\in {I}}), (\ell_i)_{i\in {I}}, ({\uptau}_i)_{i\in {I}}\Big)$ is a  $(\mathbb{Q}_x)$-\textbf{local decomposition} of $\mathtt{T}$
   if  it satisfies the following:
\begin{itemize}
\item[(LD1)] We have 
$$\mathrm{Gluing}\Big((\mathtt{T}', ({r}_i)_{i\in {I}}), ({\uptau}_i)_{i\in {I}} \Big)= \mathtt{T}, \quad \text{  in } \mathbb{T}, \text{ a.s.}$$
\item[(LD2)] Conditionally on
$\big((\mathtt{T}', ({r}_i)_{i\in {I}}), (\ell_i)_{i\in {I}}\big)$, 
the random decorated   real trees $({\uptau}_i)_{i\in  {I}}$  are independent and the conditional law of ${\uptau}_i$ is $\mathbb{Q}_{\ell_i}$ for every $i\in {I}$.
\end{itemize}
\end{definition}
We stress that the first and third components of the triplet involved in a local decomposition are random variables in spaces of equivalent classes $\mathbb{T}^{I \bullet}$ and $\mathbb{T}$. Since by Lemma \ref{L:gluingcompatible}, the equivalence class of the glued tree in (LD1) does not depend on the choice for representatives,  the above definition does makes sense. 
In this situation, we refer to $\mathtt{T}'$ as the base of the local decomposition, and to the ${\uptau}_i$'s as \textbf{subtrees dangling from} $\mathtt{T}'$. Recall also that in the construction by gluing, only the decorated real trees ${\uptau}_i$ that are not degenerate play a role. We also stress that local decompositions are only of interest if conditions \eqref{E:nullheight} and  \eqref{E:nulldeco} are fulfilled a.s.; otherwise, we will simply have $\mathtt{T} = \mathtt{0}$.

\begin{figure}[!h]
 \begin{center}
 \includegraphics[width=10cm]{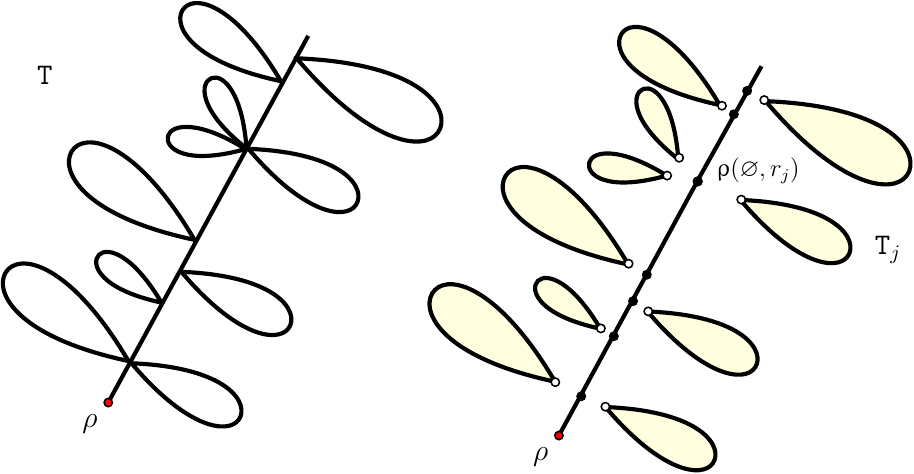}
 \caption{Illustration of the local decomposition of branching processes along the spine $T^0$: conditionally on the reproduction measure, the dangling subtrees $\texttt{T}_i$ are independent and of law $ \mathbb{Q}_{\ell_{i}}$, $i\in \N$.}  \label{fig:locdecspine}
 \end{center}
 \end{figure}

We interpret a $(\mathbb{Q}_x)_{x\geq 0}$-local decomposition  as a type of Markov property. Conditioned on the  ``present" $(\ell_i)_{i\in {I}}$, the ``past" $(\mathtt{T}', ({r}_i)_{i\in {I}})$ and the ``future" $({\uptau}_i)_{i\in  {I}}$  are independent. Furthermore, the conditional distribution of the future is determined by the probability kernel $(\mathbb{Q}_x)_{x\geq 0}$.  
This property can actually be seen as an extension of the branching property discussed in Section~\ref{sec:2.1}, within the context of random decorated trees. It can be related to the so-called \textbf{strong Markov branching property} of general branching processes  in \cite{jagers1989general}, as we shall see notably in the proof of the forthcoming Proposition \ref{prop:cutoffmarkov}.   Let us illustrate the notions that we just introduced with  a basic example involving general branching processes.
\begin{example}[Local decomposition of general branching processes along the spine] \label{Ex:locdecspine} 
Consider a decoration-reproduction kernel $(P_x)_{x>0}$ as defined in Section \ref{sec:2.1}. We write as usual $\P_x$ for the law of the family of decoration-reproduction processes $\left(f_u, \eta_u\right)_{u\in \U}$ which is induced when the ancestor has type $x$.
Assume that Property~{\hypersetup{linkcolor=black}\hyperlink{prop:P}{$(\mathcal{P})$}} holds $\P_x$-a.s. for all $x>0$. We take for $\texttt{T} \in \mathbb{T}$ the random decorated tree defined in  Theorem \ref{T:recolinfty} and denote its law under $\P_x$ by $\mathbb{Q}_x$. We also write 
$\mathbb{Q}_0$ for the Dirac point mass at the degenerate decorated real tree $\texttt{0}$ in $\mathbb{T}$. 

Let us work under $\P_y$ for some $y>0$. Recall from Notation \ref{N:subtrees}, that  $T^0=\{\uprho(\varnothing,t):~t\in [0, z_\varnothing]\}$ stands for the base subtree of $T$ induced by the ancestral individual of the general branching process. We interpret $T^0$ as the spine of $T$. We also  write $d_{T^0}$ for the restriction  of $d_T$ to $T^0$ and $g^0$ for  the associated usc-decoration -- which corresponds to the usc version of $f_\varnothing$.  Next, recall that the atoms of the reproduction process $\eta_{\varnothing}$ of the ancestor have been enumerated, say $(r_1, \ell_1), (r_2, \ell_2), \ldots$, using some deterministic rule, for instance the co-lexicographic order, and, if needed, are  completed with fictitious pairs $(\dagger, 0)$ to get an infinite sequence. For any $j\geq 1$ such that $r_j\neq \dagger$, the sub-family of decoration-reproduction processes $\left(f_{ju}, \eta_{ju}\right)_{u\in \U}$
also satisfies the Property~{\hypersetup{linkcolor=black}\hyperlink{prop:P}{$(\mathcal{P})$}} almost surely and we write  $\texttt T_j$  for the random decorated tree that it induces.
When $r_j= \dagger$, we simply decide that $\texttt T_j$ is the degenerate decorated real tree.

 On the one hand, we see from the gluing construction that 
$$\mathrm{Gluing}\Big(\big((T^0,d_{T^0}, \rho, g^0), (\uprho(\varnothing,r_i))_{i\in \N}\big), (\texttt{T}_i)_{i\in \N} \Big)= \texttt{T},$$ which is the requirement (LD1). 
  On the other hand, 
 the branching property of general branching processes discussed in Section \ref{sec:2.1} entails that conditionally on $(f_{\varnothing}, \eta_{\varnothing})$, or  equivalently  conditionally on $\big((T^0,d_{T^0}, \rho, g^0), (\uprho(\varnothing,r_i))_{i\in \N}\big)$ and $(\ell_i)_{i\in \N}$, the sub-families of decoration-reproduction processes $\left(f_{ju}, \eta_{ju}\right)_{u\in \U}$ for $j\geq 1$
are independent and the conditional law of $\left(f_{ju}, \eta_{ju}\right)_{u\in \U}$ is $\P_{\ell_j}$. Therefore 
the random decorated   real trees $(\texttt{T}_i)_{i\in  \N}$  are also conditionally independent and the conditional law of $\texttt{T}_i$ is $\mathbb{Q}_{\ell_i}$ for every $i\in \N$, which is the requirement (LD2). We conclude that 
$$\Big(\Big(\big(T^0,d_{T^0}, \rho, g^0\big), \big(\uprho(\varnothing,r_i)\big)_{i\in \N}\Big), \big(\ell_i\big)_{i\in \N},  \big(\texttt{T}_i\big)_{i\in \N}\Big)$$
  is a  $(\mathbb{Q}_x)_{x\geq 0}$ local decomposition of $\texttt{T}$; see Figure~\ref{fig:locdecspine}  for an illustration.
 \hfill $\diamond$
\end{example}

Let us now focus on our case of interest, namely when $ \mathtt{T} = (T,d_T, \rho, g)$ is built as in Proposition~\ref{P:constructionomega-} from decoration-reproduction processes $(f_{u}, \eta_{u})_{u \in \mathbb{U}}$ induced by  a characteristic quadruplet $(\sigma^{2}, \mathrm{a}, \boldsymbol{\Lambda} ; \alpha)$ satisfying Assumption \ref{A:gamma0}. Recall the notation $ \mathbb{P}_{x}$ for the law of $(f_{u}, \eta_{u})_{u \in \mathbb{U}}$ when starting from initial decoration $x>0$ and that $ \mathbb{Q}_{x}$ is the law of the equivalence class of $\mathtt{T} = (T,d_T, \rho, g)$ in $ \mathbb{T}$. Note the subtlety here: under $ \mathbb{P}_{x}$ the variable $\mathtt{T}$ is an actual and explicit decorated tree built from gluing decorated segments obtained from $(f_{u}, \eta_{u})_{u \in \mathbb{U}}$ as in Section \ref{sec:1.2}, whereas under $ \mathbb{Q}_{x}$ the variable $ \mathtt{T}$ is rather an equivalence class of such trees. Recall also that $\mathbb{Q}_0$  stands for the law of the  degenerate decorated tree $\mathtt{0}$ in $\mathbb{T}$.
  
We shall work under $ \mathbb{P} = \mathbb{P}_{1}$ and establish a decomposition of the actual decorated tree $ \mathtt{T} = (T,d_T, \rho, g)$ that will provide a local decomposition in the sense of Definition \ref{Def:locdec} once seen as variables in $ \mathbb{T}$. In this direction,  we need first to introduce  a standard procedure to decorate subtrees of $T$. Specifically, we say that $\tau\subset T$ is a subtree if it is  non-empty, connected and closed. We can endow $\tau$ with the distance $d_{\tau}$  induced by $d_T$ and  root it at the unique  point $\rho_{\tau}\in\tau$ such that $\inf_{y\in \tau} d_{T}(\rho,y)=d_T(\rho,\rho_\tau)$. Then, $({\tau}, d_{\tau}, \rho_\tau)$ is a compact rooted real tree and we further say that $\tau$ is a \textbf{base subtree} if  furthermore $\rho\in \tau$. We equip ${\tau}$  with  the decoration $g_{\tau}$ that coincides with
   $g$ on ${\tau}\setminus \{\rho_\tau\}$ and is given at the root by
    \begin{equation} \label{eq:labelatroot} g_{\tau}( \rho_{\tau}) := \limsup_{\begin{subarray}{c}y \in {\tau}\setminus \{\rho_\tau\}\\  y \to \rho_{\tau} \end{subarray}}g(y), \end{equation}
    with the convention that $ g_{\tau}( \rho_{\tau})= 0$, if $\tau=\{\rho_{\tau}\}$. We refer to  $g_{\tau}( \rho_{\tau})$ as the \textbf{germ of the decoration} on $\tau$. It is important to note that this quantity has been defined not only to ensure upper semi-continuity of the decoration, but also so that it can be evaluated by peeping only infinitesimally at the boundary point of a connected component ${\tau}\setminus \{\rho_\tau\}$, revealing the latter as little as possible. In particular, since $f_{u}(0) = \chi(u)$ in the self-similar setup, if $\tau$ is a base subtree then  we must have $g_{\tau}(\rho_\tau)=g(\rho)$. In the rest of the section we will use the standard notation ${\uptau}:=(\tau, d_{\tau}, \rho_\tau, g_\tau)$, and for definiteness we extend the definition when $\tau$ is empty by taking $\uptau:=(\{\rho\}, 0,\rho, 0)$ which is isomorphic to $\mathtt{0}$. We refer to $\uptau$ as the standard decoration of $\tau$.

   We can now explain the road map to  encode local decompositions under $\P$. 
   First,  assume that we have defined a base subtree $T'\subset T$ by means of some  algorithm or geometric definition (using the $(f_{u}, \eta_{u})_{ u \in \mathbb{U}}$). We then write $(\tau^*_i)_{i\in I}$ for the family of connected components  of $T\backslash T'$, agreeing for definitiveness 
 that some $\tau^*_i$ may be empty, notably  when the number of components is finite.    We stress that the precise choice of the indexing set $I$ is irrelevant at the stage, but will have some importance latter. Next, we write $\tau_i$ for the closure of $\tau_i^*$, in particular we have $\tau_i=\tau_i^*\cup\{\rho_{\tau_i}\}$ when $\tau_i^*$ is non empty and $\tau_i=\varnothing$ otherwise.  If, we write $r_i$ and $\ell_i$  for the root and initial usc-decoration of $\uptau_j$, that is $(r_i,\ell_i):=(\rho_{\tau_i}, g_{\tau_i}(\rho_{\tau_i}))$ if $\tau_i$ is non empty and $(r_i,\ell_i):=(\rho,0)$ otherwise, then   the requirement (LD1)  for a local decomposition is  clearly fulfilled, i.e. 
$$ \mathrm{Gluing}\Big((\mathtt{T}', ({r}_i)_{i\in {I}}), ({\uptau}_i)_{i\in {I}} \Big)= \mathtt{T}, \quad \text{  in } \mathbb{T}, \text{ a.s.}.$$
 So the remaining crucial issue  is to verify (LD2), and this requires in particular a proper choice of the indexation  of the connected components\footnote{For instance, the algorithm that ranks components in the decreasing order of their heights would not serve our purpose.}, which in our cases will use the explicit construction from $( f_{u}, \eta_{u})_{u \in \mathbb{U}}$. Under $\P$, we will say in short that  $T'\subset T$ \textbf{induces a local decomposition} of the self-similar Markov tree $\mathtt{T}$ whenever it exists such an indexation so that  $\Big((\mathtt{T}', ({r}_i)_{i\in {I}}), (\ell_i)_{i\in {I}}, ({\uptau}_i)_{i\in {I}}\Big)$ is a $ (\mathbb{Q}_{x})_{x >0}$ local decomposition of $ \mathtt{T}$. Ours proofs will consists in constructing such an indexation and checking the independence property.

\section{Genealogical Markov property}\label{sub:section:gene}
In this section we establish local decompositions induced by base subtrees that are constructed using the genealogy of Ulam's tree. Recall that we work under $\P$ and that we write $\mathtt{T}=(T,d_T,\rho, g)$ for the decorated tree  built  from  the family of decoration-reproduction processes $(f_u,\eta_u)_{u\in \mathbb{U}}$.
Recall also from Notation \ref{N:subtrees} that for every $n\geq 0$, we write $T^n$ for the base subtree of $T$ induced by individuals of the general branching process up to generation $n$ only.
The following local decomposition can be seen as an easy extension of Example \ref{Ex:locdecspine}.

\begin{proposition} \label{prop:markovsimple:gen} For every $n\geq 0$, the subtree $T^n$ induces a local decomposition of  $\mathtt{T}$ under $\P$.
\end{proposition}

\begin{proof}
We will establish the claim by induction. We start considering the case $n=0$ which has  essentially been  already discussed in Example \ref{Ex:locdecspine}, and we use the notation therein.  
The slight difference is that we decorate here the spine segment $T^0$ with the restriction of $g$ to $T^0$, whereas in Example \ref{Ex:locdecspine}, we rather used the usc-modification of  $f_{\varnothing}$.
However,  these two functions may only differ at marked points of the segment where gluing is performed. More precisely, for every $i\in \mathbb{N}$, we have 
$g(r_i)=f_{\varnothing}(r_i) \vee \max\{\ell_j: r_j=r_i\}.$
We deduce that (LD2) still holds when we decorate $T^0$ with the restriction of $g$ instead of $f_{\varnothing}$, and hence $T^0$ induces a local decomposition of  $\mathtt{T}$ under $\P$.

For $n=1$, we decompose in turn each (non-degenerate) connected component along its own spine,  using $I=\N^2$ to index the family of subtrees at the second generation.
And so on, and so forth, generation after generation. For the sake of clarity, let us spell the elements of the local decomposition out. At generation $n$, we use $I=\N^{n+1}$ as set of indices. For $u\in \N^{n+1}$, $r_u=\uprho(u)$ is the equivalence class of the root  of the segment labelled by $u$ in the construction by gluing, and $\ell_u=\chi(u)$ is the type of the individual $u$, that is the germ of the decoration $g_u(\rho_u)$ for the subtree  $\tau_u=T_u$, again with the Notation \ref{N:subtrees}. Condition (LD2) follows directly from the branching property.
 \end{proof}

In the remaining of the section, we focus  on the model of pruned tree below a fixed level, which can be though of as a variation of $T^n$ where we kill individuals with types lless than $\varepsilon$
(together with their descendance), thus only keeping finitely many individuals with types at least $\varepsilon$. We shall prove after, see Corollary \ref{lem:cutoff},  that the associated decorated subtree $ \mathtt{T}^{[ \varepsilon]}$ is an approximation of $ \mathtt{T}$ as $ \varepsilon \to 0$; 
 and this  will be notably useful in Part II, when we shall consider scaling limits of Galton-Watson branching processes with integer types.

 To define properly the pruning transformation, we fix some  threshold $\varepsilon\in(0,1)$ and introduce the set of vertices
\begin{equation}\label{E:Lepsi}
F(\varepsilon)\coloneqq\big\{u\in \U: f_{u}(0)<  \varepsilon \text{ and } f_{v}(0)\geq \varepsilon \text{ for all }v\prec u\big\},
\end{equation}
where  the notation $v\prec u$ means that $v$ is a strict prefix of $u$ in $\U$. For every  $u\in F(\varepsilon)$, we change the entire descent of  $u$ and make it fictitious by setting $(f_v^{[\varepsilon]},\eta_v^{[\varepsilon]} )=(0,0)$ for all $ v\succeq u$ ($u$ prefix of $v$).
The decoration-reproduction of the  ascendants of individuals of $F(\varepsilon)$  are unchanged, i.e.
$(f_v^{[\varepsilon]},\eta_v^{[\varepsilon]} )=(f_v,\eta_v)$ for all $v\prec u$ with $u\in F(\varepsilon)$. 
We know from  Lemma \ref{L:meanx} that the family 
$$A(\varepsilon)\coloneq \big\{v\in \U: v \prec u \text{ for some }u\in F(\varepsilon)\big\}$$  is finite $\P$-a.s.
The  tree pruned at level $\varepsilon>0$, denoted by $T^{[\varepsilon]}$, is the real tree obtained by gluing  as in Section \ref{sec:2.1} the line segments $[0,z_v]$ for $v\in A(\varepsilon)$ only, i.e.
$$T^{[\varepsilon]}:=\big\{\uprho(u,t):~u\in A(\varepsilon) \text{ and } t\in [0,z_u]\big\},$$ 
where we recall that $\uprho(u,t)$ denotes the -- equivalent class of the -- point on the segment $S_u$ at distance $t$ from the root $\rho(u)$.
In particular,  $T^{[\varepsilon]}$ is a base subtree  of $T$ built from a finite number of line segments.

\begin{proposition} \label{prop:cutoffmarkov}
For every $\varepsilon\in (0,1)$, the subtree $T^{[\varepsilon]}$ induces a  local decomposition of  $\mathtt{T}$ under~$\P$.
\end{proposition}
We stress that the heart of the proof is that the set of vertices $F(\varepsilon)$ in \eqref{E:Lepsi} is an optional line in the sense of Jagers \cite{jagers1989general} to which the strong Markov branching property applies.
Other choices of optional lines would thus yield analogous local decompositions. 

\begin{proof} First, for every $v\in F(\varepsilon)$  non-fictitious, we consider the subtree $$\tau_{v}=\big\{\uprho(vu,t):~u\in \U \text{ and } t\in [0,z_{vu}] \big\}, $$ and we write $\uptau_v$ for the induced  decorated tree. By construction, the connected components of $T\setminus T^{[\varepsilon]}$ are precisely the $\tau_{v}$, for  $v\in F(\varepsilon)$  non-fictitious. Furthermore, recalling  Notation \ref{N:subtrees}, we see that $\uptau_v$ is isomorphic to the decorated tree $\mathtt{T}_v$  associated with the family $(f_{vu}, \eta_{vu})_{u\in \U}$. We also let $r_v=\uprho(v)$ and $\ell_v=\chi(v)=f_v(0)$. For convenience,  we extend the above construction to all the Ulam tree  by taking $\uptau_v=(\{\rho\}, 0,\rho, 0)$, $r_v=\rho$ and $\ell_v=f_v(0)$,  for  $v\notin F(\varepsilon)$ or $v$ fictitious.

We have to verify (LD2) in this setting, which can be derived from Proposition \ref{prop:markovsimple:gen} by constructing $T^{[\varepsilon]}$ recursively  branch by branch. Nonetheless, since the same argument can be used to establish many other local decompositions,   it may be more instructive to rather
use the strong Markov branching property of general branching processes, which has been developed in great generality in \cite{jagers1989general}.

We first observe that the set of vertices $F(\varepsilon)$  is a  random line, that is $F(\varepsilon)$ does not contain two vertices $u$ and $v$ with $u\prec v$. 
It is furthermore optional, in the following sense. For any subset of vertices $V\subset \U$, the event $\{F(\varepsilon)\preceq V\}$  that  every $v\in V$ has some prefix in $F(\varepsilon)$,  is measurable with respect to the sigma-algebra $\mathcal F_{V}$ generated by $\left((f_w,\eta_w): w\not \succeq v \text{ for all } v\in V\right)$. In other words, the event  $\{F(\varepsilon)\preceq V\}$ does not depend on the decoration-reproduction processes indexed by vertices with a prefix in $V$.  

Then define the pre-$F(\varepsilon)$ sigma-algebra $\mathcal{F}_{F(\varepsilon)}$ of  events $A$ such that $A\cap\{F(\varepsilon)\preceq V\}\in \mathcal{F}_V$ for all $V\subset \U$. The strong Markov branching property \cite[Theorem 4.14]{jagers1989general} of general branching processes at optional stopping lines can now be stated in our framework as follows. Consider for every $u\in \U$ a measurable functional $\varphi_u: \mathbb{T}\to [0,1]$; then there is the identity
$$\E\Big( \prod_{u\in F(\varepsilon)} \varphi_u(\mathtt{T}_u)~\Big|~ \mathcal F_{F(\varepsilon)}\Big) =  \prod_{u\in F(\varepsilon)} \E_{\chi(u)}\big(\varphi_u(\mathtt{T})\big).$$
Let us now impose furthermore $\varphi_u(\mathtt{0})=1$ for a while, so that we can rewrite the preceding in the form  
$$\E\Big( \prod_{u\in \U} \varphi_u(\uptau_u) ~\Big| \mathcal F_{F(\varepsilon)}\Big) =  \prod_{u\in \U} \E_{\ell_u}\big(\varphi_u(\mathtt{T})\big).$$
It is readily seen that this identity remains valid if we drop the requirement $\varphi_u(\mathtt{0})=1$. 
Indeed, it suffices to replace $\varphi_v$ by $\indset{\mathtt{0}}$ for any given vertex $v\in \U$, use linearity, and repeat the operation for every vertex in $ \U$.
Since $(\mathtt{T}^{[\varepsilon]}, (r_v)_{v\in \U})$ and $(\ell_v)_{v\in \U}$, as variable in $\mathbb{T}^{\U\bullet}$ and $\mathbb{R}_+^{\U}$, are measurable with respect to the sigma-algebra $\mathcal{F}_{F(\varepsilon)}$, this shows  (LD2) and hence completes the proof of the proposition.
\end{proof}

We conclude this section by establishing that the subtree pruned at level $\varepsilon$ approximates the self-similar Markov tree as $\varepsilon \to 0+$.  We shall furthermore endow $T^{[ \varepsilon]}$ with approximations of the weighted length and harmonic measures constructed in Section \ref{sec:2.3}. More precisely, for any $\gamma \geq \gamma_{0}$, the tree $T^{[ \varepsilon]}$ can be equipped with the restriction $\indset{T^{[\varepsilon]}}\cdot \uplambda^\gamma$ of the length measure $ g^{\gamma-\alpha} \uplambda_{T}$ constructed in Proposition \ref{prop:lengthmeasures}. However, when Assumption \ref{A:omega-} holds, the harmonic measure $\upmu$  assigns zero mass to $T^n$ for any $n\geq 0$, and therefore also to the pruned subtree $T^{[\varepsilon]}$ for any $\varepsilon >0$.  In particular the restriction of $\upmu$ to $T^{[\varepsilon]}$ is always null and does not converge to $\upmu$ as $\varepsilon \to 0+$. For this reason, we shall then rather equip $T^{[\varepsilon]}$ with another measure, namely we let $\upmu^{[\varepsilon]}$ be  the image of the harmonic measure $\upmu$ by the canonical  projection of $T$ on $T^{[\varepsilon]}$. Specifically, we introduce
\begin{equation}\label{eq:def:nu:varepsilon}
\upmu^{[\varepsilon]} \coloneqq \sum_{u\in F(\varepsilon)} \upmu(  T_u)  \delta_{ \uprho(u)}, 
\end{equation}
where the optional line $F(\varepsilon)$ has been defined in \eqref{E:Lepsi}, and, as usual, $ \uprho(u)$ denotes the point in $T^{[\varepsilon]}$ at which the segment $S_u$ is glued. Then we have the following property.

\begin{corollary}[Convergence of cutoff approximations] \label{lem:cutoff} In the notation above, the following assertions hold $\P$-a.s. 
\begin{enumerate}
\item[(i)] Suppose Assumption \ref{A:gamma0}. Then for any $\gamma\geq \gamma_0$, we have
$$\lim_{\varepsilon\to 0+} \big( T^{[ \varepsilon]},  d_{T^{{[ \varepsilon]}}}, \rho, g_{T^{{ [\varepsilon]}}}, \indset{T^{[\varepsilon]}}\cdot \uplambda^\gamma \big) = \big( T, d_T, \rho, g,  \uplambda^\gamma \big),  \qquad \text{in }\mathbb{T}_{m}.$$
\item[(ii)] Suppose Assumption \ref{A:omega-}. Then we have
$$\lim_{\varepsilon\to 0+} \big( T^{[ \varepsilon]},  d_{T^{{[ \varepsilon]}}}, \rho, g_{T^{{ [\varepsilon]}}},  \upmu^{[ \varepsilon]}\big) = \big( T, d_T, \rho, g,   \upmu\big), \qquad \text{in }\mathbb{T}_{m}.$$
\end{enumerate}
\end{corollary}

\begin{proof} We will establish the statements with convergence in probability instead of almost surely, as then the sharper claims follow readily by an argument of monotonicity.
 To start with, recall from Proposition \ref{P:constructionomega-} and Lemma \ref{L:meanx} that
$$\E\Big(\sum_{u\in \U}\chi(u)^{\gamma_0}\Big) < \infty.$$
As a consequence, we have
$$\lim_{\varepsilon\to 0+} \E\Big(\sum_{u\in \U}\indset{\{\chi(u)<\varepsilon\}}\chi(u)^{\gamma_0}\Big) =0,$$
and \textit{a fortiori}
\begin{equation}\label{a:fortiori:T:eps}
\lim_{\varepsilon\to 0+} \E\Big(\sum_{u \in F(\varepsilon)}\chi(u)^{\gamma_0}\Big) =0.
\end{equation}
 Next, we deduce from the local decomposition  Proposition \ref{prop:cutoffmarkov}, Corollary \ref{C:CMJ} and the scaling property that there is some finite constant $c>0$ such that for every 
$\varepsilon>0$, we have
$$\E\Big(\sum_{u \in F(\varepsilon)}\mathrm{Height}(T_u)^{\gamma_0/\alpha}\Big) \leq c\cdot \E\Big(\sum_{u \in F(\varepsilon)} \chi(u)^{\gamma_0} \Big)$$
and
$$\E\Big(\sum_{u \in F(\varepsilon)}\max_{T_u} g_u^{\gamma_0}\Big) \leq c\cdot  \E\Big(\sum_{u \in F(\varepsilon)}\chi(u)^{\gamma_0}\Big).$$
We deduce that
$$\lim_{\varepsilon\to 0+} \E\Big(\sum_{u \in F(\varepsilon)}\Big(\mathrm{Height}(T_u)^{\gamma_0/\alpha} + \max_{T_u} g_u^{\gamma_0}\Big)\Big) =0.$$
Since the decorated real tree $\texttt{T}$ can be recovered by gluing on $\texttt{T}^{[\varepsilon]}$ the subtrees $T_u$ for $u\in F(\varepsilon)$, 
there is the bound
$$\d_{\TT}(\texttt{T}, \texttt{T}^{[\varepsilon]}) \leq \sup_{u\in F(\varepsilon)} \Big( \mathrm{Height}(T_u) \vee \max_{T_u} g_u\Big),
$$
and we conclude that
$$\lim_{\varepsilon \to 0+} \d_{\TT}(\texttt{T}, \texttt{T}^{[\varepsilon]})=0, \qquad\text{in probability.}$$
Note that this immediately yields the claim (ii) (again with convergence in probability in place of almost surely), since from the definition \eqref{eq:def:nu:varepsilon},  the Prokhorov distance between $\upmu^{[\varepsilon]}$ and $\upmu$  is bounded from above by the Hausdorff distance 
between $T^{[\varepsilon]}$ and $T$. 

We are thus left to deal with the weighted length measure in (i). In this direction, remark that for every $\gamma\geq \gamma_0$ we have
$$ \d_{\mathrm{Prok}}(\uplambda^{\gamma}, \indset{\texttt{T}^{[\varepsilon]}}\cdot \uplambda^{\gamma}) \leq  \uplambda^{\gamma}(T \backslash T^{[\varepsilon]}) \leq \big(\max_{T} g^{\gamma-\gamma_0} \big) \cdot \uplambda^{\gamma_0}(T \backslash T^{[\varepsilon]}).  $$
Consequently,  by Corollary \ref{C:CMJ},  to conclude it suffices to establish that $\uplambda^{\gamma_0}(T \backslash T^{[\varepsilon]})$ converges to $0$ in probability, as $\varepsilon \downarrow 0$. To this end note that, from
 Proposition \ref{prop:cutoffmarkov} and the scaling property, that for every 
$\varepsilon>0$
$$
\E\big(\sum_{u \in F(\varepsilon)}\uplambda^{\gamma_0}(T_u)\big) = \E\big(\uplambda^{\gamma_0}(T)\big) \cdot \E\big(\sum_{u \in F(\varepsilon)}\chi(u)^{\gamma_0}\big)=-\frac{1}{\kappa(\gamma_0)}\cdot \E\big(\sum_{u \in F(\varepsilon)}\chi(u)^{\gamma_0}\big),
$$
where to obtain the second equality we used \eqref{mean:uplambda:gamma:T}. The desired results follows now by \eqref{a:fortiori:T:eps}.
\end{proof}

\begin{remark}[An intrinsic cutoff] The reader may compare the subtree $ \mathtt{T}^{[ \varepsilon]}$ with the subtree $ \mathtt{T}^{( \varepsilon)}$ obtained heuristically from $ \mathtt{T}$ by starting from the root and cutting each branch as soon as the decoration drops below level $ \varepsilon$. The inclusion $ \mathtt{T}^{ ( \varepsilon)} \subset \mathtt{T}^{[ \varepsilon]}$ should be plain. The nice feature of $ \mathtt{T}^{( \varepsilon)}$ is that its equivalence class in $ \mathbb{T}$ only depends on that of $ \mathtt{T}$, that is, it is an intrinsic geometric subtree, as opposed to $ \mathtt{T}^{[ \varepsilon]}$ which uses the construction from the family of decoration-reproduction processes $(f_{u}, \eta_{u})_{ u \in \mathbb{U}}$. However, it will be technical simpler to deal with $ \mathtt{T}^{[ \varepsilon]}$ rather than $ \mathtt{T}^{( \varepsilon)}$.

\begin{figure}[!h]
 \begin{center}
 \includegraphics[width=5.2cm]{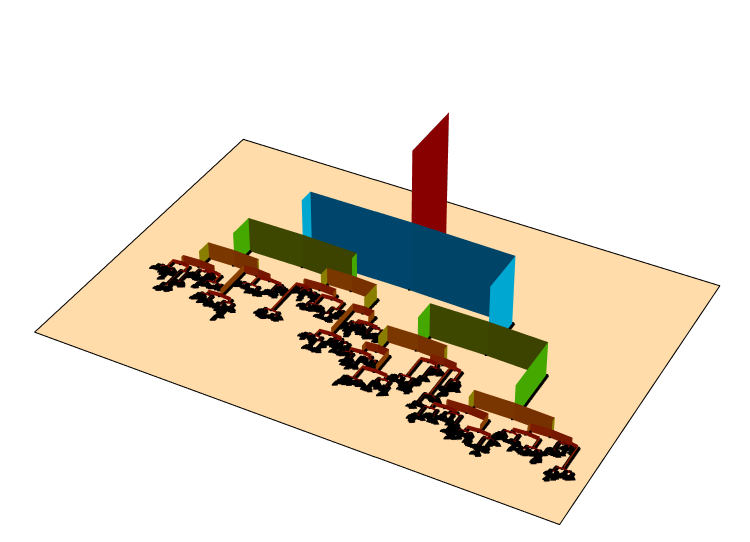}
   \includegraphics[width=5.2cm]{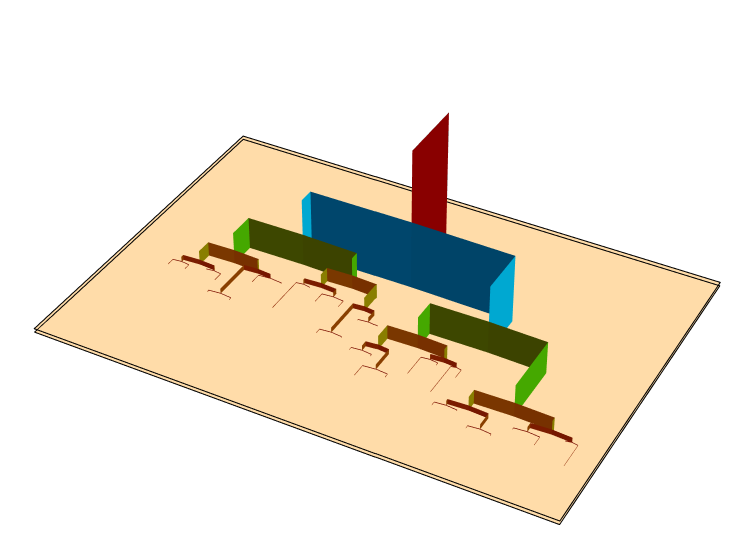}
     \includegraphics[width=5.2cm]{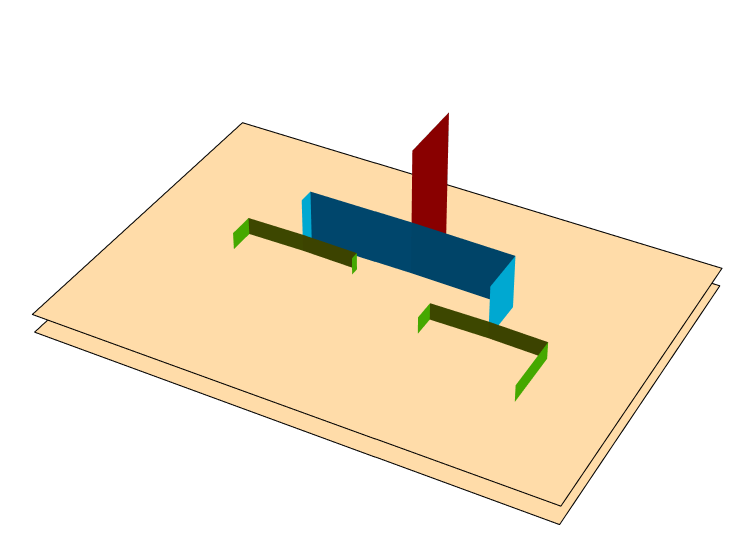}
 \caption{Illustration of the decorated subtrees $T^{{ ( \varepsilon)}}$ pruned at increasing levels from left to right. The subtrees $ T^{{[  \varepsilon]}}$ would actually be larger since they contain all branches whose initial decoration is larger than $   \varepsilon$.}
 \end{center}
 \end{figure}
\end{remark} 
 
\section{Temporal Markov property}\label{sub:section:temp}
Roughly speaking, the local decompositions described in Propositions \ref{prop:markovsimple:gen} and \ref{prop:cutoffmarkov} only rely on the genealogical branching property of self-similar Markov trees, and  could  have been stated as well for general branching processes; see Example \ref{Ex:locdecspine} and also the proof of Proposition~\ref{prop:cutoffmarkov}. In this section, we turn our attention to another type of local decompositions for self-similar Markov trees, which now rather stems
from the following temporal Markov property of decoration-reproduction processes. Recall that $(P_x)_{x\geq 0}$ denotes the self-similar decoration-reproduction kernel
with characteristic quadruplet $ (\sigma^2, \mathrm{a}, \boldsymbol{\Lambda} ; \alpha)$.
For every $t\geq 0$, we use the notation $(f,\eta)\circ \theta_t$ for the pair of shifted processes at time $t$, that is $(f,\eta)\circ \theta_t= (0,0)$ if $t\geq z$, and  otherwise,
$$f\circ \theta_t: [0,z-t]\to \R_+, \quad f\circ \theta_t(s)=f(t+s) \quad \text{for }0\leq s \leq z-t,$$
and
$$\int_{[0,\infty)\times \R_+}\varphi(s,x)\eta\circ \theta_t(\d  s, \d  x)= \int_{(t,\infty)\times \R_+}\varphi(s-t,x)\eta(\d  s, \d  x),$$
 for any measurable function  $\varphi: [0,\infty)\times \R_+ \to \R_+$.

\begin{lemma}[Markov property of the decoration-reproduction process] \label{L:Markovdecoreprod} 
Write $(\mathcal G_t)_{t\geq 0}$ for the natural filtration of the decoration-reproduction process $(f,\eta)$, i.e. $\mathcal G_t$ is the (completed) sigma-algebra generated by
$\indset{[0,t]}\cdot f$ and $\indset{[0,t]\times \R_+}\cdot \eta$. 
For every $x\geq 0$ and every $(\mathcal G_t)$-stopping time $R$, the conditional distribution under $P_x$  of the shifted decoration-reproduction process $(f,\eta)\circ \theta_R$ given $\mathcal G_R$ is $P_{f(R)}$. 
\end{lemma}

\begin{proof} The claim is an easy extension of the strong Markov property of positive self-similar Markov processes;  we merely sketch below the argument and refer to \cite[Chapter 13]{kyprianou2014fluctuations} for more details. Recall the setting of Section \ref{sec:2.2}, and in particular that $\mathbf{N}=\mathbf{N}(\d t, \d y, \d \mathbf{y})$ is a Poisson random measure on $[0, \infty) \times \mathcal{S}$ with intensity measure $\d t \boldsymbol{\Lambda}(\d y, \d \mathbf{y})$, and  $\xi$ the L\'evy process with L\'evy-It\^{o} decomposition \eqref{E:LevyIto} and possibly killed at the exponentially distributed time $\zeta$. Writing $(\mathcal H_t)_{t\geq 0}$ for the natural filtration of $( \xi, \mathbf{N})$, we have  that
 for every $t\geq 0$ and  conditionally on $t<\zeta$,  the pair $( \xi^{(t)}, \mathbf{N}^{(t)})$ given by 
$$ \xi^{(t)}(s)=\xi(s+t)-\xi(t),\quad \text{ for } s\geq 0, $$
and 
$$\int_{[0, \infty) \times \mathcal{S}} \varphi(s,y,\mathbf{y})\mathbf{N}^{(t)}(\d s, \d y, \d \mathbf{y}) = \int_{(t, \infty) \times \mathcal{S}} \varphi(s-t,y,\mathbf{y}) \mathbf{N}(\d s, \d y, \d \mathbf{y}),$$
is independent of $ \mathcal H_t$ and has the same law as $( \xi, \mathbf{N})$. This observation extends to the situation where the fixed time $t$ is replaced, first by a simple $(\mathcal H_t)$-stopping time, and then, by approximation, by any a.s. finite  $(\mathcal H_t)$-stopping time.
The statement then follows readily by inspection of the Lamperti transformation described in Section \ref{sec:2.2}.
\end{proof}

We now derive from Lemma \ref{L:Markovdecoreprod} the first two temporal local decompositions that are both induced by a stopping time $R\leq z_{\varnothing}$ in the natural filtration of the decoration-reproduction process of the ancestor $(f_\varnothing, \eta_\varnothing)$. We consider the spine truncated at distance $R$ from the root, that is the segment $T_R^0:= \{\uprho(\varnothing,r):~ r\in[0, R]\}\subset T^0$. We consider also the hull generated by the truncated spine,
$$B_R^\bullet(T)\coloneqq \big\{y\in T: d_T(\rho,p_{\varnothing}(y))\leq R\big\},$$
where $p_{\varnothing}: T\to T^0$ denotes the projection on the spine and $\rho$ stands for the root of $T$.  So roughly speaking, the hull $B_R^\bullet(T)$ is the subtree that results from killing the ancestor immediately after time $R$ in the population model. The closure of its complement $\check{B}_R^\bullet(T)$ is connected and hence also a real tree on the event $R< z_\varnothing$, and empty on the complementary event. Notice that when $R< z_\varnothing$, we have
$$\check{B}_R^\bullet(T)=\big\{y\in T: d_T(0,p_{\varnothing}(y))>R\big\}\cup \big\{\uprho(\varnothing, R)\big\}.$$

\begin{proposition} \label{prop:hullmarkov}
For every stopping time $R\leq z_{\varnothing}$ in the natural filtration of the ancestral decoration-reproduction process $(f_\varnothing, \eta_\varnothing)$, both the truncated spine $T^0_R$ and the hull $B_R^\bullet(T)$ induce a   local decomposition of  $\mathtt{T}$ under $\P$.
\end{proposition}

\begin{proof} The claim readily follows from a variation of the $(\mathbb Q_x)_{x\geq 0}$ local decomposition along the spine, i.e. Proposition \ref{prop:markovsimple:gen} when $n=0$. Recall the notation there and in Example \ref{Ex:locdecspine}, and notably 
 that the atoms of the ancestral reproduction process $\eta_{\varnothing}$ are enumerated in the  co-lexicographic order, $(r_1, \ell_1), (r_2, \ell_2), \ldots$, and that for every $j\geq 1$, 
 we denoted by $\mathtt T_j$ the random decorated tree   induced by the sub-family  $\left(f_{ju}, \eta_{ju}\right)_{u\in \U}$.
We have seen  in   the beginning of the proof of Proposition \ref{prop:markovsimple:gen}  that
\begin{equation}\label{Eq:locdecspine}
\big((\mathtt{T}^0, (\uprho(\varnothing,{r}_i))_{i\in \N}), (\ell_i)_{i\in \N}, (\mathtt{T}_i)_{i\in \N}\big) \ \text{ is a $(\mathbb{Q}_x)_{x\geq 0}$ local decomposition of } \ \mathtt{T}.
\end{equation}

 \begin{figure}[!h]
 \begin{center}
 \includegraphics[width=10cm]{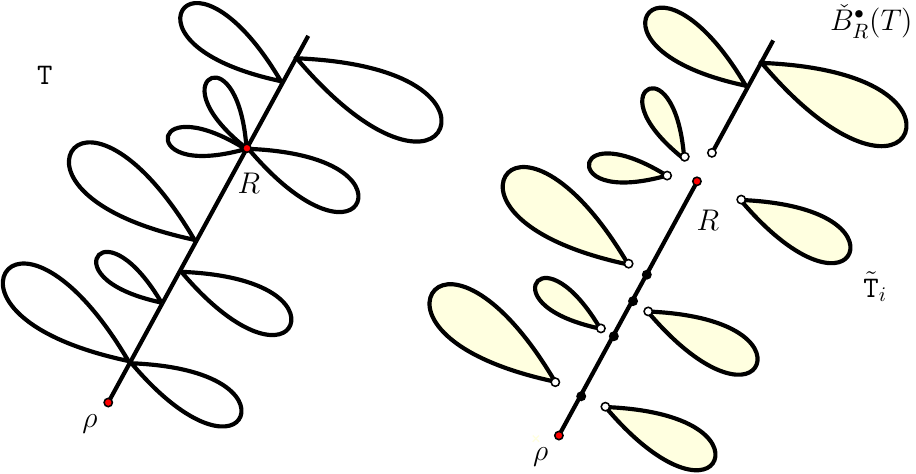}
 \caption{Illustration of the Markov property for the truncated spine: conditionally on the reproduction measure, all dangling subtrees are independent and of law $ \mathbb{Q}_{\ell_{i}}$ where $\ell_{i}$ are the germ decorations.}
 \end{center}
 \end{figure}

 We now  re-order those elements depending on their positions relative to the stopping time $R$ as follows, taking $\mathbb Z^*=\mathbb Z\backslash \{0\}$ as new set of indices.
 We construct recursively $1\leq R(1)<R(2)< \ldots$ such that for any $j\geq 1$, $R(j)=k$ where  $k$ is the rank of the $j$-th atom $(r_k, \ell_k)$  such that $r_k>R$ if any, and $R(j)=\infty$ otherwise.
 We define similarly $1\leq R(-1)<R(-2)< \ldots$ such that for any $j\geq 1$,  $R(-j)=k$ where  $k$ is the rank of the $j$-th atom $(r_k, \ell_k)$  such that $r_k\leq R$ if any, and $R(-j)=\infty$ otherwise.
 We then set for any $j\in \mathbb Z^*$, $\tilde r_j=r_{R(j)}$, $\tilde \ell_j=\ell_{R(j)}$ and $\widetilde{\mathtt T}_j=\mathtt T_{R(j)}$ when $R(j)<\infty$, and $\tilde r_j=0$, $\tilde{\ell}_j=0$ and $\widetilde{\mathtt T}_j=\mathtt{0}$ otherwise. 
It is immediate to check from  \eqref{Eq:locdecspine}, the fact that $R$ is measurable with respect to
$((\mathtt{T}^0, (\uprho(\varnothing,{r}_i))_{i\in \N}), (\ell_i)_{i\in \N})$, and the fact that the original ordering can be recovered from the new one,  that 
\begin{equation}\label{Eq:locdecspinetilde} \big((\mathtt{T}^0, (\uprho(\varnothing,\tilde r_i))_{i\in \mathbb Z^*}),(\tilde{\ell}_i)_{i\in \mathbb Z^*}, (\widetilde{\mathtt T}_i)_{i\in \mathbb Z^*}\big)
 \ \text{ is again a $(\mathbb{Q}_x)_{x\geq 0}$ local decomposition of } \ \mathtt{T}.
 \end{equation}

  We then split the spine at $R$. First, write $\widetilde{\mathtt T}_0$ for  the decorated version of $\check{B}_R^\bullet(T)$ and note that it  can be obtained by taking the set  $\{\uprho(\varnothing,r):~r\in [R,z_\varnothing]\}$ with its associated decoration, and then gluing on it the family $(\widetilde{\mathtt T}_i)_{i\geq 1}$ at the position prescribed by the points $(\uprho(\varnothing,\tilde r_i))_{i\geq 1}$. Second, write $\texttt{T}_0^R$ for the decorated truncated spine and observe that the variables $(\mathtt{T}_{R}^0, (\uprho(\varnothing,\tilde r_i))_{i\leq -1})$ and $(\tilde{\ell}_i)_{i\leq -1}$ are measurable in the  sigma-algebra  $\mathcal{G}_R$ generated by the ancestral decoration-reproduction process stopped at time $R$. By comparison with the local decomposition along the spine in Proposition~\ref{prop:markovsimple:gen}  and the temporal Markov property in Lemma~\ref{L:Markovdecoreprod}, we deduce that the conditional law of 
 $\widetilde{\mathtt T}_0$ given  $$\Big(\big(\mathtt{T}_{R}^0, \big(\uprho(\varnothing, \tilde r_i)\big)_{i\leq -1}\big),(\tilde{\ell}_i)_{i\leq -1}, (\widetilde{\mathtt T}_i)_{i\leq -1}\big)$$ is $\mathbb{Q}_{f_{\varnothing}(R)}.$
  Thus writing $\tilde r_0=R$ and $\tilde \ell_0=f_{\varnothing}(R)$, 
  we conclude from \eqref{Eq:locdecspinetilde} that 
  $$\Big(\big(\mathtt{T}_{R}^0, \big(\uprho(\varnothing, \tilde r_i)\big)_{i\leq 0}\big),(\tilde{\ell}_i)_{i\leq 0}, (\widetilde{\mathtt T}_i)_{i\leq 0}\big)$$
  is a  $(\mathbb{Q}_x)$-local decomposition of $\mathtt{T}$; we have thus proved the claim for the truncated spine. 
  
  Finally the claim for the hull follows again from \eqref{Eq:locdecspinetilde} and the observation that  the decorated version of $B_R^{\bullet}(T)$ can be obtained by gluing the sequence of decorated trees $(\widetilde{\mathtt T}_i)_{i\leq -1}$ 
on the truncated spine $\mathtt{T}_{R}^0$  at the locations $(\uprho(\varnothing,\tilde r_i))_{i\leq -1}$. 
\end{proof}

Similarly as in the previous section, one can iterate the local decomposition of the truncated spine to infer new Markov properties. We illustrate  this procedure with another important example  of base subtrees, namely the \textbf{closed ball} of radius $a$,
 $$B_a(T):= \big\{\upsilon\in T:~{d}_{T}(\rho,\upsilon)\leq a\big\}.$$

\begin{proposition}\label{prop:danglingcutoff} 
For every $a>0$, the ball $ B_a(T)$ induces a  local decomposition of  $\mathtt{T}$ under~$\P$.
\end{proposition}

\begin{proof} Recall that $T^n$  denotes the base subtree of $T$ induced by individuals of the general branching process up to generation $n$ only, and  set $B_a(T^n)\coloneqq B_a(T)\cap T^n$. 
For $n=0$, the subtree $B_a(T^0)$ is merely the spine  truncated at the fixed height $a$, and we know from Proposition \ref{prop:hullmarkov} that  $B_a(T^0)$  induces a  local decomposition of  $\mathtt{T}$ under $\P$.
Just as in Proposition \ref{prop:markovsimple:gen}, we  can then recursively decompose $\mathtt{T}$ along the truncated segments at generations $1, 2, \ldots$, and get that for any $n\geq 0$, the base tree
  $B_a(T^n)$  induces a  local decomposition of  $\mathtt{T}$ under $\P$. In this framework, 
the connected components of $T\backslash B_a(T^n)$ are indexed  by the set of vertices $u\in \U$.

 More specifically, for every $u\in \mathbb{U}$ with $|u|\leq n$,  note that there is at most one point $\uprho(u,t)$ in $\{\uprho(u,s):~s\in [0,z_u]\}$ at distance $a$ from $\rho$. When there is such a point,   we set $r_u^n:=\uprho(u,t)$,  $\ell_u^n:=f_u(t)$ and we write $\tau_u^n$ for the subtree formed by the set of points $v\in T$ such that $\uprho(u,t)\in \llbracket \rho ,v \rrbracket$. Remark that  we have  $t=a-d_T(\rho, \uprho(u))$ and   $\tau_u^{n}$ can be seen as the  closure of the  complement of the hull of radius $a-d_T(\rho, \uprho(u))$ of $T_u$.  If there no such point or if $|u|>n$, we let $r_u^n:=\rho$, $\ell_u^n:=0$ and $\tau_u^{n}:=0$. If we write $\uptau_u^n$ for the decorated version of $\tau_u^n$ and $\mathtt{B}_a(T^n)$ for the decorated version of $B_a(T^n)$, by recursively applying Proposition \ref{prop:markovsimple:gen} we infer that, under $\mathbb{P}$ the family:
 $$\Big(\big(\mathtt{B}_a(T^n), (r_u^n)_{u\in \mathbb{U}}\big), \big(\ell_u^n\big)_{u\in \mathbb{U}}, \big(\uptau_u^n\big)_{u\in\mathbb{U}}\Big) $$
 is a $(\mathbb{Q}_{x})_{x\geq 0}$ local decomposition of $\mathtt{T}$.   We stress that the labeling is consistent as we let generations increase, that is, for any vertex $u$ with $|u|\leq n$,  the mark $r_u^n$, the variable $\ell_u^n$ and the  subtree $\uptau_u^n$ for  the decomposition 
 of $T\backslash B_a(T^n)$ are identical to those obtained for the decomposition of $T\backslash B_a(T^m)$,  for any $m>n$. With the obvious notation, for $u\in \mathbb{U}$,  let us write $r_u,\ell_u$ and $\uptau_u$ for the corresponding terminal value.   Now, using  that $\big(\| f_{u}\| : u\in  \mathbb{U}\big)$  is a null family, the fact that $\sum_{u\in \mathbb{U}}z_u<\infty$ a.s. and the definition of $(\mathbb{T}^{\mathbb{U}\bullet}, \d_{\mathbb{T}^{\mathbb{U}\bullet}})$ given in \eqref{eq:def:d:bullet}, it is straightforward  to infer that $\big(\mathtt{B}_a(T^n), (r_u^n)_{u\in \mathbb{U}}\big)$ converges to    $\big(\mathtt{B}_a(T), (r_u)_{u\in \mathbb{U}}\big)$, where $\mathtt{B}_a(T)$ stands for the standard decorated version of $B_a(T)$. Therefore, taking the limit when $n\to \infty$, we deduce that
  $$\Big(\big(\mathtt{B}_a(T), (r_u)_{u\in \mathbb{U}}\big), \big(\ell_u\big)_{u\in \mathbb{U}}, \big(\uptau_u\big)_{u\in\mathbb{U}}\Big) $$
 is a $(\mathbb{Q}_{x})_{x\geq 0}$ local decomposition of $\mathtt{T}$.  This completes the proof of the proposition. 
\end{proof}
We now conclude this chapter by stressing a connection with measured-valued branching processes; we refer to \cite{li2011measure} for background and precise definition.
Consider for every $a\geq 0$, the random integer-valued measure on $(0,\infty)$
$$Z_a\coloneqq \sum_{i\in I} \delta_{\ell_i(a)},$$
where $\ell_i(a)$ denotes the germ decoration of the connected component $\tau_i^*$ of $T\backslash B_a(T)$ and we implicitly discard empty components with $\ell_i(a)=0$ in the sum. 
Plainly, $Z_a$ does not depend on the choice for indexing these connected components. We then deduce from Proposition \ref{prop:danglingcutoff} that
for ever $a\geq 0$, conditionally on $Z_a$, the shifted process $(Z_{a+b})_{b\geq 0}$ is independent of the process $(Z_{b})_{0\leq b\leq a}$. Moreover, the distribution of the shifted process 
given $Z_a$ is that of $\sum_{i\in I} Y^i$, where the $Y^i$'s are independent measured-valued processes and each $Y^i$ has the same law as $Z$ under $\mathbb Q_{\ell_i(a)}$.
In short, 
$(Z_a)_{a\geq 0}$ is measured-valued branching process. For instance,  $(Z_a)_{a\geq 0}$ may be a self-similar fragmentation process as discussed in Section \ref{sec:5.2}, or a growth-fragmentation process as in Section \ref{sec:5.3}.

\section*{Comments and bibliographical notes}

The Markov property is of course a most important concept in probability theory. Born with the theory of stochastic processes, it has been developed for more sophisticated objets such as  branching processes \cite{jagers1989general} (the stopping lines used many times in this chapter), superprocesses \cite{dynkin1991branching} (the so-called special Markov property) and random snakes, see \cite{LeG99,DLG02} and \cite{riera2024structure}. The Gaussian Free Field (GFF, in abbreviation) is essentially characterized by its domain Markov property \cite{berestycki2020characterisation} and random sets coupled with the GFF that satisfy a strong Markov property are called \textbf{local sets}, see \cite{miller2016imaginary} for their introduction and \cite{werner2021lecture} for a comprehensive survey. In particular, our formulation of the Markov property (Definition \ref{Def:locdec}) is inspired from GFF local sets and we use the terminology \textit{local decomposition} to emphasize this connection.  See also \cite{gall2023spatial,le2024peeling} for a theory of a general Markov property in Brownian geometry. We also recall from the introduction that we expect that our self-similar Markov trees (including the critical case discussed in Section \ref{sec:commentsGBP}) are essentially all random decorated trees with a positive decoration on the skeleton satisfying a Markov and self-similar property.

\chapter{Spinal decompositions and bifurcators} \label{chap:spinal:deco}

The main purpose of this chapter is to investigate so-called spinal decompositions for self-similar Markov trees governed by a characteristic quadruplet $ (\sigma^2, \mathrm{a}, \boldsymbol{\Lambda} ; \alpha)$ satisfying Assumption \ref{A:gamma0} which is fixed thorough this chapter. 
 We first take the point of view of general branching processes, and then provide an explicit description 
in the self-similar Markovian setting. We  use either a weighted length measure or the harmonic measure to
mark a point $\rho^\bullet$ at random on the tree, and describe the joint distribution of the decoration-reproduction along the segment $\llbracket \rho, \rho^{\bullet}\rrbracket$ and the family of subtrees that stem from this segment, that is the collection of -- the closure of --  the connected components of the complement $T\backslash \llbracket \rho, \rho^{\bullet}\rrbracket$. Although the decorated trees dangling from $\llbracket \rho, \rho^{\bullet}\rrbracket$ are, conditionally on their initial decorations, independent ssMt, the decoration-reproduction process along $\llbracket \rho, \rho^{\bullet}\rrbracket$ is a different self-similar Markov decoration-reproduction whose characteristics $( \sigma^2, \mathrm{a}_\gamma,  \boldsymbol{ \Lambda}_\gamma ; \alpha)$ are explicitly determined in terms of the initial ones. As a first application,  we address an important issue already mentioned in Sections \ref{sec:1.3} and \ref{sec:defssmb}, namely the fact that different characteristic quadruplets may yield  self-similar Markov trees with the same distribution,
and provide an explicit characterization of such bifurcators  in Section \ref{sec:bifurcators}.
Another application to the determination of  the Hausdorff dimension of self-similar Markov trees that fulfill the first Cramer's condition is given in Section \ref{sec:Hausdim}.

\section{Spine decompositions in the setting of general branching processes}
The notion of spinal decomposition is one of the most useful and powerful tools in the study of branching structures. The purpose of this section is to present its basic aspect, first from the point of view of general branching processes,
and then more specifically in the situation where the decoration-reproduction kernel is Markovian and self-similar. \par Let $(P_x)_{x>0}$ denote a decoration-reproduction kernel; using the notation introduced in Section \ref{sec:2.1}, we write $\P_x$ for the probability law of the associated  family of decoration-reproduction processes $(f_u, \eta_u)_{u \in \mathbb{U}}$, which results when the ancestor has type $x > 0$.
We assume the existence of  a  harmonic function  $h:(0,\infty)\to (0,\infty)$ in the sense of \eqref{E:mean-harmonic} and agree that $h(0)=0$.
Recall that $\chi(u)$ stands for the type of the individual labelled by the vertex $u$ of the Ulam tree $\U$, and that 
the process $M_n=\sum_{|u|=n}h(\chi(u))$ for $n\geq 0$ in \eqref{E:addmart} is a nonnegative martingale.
We stress that  at this stage, we  do not require Property~{\hypersetup{linkcolor=black}\hyperlink{prop:P}{$(\mathcal{P})$}} to hold, nor 
the additive martingale \eqref{E:addmart} to converge in $L^1(\P)$.
The harmonic function $h$ serves here to distinguish an infinite lineage at random, and  we will use the symbol $\star$ as an exponent  to refer to distinction. That is,   a single individual $u^\star(n)\in \N^n$ is distinguished at each generation $n\geq 0$ such that $u^\star(n)\prec u^\star(n+1)$, i.e. the  individual distinguished at generation $n$ is always the parent of the individual distinguished at generation $n+1$, and $(u^\star(n))_{n\geq 0}$ constitutes the distinguished lineage. Specifically, we introduce a probability measure $\bar \P_x^h$ 
that describes the joint law of the family of decoration-reproduction processes $\left( f_u, \eta_u \right)_{u\in \U}$  together with  $(u^\star(n))_{n\geq 0}$.
For every $n\geq 0$, every $v\in \N^n$ and every 
nonnegative functional $\Phi$ of the decoration-reproduction processes up to generation $n$, we define
  \begin{eqnarray} \label{eq:defspineh}\bar \E_x^h \big(\Phi( (f_u, \eta_u)_{|u|\leq n}) \indset{u^\star(n)=v}\big) \coloneqq h(x)^{-1}  \E_x\Big( \Phi((f_u, \eta_u)_{|u|\leq n}) h(\chi(v))
\Big).  \end{eqnarray}

The martingale property ensures the coherence of this definition. It is an easily checked and well-known fact in the field of branching processes (see for instance \cite[Chapter 4]{shi2015branching} or \cite{biggins2004measure}) that the probability measures $\bar \P_x^h$ also describe 
a generalized branching process, where individuals not only have a type $x>0$, but are furthermore either distinguished or not. More precisely, non-distinguished individuals beget only non-distinguished children and the statistics of their decoration-reproduction processes are given by the original kernel $(P_x)_{x>0}$.
Distinguished individuals always give birth to a single distinguished child, possibly together with further non-distinguished children,  the distinguished child being picked at random  in the progeny with probability proportional to the value assigned by $h$ to the type. The decoration-reproduction process of a distinguished individual follows a biased version of the original kernel, denoted by  $(\bar P^h_x)_{x>0}$.
Specifically, writing  $t^\star$ for the age of the parent when its distinguished child is born,  and $ x^\star$ for  the type of the latter, then for every nonnegative functional $\Psi$, there is the identity
\begin{equation}\label{eq:tilteddr}\bar E_x^h\big(\Psi(f,\eta,t^\star, x^\star)\big) = h(x)^{-1} E_x\Big( \int_{[0,\infty)\times(0,\infty)} \eta(\dd t, \dd y)\Psi(f, \eta,t,y) h(y)\Big)  .
\end{equation}

Let us now provide a more explicit description of the biased law $\bar P^h_x$ 
in the case when the original decoration-reproduction process is self-similar  Markovian with characteristic quadruplet $(\sigma^2, \mathrm{a}, \boldsymbol{\Lambda} ; \alpha)$ as in Section \ref{sec:2.2}. We exclude implicitly the degenerate case $\eta=0$ when there is no reproduction at all, and suppose further that the cumulant function $\kappa$ in \eqref{E:cumulant} is well-defined and vanishes at some  $\omega>0$. 
The function $h(x)=x^\omega$ is harmonic; see Lemma \ref{L:MUI} and observe that the sole harmonicity of $h$ only requires Lemma~\ref{L:verCMJ1} and not the more stringent Assumption~\ref{A:omega-}. 
In this direction, we shall now write $\bar P^\omega$ instead of $\bar P^h$ for the sake of clarity.
Recall from Section~\ref{sec:2.2} that in the self-similar Markov setting, the decoration process is denoted by $X$ rather than by $f$, so
\eqref{eq:tilteddr} reads
\begin{equation}\label{eq:tilteddr2}\bar E_x^\omega\big(\Psi(X,\eta,t^\star, x^\star)\big) = x^{-\omega } E_x\Big( \int_{[0,\infty)\times(0,\infty)} \eta(\dd t, \dd y)\Psi(X, \eta,t,y) y^\omega \Big).
\end{equation}

We now argue that a similar feature occurs also for $\gamma>0$ when $\kappa(\gamma)<0$. 
As a motivation,  we define first a probability measure $\widetilde{\P}_x^\gamma$, with associated mathematical expectation $\widetilde{\E}_x^\gamma$,  describing the joint law of the family of decoration-reproduction processes together with a marked vertex $u^\bullet$ and a marked time $t^\bullet$, 
 such that for every $v\in \U$, every measurable $\varphi: \R_+\to \R_+$  and every nonnegative functional $\Phi$,
\begin{equation} \label{E:spinedecrep}
\widetilde{\E}_x^\gamma \Big(\Phi( (f_u, \eta_u)_{u\in \U}) \varphi(t^\bullet) \indset{u^\bullet=v}\Big) \coloneqq -\kappa(\gamma) x^{-\gamma}  \E_x\Big( \Phi((f_u, \eta_u)_{u\in \U}) \int_0^{z_v} \varphi(t) f_v(t)^{\gamma-\alpha} \dd t\Big),
\end{equation}
where the assertion that $\widetilde{\P}_x^\gamma$ is probability is seen from  Proposition \ref{prop:lengthmeasures}. 
In words, the distribution of the family of decoration-reproduction processes $ (f_u, \eta_u)_{u\in \U}$ under $\widetilde{\P}_x^\gamma$ is $x^{-\gamma} \kappa(\gamma) \uplambda^\gamma (T)\cdot \P_x$, 
further conditionally on this family, the marked vertex $u^\bullet$ is picked at random in $\U$ with probability proportional to $ \int_0^{z_v}  f_v(t)^{\gamma-\alpha} \dd t$, and finally, 
the conditional law of the marked time $t^\bullet$ given $(f_u, \eta_u)_{u\in \U}$ and $u^\bullet$ is  proportional to $\mathbf{1}_{t\in [0,z_{u^\bullet})}f_{u^\bullet}(t)^{\gamma-\alpha} \dd t$.

 We then relate the law $\widetilde{\P}^\gamma_x$ to that of another generalized branching process where  individuals have a positive type and are either distinguished or not. Again non-distinguished individuals beget only non-distinguished children and the statistics of their decoration-reproduction processes are given by the original kernel $(P_x)_{x>0}$.
In turn distinguished individuals, for which an exponent $\star$ is used in the notation,  give birth to at most one distinguished child and further non-distinguished children. We stress that now, a distinguished individual may have no distinguished child (then, of course, the distinguished lineage dies out), in which case we rather distinguish an age. 
Specifically, when the distinguished parent has a distinguished child, we write $x^\star>0$ for the type of the latter and $t^\star$ for the age of the parent at the distinguished birth event.
Otherwise, i.e. when the distinguished parent has no distinguished child, we set $ x^\star=0$ and write $t^\star$ for the distinguished age. 
With this notation at hand, the decoration-reproduction kernel  $(\bar P^\gamma_x)_{x>0}$ for distinguished individuals is defined in terms of the original kernel $(P_x)_{x>0}$ by the following  variation of \eqref{eq:tilteddr}.
We set for every nonnegative functional $\Psi$, 
\begin{align}\label{eq:tilteddr3}
&\bar E_x^\gamma\big(\Psi(X,\eta,t^\star, x^\star)\big) \nonumber \\
&= x^{-\gamma} E_x\left( \int_{[0,z)\times(0,\infty)} \Psi(X, \eta,t,y) y^\gamma \eta(\dd t, \dd y)-\kappa(\gamma)\int_0^z \Psi(X, \eta,t,0) X(t)^{\gamma-\alpha}\dd t \right) ,
\end{align}
where as usual $\bar E_x^\gamma$ stands for the mathematical expectation with respect to $\bar P^\gamma_x$. 
We write $\bar{\P}^\gamma_x$ for the law of the generalized branching process  with distinguished individuals defined above , when the ancestor has type $x>0$ and is  distinguished. More precisely, $\bar{\P}^\gamma_x$  is the joint distribution of the family of decoration-reproduction processes $\left( f_u, \eta_u\right)_{u\in \U}$ together with
the (finite) sequence of distinguished individuals and distinguished ages. 
We further write  $u^\bullet$ for the ultimate distinguished individual and $t^\bullet$ for its distinguished age. 

\begin{lemma} \label{L:spinedecgbp} Suppose $\kappa(\gamma)<0$. Then the law of $\left( ( f_u, \eta_u)_{u\in \U}, t^\bullet, u^\bullet\right)$ under $\bar \P^\gamma_x$ is $\widetilde{\P}^\gamma_x$. 
\end{lemma} 
\begin{proof} To start with, we observe that for any $n\geq 1$, any vertex $v\in \N^n$ at generation $n$, and every functional $F\geq 0$, there is the identity
\begin{equation} \label{eq:jesaispasquoi}
\bar{\E}^\gamma_x\big(F\left((f_u, \eta_u)_{|u|\leq n-1}\big)\indset{v \text{ is distinguished}}\right)= \widetilde{\E}^\gamma_x\left(F\left((f_u, \eta_u)_{|u|\leq n-1}\right)\indset{ u^\bullet(n)=v}\right),
\end{equation}
where as usual $u^\bullet(n)$ denotes the forebear of $u^\bullet$ at generation $n$ if $|u^\bullet|\geq n$ and otherwise the indicator in the left-hand side is interpreted as $0$. 
Indeed, this identity is immediately checked for $n=1$ from \eqref{E:spinedecrep} and Proposition \ref{prop:lengthmeasures}. The general case  $n\geq 1$ then follows by induction, applying the branching property under $\bar \P^\gamma_x$ and under $\P_x$. 

Next, consider a measurable function $\varphi: \R_+\to \R_+$  and a nonnegative functional $\Phi$.
By the branching property under $\P_x$, there is  a functional $\Psi(f,\eta)$ of the decoration-reproduction process   such that for any $x>0$, 
$$\E_x \big(\Phi( (f_u, \eta_u)_{u\in \U}) \mid f_\varnothing, \eta_{\varnothing}\big) =\Psi(f_\varnothing, \eta_{\varnothing}).$$
Using the branching property under $\bar \P^\gamma_x$ and \eqref{eq:tilteddr3}, we get by conditioning on the decoration-reproduction of the ancestor that
\begin{align} \label{eq:jesaispasquoi2}
\bar \E_x^\gamma \left(\Phi( (f_u, \eta_u)_{u\in \U}) \varphi(t^\bullet) \indset{u^\bullet=\varnothing}\right) &= -\kappa(\gamma) x^{-\gamma}  \E_x\left( \Psi(f_\varnothing, \eta_{\varnothing}) \int_0^{z_\varnothing} \varphi(t) f_\varnothing(t)^{\gamma-\alpha} \dd t\right)\nonumber \\
& = \widetilde{\E}_x^\gamma \left(\Phi( (f_u, \eta_u)_{u\in \U}) \varphi(t^\bullet) \indset{u^\bullet=\varnothing}\right),
\end{align}
where the second equality stems from \eqref{E:spinedecrep}. The statement now readily follows from the combination of  \eqref{eq:tilteddr3} and \eqref{eq:jesaispasquoi}  with the branching property. 
\end{proof}

Our goal now is to characterize the law of $(X,\eta,t^\star, x^\star)$ under $(\bar{P}^{\gamma}_{x})_{x\geq 0}$ in both regimes.  In this direction, we fix $\gamma> 0$  such that $\kappa(\gamma)\leq 0$, and we decompose the decoration-reproduction process into three parts. The first corresponds to the time interval $[0,t^\star)$, the second to the time interval $(t^\star,\infty)$ up to a natural shift and rescaling, and the third focusses at the exact time $t^\star$ when the distinguished child is born. We now spell out these parts explicitly.
We write $ A= A(X,\eta,t^\star)$ for the pair resulting from the restriction of $(X,\eta)$ to the time-interval $[0,t^\star)$, that is
$$ A = \Big( \indset{[0,t^\star)} X, \indset{[0,t^\star)\times(0,\infty)}\cdot \eta \Big).$$
We write $ B= B(X,\eta,t^\star)$ for the pair given by the decoration-reproduction shifted at time $t^\star$ and properly rescaled, namely the image of
$ \left( \indset{[t^\star,\infty)} X, \indset{(t^\star,\infty)\times(0,\infty)}\cdot \eta \right)$ by the  map 
$$(s,y) \mapsto \Big(\frac{s-t^\star}{X(t^{\star})^{\alpha}}, \frac{y}{X(t^{\star})}\Big), \qquad s\geq t^\star, y\geq 0.$$
Finally, we consider the (possibly finite) sequence $(x_j)_{j\geq 1}$ of the types of  children (distinguished or not) which are born at time $t^\star$, as usual ranked in the non-increasing order.
In other words, the restriction of the reproduction process $\eta$ to the fiber $\{t^\star\} \times (0,\infty)$  
is 
$$\indset{\{t^\star\} \times (0,\infty)} \cdot \eta = \sum_{j\geq 1} \delta_{(t^\star,x_j)},$$
and $x^\star$ is one of the terms of the sequence $(x_j)_{j\geq 1}$.
We set
$$ C= C(X,\eta,t^\star, x^\star)= \Big( \frac{X(t^\star)}{X(t^\star-)}, \frac{x^\star}{X(t^\star-)}, \Big( \frac{x_j}{X(t^\star-)}\Big)_{j\geq 1}\Big).$$
Plainly, the variable $X(t^\star-)$ is measurable with respect to $  A$, so $X(t^\star)$ can be recovered from $ A$ and $ C$. Therefore  $ A$, $  B$, and $ C$ entirely determine the decoration-reproduction process of a distinguished individual.
   We conclude this section by an explicit description of their joint law. In this direction, recall that the function $\psi$ is the Laplace exponent given by the L\'evy-Khintchine formula \eqref{E:LKfor}, and from \eqref{E:cumulant}
 that
 $$\psi(\gamma) - \kappa(\gamma) = - \int_{\mathcal{S}} \boldsymbol{\Lambda}( \d y, \d  \mathbf y ) \left(\sum_{i=1}^{\infty} \e^{\gamma y_i} \right)\in(-\infty, 0).$$

\begin{proposition}\label{P:spinegbpgamma} Let $\gamma > 0$ such that $\kappa(\gamma)\leq 0$. 
 For every $x>0$,  the random variables $ A$, $ B$ and $ C$  are independent under $\bar P^\gamma_x$. 
Moreover we have:
\begin{itemize}
\item[(i)]
$ A$ is a self-similar Markov decoration-reproduction process with characteristic quadruplet $(\sigma^2, \bar{\mathrm{a}}_\gamma, \bar{\boldsymbol{\Lambda}}_\gamma; \alpha)$ and type $x$, where
$$ \bar{\boldsymbol{\Lambda}}_\gamma(\dd y, \dd \mathbf y)\coloneqq \e^{\gamma y}\cdot  \boldsymbol{\Lambda}(\dd y, \dd \mathbf y) - \psi(\gamma) \delta_{(-\infty, (-\infty, \ldots))}(\dd y, \dd \mathbf y),$$
and 
$$ \bar{\mathrm a}_\gamma \coloneqq {\mathrm a} +\sigma^2 \gamma+\int_{\mathcal{S}} y(\e^{  \gamma y} -1) \mathbf{1}_{|y|\leq 1} \boldsymbol{\Lambda}(\dd y, \dd \mathbf{y}).$$
\item [(ii)] $ B$ has the law $P_1$ of the initial self-similar Markov decoration-reproduction process with  characteristic quadruplet $(\sigma^2, \mathrm{a}, \boldsymbol{\Lambda} ; \alpha)$.
\item[(iii)] The law of $C$ is determined by
$$\bar E^\gamma_x\left(F(C)\right) =  |\psi(\gamma)|^{-1} \left(  - \kappa(\gamma) F(1,0,(0,...)) +  \int_{\mathcal S} \boldsymbol{\Lambda}(\dd y, \dd \mathbf{y})  \sum_{i\geq 1} \e^{\gamma y_i} F(\e^{\gamma y},\e^{\gamma y_i},(\e^{\gamma y_j})_{j\geq 1}) \right),$$
for every nonnegative functional $F$.
\end{itemize}
\end{proposition}

\begin{proof} By self-similarity, there is no loss of generality in assuming that the type of the initial distinguished individual is $x=1$. As usual, we then drop the indices $1$ in the notation for probabilities and expectations. Roughly speaking, the cornerstone of the argument consists in applying  the well-known Mecke equation for Poisson random measures and analyzing its consequences. 
For this, we need first to reformulate \eqref{eq:tilteddr2} in terms of the L\'evy process $\xi$ and the Poisson random measure $\mathbf{N}$, which essentially amounts to undoing the Lamperti transformation.

The decoration-reproduction process $(X,\eta)$ under $P$ is constructed by applying the Lamperti transformation to the L\'evy process $\xi$ and the point process $ \eta$ (see  \eqref{E:Lamperproc} and \eqref{Eq:eta:Markov}), and since $ \eta$ is defined in terms of the L\'evy process $\xi$ and the Poisson random measure $\mathbf{N}$ by \eqref{Eq:tildeeta:Markov}, we have
\begin{align*}
\bar E^\gamma\Big(\Phi(\xi,\mathbf N,\epsilon(t^\star), x^\star)\Big) = &E\Big( \int_{[0,\infty)\times \R\times \mathcal S_1} \mathbf{N}(\dd t, \dd y, \dd \mathbf y) \sum_{j\geq 1} \Phi\left(\xi, \mathbf N,t,\e^{\xi(t-)+y_j}\right ) \e^{\gamma(\xi(t-)+y_j)} \Big)\\
 &-\kappa(\gamma) E\Big( \int_{[0,\infty)} \mathrm{d} t ~\Phi\big(\xi, \mathbf N,t,0\big ) \e^{\gamma \xi(t)} \Big),
\end{align*}
where $\Phi$ denotes another generic nonnegative functional and  $\epsilon$ is given by \eqref{E:espilonlamperti}.

 We can now apply the Mecke equation (see for instance \cite[Section 4.1]{last2017lectures}), which involves adding a Dirac mass at $(t,y,\mathbf y)$, under $\mathbf{1}_{t>0} \d t  \boldsymbol{\Lambda}( \dd y, \dd \mathbf y) $, 
to the Poisson random measure $\mathbf N$. Note from the L\'evy-It\^{o} decomposition \eqref{E:LevyIto} that this addition also changes $\xi$ into $\xi + y\indset{[t,\infty)}$. We rewrite the first expectation of the right-hand side above as
\begin{align*}
&\int_0^{\infty} \dd t \int_{\R\times \mathcal S_1} \boldsymbol{\Lambda}( \dd y, \dd \mathbf y) E\Big(  \sum_{j\geq 1} \Phi\left(\xi + y\indset{[t,\infty)}, \mathbf N + \delta_{(t,y,\mathbf y)},t,\e^{\xi(t-)+y_j}  \right)\e^{\gamma(\xi(t-)+y_j)}\Big) \\
&=  \int_{\R\times \mathcal S_1} \boldsymbol{\Lambda}( \dd y, \dd \mathbf y)  \sum_{j\geq 1} \e^{\gamma y_j} \int_0^{\infty}  E\Big( \e^{\gamma \xi(t-)} \Phi\Big(\xi + y\indset{[t,\infty)}, \mathbf N + \delta_{(t,y,\mathbf y)},t,\e^{\xi(t-)+y_j} \Big) \Big) \dd t.
\end{align*}

Recall that variables $A$, $B$, and $C$ in the statement essentially correspond to splitting the decoration-reproduction process of a distinguished individual  at time $t^\star$. Aiming similarly at undoing the Lamperti transformation for the parts before and after $t^\star$, this leads us to decompose $(\xi, \mathbf N)$ under $\bar P^\gamma$ into three parts corresponding to times strictly before, strictly after, and exactly at $\epsilon(t^\star)$. 
Specifically, we write $A'$ for the restriction of $(\xi, \mathbf N)$ to $[0, \epsilon(t^\star))$, and $B'$ for the pair 
$$\Big(\indset{(\epsilon(t^\star),\infty)}(\xi-\xi(\epsilon(t^\star)), \indset{(\epsilon(t^\star),\infty)\times \mathcal S}\cdot \mathbf N\Big)$$ further shifted in time by $\epsilon(t^\star)$. Last, we consider the atom of $\mathbf N$ at time $\epsilon(t^\star)$, say $(\epsilon(t^\star), w, (w_i)_{i\geq 1})$. We denote by $w^\star$  the distinguished element of the sequence 
$(w_i)_{i\geq 1}$ and set $C'=(w,w^\star,(w_i)_{i\geq 1})$. In particular, $w$ coincides with the jump $\Delta \xi(\epsilon(t^\star))$ of $\xi$ at time $\epsilon(t^\star)$ and 
$w^\star = \log x^\star - \xi(\epsilon(t^\star)-)$. Then $A$ and $B$ are recovered by applying the Lamperti transformation to $A'$ and to $B'$ respectively, and $C$ by exponentiating the elements of $C'$.

Specializing the previous discussion to functionals $\Phi$ of the form
$$\Phi(\xi,\mathbf N,\epsilon(t^\star), x^\star)=\Phi_1(A')\Phi_2(B')\Phi_3(C'), $$
and using that for  Lebesgue almost every $t>0$ we have $\xi(t-)=\xi(t)$, we obtain that the quantity $\bar E^\gamma\left( \Phi_1(A')\Phi_2(B')\Phi_3(C') \right)$ equals
\begin{align*}
\Big(&\int_{\R\times \mathcal S_1} \boldsymbol{\Lambda}( \dd y, \dd \mathbf y)  \sum_{j\geq 1} \e^{\gamma y_j}  \Phi_3\big(y,y_j, \mathbf{y} \big)- \kappa(\gamma) \Phi_3\big(0,-\infty,(-\infty,...)\big) \Big)\\
&\cdot \int_0^{\infty}  E\left( \e^{\gamma \xi(t)} \Phi_1(\indset{[0,t)}\xi, \indset{[0,t)\times \mathcal S}\cdot \mathbf N))\Phi_2(\xi^{(t)}, \mathbf{N}^{(t)}) \right) \dd t,
\end{align*}
where $\xi^{(t)}(s):=\xi(t+s)-\xi(t)$, for $s\geq 0$, and $\mathbf{N}^{(t)}(\dd s, \dd y', \dd{\mathbf y}') := \indset{\{s>0\}\times \mathcal S} \cdot \mathbf{N}(t+\dd s, \dd y', \dd{\mathbf y}')$ are the shifted versions of $\xi$ and $\mathbf N$. 
It is now elementary to check that for any $t\geq 0$, the expectation in the right-hand side can be expressed as the product
$$ E\Big( \e^{\gamma \xi(t)} \Phi_1(\indset{[0,t)}\xi, \indset{[0,t)\times \mathcal S}\cdot \mathbf N))\Big) \cdot E\big(\Phi_2(\xi, \mathbf{N}) \big);$$
see e.g. the proof of Proposition \ref{prop:hullmarkov} which is closely related.

This proves the independence claim in the statement, as well as (ii) and (iii), and we are left the computation of
$$\int_0^{\infty}  E\left( \e^{\gamma \xi(t)} \Phi_1(\indset{[0,t)}\xi, \indset{[0,t)\times \mathcal S}\cdot \mathbf N))\right) \dd t.$$
For this, recall that $\psi(\gamma)\in(-\infty,0)$ is the value taken at $\gamma$ by the Laplace exponent $\psi$ of the L\'evy process $\xi$, 
and rewrite the preceding quantity as
$$\int_0^{\infty} \e^{\psi(\gamma)t} E\left( \e^{\gamma \xi(t)-t\psi(\gamma)} \Phi_1(\indset{[0,t)}\xi, \indset{[0,t)\times \mathcal S}\cdot \mathbf N))\right) \dd t.$$
The process $\exp(\gamma \xi(t)-\psi(\gamma)t)$ is a martingale under $P$, which can be use as a density to define  a locally equivalent probability measure. This procedure is known  as an Esscher transform. 
In the current  setting, if we write $(P^\prime_x)_{x>0}$ for the self-similar Markov decoration-reproduction kernel with characteristic quadruplet $(\sigma^2, \bar{\mathrm{a}}_\gamma, \bar{\boldsymbol{\Lambda}}_\gamma ; \alpha)$ defined in the statement,
then an elementary calculation as in the proof of \cite[Theorem 3.9]{kyprianou2014fluctuations} yields
$$|\psi(\gamma)| \int_0^{\infty} \e^{\psi(\gamma)t} E\left( \e^{\gamma \xi(t)-t\psi(\gamma)} \Phi_1(\indset{[0,t)}\xi, \indset{[0,t)\times \mathcal S}\cdot \mathbf N))\right) \dd t= E^\prime \left( \Phi_1(\indset{[0,\zeta)}\xi, \indset{[0,\zeta)\times \mathcal S}\cdot \mathbf N))\right),$$
where we recall that $\zeta$ stands for the lifetime of the L\'evy process, which is exponentially distributed with parameter $|\psi(\gamma)|$ under $P^\prime$, and $E^\prime $ corresponds to the mathematical expectation with respect to $P^\prime$. 
This completes the proof.
\end{proof}

\section{Size-biased spine decompositions for self-similar Markov trees}
\label{Section:spine}
In the setting of  self-similar Markov trees, a spinal decomposition  can be viewed as a remarkable local decomposition  in the sense of Definition \ref{Def:locdec} for some size-biased version of a self-similar Markov tree under $\P_x$. Slightly more precisely, the base subtree that induces this decomposition is merely the segment $\llbracket \rho, \rho^{\bullet}\rrbracket$ from the root $\rho$ to a marked point $\rho^{\bullet}$ which is picked at random according to some natural probability measure on the tree. We may work either with a weighted length measure or the harmonic measure on $T$, and to avoid many repetitions, we will use the same notation for both cases. In this setting, a spinal decomposition further describes explicitly the joint law of the decoration on the segment  $\llbracket \rho, \rho^{\bullet}\rrbracket$
and of the point measure induced by the germs of the decorations of the subtrees dandling from $\llbracket \rho, \rho^{\bullet}\rrbracket$. 

We consider a characteristic quadruplet $ (\sigma^2, \mathrm{a}, \boldsymbol{\Lambda} ; \alpha)$ satisfying Assumption \ref{A:gamma0}.
Either under this sole requirement, we fix 
 $\gamma>0$ such that $\kappa(\gamma)<0$ and write $\upnu\coloneqq -\kappa(\gamma) \uplambda^{\gamma}$, or under the stronger  Assumption \ref{A:omega-},
 we take  $\gamma= \omega_-$, so $\kappa(\gamma)=0$,  and write $\upnu:=\upmu$ for the harmonic measure.  
  Recall  from \eqref{mean:uplambda:gamma:T}  and \eqref{eq:expec:mu:T:x} that  the total mass $\upnu(T)$ has expectation
\begin{equation} \label{eq:tildeproba}
\E_{x}\left( \upnu(T)\right)=  x^{\gamma}, \quad \text{ for } x>0.
\end{equation}

  The characteristic quadruplet $ (\sigma^2, \mathrm{a}, \boldsymbol{\Lambda} ; \alpha)$ induces self-similar laws $(\P_x)_{x>0}$ for families of reproduction processes $\left( f_u, \eta_u\right)_{u\in \U}$, and in turn the latter  yield 
  distributions $(\mathbb{Q}_x)_{x>0}$ on the space $\TT$ of (equivalence classes up to isomorphisms of) measured decorated real trees $\mathbf{T}=(T,d_T,\rho,g,\upnu)$. 
  We next introduce for every $x>0$ a probability measure $\widetilde{\mathbb{Q}}_x$ on the space $\TT^\bullet$ of (equivalence classes up to isomorphisms of) non-measured decorated real trees with a single marked point, $\mathtt{T}^\bullet=(\mathtt{T},\rho^\bullet)$, such that for any positive measurable function $F : \mathbb{T}^{\bullet} \to \mathbb{R}_{+}$,
   \begin{eqnarray} \label{eq:defspinebiaised} \widetilde{\mathbb{Q}}_{x}^{\gamma}\Big(F\big( \mathtt{T}^\bullet\big)\Big) :=   x^{-\gamma} {\mathbb{Q}}_{x} \left( \int_{T} F\big( \mathtt{T}, r\big)  \upnu (\mathrm d r)  \right), \end{eqnarray}
 where for the sake of simplicity, we use the same notation  for the probability measures $\widetilde{\mathbb{Q}}_{x}$ and $ {\mathbb{Q}}_{x}$ as for their corresponding mathematical expectations. In words, $\widetilde{\mathbb{Q}}_{x}^{\gamma}$ is obtained from ${\mathbb{Q}}_{x}$ by first biasing the latter with the total mass  $x^{-\gamma} \upnu(T)$,
and then marking a point $\rho^\bullet$  at random according to the normalized law $\upnu(\dd r)/ \upnu(T)$ on $T$. 
The assertion that $\widetilde{\mathbb{Q}}_{x}^{\gamma}$ is a probability measure on  $\TT^\bullet$ is seen from \eqref{eq:tildeproba}.

\begin{figure}[!h]
 \begin{center}
 \includegraphics[width=12cm]{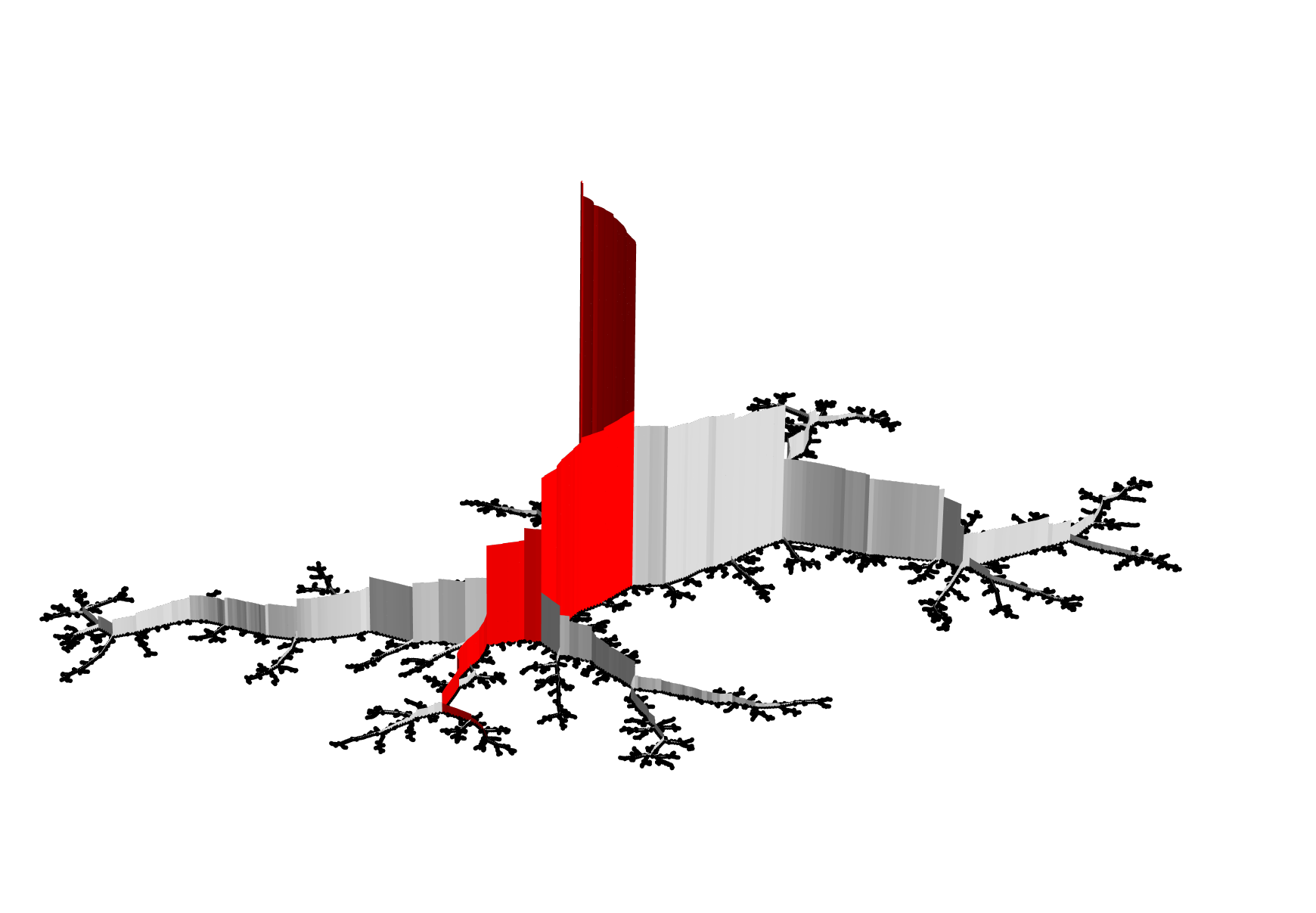}
 \caption{Illustration of  the spinal decomposition where a branch (in red above) has been distinguished and along which the decoration-reproduction process evolves according to tilted characteristics.}
 \end{center}
 \end{figure}
It will often be convenient to work with a specific realization of the marked decorated tree $\mathtt{T}^\bullet=(\mathtt{T},\rho^\bullet)$ under $\widetilde{\mathbb{Q}}_{x}^{\gamma}$
in terms of a certain general branching process. In this direction, let us first make the connexion with  the preceding section in the case when Assumption \ref{A:omega-} holds. When $\kappa(\gamma)<0$ this connection is transparent. Specifically,  if we set $\varrho^\bullet\coloneqq \varrho(u^\bullet, t^\bullet)$, then by Proposition \ref{P:constructionomega-} and Lemma \ref{L:spinedecgbp} , under $\bar{\P}_x^\gamma$, 
 the family of decoration-reproduction processes $(f_u,\eta_u)_{u\in \U}$ satisfies Property~{\hypersetup{linkcolor=black}\hyperlink{prop:P}{$(\mathcal{P})$}}, and 
  the construction of decorated trees by gluing  building blocks immediately yields that
\begin{equation} \label{eq:spinevsspine2}
\text{  the distribution of  $(\mathtt{T}, \varrho^\bullet)$ under $\bar{\P}_x^\gamma$  is $\widetilde{\mathbb Q}^\gamma_x$. }
\end{equation}
Let us now explain why the analog result also holds for the harmonic measure. More precisely, recall that the probability measure  $\bar{\P}_x^\omega$ stands for  the joint law of the family of decoration-reproduction processes $\left( f_u, \eta_u\right)_{u\in \U}$  together with the distinguished lineage  $(u^\star(n))_{n\geq 0}$.

\begin{proposition} \label{P:spinebar} Let Assumption \ref{A:omega-} be satisfied and take $\gamma=\omega_-=\omega$. Then the following assertions hold under $\bar \P^{\omega_-}_x$ for any $x>0$:
\begin{enumerate}
\item[(i)]
The law of $(f_u,\eta_u)_{u\in \U}$    is  $x^{-\omega_-} \upnu(T) \cdot \P_x$ and satisfy 
 Property~{\hypersetup{linkcolor=black}\hyperlink{prop:P}{$(\mathcal{P})$}}. 
We write as usual $\mathtt{T}$ for the resulting  decorated tree.

\item[(ii)] The sequence   $(\varrho(u^\star(n)))_{n\geq 0}$ of the locations on $T$ induced by the births of the distinguished individuals $u^\star(n)$ converges a.s. We write $\varrho^\bullet$ for its limit.

\item[(iii)] The law of the marked decorated tree
$(\mathtt{T}, \varrho^\bullet)$  is $\widetilde{\mathbb Q}^{\omega_-}_x$. 
\end{enumerate}
\end{proposition}

\begin{proof} (i) By definition, for any $n\geq 0$, the distribution of $(f_u,\eta_u)_{|u|\leq n}$  under $\bar \P^{\omega_-}_x$  is the same as under $x^{-\omega_-} M_{n+1}(\omega_-) \cdot \P_x$,
where $(M_n(\omega_-))_{n\geq 0}$ is the intrinsic martingale.  The remaining assertions are
then immediate from Proposition \ref {P:constructionomega-}.

(ii) The sequence   $(\varrho(u^\star(n)))_{n\geq 0}$ is monotone increasing, in the sense that $\varrho(u^\star(n))$ is te parent of $ \varrho(u^\star(n+1))$ for all $n\geq 0$. 
Since $T$ is compact, this sequence converges in $T$.

(iii) Again it is seen from the very construction that for every $n\geq 0$, the conditional distribution under $\bar{\P}_x^{\omega_-}$ of $\varrho(u^\star(n))$ given  $\mathtt{T}$ is $\upmu^n(\dd v) /\upmu(T)$, where $\upmu^n$ stands for the projection of the harmonic measure on $T^n$ (recall that the latter denotes the subtree induced by the individuals up to the generation $n$ only). Since we know from Proposition \ref{P:newPmass} that $\upmu^n$ converges towards $\upmu$ as $n\to \infty$ in the sense of Prokhorov, $\bar{\P}_x^{\omega_-}$-a.s. given $\mathtt{T}$, we conclude that  the conditional law under $\bar{\P}_x^{\omega_-}$ of $\varrho^\bullet$ given  $\mathtt{T}$ is indeed $\upmu(\dd v)/\upmu(T)=\upnu(\dd v)/\upnu(T)$. 
\end{proof}

Now that we have defined the $\upnu$-marked version of a ssMt and explained how it can be constructed using an explicit family of decoration-reproduction processes, let us turn our attention to the spinal decomposition, which involves decomposing the latter along the marked segment. Informally, our goal is to explicitly describe the joint law of the decoration on the segment $\llbracket \rho, \rho^{\bullet}\rrbracket$,  the standard decorated versions of the subtrees dandling from $\llbracket \rho, \rho^{\bullet}\rrbracket$, and the point measure induced by the germs of the decorations of these dandling subtrees. Let us explain the road map that we will follow.
First, we introduce the candidate for the law of the decoration-reproduction process encoding the decoration and germs on the segment $\llbracket \rho, \rho^{\bullet}\rrbracket$
by means of a new characteristic quadruplet. As in Section \ref{sec:1.2}, this new decoration-reproduction process allows us to define a decorated segment, and by analogy with Section \ref{sub:section:gene} we refer to it as the marked spine. Next, we demonstrate in Proposition  \ref{P:spinegeneralnew} that we can construct a decorated tree by gluing self-similar Markov trees, with the original characteristic quadruplet  $ (\sigma^2, \mathrm{a}, \boldsymbol{\Lambda} ; \alpha)$, onto this marked  spine. Finally, we show that the resulting decorated tree, marked at the endpoint of the marked  spine,  is distributed according to $ \widetilde{\mathbb{Q}}_{x}^{\gamma}$, see Theorem \ref{T:Spindec:gamma}. The proof will rely on the characterization given in Proposition~\ref{P:spinegbpgamma}.

Let us proceed by introducing our candidate for the decoration-reproduction process along the marked spine.  In this direction, we introduce  a measure $ \boldsymbol{\Lambda}_\gamma$ on $\mathcal S$ which is derived from the generalized L\'evy measure  $ \boldsymbol{\Lambda}$ as follows. Fix $i\geq 1$ and consider for any pair $(y, \mathbf y)$ in $\mathcal{S}=\R\times \mathcal {S}_1$ the pair $(y, \mathbf y)^{\backsim i}$ in $\mathcal S$ that results by swapping\footnote{Observe that we often use the symbol $\sim$ to indicate a tilting transformation of probability measures; this should not be confused with $\backsim$ which rather refers to swapping two elements.}
$y$ and $y_i$. 
Specifically, the first element of $(y, \mathbf y)^{\backsim i}$ is given by the $i$-th term $y_i$ of the sequence $\mathbf y =(y_j)_{j\geq 1}$,
and its second element  is obtained from the sequence $(y_1, \ldots, y_{i-1}, y, y_{i+1}, \ldots)$  (i.e. we replace the $i$-th term $y_i$ in $\mathbf y $ by $y$) after re-ordering terms in the non-increasing order. 
Then $ \boldsymbol{\Lambda}^{\backsim i}$  is simply the push-forward measure of  $ \boldsymbol{\Lambda}$ by the map $(y, \mathbf y)\mapsto (y, \mathbf y)^{\backsim i}$. 
We next define the measure $\boldsymbol{\Lambda}^\backsim_\gamma$
 on  $\mathcal{S}$ by
$$ \boldsymbol{\Lambda}^{\backsim}_{\gamma}(\dd y, \dd \mathbf{y})\coloneqq
  \e^{\gamma y}\cdot \left( \sum \limits_{i\geq 1}  \boldsymbol{\Lambda}^{\backsim i}(\dd y, \dd \mathbf{y})\right),$$
 and observe from the finiteness of $\kappa(\gamma)$  and \eqref{E:cumulant} that
 \begin{equation} \label{E:massbacksim} \boldsymbol{\Lambda}^{\backsim}_\gamma({\mathcal S} )
 =   \int_{\mathcal S}  \left(  \sum_{i\geq 1} \e^{\gamma y_i}\right) \boldsymbol{\Lambda}(\dd y, \dd \mathbf{y})= \kappa(\gamma)-\psi(\gamma)<\infty.
 \end{equation}
Finally, we define
\begin{equation}\label{Lambda^*_def_2:gamma}
 \boldsymbol{\Lambda}_\gamma(\dd y, \dd \mathbf{y}):=\e^{\gamma y}\cdot  \boldsymbol{\Lambda}(\dd y, \dd \mathbf{y})+\ \boldsymbol{\Lambda}^{\backsim}_{\gamma}(\dd y, \dd \mathbf{y})-\kappa(\gamma) \delta_{(-\infty, (0, -\infty,\cdots))}(\dd y, \dd \mathbf{y}) .
 \end{equation}

  We stress from the finiteness of $\psi(\gamma)$ that $\int_{\mathcal S}   (1\wedge y^2 ) \e^{\gamma y}\boldsymbol{\Lambda}(\dd y, \dd \mathbf{y})<\infty, 
$ 
and therefore  $\boldsymbol{\Lambda}_\gamma$ is a generalized L\'evy measure, since
$$
\int_{\mathcal S}   (1\wedge y^2 ) \boldsymbol{\Lambda}_\gamma(\dd y, \dd \mathbf{y})<\infty.
 $$
Of course \eqref{Lambda^*_def_2:gamma} is reminiscent of the definition of  $\bar{\boldsymbol{\Lambda}}_\gamma$ given in Proposition \ref{P:spinegbpgamma}.   One of our motivations for introducing  $ \boldsymbol{\Lambda}_\gamma$ stems from the following observation. 

  \begin{lemma} \label{L:LKtilde}   The function 
 \begin{equation} \label{eq:kappaspi}\psi_{\gamma}(q) \coloneqq \kappa(\gamma+q), \qquad q\geq 0 \end{equation}
can then be expressed in the L\'evy-Khintchine form
 $$\psi_{\gamma}(q)= \frac{1}{2} \sigma^2 q^2 + \mathrm{a}_{\gamma}q + \int_{\mathcal S}   \left( \mathrm{e}^{q y} -1-  q y \mathbf{1}_{|y|\leq 1} \right) \boldsymbol{\Lambda}_\gamma(\dd y, \dd \mathbf{y}),$$
where  the drift coefficient $\mathrm{a}_{\gamma}$ is given by
  \begin{equation} \label{eq:driftspine}{\mathrm a}_\gamma \coloneqq {\mathrm a} +\sigma^2 \gamma+\int_{\mathcal{S}}\Big( y(\e^{  \gamma y} -1) \mathbf{1}_{|y|\leq 1}+\sum_{i=1}^{\infty} y_i\e^{ \gamma y_i} \mathbf{1}_{|y_i|\leq1}\Big) \boldsymbol{\Lambda}(\dd y, \dd \mathbf{y}),
  \end{equation}
\end{lemma}

\begin{proof} 
The  claim should be viewed as a variation of a well-known result  related to the Esscher transform, see e.g. \cite[Theorem 3.9]{kyprianou2014fluctuations}. It states that if $\psi$ is the Laplace exponent of a real L\'evy process with characteristic triplet  $(\sigma^2, \mathrm a, \Lambda)$ (recall our convention that the killing rate $\mathrm k=\Lambda(\{-\infty\})$ is specified by the mass of the L\'evy measure at $-\infty$), and if $\psi(\gamma)\leq 0$, then 
the shifted function $\psi(\gamma+\cdot)$ can be expressed in the L\'evy-Khintchine form \eqref{E:LKfor} for the same Gaussian coefficient $\sigma^2$, the tilted drift coefficient 
$\mathrm a+\sigma^2 \gamma + \int(\e^{\gamma y}-1)y \indset{|y|\leq 1}\Lambda(\dd y)$, and the  tilted L\'evy measure $\e^{\gamma y} \Lambda(\dd y)-\psi(\gamma)\delta_{-\infty}$.

Let us now proceed with the proof of the lemma. Since   $\boldsymbol{\Lambda}_\gamma$ is a generalized L\'evy measure with killing rate $-\kappa(\gamma)\geq 0$, we have 
\begin{align*}
&\int_{\mathcal S}   \left( \mathrm{e}^{q y} -1-  q y \mathbf{1}_{|y|\leq 1} \right) \boldsymbol{\Lambda}_\gamma(\dd y, \dd \mathbf{y})\\
&= \int_{\mathcal{S}}   \left( \mathrm{e}^{q y} -1-  q y \mathbf{1}_{|y|\leq 1} \right) \e^{\gamma y}  \boldsymbol{\Lambda}(\dd y, \dd \mathbf{y}) + 
\int_{\mathcal{S}}\Big( \sum_{i=1}^{\infty} \e^{ \gamma y_i}\left( \mathrm{e}^{q y_i} -1-  q y_i \mathbf{1}_{|y_i|\leq 1} \right)\Big) \boldsymbol{\Lambda}(\dd y, \dd \mathbf{y})
+\kappa(\gamma).
\end{align*}
Now, by the Esscher transformation, the first term in the sum on the right-hand side can be expressed as 
\begin{align*}
& \int_{\mathcal{S}}  \left( \mathrm{e}^{q y} -1-  q y \mathbf{1}_{|y|\leq 1} \right) \e^{\gamma y}  \boldsymbol{\Lambda}(\dd y, \dd \mathbf{y})  \\
&=\psi(\gamma+q) - \psi(\gamma) -\frac{1}{2}\sigma^2 q^2 - q \left( \mathrm a+\sigma^2 \gamma + \int_{\mathcal{S}} (\e^{\gamma y}-1)y  \mathbf{1}_{|y|\leq 1} \boldsymbol{\Lambda}(\dd y, \dd \mathbf{y})  \right).
\end{align*}
Recall also from \eqref{E:massbacksim} that
$\kappa(\gamma) =\psi(\gamma) +\boldsymbol{\Lambda}^{\backsim}_\gamma({\mathcal S} ) $. 
Therefore, if we set 
  \begin{equation*}{\mathrm a}_\gamma \coloneqq {\mathrm a} +\sigma^2 \gamma+\int_{\mathcal{S}}\Big( y(\e^{  \gamma y} -1) \mathbf{1}_{|y|\leq 1}+\sum_{i=1}^{\infty} y_i\e^{ \gamma y_i} \mathbf{1}_{|y_i|\leq1}\Big) \boldsymbol{\Lambda}(\dd y, \dd \mathbf{y}),
  \end{equation*}
and use  \eqref{E:cumulant}, then we arrive at the identity
 $$\kappa(\gamma+q)= \frac{1}{2} \sigma^2 q^2 + \mathrm{a}_{\gamma}q + \int_{\mathcal S}   \left( \mathrm{e}^{q y} -1-  q y \mathbf{1}_{|y|\leq 1} \right) \boldsymbol{\Lambda}_\gamma(\dd y, \dd \mathbf{y}),$$
 and our claim is checked. 
 \end{proof}

 Lemma \ref{L:LKtilde} allows to apply the construction devised in Section \ref{sec:2.2} by keeping the same Gaussian coefficient $\sigma^2$ and the same exponent of self-similarity $\alpha$, but  
 replacing the Poisson random measure $\mathbf N$  on $[0,\infty)\times  \mathcal{S}$  with intensity $ \dd  t \boldsymbol{\Lambda}(\dd y, \dd \mathbf{y})$ there by  a Poisson random measure   $\mathbf N_{\gamma}$ with intensity $ \dd  t \boldsymbol{\Lambda}_{\gamma}( \dd y ,\dd \mathbf{y})$ 
 and the drift coefficient $\mathrm a$ by $\mathrm a_{\gamma}$ given in \eqref{eq:driftspine}. We write $(P^{\gamma}_x)_{x>0}$ for the self-similar  kernel of decoration-reproduction laws induced there by the characteristic quadruplet 
 $( \sigma^2,\mathrm{a}_{\gamma},   \boldsymbol{\Lambda}_{\gamma};  \alpha)$. In particular, the decoration process under this kernel is a self-similar Markov process $X_\gamma$ with exponent $\alpha$ that is associated by the Lamperti transformation to a L\'evy process $\xi_{\gamma}$ with Laplace exponent $\psi_\gamma$. 
 We also write $\eta_\gamma$ for the reproduction process. We stress that in the case where $\kappa(\gamma)<0$, the decoration is strictly positif immediately before the deathtime, $X_\gamma(z-)>0$, and the reproduction process $\eta_\gamma$ has an atom at $(z,X_\gamma(z-))$ (and there are no further atoms at time $z$). 
 
Our goal now is to glue self-similar Markov trees, with characteristics $(\sigma^2, \mathrm{a}, \boldsymbol{\Lambda}; \alpha)$, onto the decorated segment induced by $(X_\gamma, \eta_\gamma)$. To this end, we rely on the following technical result.

   \begin{proposition}\label{P:spinegeneralnew} Fix $x>0$. Under $P_x^\gamma$, consider $(t_1,y_1), (t_2,y_2), \ldots$ the atoms of $\eta_\gamma$ in co-lexicographical order.
   Let $(\normalfont{\texttt T}_j=(T_j, d_{T_j}, \rho_j, g_j))_{j\geq 1}$ denote a sequence of independent ssMt with
 characteristic quadruplet $(\sigma^2, \mathrm{a}, \boldsymbol{\Lambda} ; \alpha)$ 
  and such that each $\normalfont{\texttt T}_j$ has the law $\P_{y_j}$.  Then, almost surely, the two families 
 $$\big( \mathrm{Height}(T_j)\big)_{j\geq 1} \ \text{and} \ \big( \max_{T_j} g_j\big)_{j\geq 1}$$
 are  null.
    \end{proposition}
Before proving the proposition, let us discuss some of its implications. First, with the notation above, it allows us to use Lemma \ref{L:newL2} and  consider the non-degenerated  decorated tree    
$$\mathrm{Gluing}\Big( \big([0, z], d,0,X_\gamma \big), (y_j)_{j\geq 1}, \big(\texttt{T}_i\big)_{j\geq 1}\Big), $$
where  $z$ stands for the lifetime of $X_\gamma$ and    $d$ for the usual distance on segments. We  then let $ \widehat{\mathbb{Q}}_{x}^{\gamma}$ be the resulting law on the space  $\TT^\bullet$ of decorated real trees with a single marked point, where the marked point is induced by the right-extremity $z$ of the ancestral segment. 
Heuristically speaking, under $ \widehat{\mathbb{Q}}_{x}^{\gamma}$,  the evolution along the spine is governed by the characteristic $(  \sigma^{2}, \mathrm{a}_{\gamma}, \boldsymbol{ \Lambda_{\gamma}}; \alpha)$,  while the other branches are governed by $(  \sigma^{2}, \mathrm{a}, \boldsymbol{ \Lambda} ; \alpha)$. We can now provide a formal statement for the spinal decomposition. Recall the definition \eqref{eq:defspinebiaised} of the law  $ \widetilde{ \mathbb{Q}}_{x}^{\gamma}$.

\begin{theorem} [Spinal decomposition] \label{T:Spindec:gamma} Either let Assumption \ref{A:gamma0} hold and pick any 
 $\gamma>0$ such that $\kappa(\gamma)<0$, or let the stronger  Assumption \ref{A:omega-} hold and take  $\gamma= \omega_-$. 
 For every $x>0$,  the probability measures  $ \widehat{\mathbb{Q}}_{x}^{\gamma}$ and  $ \widetilde{ \mathbb{Q}}_{x}^{\gamma}$  on $\TT^\bullet$ are identical.
\end{theorem}
\begin{figure}[!h]
 \begin{center}
 \includegraphics[width=12cm,angle=1]{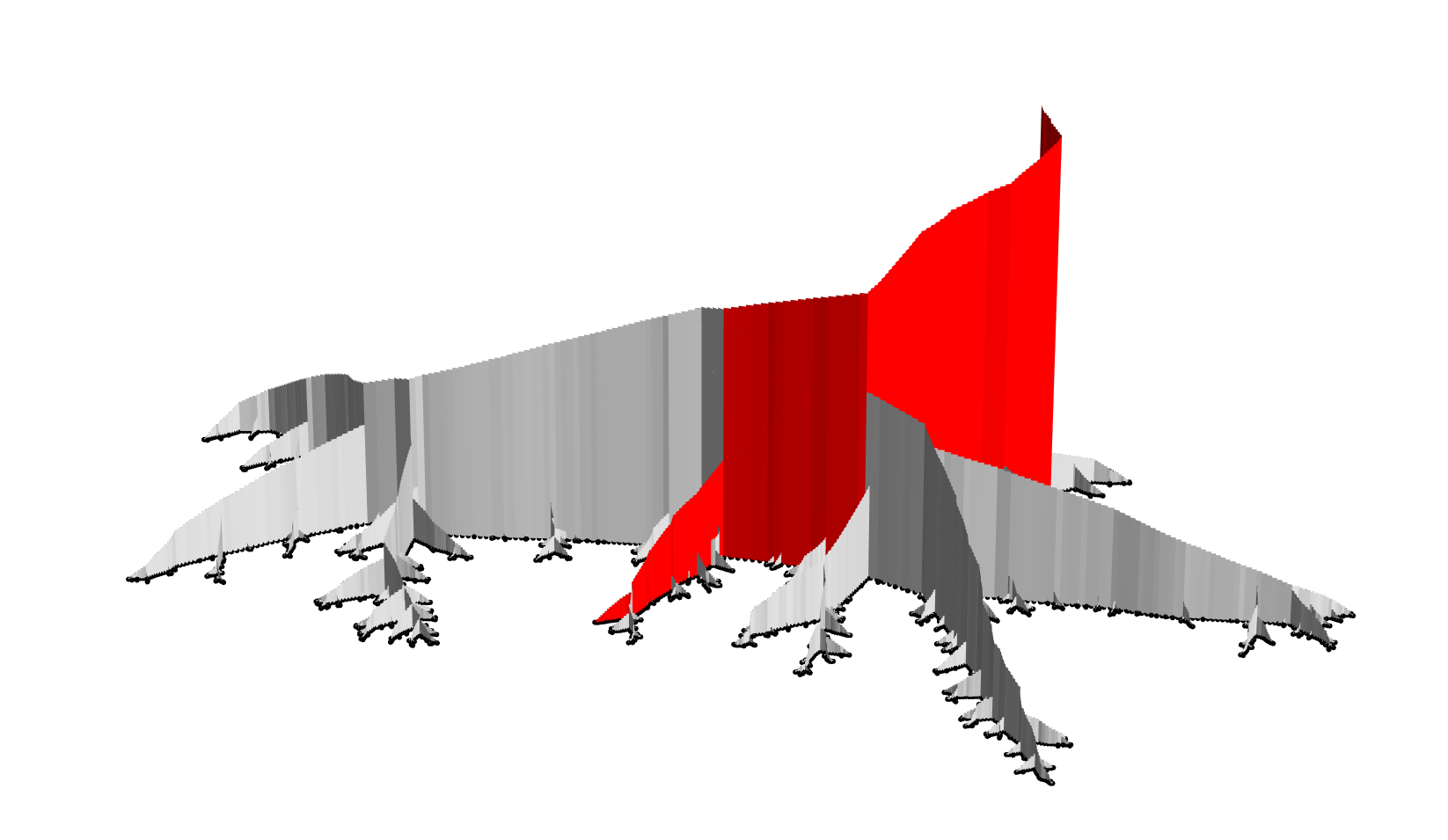}
 \caption{Illustration of Theorem \ref{T:Spindec:gamma}.}
 \end{center}
 \end{figure}
 
 The proofs of Proposition  \ref{P:spinegeneralnew} and Theorem \ref{T:Spindec:gamma} rely on  two further technical lemmas that connect the distribution of the decoration-reproduction process $(X_\gamma,\eta_\gamma)$ under $P_1^\gamma$ to $(X,\eta)$ under $P_1$. These two lemmas will allow us to prove Proposition \ref{P:spinegeneralnew} and Theorem \ref{T:Spindec:gamma}   simultaneously by comparing these results with the characterization of $\bar{P}_1^\gamma$ given in Proposition \ref{P:spinegbpgamma}.  To start with,  \eqref{Lambda^*_def_2:gamma} suggests 
to decompose $\mathbf{N}_\gamma$ 
under $P^\gamma_1$ as the sum of two independent Poisson point processes, $\mathbf{N}^{\prime}$ and $\mathbf{N}^{\prime\prime}$, with respective intensities
$$\e^{\gamma y}\cdot \dd t \boldsymbol{\Lambda}(\dd y, \dd \mathbf{y})\quad \text{and}\quad  \dd t \left(\boldsymbol{\Lambda}^{\backsim}_\gamma-\kappa(\gamma) \delta_{(-\infty, (0, -\infty,\cdots))}\right) (\dd y, \dd \mathbf{y}).$$
 Since  $\boldsymbol{\Lambda}^{\backsim}_\gamma$ is a finite measure, the set of times at which $\mathbf{N}^{\prime\prime}$ has an atom is discrete $P^\gamma_1$-a.s.,
 and we write $r^{\prime\prime}$  for the first one, 
 $$r^{\prime\prime}\coloneqq \sup\big\{r\geq 0: \mathbf{N}^{\prime\prime}([0,r]\times \mathcal S)=0\big\}.$$
 In particular,  by \eqref{E:massbacksim}, $r^{\prime\prime}$ is exponentially distributed with parameter 
 $\boldsymbol{\Lambda}^{\backsim}_\gamma(\mathcal S)-\kappa(\gamma)= -\psi(\gamma)$.

 We next write $\epsilon_\gamma(t)=\int_0^t\exp(\alpha \xi_{\gamma}(s))\dd s$ for the exponential functional  \eqref{E:espilonlamperti} of the L\'evy process  $\xi_\gamma$
 that appears in the Lamperti transformation.  So $t^{\prime\prime}\coloneqq \epsilon_\gamma(r^{\prime\prime})$ is the first time  at which the reproduction process $\eta_\gamma$ has atoms originating from $\mathbf{N}^{\prime\prime}$. Write furthermore $(t^{\prime\prime},x^{\prime\prime}_j)_{j\geq 1}$ for the sequence of atoms of $\eta_\gamma$ 
 at time  $t^{\prime\prime}$, in non-increasing order of the $x_j$'s.
 
 \begin{lemma} \label{lemma:xi^*prime:gamma} Under $P^\gamma_1$, and in the notation above,  the pair
 $$\left(  X_\gamma(t^{\prime\prime})/X_\gamma(t^{\prime\prime}-), \left(x_j^{\prime\prime}/X_\gamma(t^{\prime\prime}-)\right)_{j\geq 1}\right),$$
  is independent of  the restriction of the decoration-reproduction process to the time-interval $[0,t^{\prime\prime})$,
 $\left( \indset{[0,t^{\prime\prime})} X_\gamma, \indset{[0,t^{\prime\prime})\times (0,\infty)}\cdot \eta_\gamma\right)$,
and the distribution of the former is the push-forward of the law  on $\mathcal S$
 $$|\psi(\gamma)|^{-1}\left(\boldsymbol{\Lambda}^{\backsim}_\gamma-\kappa(\gamma) \delta_{(-\infty, (0, -\infty,\cdots)}\right)$$
 by the exponential map $(y, (y_j)_{j\geq 1})\mapsto (\e^y, (\e^{y_j})_{j\geq 1})$. 
 Moreover, for every functional $F\geq 0$, there is the identity 
 $$E^\gamma_1\Big( F\big( \indset{[0,t^{\prime\prime})} X_\gamma, \indset{[0,t^{\prime\prime})\times (0,\infty)}\cdot \eta_\gamma\big)\Big) = 
|\psi(\gamma)| \int_0^{\infty} \dd t \,
 E_1\Big( X(t)^{\gamma-\alpha} F\big( \indset{[0,t)} X, \indset{[0,t)\times (0,\infty)}\cdot \eta\big)\Big).$$
 \end{lemma}

\begin{proof} Let $(r^{\prime\prime},(y^{\prime\prime}, \mathbf y^{\prime\prime}))$ denote the first atom of $\mathbf{N}^{\prime\prime}$. By basic properties of Poisson point measures, 
 $(y^{\prime\prime},\mathbf y^{\prime\prime})$ is a variable in $\mathcal S$ distributed according to the normalized intensity
 $$|\psi(\gamma)|^{-1}\left(\boldsymbol{\Lambda}^{\backsim}_\gamma-\kappa(\gamma) \delta_{(-\infty, (0, -\infty,\cdots)}\right).$$
 Moreover  $r^{\prime\prime}$, $(y^{\prime\prime}, \mathbf y^{\prime\prime})$, and  $\mathbf N^{\prime}$ are independent. 
 
Next,  let $\xi^\prime$  denote the process derived from $\xi_\gamma$ by suppressing the jumps of the latter coming from $\mathbf{N}^{\prime\prime}$. 
Then  $\xi^{\prime}$ is a L\'evy process whose L\'evy-It\^{o} decomposition \eqref{E:LevyIto} uses the same  Brownian component as $\xi_\gamma$, the Poisson random measure $\mathbf{N}^{\prime}$ instead of $\mathbf{N}_\gamma$, and finally 
the drift coefficient 
 $${\mathrm a}^\prime \coloneqq {\mathrm a}_\gamma-\int_{\mathcal{S}}y\indset{|y|\leq 1} \boldsymbol{\Lambda}^\backsim_\gamma(\dd y, \dd \mathbf{y})={\mathrm a} +\sigma^2 \gamma+\int_{\mathcal{S}}  y(\e^{  \gamma y} -1) \mathbf{1}_{|y|\leq 1} \boldsymbol{\Lambda}(\dd y, \dd \mathbf{y}).$$
We stress that the latter quantity has been tuned up to take into account the compensation in Poissonian integrals.
By the L\'evy-Khintchine formula, the Laplace exponent of $\xi^{\prime}$ is $q\mapsto \psi(\gamma+q)-\psi(\gamma)$. Obviously, $r^{\prime\prime}$, $(y^{\prime\prime}, \mathbf y^{\prime\prime})$ and  $(\xi^{\prime},\mathbf N^{\prime})$ are independent. 

Since $\xi_\gamma$ and $\xi^\prime$ coincide on the  time-interval  $[0,r^{\prime\prime})$, we have
\begin{equation}\label{eq:time:change:gamma:eps}
 \epsilon_\gamma(t)=\epsilon^{\prime}(t)\coloneqq \int_0^t\exp(\alpha \xi^{\prime}(s))\dd s\qquad \text{ for all }t\leq r^{\prime\prime},
 \end{equation}
 and  by the Lamperti transformation, $(y^{\prime\prime}, \mathbf y^{\prime\prime})$ and  $\left( \indset{[0,t^{\prime\prime})} X_\gamma, \indset{[0,t^{\prime\prime})\times (0,\infty)}\cdot \eta_\gamma\right)$ are independent. By construction, there are the identities
 $$ X_\gamma(\epsilon_\gamma(r^{\prime\prime}))/ X_\gamma(\epsilon_\gamma(r^{\prime\prime})-)= \exp(y^{\prime\prime})  \quad \text{and} \quad 
x^{\prime\prime}_j / X_\gamma(\epsilon_\gamma(r^{\prime\prime})-) = \exp(y_j^{\prime\prime}) ,$$
and the first two claims of the statement follow.

We turn our attention to the third claim.  We deduce by the discussion  above that for every functional $F\geq 0$, there is the identity
\begin{align*}
&E_1^\gamma\Big( F\left( \indset{[0,t^{\prime\prime})} X_\gamma, \indset{[0,t^{\prime\prime})\times (0,\infty)}\cdot \eta_\gamma\right)\Big) \\
&= |\psi(\gamma)| \int_0^{\infty} \dd t \, E^{\gamma}_1 \left( \e^{\psi(\gamma)t} 
F\left( \indset{[0,\epsilon^{\prime}(t))} X',
 \indset{[0,\epsilon^{\prime}(t))\times (0,\infty)}
 \cdot \eta^\prime\right)\right),
\end{align*}
where $(X^\prime, \eta^\prime)$ denotes the decoration-reproduction process derived from $(\xi^\prime, \mathbf N^\prime)$ by the Lamperti transformation.
On the other hand, again from a version of the Esscher transformation (\cite[Theorem 3.9]{kyprianou2014fluctuations}), we know that
 the process $(\exp(\gamma \xi(t)-t \psi(\gamma)))_{t\geq 0}$ is a martingale under $P_1$, and for every  $t>0$, the distribution of 
 the pair
 $ \left( \indset{[0,t]} \xi,\indset{[0,t]\times \mathcal S} \cdot \mathbf{N}\right)$
  under the tilted law $\exp(\gamma\xi(t)-t\psi(\gamma))\cdot P_{1}$  is the same as that of the pair  
  $\left( \indset{[0,t]}\xi^{\prime}, \indset{[0,t]\times \mathcal S} \cdot \mathbf{N}^{\prime}\right)$. 
A final application of the Lamperti transformation yields
\begin{align*}
&\int_0^{\infty} \dd t \, E^{\gamma}_1 \left( \e^{\psi(\gamma)t} 
F\left( \indset{[0,\epsilon^{\prime}(t))} X', \indset{[0,\epsilon^{\prime}(t))\times (0,\infty)} \cdot \eta^\prime\right)\right) \\ 
 &= \int_0^{\infty} \dd s \, E_1\Big( X(\epsilon_\gamma(s))^{\gamma} F\left( \indset{[0,\epsilon_\gamma(s))} X, \indset{[0,\epsilon_\gamma(s))\times (0,\infty)}\cdot \eta\right)\Big),\end{align*}
 and   combining Tonelli's theorem with a time change using \eqref{eq:time:change:gamma:eps} we get that the previous displays equals
 $$\int_0^{\infty} \dd s \, E_1\left( X(s)^{\gamma-\alpha} F\left( \indset{[0,s)} X, \indset{[0,s)\times (0,\infty)}\cdot \eta\right)\right), $$
which completes the proof. 
\end{proof}

 The second technical lemma shows that the distributions in Lemma \ref{lemma:xi^*prime:gamma}  also naturally appear under a tilted version of $P_1$; we need to introduce some notation in this direction.
 Consider a realization of the decoration-reproduction process $(X,\eta)$ under the law $P_1$, and assume that $(t,x)$ is an atom of $\eta$, that is a child with type $x$ is born at time $t$.
Let $\mathbf x(t)=(x_j)_{j\geq 1}$ be the sequence of the types of the children born at that  time, repeated according to their multiplicities and ranked in the non-increasing order. In particular, $x$ is one of the terms of the sequence $\mathbf x$,
say $x=x_i$. We then write $(X(t),\mathbf x)^{x \backsim}$ for the pair whose first element is $x$ and second element is derived from the sequence $(x_1, \ldots, x_{i-1}, X(t), x_{i+1}, \ldots)$  (i.e. we replace a term $x$ in $\mathbf x$ by $X(t)$) after re-ordering in the non-increasing order.

\begin{lemma}\label{lemma:spine:branching:gamma} For every nonnegative functionals $F,G$, there is the identity
\begin{align*}
&E^\gamma_1\left( F\left( \indset{[0,t^{\prime\prime})} X_\gamma, \indset{[0,t^{\prime\prime})\times (0,\infty)}\cdot \eta_\gamma\right) G(X_\gamma(t^{\prime\prime}),(x^{\prime\prime}_j)_{j\geq 1})\right) \\
&= E_1\left( \int_{[0,z)\times (0,\infty)} \eta(\dd t, \dd x) x^\gamma F\left( \indset{[0,t)} X, \indset{[0,t)\times (0,\infty)}\cdot \eta\right) G\big(\left( X(t),\mathbf x(t)\right)^{x\backsim}\big) \right) \\
&-\kappa(\gamma)\cdot E_1\Big(\int_0^z \dd  t~ X(t)^{\gamma-\alpha}F\left( \indset{[0,t)} X, \indset{[0,t)\times (0,\infty)}\cdot \eta\right) G\big(0, (X(t), 0, \dots)\big)\Big).\end{align*}
\end{lemma}

\begin{proof} The notion of compensators of optional processes and random measures (see \cite[Sections I.3 and  II.1]{JS03}) lies at the heart of the proof. Roughly speaking it plays a role similar to that of the Mecke equation in  Proposition \ref{P:spinegbpgamma}.
To start with, recall from Lemma \ref{lemma:xi^*prime:gamma} that 
\begin{align*}
&E^\gamma_1\left( F\left( \indset{[0,t^{\prime\prime})} X_\gamma, \indset{[0,t^{\prime\prime})\times (0,\infty)}\cdot \eta_\gamma\right) G(X_\gamma(t^{\prime\prime}),(x^{\prime\prime}_j)_{j\geq 1})\right) \\
&= E_1\left(\int_0^{z} \dd t \,
  X(t)^{\gamma-\alpha } F\left( \indset{[0,t)} X, \indset{[0,t)\times (0,\infty)}\cdot \eta\right)
  \int_{\mathcal S} \boldsymbol{\Lambda}^{\backsim}_\gamma (\dd y, \dd \mathbf y) G(X(t)\e^y,(X(t)\e^{y_j})_{j\geq 1})\right)   \\
&  - \kappa(\gamma) E_1\left(  \int_0^{\infty} \dd t \,
  X(t)^{\gamma-\alpha } F\left( \indset{[0,t)} X, \indset{[0,t)\times (0,\infty)}\cdot \eta\right)G(0,(X(t),0, \ldots))\right) 
.\end{align*}
Comparing this identity with that of the statement, it suffices to identify the first terms of the differences in the respective right-hand sides.

We work under $P_1$ and  may assume without loss of generality that the functionals $F$ and $G$ are bounded.  In the natural filtration of $(X, \eta)$, we consider the optional  increasing process
$$A(t)\coloneqq\int_{[0,t]\times (0,\infty)} \eta(\dd s, \dd x)   x^\gamma G((X(s),\mathbf x(s))^{x \backsim}) .$$
We claim that its compensator $A^{(\mathrm p)}$ is the predictable increasing process
given by
\begin{equation}\label{Eq:compensA}
A^{(\mathrm p)}(t) = \int_0^t  \dd s X(s)^{\gamma-\alpha}   \int_{\mathcal S} \boldsymbol{\Lambda}^{\backsim}_\gamma (\dd y, \dd \mathbf y) G(X(s)\e^y,(X(s)\e^{y_j})_{j\geq 1}) ,\end{equation}
in the sense that $A(t)-A^{(\mathrm p)}(t) $ is a martingale. Since the process $F\left( \indset{[0,t)} X, \indset{[0,t)\times (0,\infty)}\cdot \eta\right)$ is predictable and bounded, it follows that
\begin{align*}
&E_1\Big( \int_{[0,z)\times (0,\infty)} F\left( \indset{[0,t)} X, \indset{[0,t)\times (0,\infty)}\cdot \eta\right) \dd A(t)\Big) \\
&= E_1\Big( \int_{[0,z)\times (0,\infty)} F\left( \indset{[0,t)} X, \indset{[0,t)\times (0,\infty)}\cdot \eta\right) \dd A^{(\mathrm p)}(t)\Big),
\end{align*}
which proves of the statement. 

We still have to check \eqref{Eq:compensA}, for which we need to return to the construction of the reproduction process  $\eta$ in terms of the Poisson point measure $\mathbf N$ and Lamperti's  transformation in Section \ref{sec:2.2}. 
The expression for the optional increasing process $A$ leads us to introduce the optional process (now in the natural filtration of $(\xi, \mathbf N)$) 
$$D(t)\coloneqq \int_{[0,t]\times \mathcal S} \mathbf N(\dd s, \dd y, \dd \mathbf{y}) \sum_{i\geq 1} \exp(\gamma(\xi(s-)+y_i)) H(\xi(s-), (y, \mathbf y)^{\backsim i}),
$$
for the functional 
$$H(r, (y, \mathbf y))\coloneqq G(\e^{r+y}, (\e^{r+y_i})_{i\geq 1}).$$ 
Observe that if we write $x=\e^{r+y}$ and $\mathbf{x}=(\e^{r+y_i})_{i\geq 1}$, then 
$$H(r, (y, \mathbf y)^{\backsim i})= G((x,\mathbf x)^{x_i\backsim}).$$
Finally,  by Poissonian calculus, the compensator $D^{(\mathrm p)}$ of $D$ is
\begin{align*}
D^{(\mathrm p)}(t)&= \int_0^t \dd s \int_{\mathcal S} \boldsymbol{\Lambda} (\dd y, \dd \mathbf{y})  \sum_{i\geq 1} \exp(\gamma(\xi(s-)+y_i)) H(\xi(s-), (y, \mathbf y)^{\backsim i})\\
&= \int_0^t \dd s \exp(\gamma \xi(s-)) \int_{\mathcal S} \boldsymbol{\Lambda}_\gamma^{\backsim} (\dd y, \dd \mathbf{y}) H(\xi(s-), (y, \mathbf y)),
\end{align*}
and \eqref{Eq:compensA} follows from an application of the Lamperti time-substitution. 
\end{proof}

 We can now tackle  the proofs of Proposition  \ref{P:spinegeneralnew} and Theorem \ref{T:Spindec:gamma}.
 \begin{proof}[Proofs of Proposition  \ref{P:spinegeneralnew} and Theorem \ref{T:Spindec:gamma}]
 Depending on whether $\gamma=\omega_-$ or $\kappa(\gamma)<0$, we use the realization of $\widetilde{\mathbb Q}^\gamma_x$ in terms of a general branching process with a distinguished lineage,  provided in Proposition \ref{P:spinebar}  and the discussion above it. We henceforth work under the associated law $\bar \P^\gamma_x$, and by scaling, we may assume  without loss of generality that $x=1$. This enables us to  explore the segment $\llbracket \rho, \rho^{\bullet}\rrbracket$ by following the trajectory of distinguished individuals, switching from parent to child at each distinguished birth event, in the sense that we then stop following the distinguished parent and rather follow its distinguished child. This allows us to analyze the decoration-reproduction process along the spine $\llbracket \rho, \rho^\bullet \rrbracket$, say $(f^\bullet, \eta^\bullet)$, 
which is defined rigorously in terms of the distinguished lineage  as follows. 
In the case $\kappa(\gamma)<0$, we write $n^\bullet\coloneqq |u^\bullet|$ for the generation of the ultimate distinguished individual, 
so $n^\bullet+1$ has the geometric law with success probability $\kappa(\gamma)/\psi(\gamma)$, as it can readily be checked from  \eqref{E:spinedecrep} and using the first lines of the proof of Proposition \ref{prop:lengthmeasures}. In the case $\gamma=\omega_-$, we set $n^\bullet=\infty$.
Next for every generation $0\leq n < n^\bullet$, we  write $t^\star(n)$ for the distinguished  age at which the distinguished individual $u^\star(n)$ begets its distinguished child, and further set 
$t^\star(n^\bullet)=t^\bullet$ in the notation of  Lemma \ref{L:spinedecgbp} when $\kappa(\gamma)<0$.
In this framework, the segment $\llbracket \rho, \rho^\bullet \rrbracket$ is realized by concatenating one after the other the segments corresponding to distinguished individuals and truncated at their distinguished  ages.
We set $s^\star(n)\coloneqq \sum_{k=0}^n t^\star(k)$, in particular $s^\star(-1)\coloneqq 0$ and $s^\star(n^\bullet)\coloneqq z^\bullet$ is the length of $\llbracket \rho, \rho^\bullet \rrbracket$.
The (rcll) decoration $f^\bullet: [0, z^\bullet)\to \R_+$ is simply obtained by concatenating the truncated decorations of distinguished individuals, 
$$f^\bullet(t)\coloneqq  f_{u^\star(n)}\big(t-s^\star(n-1)\big)\quad \text{for any }s^\star(n-1)\leq t <s^\star(n).$$

In turn, the reproduction process $\eta^\bullet$ only  partly results from the concatenation of the truncated reproduction processes of distinguished individuals; a special attention must be given to the birth-times of distinguished children. 
Discarding the latter at first, we define the point process $\eta^\circ$ by
$$\int_{[0, z^\bullet)\times (0,\infty)} \eta^\circ (\dd t, \dd x)\varphi(t,x) \coloneqq \sum_{k=0}^{n^\bullet} \int_{[0, t^\star(k))\times (0,\infty)} \eta_{u^\star(k)} (\dd t, \dd x)\varphi(t+s^\star(k-1), x),$$
where $\varphi:\R_+\times(0,\infty)\to \R_+$ stands for a generic measurable function. 
Next,  for any $0\leq k < n^\bullet$, the distinguished individual $u^\star(k)$ begets children at its distinguished age $t^\star(k)$; for simplicity just write $(x_j)_{j\geq 1}$ for  the sequence of the types of those children, repeated according to their multiplicities and, say, listed in the non-increasing order. The type  of the distinguished child is an element of this sequence, and we then write $(x^\backsim_j)_{j\geq 1}$ for the sequence obtained from $(x_j)_{j\geq 1}$ after replacing the type of the distinguished child by $f_{u^\star(k)}(t^\star(k))$, the value of the decoration immediately after this birth event. This transformation reflects the fact that at distinguished birth events, we cease to follow the distinguished parent 
and rather switch to its distinguished child. We then write 
$$\eta^{\star\backsim}_k \coloneqq \sum_{j\geq 1} \delta_{(s^\star(k), x^\backsim_j)}.$$
Finally, in the case $\kappa(\gamma)<0$, the ultimate distinguished individual $u^\bullet$ at generation $n^\bullet$ does not beget any distinguished child,
and we set 
$$\eta^{\star\backsim}_{n^\bullet} \coloneqq \delta_{(z^\bullet, f_{u^\bullet}(t^\bullet))}.$$
We now have all the ingredients to define the reproduction process on $\llbracket \rho, \rho^\bullet \rrbracket$ by 
$$\eta^\bullet\coloneqq \eta^\circ + \sum_{k=0}^{n^\bullet} \eta^{\star\backsim}_k .$$
 We claim that
\begin{itemize}
\item[(i)] The distribution of $(f^\bullet, \eta^\bullet)$  is  $P_1^\gamma$.
\item[(ii)] If we write $\eta^\bullet=\sum_{i\geq 1}\delta_{(t_i,x_i)}$ where indices are chosen in co-lexicographical order, then  conditionally on $(f^\bullet,\eta^\bullet)$, the associated standard decorated subtrees $(\mathtt{T}_i)_{i\geq 1}$ dangling from $\llbracket \rho, \rho^\bullet \rrbracket$ in  $  \mathtt{T}^{\bullet} = (\mathtt{T},\rho^\bullet)$ are independent and,  for every $i\geq 1$, the law of $\mathtt{T}_i$ is $\P_{x_i}$.
\end{itemize}
Before proving the claim let us explain why  Proposition  \ref{P:spinegeneralnew} and Theorem \ref{T:Spindec:gamma} follows directly from it. First, by Proposition  \ref{P:constructionomega-} and Lemma  \ref{L:spinedecgbp},  combined with the definition of  $\bar{\P}_1^\gamma$, we infer that under $\bar{\P}_1^\gamma$ the family $(f_u,\eta_u)_{u\in \mathbb{U}}$ satisfies Property~{\hypersetup{linkcolor=black}\hyperlink{prop:P}{$(\mathcal{P})$}}. It follows that the  families $$\big( \mathrm{Height}(T_i)\big)_{i\geq 1} \ \text{and} \ \big( \max_{T_i} g_i\big)_{i\geq 1}$$
 are  both null. Therefore Point (i) implies Proposition \ref{P:spinegeneralnew}. Moreover, by definition of the dangling subtrees and Point (ii), we must have:
 $$\mathtt{T}=\mathrm{Gluing}\Big( \big([0, z^\bullet], d,0,f^\bullet \big), (x_i)_{i\geq 1}, \big(\texttt{T}_i\big)_{j\geq 1}\Big),$$
 where as usual     $d$ stands for the usual distance on segments. 
This entails Theorem \ref{T:Spindec:gamma}, since by Point (i) the right-side hand of the previous display is distributed according to $ \widehat{\mathbb{Q}}_{x}^{\gamma}$. Let us now conclude by proving the claim. We focus on the case $\kappa(\gamma)<0$, as the case $\gamma=\omega_-$ follows from  similar (actually, slightly simpler) arguments and we leave the extension the reader -- it can also be  directly deduced from the case $\kappa(\gamma)<0$ by taking the limit $\gamma\downarrow \omega_-$  and using Proposition \ref{thm:cov:wei:leng}. In order to verify (i), recall the setting of Lemma \ref{lemma:xi^*prime:gamma}. Imagine that under the law $P^\gamma_1$, we tag the times at which the reproduction process $\eta_\gamma$ has atoms originating from the Poisson random measure $\mathbf N^{\prime\prime}$. In particular the first tagged time is $t^{\prime \prime}$, and the ultimate one is related via the Lamperti time-change to the first time when $\mathbf N^{\prime\prime}$ has an atom on the fiber $[0,\infty)\times \{(-\infty, (0, -\infty, \ldots))\}$.  It follows that the
total number of tagged times has the geometric distribution with success parameter $-\kappa(\gamma)/\left(\boldsymbol{\Lambda}^{\backsim}_\gamma(\mathcal S)-\kappa(\gamma)\right)= \kappa(\gamma)/\psi(\gamma)$, that is the same distribution as the number $n^\bullet +1$ of distinguished individuals under $\bar \P^\gamma_1$.
This is of course not a mere coincidence, and we shall actually see that we can couple $\P_1^\gamma$ and $\bar \P^\gamma_1$ such that the tagged times for $\eta_\gamma$ correspond to the times $t^\star(k)$ for $k=0, \ldots , n^\bullet$
at which a new distinguished individual is born under $\bar \P^\gamma_1$.

Indeed, consider the restriction of the decoration-reproduction process under $\bar \P^\gamma_1$ to the closed time interval $[0,t^\star(0)]$, 
$$\left( \indset{[0,t^\star(0)]} f^\bullet,  \indset{[0,t^\star(0)]\times \mathcal S}\cdot \eta^\bullet\right).$$
By the construction of $(f^\bullet, \eta^\bullet)$ and the very
definition of the general branching process with law $\bar \P^\gamma_1$, the former has the same law as 
$$\Big( \indset{[0,t^\star)} f + \indset{\{t^\star\}}x^\star,  \indset{[0,t^\star)\times \mathcal S}\cdot \eta +  \sum_{j\geq 1}\delta_{(t^\star, x^\backsim_j)} \Big)$$
under $\bar P^\gamma_1$, where the sequence $(x^\backsim_j)_{j\geq 1}$ is reduced to the single term $f(t^\star)$ if $x^\star=0$, and otherwise
is obtained by replacing  $x^\star$  in the sequence of the types of the children born at time $t^\star$ by $f(t^\star)$. 

The comparison of Proposition  \ref{P:spinegbpgamma} with Lemma \ref{lemma:spine:branching:gamma} now confirms
that the law under $\bar \P^\gamma_1$ of the restriction $\left( \indset{[0,t^\star(0)]} f^\bullet,  \indset{[0,t^\star(0)]\times \mathcal S}\cdot \eta^\bullet\right)$ 
is indeed the same as that of the restriction $\left(\indset{[0, t^{\prime\prime}]}f, \indset{[0, t^{\prime\prime}]\times \mathcal S}\cdot \eta\right)$ under $P^\gamma_1$. 
The verification of (i) can then be completed by an application of the strong Markov property, conditionally on $f^\bullet(t^\star(0))$ under $\bar \P^\gamma_1$,
respectively conditionally on $f(t^{\prime\prime})$ under $P^\gamma_1$, iteratively as long as these quantities remain non-zero.

It remains to check Point  (ii) about the dangling subtrees. Since the non-distinguished individuals under $\bar \P^\gamma_1$ have the same evolution as under $\P_1$,
we only need to consider the subtrees induced by distinguished individuals strictly after their distinguished age $t^\star$. Recall then from   Proposition \ref{P:spinegbpgamma}(ii) that the decoration-reproduction shifted at time $t$ and properly rescaled, denote by $B$, has law $P_1$. Since by definition, a distinguished parent does not have any further distinguished children strictly after time $t^\star$, we conclude, by an application of the branching property and Lemma \ref{L:Markovdecoreprod}, that the subtree induced by the $k$-th distinguished individual strictly after time $t^\star(k)$ (i.e. strictly after the distinguished individual has reached the age $t^\star(k)$) and properly rescaled has indeed the law $\mathbb{Q}_1$, independently of  all the others dangling subtrees. This completes the verification of (ii). 

\end{proof}

\section{Bifurcators} \label{sec:bifurcators}

In this section, we use the spinal decomposition to determine all the characteristic quadruplets satisfying Assumption \ref{A:gamma0} which yield  the same law on (unmarked, non-measured) decorated trees. To this end, we  introduce the map 
$$\mathrm{ord}: \mathcal{S}\to \mathcal{S}_1, \quad \mathrm{ord}\big(y_0,(y_i)_{i\geq 1}\big) =  (y^\downarrow_i)_{i\geq 1},$$
where the right-hand side denotes the sequence obtained by  ranking the collection of $(y_i)_{i\geq 0}$ in non-increasing order. 
 Borrowing the terminology from \cite{pitman2015regenerative,shi2017growth}, we say that two characteristic quadruplets $(\sigma^2, {\mathrm a},  \boldsymbol{\Lambda};\alpha)$ and $ (\sigma_{\Yleft}^2, {\mathrm a}_{\Yleft}, { \boldsymbol{\Lambda}}_{\Yleft};\alpha_{\Yleft})$ are \textbf{bifurcators} of one another and we write  $$(\sigma^2, {\mathrm a},   \boldsymbol{\Lambda};\alpha)\approx (\sigma_{\Yleft}^2, {\mathrm a}_\Yleft,  \boldsymbol{\Lambda}_{\Yleft};\alpha_{\Yleft})$$ 
 if and only if 
\begin{equation}\label{equi:bif}
\sigma^2=\sigma^2_{\Yleft}\quad \mathrm{,} \quad \boldsymbol{\Lambda}\circ  \mathrm{ord}^{-1}=\boldsymbol{\Lambda}_{\Yleft} \circ  \mathrm{ord}^{-1}\quad,\quad\alpha=\alpha_{\Yleft},
\end{equation} 
and 
\begin{equation}\label{eq:rela:a:Yleft}
{\mathrm a}-{\mathrm a}_{\Yleft}=\lim \limits_{\eps\to 0+}\Big(\int_{\esp<|y|\leq 1} \boldsymbol{\Lambda}(\dd y,\dd \mathbf{y})~y-\int_{\esp<|y|\leq1} \boldsymbol{\Lambda}_{\Yleft}(\dd y,\dd \mathbf{y})~y\Big).
\end{equation}
\begin{remark} If \eqref{equi:bif} holds for two characteristic quadruplets  that both fulfill Assumption \ref{A:gamma0}, then we will see in the proof of Theorem \ref{theo:bif} below that $ \boldsymbol{\Lambda}- \boldsymbol{\Lambda}_\Yleft$ is always a finite signed measure on $\mathcal S$. As a consequence, the condition \eqref{eq:rela:a:Yleft} can then be re-expressed in the simpler form
$${\mathrm a}_{\Yleft}={\mathrm a} + \int_{|y|\leq 1} y \left(\boldsymbol{\Lambda}_\Yleft - \boldsymbol{\Lambda} \right) (\dd y,\dd \mathbf{y}).$$
\end{remark} 

Plainly $\approx$ is an equivalence relation; the notation is also meant to suggest that when both quadruplet fulfill Assumption \ref{A:gamma0},
 the self-similar Markov trees, say $\mathtt T$  with  characteristic quadruplet $(\sigma^2, {\mathrm a},   \boldsymbol{\Lambda};\alpha)$ and $\mathtt T_\Yleft$  with  characteristic quadruplet $ (\sigma_{\Yleft}^2, {\mathrm a}_{\Yleft}, { \boldsymbol{\Lambda}}_{\Yleft};\alpha_{\Yleft})$, should be thought of as  isomorphic; see Chapter \ref{chap:topology}. 
We immediately see from  \eqref{E:cumulant} that if $(\sigma^2, {\mathrm a},  \boldsymbol{\Lambda};\alpha)\approx (\sigma_{\Yleft}^2, {\mathrm a}_{\Yleft}, \boldsymbol{\Lambda}_{\Yleft};\alpha_{\Yleft})$, then the cumulant function $\kappa_{\Yleft}$  associated to  $ (\sigma_{\Yleft}^2, {\mathrm a}_{\Yleft}, { \boldsymbol{\Lambda}}_{\Yleft})$ 
is identical to $\kappa$. 
Therefore, if $\kappa(\gamma)<0$ for some $\gamma >0$, then $\kappa_\Yleft(\gamma)<0$ as well, and 
Assumption \ref{A:gamma0} holds for both characteristic quadruplets. 

\begin{example}\label{ex:loclargest}
 An important example of bifurcator is obtained by taking $ { \boldsymbol{\Lambda}}_{\Yleft}=\boldsymbol{\Lambda}_\ast$ to be the push-forward of $ \boldsymbol{ \Lambda}$ by the transformation on $\mathcal S$ that swaps $y$ and $y_1$ if $y_1>y$ and leaves $(y,\mathbf y)$ unchanged otherwise, and then choosing $ \mathrm{a}_\ast$ so that \eqref{eq:rela:a:Yleft} is verified. See  the proof of Theorem \ref{theo:bif} below for details.
 The bifurcator $(\sigma^2, {\mathrm a}_\ast,  \boldsymbol{\Lambda}_{\ast};\alpha)$ is called the \textbf{locally largest} bifurcator since, in the genealogical interpretation, all children have a type smaller than the value of the decoration of their parent immediately after the birth event. This allows us to distinguish a canonical element in every equivalence class of bifurcators. 
 \end{example}

We stress that it is implicitly assumed in the definition of bifurcators that both $ \boldsymbol{\Lambda}$ and $\boldsymbol{\Lambda}_{\Yleft}$ are generalized L\'evy measures, and in particular that their images by the first projection $\mathcal S\to \R$, $(y,\mathbf y)\mapsto y$,  are standard L\'evy measures. For instance, if ${ \boldsymbol{\Lambda}}^\backsim$ denotes the measure on $\mathcal S$ obtained from $\boldsymbol{\Lambda}$ by swapping the first coordinate $y$ and the first term $y_1$ of the sequence $\mathbf y$, no matter whether $y_1>y$ or not,  then plainly $ \boldsymbol{\Lambda}\circ  \mathrm{ord}^{-1}=\boldsymbol{\Lambda}^\backsim\circ  \mathrm{ord}^{-1}$. However generically $ \boldsymbol{\Lambda}^\backsim((-\infty,-1]\times \mathcal S_1)=\infty$ and therefore $\boldsymbol{\Lambda}^\backsim$ is \textbf{not} a generalized L\'evy measure. 

In the sequel, we fix  $(\sigma^2, {\mathrm a}, \boldsymbol{\Lambda}; \alpha)$ and  $(\sigma_{\Yleft}^2, {\mathrm a}_{\Yleft} \boldsymbol{\Lambda}_{\Yleft}; \alpha_{\Yleft})$  two characteristic  quadruplets and   use the obvious notation $P_x$ and $P_x^{\Yleft}$ for the laws of the decoration-reproduction process of an individual with type $x>0$,   $\P_x$ and  $\P_x^{\Yleft}$
for the law of the family of decoration-reproduction processes indexed by the Ulam tree, and finally $\mathbb Q_x$ and $\mathbb Q^{\Yleft}_x$ for the laws of the self-similar Markov trees endowed with the zero measure -- when  Assumption  \ref{A:gamma0} is verified.

\begin{theorem}[Bifurcators] \label{theo:bif} Assume that  the two characteristic  quadruplets,  $(\sigma^2, {\mathrm a}, \boldsymbol{\Lambda}; \alpha)$ and  $(\sigma_{\Yleft}^2, {\mathrm a}_{\Yleft}, \boldsymbol{\Lambda}_{\Yleft}; \alpha_{\Yleft})$, satisfy Assumption  \ref{A:gamma0}. Then, $\mathbb Q_x=\mathbb Q^{\Yleft}_x$, for  all $x>0$, if and only  if $(\sigma^2, {\mathrm a},  \boldsymbol{\Lambda};\alpha)$ and $ (\sigma_{\Yleft}^2, {\mathrm a}_{\Yleft}, \boldsymbol{\Lambda}_{\Yleft}; \alpha_{\Yleft})$ are bifurcators of one another.
\end{theorem}

Before establishing Theorem \ref{theo:bif}, we shall first illustrate the statement by establishing a special case, which will turn out later to be a key step of the proof.
In this direction, let  $\boldsymbol{\Lambda}^\prime$  and $\boldsymbol{\Lambda}^{\prime\prime}$ two measures on $\mathcal{S}$ such that $\boldsymbol{\Lambda}=\boldsymbol{\Lambda}^{\prime}+\boldsymbol{\Lambda}^{\prime\prime}$, and satisfying $\boldsymbol{\Lambda}^{\prime}(\{-\infty\}\times \mathcal{S}_1)=0$
 and $\boldsymbol{\Lambda}^{\prime\prime}( \mathcal{S})<\infty$. We then write  $\boldsymbol{\Lambda}^{\prime\prime\prime}$ for the push-forward of $\boldsymbol{\Lambda}^{\prime \prime}$ by the transformation on $\mathcal S$
$$(y, \mathbf y)\mapsto  (y, \mathbf y)^\dagger\coloneqq  (-\infty, \mathrm{ord} (y, \mathbf y)).
$$
We set $\boldsymbol{\Lambda}^\dagger\coloneqq \boldsymbol{\Lambda}^\prime +\boldsymbol{\Lambda}^{\prime\prime\prime}$ and 
$$\mathrm{a}^\dagger\coloneqq \mathrm{a} - \int_{\mathcal S} y \indset{|y|\leq 1}\boldsymbol{\Lambda}^{\prime\prime}(\dd y, \dd \mathbf y).$$
Clearly, $(\sigma^2, {\mathrm a}^\dagger, \boldsymbol{\Lambda}^\dagger; \alpha)$ is a characteristic  quadruplet and we  have $(\sigma^2, {\mathrm a}^\dagger, \boldsymbol{\Lambda}^\dagger; \alpha)\approx (\sigma^2, {\mathrm a}, \boldsymbol{\Lambda}; \alpha)$,
and the next statement is thus a version of Theorem \ref{theo:bif} in this setting. With transparent notations, we write  $\mathbb Q^{\dagger}_x$ for the laws of the self-similar Markov trees  associated with  the characteristic quadruplet   $(\sigma^2, {\mathrm a}^\dagger, \boldsymbol{\Lambda}^\dagger; \alpha)$, when the latter is well defined.

\begin{lemma} \label{L:exbifurc} If  $(\sigma^2, {\mathrm a}, \boldsymbol{\Lambda}; \alpha)$ satisfies Assumption  \ref{A:gamma0}, then so does $(\sigma^2, {\mathrm a}^\dagger, \boldsymbol{\Lambda}^\dagger; \alpha)$ and  the law of the self-similar Markov tree with characteristic  quadruplet $(\sigma^2, {\mathrm a}^\dagger, \boldsymbol{\Lambda}^\dagger; \alpha)$ coincides with $(\mathbb Q_x)_{x>0}$.
\end{lemma}

\begin{proof} The strategy of the proof is similar to that of Theorem \ref{T:Spindec:gamma}. Let us proceed. First, notice that by scaling it suffices to treat the case $x=1$, and remark that $(\sigma^2, {\mathrm a}^\dagger, \boldsymbol{\Lambda}^\dagger; \alpha)$  directly satisfies Assumption  \ref{A:gamma0} since   $(\sigma^2, {\mathrm a}^\dagger, \boldsymbol{\Lambda}^\dagger; \alpha)\approx (\sigma^2, {\mathrm a}, \boldsymbol{\Lambda}; \alpha)$. We shall use a coupling argument, and for this recall  the construction of self-similar Markov decoration-reproduction processes in  Section \ref{sec:2.2}. 
Let $B$ be a standard Brownian motion,  and $ \mathbf{N}^{\prime}$ and $  \mathbf{N}^{\prime\prime} $,  two Poisson point measures  on $[0,\infty)\times  \mathcal{S}$ with respective intensities $\dd t   \boldsymbol{\Lambda}^\prime(\dd  y , \dd \mathbf{y}) $ and $ \dd t   \boldsymbol{\Lambda}^{\prime\prime}(\dd  y , \dd \mathbf{y}) $. 
 We assume that  $B$,  $\mathbf{N}^{\prime}$, and  $ \mathbf{N}^{\prime\prime} $ are independent. Let $\mathbf N^{\prime\prime\prime}$ denote the image of $ \mathbf{N}^{\prime\prime} $ by the map $(t,(y,\mathbf y))\mapsto  (t, (y, \mathbf y)^\dagger)$. 
Then $\mathbf{N}= \mathbf{N}^{\prime} +  \mathbf{N}^{\prime\prime}$ is a Poisson point measure with intensity $\dd t \boldsymbol{\Lambda}(\dd y, \dd \mathbf{y})$, whereas 
$\mathbf{N}^\dagger= \mathbf{N}^{\prime} + \mathbf{N}^{\prime\prime\prime}$ is a Poisson point measure with intensity $\dd t \boldsymbol{\Lambda}^\dagger(\dd y, \dd \mathbf{y})$, and each of these point measures is independent of $B$. 
 We then construct the L\'evy processes $\xi$ and $\xi^\dagger$, and the  decoration-reproduction processes $(X,\eta)$ and $(X^\dagger, \eta^\dagger)$ with law $P$ and $ P^\dagger$, by using the same Brownian motion $B$,  the point measures  $ \mathbf{N}$ and  $ \mathbf{N}^\dagger$, and the drift coefficients $\mathrm a$ and $\mathrm a^\dagger$,  respectively. 
  
 Since  $0<\boldsymbol{\Lambda}^{\prime\prime}(\mathcal S)<\infty$, we may consider the first time $\uptau_1$ at which $ \mathbf{N}^{\prime\prime}$ has an atom.  The L\'evy processes $\xi$ and $\xi^\dagger$ coincide on the time interval $[0, \uptau_1)$, and so do the point processes $\mathbf{N}$ and $\mathbf{N}^\dagger$. Performing the Lamperti transformation with the exponential functional $\upepsilon$ in \eqref{E:espilonlamperti}, the restriction of the decoration processes $X$ and $X^\dagger$ to the time interval $[0,\epsilon(\uptau_1))$ coincide.
 Since $\boldsymbol{\Lambda}^{\prime}(\{-\infty\}\times \mathcal{S}_1)=0$, the measure  $\mathbf N^\prime$ has no atoms on $\R_+\times \{-\infty\}\times \mathcal S_1$ and 
 $$\uptau_1=\inf\big\{t\geq 0:  \mathbf N^\dagger([0,t]\times \{-\infty\}\times \mathcal S_1)>0\big\}.$$
 Hence $\upepsilon(\uptau_1)=z^\dagger$ is the lifetime of $X^\dagger$, and is  smaller than the lifetime $z$ of $X$. Moreover, from the very definition of the transformation $(y, \mathbf y)\mapsto  (y, \mathbf y)^\dagger$,
 we can express the entire reproduction process $\eta^\dagger$ as
 $$\eta^\dagger=  \indset{[0,\epsilon(\uptau_1)]\times(0,\infty)}\cdot \eta + \mathbf{1}_{X(\upepsilon(\uptau_1))>0}\cdot\delta_{(\upepsilon(\uptau_1), X(\upepsilon(\uptau_1)))}.$$

 The interpretation in terms of the evolution of populations  is that  $(X^\dagger, \eta^\dagger)$ results from $(X,\eta)$ by killing at time $z^\dagger$, the reproduction event for $\eta$ occurring at time $z^\dagger$ being still taken into account in $\eta^\dagger$,  and adding an extra child with type $X(z^\dagger)$ at that time (which can be of type $0$ if $X$ is also killed at $z^\dagger$). Let us tag this extra child to distinguish it from the progeny of the $\dagger$-parent that stems from $\eta$. It should now be intuitively clear from the Markov property for $(X,\eta)$ (see Lemma \ref{L:Markovdecoreprod})  that such a killing combined with the addition of an extra child at the killing time having precisely the type given by the decoration of the parent when it is killed is essentially a neutral operation for the population model (even though it clearly impacts the genealogical structure). More precisely, when $X(z^\dagger)>0$, we write $\uptau_2$ for the second time at which $ \mathbf{N}^{\prime\prime}$ has an atom. We can the use for  the decoration-reproduction process of the tagged child the restriction of $(X,\eta)$ to the time-interval $(\uptau_1, \uptau_2]$ and shifted by $-\uptau_1$ backward in time.
  The concatenation of the decoration-reproduction processes of the $\dagger$-parent and its tagged child is then given by the restriction of $(X,\eta)$ to the time-interval $[0, \uptau_2]$. Concatenating iteratively along the tagged lineage, we then recover the entire process $(X,\eta)$, and hence in this respect, the claimed identity $\mathbb Q_1 = \mathbb Q^\dagger_1$ thus should not come as a surprise -- here we use the transparent notation $ \mathbb Q^\dagger_1$ for the law of the self-similar Markov tree with initial decoration $1$ associated with  $(\sigma^2, {\mathrm a}^\dagger, \boldsymbol{\Lambda}^\dagger; \alpha)$.
 
 More formally, recalling the gluing construction of self-similar Markov trees, the analysis above enables us to construct recursively two decorated trees $\texttt T=(T,d_T,\rho,g)$ and $\texttt T^\dagger=(T^\dagger,d_{T^\dagger},\rho^\dagger,g^\dagger)$, such that the distributions induced on $\mathbb{T}$ are respectively $\mathbb Q_1$ and $\mathbb Q_1^\dagger$, and $T^\dagger$ is a subtree of $T$, $\rho^\dagger=\rho$, $d_{T^\dagger}$ corresponds to the restriction distance of $d_T$ to $T^\dagger$, and $g^\dagger(\upsilon)\leq g(\upsilon)$, for every $\upsilon \in T^\dagger$. To conclude the proof, we need to establish that $T=T^\dagger$ and $g=g^\dagger$. To this end, fix $\gamma>0$ such that $\kappa(\gamma)<0$.  Recalling that the functions used in the gluing construction of Section \ref{sec:1.2} are the usc-modification of rcll functions that can vanish only at the end of their lifetimes, we infer that to obtain the desired result, it suffices to show that:
  $$\int_T g(\upsilon)^{\gamma-\alpha} ~\uplambda_T(\d \upsilon)= \int_{T^\dagger} \big(g(\upsilon)^{\dagger}\big)^{\gamma-\alpha} ~\uplambda_{T^\dagger}(\d \upsilon).$$
We already know that the left-hand term is greater than or equal to the right-hand side. Moreover, by Proposition \ref{prop:lengthmeasures} and since $(\sigma^2, {\mathrm a}, \boldsymbol{\Lambda}; \alpha)$ and $(\sigma^2, {\mathrm a}^\dagger, \boldsymbol{\Lambda}^\dagger; \alpha)$ have the same cumulant function, both quantities have the same expectation and therefore they must be equal a.s. for every $g$. This completes the proof of the lemma.
 \end{proof}

We now prove  Theorem \ref{theo:bif}.

\begin{proof}[Proof of  Theorem \ref{theo:bif}]
We first establish the sufficiency part.
Consider a characteristic quadruplet $(\sigma^2, {\mathrm a},  \boldsymbol{\Lambda};\alpha)$ that fulfills Assumption  \ref{A:gamma0} and fix any $\gamma>0$  such that $\kappa(\gamma)<0$. 
As we suggested in Example \ref{ex:loclargest}, it is natural to consider the locally largest bifurcator. In order to introduce the latter  rigorously, we consider the partition $\mathcal{S}=\mathcal{S}_{< } \sqcup \mathcal{S}_{\geq }$,
where $\mathcal{S}_{<}\coloneqq \{(y,\mathbf y)\in \mathcal{S}: y< y_1\}$ and $\mathcal{S}_{\geq } \coloneqq \mathcal{S} \backslash \mathcal{S}_{<}$,
and observe, from the fact that $\Lambda_0(\d y)= \boldsymbol{\Lambda}(\dd y, \mathcal S_1)$ is a L\'evy measure and the finiteness of $\kappa(\gamma)$, that
\begin{align*}
\boldsymbol{\Lambda}(\mathcal{S}_{<}) &\leq \int_{\mathcal{S}}\boldsymbol{\Lambda}(\dd  y, \dd \mathbf{y})\mathbf{1}_{y\leq -1}+\int_{\mathcal{S}}\boldsymbol{\Lambda}(\dd  y, \dd \mathbf{y})\mathbf{1}_{y_1\geq -1}\\
&\leq \Lambda_0([-\infty,-1])+\e^{\gamma}\int_{\mathcal{S}}\boldsymbol{\Lambda}(\dd  y, \dd \mathbf{y})\sum \limits_{i\geq 1} \e^{\gamma y_i}<\infty.
\end{align*}

Write $\boldsymbol{\Lambda}_\ast$ for  the push-forward of  the generalized L\'evy measure  $\boldsymbol{\Lambda}$ by the transformation $ (y, \mathbf y)\mapsto  \mathrm{ord}(y, \mathbf y)$, 
so  $\boldsymbol{\Lambda}_\ast- \boldsymbol{\Lambda}$ is a finite signed measure on $\mathcal S$. We then set
$$
{\mathrm a}_\ast\coloneqq {\mathrm a}+\int_{|y|\leq 1} y \left(\boldsymbol{\Lambda}_\ast- \boldsymbol{\Lambda}\right)(\dd y, \dd \mathbf{y}).$$
By construction, $(\sigma^2, {\mathrm a},  \boldsymbol{\Lambda};\alpha)$ and $(\sigma^2, {\mathrm a}_\ast,  \boldsymbol{\Lambda}_\ast;\alpha)$ are bifurcators one another. Now remark that $\boldsymbol{\Lambda}_\ast- \mathbf{1}_{\mathcal{S}_{\geq }}\boldsymbol{\Lambda}$ is a finite measure, and let us write $\boldsymbol{\Lambda}^{\prime \prime\prime}_*$ for the push-forward measure of  $\boldsymbol{\Lambda}_\ast- \mathbf{1}_{\mathcal{S}_{\geq }}\boldsymbol{\Lambda}$ by $(y, \mathbf y)\mapsto  (y, \mathbf y)^\dagger$. We set $\boldsymbol{\Lambda}^\dagger_\ast= \mathbf{1}_{\mathcal{S}_{\geq }}\boldsymbol{\Lambda}+ \boldsymbol{\Lambda}^{\prime \prime\prime}_*$ and 
$${\mathrm a}_\ast^\dagger \coloneqq {\mathrm a}_\ast+\int_{|y|\leq 1} y \left(\boldsymbol{\Lambda}^\dagger_\ast- \boldsymbol{\Lambda}_\ast\right)(\dd y, \dd \mathbf{y}).$$
We define similarly $\boldsymbol{\Lambda}^\dagger$ and ${\mathrm a}^\dagger$ replacing $\boldsymbol{\Lambda}$ by $\boldsymbol{\Lambda}^\dagger$, $\boldsymbol{\Lambda}_\ast- \mathbf{1}_{\mathcal{S}_{\geq }}\boldsymbol{\Lambda}$ by $\mathbf{1}_{\mathcal{S}_{<}}\boldsymbol{\Lambda}=\boldsymbol{\Lambda}- \mathbf{1}_{\mathcal{S}_{\geq }}\boldsymbol{\Lambda}$  and ${\mathrm a}_\ast$ by $\mathrm{a}$  in the construction above. By definiton, we have the identity  
$$(\sigma^2, {\mathrm a}^\dagger,  \boldsymbol{\Lambda}^\dagger;\alpha) = (\sigma^2, {\mathrm a}_\ast^\dagger,  \boldsymbol{\Lambda}_\ast^\dagger;\alpha),$$
and we deduce from Lemma \ref{L:exbifurc} that
$$(\mathbb Q_x)_{x>0}=(\mathbb Q^\ast_x)_{x>0},$$
where as usual $ \mathbb Q^\ast_x$ stands for the law of the self-similar Markov tree with initial decoration $x$ associated with  $(\sigma^2, {\mathrm a}_\ast, \boldsymbol{\Lambda}_\ast; \alpha)$. Finally  consider another characteristic quadruplet 
$ (\sigma^2, {\mathrm a}_{\Yleft},  \boldsymbol{\Lambda}_{\Yleft};\alpha) \approx 
 (\sigma^2, {\mathrm a},  \boldsymbol{\Lambda};\alpha)$. 
 Then it is readily checked that the locally largest bifurcator for $(\sigma^2, {\mathrm a}_{\Yleft},  \boldsymbol{\Lambda}_{\Yleft};\alpha)$ is again
 $(\sigma^2, {\mathrm a}_\ast,  \boldsymbol{\Lambda}_\ast;\alpha)$, from which we conclude that $(\mathbb Q_x)_{x>0}=(\mathbb Q^\Yleft_x)_{x>0}$.

We turn our attention to the necessary condition.
 Suppose that the laws $\mathbb Q_x$ and $\mathbb Q^{\Yleft}_x$ coincide for all $x>0$.
 Comparing for instance the distribution of the height of $T$ under $\mathbb Q_x$, $\mathbb Q^{\Yleft}_x$, $\mathbb Q_1$ and $\mathbb Q^{\Yleft}_1$,
 we immediately see that the exponents of self-similarity $\alpha$ and $\alpha_\Yleft$ must be the same. 
 We can now focus on the case $x=1$ and then drop the index $x$ from the notation, i.e. we write as often $\mathbb Q=\mathbb Q_1$ and $\mathbb Q_{\Yleft}=\mathbb Q^{\Yleft}_1$.

Next, consider the random variable $\int_T g(v)^{\gamma-\alpha} \uplambda_T(\dd v)$,
 where we recall that $\uplambda_T$ stands for the length (or Lebesgue) measure on $T$, and we write $\kappa$ and $\kappa_{\Yleft}$ for the cumulant function associated respectively with $ (\sigma^2, {\mathrm a},  \boldsymbol{\Lambda};\alpha)$ and $ (\sigma^2, {\mathrm a}_{\Yleft},  \boldsymbol{\Lambda}_{\Yleft};\alpha) $. 
According to   Proposition \ref{prop:lengthmeasures}, the variable $\int_T g(v)^{\gamma-\alpha} \uplambda_T(\dd v)$ has expectation $-1/\kappa(\gamma)<\infty$ under $\mathbb Q$ and $-1/\kappa_\Yleft(\gamma)<\infty$ under  $\mathbb Q_{\Yleft}$.
This forces $\kappa_\Yleft(\gamma)= \kappa(\gamma)<0$.  So we can  equip the decorated real tree $\texttt T$    with
  the (finite) weighted length measure $\upnu=\uplambda^{\gamma}$ under $\mathbb Q$  as well as  under $\mathbb Q_{\Yleft}$. The identity $\mathbb Q=\mathbb Q_{\Yleft}$ then extends to the framework of decorated real trees with a single marked point of Section \ref{Section:spine}.
  More precisely, we define   
   the law  $\widetilde{\mathbb Q}^\gamma$ (respectively, $\widetilde{\mathbb Q}^\gamma_{\Yleft}$) on $\TT^\bullet$ as in \eqref{eq:defspinebiaised}, that is  by first biasing ${\mathbb Q}$ (respectively, $\mathbb Q_{\Yleft}$) with the variable $\kappa(\gamma) \upnu(T)$ and then picking a point $\rho^\bullet$ in $T$ at random according to the normalized probability measure $\upnu(\dd v)/\upnu(T)$. Obviously, we have again $\widetilde{\mathbb Q}^\gamma=\widetilde{\mathbb Q}^\gamma_{\Yleft}$. 
   
   Recall from the preceding section that in this setting,  $f^\bullet$ denotes the decoration on the marked segment $\llbracket \rho, \rho^{\bullet}\rrbracket$
   and $\eta^\bullet$ the point process on $\llbracket \rho, \rho^{\bullet}\rrbracket\times (0,\infty)$ that records the germs of the decorations of the subtrees dandling from 
$\llbracket \rho, \rho^{\bullet}\rrbracket$, see the proof of Theorem \ref{T:Spindec:gamma} for a formal definition.
     We know from the spinal decomposition, i.e. Theorem \ref{T:Spindec:gamma},  that the law of $(f^\bullet, \eta^\bullet)$ 
   under  $\widetilde{\mathbb Q}^\gamma$ is $P^\gamma$,  that is that of the self-similar Markov decoration-reproduction process with tilted characteristic quadruplet
   $( \sigma^2,\mathrm{a}_{\gamma}, \boldsymbol{\Lambda}_{\gamma};  \alpha)$, where the drift  $\mathrm{a}_{\gamma}$ is given by \eqref{eq:driftspine}
   and the generalized L\'evy measure $\boldsymbol{\Lambda}_\gamma$ by \eqref{Lambda^*_def_2:gamma}. Using an obvious notation under $\widetilde{\mathbb{Q}}^\gamma_\Yleft$, we arrive at the identity 
   $P^\gamma=P^\gamma_\Yleft$ where the right-hand side is the law of the self-similar Markov decoration-reproduction process with tilted characteristic quadruplet
   $( \sigma_{\Yleft}^2,\mathrm{a}_{\gamma}^{\Yleft}, \boldsymbol{\Lambda}_{\gamma}^{\Yleft}; \alpha)$.
   
      Just in the same way as the distribution of a L\'evy process determines its characteristic triplet, one readily sees by
  undoing the Lamperti transformation that the law of a self-similar Markov decoration-reproduction process with a given exponent of self-similarity $\alpha$ entirely determines its characteristic quadruplet. We infer from above that there is the identity
 $$( \sigma^2,\mathrm{a}_{\gamma}, \boldsymbol{\Lambda}_{\gamma}) = ( \sigma_{\Yleft}^2,\mathrm{a}_{\gamma}^{\Yleft}, \boldsymbol{\Lambda}_{\gamma}^{\Yleft}).$$

  To conclude the proof of the necessary part, it now suffices to observe from  \eqref{Lambda^*_def_2:gamma} that the push-forward of the generalized L\'evy measure $\boldsymbol{\Lambda}$
   by the function $\mathrm{ord}: \mathcal{S}\to \mathcal{S}$ that ranks all the terms of $(y, \mathbf{y})$ in the non-increasing order can be expressed in terms of $\boldsymbol{\Lambda}_{\gamma}$ and $\kappa(\gamma)$ via the identity
 $$\left( \e^{\gamma y}+\sum \limits_{i\geq 1}\e^{\gamma y_i}\right) \cdot \left(\boldsymbol{\Lambda}\circ \mathrm{ord}^{-1}\right) (\dd y,\dd \mathbf{y}) 
 =  \left(\left( \boldsymbol{\Lambda}_{\gamma}+
 \kappa(\gamma) \delta_{\{-\infty\}\times \{0, -\infty,\cdots\}}\right)\circ \mathrm{ord}^{-1}\right)(\dd y,\dd \mathbf{y}).
 $$
Recalling that $\kappa(\gamma)= \kappa_\Yleft(\gamma)$, we now see that $\boldsymbol{\Lambda}\circ \mathrm{ord}^{-1}=\boldsymbol{\Lambda}_\Yleft \circ \mathrm{ord}^{-1}$. Finally, we check \eqref{eq:rela:a:Yleft}
using again $\kappa(\gamma)= \kappa_\Yleft(\gamma)$, \eqref{E:cumulant}, and the L\'evy-Khintchine formula \eqref{E:LKfor}. 
\end{proof}

We saw in the above proof that the tilted characteristics $( \sigma^2, \mathrm{a}_\gamma, \boldsymbol{ \Lambda}_\gamma ; \alpha)$ can be recovered from the laws $( \mathbb{Q}_x)_{x >0}$. Those tilted characteristics are in fact  more intrinsic that the initial ones since they generically uniquely characterize the law of the ssMt. We say that a generalized L\'evy measure $ \boldsymbol{\Lambda}$ is \textbf{asymmetric} if there exists a point $\{y_{0},(y_{1}, y_{2}, ...)\} \in \mathcal{S}$ in the support of $ \boldsymbol{ \Lambda}$ so that $y_{0}>y_{1}>-\infty$. Then we have

\begin{corollary} \label{C:tiltedunique} Suppose that $ \boldsymbol{\Lambda}$ is asymmetric. The law $(\mathbb{Q}_x)_{x >0}$ is uniquely characterized by the data $( \sigma^2, \mathrm{a}_\gamma, \boldsymbol{ \Lambda}_\gamma ; \alpha)$.
\end{corollary}
Remark that in the binary conservative case, the tilted characteristics are completely described by the sole L\'evy--Khintchine exponent $\psi_\gamma$ of the pssMp decoration $X_{\gamma}$ along the tagged branch since the decoration-reproduction process $\eta_{\gamma}$ is given in terms of $X_{\gamma}$ by \eqref{Eq:reprodbinary}.  See the end of the next section for applications.

\begin{proof} We saw in the previous proof that $ ( \mathbb{Q}_x)_{x>0}$ is characterized by the data of $(\kappa(\gamma), \gamma)$ together with the tilted characteristic quadruplet $( \sigma^2, \mathrm{a}_\gamma, \boldsymbol{ \Lambda}_\gamma ; \alpha)$. Using \eqref{eq:kappaspi}, we deduce that the function $\kappa$ can itself be recovered from $( \sigma^2, \mathrm{a}_\gamma, \boldsymbol{ \Lambda}_\gamma ; \alpha)$ together with the parameter $\gamma$. To complete the proof it suffices to show that $\gamma>0$ can be recovered from the generalized L\'evy measure $\boldsymbol{ \Lambda}_\gamma$ only. To prove this, let us denote by $ \boldsymbol{ \Lambda}_\ast$ the generalized L\'evy measure of the locally largest bifurcator so that we have 
$$ \boldsymbol{\Lambda}_\gamma(\dd y, \dd \mathbf{y}) = \mathrm{e}^{\gamma y} \boldsymbol{\Lambda}_\ast(\dd y, \dd \mathbf{y}) +
  \e^{\gamma y}\cdot \left( \sum \limits_{i\geq 1}  \boldsymbol{\Lambda}_\ast^{\backsim i}(\dd y, \dd \mathbf{y})\right).$$
Since $  \boldsymbol{\Lambda}_{\ast}$ is asymmetric (this is equivalent to asking that a bifurcator $ \boldsymbol{\Lambda}$ or even $ \boldsymbol{ \Lambda}_{\gamma}$ is asymmetric), then there exists a point $( y_0, y_1, y_2, ... ) \in \mathcal{S}_1$ with $y_0>y_1=y_2= \cdots = y_k> y_{k+1} \geq  .... \geq -\infty$ in the support of $\boldsymbol{ \Lambda}_\ast$. Let us now consider the measure $ \boldsymbol{\Lambda}_\gamma$ restricted to the vicinity of the points  $( y_0, (y_1, y_2, ... ))$ and $(y_1,(y_0,y_2,...))$ in $ \mathcal{S}$. Denote them respectively by $\boldsymbol{\Lambda}^0_\gamma$ and $\boldsymbol{\Lambda}^1_\gamma$. Using the previous display, we deduce that $\boldsymbol{\Lambda}^1_\gamma \circ \mathrm{Ord}^{-1}$ and $\boldsymbol{\Lambda}^0_\gamma$ are absolutely continuous with respect to each-other in the vicitinity of $( y_0, (y_1, y_2, ... ))$ with Radon--Nikodym derivative equal to 

$$ \frac{\mathrm{d} \big( \boldsymbol{\Lambda}^1_\gamma \circ \mathrm{Ord}^{-1} \big)}{\mathrm{d} \boldsymbol{\Lambda}^0_\gamma}( y_0, (y_1, y_2, ... )) = \frac{ c\cdot \mathrm{e}^{\gamma y_1}}{\mathrm{e}^{\gamma y_0}},$$
where $c>0$ is an explicit constant.  Since $y_{1} \ne y_{0}$  are known, this formula thus enables us to recover $\gamma$ from the knowledge of $ \boldsymbol{ \Lambda}_\gamma$.\end{proof}

In the symmetric case, the previous corollary may not hold. Consider for example the two finite generalized L\'evy measures
 \begin{eqnarray*} \boldsymbol{\Lambda} &=&  \delta_{(-\log 3, (-\log 3, -\log 3, -\infty, ...))} + 3\cdot  \delta_{(\log 3, ( -\infty, -\infty, ...))}+7 \cdot  \delta_{(-\infty, ( -\infty, -\infty, ...))},\\
 \tilde{ \boldsymbol{\Lambda}}& =& 3\cdot \delta_{(-\log 3, (-\log 3, -\log 3, -\infty, ...))} + \delta_{(\log 3, ( -\infty, -\infty, ...))}+7\cdot  \delta_{(-\infty, ( -\infty, -\infty, ...))}.  \end{eqnarray*}
Then a straightforward calculation shows that the $(\gamma=1)$-tilted characteristics of $(0,0, \boldsymbol{ \Lambda} ; \alpha)$ coincide with those of the $(\tilde{\gamma}=2)$-tilted characteristics of $(0,0, \tilde{\boldsymbol{ \Lambda}} ; \alpha)$.

\section{Hausdorff dimensions} \label{sec:Hausdim}
In this section, we use  the spinal decomposition of Theorem \ref{T:Spindec:gamma} to determine the Hausdorff dimension of some random sets that appear naturally for self-similar Markov trees satisfying the first Cramer's condition, completing  Lemma \ref{Lem:haus:upper} in this context.

\begin{proposition}\label{Prop:spine:Haus:bis} Fix $ (\sigma^2, \mathrm{a},  \boldsymbol{\Lambda} ; \alpha)$ satisfying Assumption \ref{A:omega-} for some $\omega_->0$. Then  $\mathbb P_1$-a.s., the Hausdorff dimensions of $\partial_0 T$, $T$ and 
$\mathrm{Hyp}(g)$ are $\omega_-/\alpha$, $1\vee (\omega_-/\alpha)$ and  $2\vee (\omega_-/\alpha)$ respectively.
\end{proposition}
\begin{proof}
Thanks to Lemma \ref{Lem:haus:upper}, we only need to establish that $\mathrm{dim}_{H}( \partial_0 T) \geq \omega_-/\alpha$,  $\mathbb P_1$-a.s. In this direction, we claim that it suffices to establish that 
 \begin{eqnarray} \label{eq:goalHDlower} \limsup_{r \to 0} \frac{\upmu\big(B_{r}(\rho^{\bullet})\big)}{r^{\frac{\omega_{-}}{\alpha}-\delta}} = 0, \quad \widehat {\mathbb Q}_1^{\omega_-}\text{-a.s.}, \end{eqnarray}
where we write $\rho^\bullet$ for the point corresponding in $\mathtt{T}$ to the extremity $z_\varnothing$ of the ancestral individual $(f_\varnothing, \eta_\varnothing)$, and $B_{r}(\rho^{\bullet})$ for the closed  ball of radius $r$ centered at $\rho^\bullet$. Indeed, recall from Proposition~\ref{P:newPleaves} that  the harmonic measure is supported on $\partial_0 T$, $\mathbb{P}_1$- a.s. Since by   \eqref{eq:defspinebiaised} and Theorem \ref{T:Spindec:gamma}, under  $\widehat {\mathbb Q}_1^{\omega_-}$,  the marked point $\rho^{\bullet}$ is distributed according to  $\upmu(\dd v)/\upmu(T)$, it follows from standard density theorems for Hausdorff measures  that $ \mathrm{dim}_{H}(\partial_0 T) \geq \omega_-/\alpha$, $\mathbb P_1$-a.s.

To prove the claim \eqref{eq:goalHDlower}  recall that, under $\widehat {\mathbb Q}_1^{\omega_-}$, the decoration-reproduction process  $(f_\varnothing, \eta_\varnothing)$ of the spine is distributed according to the biased decoration-reproduction kernels $(P^{\omega_-}_x)_{x>0}$. Recall also from Section \ref{sub:section:temp} the notation $\check{B}_a^{\bullet}(T)$   for  the closure of the  complement of  the hull of radius $a$ when $d_T(\rho,\rho^\bullet)>a$. By definition, $\check{B}_a^{\bullet}(T)$ contains the open ball $B_{r}(\rho^{\bullet})$ for the radius  $r=d_T(\rho,\rho^\bullet)-a$. Finally, fix $\delta >0$ and for every $\varepsilon>0$, set $\vartheta_{ \varepsilon} := \inf \{ t \geq 0 : f_\varnothing(t) < \varepsilon\}$. We shall show that for every $\delta_1,\delta_2>0$  such that 
\begin{equation}\label{eq:tune:delta:hauss}
\omega_--\delta_1>(\alpha+\delta_2)(\omega_-/\alpha -\delta), 
\end{equation}
we have $\widehat {\mathbb Q}_1^{\omega_-}$-a.s. that for every $k\in \mathbb{N}$ large enough,
 \begin{equation}\label{eq:haus:even}
\upmu \big ( \check{B}_{\vartheta_{ 2^{-k}}}^{\bullet}(T)\big)\leq 2^{-k(\omega_{-}-\delta_1)} \quad \text{ and }\quad   2^{-k (\alpha+ \delta_2)}\leq d_{T}( \rho, \rho^{\bullet}) - \vartheta_{2^{-k}}.
\end{equation}
The desired result \eqref{eq:goalHDlower} then follows, since  we deduce from the inclusions 
$$B_{2^{-k (\alpha+ \delta_2)}}(\rho^\bullet)\subset B_{d_{T}( \rho, \rho^{\bullet}) - \vartheta_{2^{-k}}}(\rho^\bullet)\subset  \check{B}_{\vartheta_{ 2^{-k}}}^{\bullet}(T)$$
that
$$\upmu\big( B_{2^{-k (\alpha+ \delta_2)}}(\rho^\bullet)\big) \leq 2^{-k(\omega_{-}-\delta_1)}, \quad \text{ for $k$ sufficiently large,}$$
 and we conclude from \eqref{eq:tune:delta:hauss}.

Let us proceed with the proof of \eqref{eq:haus:even}. First remark that by the Markov property of $(f_\varnothing,\eta_{\varnothing})$ of Lemma \ref{L:Markovdecoreprod}  and the branching property, conditionally on $f_{\varnothing}( \vartheta_{ \varepsilon})$, the variable
$ \upmu(  \check{B}_{\vartheta_{ 2^{-k}}}^{\bullet}(T)\big)$ is distributed as $\upmu(T)$ under $\widehat {\mathbb Q}_{f_{\varnothing}( \vartheta_{ \varepsilon})}^{\omega_-}$. Therefore, by the scaling property combined with the fact that $f_{\varnothing}( \vartheta_{ \varepsilon})\leq \eps$,   we get 
$$\widehat {\mathbb Q}_1^{\omega_-}\Big( \upmu\big(  \check{B}_{\vartheta_{ 2^{-k}}}^{\bullet}(T)\big)  \geq  \varepsilon^{\omega_{-}-\delta_1}\Big) 
\leq \widehat {\mathbb Q}_1^{\omega_-}\Big(  \upmu(  T)   \geq \varepsilon^{-\delta_1}\Big).$$ 
Since $\E_1(\upmu(T)^{p})<\infty$ for some $p>1$ appearing in Assumption \ref{A:omega-}, the variable $\upmu(T)^{p-1}$ has a finite mean under $\widehat {\mathbb Q}_1^{\omega_-}$. 
By the Markov inequality, the right-hand side of the last display is  $ O(\varepsilon^{{\delta_1(p-1)}})$. The first inequality in \eqref{eq:haus:even} follows by taking $ \varepsilon= 2^{-k}$ and applying the Borel--Cantelli lemma.
On the other hand, recall 
 that under $\widehat {\mathbb Q}_1^{\omega_-}$,  the decoration process $f_\varnothing= X_{\omega_-}$ along the spine is obtained by performing the Lamperti  transformation to a L\'evy process $\xi_{\omega_-}$, and from Lemma \ref{L:LKtilde} that the latter has Laplace exponent $q\mapsto\kappa(\omega_-+q)$. 
From basic results on L\'evy processes, we have  $\xi_{\omega_-}(t) \sim t\cdot  \widehat{\mathbb{E}}_1^{\omega_-}(\xi_{\omega_-}(1))$  almost surely as $t \to \infty$, and we infer from the  Lamperti transformation that a.s.~we eventually have 
  $$  2^{-k (\alpha+ \delta_2)}\leq d_{T}( \rho, \rho^{\bullet}) - \vartheta_{2^{-k}}.$$ 
  This completes the verification of  \eqref{eq:haus:even} and hence the proof. 
\end{proof}

\section{Back to Examples} \label{sec:examplespinal} \label{sec:spinalex}

We now revisit the Examples of Chapter \ref{chap:example} and explicit their spinal decomposition. We also recall some background about stable L\'evy processes conditioned to die continuously at zero, since they appear in the $\omega_{-}$-spinal decomposition of  many natural examples. Interestingly, we will see that the only spectrally negative $\beta$-stable L\'evy processes conditioned to die continuously at zero that may appear in binary conservative ssMt are present in Examples  \ref{ex:brownian} and \ref{ex:stablefamily},  and satisfy 
$$ \beta \in (0,1/2] \cup (1,3/2].$$

Recall from   \eqref{Lambda^*_def_2:gamma} and \eqref{eq:driftspine} the definition of the tilted characteristics $( \sigma, \mathrm{a}_\gamma, \boldsymbol{ \Lambda}_\gamma ; \alpha)$ in the spinal decomposition (Theorem \ref{T:Spindec:gamma}) for $\kappa(\gamma)<0$ or $\gamma= \omega_-$ together with Assumption \ref{A:omega-}. Recall also from \eqref{eq:kappaspi} that the L\'evy--Khintchine exponent $\psi_\gamma$ of the L\'evy process underlying the pssMp $X_\gamma$ of the decoration along the distinguished tagged branch under $ \widetilde{\mathbb{Q}}^\gamma$ is given by $$\psi_{\gamma}(q) = \kappa(\gamma+q).$$
Similarly to what we did in the opening of Chapter \ref{chap:example}, in the case when the image $ { \Lambda}_{\gamma,0}$ of $ \boldsymbol{\Lambda}_\gamma$ by the application $(y, \mathbf{y}) \mapsto y$ integrates $1 \wedge |y|$ we define the canonical drift coefficient 
$$ \mathrm{ \mathrm{a}_{\gamma}^ {\mathrm{can}}} \coloneqq {\mathrm a}_\gamma - \int \Lambda_{\gamma,0}(\dd y)\  y {\mathbf 1}_{|y|\leq 1}.$$ It is easy to see that when $\kappa(\gamma)\leq 0$ then $ {\Lambda}_{\gamma,0}$ integrates $ 1 \wedge |y|$ if and only if ${\Lambda}_0$ does and in this case we have   \begin{eqnarray} \label{eq:acangamma} \mathrm{a^{can}_\gamma} = \mathrm{a_{ can}} + \sigma^2\gamma.  \end{eqnarray}

It is rather straightforward to compute the tilted characteristics $( \sigma^2, \mathrm{a}_\gamma, \boldsymbol{ \Lambda}_\gamma ; \alpha)$ in the finite branching activity case, see Examples \ref{Ex:stat}, \ref{ex:reduced} and \ref{Ex:bbess}. We shall do so only in the case $\gamma = \omega_-$ to lighten the prose and because there is no killing involved in this case. 
\begin{itemize}
\item In Example \ref{Ex:stat} we have $\omega_-=1$ and Assumption \ref{A:omega-} holds. After performing the spinal decomposition with $\gamma=\omega_-$  we still have $ \mathrm{a_{\gamma}^{can}}=0$, $\sigma^2=0$ and the tilted generalized L\'evy measure is merely equal to the initial one, more precisely $ \boldsymbol{\Lambda}_{ \mathrm{half}, \omega_-} = \boldsymbol{\Lambda}_{ \mathrm{half}} \circ \mathrm{Ord}^{-1}$. The tagged branch thus evolves as a standard branch in this model.
\item In Example \ref{ex:reduced} we have $\omega_-=1$ and Assumption \ref{A:omega-} holds. After performing the spinal decomposition with $\gamma=\omega_-$  we still have $ \mathrm{a_{\gamma}^{can}}=-1$, $\sigma^2=0$ and the tilted generalized L\'evy measure is twice the original one, more precisely $ \boldsymbol{\Lambda}_{ \mathrm{two}, \omega_-} = 2 \times \boldsymbol{\Lambda}_{ \mathrm{two}} \circ \mathrm{Ord}^{-1}$. Along the tagged branch, the intensity of splittings is twice that of a standard branch. Undoing the Lamperti transformation we recover the famous spine decomposition of Yule trees, see e.g. \cite[Proposition 5]{CLGharmo}.
\item In Example \ref{Ex:bbess} we have $\omega_- = - \mathrm{a_{can}} - \sqrt{\mathrm{a_{can}}^2 -2}$ and Assumption \ref{A:omega-} holds as long as $\mathrm{a_{can}} < - \sqrt{2}$.  After performing the spinal decomposition with $\gamma=\omega_-$ we still have $ \sigma^2=1$, the tilted generalized L\'evy measure is again $ \boldsymbol{\Lambda}_{ \mathrm{two}, \omega_-} = 2 \times \boldsymbol{\Lambda}_{ \mathrm{two}} \circ \mathrm{Ord}^{-1}$ and the canonical drift gets changed to $ \mathrm{a_\gamma^{ can}} =  - \sqrt{\mathrm{a_{can}}^2 -2}$ using \eqref{eq:acangamma}. In particular, the intensity of splittings along the tagged branch is multiplied by two as in the previous case, and the decoration $X_{\omega_{-}}$ evolves as a Bessel process with dimension $2 \mathrm{a_{\gamma}^{ can}}+2$.
\end{itemize}

Let us now move to examples with an infinite branching activity, starting with Example \ref{ex:brownianheight}. We have $\omega_-=2$ and Assumption \ref{A:omega-} holds. After performing the spinal decomposition with $\gamma=\omega_-$,  we still have $ \mathrm{a_{\gamma}^{can}}=-1$, $\sigma^2=0$ and the tilted generalized L\'evy measure is the sum of the original one $\boldsymbol{\Lambda}$ and a finite measure $\boldsymbol{\Lambda}'$ given by 
$$ \int_{ \mathcal{S}}F\big( \mathrm{e}^{y_0},(  \mathrm{e}^{y_1}, \mathrm{e}^{y_2}, ...)\big) \boldsymbol{\Lambda}' (\dd \mathbf{y})= 2\int_{0}^{1} \mathrm{d}u\ F(u,(1, 0,0,...)).$$ In probabilistic terms, this means that the decoration $X_{\omega_{-}}$ along the tagged branch is a pssMp which evolves as follows: starting from a value $x$, it decreases until $x \cdot (1- \sqrt{U})$ as a standard branch, and then jumps to $x \cdot (1- \sqrt{U}) \cdot V$ where $U,V$ are independent and uniform on $[0,1]$. Afterwards, the evolution iterates the same dynamic.   

\medskip 
The next, and perhaps one of the most important, example is given by the Brownian CRT with mass $1$, see Example \ref{ex:brownian}. It is an example of conservative fragmentation so we have $\omega_-=1$ and Assumption \ref{A:omega-} holds. The spinal decomposition with $\gamma = \omega_-$ has characteristics $\sigma^2=0$, canonical drift $ \mathrm{a_{can}}=0$ and the tilted generalized L\'evy measure is 
$$ \int_{ \mathcal{S}} F\big(   \mathrm{e}^{y_0}, ( \mathrm{e}^{y_1},  \ldots )\big)  \  \boldsymbol{ \Lambda}_{ \mathrm{Bro},\omega_-}( \mathrm{d} y_0, \mathrm{d} \mathbf{y}) =    \sqrt{\frac{2}{\pi}} \int_{0}^{1} F(x,1-x,0,0, ...)\cdot x \cdot  \frac{ \mathrm{d}x}{(x(1-x))^{3/2}}.$$ In particular, the L\'evy--Khintchine exponent of the underlying L\'evy process of the pssMp decoration $X_{\omega_{-}}$ along the tagged branch is 
$$\psi_{\omega_-}(z) = \kappa_{ \mathrm{Bro}}(z+ \omega_-) = -2 \sqrt{2} \cdot \frac{\Gamma(z+ \frac{1}{2})}{\Gamma(z)}.$$ We deduce 
from \cite[Proposition 1]{uribe2009falling} (see \cite[Chapter 5]{kyprianou2022stable} for more general results that we shall use below) that the decoration along the tagged branch $X_{\omega_-}$ has the law of the opposite of a stable subordinator of index $1/2$ starting from $1$, conditioned to visit $0$, and finally killed when hitting $0$. See below for details. Notice that thanks to the conservative and binary properties, the decoration-reproduction $\eta_{\omega_{-}}$ is recovered from $X_{\omega_{-}}$ using \eqref{Eq:reprodbinary}. Interestingly, the reproduction process  bears a close relation to the Poisson--Dirichlet distribution with parameters $(1/2,1/2)$; see \cite[Theorem 1]{uribe2009falling}.
A similar phenomenon occurs in the spinal decomposition of the stable trees of Example \ref{ex:stable}. Again, since the fragmentation is conservative we have $\omega_{-}=1$ and Assumption \ref{A:omega-} holds. As above, the decoration along the tagged branch has the same law as the opposite of a stable subordinator of
index $1-1/\beta$ starting from $1$, conditioned to visit $0$, and finally killed when hitting $0$; see \cite[Proposition 1]{uribe2009falling} and also \cite[Proposition 1]{Mie03}. 
See also below for details. Since this case is not binary, the decoration-reproduction process is more involved, as above it bears close relation with the Poisson--Dirichlet 
distribution with parameter $( 1-\frac{1}{\beta}, 1 - \frac{1}{\beta})$, see \cite{Mie03,uribe2009falling} for details.

\medskip 
Finally, let us consider the family of Examples \ref{ex:stablefamily}. Recall that for $ \texttt{a} \in (0,1]$ and $\texttt{b} \in (0,1/2]$ the characteristics $(\sigma^{2}=0, \mathrm{a}_{ \texttt{a},\texttt{b}}, \boldsymbol{ \Lambda}_{ \texttt{a},\texttt{b}} ; \alpha)$ yield a ssMt for which $\omega_{-}= { \texttt{a}+2\texttt{b}}$ and $\omega_{+} ={ \texttt{a}+2\texttt{b}}+1$ and for which Assumption \ref{A:omega-} always holds. We introduce the notation 
$$ \beta = { \texttt{a}+\texttt{b}}, \quad \mbox{ and } \quad \varrho =  \frac{ \texttt{a}}{ \texttt{a}+\texttt{b}}.$$
Then using  \eqref{eq:kappastable} and Lemma \ref{L:LKtilde} we see that the L\'evy--Khintchine exponent of the pssMp decoration along the $\omega_-$-tagged branch is given by 
$$ \psi_{\omega_-}(z) = - \Gamma(1 +  \texttt{a} - z) \Gamma( \texttt{b} + z) \frac{ \sin(\pi z)}{\pi}.$$ Using \cite[Theorem 5.15]{kyprianou2022stable} we deduce that in the $ \omega_{-}$- spine decomposition of those ssMt, the decoration $X_{\omega_{-}}$ along the tagged branch evolves as a stable L\'evy process with parameters $(\beta, \varrho)$ conditioned to die continuously at $0$. Since the generalized L\'evy measure is conservative \eqref{def:conservativeLM} and binary \eqref{Eq:binaryLM},  the decoration-reproduction process is recovered from $X_{\omega_{-}}$ using \eqref{Eq:reprodbinary}. In the case of the overlay Example \ref{ex:overlay}, the decoration $X_{\omega_{-}}$ in the $\omega_{-}$-spinal decomposition is a so-called ricocheted stable process, introduced by Budd \cite{BudOn} and recently studied in \cite{kyprianou2021double,watson2023growth}. The critical case Example \ref{ex:ADS} is left aside since Assumption \ref{A:omega-} does not hold, see Section \ref{sec:commentsGBP} for a discussion.

\bigskip

We saw in Examples \ref{ex:brownian}, \ref{ex:stable} and \ref{ex:stablefamily} the appearance of conditioned stable processes as the pssMp decoration along the tagged branch of a ssMt. We will actually show that not all such processes can appear in spine decomposition of binary conservative ssMt. Let us first give some background on stable L\'evy processes, their conditioned versions and their relations with hypergeometric L\'evy processes to unify the results. Recall that a  L\'evy process $(\xi_{t})_{t \geq 0}$ is \textbf{stable} with index $\alpha \in (0,2]$ if it satisfies the scaling relation $c^{ -1}\cdot (\xi({c^{\alpha} t}))_{t \geq 0} = (\xi({t}))_{t \geq 0}$ in law. Up to dilation, they can be classified by their index index of similarly $\alpha \in (0,2]$ together with the positivity parameter $\varrho \geq 0$ given by $\varrho = \mathbb{P}( \xi(t) \geq 0)$. Specifically, if 
$$(\alpha, \varrho) \in \big\{ \alpha \in (0,1), \varrho \in [0,1]\big\} \cup \big\{ \alpha =1, \varrho = 1/2\big\} \cup \big\{ \alpha \in (1,2), \varrho \in [ 1- \frac{1}{\alpha}, \frac{1}{\alpha}]\big\},$$ their L\'evy measure is given by 
 \begin{eqnarray} \label{eq:levymeasurestable} \Pi( \mathrm{d}x) =  \frac{\mathrm{d}x}{|x|^{\alpha+1}} \left(  \Gamma(1+ \alpha) \frac{\sin( \pi \alpha \varrho)}{\pi}\mathbf{1}_{ x >0} + \Gamma(1+ \alpha) \frac{\sin( \pi \alpha (1-\varrho))}{\pi} \mathbf{1}_{x <0}\right). \end{eqnarray} The case $\alpha=2$ is the case of Brownian motion (no jumps). We refer to \cite[Chapter VIII]{Ber96} or \cite[Chapter 4]{kyprianou2022stable} for details. A stable L\'evy process naturally gives rise to a pssMp without performing the Lamperti transformation: if $ \xi$ is an $\alpha$-stable L\'evy process starting from $1$, then the \textbf{censored} process  $\xi^{\dagger}$ defined by 
$$ \xi^{\dagger}(t) := \xi(t) \mathbf{1}_{t \leq T_{ \mathbb{R}_{-}}}, \quad \mbox{ with }T_{ \mathbb{R}_{-}} := \inf \left\{ t \geq0 : \xi(t) <0 \right\},$$ is a pssMp. There are other ways to build pssMp from $\xi$ using $h$-transformations, see \cite{caballero2006conditioned}. More precisely, recall that $h : \mathbb{R}_{+} \to \mathbb{R}_{+}$ is a positive harmonic function for $\xi^{\dagger}$ if for each $x>0$ we have 
$${h}(x) = \mathbb{E}_{x}\left( h( \xi^{\dagger}(t)) \mathbf{1}_{t < T_{ \mathbb{R}_{-}}}\right).$$
Given such a harmonic function, it is possible to define a new process $\xi^{h}$ by the formula
$$ \mathbb{P}_{x}( {\xi^{h} \in A}) = \frac{1}{h(x)} \mathbb{P}_{x}\left( {\xi^{\dagger} \in A} \cdot h(\xi^{\dagger}(t)) \mathbf{1}_{t < T_{ \mathbb{R}_{-}}}\right),$$ where $A$ is measurable with respect to $ \mathcal{F}_{t}$ (a priori it is unclear whether the obtained process is conservative or not, but it will be the case in what follows). It turns out that for stable L\'evy processes, any positive harmonic function for $\xi^{\dagger}$ is a linear combination of the two functions $$h^{\uparrow}(x) = x^{\alpha (1-\varrho)}, \quad \mbox{ and }\quad h^{\downarrow}(x) = x^{\alpha (1-\varrho) -1},$$ see \cite{silverstein1980classification}. The two processes $\xi^{\uparrow}$ and $\xi^{\downarrow}$ obtained this way  are respectively called the $(\alpha, \varrho)$-stable L\'evy processes conditioned to survive, resp. to die continuously at $0$, since they can be alternatively obtained by a limiting conditioning procedure associated to their names. See \cite[Chapter 5]{kyprianou2022stable} or \cite{Cha96} for details. Then the three processes $ \xi^{\dagger}, \xi^{\downarrow}, \xi^{\uparrow}$ are positive self-similar Markov processes, for which it is possible to compute the characteristics of the underlying  L\'evy process in the Lamperti transformation \cite{caballero2006conditioned} (beware, those are not stable processes anymore). To present them in a unified way, it is convenient to introduce the formalism of hypergeometric pssMp and L\'evy processes. Recall that a pssMp is called hypergeometric if the underlying L\'evy-Khintchine exponent is of the form 
\begin{equation}\label{eq:hypergeo}
\psi(z)=- \frac{\Gamma(1-\varsigma+\upsilon-z)}{\Gamma(1-\zeta-z)}\cdot \frac{\Gamma(\hat\varsigma+\hat\upsilon+z)}{\Gamma(\hat\varsigma+z)},
\end{equation}
where $(\varsigma,\upsilon,\hat \varsigma,\hat \upsilon)$ belongs to the admissible set of parameters 
$\{ \varsigma \le 1, \; \upsilon \in [0,1), \; \hat \varsigma \ge 0, \; \hat \upsilon \in (0,1) \}$.\footnote{To be precise, the borderline case $\upsilon=0$ is often excluded in the definition. Nonetheless, it can be added to the family of hypergeometric processes thanks to Proposition 4.1 and Theorem 4.4 in \cite{kyprianou2022stable}.}  We refer and use the same notation as in \cite{kuznetsov2013fluctuations} where we replaced $\beta$ by $\varsigma$ and $\gamma$ by $\upsilon$ to avoid conflict in the notation. In particular, it is proved there that the pssMp $ \xi^{\dagger}, \xi^{\downarrow}$ and $ \xi^{\uparrow}$ are hypergeometric  with parameters
\begin{center}
\begin{tabular}{|c||c|c|c|c|}
\hline \rule[0pt]{0pt}{15pt}  
    & \;$\varsigma$ \;  & \; $\upsilon$ \; & \; $\hat \varsigma$ \; & \; $\hat \upsilon$ \;  \\ [0.5ex] \hline \hline 
   \rule[2pt]{0pt}{11pt}   $\xi^\dagger$  & $1-\alpha (1-\varrho)$ & $\alpha \varrho$ & $1-\alpha (1-\varrho)$ & $\alpha (1-\varrho)$ \\ [0.5ex] \hline
   \rule[2pt]{0pt}{11pt}  $\xi^{\uparrow}$  & $1$ & $\alpha\varrho$ & $1$ & $\alpha (1-\varrho)$ \\ [0.5ex] \hline
   \rule[2pt]{0pt}{11pt}  $\xi^{\downarrow}$  & $0$ & $\alpha\varrho$ & $0$ & $\alpha (1-\varrho)$ \\ [0.5ex] \hline	
\end{tabular}
\end{center}

Let us now examine which of the pssMp $ \xi^{\downarrow}, \xi^{\dagger}$ can arise as the decoration along a tagged branch in a conservative and binary ssMt. We shall actually consider the case where the generalized L\'evy measure $ \boldsymbol{ \Lambda}$ is binary, and \textbf{almost conservative} meaning that the only possible loss of mass during splitting event is due to a killing:
 \begin{eqnarray} \boldsymbol{\Lambda} \Big( \Big\{(y_{0},(y_{1}, ... )) :  \sum_{j=0}^{\infty} \e^{y_j}  \ne 1 \Big\} \Big) = \boldsymbol{\Lambda} \big( \left\{ (-\infty, (-\infty, -\infty...) \right\} \big),  \end{eqnarray} note that in this case the decoration-reproduction process $\eta$ is still recovered from the pssMp evolution as in \eqref{Eq:reprodbinary} except that the final jump does not yield to an atom in the decoration-reproduction.  Since we are in the spectrally negative case,  we shall treat separately the subordinator case $\alpha \in (0,1)$ and the case $\alpha \in (1,2)$. In the case $\alpha \in (0,1)$, consider the opposite of a $\alpha$-stable subordinator started from $1$ either killed $\xi^{\dagger}$ when reaching $ \mathbb{R}_{-}$ or conditioned $\xi^{\downarrow}$ to die continuously at $0$. Then we have $\varrho=0$. By \cite[Theorem 5.10 and 5.15]{kyprianou2022stable} the L\'evy-Khintchine exponents in the Lamperti representation is explicit and given by 
 \begin{equation}\label{eq:psi:dagger:downarrow}
 \psi^{\dagger}_{\alpha}(z) = -\frac{\Gamma(1 + z)}{\Gamma(1 - \alpha + z)}, \quad \mbox{ and }\quad \psi_{\alpha}^{\downarrow}(z) = -\frac{\Gamma(\alpha + z)}{\Gamma(z)}.
 \end{equation}
 
  By \cite[Theorem 4.6 (ii)]{kyprianou2022stable}, their  L\'evy measures are,  after a push-forward by $ x \mapsto \mathrm{e}^{x}$, respectively given by 
$$ \pi^{\dagger}_{\alpha}( \mathrm{d}x) = \frac{-1}{\Gamma(-\alpha)}\cdot  \frac{ \mathrm{d}x}{(1-x)^{1+\alpha}} \mathbf{1}_{x \in [0,1]}, \quad \mbox{ and }\quad  \pi^{\downarrow}_{\alpha}( \mathrm{d}x) =\quad \frac{-1}{\Gamma(-\alpha)}\cdot  \frac{ \mathrm{d}x}{x^{1- \alpha}(1-x)^{1+\alpha}}\mathbf{1}_{x \in [0,1]}.$$
Moreover, we infer from \eqref{eq:psi:dagger:downarrow} that  their killing rate are $ \frac{1}{\Gamma(1-\alpha)}$ and $0$ respectively. 
In order to obtain $\xi^\dagger$ and $\xi^{\downarrow}$ as tilted versions of a characteristic quadruplet $(\mathrm{a}_{*}, \sigma_*^2,\boldsymbol{\Lambda_{\ast}}; \alpha_* )$, we must have $\mathrm{a}_{*}^{\text{can}}=\sigma_*=0$ and $\alpha_*=\alpha$ ; where we use the transparent notation $\mathrm{a}_{*}^{\text{can}}$ for the associated canonical drift. Furthermore, without loss of generality we may assume that  $\boldsymbol{\Lambda_{\ast}}$ is a locally  largest generalized L\'evy measure, and by \eqref{eq:pi:dagger:downarrow} and the almost conservative assumption it must be of the form:
 \begin{align}\label{eq:pi:dagger:downarrow}
\int_{ \mathcal{S}} F( \mathrm{e}^{y_{0}}, &( \mathrm{e}^{y_{1}}, \mathrm{e}^{y_{2}}, ...) \boldsymbol{\Lambda_{\ast}}( \mathrm{d}y_{0},  \mathrm{d} \mathbf{y}) \nonumber\\
&:= \int_{1/2}^{1} \frac{-1}{\Gamma(-\alpha)}\cdot  \frac{ \mathrm{d}x}{(x(1-x))^{1+\alpha}} F\big(x,(1-x,0,0, ...)\big)+ \mathrm{k}_*\cdot F\big(0,(0,\dots)\big),
\end{align}
for some constant $\mathrm{k}_*\geq 0$ corresponding to the killing rate. Recalling \eqref{E:massbacksim}, we see that we must tilt the latter by   $\gamma = \alpha+1$ in the case of $\xi^{\dagger}$ and  $\gamma = 2 \alpha$ in the case of $\xi^{\downarrow}$. Finally using the \eqref{eq:kappaspi} combined with the almost conservative assumption,  we infer that we have to take $\mathrm{k}_*= -\psi^{\dagger}_{\alpha}(-\alpha)$ for $\xi^{\dagger}$  and $\mathrm{k}_*= -\psi^{\downarrow}_{\alpha}(1-2\alpha)$ for $\xi^{\downarrow}$. 
Performing the calculation, it turns out that 
$$\psi^{\dagger}_{\alpha}(-\alpha)  =  \psi^{\downarrow}_{\alpha}(1-2 \alpha)  = \frac{ -\Gamma(-\alpha)}{2\Gamma(-2\alpha)}= \frac{-4^\alpha \sqrt{\pi}}{\Gamma(1/2-\alpha)},$$ and the above display is non-positive if and only if $\alpha \in (0,1/2]$. So, to summarize, the processes $\xi^{\dagger}$ and $\xi^{\downarrow}$ can be obtained as the decoration processes along the tagged branch in a ssMt only when $\alpha \in (0,1/2]$, and in this case we can  even use the same ssMt. Also, note that the associated generalized Lévy measure has killing in all cases except for $\alpha=1/2$.

The case when $\alpha \in (1,2)$ is very similar except that some care is needed to deal correctly with compensation. The spectrally negative case correspond to $\varrho = \frac{1}{\alpha}$.  By \cite[Theorem 5.10 and 5.15]{kyprianou2022stable} the L\'evy--Khintchine exponents are now given by 
$$ \psi^{\dagger}_{\alpha}(z) =  \frac{1}{\pi}\left(\Gamma(\alpha - z) \Gamma(1 + z) \sin(\pi (\alpha-z))\right) = (z-\alpha+1) \frac{\Gamma(1+z)}{\Gamma(2-\alpha+z)}$$
$$ \:\:\mbox{ and }\:\: \psi_{\alpha}^{\downarrow}(z) = \frac{1}{\pi}\left(-\Gamma(2 - z) \Gamma(\alpha-1 + z) \sin(\pi z)\right),$$ and the underlying L\'evy measures are, after a push-forward by $ x \mapsto \mathrm{e}^{x}$,  similarly given by 
$$ \pi^{\dagger}_{\alpha}( \mathrm{d}x) = \frac{-\Gamma(\alpha+1)\sin(\alpha \pi)}{\pi} \frac{ \mathrm{d}x}{(1-x)^{1+\alpha}} \quad \mbox{ and }\quad \pi^{\downarrow}_{\alpha}( \mathrm{d}x) = \frac{-\Gamma(\alpha+1)\sin(\alpha \pi)}{\pi} \frac{ \mathrm{d}x}{x^{2-\alpha}(1-x)^{1+\alpha}}.$$ 
As above, those measures are necessarily obtained by starting from a locally largest of the form
\begin{align*} 
\int_{ \mathcal{S}} F( \mathrm{e}^{y_{0}}, &( \mathrm{e}^{y_{1}}, \mathrm{e}^{y_{2}}, ...) \boldsymbol{\Lambda_{\ast}}( \mathrm{d}y_{0},  \mathrm{d} \mathbf{y}) \\
&:= \int_{1/2}^{1} \frac{-\Gamma(\alpha+1)\sin(\alpha \pi)}{\pi}  \frac{ \mathrm{d}x}{(x(1-x))^{1+\alpha}} F(x,(1-x,0,0, ...))+ \mathrm{k}_*\cdot F\big(0,(0,\dots)\big),
\end{align*} after tilting  by $\gamma = \alpha+1$ in the case of $\xi^{\dagger}$ and by $\gamma = 2 \alpha-1$ in the case of $\xi^{\downarrow}$. The drift is then adjusted in the characteristics of the ssMt to match the one obtained in $\psi^{\dagger}_{\alpha}$ or $\psi^{\downarrow}_{\alpha}$. As  above, the only condition to check is that the killing rate is non-negative. Using a formal calculation software, it is possible to obtain a closed formula for the cumulant function in terms of the killing rate, and we get the condition 
$\mathrm{k}_*= \frac{-4^\alpha \sqrt{\pi}}{\Gamma(1/2-\alpha)},$
 in both regimes. The latter is non-negative if and only if $\alpha \in (1,3/2]$.
Notice in particular, the perhaps striking fact that the only stable processes conditioned to die continuously at $0$ that can appear as the decoration $X_{\gamma}$ in a binary conservative ssMt without killing are the $1/2$ and $3/2$-stable spectrally negative cases.

\section*{Comments and bibliographical notes}

The notion of spinal decomposition is one of the most useful and powerful tools in the study of branching structures. It can be traced back at least to Kahane \& Peyri\`ere \cite{kahane1976certaines}, while the first geometric formulation on trees is due to  Chauvin \& Rouault \cite{chauvin1988kpp}. It  has been famously popularized in the 90's by Lyons, Pemantle and Peres \cite{LPP95b,lyons1997simple} and  has found numerous applications, notably in branching random walk theory, see \cite{shi2015branching}. In a context close to ours, these tools have been applied also for growth-fragmentation in \cite{BBCK18} and branching L\'evy processes \cite{bertoin2018biggins}.

The notion of bifurcators has already been introduced by Pitman \& Winkel in the setting of fragmentations trees by \cite{pitman2015regenerative} and by Shi \cite{shi2017growth} for growth-fragmentations processes. Section \ref{sec:bifurcators} build upon these works. See \cite[Section 5]{BBCK18} for a similar use of the tilted characteristics in the growth-fragmentation case.

The \textbf{profile} of a random measured (decorated) tree $ \mathbf{T} = (T, d_{T}, \rho, g, \nu)$ is the push forward of $\nu$ by the distance to the origin  $ x \in T \mapsto d_{T}(\rho,x)$. After size-biaising by the total mass, it is related to the law of $d_{T}(\rho, \rho^{\bullet})$ defined in \eqref{eq:defspinebiaised}. The profile has been studied in details in the fragmentation \cite{haas2004regularity} and the growth-fragmentation \cite{ged2019profile}  cases. In particular, it is proved there that in the case of the harmonic measure, the profile has a continuous density if and only if the self-similarity exponent $\alpha$ is strictly larger than $\omega_{-}$. In the case of the Brownian CRT of Example~\ref{ex:brownian}, the profile is famously linked to Ray-Knight theorems and a similar phenomenon, though more complicated, has been established in the case of the Brownian growth-fragmentation tree of Example \ref{ex:3/2stable}, see  \cite{le2024markov,gall2023stochastic,chapuy2024note}. We wonder whether a similar Markov property of local times holds in a greater generality.

The results in Section \ref{sec:examplespinal} show that the stable family introduced in Example \ref{ex:stablefamily} is uniquely characterized as being ssMt, with the evolution along tagged branches closely related to (a version of) stable Lévy processes. For a related perspective from the growth-fragmentation process, see \cite[Section 5]{BBCK18}. This line of analysis was originally inspired by a computation of Miller and Sheffield \cite[Section 4.6]{MS15}, as well as by the work of Miermont and Schweinsberg \cite{miermont2003self} on fragmentation processes, which respectively distinguished the Brownian fragmentation and the Brownian growth-fragmentation. More recently, these insights have played a role in establishing the universality of polynomial exponents in catalytic equations; see \cite{ConCur25}.

\part{\sc{Discrete Markov branching trees}}\label{part:discrete}
\chapter{Galton--Watson processes with integer types and their scaling limits} \label{chap:6}

 Any scaling limit is intrinsically self-similar.
This important observation, which goes back to Lamperti \cite[Theorem 2]{LampertiSSSP}, can be stated a bit loosely as follows.
Consider a random process $X=(X_t)_{t\geq 0}$  with values in an  Euclidean space. Suppose that there is a positive function $c$ with $\lim_{r\to \infty} c(r)=\infty$, such that the rescaled process $(c(r)^{-1} X_{rt})_{t\geq 0}$ converges in law as $r\to \infty$ to some non-degenerate process $Y=(Y_t)_{t\geq0}$. Then the limiting process $Y$ must be self-similar; furthermore  the function $c$ must be regularly varying at $\infty$.
This leads to the fundamental problem of identifying classes of processes $X$ that are attracted after a proper rescaling to a given self-similar process $Y$. For instance, when $Y$ is a Brownian motion, Donsker's invariance principe enables us to choose for $X$ any centered random walk with finite variance. In this second part of the text, the role of $Y$ will be played by the self-similar Markov trees that were constructed and studied in the first part.

As it was already mentioned in the introductory chapter,  this monograph has been motivated by the breakthrough by Haas and Miermont \cite{HM12} on scaling limits of Markov branching trees. The latter can be viewed as a deep extension   to the branching framework of their scaling limit theorem in  \cite{haas2011self} for non-increasing Markov chains with positive integer values. On the other hand, the assumption in  \cite{haas2011self} that the Markov chain is non-increasing
has been relaxed  afterwards by Bertoin and Kortchemski \cite{BK14}. This suggests that  \cite{HM12}   could be considerably generalized by removing a condition of conservativeness there (which entails that the Markov chains along the branches of the discrete trees are non-increasing). Our goal will be achieved in the present and the next chapters; we also refer to Chapter \ref{chap:applications} for a number of illustrations.

Specifically,  we consider here Galton--Watson processes with integer types.  The system  started from a single particle  is naturally represented in terms of a genealogical tree decorated by the types of the particles and equipped with a measure that counts particles with some weights depending on types. 
We are interested in the asymptotic behavior of the latter when the type of the ancestral particle goes to infinity. 
For this purpose, we first view lineages of particles induced by some selection rule as individuals. Individuals are further equipped with a decoration that records the types of particles along the lineage, so that the whole Galton--Watson process can then be seen as the evolution of a population modeled by a general branching process with a decoration, just as in Chapter \ref{chap:generalBP}. 
The decoration and reproduction along a lineage are Markovian. One naturally expects that requesting the existence of a scaling limit for the latter should be the key to 
a scaling limit theorem for the whole  decorated  Galton--Watson tree, and that  the  limit should be a self-similar  Markov tree of Chapter \ref{chap:generalBP}. 
However this  is not sufficient and we shall need to impose further technical requirements which will be further discussed in the next chapter.

\section{Galton--Watson processes as general branching processes} \label{sec:4.1}
An integer type Galton--Watson process is a branching particle system in discrete time, where particles have  types in $\N=\{1,2,\ldots\}$. We suppose that at the initial time $n=0$, the system consists of a single particle with some given type. The state of the system at an arbitrary  time $n\geq 0$ 
is represented by $\Z(n)=(Z_k(n))_{k\geq 1}$, where $Z_k(n)$ counts the number of particles of type $k$ at time $n$ in the system. We will often view the sequence $\Z(n)$ as an integer-valued measure on $\N$.
The evolution of particles satisfies the branching property, meaning that 
 at each step, every parent particle is replaced by a random family of children particles of different types, independently of the other particles present in the system and also independently of the preceding steps. Moreover,  the law of this family of children particles only depends on the 
type of the parent particle.

\begin{figure}[!h]
 \begin{center}
 \includegraphics[width=13cm]{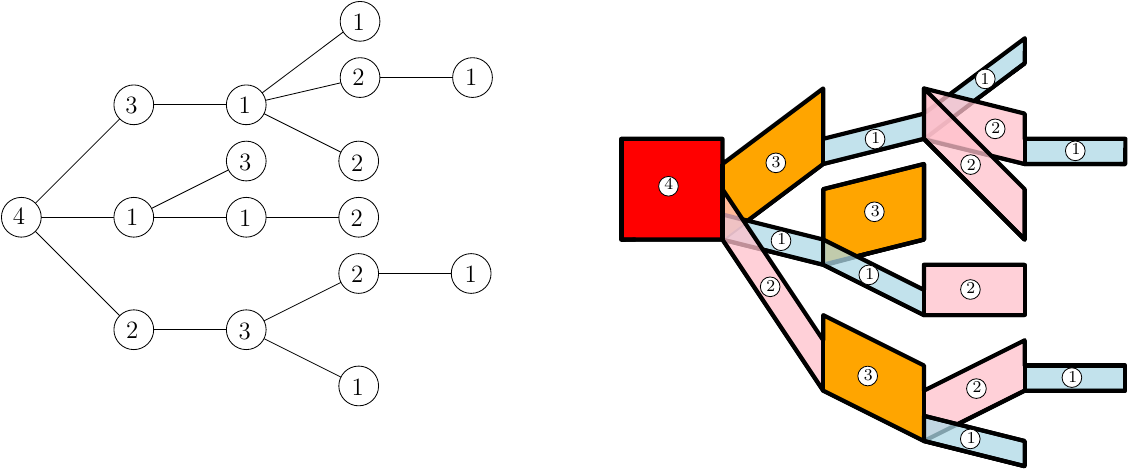}
 \caption{A representation of a Galton--Watson process with integer types as a random decorated tree. \label{fig:GWtoBPlight2}}
 \end{center}
 \end{figure}

Usually,  a Galton--Watson process is conveniently encoded by a discrete genealogical tree, with vertices representing particles and oriented edges connecting parent to children; furthermore each vertex is  tagged with the type of the particle it represents. See the left picture in Figure \ref{fig:GWtoBPlight2} for an illustration. Nonetheless, 
in order to match the framework developed in Chapter \ref{chap:topology}, we will rather use here a slightly different graphical encoding, where now particles are represented by line segments of unit length, and the type of a particle by a constant function on the segment. The genealogical tree $T_\GW$, where the subscript $\GW$ refers of course to Galton--Watson, is then obtained by gluing the left extremity of the segment representing a child particle to the right extremity of that representing its parent particle; 
see  the right picture in Figure \ref{fig:GWtoBPlight2}.  
The root $\rho_{\GW}$ of $T_\GW$ is chosen as the left-extremity of the line segment representing the initial particle. Hence $T_\GW$ is a rooted real tree with a discrete combinatorial structure, such that the distance $d_{T_{\GW}}$ on $T_\GW$  measures  time durations for the particle system. As usual, we abuse notation and simply write $T_{\GW}$ for $( T_{\GW}, d_{T_{\GW}}, \rho_{\GW})$. The types of particles induce a decoration function $g_{\GW}: T_\GW\to \N$; for instance, $Z_k(n)$ can then be recovered as the number of points  having the decoration $k$ at distance $n+1/2$ from the root.  We write 
$$\texttt{T}_{\GW}:=(T_\GW,  g_{\GW})$$
for the resulting decorated real tree, and $\Q^\GW_j$ for its distribution  when $j\in \N$ is the type of the initial particle. We merely ignore measures at this stage for the sake of simplicity, and will equip $T_\GW$ with some natural measures later on.
In the remainder of this section, we assume that the Galton--Watson process eventually becomes extinct, i.e., $\Z(n)=0$ for some sufficiently large $n$, $\Q_j^{\GW}$-a.s., for any $j\geq 1$. As a consequence, $T_{\GW}$ is always bounded a.s.

Our main purpose is to describe simple conditions that ensure the convergence in distribution of a properly rescaled version of the decorated tree $\texttt{T}_{\GW}$ under $\Q^\GW_j$
as $j\to\infty$, where the limit is a self-similar Markov tree with some given characteristic quadruplet. Let us briefly sketch our strategy. We will argue that
$\texttt{T}_{\GW}$ can also be viewed as the decorated genealogical tree  constructed from a family of decoration-reproduction processes representing a general branching process as in Chapter \ref{chap:generalBP}. It is then natural to expect that if these decoration-reproduction processes possess a scaling limit, then this limit gives rise to a self-similar  decoration-reproduction kernel in the sense of Definition \ref{D:SSGBP} and hence  the decorated Galton--Watson tree should have a scaling limit given by a self-similar Markov tree. To make this strategy rigorous, we shall also need some technical conditions which will be introduced later on.

\begin{figure}[!h]
 \begin{center}
 \includegraphics[width=15cm]{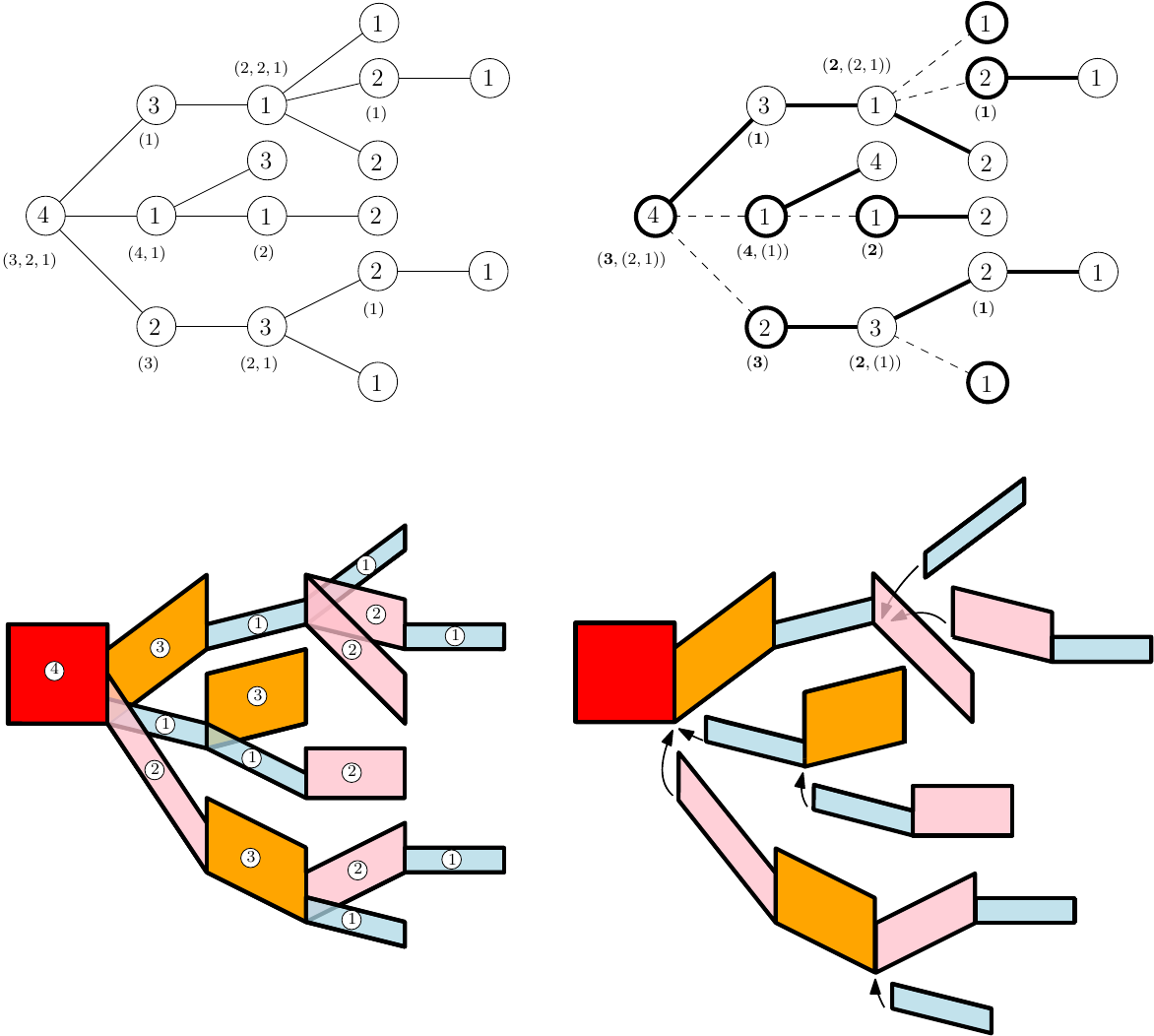}
 \caption{Illustration of the transformation of a Galton--Watson process with integer types into a general branching process and decorated real tree. First line, in the left-hand side, the integer marks assigned to vertices are the types of the particles, and the ranked sequence tagged below a vertex is that of the types of its children. First line, in the right-hand side, the largest type is selected (in bold) for every progeny, and 
 each lineage made of thick lines  is then viewed as a single individual. {The particles that founded lineages are represented by bold circles}. The second line illustrates the transformation into a decorated tree, itself decomposed into decoration-reproduction processes as in Chapter \ref{chap:generalBP}. 
 \label{fig:GWtoBP}}
 \end{center}
 \end{figure}

We shall first explain how a Galton--Watson particle system with integer types can also be viewed as a general branching process\footnote{Of course, Galton--Watson process form the subclass of general branching processes where all individuals  have a unit lifetime and beget children at their death only; however this trivial observation does not bring any insight.} with a decoration, as this has been presented in Chapter~\ref{chap:generalBP}. For this purpose, we stress that the notions of particle and of individual, which are often used indifferently in the literature on branching processes, have here a distinct interpretation. More precisely, an \textbf{individual} will represent a descending lineage of \textbf{particles}, 
the type of an individual being defined as the type of the particle that founds the lineage; see Figure \ref{fig:GWtoBP}.  To distinguish between these two aspects, we will use the name  Galton--Watson process to refer to the evolution of the \textbf{particle system}, and the name \textbf{population model} when we rather consider individuals of the general branching process. The consequences of changing the perspectives from particle system to population model, or vice-versa,
  are easy to analyze, and it is straightforward to switch from one point of view to the other. Essentially, the particle system 
evolves in discrete times, and dealing with individuals living in continuous times  instead implies an obvious interpolation. Last but not least, we stress that the population model representing the Galton--Watson process is not unique; more precisely it depends of some  selection rule that we now introduce\footnote{This is obviously reminiscent of the bifurcation phenomenon identified in Figure \ref{fig:bifurcation} and studied in depth in Chapter \ref{chap:spinal:deco}.}.

We call \textbf{selection rule} a procedure $\varsigma$ for distinguishing at most one atom of any finite integer valued measure $\nu$ on $\N$. For instance, we may choose the largest atom of $\nu$ (as it often will be the case),  the closest atom to some given integer, a randomly selected atom, ..., or even decide to distinguish no atom at all. In the case when one atom of $\nu$ has been selected, say $v_0\geq 1$,  
then we write
$$\varsigma(\nu)=(v_0, v),$$
where   $v=(v_1, \ldots, v_{\ell})$ the sequence of the remaining atoms of $\nu$ (if any), repeated according to their multiplicities and ranked in the
non-increasing  order. 
If no atom has been selected, then we use the same notation $\varsigma(\nu)=(v_0, v)$, now with the conventions  that $v_0=0$ and that $v$ is the ranked sequence of the atoms of $\nu$, repeated according to their multiplicities.
In both cases, the second component $v$ of $\varsigma(\nu)$ belongs to   the subspace  $\U^{\downarrow}$  of non-increasing (possibly empty) sequences in $ \U$ (here, we should merely think of $v\in\U^{\downarrow}$ as a finite family of  types and disregard any genealogical  interpretation). Note also that the identity $\nu = \sum_{i=0}^\ell \delta_{v_i}$  holds in the sense of measures on $\N$ even when $v_0=0$, since $0\not\in \N$. 

Let us now explain how to construct a general branching process and get another perspective on the decorated  tree $\texttt{T}_{\GW}$ that encodes  the particle system.  We start by applying the selection rule  to every progeny in the particle system. That is, we represent generically the family of types in the progeny of a particle by an integer valued measure $\nu$ on $\N$,
 and then the selection rule yields $\varsigma(\nu)=(v_0, v)$.  
Each time, we think of the selected particle with type $v_0$ as the sole, if any, legitimate child of the parent particle, and of all the other particles with types forming the sequence $v$  as illegitimate children particles.\footnote{Formally, the selection rule only specifies the type of the legitimate child. In cases where multiple particles share this type, one is chosen randomly as the legitimate child.} The legitimate child particle continues the lineage and  bears  the same \textbf{name} (label), whereas
each illegitimate child particle founds its own lineage and  bears  a new name (label). More precisely, the name of an illegitimate child is formed as usual by a prefix inherited from the parent and an integer suffix that specifies the rank of this illegitimate child among all illegitimate children begot on the whole lineage. The legitimate lineage terminates when no legitimate child particle has been  selected from the progeny of a legitimate parent particle, and of course in particular when the parent particle has no progeny at all. 

We view legitimate lineages as individuals in a population model, and agree that the \textbf{type} of an individual is simply that of the first particle of the lineage.
The lifetime $z_{u}$ of an individual $u$  is  the number of particles on this lineage, that is the total length of the line segments representing particles of the lineage.
 Since we have assumed that the Galton--Watson process is eventually extinct, $z_{u}$ is finite a.s.
The decoration is induced by the sequence of types of particles on the lineage. 
In particular, the ancestral individual labeled by $\varnothing$ consists of the initial particle and its legitimate descent, 
and its type is  $\chi(\varnothing)=j$ under $\Q^\GW_j$. Its decoration  is defined as the  rcll function $f_{\varnothing}: [0,z_{\varnothing}]\to \mathbb{Z}_+= \N\cup\{0\}$ such that for every integer $0\leq n < z_{\varnothing}$ and $t\in [n, n+1)$, the value $f_{\varnothing}(t)$ is the type of the legitimate descent  at time $n$ of the ancestral particle and $f_\varnothing(z_{\varnothing})=0$.  Finally, the reproduction process $\eta_{\varnothing}$ of the ancestor for the population model can be expressed as 
\begin{equation}\label{E:GWeta}
\eta_{\varnothing} \coloneqq \sum_{n=1}^{z_{\varnothing}} \sum_{\ell}\delta_{(n, v_\ell^{(n)})}, 
\end{equation}
where  $v^{(n)}=(v^{(n)}_{\ell})\in \U^{\downarrow}$ stands for the ranked sequence of the types of illegitimate children particles which the legitimate descent begets  at time $n$.
For example, for the realization of the Galton--Watson process illustrated in Figure \ref{fig:GWtoBP}, the reproduction process 
$\eta_{\varnothing}$ has four atoms at $(1,2), (1,1), (3,2)$ and $(3,1)$. 

Recall that we interpret an atom $(n, v_\ell^{(n)})$   in \eqref{E:GWeta} as   the arrival at time $n$ in the system of an illegitimate particle of type $v^{(n)}_{\ell}$ from the ancestral lineage, that is,  as the  birth of a child individual  with type $v^{(n)}_{\ell}$. More precisely, this child individual is identified in turn  with the legitimate lineage descending from that particle; the decoration is given as preceedingly by the sequence of the types of the particles on this lineage.  As usual, each individual of the first generation receives a different integer label compatible with  the co-lexicographical order as in Section~\ref{sec:2.1}, i.e.\ is indexed by a vertex at the first generation of the Ulam tree, and all the particles in a legitimate lineage then bear the same label. 
The family formed by all the $v^{(n)}_{\ell}$ in \eqref{E:GWeta} is that of the types $\chi(i)$  of the  individuals $i\in \N$ at the first generation (recall that the type of a lineage has been defined as that of the particle which founded this lineage and that the type $0$ is assigned to fictitious individuals). 
And so on, and so forth, generation after generation.
 It should be plain from the branching property of Galton--Watson processes that the resulting population model is a general branching Markov process\footnote{Actually, types are positive real numbers in Chapter 2, whereas here we only deal with integer types. Rather than extending artificially the Galton--Watson process to real types, we simply overlook this slight difference. }  endowed with a decoration in the sense of Section \ref{sec:2.1}. 
  We write $\P^\GW_j$ for the  law of the family of decoration-reproduction processes $(f_u, \eta_u)_{u\in \U}$ when
 the system starts at time $n=0$ from a single particle with type $j$. We stress that the notation could be sightly misleading, because this law  also depends on the choice of the selection rule $\varsigma$, even though the latter does not appear in the notation for the sake of simplicity.

After this long qualitative discussion, we turn our attention to quantitative elements that are needed for our analysis.
The most basic notion in this setting is the \textbf{mean matrix} $\boldsymbol{m}_{\GW}$ of the Galton--Watson process, see e.g. \cite[Section V.2]{AN72}, which is defined by
  \begin{eqnarray} \label{eq:meanrepoductionmatrix}{m}_{\GW}(i,j)\coloneqq \E_i^{\GW}\big(Z_1(j)\big), \qquad i,j\in \N.  \end{eqnarray}
At this point, this matrix may have infinite entries, although this situation will be excluded by forthcoming assumptions.  We also introduce  the  kernel $(\boldsymbol{\pi}^{\GW}_j)_{j\geq 1}$ of probability measures on $\mathbb{Z}_+\times \U^{\downarrow}$ 
that stems from the reproduction laws of the Galton--Watson process
and the selection rule.  We set
 \begin{equation} \label{Eq:reprodkernelGW}
 \boldsymbol{\pi}^{\GW}_j\big(v_0,v\big)\coloneqq \Q^\GW_j\big( \varsigma(\Z(1))=(v_0,v) \big),\qquad \text{for all }v_0\geq 0 \text{ and } v\in \U^{\downarrow}.
 \end{equation}
We will often refer to $(\boldsymbol{\pi}^{\GW}_j)_{j\geq 1}$ as the \textbf{reproduction  kernel} of the Galton--Watson process, omitting  the role of the selection rule for the sake of simplicity.
Needless to say, the reproduction  kernel $(\boldsymbol{\pi}^{\GW}_j)_{j\geq 1}$ determines all the measures  $ (\mathbb{Q}_{j}^{\GW})_{j \geq 1}$ and the mean matrix $\boldsymbol{m}_{\GW}$. 

The formal Markovian description of  the law $P^{\GW}_j$ of the decoration-reproduction process $(f_{\varnothing},\eta_{\varnothing})$ of the ancestor for the population model under $\Q^\GW_j$, can now be stated as follows.

\begin{proposition} \label{P:reprodkernGW} Fix arbitrary $k \geq 0$, $v^{(k+1)}\in \U^{\downarrow}$ and $(j^{(n)}, v^{(n)}) \in \N\times \U^{\downarrow}$  for each $1\leq n \leq k$.
 For every $j\geq 1$, the probability under $\Q^\GW_j$ that $z=k+1$, $f_{\varnothing}(t)= j^{(n)}$ for all $t\in (n, n+1)$ and $1\leq n \leq k$, and 
 $\eta_{\varnothing}$ is given by the right-hand side of \eqref{E:GWeta}, equals
 $$\boldsymbol{\pi}^{\GW}_j(j^{(1)}, v^{(1)}) \times \boldsymbol{\pi}^{\GW}_{j^{(1)}}(j^{(2)}, v^{(2)}) \times \cdots \times \boldsymbol{\pi}^{\GW}_{j^{(k-1)}}(j^{(k)}, v^{(k)}) \times \boldsymbol{\pi}^{\GW}_{j^{(k)}}(0,v^{(k+1)}). $$
 \end{proposition}
\noindent We now conclude this section with the  important observation analog to Remark \ref{rek:PQ}. 

\begin{remark}[$\mathbb{P}^{\GW}$ and $\mathbb{Q}^{\GW}$] 
As in the continuous setting, we made a distinction between the underlying Galton--Watson trees (and its associated decorated trees), whose law is denoted by $\mathbb{Q}^{\GW}_{j}$ when started from $j$, and the induced decoration-reproduction labeling of Ulam's tree after choosing a selection rule, whose law is denoted by $\mathbb{P}^{\GW}_{j}$. From the Galton--Watson point of view, the decorated genealogical tree $\texttt{T}_{\GW}$ is obtained by gluing the left extremity of the segment representing a child particle to the right extremity of that representing its parent particle, generation after generation. It should be evident that $\texttt{T}_{\GW}$ is isomorphic to the decorated genealogical tree that arises when we adopt the population model perspective. Indeed, the latter simply corresponds to choosing a different order for the gluing, specified by the selection rule $\varsigma$. Of course, in the end, the result of the gluing is the same up to an isomorphism; see once again Figure \ref{fig:GWtoBP}. In this regard, note that the Galton--Watson particle system is actually identical to the population model under the trivial selection rule, where no atom is ever distinguished. As in the continuous setting, this subtlety may be overlooked, but in the following, any random variable defined using a selection rule will be considered under $\mathbb{P}^{\GW}_{j}$, whereas any geometric event depending only on the labeled tree structure will be considered under $\mathbb{Q}^{\GW}_{j}$.
\end{remark}

\section{Self-similar Markov trees as scaling limits}\label{sec:4.2}
In this section, we present three versions of the invariance principle (or scaling limit) that lie at the heart of this monograph; they will be established in Section~\ref{sec:invprinc}.  The first, Theorem~\ref{T:mainuncondw}  below, is a weaker version discarding  
measures, whereas the last two, Theorem \ref{T:mainunconds-length} and \ref{T:mainunconds-mass}  below, are stronger versions taking measures into account. Roughly speaking, they state that under appropriate assumptions,  the genealogical tree of the Galton--Watson process started from a single particle with type $n$, its decoration which assigns types to edges on the tree, and eventually  a measure that counts particles weighted by a function of their types, have a joint scaling limit in distribution as $n\to \infty$. The normalization is induced by the linear map 
\begin{equation}\label{eq:linmap} \underbrace{\R_+}_{ \mathrm{indexing}} \times \underbrace{(0,\infty)}_{ \mathrm{decoration}}\to \R_+\times (0,\infty), \qquad (s,y)\mapsto (n^{-\alpha}s, n^{-1}y), \end{equation}
and the 
 limit is given by one of the self-similar Markov trees which have been constructed in Chapter  \ref{chap:generalBP}, endowed with one of the weighted length measures  $\uplambda^\gamma$ (cf. Section \ref{sec:2.3.1}), or with the harmonic measure  $\upmu$ (cf. Section \ref{sec:cramer:cond:1}),
depending on the asymptotic behavior of the weigh function.
In order to help the reader, let us describe the content of the next two chapters. The rest of the current chapter will be devoted to the proof of our main scaling limits results under the forthcoming Assumptions \ref{assum:main:gen:deco:repro} (scaling limit of decoration-reproduction processes) and \ref{A:all:existance:PHI} (existence of superharmonic function with prescribed growth). Although Assumption \ref{assum:main:gen:deco:repro} seems natural, its verification may not be obvious to the reader. Similarly, Assumption \ref{A:all:existance:PHI} may look like an ad-hoc condition. Indeed, these two assumptions are tailor made for the proof of our scaling limit results. We shall provided more natural analytical conditions in Chapter \ref{C:scaling} which will enforce Assumptions \ref{assum:main:gen:deco:repro} and \ref{A:all:existance:PHI}. Notice however an important qualitative difference: one the one hand, Assumption \ref{assum:main:gen:deco:repro} is only \textit{asymptotic} in the type, whereas Assumption \ref{A:all:existance:PHI} assumes a \textit{global} control and thus depends very much of the behavior of small types particles. A similar dichotomy will be present in the assumptions of Chapter \ref{C:scaling}.

\begin{figure}[!h]
 \begin{center}
 \includegraphics[width=15cm]{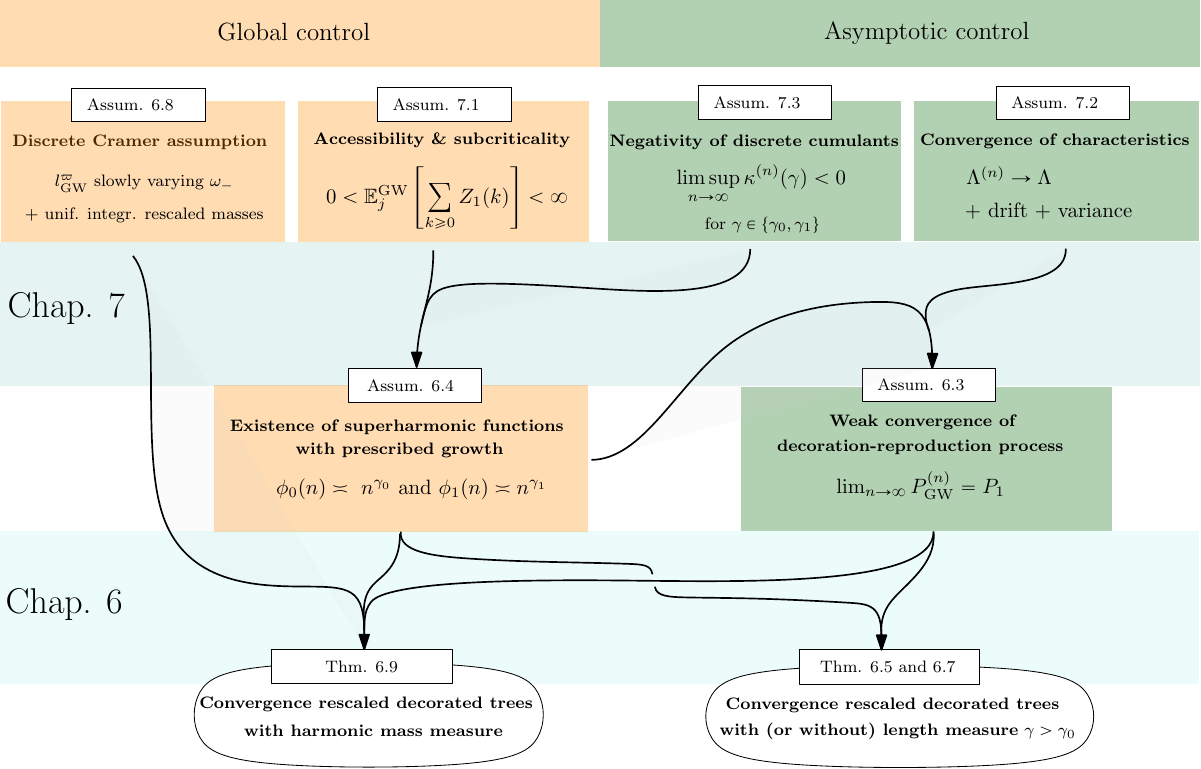}
 \caption{Diagram of the implications between assumptions of Chapter \ref{chap:6} and Chapter \ref{C:scaling}. The rest of the current chapter is devoted to the proof of the main scaling limits results (Theorems \ref{T:mainuncondw}, \ref{T:mainunconds-length} and \ref{T:mainunconds-mass}) from Assumptions \ref{assum:main:gen:deco:repro}, \ref{A:all:existance:PHI} and \ref{A:hGWrv}. The next chapter provides more natural Assumptions \ref{A:all:GW}, \ref{A:BK} and \ref{A:typex} which will enforce the latter. \label{fig:diagram}}
 \end{center}
 \end{figure}

 Throughout the rest of this chapter,  we fix some  characteristic quadruplet  $(\sigma^{2}, \mathrm a,  \boldsymbol{\Lambda}; \alpha)$ and assume without further mentions that its cumulant $\kappa$ satisfies  Assumption \ref{A:gamma0}.  Recall that we also write $P_1$ for the law of the decoration-reproduction process $(X,\eta)$ which is constructed in Section \ref{sec:2.2} by the Lamperti transformation.

Let us next present the assumptions on the Galton--Watson process that we shall need to state our scaling limit theorems.
 Our first assumption should not come as a surprise. It requests that the law $P_1$ arises as the scaling limit 
 of the laws $P^{\GW}_n$ of  the decoration-reproduction processes for the population model when the Galton--Watson process starts from a particle with a large type $n$ (recall Proposition \ref{P:reprodkernGW}). 
 Specifically, the linear map \eqref{eq:linmap}  induces a scaling transformation on decoration-reproduction pairs, 
 \begin{equation}\label{Eq:scalefeta}(f,\eta)\mapsto ( f^{(n)}, \eta^{(n)}),
 \end{equation}
where 
  $f^{(n)}(\cdot )\coloneqq  n^{-1} f(n^{\alpha}\cdot )$ and $\eta^{(n)}$ is the push-forward of $\eta$ by \eqref{eq:linmap}. 
  Here the decoration $f: [0,z]\to \R_+$ is a generic rcll function, the reproduction $\eta=\eta(\d t, \d x)$ a generic integer-valued Radon measure on $[0,z]\times (0,\infty)$, and the lifetime $z\geq 0$ is arbitrary. 
We write $P^{(n)}_{\GW}$ for the image of $P_n^{\GW}$ by this scaling transformation.  
Last but not least, the space of decoration-reproduction pairs  is naturally equipped with the product topology, where  we use the Skorokhod's distance between rcll functions
on compact intervals 
 for the decoration component, and the vague topology for Radon measures on $\R_+\times (0,\infty)$  for the reproduction component.
Formally,  we endow the space  
$$\mathbb D^{\dagger}\coloneqq \bigsqcup_{z\geq 0}\mathbb D([0,z],\R_+)$$  of  rcll functions $\omega: [0,z]\to \R_+$ with a finite lifetime $z$, with the Skorokhod's topology; see \cite[Chapter VI]{JS03} for background. So a sequence $(\omega_n)_{n\geq 1}$ in $\mathbb D^{\dagger}$ converges to some $\omega \in \mathbb D^{\dagger}$ if and only if  for every $n\geq 1$, we can find an increasing bijection 
$\beta_n: [0,z]\to [0,z_n]$ such that, as $n\to \infty$, the bijection $\beta_n$ converges to the identity function on $[0,z]$ and $\omega_n\circ \beta_n$ to $\omega$, both uniformly on $[0,z]$. Plainly, this forces $\lim_{n\to \infty} z_n=z$.
 
 \begin{assumption}[Scaling limit of decoration-reproduction process] \label{assum:main:gen:deco:repro} There is the weak convergence of probability measures
\begin{equation}\label{eq:assum:main:gen:deco:repro}
\lim_{n\to \infty} P^{(n)}_{\GW} = P_1.
\end{equation}
\end{assumption}
See also the forthcoming Lemma \ref{lem:convreprodecolight} for a weaker hypothesis. Due to the lack of continuity of the function that maps a family of decoration-reproduction pairs to a decorated tree, 
the existence of a scaling limit for the Galton--Watson decoration-reproduction process is not sufficient to ensure the convergence of the rescaled decorated genealogical trees, and we will also need some uniform control for the reproduction kernel of the Galton--Watson process. 

Recall the notation $\boldsymbol{m}_{\GW}$ for the mean matrix  of the Galton--Watson process. 
Taking the point of view of linear operators,  we define every function $\varphi:  \mathbb{N} \to [0,\infty]$ the function $\boldsymbol{m}_{\GW} \varphi$ by 
  \begin{eqnarray} \boldsymbol{m}_{\GW} \varphi(i) \coloneqq \sum_{j\geq 1} {m}_{\GW}(i,j) \varphi(j)=\mathbb{E}^{\GW}_{i}\Big(\sum_{j \geq 1} \varphi(j) Z_{j}(1)\Big), \qquad i\geq 1. \label{def:pinbold}\end{eqnarray}
Then the generator $ \boldsymbol{G}_{\GW}\coloneqq  \boldsymbol{m}_{\GW} -\boldsymbol{I}$ is the operator given for every nonnegative and finite function $\varphi$ by
\begin{equation} \label{Eq:Laplacediscret} \boldsymbol{G}_{\GW} \varphi \coloneqq   \boldsymbol{m}_{\GW} \varphi-\varphi.
\end{equation} 
We call $\varphi$ \textbf{superharmonic} for the Galton--Watson process  on a subset $C\subset \mathbb{N}$ if $\boldsymbol{G}_{\GW} \varphi(n) \leq 0$  for every $n\in C$.
When $C=\mathbb{N}$, we simply say that $\varphi$ is superharmonic and add the adverb \textbf{strictly} when the preceding inequality is strict. 

Finally, we write $\varphi(n)\asymp n^{\gamma}$ for some $\gamma >0$ when the ration  $\varphi(n)/n^{\gamma}$ remains bounded away from $0$ and $\infty$ for $n\in \N$.
Recall also from Assumption \ref{A:gamma0} that the cumulant $\kappa$ is subcritical and that $\gamma_0>0$ has been chosen such that $\kappa(\gamma_0)<0$. 

  \begin{assumption}[Existence of superharmonic functions with prescribed growth] \label{A:all:existance:PHI} 
There exist   $\phi_0,\phi_1:\mathbb{N}\to \mathbb{R}_+^*$  and $\gamma_1\in(0,\gamma_0)$ such that
\begin{itemize}
\item $\phi_0$   is  strictly superharmonic for the Galton--Watson process,    $\phi_0(n) \asymp n^{\gamma_0}$, and
$$ \limsup \limits_{n\to \infty} n^{\alpha-\gamma_0}\boldsymbol{G}_{\GW} \phi_0(n)<0;$$

\item $\phi_1$ is superharmonic for the Galton--Watson process and  $\phi_1(n) \asymp n^{\gamma_1}$.

\end{itemize}
 \end{assumption}
The strictly superharmonic function $\phi_0$ of Assumption \ref{A:all:existance:PHI}  serves a similar purpose as the function in Equation \eqref{E: AssumpH1} in the continuous setting. We will also see that Assumption \ref{A:all:existance:PHI}  forces the a.s. extinction of the Galton--Watson system, which was part of our initial assumptions.

In order to state our first scaling limit theorem,  we recall the notation $(\mathbb{T}, \d_{\mathbb{T}})$ for the Polish space of equivalent classes of decorated trees,  see Theorem \ref{theo:Polish} and Corollary~\ref{C:dmrtclosed}, and that we see $\texttt{T}$ as a random element of $\mathbb{T}$.  We still have to specify the notion of scaling for decorated trees. The linear transformation \eqref{eq:linmap} suggests that   for every
 $n\geq 1$,  the  rescaled version of a decorated tree $(T_\GW, d_{T_\GW}, \rho_{\GW}, g_\GW)$  should be  given by
 keeping the same space and root, and by normalizing distance  with a factor $1/n^{\alpha}$ and the decoration  with a factor $1/n$. 

\begin{theorem}[Without measures] \label{T:mainuncondw} Let  Assumptions  \ref{assum:main:gen:deco:repro} and  \ref{A:all:existance:PHI}  hold.
As $n\to \infty$, the sequence  of the distributions  of the rescaled decorated genealogical  trees 
 $$\big(T_{\GW}, n^{-\alpha} d_{T_{\GW}}, \rho_{\GW}, n^{-1} g_{\GW}\big), \quad \text{ under }\Q_n^\GW,$$ 
 converges to the law  $\Q_1$ of 
the self-similar Markov tree $\normalfont{\texttt{T}}$ with characteristic quadruplet $(\sigma^{2}, \mathrm a,  \boldsymbol{\Lambda}; \alpha)$.
\end{theorem}

\begin{remark}\label{R:roleselection} It may be useful to briefly comment on the role of the selection rule $\varsigma$ which might be overlooked otherwise, as we chose to drop it from the notation.  
First, whenever it exists, a scaling limit for the decorated genealogical tree ${\texttt{T}}_{\GW}$ of the Galton--Watson process shall of course not depend on any selection rule.
Observe in this direction that although  Assumption \ref{A:all:existance:PHI} is independent of $\varsigma$, this is not the case for Assumption  \ref{assum:main:gen:deco:repro}.
There are examples for which Assumption  \ref{assum:main:gen:deco:repro} is fulfilled for two different selection rules, say $\varsigma$ and $\varsigma_{\Yleft}$, and two 
different characteristic quadruplets, say  $(\sigma^2, {\mathrm a},  \boldsymbol{\Lambda};\alpha)$ and $ (\sigma_{\Yleft}^2, {\mathrm a}_{\Yleft}, { \boldsymbol{\Lambda}}_{\Yleft};\alpha_{\Yleft})$.
Nonetheless, since the limiting laws of self-similar Markov trees must be the same, i.e. $\Q_1=\Q^{\Yleft}_1$, the two characteristic quadruplets must be bifurcators of one another. See Section \ref{sec:bifurcators},
and notably Theorem \ref{theo:bif} there.
\end{remark}

We next turn our attention to stronger scaling limit results involving measures, and in this direction, we naturally endow $T_\GW$  with a weighted length measure. Specifically, we write $\uplambda_{\GW}$ for the length (i.e. Lebesgue) measure
on $T_{\GW}$, which is finite a.s. since the Galton--Watson process is eventually extinct. Given a weight function $\varpi: \mathbb{N} \to \mathbb{R}_+$, we define $\uplambda^{\varpi}_{\GW}\coloneqq \varpi\circ g_{\GW} \cdot \uplambda_{\GW}$, i.e. 
 $\uplambda^{\varpi}_{\GW}$ is absolutely continuous with respect to $\uplambda_{\GW}$
 with Radon-Nikodym derivative (or density) $\varpi\circ g_{\GW}$. Roughly speaking, $\uplambda^{\varpi}_{\GW}$ essentially amounts to counting particles weighted by a function of their types; in particular the total mass of $\uplambda^{\varpi}_{\GW}$ is  given by $\sum_k \varpi(k)\bar Z_k$,
where $\bar Z_k=\sum_{n=0}^{\infty}Z_k(n)$ denotes the total number of particles of type $k$ that ever existed in the system. 

We will obtain different scaling limit theorems depending on the asymptotic behavior of the weight function. Taking Assumptions  \ref{assum:main:gen:deco:repro} and  \ref{A:all:existance:PHI} for granted, we first present the case when the limiting measure is one of the weighted length measure $\uplambda^{\gamma}= g^{\gamma-\alpha}\cdot \uplambda$ for some $\gamma>\gamma_0$. Before stating the result, a quick discussion about  length measures and rescaling may be appropriate. Following Theorem \ref{T:mainuncondw}, the length measure of the rescaled tree $\big(T_{\GW}, n^{-\alpha} d_{T_{\GW}}\big)$
is $n^{-\alpha}\cdot  \uplambda_{\GW}.$ Heuristically, since the length measure of a segment is given by the distance between its extremities, 
we may hope from Theorem \ref{T:mainuncondw} that for $n$ large, $n^{-\alpha}\cdot  \uplambda_{\GW}$ should be close, in some loose sense\footnote{This is only informal, since in many cases, the length measure on a real tree is not a Radon measure.
Beware also that the assertion can be wrong even in situations where total masses remain bounded.}, to the length measure $\uplambda$ on the limiting self-similar Markov tree. 
Since in turn the rescaled decoration $n^{-1}g_{\GW}$  is close to $g$, if we further assume that the weight function $\varpi$  is regularly varying with exponent 
$\gamma-\alpha$, meaning the  weight function  satisfies 
$$\frac{\varpi(\lfloor nx\rfloor)}{\varpi(n)} \sim x^{\gamma-\alpha}\qquad \text{ as $n\to \infty$,   for any $x>0$,}$$ 
 then we should have, at least informally, that
$$\frac{\uplambda_{\GW}^{\varpi}}{n^\alpha\varpi(n)}\approx g^{\gamma-\alpha} \cdot \uplambda =\uplambda^{\gamma}.$$
 In this direction, it may be worth stressing that, although weighted lengths measures on a decorated tree are functionals of the latter, it is easily seen that these functionals are everywhere discontinuous. As a consequence, we cannot simply deduce measured versions  of Theorem  \ref{T:mainuncondw} from an application of the continuous mapping theorem\footnote{In the converse direction,  the projection $(\texttt{T}, \upmu)\mapsto (\texttt{T},0)$ on $\mathbb{T}_m$ is continuous for the distance $\mathrm{d}_{\mathbb{T}_m}$,
and  Theorem~\ref{T:mainuncondw} hence follows from a stronger measure version.}.

Our next result confirms that the intuition gained above is actually correct.

\begin{theorem}[With a weighted length measure]\label{T:mainunconds-length}
Suppose that Assumptions  \ref{assum:main:gen:deco:repro} and \ref{A:all:existance:PHI}  hold, and that $\varpi$ is regularly varying at $\infty$ with exponent $\gamma-\alpha$ for some $\gamma>\gamma_0$. In the notation above,  the sequence of the distributions  of the rescaled measured decorated trees 
 $$\Big(T_{\GW}, n^{-\alpha} d_{T_{\GW}}, \rho_{\GW}, n^{-1} g_{\GW}, n^{-\alpha}\varpi(n)^{-1} \uplambda_{\GW}^{\varpi}\Big), \quad \text{ under }\Q_n^\GW,$$ 
 converges as $n\to \infty$,  to the law  of 
the measured self-similar Markov tree $(\normalfont{\texttt{T}}, \uplambda^\gamma)$ with characteristic quadruplet $(\sigma^{2}, \mathrm a,  \boldsymbol{\Lambda}; \alpha)$ under $\Q_1$.\end{theorem}

Finally, we present a scaling limit theorem involving the harmonic measure. This case is more challenging, and we need some additional  requirements. The first is that the characteristic quadruplet and its cumulant $\kappa$ verify Assumption \ref{A:omega-}, as we have to define the harmonic measure for the self-similar Markov tree. Recall in particular that $\omega_-$ is the first  root of cumulant $\kappa$.
The second should be viewed as a version of Assumption \ref{A:omega-}
for the Galton--Watson process; it is given in terms of  the mean of the total mass of $\uplambda^{\varpi}_{\GW}$, 
 \begin{eqnarray} \label{def:lA} l^{\varpi}_{\GW}(n) \coloneqq  \mathbb{E}^{\GW}_{n} \Big(  \uplambda^{\varpi}_{\GW}\big( T_{\GW}\big)\Big)=\mathbb{E}^{\GW}_n \Big( \sum_{j \geq 0} \sum_{k \geq 1}\varpi(k)Z_{k}(j)  \Big), \quad n\geq 1.\end{eqnarray} 
 \textit{A priori}, the function $l^{ \varpi}_{\GW}$ can take infinite values; however we will see  in Proposition \ref{prop:massepresqueomega-}  that 
this is not the case provided that the weight $\varpi$ does not grow too rapidly. 

\begin{assumption}[Discrete Cramer condition]\label{A:hGWrv} We suppose the weight function $\varpi$ is such that:
\begin{enumerate}
\item[(i)] the function $l^{ \varpi}_\GW$  is regularly varying with exponent $\omega_-$; 
\item[(ii)]   the family $\big( l^{ \varpi}_\GW(n)^{-1} \cdot  \uplambda^{\varpi}_{\GW}(T_\GW)\ \mbox{ under } \mathbb{Q}^\GW_n\big)_{n\geq 1}$ is uniformly integrable.
\end{enumerate}
\end{assumption}
 An application of the branching property for the Galton--Watson process enables us to interpret  $l^{\varpi}_{\GW}$ as the  potential\footnote{Of course, if the spectral radius of the mean matrix is strictly less than $1$ (which we do not impose)
then the series $\sum_{k= 0}^{\infty} \boldsymbol{m}_{\GW}^k$ converges; its sum is given by the inverse matrix of $\boldsymbol{I}-\boldsymbol{m}_{\GW}= -\boldsymbol{G}_{\GW}$ and is 
known as the Green matrix. We then have $l^{ \varpi}_\GW= \left(\boldsymbol{I}-\boldsymbol{m}_{\GW}\right)^{-1} \varpi$, or equivalently
$\boldsymbol{G}_{\GW} l^{ \varpi}_\GW=-\varpi$ for any weight functions with finite support.
This can sometimes be used to compute $l_\GW$ explicitly, for instance in the case when for every $j\geq 1$, a particle with type $j$ has no offsprings with a greater type, since then the mean  matrix
$\boldsymbol{m}_\GW$ is triangular. } 
of  $\varpi$ for 
  the transition operator \eqref {def:pinbold}. Namely, writing $\boldsymbol{m}_{\GW}^k$ for the $k$-th power of the mean matrix 
  $\boldsymbol{m}_{\GW}$, we have
 $$l^{\varpi}_{\GW}=\sum_{k= 0}^{\infty} \boldsymbol{m}_{\GW}^k\varpi.$$
Clearly,  there is the equation
  \begin{equation} \label{eq:ellsuperh}
   l_{\GW}^{ \varpi} =  \boldsymbol{m}_{\GW} l_{\GW}^{ \varpi} +{ \varpi}.
   \end{equation}
   In particular, $l_{\GW}^{ \varpi}$ is superharmonic whenever it is finite everywhere. Observe also that, since the Galton--Watson process is eventually extinct, there are no non-trivial mean harmonic functions and the construction of Section \ref{sec:cramer:cond:1} does not apply.

We can now state our third and last scaling limit theorem.

\begin{theorem}[With the harmonic measure]\label{T:mainunconds-mass} 
Suppose that Assumptions \ref{A:omega-},  \ref{assum:main:gen:deco:repro} and \ref{A:all:existance:PHI}  hold, and that the weight function $\varpi$ verifies Assumption \ref{A:hGWrv}.
The sequence of the distributions  of the rescaled measured decorated trees 
 $$ \left(T_{\GW}, n^{-\alpha} d_{T_{\GW}}, \rho_{\GW}, n^{-1} g_{\GW},  l^{\varpi}_\GW(n)^{-1} \,  \uplambda^{\varpi}_{\GW}\right), \quad \text{ under }\Q_n^\GW,$$ 
 converges as $n\to \infty$,  to the law   under $\Q_1$ of 
the self-similar Markov tree $(\normalfont{\texttt{T}}, \upmu)$ with characteristic quadruplet $(\sigma^{2}, \mathrm a,  \boldsymbol{\Lambda}; \alpha)$
endowed with its harmonic measure. 
\end{theorem}

These invariance principles will be proved in Section \ref{sec:invprinc}; beforehand we need to establish some uniform controls in the discrete setting, which motivates the next technical section.

 \section{Uniform controls}\label{sect:uniform:controls}
 Our main purpose here is to derive some uniform controls for the Galton--Watson process from 
Assumption \ref{A:all:existance:PHI}. We start by establishing that the system becomes eventually extinct, and by controlling the height and the decoration of the genealogical tree.  The following bounds are the discrete counterparts of those in Theorem \ref{T:CMJ}.

 \begin{proposition}\label{L:Aspandii} Let  Assumption  \ref{A:all:existance:PHI}  hold. Then we have
 \begin{equation}\label{eq:L:Aspandii:1}
 \sup_{x >0} \sup_{n \geq 1}  \left(x^{\gamma_0}   \Q_n^\GW\left(\sup_{T_{\GW}} g_{\GW} \geq x n\right)\right) < \infty,
 \end{equation}
and 
 \begin{equation}\label{eq:L:Aspandii:2} \sup_{x > 0}\sup_{n\geq 1}   \left(x^{ \gamma_0/\alpha}  \Q_n^\GW\big(  \mathrm{Height}(T_{\GW}) \geq x n^{\alpha}\big) \right)< \infty.
 \end{equation}
 In particular,  the Galton--Watson process becomes eventually extinct; its  genealogical  tree  $T_{\GW}$ and also a fortiori its decoration $g_{\GW}$ are bounded, $\Q_n^\GW$-a.s.
 for all $n\geq 1$.
\end{proposition}

\begin{proof} Let   $\phi_0$  given by Assumption  \ref{A:all:existance:PHI} and let $C>0$ large enough so that 
$$C^{-1} n^{\gamma_0} \leq \phi_0(n) \leq C n^{\gamma_0}.$$ 
We start with the bounds 
$$
\sup\limits_{v\in T_\GW} g_{\GW}(v)^{\gamma_0}\leq  C  \sup\limits_{v\in T_\GW} \phi_0(g_{\GW}(v)) \leq  C \sup\limits_{\ell \geq 0} \phi_0\big(\mathbf{Z}(\ell)\big),
$$
where 
$$\phi_0\big( \mathbf{Z}(\ell)\big)\coloneqq \sum_{j\geq 1} \phi_0(j) Z_j(\ell).$$
In particular, for any $y>0$, we have
$$  \Q_n^\GW\big(\sup_{T_{\GW}} g_{\GW} \geq y^{1/\gamma_0}\big)  \leq \Q^\GW_n\Big(\sup\limits_{\ell\geq 0} \phi_0( \mathbf{Z}(\ell))\geq y/C\Big).$$
We readily check from the superharmonicity of  $\phi_0$ that the nonnegative process $\phi_0( \mathbf{Z}(\cdot))$
is a supermartingale, hence
 the maximal inequality yields 
$$\Q^\GW_n\Big(\sup\limits_{\ell\geq 0} \phi_0( \mathbf{Z}(\ell))\geq y/C\Big) \leq Cy^{-1}\phi_0(n) \leq C^2 y^{-1} n^{\gamma_0} .$$
The substitution $y=(xn)^{\gamma_0}$ now entails \eqref{eq:L:Aspandii:1}.

Next, we claim that
\begin{equation}\label{eq:lem:Heigh:GW:bis}
\sup_{x > 0}\sup_{n\geq 1}   \left(x^{ \gamma_0/\alpha} \cdot \Q_n^\GW\big(\sup\big\{\ell\geq 0:~\mathbf{Z}(\ell)\neq 0\big\}  \geq x n^{\alpha}\big) \right)< \infty.
\end{equation}
Before proving \eqref{eq:lem:Heigh:GW:bis}, let us explain why the remaining assertions of the statement follow. First, \eqref{eq:lem:Heigh:GW:bis} plainly ensures that the Galton--Watson process becomes eventually extinct,
$\Q_n^\GW$-a.s. for all $n\geq 1$,  and \textit{a fortiori}  $(T_{\GW},g_{\GW})$ is bounded a.s. 
Moreover, since
$$
 \mathrm{Height}(T_{\GW})=1+\sup\big\{\ell\geq 0:~\mathbf{Z}(\ell)\neq 0\big\},
$$
the identity \eqref{eq:L:Aspandii:2} is checked. 

 Let us finally prove  \eqref{eq:lem:Heigh:GW:bis}.  To this end,  we get from the branching property and the notation \eqref{def:pinbold} and \eqref{Eq:Laplacediscret} that
$$\mathbb{E}_{n}^{\GW}\Big(\phi_0\big( \mathbf{Z}(\ell+1)\big) \Big | \mathbf{Z}(1), ... , \mathbf{Z}(\ell)\Big)=\sum\limits_{k\geq 1} Z_{k}(\ell) \boldsymbol{m}_{\GW}\phi_0(k)=\sum\limits_{k\geq 1} Z_{k}(\ell) \phi_0(k) \Big(1+\frac{\boldsymbol{G}_{\GW} \phi_0(k)}{\phi_0(k)}\Big).$$ 
Using again that $\phi_0(n) \asymp n^{\gamma_0}$ and $ \limsup_{n\to \infty} \boldsymbol{G}_{\GW} \phi_0(n)/n^{\gamma_0-\alpha}<0$, we infer the existence of  a constant $c_{\GW}>0$ such that $\boldsymbol{G}_{\GW} \phi_0(k)/\phi_0(k)\leq -c_{\GW} \phi_0(k)^{-\alpha/\gamma_0}$ for every $k\geq 1$. This implies that
 \begin{eqnarray*} \mathbb{E}_{n}^{\GW}\Big(\phi_0\big( \mathbf{Z}(\ell+1)\big)\, \Big |\, \mathbf{Z}(1), ... , \mathbf{Z}(\ell)\Big)&\leq& \sum_{k \geq 1} Z_{k}(\ell) \phi_0(k)\big(1-  c_{ \GW}\phi_0(k)^{-\alpha/\gamma_0}\big)\\ & \leq &\sum_{k \geq 1} Z_{k}(\ell) \phi_0(k)\big(1-   c_{\GW} \phi_0( \mathbf{Z}(\ell))^{-\alpha/\gamma_0}\big)\\ & =& \phi_0( \mathbf{Z}(\ell))\big(1- c_{\GW} \phi_0( \mathbf{Z}(\ell))^{-\alpha/\gamma_0}\big),  \end{eqnarray*}
 as long as $\mathbf{Z}(\ell)\neq 0$. The desired result  \eqref{eq:lem:Heigh:GW:bis} now follows from general martingale techniques; see Lemma~\ref{lem:supermartingaleextinct} in the Appendix for a general statement which is of independent interest. 
 \end{proof}
 
Assumption \ref{A:all:existance:PHI} also enables us to control the function $l^{\varpi}_{\GW}$ defined in \eqref{def:lA}.
\begin{proposition}\label{prop:massepresqueomega-}  
Let Assumption \ref{A:all:existance:PHI}  hold. For every $\gamma\in (\gamma_1,\gamma_0]$ and every weight function $\varpi$ such that $\varpi(n)=\mathcal{O}(n^{\gamma-\alpha})$, the function $l_{\GW}^{\varpi}$ takes finite values and satisfies  $l_{\GW}^\varpi(n)=\mathcal{O}(n^{\gamma})$.
\end{proposition}
The proof of Proposition \ref{prop:massepresqueomega-} relies on the following elementary comparison result. Recall the notation \eqref{Eq:Laplacediscret} for the generator.
  \begin{lemma}[Upper bound on superharmonic functions]  \label{lem:minimal} Let $\varpi: \mathbb{N}\to \mathbb{R}_+$ and $B\subset \mathbb{N}$. Assume that there exists $\phi:\mathbb{N}\to \mathbb{R}_+$  so that  $\phi \geq \varpi$ and 
\begin{equation}
\left\{
\begin{aligned}
    \phi &\geq  l_{\GW}^{ \varpi} && \text{on } B, \\
    \boldsymbol{G}_{\GW} \phi &\leq -\varpi && \text{on } B^{c}.
\end{aligned}
\right.
\label{eq:deltaplus}
\end{equation}
  Then  $\phi \geq l_{\GW}^{ \varpi}$ on $\N$, and consequently, $l_{\GW}^{ \varpi}$ takes finite values. 
  \end{lemma}

\begin{proof} For every  $j\geq 0$, we define the function
$$ \ell_j \coloneqq \sum_{i=0}^j \boldsymbol{m}_{\GW}^i\varpi,$$ 
which corresponds  to the expected mass of $\uplambda^{\varpi}_{\GW}$ assigned to  the first $j$ generations of the tree $T_\GW$.  Remark that $ \ell_0 = \varpi \leq \phi$, and  that $ \ell_j =   \boldsymbol{m}_{\GW} \ell_{j-1} + \varpi$ for $j\geq 1$.

Next we note that $ \ell_j \leq  \phi$ on $B$, as a consequence of the upperbound $ \ell_j \leq  l_{\GW}^{ \varpi}$ and the first inequality in \eqref{eq:deltaplus}.
We then argue by induction that for all $j\geq 0$, the bound $ \ell_j \leq \phi $ holds also on $B^c$. First, this assertion is obviously true for $j=0$, since $\ell_0=  \varpi \leq \phi$ by assumption. 
Next, suppose that the assertion holds for some $j\geq 0$. We have for any $n \in B^c$ 
 $$ \ell_{j+1}(n) =  \boldsymbol{m}_{\GW} \ell_j(n) + \varpi(n) \leq  \boldsymbol{m}_{\GW} \phi(n) + \varpi(n) \leq  \phi(n),$$
 where the first  inequality stems from our assertion at $j$ and the positivity of the operator $ \boldsymbol{m}_{\GW}$, and the second from \eqref{eq:deltaplus}. 
We conclude that $ l^{ \varpi}_{\GW} = \lim_{j\to \infty} \ell_j  \leq \phi$. 
 \end{proof}
 
 We are now in position to establish Proposition \ref{prop:massepresqueomega-}.

 \begin{proof}[Proof of Proposition \ref{prop:massepresqueomega-}]  Recall Assumption \ref{A:all:existance:PHI}, take $\gamma\in (\gamma_1,\gamma_0]$, and consider a weight function $\varpi:\mathbb{N}\to \mathbb{R}_+$ such that $\varpi(n)=\mathcal{O}(n^{\gamma-\alpha})$. In order to control $l_\GW^\varpi$, we will introduce a suitable superharmonic function that stems from Lemma \ref{lem:minimal}. We write $\gamma=(1-q)\gamma_1+q\gamma_0$ for some   $q\in(0,1]$, and introduce  $\phi\coloneqq \phi^{1-q}_1\cdot \phi_0^q$. Notice that   $\phi$ takes positive values and $\phi(n)\asymp n^{\gamma}$. Furthermore, by the inequality of arithmetic and geometric means, we have
\begin{align*}
\frac{\phi(k)}{\phi(n)}= \left(\frac{\phi_1(k)}{\phi_1(n)}\right)^{1-q}\cdot \left(\frac{\phi_0(k)}{\phi_0(n)}\right)^q\leq (1-q) \frac{\phi_1(k)}{\phi_1(n)}+ q \frac{\phi_0(k)}{\phi_0(n)},
\end{align*} 
for every $n,k\in \mathbb{N}$.  As a consequence, we get
$$\frac{\boldsymbol{G}_{\GW} \phi(n)}{\phi(n)}\leq (1-q) \frac{\boldsymbol{G}_{\GW} \phi_1(n)}{\phi_1(n)}+q \frac{\boldsymbol{G}_{\GW} \phi_0(n)}{\phi_0(n)},\quad n\in \mathbb{N}.$$
We now deduce from Assumption \ref{A:all:existance:PHI} that $\phi$ is strictly superharmonic with
$$ \limsup \limits_{n\to \infty} {n^{\alpha-\gamma}} \boldsymbol{G}_{\GW} \phi(n)<0,$$
Therefore,   we can always find $K>0$, such that both function $K \cdot \phi$ and $-K\cdot \boldsymbol{G}_{\GW} \phi$ are bounded from below by $\varpi$. We can then apply Lemma~\ref{lem:minimal}, with $B=\varnothing$, to deduce that $l_{\GW}^\varpi\leq K  \cdot\phi$ everywhere. In particular, $l_{\GW}^\varpi$ takes finite values and satisfies  $l_{\GW}^\varpi(n)=\mathcal{O}(n^{\gamma})$.
\end{proof}
 
We now conclude this section with another useful bound, which should be viewed as a version of the optional sampling theorem for supermartingales in a branching setting.
Recall that, following Jagers \cite{jagers1989general}, we call \textbf{line}  a subset of individuals $L\subset \U$ with the property that for any $u\in L$, the line $L$ does not contain any strict forebear $v\prec u$ of $u$.

\begin{lemma} \label{L:optsampsuperhar}
Let $\phi: \N\to \R_+$ be a superharmonic function for the Galton-Watson process. Consider any random line $L\subset \U$ such that
for any $u\in \U$, the event $\{u\in L\}$ is measurable with respect to  the sub-family of decoration-reproduction processes $(f_v, \eta_v)_{v\prec u}$. Then for any $n\geq 1$, there is the inequality
$$\E_n^{\GW}\left(\sum_{u\in L} \phi(\chi(u))\right) \leq \phi(n).
$$
\end{lemma} 
\begin{proof} The superharmonicity of $\phi$ for the particle system is transmitted to the population model 
by following the legitimate lineage of the ancestral particle. Indeed, we easily deduce from Proposition \ref{P:reprodkernGW} that 
$$\E_n^{\GW}\left(\sum_{j\geq 1} \phi(\chi(j))\right) \leq \phi(n).
$$
This is our claim for the line $L=\N$ formed by the first generation. The general case of a predictable random line follows; see the argument in the proof of \cite[Theorem 5.1]{jagers1989general}.
\end{proof}

\section{Invariance principles via a cutoff procedure}\label{sec:invprinc}

The three scaling limit theorems  stated in Section \ref{sec:4.2} will be established in this section.
 Let us first provide a brief overview of the strategy, which builds upon the pruning transformation introduced in Section \ref{sub:section:gene}. 
Recall that  $\P^\GW_n$ denotes the  law of the family of decoration-reproduction processes $(f_u, \eta_u)_{u\in \U}$ when
 the Galton--Watson system starts from a single particle of type $n$.  We then write ${\P}^{(n)}_{\GW}$ for its image by the scaling transformation induced by \eqref{eq:linmap}, which maps $(f_u, \eta_u)$  to $( f^{(n)}_u, \eta^{(n)}_u)$ for each $u\in \U$. In particular, every type takes values in $n^{-1}\N$ under ${\P}^{(n)}_{\GW}$.  For simplicity, in this section we systematically omit the subscript $\text{GW}$ on our variables under ${\P}^{(n)}_{\GW}$ and  $\P^\GW_n$,    when no ambiguity arises.  In this direction, we systematically  write $\texttt{T}$ for  the decorated tree constructed from $(f_u, \eta_u)_{u\in \mathbb{U}}$.  
 
Under ${\P}^{(n)}_{\GW}$, ${\P}_n^{\GW}$,  and  $\P_1$,   we consider  for every $r>0$ the stopping line
\begin{equation}
F(r)\coloneqq\big\{u\in \U: f_u(0)< r\text{ and } f_v(0)\geq r \text{ for all }v\prec u\big\}.
\end{equation}
That is, an individual $u$ belongs to $F(r)$ and is then viewed as frozen  if and only if it is the first individual of its ancestral lineage with type smaller than $r$. We write  $\texttt{T}^{[r]}:=(T^{[r]},d_{T^{[r]}},\rho^{[r]}, g^{[r]})$ for the decorated tree that is constructed by gluing from the subfamily $(f_u, \eta_u)_{u\in A(r)}$, where $A(r)\subset \U$ stands for the set of strict ascendents of individuals in $F(r)$. In words,  $\texttt{T}^{[r]}$  is the genealogical subtree resulting from freezing any individual  with type less than $r$ in the population model. By construction, we might and will see  $T^{[r]}$ as a subset of $T$ equipped with the restriction distance such that $\rho^{[r]}=\rho$ and $g^{[r]}\leq g$ on $T^{[r]}$. We refer to $\texttt{T}^{[r]}$  as the subtree pruned below the threshold $r$, and will also later on equip $\texttt{T}^{[r]}$ with different measures, depending on the situation. 

Our first task is to establish the limit theorems for sequences of decorated trees pruned below a fixed level $r>0$. This partly relies on a (deterministic) convergence result involving the gluing operation on line segments which we now state. Consider a measured decorated line segment  $\mathbf{T}_0=([0,z],d,t_0,f_0,\upnu_0)$, where $d$ stands for the usual distance on an interval,  
and  $t_0\in [0,z]$ serves as root. For some fixed integer  $k\geq 1$, consider also $t_1, \ldots , t_k\in [0,z]$
 and (rooted) measured decorated trees $\mathbf{T}_1=(\texttt{T}_1, \upnu_1), \ldots, \mathbf{T}_k=(\texttt{T}_k, \upnu_k)$. Perform the gluing operator and write  $\mathbf{T}=(\texttt{T}, \upnu)$ 
  where, in the notation of Section \ref{sec:1.1},
  $$
  \texttt{T}=\mathrm{Gluing}\left(\texttt{T}_0,(t_i)_{1\leq i \leq k}, (\texttt{T}_i)_{1\leq i \leq k}\right) ,$$
  and $\upnu$ is obtained by adding the measures $\upnu_j$ for $j=0, \ldots, k$ on each component $T_j$ (so the glue-points $t_j$ may eventually receive masses from several components).
 The following continuity claim should be intuitively obvious; we will leave its elementary proof (which bears some similarities with that of Lemma \ref{lem:topo}) to the interested reader.
 
 \begin{lemma}\label{L:convgluing} Consider further measured decorated line segments  $\mathbf{T}'_0=([0,z'],d,t'_0,f'_0,\upnu'_0)$,  locations  $t'_1, \ldots , t'_k\in [0,z']$,  and rooted measured decorated trees $\mathbf{T}'_1, \ldots, \mathbf{T}'_k$. Write  $\mathbf{T}'=(\normalfont{\texttt{T}}', \upnu')$ 
  where as, above,
 $$
 \normalfont{\texttt{T}}'=\mathrm{Gluing}\left( \normalfont{\texttt{T}}'_0,(t'_i)_{1\leq i \leq k}, ( \normalfont{\texttt{T}}'_i)_{1\leq i \leq k}\right),
 $$
 and $\upnu'$ results from adding the measures $\upnu'_j$ on each component. 
 
 If  $t'_j\to t_j$ and $\mathbf{T}'_j \to\mathbf{T}_j$ in $\TT_m$  for each $j=0, \ldots, k$, then $\mathbf{T}'$ converges to $\mathbf{T}$ in $\TT_m$.
 \end{lemma}

Once the limit theorem for the sequence of rescaled decorated trees pruned below a fixed threshold $r>0$ has been established, we only need to check that as the threshold goes to $0$, the pruned trees converge uniformly toward the original ones. This program is relatively easy to implement for Theorems \ref{T:mainuncondw} and \ref{T:mainunconds-length}, and we shall start with those. The case of Theorem~\ref{T:mainunconds-mass} will be more involved, and requires some further technical efforts.

\subsection{Proof of Theorem \ref{T:mainuncondw} }\label{Sub:sec:prt1}

We write  $\Upomega$ for the space of families of decoration-reproduction pairs $(f_u, \eta_u)_{u\in \U}$ indexed by the Ulam tree, and endow $\Upomega$ with the product topology.
The laws $\P_1$ and ${\P}^{(n)}_{\GW}$ for $n\geq 1$ are thus probability measures  on $\Upomega$, and in this setting, weak convergence of probability measures on $\Upomega$
is the same as  convergence in the sense of finite dimensional distributions.
We start by observing that Assumption \ref{assum:main:gen:deco:repro} about the  scaling behavior of the decoration-reproduction process  of the ancestor propagates to the entire family of decoration-reproduction processes $(f_u, \eta_u)_{u\in \U}$. 
  
 \begin{lemma} \label{Lem:confindim} Let Assumption \ref{assum:main:gen:deco:repro} hold. As $n\to \infty$, ${\P}^{(n)}_{\GW}$  converges towards $\P_1$
weakly on $\Upomega$.
 \end{lemma} 
 
 \begin{proof} We have to check that for every finite subset $V\subset \U$, the law of the rescaled subfamily $( f^{(n)}_v, \eta^{(n)}_v)_{v\in V}$ under $\P^\GW_n$  converges as $n\to \infty$ towards that of $(f_v,\eta_v)_{v\in V}$ under $\P_1$. Without loss of generality, we may suppose that $V$ is a subtree of $\U$ that contains the root $\varnothing$, and 
 we  argue by induction on the height of $V$,  that is $\max_V|v|$.

   For $V=\{\varnothing\}$, the claim reduces to  Assumption \ref{assum:main:gen:deco:repro}.
 Suppose next that $\max_V|v|=1$.  Skorokhod representation theorem allows us to assume that the convergence of  the reproduction processes $\eta^{(n)}_{\varnothing}$ towards $\eta_{\varnothing}$ holds almost-surely for the vague topology (which is metrizable).  We see that for any $j\geq 1$,  if we write $\chi^{(n)}(j)=f^{(n)}_j(0)$ for the  type  of the $j$-th child of the ancestor under $\P_{\GW}^{(n)}$, then $\lim_{n\to \infty} \chi^{(n)}(j)=\chi(j)$; where we recall that the atoms $(t_j,f_j(0))_{j\in \mathbb{N}}$ of $\eta_\varnothing$ are indexed in co-lexicographical order. It now readily follows from  Assumption \ref{assum:main:gen:deco:repro}, self-similarity
 and the branching property, that  the law of  $( f^{(n)}_j, \eta^{(n)}_j)_{j\geq 1}$ under $\P^\GW_n$ tends indeed to that of $(f_j,\eta_j)_{j\geq 1}$ under $\P_1$ as $n\to \infty$.
 \textit{A fortiori} the assertion holds for any subtree $V$ with height $\max_V|v|=1$. 
 An iteration of this argument using the branching property completes the proof. 
  \end{proof}

 The sought convergence for the sequence of pruned decorated trees is now easy to establish. 
  
  \begin{corollary} \label{Cor:convpruned} Let Assumption \ref{assum:main:gen:deco:repro} hold. For every threshold $0<r <1$ such that
  $$\P_1(\chi(u)=r)=0 \qquad \text{for any  }u\in \U,$$
  the law under  ${\P}^{(n)}_{\GW}$ of the decorated  subtree 
  $\normalfont{\texttt{T}}^{[r]}$  pruned below the threshold $r$  converges  as $n\to \infty$ towards that under $\P_1$.
   \end{corollary} 
 Later on, we will use the fact that there are at most countably many $r>0$ such that $\P_1(\chi(u)=r)>0$ for some $u\in \U$.
 
 \begin{proof} The claim might look quite intuitive from Lemma \ref{Lem:confindim}; nonetheless some care is needed to deal with technicalities,  and we present a  
 detailed argument.  

 Consider any  finite subtree $V\subset \U$ which contains the root and write $\partial V$ for the subset of leaves of $V$. We introduce the event 
 $$\mathcal A(V,r)\coloneqq \{(f_u, \eta_u)_{u\in \U}\in \Upomega: f_u(0)>r\text{ for all }u\in V\, \&\  \eta_v(\R_+\times [r,\infty))=0 \text{ for all }v\in \partial V\}.$$
 Observe  that $\mathcal A(V,r)$  is 
  open in $\Upomega$, and actually  measurable with respect to the sub-family $(f_u,\eta_u)_{u\in V}$. We also point out that $A(r)= V$ on the event $\mathcal A(V,r)$, where we recall that  $A(r)$ denotes the part of $\U$ strictly below the stopping line $F(r)$, that is the subset of  individuals $u\subset \U$ such that $u$ and all its forebears have types at least $r$.
 
 Next, introduce the function ${\texttt{T}}^{V}: \Upomega \to \TT$ that maps $(f_u, \eta_u)_{u\in \U}$ to
 the decorated tree obtained by gluing from the sub-family  $(f_u, \eta_u)_{u\in V}$. We readily check from Lemma \ref{L:convgluing} using an induction on the height of $V$, that ${\texttt{T}}^{V}$
  is continuous on $\mathcal A(V,r)$ and coincides there with the decorated tree pruned below level $r$, i.e.~
  ${\texttt{T}}^{V}={\texttt{T}}^{[r]}$ on $\mathcal A(V,r)$. We can now deduce from Lemma \ref{Lem:confindim} and the Portemanteau theorem
  that for every continuous functional   $\varPhi: \TT\to [0,1]$
  $$\liminf_{n\to \infty} \E^{(n)}_{\GW}\left(\varPhi({\texttt{T}}^{[r]}) \indset{\mathcal A(V,r)}\right) \geq  \E_1\left(  \varPhi({\texttt{T}}^{[r]}) \indset{\mathcal A(V,r)}\right).$$
  The events $\mathcal A(V,r)$ are pairwise disjoint for different  finite subtrees $V\subset \U$ that contains the root, and if we write
  $ \mathcal A(r)\coloneqq \bigsqcup \mathcal A(V,r)$ for their union, then
    Fatou's lemma entails
   $$\liminf_{n\to \infty} \E^{(n)}_{\GW}\left(   \varPhi({\texttt{T}}^{[r]}) \right) \geq  \E_1\left(   \varPhi({\texttt{T}}^{[r]}) \indset{\mathcal A(r)}\right). $$

The assumption that $\P_1( \chi(u)=r)=0$ for all $u\in \U$ is equivalent to $\eta_v(\R_+\times\{r\})=0$, $\P_1$-a.s. for all $v\in \U$, 
from which it is easily deduced that $\P_1(\mathcal A(r))=1$. Thus we have
$$\liminf_{n\to \infty} \E^{(n)}_{\GW}\left(   \varPhi({\texttt{T}}^{[r]}) \right) \geq  \E_1\left(   \varPhi({\texttt{T}}^{[r]}) \right).$$
Finally replacing $\Phi$ by $1-\Phi$, we conclude  that
$$\lim_{n\to \infty} \E^{(n)}_{\GW}\left(   \varPhi({\texttt{T}}^{[r]}) \right) =  \E_1\left(   \varPhi({\texttt{T}}^{[r]}) \right),$$
and the proof is complete. 
 \end{proof}
 
 We now have all the ingredients needed to establish  Theorem \ref{T:mainuncondw}. We stress that Assumption \ref{A:all:existance:PHI}, that was not needed previously, now plays a crucial role.

\begin{proof}[Proof of Theorem \ref{T:mainuncondw}]
The heart of the matter is to get a uniform version of Corollary \ref{Cor:convpruned}, since Corollary \ref{lem:cutoff} already ensures the convergence of pruned self-similar Markov trees as the threshold goes to $0$. 
 More precisely, we will show that for any $\varepsilon >0$, there exists some finite constant $c(\varepsilon)$ such that, for all $0< r <1$ and all $n\geq 1$,
 \begin{equation} \label{Eq:borneunif}
 \P^{(n)}_{\GW}\left(\d_{\mathrm {Hyp}}(g,g^{[r]}) >\varepsilon\right)\leq c(\varepsilon)  r^{\gamma_0-\gamma_1}.
 \end{equation}
In this direction, we bound the distance between the pruned trees and the initial one by analyzing the subtrees dangling from the former.
For any decorated tree $(T^\prime,g^\prime)$,   the hypograph distance between $\texttt{T}$ and the degenerate tree is given by definition as $ \mathrm{Height}(T^\prime) \vee \max_{T^\prime} g^\prime$,
  and this
 entails the bound
 $$\d_{\mathrm {Hyp}}(g,g^{[r]})  \leq \sup_{u\in F(r)} \left\{ \mathrm{Height}(T_u) \vee \max_{T_u}g_u\right\},
 $$
where for every $u\in \U$, $(T_u, g_u)$ denotes the decorated subtree that stems from $u$. 
Specializing this to rescaled Galton--Watson trees, we get that for any $\varepsilon >0$,
\begin{align*}
& \P^{(n)}_{\GW}\left(\d_{\mathrm {Hyp}}(g,g^{[r]}) >\varepsilon\right)\\
& \leq 
\P^{(n)}_{\GW}\big(\exists u\in F(r): \mathrm{Height}(T_u)>\varepsilon  \text{ or } \max_{T_u}g_u>\varepsilon  \big)\\
& \leq \P^{\GW}_n\big(\exists u\in F(rn): \mathrm{Height}(T_u)>\varepsilon n^{\alpha}\big)+ \P^{\GW}_n\big(\exists u\in F(rn): \max_{T_u} g_u>\varepsilon n\big),
\end{align*}
where we used the scaling transformation for the second inequality.

We detail the argument to bound the first term of the sum above only, that for the second being similar. We use $c_1, c_2, \ldots$ to denote some positive finite constants
that will appear in the calculation. 
By the branching property  under  $\P^{\GW}_n$, conditionally given that
the type of the individual  $u$ is $\chi(u)=\ell$,  the sub-genealogical decorated tree $\texttt{T}_u$ is independent of the types of the forebears of $u$ and has the law $\Q^{\GW}_{\ell}$. Recalling from Proposition \ref{L:Aspandii} that 
for  any $\ell, n\geq 1$,
$$\Q_{\ell}^{\GW}\left(\mathrm{Height}(T_{\GW})\geq \varepsilon n^{\alpha}\right) \leq c_1 \left( \frac{\ell}{\varepsilon^{1/\alpha} n} \right)^{\gamma_0}$$
and applying the union bound, we arrive at
$$ \P^{\GW}_n\big(\exists u\in F(rn): \mathrm{Height}(T_u)>\varepsilon n^{\alpha}\big) 
\leq c_1 \varepsilon^{-\gamma_0/\alpha} n^{-\gamma_0} \E^{\GW}_n\Big( \sum_{u\in F(rn)} \chi(u)^{\gamma_0} \Big).
$$
We then use Assumption \ref{A:all:existance:PHI} and Lemma \ref{L:optsampsuperhar} to deduce from the superharmonicity of $\phi_1$  that
 $$\E^{\GW}_n\Big( \sum_{u\in F(rn)} \phi_1(\chi(u)) \Big) \leq \phi_1(n)\leq c_2 n^{\gamma_1}.$$
Since for any $ u\in F(rn)$,
$$\chi(u)^{\gamma_0} \leq c_3 \chi(u)^{\gamma_0-\gamma_1}\phi_1(\chi(u)) \leq c_3(rn)^{\gamma_0-\gamma_1}\phi_1(\chi(u)), $$
we conclude that
$$ \P^{\GW}_n\big(\exists u\in F(rn): \mathrm{Height}(T_u)>\varepsilon n^{\alpha}\big)
\leq c_4 \varepsilon^{-\gamma_0/\alpha} r^{\gamma_0-\gamma_1}.$$
A very similar calculation gives
$$ \P^{\GW}_n\big(\exists u\in F(rn): \max_{T_u} g_u>\varepsilon n\big)
\leq c_5 \varepsilon^{-\gamma_0} r^{\gamma_0-\gamma_1},$$
so that putting pieces together, we have established \eqref{Eq:borneunif}.

The rest of the proof is now straightforward. Consider any $1$-Lipschitz continuous functional $\Phi: \TT\to [0,1]$, so that thanks to
\eqref{Eq:borneunif}, for all $0< r <1$ and all $n\geq 1$,
$$\E^{(n)}_{\GW} \left( | \Phi(\texttt{T}) - \Phi(\texttt{T}^{[r]})|\right) \leq \varepsilon + c(\varepsilon)  r^{\gamma_0-\gamma_1}, 
$$
where to apply \eqref{Eq:borneunif} we used that by definition $\d_{\mathbb{T}}(\texttt{T}, \texttt{T}^{[r]})\leq \d_{\mathrm {Hyp}}(g,g^{[r]})$.  Corollary \ref{lem:cutoff} now enables us to choose $r>0$ small enough such that  $c(\varepsilon)  r^{\gamma_0-\gamma_1}\leq \varepsilon$, 
$$\E_1\left( | \Phi(\texttt{T}) - \Phi(\texttt{T}^{[r]})|\right) \leq \varepsilon, $$
and furthermore $\P_1(\chi(u)=r)=0$ for all $u\in \U$.  
Invoking Corollary \ref{Cor:convpruned}, we conclude that
$$\limsup_{n\to \infty} \left| \E_1\left(  \Phi(\texttt{T})\right) -\E^{(n)}_{\GW} \left(  \Phi(\texttt{T}) \right)
\right| \leq 3 \varepsilon.$$
Since $\varepsilon$ can be chosen arbitrarily small, this proves Theorem \ref{T:mainuncondw}.
\end{proof}
 
\subsection{Proof of Theorem \ref{T:mainunconds-length}}\label{Sub:sec:scali:dec:2}

Theorem \ref{T:mainunconds-length} is a stronger version of Theorem \ref{T:mainuncondw} which takes further into account weighted length measures; its proof via pruning follows the same guiding line. Throughout this section, we will implicitly assume that Assumptions \ref{assum:main:gen:deco:repro} and \ref{A:all:existance:PHI} hold, and that
the weight function $\varpi$ is regularly varying at $\infty$ with exponent $\gamma-\alpha$ for some $\gamma>\gamma_0$. Beware that these assumptions will not be repeated  in the statements!
Of course, they will also cease to be enforced in the next section.

Let us start introducing some useful notation. In this direction, we write
$$ \varpi^{(n)}:= \varpi(\lfloor n \cdot\rfloor)/\varpi(n),\quad n\geq 1,$$
for the rescaled weight function. For a generic decoration-reproduction pair $(f,\eta)$ with domain $[0,z]$, we denote by  $\upnu$ the measure on $[0,z]$ which has the density $\varpi\circ f$ with respect to the Lebesgue measure. For every $n\geq 1$, we write  $(f^{(n)}, \eta^{(n)}, \upnu^{(n)})$ for its rescaled version, where the first two components have been defined in \eqref{Eq:scalefeta}, 
and the third is the image of $ n^{-\alpha} \varpi(n)^{-1}\upnu$ by the map $t\mapsto n^{-\alpha}t$ on $[0,z]$. 
We also write $(f_u^{(n)}, \eta_u^{(n)}, \upnu_u^{(n)})$ when the initial decoration-reproduction pair $(f_u,\eta_u)$ is associated to some individual $u$.
We now point at the following variation of Lemma  \ref{Lem:confindim}.

\begin{lemma} \label{Lem:confindimmeas} 
As $n\to \infty$, the law of the rescaled family $(f_u^{(n)}, \eta_u^{(n)}, \upnu_u^{(n)})_{u\in \U}$ under ${\P}_n^{\GW}$  converges in the sense of finite-dimensional distributions towards 
that of $(f_u, \eta_u, \uplambda^{\gamma}_u)_{u\in \U}$ under $\P_1$, where 
$$\uplambda^{\gamma}_u(\dd t) = f_u(t)^{\gamma-\alpha}\dd t, \qquad t\in[0,z_u].$$
 \end{lemma} 
\begin{proof} Combining the Skorokhod representation and Lemma  \ref{Lem:confindim}, and by a slight abuse of notation, we may assume that we are given a
family of decoration-reproduction processes $(f_u, \eta_u)_{u\in \U}$ with law $\P_1$ and 
a rescaled family of processes $(f_u^{(n)}, \eta_u^{(n)})_{u\in \U}$ with law ${\P}_{\GW}^{(n)}$, such that  $\lim_{n\to\infty} (f_u^{(n)}, \eta_u^{(n)})=(f_u, \eta_u)$ for every $u\in \U$, almost-surely. 

Now remark that, with the notation   $\varpi^{(n)}$ introduced above, we must have $\lim_{n\to \infty}\varpi^{(n)}(x)= x^{\gamma-\alpha}$. Therefore, for Lebesgue-almost every $t\geq 0$,  it follows that
\begin{equation}\label{Eq:avantui}
\lim_{n\to\infty}\indset{[0,z^{(n)}_u]}(t)\varpi^{(n)}\circ f_u^{(n)}(t) =\indset{[0,z_u]}(t)f_u(t)^{\gamma-\alpha}, \qquad \text{a.s.}
\end{equation}
The heart of the proof consists of verifying that the converge above holds also in $L^1(\dd t)$, in probability, that is, equivalently, that each sequence 
$\upnu_u^{(n)}$ converges in total variation as $n\to \infty$ towards $ \uplambda^{\gamma}_u$, in probability. By the branching property, we may focus on the case $u=\varnothing$ without loss of generality, and for the sake of notational simplicity, we now drop the indices $u$ from the notation.

The small values of $f^{(n)}$  require a special attention when $\gamma-\alpha\leq 0$, and for this we introduce some continuous function $\varphi: \R_+\to [0,1]$ such that $\varphi(x)=1$ for $0\leq x \leq 1$ and $\varphi(x)=0$ for $x\geq 2$.
On the one hand, it is plain from the convergence of $f^{(n)}$ to $f$ in the sense of Skorokhod that 
$$\lim_{n\to\infty} \int_0^{\infty} \dd t \left|\indset{[0,z^{(n)}]}(t)\big(1-\varphi(f^{(n)}(t))\big) \varpi^{(n)}(f^{(n)}(t))-\indset{[0,z]}(t)\big(1-\varphi(f(t))\big)f(t)^{\gamma-\alpha}\right|= 0, \text{ a.s.}$$
On the other hand, since the weight function $\varpi$ is regularly varying at $\infty$ with exponent 
$\gamma-\alpha$, we know from the Potter bound \cite[Theorem 1.5.6]{BGT89}
 that for any $\gamma'\in (\gamma_0,\gamma)$
 \begin{equation}\label{Eq:Potterbound}c(\gamma')\coloneqq \sup_{x\geq 1/n} \varpi^{(n)} (x)\varphi(x)x^{\alpha-\gamma'} <\infty.
 \end{equation}
We now choose $p>1$ sufficiently close to $1$ such that $p\gamma'-(p-1)\alpha>\gamma_0$ and write
$$ \int_0^{z^{(n)}} \dd t \left|\varpi^{(n)}(f_u^{(n)}(t))\varphi(f^{(n)}(t))\right|^p\leq c(\gamma')^p \int_0^{z^{(n)}} |f^{(n)}(t)|^{p(\gamma'-\alpha)} \dd t .$$
Keeping track of the scaling transformations, the expectation of the integral in the right-hand side is given in terms of the Galton--Watson process as
$$\E_n^{\GW} \left(n^{-\alpha-p(\gamma'-\alpha)} \int_0^z f(t)^{p(\gamma'-\alpha)}\dd t \right) \leq n^{(p-1)\alpha-p\gamma'}  \mathbb{E}^{\GW}_{n} \Big(  \uplambda^{\upomega}_{\GW}\big( T_{\GW}\big)\Big),
$$
where in the right-hand side, $T_{\GW}$ is measured with the weight function  $\upomega(j)=j^{p(\gamma'-\alpha)}$.
Proposition \ref{prop:massepresqueomega-}  enables us to bound the right-hand side uniformly for all $n\geq 0$, so that
$$\sup_{n\geq 1} \E\left( \int_0^{z^{(n)}} \dd t \left|\varpi^{(n)}(f_u^{(n)}(t))\varphi(f^{(n)}(t))\right|^p \right)<\infty.$$

Returning to \eqref{Eq:avantui},  we conclude from uniform integrability that
$$\lim_{n\to\infty} \E\left(\int_0^{\infty} \dd t \left|\indset{[0,z^{(n)}]}(t)\varpi^{(n)}(f_u^{(n)}(t))\varphi(f^{(n)}(t)) -\indset{[0,z_u]}(t)f(t)^{\gamma-\alpha}\varphi(f(t))\right|\right)= 0.$$
Putting the pieces together, for every $u\in \U$, the measure $\upnu_u^{(n)}$ converges as $n\to \infty$ in total variation to $\uplambda_u^{\gamma}$, in probability; which is more than needed for the proof of the statement.
\end{proof}

Lemma \ref{Lem:confindimmeas} enables us to incorporate weighted length measures in the convergence of pruned trees stated in Corollary \ref{Cor:convpruned}.
To make the notation less cluttered, we just write $\uplambda$ instead of the more specific $\uplambda_T$ for the length measure on a real $T$ when the latter is clear from the context.

  \begin{corollary} \label{Cor:convprunedmes} 
  For every threshold $0<r <1$ such that
  $$\P_1(\chi(u)=r)=0 \qquad \text{for any  }u\in \U,$$
  the law under  ${\P}^{(n)}_{\GW}$ of the decorated  subtree $\normalfont{\texttt{T}}^{[r]}$ 
   pruned below the threshold $r$,   further  endowed with the weighted length measure 
   $$\uplambda^{\varpi^{(n)}}_{T^{[r]}}\coloneqq \indset{T^{[r]}}\cdot \uplambda^{\varpi^{(n)}} =  \indset{T^{[r]}}\cdot \left( \varpi^{(n)}\circ g\cdot \uplambda\right) ,$$  converges  as $n\to \infty$, towards the  law under  ${\P}_1$ of the decorated  subtree 
  $\normalfont{\texttt{T}}^{[r]}$ equipped with the weighted length measure 
  $$\uplambda^{\gamma}_{T^{[r]}}\coloneqq\indset{T^{[r]}}\cdot  \uplambda^\gamma= \indset{T^{[r]}}\cdot  \left( g^{\lambda-\alpha} \cdot \uplambda\right) .$$
   \end{corollary} 

\begin{proof} The argument is the same as for Corollary \ref{Cor:convpruned}, except that we incorporate weighted length measures in the functionals using Lemma \ref{Lem:confindimmeas}.
\end{proof}

We can now check the version of invariance principle including weighted length measures.
\begin{proof}[Proof of Theorem \ref{T:mainunconds-length}]

Recall  that $\varpi$ is is regularly varying at $\infty$ with exponent $\gamma-\alpha$, for some $\gamma>\gamma_0$, and  from Lemma \ref{lem:cutoff} the convergence of pruned self-similar Markov trees endowed with a weighted length measure as the threshold tends to $0$. 
The only addition needed to the proof of Theorem \ref{T:mainuncondw} in the preceding section, is to verify that 
the Prokhorov distance between the weighted length measure $\uplambda^{\varpi^{(n)}}$ on $T$ and its restriction $\uplambda^{\varpi^{(n)}}_{T^{[r]}}$ to the pruned tree $T^{[r]}$ 
 converge in probability to $0$ as $r\to 0$  under $\P_{\GW}^{(n)}$, uniformly in $n$.

To start with, introduce the weight-function $\upomega_0(x)= x^{\gamma_0-\alpha}$. We pick $\gamma'\in(\gamma_0, \gamma)$ and  for any $0<a<1$, rewrite the Potter's bound \eqref{Eq:Potterbound}  in the form 
$$\varpi^{(n)}(x) \leq c(\gamma') a^{\gamma'-\gamma_0} \upomega_0(x), \qquad \text{for all } n\geq 1 \text{ and } 1/n\leq x \leq a.$$
We then have 
$$\uplambda^{\varpi^{(n)}}(T\backslash T^{[r]})\leq c(\gamma') a^{\gamma'-\gamma_0} \uplambda^{\upomega_0}(T\backslash T^{[r]}) \text{ on the event }\sup_{T\backslash T^{[r]}} g\leq a.$$
Next, keeping track of the rescaling, we write in terms of the dandling subtrees
$$\E_{\GW}^{(n)}\left( \uplambda^{\upomega_0}(T\backslash T^{[r]})\right) = \E_{\GW}^{(n)}\Big( \sum_{u\in F(r)} \uplambda^{\upomega_0}(T_u)\Big) = n^{-\gamma_0} \E^{\GW}_n\Big( \sum_{u\in F(rn)} \uplambda^{\upomega_0}(T_u)\Big).$$
Recall from the branching property  under  $\P^{\GW}_n$ that conditionally given that
the type of the individual  $u$ is $\chi(u)=\ell$,  the decorated genealogical subtree $\texttt{T}_u$ is independent of the types of the forebears of $u$ and has the law $\Q^{\GW}_{\ell}$. As a consequence, recalling the notation \eqref{def:lA},  there is the identity
$$\E^{\GW}_n\left( \indset{u\in F(rn)} \uplambda^{\upomega_0}(T_u)\right)= \E^{\GW}_n\left( \indset{u\in F(rn)} l^{\upomega_0}_{\GW}(\chi(u)) \right) .$$
We then deduce from Proposition \ref{prop:massepresqueomega-} and the assumption $\phi_0(n)\asymp n^{\gamma_0}$, that  there is some finite constant $c_1(\gamma')$ such that 
$$
\E_{\GW}^{(n)}\Big( \uplambda^{\varpi^{(n)}}(T\backslash T^{[r]}), \sup_{T\backslash T^{[r]}} g\leq a\Big) 
\leq c_1(\gamma') a^{\gamma'-\gamma} n^{-\gamma_0} \E^{\GW}_n\Big( \sum_{u\in F(rn)} \phi_0(\chi(u))\Big).
$$
Recalling (from the superharmonicity of $\phi_0$ and Lemma \ref{L:optsampsuperhar}) that
 $$\E^{\GW}_n\Big( \sum_{u\in F(rn)} \phi_0(\chi(u)) \Big) \leq \phi_0(n)\asymp n^{\gamma_0},$$
we conclude that for any $a\in (0,1)$ and $r>0$, there is the bound
\begin{equation} \label{Eq:bornunifmeas}
\sup_{n\geq 1}\E_{\GW}^{(n)}\Big( \uplambda^{\varpi^{(n)}}(T\backslash T^{[r]}), \sup_{T\backslash T^{[r]}} g\leq a\Big) \leq c_2(\gamma') a^{\gamma'-\gamma} .
\end{equation}

The rest of the proof is now easy. We know from \eqref{Eq:borneunif} that for any $a>0$,
$$\lim_{r\to 0} \sup_{n\geq 1}\P^{(n)}_{\GW}\Big (\sup_{T\backslash T^{[r]}} g >a\Big) =0.$$
We deduce from \eqref{Eq:bornunifmeas} by Markov's inequality that for any $b>0$,
$$\lim_{r\to 0} \sup_{n\geq 1}\P_{\GW}^{(n)}\left( \uplambda^{\varpi^{(n)}}(T\backslash T^{[r]})>b\right) =0.$$
As a consequence, we have in terms of the Prokhorov distance that
$$\lim_{r\to 0} \sup_{n\geq 1}\P_{\GW}^{(n)}\left( \mathrm{d}_{\mathrm{Prok}}(\uplambda^{\varpi^{(n)}}, \indset{ T^{[r]}}\cdot \uplambda^{\varpi^{(n)}})>b\right) =0.$$
Then,  by combination  with  \eqref{Eq:borneunif},  we get that
$$\lim_{r\to 0} \sup_{n\geq 1}\P_{\GW}^{(n)}\left( \mathrm{d}_{\TT_m}\left((\texttt{T},\uplambda^{\varpi^{(n)}}), (\texttt{T}^{[r]}, \indset{ T^{[r]}}\cdot \uplambda^{\varpi^{(n)}})\right)>b\right) =0.$$
Just as in the proof of Theorem \ref{T:mainuncondw}, we can now deduce from  Lemma \ref{lem:cutoff} and Corollary \ref{Cor:convprunedmes} that for 
any $1$-Lipschitz continuous functional $\Phi: \TT_m\to [0,1]$, 
$$\lim_{n\to \infty} \left| \E_1\left( \Phi(\mathbf{T})\right) -\E^{(n)}_{\GW} \left( \Phi(\mathbf{T}) \right)
\right| =0,$$
which proves Theorem \ref{T:mainunconds-length}.
\end{proof}

\subsection{Scaling limit of the total mass and proof of Theorem \ref{T:mainunconds-mass}} \label{sec:4.6} 

Our  goal here is to prove Theorem \ref{T:mainunconds-mass}; our approach using pruned trees emulates that implemented the preceding sections. 
Throughout this section, we will implicitly assume that  Assumptions \ref{A:omega-},  \ref{assum:main:gen:deco:repro} and \ref{A:all:existance:PHI}  hold, and that the weight function $\varpi$ verifies Assumption \ref{A:hGWrv}. 
We also  write $\Q_1$ for the law of a self-similar Markov tree  with characteristic quadruplet $(\sigma^{2}, \mathrm a,  \boldsymbol{\Lambda}; \alpha)$ endowed with its harmonic measure, $\mathbf{T}=(\texttt{T}, \upmu)$. 
Beware again that these assumptions and notation will not be repeated in the statements!

We start  by pointing out
that when the initial type $n$ of the Galton--Watson process is large, the mass of the legitimate lineage of the ancestral particle becomes asymptotically negligible compared to  the mass of the whole tree. As a consequence of the branching property, the measure  $\uplambda^{\varpi}_{\GW}$ concentrates   as $n\to \infty$ around the leaves of $T_{\GW}$, which contrasts with the preceding section. 
\begin{lemma} \label{L:asympneglbr} We have
$$\E_n^{\GW}\left( \int_0^{z_{\varnothing}} \varpi(f_{\varnothing}(t)) \dd t\right) =  o\left(  l^{\varpi}_\GW(n)\right), \qquad \text{as }n\to \infty.$$
\end{lemma}

\begin{proof} By decomposing the genealogical tree of the Galton--Watson process along the ancestral branch (i.e the legitimate lineage of the ancestral particle), and applying the branching property, we get with the notation \eqref{def:lA}
that 
$$ l^{\varpi}_\GW(n) = \E_n^{\GW}\left( \int_0^{z_{\varnothing}} \varpi(f_{\varnothing}(t)) \dd t\right)  +  \E_n^{\GW} \left( \sum_{j\geq 1} l^{\varpi}_\GW (f_j(0))\right).$$
Since $ l^{\varpi}_\GW$ is regularly varying with exponent $\omega_-$ and for every $j\geq 1$,  the law of $f_j(0)/n$ under $\P_n^{\GW}$ converges to that of 
$f_j(0)=\chi(j)$ under $\P_1$, we conclude from Fatou's Lemma that
$$\limsup_{n\to \infty} l^{\varpi}_\GW(n)^{-1} \E_n^{\GW}\left( \int_0^{z_{\varnothing}} \varpi(f_{\varnothing}(t)) \dd t\right) 
\leq 1-\E_1\left( \sum_{j\geq 1}\chi(j)^{\omega_-}\right).$$
It suffices now to recall from Lemma \ref{L:MUI} that the right-hand side equals $0$.
\end{proof}

We next point out
that the Assumption \ref{A:hGWrv} about the behavior of mean 
 total masses of the Galton-Watson trees can be reinforced into a strong limit theorem. In this direction, it is again convenient to
combine the Skorokhod representation and Lemma \ref{Lem:confindim}. In this direction, and up to enlarge the underlying the probability space, under $\mathbb{P}_1$,  assume that we 
 we are  a rescaled family of processes $(f_u^{(n)}, \eta_u^{(n)})_{u\in \U}$ with law ${\P}_{\GW}^{(n)}$, such that $\lim_{n\to\infty} (f_u^{(n)}, \eta_u^{(n)})=(f_u, \eta_u)$ for every $u\in \U$, almost-surely.  For each $n\geq 1$, we write $\texttt{T}^{(n)}$ for the decorated tree constructed from $(f_u^{(n)}, \eta_u^{(n)})_{u\in \U}$
and  set 
$$\upnu^{(n)}(T^{(n)})\coloneqq l^{\varpi}_\GW(n)^{-1} \sum_{u\in \U} \int_0^{ n^{\alpha}  z^{(n)}_u} \varpi\big(n f^{(n)}_u(t n^{-\alpha}  )\big) \dd t.$$ 
Inverting the scaling transformation confirms that $ \upnu^{(n)}(T^{(n)})$ has the same  law  as  $l^{\varpi}_\GW(n)^{-1} \uplambda^{\varpi}_{\GW}(T_{\GW})$ under $\Q_n^\GW$.

\begin{proposition} \label{P:cvnombrepar} In the notation and coupling above, there is the convergence
$$\lim_{n\to \infty} \upnu^{(n)}(T^{(n)}) = \upmu(T) \qquad \text{ in } L^1(\P_1).$$
\end{proposition}

The proof of Proposition \ref{P:cvnombrepar} uses an elementary perturbative result, see Lemma \ref{L:tailored} in the Appendix, akin to a law of large number for uniformly integrable variables; its setting is tailored for our purpose.

\begin{proof}
We will apply Urysohn's subsequence principle, and for this, fix an arbitrary subsequence. The extraction $(n(i))_{i\geq 1}$ that we will construct below is implicitly taken from this subsequence, even though for the sake of simplicity, no further mention to the latter will be made.

For each generation $i\geq 1$, we pick  a finite subpopulation $G_i\subset \N^i$  of individuals at the $i$-th generation, such that
$$\P_1\left( \sum_{u\in \N^i \backslash G_i} \chi(u)^{\omega_-}>2^{-i}\right) \leq 2^{-i}.$$
In particular, recalling Lemma \ref{L:MUI}, we have  by the Borel-Cantelli lemma that
$$\lim_{i\to \infty} \sum_{u\in  G_i} \chi(u)^{\omega_-} = \upmu(T), \qquad \text{a.s. and in }L^1(\P_1),$$
where the convergence in $L^1(\P_1)$ follows again from the Scheff\'e lemma.

Since $l^{ \varpi}_\GW$ is regularly varying with exponent $\omega_-$,  we have $\lim_{n\to \infty} f^{(n)}_u(0)=\chi(u)$ for every $u\in \U$, and since $G_i$ is finite, 
we can extract an increasing subsequence $(n(i))_{i\geq 1}$ such that, for each $i$, 
$$\P_1\left( \sum_{u\in  G_i} \left| \chi(u)^{\omega_-}- a_{i,u}\right|>2^{-i}\right) \leq 2^{-i},$$
where 
$$a_{i,u}\coloneqq \frac{l^{ \varpi}_\GW(n(i)f^{(n(i))}_u(0))}{l^{ \varpi}_\GW(n(i))}, \qquad \text{ for }u\in G_i.$$
 Observe  that  $(a_{i,u}: u\in G_i \text{ and }i\geq 1)$ is a null array of nonnegative random variables with (again by the Borel-Cantelli lemma)
$$\lim_{i\to \infty} \sum_{u\in G_i} a_{i,u}= \upmu(T) \qquad \text{a.s.}$$
Recalling from \eqref{eq:ellsuperh}
that $l^{ \varpi}_\GW$ is superharmonic, we must have $\E\left( \sum_{u\in G_i} a_{i,u} \right)\leq 1$ for all $i\geq 1$, and since $\E_1(\upmu(T))=1$, another application of the Scheff\'e lemma entails that
the convergence above also holds in $L^1(\P_1)$. 

Next, for every $i\in \N$ and every individual $u\in G_i$, we consider the subtree of $T^{(n(i))}$ that stems from $u$, and denote by $b_{i,u}$ its normalized mass, where the normalization is chosen such that $b_{i,u}$ has a unit mean. That is explicitly, writing for simplicity  $k=n(i)f_u^{(n(i))}(0)$ and $n=n(i)$,
$$ b_{i,u}= l^{\varpi}_\GW(k )^{-1} \sum_{u\preceq v } \int_0^{ n^{\alpha}  z^{(n)}_v} \varpi(n f^{(n)}_v(t n^{-\alpha}  )) \dd t.$$
By the branching property of the Galton--Watson process,  if the (rescaled) type of $u$ for $\texttt{T}^{(n(i))}$ is, say $f_u^{(n(i))}(0)=a$ with $a=k/n(i)$ for some $k\geq 1$, then the conditional distribution of $b_{i,u}$ given the types of all individuals in $G_i$
has the same law as  $\upbeta_{i,a}$  where, 
$$\upbeta_{i,a} \text{ has the law of }\uplambda^{\varpi}_{\GW}(T_{\GW})/l^{ \varpi}_\GW(n(i)a) \text{ under }\P^{\GW}_{n(i)a}.$$
We further agree for definitiveness that  $\upbeta_{i,a}\equiv 1$  if $n(i)a\not \in \N$.

By Assumption \ref{A:hGWrv} and the branching property at generation $i$ for the Galton--Watson process,  we are now in the framework of Lemma \ref{L:tailored}, from which we get that 
$$\lim_{i\to \infty} \sum_{u\in G_i} a_{i,u}b_{i,u}= \upmu(T) \qquad \text{in } L^1(\P_1).$$
Since by construction, $ \sum_{u\in G_i} a_{i,u}b_{i,u}\leq \upnu^{(n(i))}(T^{(n(i))})$,  and since further 
$$\E_1(\upnu^{(n(i))}(T^{(n(i))}))= \E_1(\upmu(T))=1,$$ we conclude from a version of Scheff\'e's lemma that 
the sequence $\upnu^{(n(i))}(T^{(n(i))})$ converges as $i\to \infty$ towards  $\upmu(T)$, in $L^1(\P_1)$. 

As a summary, we have extracted  from an arbitrary subsequence a further subsequence $(n(i))_{i\geq 1}$ along which the claim holds. Thanks to  Urysohn's subsequence principle, this  is all what was needed. 
\end{proof}

We are now in a position to tailor the strategy from Section  \ref{Sub:sec:prt1}. The key step is to establish the convergence of the pruned rescaled  trees equipped with adequate measures. We use the same notation and setting as for Proposition \ref{P:cvnombrepar}, and recall that the rescaled Galton--Watson tree $T^{(n)}$ is equipped with the measure
$$\upnu^{(n)}\coloneqq l^{\varpi}_{\GW}(n)^{-1}n^{\alpha}  \varpi(n g)\cdot \uplambda_{T^{(n)}}.$$

Next, recall the convergence of cutoff approximations in Corollary \ref{lem:cutoff}(ii) for self-similar Markov trees,  the notation 
 $\normalfont{\texttt{T}}^{[r]}$ for the decorated  subtree pruned below the threshold $r$ and the definition \eqref{eq:def:nu:varepsilon} given there  for the projection $\upmu^{[r]}$ of the harmonic measure $\upmu$ on $\texttt{T}^{[r]}$. This  incites us to write similarly $\texttt{T}_n^{[r]}$ for the decorated  subtree of $\texttt{T}^{(n)}$ that results from pruning below the threshold $r$, and then 
 set
 $$\upmu^{[ r]} _n\coloneqq   \sum_{\uptau} \upnu^{(n)}(\uptau)\delta_{\uprho_{\uptau}},$$
where in the right-hand side,   the sum is over the family of subtrees $\uptau$ of $T^{(n)}$ dandling from $T_n^{[r]}$ and $\uprho_{\uptau}$ stands for the root of $\uptau$.
Plainly, the roots  of the dandling trees belong to the truncated tree $T_n^{[r]}$, and $\upmu^{[ r]}_n$ is a measure on $T^{[r]}_n$ with total mass at most $\upnu^{(n)}(T^{(n)})$. 

\begin{corollary}\label{Cor:convpruned:2}   Fix a threshold $0<r <1$ such that
  $$\P_1(\chi(u)=r)=0 \qquad \text{for any  }u\in \U.$$
  In the notation above, the sequence of  measured  truncated decorated trees $(\normalfont{\texttt{T}}^{[r]}_n, \upmu_n^{[ r ]})$  converges as $n\to \infty$ to  
 $(\normalfont{\texttt{T}}^{[r]}, \upmu^{[ r]})$,  where $\upmu^{[r]}$ is given by \eqref{eq:def:nu:varepsilon} in terms of the harmonic measure $\upmu$ of self-similar Markov tree.  
\end{corollary}

\begin{proof} The argument is similar to that for Corollary \ref{Cor:convpruned},  with Proposition \ref{P:cvnombrepar} being the key to the convergence of the measure components. 
To avoid lengthy repetitions, we merely indicate the main lines and leave a few technical details to the interested reader.

Recall from Notation \ref{N:subtrees}(i) that $\texttt{T}_u=(T_u, d_{T_u}, \uprho(u), g_u)$ is the decorated real subtree generated by an individual $u\in \U$,
so that defining $\uptau_u=\texttt{T}_u$ for $u\in F(r)$ and agreeing that $\uptau_u$ is degenerate otherwise yields a natural indexation by $\U$ of the family of the subtrees dandling from $T^{[r]}$ 
(this appeared already in the proof of Proposition \ref {prop:cutoffmarkov}). We view  the roots  of the dandling trees as marks (possibly fictitious) indexed by $\U$ on $T^{[r]}$, and  further record the germ of the decoration of dandling trees, that is the type $\chi(u)$ of the individual $u\in F(r)$; see \eqref{eq:labelatroot}. We shall first extend
Corollary \ref{Cor:convpruned}  to include those marks and germs of decoration. 
  
Recall that we are given a family of decoration-reproduction processes $(f_u, \eta_u)_{u\in \U}$ with law $\P_1$ and 
a rescaled family of processes $(f_u^{(n)}, \eta_u^{(n)})_{u\in \U}$ with law ${\P}_{\GW}^{(n)}$, such that  $\lim_{n\to\infty} (f_u^{(n)}, \eta_u^{(n)})=(f_u, \eta_u)$ for every $u\in \U$, almost-surely. We write respectively $F^{(n)}(r)$ and $F(r)$ for the lines of frozen individuals at  threshold $r$, and stress that  as $n\to \infty$, the sequence of indicator functions
$\indset{F^{(n)}(r)}$ on $\U$ converges pointwise to $\indset{F(r)}$, almost surely.
Then let  $\texttt{T}^{[r]\bullet}_n$ denote the decorated tree constructed from $(f_u^{(n)}, \eta_u^{(n)})_{u\in \U}$, 
 pruned below the threshold $r$, and with marks at the roots of the dandling trees. Write similarly $\texttt{T}^{[r]\bullet}$ for the same object constructed from the family $(f_u, \eta_u)_{u\in \U}$. Then as $n\to \infty$, $\texttt{T}^{[r]\bullet}_n$ converges a.s. to $\texttt{T}^{[r]\bullet}$ on  the space of  decorated real trees with marks (see the end of Section \ref{sec:1.3} for the distance on that space). The argument to justify this assertion is similar to that in the proof of Corollary \ref{Cor:convpruned}, except that it uses an immediate extension of Lemma \ref{L:convgluing} where the trees $T_i$ for $i=0,1, \ldots,  k$ further carry marks.
 Observe that there is also the a.s. convergence of the germs of the decoration of the dandling trees, since 
 $\lim_{n\to \infty} f_u^{(n)}=f_u$ a.s. for all $u\in \U$.

It now suffices to verify that
\begin{equation}\label{Eq:prescheff}
\lim_{n\to \infty} \E_1\left( \sum_{u\in \U}\left| \indset{u\in F(r)} \upmu(T_u)- \indset{u\in F^{(n)}(r)} \upnu^{(n)}(T^{(n)}_u)\right| \right) =0,
\end{equation}
as indeed this then implies the convergence of $\upmu_n^{[ r ]}$ towards $\upmu^{[ r ]}$ in the Gromov-Prokhorov sense as $n\to \infty$. 
We readily deduce from above and Proposition \ref{P:cvnombrepar}  that for each $u\in \U$,
$$\lim_{n\to \infty} \E_1\left( \left| \indset{u\in F(r)} \upmu(T_u)- \indset{u\in F^{(n)}(r)} \upnu^{(n)}(T^{(n)}_u)\right| \right) =0.$$
On the one hand, since the harmonic measure $\upmu$ gives no mass to the skeleton of $T$, and \textit{a fortiori} to the pruned tree $T^{[r]}$, we have 
$$ \sum_{u\in \U}\indset{u\in F(r)} \upmu(T_u)= \upmu(T\backslash T^{[r]}) = \upmu(T),$$
and this quantity has expectation $1$.
On the other hand, we have also
$$ \sum_{u\in \U}\indset{u\in F^{(n)}(r)} \upnu^{(n)}(T^{(n)}_u)\leq  \upnu^{(n)}(T^{(n)}),$$
and, by definition, the expectation of right-hand side equals $1$ for every $n\geq0$. A final application of the Scheff\'e lemma yields \eqref{Eq:prescheff}, and hence completes the proof.
\end{proof}

We can now combine Corollary \ref{Cor:convpruned:2} with the convergence \eqref{Eq:borneunif}  to deduce Theorem  \ref{T:mainunconds-mass}.
\begin{proof}[Proof of Theorem  \ref{T:mainunconds-mass}]
We will use the following elementary fact that is seen from the very definition of the Prokhorov distance. Consider a real tree $\tau$ equipped with a finite measure $\nu$, and
$\tau'\subset \tau$ a subtree of $\tau$ that we equip with the push-forward image $\nu'$ of $\nu$ by the projection from $\tau$ to $\tau'$. Then the Prokhorov distance $\d_{\mathrm {Prok}}(\nu, \nu')$ is bounded from above by the Hausdorff distance $\d_{\mathrm {Haus}}(\tau, \tau')$.

Fix some arbitrarily small $\varepsilon >0$. Thanks to Corollary \ref{lem:cutoff} and
 \eqref{Eq:borneunif}, we can take $r>0$ sufficiently small with  $\P_1(\chi(u)=r)=0 $ for any  $u\in \U$, such that 
  both
 $$\P_1\left( \d_{\TT_m}((\texttt{T}, \upmu), ( \texttt{T}^{[r]}, \upmu^{[r]}))>\varepsilon \right)< \varepsilon.$$
 and
$$\P_1\left( \d_{\TT_m}(\texttt{T}^{(n)}, \texttt{T}^{[r]}_n)>\varepsilon \right)< \varepsilon.$$
Recall that $T^{(n)}$ is endowed with the measure $\upnu^{(n)}$, and observe that its image by the projection from $T^{(n)}$ to $T^{[r]}_n$ is given, in the notation of 
Corollary \ref{Cor:convpruned:2}, by $\upmu^{[r]}_n + \indset{T^{[r]}_n} \cdot \upnu^{(n)}$. 
Since we know from Lemma \ref{L:asympneglbr} and the branching property that the mass assigns by $\upnu^{(n)}$ to any segment $S_u$ becomes asymptotically negligible 
as $n\to \infty$, we deduce that $\lim_{n\to \infty} \upnu^{(n)}(T^{[r]}_n)=0$ in probability.

We thus have
$$\P_1\left(\upnu^{(n)}(T^{[r]}_n) > \varepsilon \right) < \varepsilon \qquad \text{ for all sufficiently large }n,$$
and hence also (from the elementary fact observed at the beginning of the proof),
 $$\P_1\left( \d_{\TT_m}(( \texttt{T}^{(n)}, \upnu^{(n)}), (\texttt{T}^{[r]}_n, \upmu^{(n)}))>2\varepsilon \right)< 2\varepsilon \qquad \text{ for all sufficiently large }n.$$
We then conclude from Corollary \ref{Cor:convpruned:2} that
 $$\P_1\left( \d_{\TT_m}(( \texttt{T}, \upmu), (\texttt{T}^{(n)}, \upnu^{(n)}))>3\varepsilon \right)< 3\varepsilon \qquad \text{ for all sufficiently large }n,$$
and the third version of our invariance principle is proven.
\end{proof}

We finish this chapter by stating a weaker version of Assumption \ref{assum:main:gen:deco:repro} in presence of Assumption \ref{A:all:existance:PHI}. This will be in particular used in the proof of Theorem \ref{Th:cvdecoreprod} and Corollary \ref{cor:cactus}. Recall the setup of Assumption \ref{assum:main:gen:deco:repro}. For any $  \varepsilon>0$, if $ {P}$ is a law on decoration-reproduction process $(f, \eta)$ started from $F(0)=1$, we denote by  $ {P}^{ \varepsilon}$ the law of the decoration-reproduction process killed at the first instant it drops below level $ \varepsilon$ that is the law under $ {P}$ of 
$$ ( f|_{[0 , \theta^{ \varepsilon}]}, \eta|_{[0, \theta^{ \varepsilon}]}), \quad \mbox{ where } \theta^{ \varepsilon} = \inf\{u \geq 0 : f(u) \leq \varepsilon\}.$$
\begin{lemma}  \label{lem:convreprodecolight} Suppose Assumption \ref{A:all:existance:PHI}. If for every $ \varepsilon>0$ we have 
$$ P^{(n), \varepsilon}_{\GW} \xrightarrow[n\to\infty]{(d)} P^{ \varepsilon}_{1},$$ in the sense of the topology described just before Assumption \ref{assum:main:gen:deco:repro}, then   Assumption \ref{assum:main:gen:deco:repro} holds.
\end{lemma}
\begin{proof} In order to pass from the convergence of $P^{(n), \varepsilon}_{\GW}$ to $P^{ \varepsilon}_{1}$ to the convergence of $P^{(n)}_{\GW}$ to $P_{1}$, one needs to prove that after dropping below level $ \varepsilon$, the decoration die quickly while staying small and do not give rise to non-small reproduction. Both are ensured uniformly in $n$ (and so at the limit) by Lemma \ref{L:Aspandii}. We leave the easy adaptation  to the reader.
\end{proof}

\section*{Comments and bibliographical notes}
Seeing a (multi-type) Galton--Watson tree as a Crump-Mode-Jagers branching system where individuals represent lineages of particles is a fruitful idea that already appeared in different forms in the literature. For example, a version of this transformation already appears in \cite{minami2005number} when relating a monotype Galton--Watson tree conditioned by the number of leaves to another Galton--Watson process conditioned on the total progeny. The proofs of our invariance principles build upon many previous works including \cite{BCK18,BBCK18} and \cite{dadoun2019self}. In \cite{dadoun2019self}, Dadoun circumvented the use of global superharmonic functions to the cost of pruning the particles dropping below a certain level. Notice also, that none of the previous works considered the measure or the decoration.

\chapter{Analytical criteria for the invariance principle}\label{C:scaling}
 
 \newcommand{\kappan}{\kappa^{(n)}}

 The three versions of the invariance principle established in the preceding chapter  depend crucially on a scaling limit for the decoration-reproduction induced by the Galton--Watson process and a selection rule, namely Assumption  \ref{assum:main:gen:deco:repro}, as well as further technical conditions
 (Assumptions  \ref{A:all:existance:PHI} and \ref{A:hGWrv}). Our purpose here is to provide more explicit analytic criteria in terms of the reproduction kernel of the  Galton--Watson process that ensure these assumptions. In a nutshell, as shown in the diagram of Figure \ref{fig:diagram}, we demonstrate that the probabilistic Assumptions \ref{assum:main:gen:deco:repro} and \ref{A:all:existance:PHI} are implied by the analytical Assumptions \ref{A:all:GW}, \ref{A:typex} and \ref{A:BK}.  Just as in the preceding chapter,  we fix some  characteristic quadruplet  $(\sigma^{2}, \mathrm a,  \boldsymbol{\Lambda}; \alpha)$ and assume without further mentions that its cumulant $\kappa$ satisfies  Assumption \ref{A:gamma0} for some $\gamma_0>0$. 
We suppose also for simplicity that
 \begin{equation} \label{eq:tech}  \Lambda_0(\{-1,1\}) =0.  \end{equation} 
  This is a very mild restriction  that stems  from the discontinuities at $\pm 1$  of the indicator function $\indset{|y|<1}$ appearing in the L\'evy-Khintchin formula \eqref{E:LKfor}. It is made for the sake of simplicity only and can always be easily circumvented.

\section{Scaling limit for the decoration-reproduction process} \label{sec:7.1}

Throughout this chapter, we will assume that, no matter the type of the ancestral particle, the expected number of particles of type $1$ that ever exist is both strictly positive and finite.
This can be interpreted as a requirement of  accessibility  and sub-criticality of type $1$ for the Galton--Watson process.
  \begin{assumption}[Type $1$ is accessible and sub-critical] \label{A:all:GW} For any $j\geq 1$, we have 
 \begin{equation}\label{eq:non:explo:1}
0< \E_j^{\GW}\Big(\sum \limits_{k\geq 0} Z_1(k)\Big)= \sum_{k=0}^{\infty} \boldsymbol{m}^k_{\GW}(j,1)<\infty.
\end{equation}
 \end{assumption}
In other words, the above ensures the finiteness and the positivity of the ``smallest possible'' super-harmonic function $ l_{\GW}^{ \mathbf{1}_{ \{ 1\}}}$ for the Galton--Watson process. Checking Assumption \ref{A:all:GW}  for multi-type Galton--Watson case is in general a hard task. However, in many situations of interest, it is granted from combinatorial arguments. The next   assumptions  are finer but  only asymptotic when the type of the ancestral particle goes to infinity. 

To state these assumptions, recall the notation \eqref{Eq:reprodkernelGW} for  the reproduction kernel of the Galton--Watson process $\boldsymbol{\pi}^{\GW}_n$,  and let $\boldsymbol{\Lambda}^{(n)}$ denote the $n^{\alpha}$ multiple of its push-forward image of   by 
 the logarithmic transformation 
 $$\log(n^{-1} \cdot): \mathbb{Z}_+\times \U^{\downarrow} \to \mathcal{S}=[-\infty, \infty)\times \mathcal S_1$$  defined by 
$$ \big(n_0,(n_1, \ldots, n_k, )\big) \mapsto \left(\log  \frac{n_{0}}{n} ,\left(\log \frac{n_1}{n}, ...  , \log \frac{n_k}{n}\right)\right).$$
Specifically, for every measurable function   $ {f} :  \mathcal S \to \R$, there is the identity 
\begin{equation} \label{eq:deflambdan} 
   \int_{  \mathcal{S}}\boldsymbol{\Lambda}^{(n)} (\mathrm{d}  \mathbf{y}) {f}(\mathbf{y}) = n^{\alpha } \cdot \E^\GW_n\Big( {f}\big(\log\big(n^{-1}\varsigma(\mathbf{Z}(1))\big)\big)\Big).
      \end{equation} We also write $\Lambda^{(n)}_0$ for the image of   $ \boldsymbol{\Lambda}^{(n)}$
  by the first projection $[-\infty, \infty)\times \mathcal S_1\to [-\infty, \infty)$. Last, recall that $(-\infty, -\infty, \ldots)\in \mathcal{S}_1$ is interpreted as the empty sequence and $\mathbf 0\coloneqq (0,(-\infty, -\infty, \ldots))$ as the neutral element in $\mathcal{S}$.  
   Our first asymptotic  assumption  requests the convergence of a rescaled version of the reproduction kernel of the Galton--Watson process towards the 
  infinitesimal characteristics of the self similar Markov tree, and  is crucial for establishing the convergence in distribution of the rescaled decoration-reproduction processes.
\begin{assumption}[Convergence of characteristics] \label{A:BK} 
There is the convergence of the generalized L\'evy measures
$$\lim_{n\to \infty}
\int_{ \mathcal{S}}\boldsymbol{\Lambda}^{(n)} ( \mathrm{d}  \mathbf{y}) {f}( \mathbf{y})  = \int_{ \mathcal{S}}\boldsymbol{\Lambda} ( \mathrm{d}  \mathbf{y}) {f}(\mathbf{y}),$$
where $ {f} : \mathcal{S}\to \mathbb{R}$ denotes a generic bounded continuous function that vanishes on some neighborhood of $\mathbf 0$ and admits a limit at $(-\infty,(-\infty, -\infty, ...))$. There is furthermore  convergence 
 of the drift and variance,
 $$\lim_{n\to \infty}  \int_{-1}^1 \Lambda_0^{(n)}( \mathrm{d}y) y= \mathrm a \quad\text{and}
 \quad \lim_{n\to \infty}   \int_{-1}^1  \Lambda_0^{(n)}( \mathrm{d}y) y^2=  \sigma^2 + \int_{-1}^1  \Lambda_0( \mathrm{d}y) y^2,
 $$
 where $(\sigma^{2}, \mathrm a,\boldsymbol{\Lambda}; \alpha)$ is a characteristic quadruplet  satisfying Assumption \ref{A:gamma0}.
 \end{assumption}
  
Our second asymptotic assumption can be interpreted as a discrete analogue—and strengthening—of Assumption~\ref{A:gamma0}. To formalize it,  by analogy with the continuous setting, we define the \textbf{discrete cumulant function} as
  \begin{eqnarray}\label{E:discretecumul}  \kappan(\gamma) &\coloneqq & n^{\alpha-\gamma} \boldsymbol{G}_{\GW} ( n\mapsto n^{\gamma}) \\ &=&    n^{\alpha} \sum_{(k_{0},(k_{1}, k_{2},...)) \in \mathbb{Z}_+ \times \mathbb{U}^{\downarrow}}  \boldsymbol{\pi}^{\GW}_{n}\big(k_{0},( k_{1}, k_{2}, \dots)\big) \left( \sum_{i=0}^{\infty} \left(\frac{k_{i}}{n}\right)^{\gamma} -1 \right)\nonumber \\
   &=&  \int_{ \mathcal{S}} \boldsymbol{\Lambda}^{(n)} ( \mathrm{d} \mathbf{y}) \left( \sum_{i \geq 0} \mathrm{e}^{\gamma y_i}-1\right).  \end{eqnarray} 
   Recall from Assumption \ref{A:gamma0} that $\gamma_0>0$ has been chosen such that $\kappa(\gamma_0)<0$. 
\begin{assumption}[Negativity of discrete cumulants]\label{A:typex} There exists $0<\gamma_1<\gamma_0$ such that 
$$ \limsup \limits_{n\to \infty}\kappan(\gamma_1)<0\quad \text{ and } \quad \limsup \limits_{n\to \infty}\kappan(\gamma_0)<0.$$
 \end{assumption} 
We now outline the structure of the remainder of the chapter and explain how we combine the assumptions to derive the key  Assumptions \ref{assum:main:gen:deco:repro} and  \ref{A:all:existance:PHI}. In Section \ref{sec:superharmonic}, we begin by studying superharmonic functions and establish the following result:
 \begin{proposition}\label{prop:PHI:ex}
  Assumption \ref{A:all:existance:PHI} follows from Assumptions \ref{A:all:GW} and \ref{A:typex}. 
  \end{proposition}

Next, in Section \ref{sec:proofscaldecrep}, we turn to the main result of the chapter. This result can be viewed as an extension of a theorem by Bertoin and Kortchemski \cite{BK14} concerning scaling limits of Markov
 chains with values in $\mathbb{N}$, which itself generalizes earlier work by Haas and Miermont \cite{haas2011self} on decreasing Markov chains. In our context, it takes the following form:
\begin{theorem} \label{Th:cvdecoreprod} 
Assumption \ref{assum:main:gen:deco:repro} follows from Assumptions \ref{A:all:GW}, \ref{A:BK}, and \ref{A:typex}.
 \end{theorem}
Finally, in Section \ref{sec:minimal:reverse}, we conclude our analysis of superharmonic functions by proving lower bounds, which are of independent interest. Throughout the remainder of the chapter, we assume that Assumption \ref{A:all:GW} holds, as it underpins all subsequent arguments.

\section{Upper bounds on superharmonic functions}
\label{sec:superharmonic}
Our main purpose is to prove  Proposition \ref{prop:PHI:ex} (i.e. that Assumption \ref{A:all:existance:PHI} is implied by Assumptions~\ref{A:all:GW}~and~\ref{A:typex})  and construct strictly superharmonic functions with polynomial growth.  
Of course, recall not only Assumption \ref{A:all:existance:PHI}, but also the notation and definition given right before,
notably \eqref{Eq:Laplacediscret} for the generator.  We assume throughout this section and without further mentions  that  Assumptions \ref{A:all:GW} is satisfied.

We  have already observed in \eqref{eq:ellsuperh} that the   mean size functions $ l^{\varpi}_{\GW}$ defined  by \eqref{def:lA} for some weight function $\varpi$, are 
naturally superharmonic provided that they are finite everywhere. We now point out that this is always the case when $\varpi$ has finite support. 

\begin{lemma} 
 \label{eq:bis:finite}
For every $\varpi:\mathbb{N}\to \mathbb{R}_+$ with finite support, the function $l_{\GW}^{\varpi}$ takes finite values, and as a consequence, is superharmonic. Moreover $l_{\GW}^{\varpi}(n)=O(n^{\gamma})$ for any $\gamma>0$ with 
$ \limsup \limits_{n\to \infty}\kappan(\gamma)<0$.
\end{lemma}
\begin{proof} It suffices  to check the first assertion for the weight function $\varpi= \indset{\{j\}}$,  where  $j\geq 1$ is arbitrary. 
In this direction, write $\mathcal Z_j(k)$ for the number of particles of type $j$ in the entire Galton--Watson process such that exactly $k$ of its strict predecessors  are also particles of type $j$. Obviously, there is the identity
$$\sum_{k\geq 0} \mathcal Z_j(k) = \bar Z_j\coloneqq \sum_{n=0}^{\infty} Z_j(n).$$
The  strong branching property \cite[Theorem 4.14]{jagers1989general} of Galton--Watson processes at stopping lines  readily entails 
that under $\Q_\ell^{\GW}$, the process $(\mathcal Z_j(k))_{k\geq 0}$ is a monotype Galton--Watson process
started from the random variable $\mathcal Z_j(0)$ (of course $\mathcal Z_j(0)=1$, $\Q_j^{\GW}$-a.s.). 

We first focus on the type  $j=1$.  The sub-criticality and accessibility Assumption \ref{A:all:GW} yields (again by  the branching property of Galton--Watson processes at stopping lines) that
 $$\E^{\GW}_1( \mathcal Z_1(0) ) < 1  \text{ and } 
0< \E^{\GW}_j( \mathcal Z_1(0) ) < \infty \text{ for all }j\geq 2.$$
On the other hand, by decomposing the Galton--Watson process at the stopping line formed by particles of type $j$ that have exactly $k$ strict predecessors of type $j$ 
and applying once again  the branching property  at stopping lines, we get that 
$$\E_\ell^{\GW}\Big(\sum_{n=0}^{\infty} Z_1(n) \Big) \geq \E_\ell^{\GW}\Big(\sum_{k=0}^{\infty} \mathcal{Z}_j(k) \Big) \cdot \E^{\GW}_j\big( \mathcal Z_1(0) \big),$$
for any $j\geq 2$ and $\ell\geq 1$.\footnote{The previous inequality is actually an equality when $\ell = j$.} We know from Assumption \ref{A:all:GW} that the left-hand side is finite, and we have observed above that $\E^{\GW}_j( \mathcal Z_1(0) )>0$. 
 We conclude that
 $$\E_\ell^{\GW}\Big(\sum_{k=0}^{\infty} \mathcal{Z}_j(k) \Big)  = \E_\ell^{\GW}\Big(\sum_{n=0}^{\infty} Z_j(n) \Big) <\infty,$$
which is what we wanted to verify.

  To prove the second claim of the lemma,  introduce the function $\Phi(n)=n^{\gamma}$ for $n\geq 1$ and observe from \eqref{E:discretecumul}  that 
$$\limsup \limits_{n\to \infty} n^{\alpha-\gamma} \boldsymbol{G}_{\GW} \Phi(n)=\limsup \limits_{n\to \infty} \kappan(\gamma)<0.$$
Since $\varpi$ has finite support and $l_{\GW}^{\varpi}$ takes finite values, we can find $N\geq 1$  and $c>0$ large enough, such that $\boldsymbol{G}_{\GW} \Phi(n)< -\varpi(n)$ for every $n\geq N$, and $c  \cdot \Phi(n)\geq l_{\GW}^{\varpi}(n)$ for $n\leq N$. An application of Lemma~\ref{lem:minimal} with $B=\{1,\dots, N\}$ gives that $ l_{\GW}^{\varpi} \leq c \cdot \Phi$ everywhere. This completes the proof of the lemma.
\end{proof}

We are now able not only to verify Proposition \ref{prop:PHI:ex} but also to upper-bound the growth of fairly general functions $l_{\GW}^{\varpi}$.

\begin{proposition}\label{prop:APHI}  For any  $\gamma>0$ with $\limsup_{n\to \infty}\kappan(\gamma)<0$, there exists $\phi:\mathbb{N}\to \mathbb{R}_+^*$ strictly superharmonic such that $\phi(n) \asymp n^{\gamma}$ and 
$$ \limsup \limits_{n\to \infty} n^{\alpha-\gamma}\boldsymbol{G}_{\GW} \phi(n)<0.$$ Moreover, for every weight function $\varpi:\mathbb{N}\to \mathbb{R}_+$ such that $\varpi(n)=\mathcal{O}(n^{\gamma-\alpha})$, the function $l_{\GW}^{\varpi}$ is finite everywhere (and hence superharmonic), and satisfies  $l_{\GW}^\varpi(n)=\mathcal{O}(n^{\gamma})$.
\end{proposition}
Before proving the proposition, let us explain why it implies Proposition \ref{prop:PHI:ex}. If we further assume that Assumption \ref{A:typex} is verified in addition to Assumption \ref{A:all:GW}, then we can simply apply Proposition \ref{prop:APHI} for $\gamma_0$ and $\gamma_1$, which gives the desired superharmonic functions of Assumption~\ref{A:all:existance:PHI}.
\begin{proof}
Fix $\gamma >0$ such that $\limsup_{n\to \infty}\kappan(\gamma)<0$. Let again $\Phi(n)=n^{\gamma}$, and note that  $\limsup_{n\to \infty} n^{\alpha-\gamma}\boldsymbol{G}_{\GW} \Phi(n) <0$. Hence $\Phi$ satisfies  the desired properties on a co-finite subset of $\N$, and we shall  construct the function $\phi$ from $\Phi$. To this end, we consider $N\geq 1$ such that $\boldsymbol{G}_{\GW} \Phi(n)<0$ for all $n\geq N$, and let $\varpi^{\circ}=\mathbf{1}_{\{1, \dots, N\}}$. Since $l_{\GW}^{\varpi^\circ}$ is finite by Lemma \ref{eq:bis:finite}, we have from  \eqref{eq:ellsuperh} that $\boldsymbol{G}_{\GW} l_{\GW}^{\varpi^\circ}=-\varpi^\circ$. Thus, we can find $c>0$ sufficiently large such that 
$$\boldsymbol{G}_{\GW}(\Phi+ c   l_{\GW}^{\varpi^\circ})(n)= \boldsymbol{G}_{\GW} \Phi(n) - c  \varpi^\circ(n)<0,\quad \text{ for } n\geq 1. $$
In particular, the function  $\phi= \Phi+ c   l_{\GW}^{\varpi^\circ}$ is superharmonic on $\mathbb{N}$ and  satisfies
$$ \limsup \limits_{n\to \infty} n^{\alpha-\gamma} \boldsymbol{G}_{\GW} \phi(n)= \limsup \limits_{n\to \infty} n^{\alpha-\gamma} \boldsymbol{G}_{\GW} \Phi(n)<0.$$
Lemma \ref{eq:bis:finite} now yields $l_{\GW}^{\varpi^\circ}(n)=\mathcal{O}(n^{\gamma})$, and hence $\phi(n)\asymp n^{\gamma}$. 

Finally consider a weight function $\varpi:\mathbb{N}\to \mathbb{R}_+$ with $\varpi(n)=\mathcal{O}(n^{\gamma-\alpha})$.  We can always find $K>0$, such that both  $K  \phi$ and $-K \boldsymbol{G}_{\GW} \phi$ are bounded from below by $\varpi$. Therefore, we can again apply Lemma  \ref{lem:minimal}, this time with $B=\varnothing$, to deduce that $l_{\GW}^\varpi\leq K \phi$ everywhere. In particular, $l_{\GW}^\varpi$ takes finite values and satisfies  $l_{\GW}^\varpi(n)=\mathcal{O}(n^{\gamma})$.
\end{proof}

\section{Proof of Theorem \ref{Th:cvdecoreprod}}\label{sec:proofscaldecrep}

The title explains the purpose of this section, but hides the fact that the task has rather technical aspects. For the sake of simplicity, we will focus on the case when there is no killing, $\texttt{k}=0$, and will briefly indicate at the end the simple changes needed to deal with the case $\texttt{k}>0$. 
We suppose throughout this section that Assumptions \ref{A:BK} and  \ref{A:typex} hold, in addition to  Assumption~\ref{A:all:GW}.

Let $(X, \eta)$ denote a self-similar Markov decoration-reproduction process with law $P=P_1$ as constructed in Section \ref{sec:2.2}. We suppose that its characteristic quadruplet $(\sigma^2, \mathrm{a}, \boldsymbol{\Lambda} ; \alpha)$  satisfies Assumption \ref{A:gamma0} and has killing rate $\texttt{k}=0$. Recall the notation $P^{(n)}_{\GW}$ for the image by the scaling transformation \eqref{Eq:scalefeta}
of the law $P_n^{\GW}$ of  the decoration-reproduction process induced by a selection rule for a Galton--Watson process 
when the ancestor has  type $n\ge 1$. Up to Skorokhod's embedding, Theorem \ref{Th:cvdecoreprod} amounts to the following strong convergence claim.

\begin{proposition} \label{P:cvdecoreprod} For every $n\geq 1$, there exist a rcll process
$f^{(n)}: [0,z^{(n)}]\to \R_+$ and a point process $\eta^{(n)}$ on $[0,z^{(n)}]\times (0,\infty)$  
with joint law  $P_\GW^{(n)}$, such that the following holds as $n\to \infty$:

\begin{enumerate}
\item[(i)] $f^{(n)}$ converges a.s.  in $\mathbb D^{\dagger}$ to the  self-similar Markov process $X$,

\item[(ii)]  $ \eta^{(n)}$ converges a.s.  in the sense of vague convergence of Radon measures on $[0,\infty)\times (0,\infty]$ to the reproduction process $\eta$. 
\end{enumerate}
\end{proposition}

The main idea for the proof of Proposition \ref{P:cvdecoreprod} is borrowed from \cite{BK14}.
Rather than dealing directly with decoration-reproduction processes, we consider the ratios of the types of children particles to their parents. We then perform a logarithmic transformation and a convenient  time-change such  that  as the type of the ancestral particle goes to infinity, this transformed version converges in distribution to a L\'evy process endowed with a Poisson point process on $\mathcal S_1$ as in Section \ref{sec:2.2}. We are left to check that the inverse transformation corresponds asymptotically to the Lamperti transformation.

Specifically, working implicitly under the law $P^{(n)}_{\GW}$  for an arbitrary $n\geq 1$, we introduce a real-valued process $L^{(n)}$ and a measure-valued process
$A^{(n)}$ on $ \R$ defined gradually from $(f^{(n)},\eta^{(n)})$ as follows.
First, given the decoration process $f^{(n)}$, we replace the holding times with durations $n^{-\alpha}$ by independent exponential times with parameter $k^{\alpha}$ when $f^{(n)}$ is at state $k/n$ with $k\in \N$, and then obtain $L^{(n)}$ by applying the logarithm transformation. That is, let $(e_j)_{j\geq 1}$ denote an independent sequence of i.i.d. standard exponential variables,  introduce recursively the sequence of jump times
$$t_0= 0, \quad t_j=t_{j-1} + (nf^{(n)}(j-1))^{-\alpha} e_j \quad \text{for }j\geq 1,$$
and set 
$$L^{(n)}_t=\log (f^{(n)}(t_{j-1})) \quad \text{ for }t_{j-1} \leq t < t_j.$$ 
Hence, 
$L^{(n)}$ is a Markov jump process on the state space $\log(n^{-1}\N)$,  its jump chain is 
given by the successive values of $\log f^{(n)}$. More precisely,
$L^{(n)}$ 
starts from $0$ and is absorbed at $-\infty$, and has jump rate $k^{\alpha} \boldsymbol{\pi}^{\GW}_{k}(\ell, \U^{\downarrow})$ from $\log(k/n)$ to $\log(\ell/n)$. 

Second,  introduce the  integer-valued measure $\Updelta^{(n)}(\dd t, \dd y)$ on $\R_+\times \R$ such that, 
for every $1\leq j <z$, if $k=nf^{(n)}(n^{\alpha}(j-1))$ denotes the type of the legitimate particle at generation $j-1$ and  $v=(v_1, \ldots,v_\ell)\in \U^{\downarrow}$  the ranked sequence of the types of its illegitimate progeny (if any),  then  $\Updelta^{(n)}$  has $\ell$  atoms located at 
$$(t_j, \log(v_i/k))= \left (t_j, \log(v_i/n)-L^{(n)}_{t_j-} \right), \qquad \text{for }i=1, \ldots, \ell,$$ and there are no other atoms of $\Updelta^{(n)}$ aside these ones (for all $j$'s).
Integrating over times naturally yields the measure-valued process
$$A^{(n)}_s(\dd y) \coloneqq  \int_{0\leq r \leq s}  \Updelta^{(n)}(\dd r, \dd y), \qquad s\geq 0.$$

By construction, the pair  $(L^{(n)},A^{(n)})$ is a Markov additive process\footnote{  For the sake of coherence with our earlier notation, the additive component
is the second coordinate of the pair, at the opposite of the convention in the literature on Markov additive processes.}, in the sense that for any $t\geq 0$, 
the conditional distribution of the shifted process $ (L^{(n)}_{t+\cdot},A^{(n)}_{t+\cdot}-A^{(n)}_t)$
given $(L^{(n)}_t,A^{(n)}_t)$ only depend on $L^{(n)}_t$. In other words,  its jump rates only depend on the first component, and more precisely, for every  $(v_0,v)\in\N\times \U^{\downarrow}$, 
 we have for every $t\geq 0$ that 
\begin{align}\label{E:addMarkp}
& \lim_{\varepsilon\to 0+} \varepsilon^{-1} P_\GW^{(n)}\left( 
 L^{(n)}_{t+\varepsilon}=\log(v_0/n),
A^{(n)}_{t+\varepsilon}-A^{(n)}_{t}
= \sum_{j}\delta_{\log(v_j/k)}
\Big | ~L^{(n)}_t=\log(k/n)\right) \nonumber \\
&=k^{\alpha}\boldsymbol{\pi}^{\GW}_k(v_0,v).
\end{align} 

The next result is the core of the proof of Proposition \ref{P:cvdecoreprod}.
Given a characteristic quadruplet $(\sigma^{2}, \mathrm a, \boldsymbol{\Lambda}, \alpha)$ with no killing, let $B$ be
a Brownian motion  and  $\mathbf N$ an independent Poisson random measure on 
 $[0,\infty)\times  \mathcal{S}$ with intensity $ \d  t \boldsymbol{\Lambda}( \d y ,\d \mathbf y )$.
 Recall from  Section \ref{sec:2.2} that $\xi$ then denotes the L\'evy process whose jump process is described by the first projection $N_0$ of $\mathbf N$ on $[0,\infty)\times \R$, and that
$\mathbf N_1$ stands for the second projection  of $\mathbf N$ on $[0,\infty)\times \mathcal S_1$. 
We write $\mathcal{M}(\R)$ for the space of Radon (i.e. locally finite) measures on $\R$ endowed with a distance that metrizes the vague topology (see \cite[Lemma A.5.5]{Kal07})
and introduce the measure valued process $(A_s)_{s\geq 0}$ on $\R$ such that for every 
continuous function with compact support $\varphi: [0,\infty)\to \R$, 
$$\int_{\R} \varphi(y) A_s(\dd y) = \int_{[0,s]\times \mathcal S_1}  \sum_{j}\varphi(y_j) \mathbf N_1(\dd r, \dd \mathbf y).$$

For any $d\geq 1$, we endow the space  
$$\mathbb D^{\dagger}(\R^d)\coloneqq \bigsqcup_{\zeta\geq 0}\mathbb D([0,\zeta],\R^d)$$  of  rcll functions $\omega: [0,\zeta]\to \R^d$ with a finite lifetime $\zeta$, with the Skorokhod's topology; see \cite[Chapter VI]{JS03} for background. Recall that in general,  the convergence  in $\mathbb D^{\dagger}(\R^d)$ for a sequence $(\omega_n)_{n\geq 1}$ is a stronger property than the mere 
convergence in $\mathbb D^{\dagger}(\R)$ of each sequence of coordinates, say $(\omega^i_n)_{n\geq 1}$ for $i=1, \dots, d$.

\begin{lemma} \label{L:cvdecoreprod} Fix an arbitrary integer $d\geq 1$, and for every $j=1, \ldots, d$, let $\varphi_j: [0,\infty)\to \R$ be 
a continuous function with compact support.
Then  the sequence of the $(d+1)$-tuple of processes
$$\left(L^{(n)}_s, \int_{[0,s]}  \varphi_1(y) A_s^{(n)}(\dd y), \ldots, \int_{[0,s]}  \varphi_d(y) A_s^{(n)}(\dd y)
\right)_{s\geq 0}$$
converges in distribution on $\mathbb D^{\dagger}(\R^{d+1})$  as $n\to \infty$ towards
the pair
$$\left(\xi(s), \int_{\R} \varphi_1(y) A_s(\dd y), \ldots, \int_{\R} \varphi_d(y) A_s(\dd y)
\right)_{s\geq 0}.$$
\end{lemma}

\begin{proof}  We aim at applying Theorem VIII.2.17
in Jacod and Shiryaev \cite{JS03}, which provides a general framework for the functional convergence in distribution of a sequence of semi-martingales to a process with independent increments and no fixed discontinuities in terms of the so-called modified characteristics. Readers may wish to have \cite{JS03} at hand for certain definition and results that will be quoted here. 
For the sake of notational simplicity, we assume that $d=1$ and write $\varphi=\varphi_1$; the calculations for general $d\geq 2$ are heavier but similar. 

To start with, the process
$$\left(\xi(s), \int_{[0,s]\times \mathcal S_1}  \sum_{j}\varphi(y_j) \mathbf N_1(\dd r, \dd \mathbf y_1)
\right)_{s\geq 0}$$ is a L\'evy process, and choosing 
for truncation function\footnote{The second coordinate has bounded variation, so there is no need to truncate the second variable. } $(y_0,y)\mapsto y_0\indset{|y_0| \leq 1}$, we get that 
its first characteristic (see \cite[Definition II.2.6]{JS03}) is 
$$t\mapsto (\mathrm a t,0),$$
its modified second characteristic (see \cite[Definition II.2.16, and Proposition II.2.17]{JS03}) 
$$t\mapsto t\begin{pmatrix} \sigma^2  +  \int_{(-1,1)}y_0^2 \Lambda_0(\dd y_0) & 0\\0&0\end{pmatrix},$$
and finally, the compensator $\nu(\dd t, \dd y_0, \dd y)$ of the jump measure is given for a generic measurable function $F:\R_+\times \R\times \R\to \R_+$  by
$$\int_{\R_+\times \R\times \R}F(t,y_0, y)\nu(\dd t, \dd y_0, \dd y)=\int_{\R_+\times \R\times \mathrm S_1}F\left (t,y_0, \sum_{j\geq 1}\varphi(y_j)\right )\dd t \boldsymbol{\Lambda}(\dd y_0, \dd \mathbf y). $$

We next turn our attention  to 
 the processes
$$\left(L^{(n)}_s, \int_{[0,s]}  \varphi(y) A_s^{(n)}(\dd y)
\right)_{s\geq 0}.$$
These 
 are semi-martingales with modified characteristic triplets given as follows.
 Recalling that  $(L^{(n)},A^{(n)})$ is a Markov additive process with  jump rates  \eqref{E:addMarkp}, 
  noting that if $r=\log(k/n)$, then $k=n\e^r$ and $k^{\alpha}=n^{\alpha}\e^{\alpha r}$, and finally 
  using   the notation \eqref{eq:deflambdan} for $\boldsymbol{\Lambda}^{(k)}$,  the first characteristic is
$$t\mapsto \left( \int_0^t  \dd s   \int_{-1}^1 y 
 \Lambda_0^{(n\exp(L^{(n)}_s))}( \mathrm{d}y), 0\right), $$
the second (modified)
$$t\mapsto \begin{pmatrix}  \int_0^t  \dd s  \int_{-1}^1 y^2 
\Lambda_0^{(n\exp(L^{(n)}_s))}( \mathrm{d}y) & 0\\0&0\end{pmatrix},$$
and finally, the compensator $\nu^{(n)}(\dd t, \dd y_0, \dd y_1)$ of the jump measure is given for a generic measurable function $F:\R_+\times \R\times \R\to \R_+$  by
\begin{align*}
&\int_{\R_+\times \R\times \R}F(s,y_0, y)\nu^{(n)}(\dd s, \dd y_0, \dd y)\\
&=\int_{\R_+}\dd s  \int_{\R\times \mathrm S_1}  F\left (s,y_0, \sum_j\varphi(y_j)\right ) 
\boldsymbol{\Lambda}^{(n\exp(L^{(n)}_s))}(\dd y_0, \dd \mathbf y). 
\end{align*}

With these expressions for the modified characteristics at hand, it is easily to verify the conditions of  \cite[Theorem VIII.2.17]{JS03} about the convergence of the modified characteristics follow  from Assumption \ref{A:BK}, which entails the claim. 
\end{proof}

We next strengthen the preceding lemma to convergence of measured-valued processes.

\begin{corollary} \label{C:cvdecoreprod} 
There exists a version of the sequence of Markov additive processes $(L^{(n)},A^{(n)})$ 
which converges almost surely on $\mathbb D(\R\times \mathcal{M}(\R) )$  as $n\to \infty$ towards
$(\xi, A)$.
\end{corollary}

\begin{proof} We first infer from Lemma \ref{L:cvdecoreprod} and \cite[Theorem A.5.7]{Kal07} that the sequence of the laws of the processes $(L^{(n)},A^{(n)})$ is relatively compact in 
the space of probability measures on $\mathbb D(\R\times \mathcal{M}(\R) )$. Next, we deduce from \cite[Theorem 23.16]{Kal07} and again Lemma \ref{L:cvdecoreprod} that the finite-dimensional  distributions (obtained by evaluating these processes at finitely many times)  of any  sub-sequential limit coincide with those of $(\xi,A)$. Therefore, the sequence $(L^{(n)},A^{(n)})$ converges in distribution  on $\mathbb D(\R\times \mathcal{M}(\R) )$ towards $(\xi,A)$  as $n\to \infty$.
An appeal to the Skorokhod's coupling theorem \cite[Theorem 5.31]{Kal07} complete the proof.
\end{proof}

The final ingredient that will be needed for the proof of Proposition \ref{P:cvdecoreprod} is the following technical observation.

\begin{lemma} \label{L:cvintskor} For every $n\geq 1$, let $\omega_n: [0, \zeta_n]\to \R$ be a function  in $\mathbb D^{\dagger}(\R)$ and $v_n(\dd t, \dd y)$ a locally finite integer-valued measure on
$[0,\zeta_n]\times \R$.  For every $\varphi\in \mathcal C_c(\R)$, that is $\varphi: \R\to \R$
is a  continuous function with compact support, define the function $v_n^{\varphi}$ in $ \mathbb D^{\dagger}(\R)$ by
$$v_n^{\varphi} (t)\coloneqq \int_{[0,t]\times \R} \varphi(y) v_n(\dd s, \dd y) \quad\text{for }t\in [0,\zeta_n].$$
Assume  that there are  $\omega: [0, \zeta]\to \R$  in $\mathbb D^{\dagger}(\R)$ and  a locally finite integer-valued measure $v(\dd t, \dd y)$ on
$[0,\zeta]\times \R$, such that for every $\varphi\in \mathcal C_c(\R)$,
the pair
$(\omega_n, v_n^{\varphi})$ converges to  $(\omega,v^{\varphi})$ in $\mathbb D^{\dagger}(\R^2)$ as $n\to \infty$. 

Then for every  $\Phi\in \mathcal C_c(\R_+\times \R\times \R)$, 
we have
$$\lim_{n\to \infty} \int_{[0,\zeta_n]\times \R} \Phi(s,\omega_n(s-), y) v_n(\dd s, \dd y) = \int_{[0,\zeta]\times \R} \Phi(s,\omega(s-), y) v(\dd s, \dd y).$$
\end{lemma}
\begin{proof} Consider for every $n\geq 1$ the integer-valued measure $\bar v_n$  on
$[0,\zeta_n]\times \R \times \R$ defined by 
$$\int_{[0,\zeta_n]\times \R\times \R} \Phi(s,w, y) \bar v_n(\dd s, \dd w, \dd y)= \int_{[0,\zeta_n]\times \R} \Phi(s,\omega_n(s-), y) v_n(\dd s, \dd y),$$
and define $\bar v$ similarly. We have to check that $\bar  v_n$ converges vaguely to $\bar v$ as $n\to \infty$, and
in this direction, it suffices to consider the case where $\Phi(s,w, y)=\gamma(s,w)\varphi(y)$, with $\gamma: \R_+\times \R\to [0,1]$  and $\varphi: \R\to [0,1]$ two continuous functions with compact support.
 Using the notation $\dd v_n^{\varphi}(s)$ for the Stieltjes measure associated to 
the increasing function $s \mapsto v_n^{\varphi}(s)$, we can write
$$\int_{[0,\zeta_n]\times \R\times \R} \Phi(s,w, y) \bar v_n(\dd s, \dd w, \dd y)=\int_{[0,\zeta_n]} \gamma(s,\omega_n(s-)) \dd v_n^{\varphi}(s).$$

We next choose a sequence of increasing bijections 
$\beta_n: [0,\zeta]\to [0,\zeta_n]$ such that as $n\to \infty$, $\beta_n$ converges to the identity function on $[0,\zeta]$ and $\omega_n\circ \beta_n$ to $\omega$, both uniformly on $[0,\zeta]$.
We write $0< s_1< \ldots < s_k\leq \zeta$ for the sequence of the jump times of the step function $v^{\varphi}$, and 
 have plainly
 \begin{align*}& \liminf_{n\to \infty}  \int_{[0,\zeta_n]} \gamma(s,\omega_n(s-)) \dd v_n^{\varphi}(s) \\
 & \geq
 \liminf_{n\to \infty} \sum_{i=1}^k \gamma(\beta_n(s_i), \omega_n(\beta_n(s_i)-))\left( v_n^{\varphi}(\beta_n(s_i))-v_n^{\varphi}(\beta_n(s_i)-)\right)\\
& \geq    \sum_{i=1}^k \gamma(s_i, \omega(s_i-))\left( v^{\varphi}(s_i)-v^{\varphi}(s_i-)\right)\\
& \geq \int_{[0,\zeta]}\gamma(s,\omega(s-)) \dd v^{\varphi}(s)
 \end{align*}
 (the last two inequalities are actually equalities). 
 Replacing $\gamma$ by $1-\gamma$, we conclude that
 $$\lim_{n\to \infty} \int_{[0,\zeta]}\gamma(s,\omega_n(s-)) \dd v_n^{\varphi}(s)=  \int_{[0,\zeta]}\gamma(s,\omega(s-)) \dd v^{\varphi}(s),$$
 and  the proof is complete. 
\end{proof}

We can now establish Proposition \ref{P:cvdecoreprod}.

\begin{proof}[Proof of Proposition \ref{P:cvdecoreprod}]
From Corollary \ref{C:cvdecoreprod}, we can construct for each $n\geq 1$ a step process
$f^{(n)}: [0,z^{(n)}]\to \R_+$ and a point process $\eta^{(n)}$ on $[0,z^{(n)}]\times \R_+$  
such that the pair $(f^{(n)},\eta^{(n)})$ has the  law $P_\GW^{(n)}$,
and as $n\to \infty$ 
the sequence of Markov additive processes $(L^{(n)},A^{(n)})$ derived from $(f^{(n)},\eta^{(n)})$
 converges almost surely $(\xi,A)$ to  on $\mathbb D(\R\times \mathcal{M}(\R) )$.
 Roughly speaking, we now have to check that, in the limit as $n\to \infty$, inverting this transformation to recover $(f^{(n)},\eta^{(n)})$ from 
 $(L^{(n)},A^{(n)})$  amounts to  performing the Lamperti transformation on $(L^{(n)},\Updelta^{(n)})$, where we recall that  the additive component $A^{(n)}$ has been defined as the antiderivative of $\Updelta^{(n)}$.

 Essentially, the inverse transformation has two main components, namely we  apply an exponential map $y\mapsto \e^y$ to $L^{(n)}$ and simultaneously shift each atom of $A^{(n)}$, say at $(s,y)$, at $(s, \e^y \exp(L^{(n)}_{s-}))$, and 
second we replace the total rate of jumps $k^{\alpha}$ when $\exp(L^{(n)})$
 is at $k/n$, by $n^{\alpha}$. 
   More precisely,
  we define first  the time change 
 $\tau^{(n)}: \R_+\to \R_+$  implicitly by
 $$ \int_0^{\tau^{(n)}(t)} \exp(\alpha L^{(n)}_s) \dd s = t \quad\text{if}\quad t\leq  \int_0^{\infty} \exp(\alpha L^{(n)}_s) \dd s,  \quad\text{and} \quad  \tau^{(n)}(t)= \infty \text{ otherwise.}$$
We then define the   
  process $\check{f}^{(n)}$ by 
  $$\check{f}^{(n)}(t) \coloneqq \exp\left( L^{(n)}_{\tau^{(n)}(t)}\right), \qquad 0\leq t \leq \int_0^{\infty} \exp(\alpha L^{(n)}_s) \dd s,$$
and,    using for convenience  a generic continuous function with compact support,  $\gamma: \R_+ \times \R_+\to \R$,
 a  point process $\check{\eta}^{(n)}(\dd t, \dd x)$ by
  $$\int_{\R_+\times \R_+} \gamma(t,x) \check{\eta}^{(n)}(\dd t, \dd x) =
  \int_{\R_+\times \R_+} \gamma \left (\tau^{(n)}(s),\exp\left( L^{(n)}_{\tau^{(n)}( s)-} \right)y\right)\Updelta^{(n)}(\dd s, \dd y).$$
  The law of the  pair $(\check f^{(n)}, \check \eta^{(n)})$ is not precisely that $(f^{(n)},\eta^{(n)})$, but rather that of a minor modification of the latter in which independent exponential  
   waiting times with mean $n^{-\alpha}$ between consecutive steps replace the deterministic waiting times $n^{-\alpha}$. By the law of large numbers,   has no effects on the limit behaviors that we want to verify.

Since $L^{(n)}$ converges to $\xi$ as $n\to \infty$, we naturally expect $\tau^{(n)}$ to converge to 
   the Lamperti time-change $\tau$ defined by \eqref{E:Lampertime}; beware that this has to include the convergence of the explosion times, viz.
   $$\lim_{n\to \infty} \int_0^{\infty} \exp(\alpha L^{(n)}_s) \dd s  = \int_0^{\infty} \exp(\alpha \xi(s)) \dd s.$$
    This is where Assumption \ref{A:typex} 
 is needed, and that the above indeed holds can now be seen from \cite[Theorem 3(i)]{BK14} (alternatively, one could derive also this from Proposition  \ref{L:Aspandii}).
 We can now conclude from Corollary \ref{C:cvdecoreprod},  Proposition  \ref{L:Aspandii} and Lemma \ref{L:cvintskor} that almost surely as $n\to \infty$, 
 $\check f^{(n)}$ converges   in $\mathbb D^{\dagger}$ to the  self-similar Markov process $(X(t))_{0\leq t \leq z}$
 and $\check \eta^{(n)}$    in the sense of vague convergence of Radon measures on $[0,\infty)\times (0,\infty]$ to the reproduction process $\eta$  constructed in  Section \ref{sec:2.2}. 
\end{proof}

Theorem \ref{Th:cvdecoreprod} has been now proved in the case without killing, and we still need to address the situation where $\texttt{k}>0$. In short, we shall reincorporate the killing a posteriori and use Lemma \ref{lem:convreprodecolight} to control its effect. Fix $A>1$ very large and suppose that $A$ is not an atom of $\Lambda_{0}( \mathrm{d}y)$. We shall modify both the discrete and continuous generalized L\'evy measure by sterilizing the particles issued from jumps lower than $-A$, that is we consider the push forwards $ {{ \boldsymbol{\Lambda}}}^{{(n),[A]}}$ and ${ \boldsymbol{\Lambda}}^{[A]}$ of  $ { \boldsymbol{\Lambda}}^{{(n)}}$ and ${ \boldsymbol{\Lambda}}$ under the mapping $$ \mathbf{y} = \big(y_{0},(y_{1}, y_{2}, ... ) \big) \mapsto \big(y_0-\infty  \mathbf{1}_{y_{0} \leq -A} , (y_{1}, y_{2}, ... )\big).$$
We shall write $ P_{\GW}^{(n),[A]}$ for the law of the rescaled discrete decoration-reproduction process associated to $ {{ \boldsymbol{\Lambda}}}^{{(n),[A]}}$ and $ P_{1}^{[A]}$ for the continuous one associated to the characteristics $(  \sigma^{2},\mathrm{a}, { \boldsymbol{ \Lambda}^{[A]}} ; \alpha)$. Since we sterilized some particles, those new generalized L\'evy measure still satisfy Assumptions \ref{A:all:GW}, \ref{A:BK} and  \ref{A:typex} (the drift and variances being untouched) and moreover we have a ``pure'' convergence of the killing term $$ \mathrm{k}^{{[A]}} =  { \boldsymbol{ \Lambda}^{{[A]}}}(\{-\infty\} \times \mathcal{S}_{1}) =  \boldsymbol{ \Lambda}((-\infty,A] \times \mathcal{S}_{1}),$$ that is
\begin{equation} \label{E:killingalakon} \lim_{n\to \infty}   \int_{ \mathcal{S}_{1}} {\boldsymbol{\Lambda}}^{(n),[A]} ( \{-\infty\} \times \mathrm{d}  \mathbf{y}) g( \mathbf{y})  = \int_{ \mathcal{S}_{1}} {\boldsymbol{\Lambda}^{[A]}} ( \{-\infty\} \times \mathrm{d}  \mathbf{y}) {g}(\mathbf{y}), \end{equation}
for any bounded continuous function $ g : \mathcal{S}_{1} \to \mathbb{R}$.
In this case, the killing can be interpreted as an extra decoration and one can go over the proof of Theorem \ref{Th:cvdecoreprod} and get the convergence of $ P_{\GW}^{{(n),[A]}}$ towards $P_{1}^{[A]}$ in the sense of $ \mathbb{D}^{\dagger}$. In particular, since $ \mathrm{k}_{A} >0$ the continuous decoration process dies with a jump, that is under $ P_{1}^{[A]}$ we have $f(z-) > 0$. In fact, the law $P_{1}^{[A]}$ can simply be obtained by killing the decoration-reproduction process under $P_{1}= P_{1}^{[-\infty]}$ at the first time the decoration process makes a negative jump of size at least $ 1-\mathrm{e}^{-A}$ times its current position. In particular, for any $ \varepsilon>0$ one can assume that $A$ is taken large enough so that 
$$ P_{1}^{[A]}( f(z-)\cdot \mathrm{e}^{-A} > \varepsilon) \leq  \varepsilon.$$
On the event above, we have convergence in distribution of the rescaled discrete decoration-reproduction process freezed at the first time it drops below $ \varepsilon$ towards the continuous one (also freezed below level $  \varepsilon$). For any $ \varepsilon>0$, letting  $A \to \infty$ we deduce in the notation of Lemma \ref{lem:convreprodecolight} that $ P_{\GW}^{{(n), \varepsilon}} \to P_{1}^{ \varepsilon}$ and deduce  Theorem \ref{Th:cvdecoreprod} from it.

\section{Lower bounds on superharmonic functions}\label{sec:minimal:reverse}

We have not been able to get a simple explicit condition in terms of the reproduction kernel of the Galton--Watson process only that would ensure the discrete Cramer condition of Assumption~\ref{A:hGWrv}. Nonetheless, in this direction,  we establish here interesting lower bounds for the growth of the superharmonic function $ l^{\varpi}$ that complete the upper bounds of Section \ref{sec:superharmonic}. If we further assume convergence of the discrete cumulant functions, we are able to give an equivalent for $l^{\varpi}_\GW$  when $\varpi(n) = n^{\gamma-\alpha}$; see Proposition \ref{prop:si:tout:va:bien:pour:kappa:n}. Since the results are not needed in the rest of the section, this part may be skipped at first reading.

Again, in this section,  Assumption \ref{A:all:GW} is enforced.  The definition of the generator \eqref{Eq:Laplacediscret}  extends by linearity
to functions $\varphi :  \mathbb{N} \to \mathbb{R}$ such that  $ \boldsymbol{m}_{\GW}|\varphi|$ is everywhere finite; we then call $\varphi$
superharmonic on a subset $B\subset \N$ if $ \boldsymbol{m}_{\GW}\varphi \leq \varphi$ on $B$, without requesting $\varphi\geq 0$.  
 We point at the following version of the minimum principle of superharmonic functions in this setting.
 
  \begin{lemma}[Nonnegativity of superharmonic functions] \label{L:shgeq0}
   Let 
 $\varphi :  \mathbb{N} \to \mathbb{R}$ be a real-valued  function and $\phi :  \mathbb{N} \to \mathbb{R}_+$ a nonnegative function.
 Write $\varphi^-$ for the negative part of $\varphi$ and assume that
   with $\varphi^-(n) = o(\phi(n))$ as $n\to \infty$.
   Suppose   that $\varphi$ and $\phi$ are both superharmonic on some subset $B\subset \N$. 
   
   If furthermore $\varphi\geq 0$ on $B^c$, then $\varphi \geq 0$ everywhere. 
  \end{lemma}
 \begin{proof} We first treat the case $B=\N$ and will then reduce the case $B\subset \N$ to the former.
 
 So suppose that $\varphi$ and $\phi$ are superharmonic, and write 
 $$\varphi(\Z(n))\coloneqq \sum_{j=1}^{\infty}\varphi(Z_j(n)).$$ 
 The process $\left(\varphi(\Z(n))\right)_{n\geq 0}$ is a $\Q_k^{\GW}$-supermartingale for any $k\in \N$, in particular we
 have the lower-bound
 $$\varphi(k)\geq \E_k^{\GW}\left(\varphi(\Z(n))\right)\qquad \text{ for all }n\geq 0.$$
Then take $\varepsilon>0$, and choose $N\geq 1$ sufficiently large so that $\varphi^-(n)\leq \varepsilon \phi(n)$ for all $n> N$.
There is the upper bound bound
$$\varphi^-(\Z(n)) \leq \sum_{j=1}^{N} \varphi^-(j) Z_j(n) + \varepsilon \phi(\Z(n)).$$
Recall from Lemma \ref{eq:bis:finite} that if we write $\bar Z_j\coloneqq \sum_{i=0}^{\infty}Z_j(i)$ for the total number of particles of type $j$ produced by the Galton--Watson process, then $\E_k^{\GW}( \bar Z_j)<\infty$. This entails
$$\lim_{n\to \infty} \E_k^{\GW}\left(\sum_{j=1}^{N} \varphi^-(j) Z_j(n) \right)=0.$$
On the other hand, $\E_k^{\GW}(\phi(\Z(n)))\leq \phi(k)$, since $\phi$ is superharmonic, and we conclude that
$$\limsup_{n\to \infty} \E_k^{\GW}\left(\varphi^-(\Z(n))\right)\leq \varepsilon.$$
Since $\varepsilon$ can be chosen arbitrarily small, we have $$\liminf_{n\to \infty} \E_k^{\GW}\left(\varphi(\Z(n))\right) \geq 0,$$
which proves that indeed $\varphi(k)\geq 0$ for any $k\geq 1$.

Now consider any $B\subset \N$ and assume that $\varphi$ and $\phi$ are superharmonic on $B$ only, and that $\varphi \geq 0$ on $B^c$. 
Imagine that we modify the Galton--Watson process by sterilizing particles with types in $B^c$, leaving unchanged the reproduction of particles with types in $B$.
In other words, we consider the Galton--Watson process $\left(\Z'(n)\right)_{n\geq 0}$ with mean measure $\boldsymbol{m}_{\GW}'$ given by 
$m_{\GW}'(i,j)=m_{\GW}(i,j)$ for $i\in B$ and $m_{\GW}'(i,j)=0$ for $i\in B^c$. Since $\varphi$ and  $\phi$ are both nonnegative on $B^c$, 
they are both superharmonic everywhere for the sterilized process. Plainly, the sterilized Galton--Watson process produces fewer particles than the original one,
so the expectation of the total number of particles of any given type $j$ remains finite, no matter the type of the ancestral particle\footnote{This observation is needed as the accessibility requirement in Assumption \ref{A:all:GW} may fail for the sterilized Galton--Watson process.}. Hence, by the first part of the proof, $\varphi\geq 0$ on $\N$. \end{proof}
 
 We immediately derive the following reverse version of Lemma \ref{lem:minimal}

 \begin{lemma} [Lower bound on superharmonic functions] \label{lem:minimal:reverse}
Let $\varphi: \N\to \R$ and $\varpi: \mathbb{N}\to \mathbb{R}_+$ such that $l_{\GW}^{\varpi}$ is everywhere finite. 
Suppose that  there exist a  subset $B\subset \N$  and $\phi:\mathbb{N}\to \mathbb{R}^*_+$ superharmonic on $B$,  such that $\varphi=o(\phi)$ (hence    $ \boldsymbol{m}_{\GW}(|\varphi|)$ is everywhere finite). If
\begin{equation}
\left\{
\begin{aligned}
    \varphi &\leq  l_{\GW}^{ \varpi} && \text{on }B^c, \\
    -\boldsymbol{G}_{\GW} \varphi &\leq \varpi && \text{on } B,
\end{aligned}
\right.
\label{eq:deltaminus}
\end{equation}
then $\varphi\leq l_{\GW}^{\varpi}$ everywhere.
 \end{lemma}
\begin{proof}
It suffices to apply Lemma \ref{L:shgeq0} to the function $\ell_{\GW}^{\varpi}-\varphi$. 
\end{proof}

We conclude this section by presenting a sharper estimate for the super-harmonic functions $l^{\varpi}_{\GW}$. Although it will not be used in the sequel, it provides additional intuition and helps justify Assumption \ref{A:hGWrv}, making it more natural. Since proving this lemma requires the convergence of discrete cumulants, we begin by showing that under mild conditions, this convergence follows from the convergence of characteristics (Assumption \ref{A:BK}). We use the notation $ \boldsymbol{\Lambda}_1^{(n)}$  for the image of   $ \boldsymbol{\Lambda}^{(n)}$ by the  projection $[-\infty, \infty)\times \mathcal S_1\to S_1$.

\begin{lemma} \label{lem:convkappacond}Suppose that Assumption \ref{A:BK} holds. Then for every $\gamma >0$ we have 
$$ \liminf_{n \to \infty} \kappan(\gamma) \geq \kappa( \gamma) \quad  \in  \mathbb{R} \cup \{+ \infty\}.$$
If furthermore, for some $\gamma,\delta >0$, we have 
$$ \left\{\begin{array}{l} \limsup_{n \to \infty }  \int_{ \mathcal{S}_1} \boldsymbol{\Lambda}^{(n)}_1 ( \mathrm{d} \mathbf{y}) \sum_{i \geq 1} \mathrm{e}^{(\gamma -\delta) y_i} < \infty\\
\limsup_{n \to \infty }  \int_{ \mathcal{S}_1} \boldsymbol{\Lambda}^{(n)}_1 ( \mathrm{d} \mathbf{y}) \left( \sum_{i \geq 1} \mathrm{e}^{\gamma  y_i}\right)^{1+\delta} < \infty \end{array} \right.
 \quad \mbox{ and } \quad \limsup_{n \to \infty }  \int_{(1,\infty)} \Lambda^{(n)}_{0} ( \mathrm{d} y_{0}) \mathrm{e}^{(\gamma +\delta) y_0}< \infty,$$ then we have $\kappan(\gamma) \to \kappa(\gamma) \in \mathbb{R}$ as $n \to \infty$. 
\end{lemma}

\begin{proof}
 First, for every $n\geq 1$, we decompose $\kappan$ in the form $ \kappan(\gamma) = I_0^{(n)}+I_0^{(n),\prime} +I_1^{(n)},$ 
where $$I_0^{(n)}:=\int \Lambda^{(n)}_0 (\dd y_0 ) \big(\exp(\gamma y_0)-1-\big(\gamma y_0+\gamma^{2}y_0^2/2)\mathbf{1}_{ |y_{0}| < 1}\big)\big),$$
$$I_0^{(n),\prime}:= \int_{(-1,1)} \Lambda^{(n)}_0 (\dd y_0 )  \big(\gamma y_0+\gamma^{2}y_0^2/2\big)\quad \text{ and }\quad I_1^{(n)}:=\int_{\mathcal{S}_1}  \boldsymbol{\Lambda}^{(n)}_1 (\dd \mathbf{y} )\sum_{i\geq 1}\exp(\gamma y_i).$$ We are going to study each term separately. By Assumption \ref{A:BK} we already know that 
$$\lim \limits_{n\to \infty}  I_0^{(n),\prime}= a \gamma + \frac{1}{2}\sigma^{2} \gamma^{2}. $$
Let us now study the integral $I_0^{(n)}$. In this direction, note that since
$$\exp(\gamma y_0)-1-\gamma y_0-\gamma^{2}y_0^2/2= o(y_0^{2}), \quad \text{ as } y_0\to 0,$$
 by Assumption \ref{A:BK} it plainly follows that
 \begin{align*}
\lim \limits_{n\to \infty}  \int \Lambda^{(n)}_0 (\dd y_0 ) \big(\exp(\gamma y_0) \mathbf{1}_{y_0<A} -1&-\big(\gamma y_0+\gamma^{2}y_0^2/2)\mathbf{1}_{ |y_{0}| < 1}\big)\big)\\
&=\int \Lambda_0 (\dd y_0 ) \big(\exp(\gamma y_0) \mathbf{1}_{y_0<A} -1-\big(\gamma y_0+\gamma^{2}y_0^2/2)\mathbf{1}_{ |y_{0}| < 1}\big)\big)
\end{align*}
for every $A>1$ with $\Lambda_{0}(\{A\})=0$. Hence,  we infer by monotone convergence taking the limit along $A\to\infty$, with $\Lambda_{0}(\{A\})=0$, that:
$$\liminf\limits_{n\to\infty} I_{0}^{(n)}\geq \int \Lambda_0 (\dd y_0 ) \big(\exp(\gamma y_0)-1-\big(\gamma y_0+\gamma^{2}y_0^2/2)\mathbf{1}_{ |y_{0}| < 1}\big)\big)=\psi(\gamma)-a\gamma -\frac{1}{2}\sigma^{2}\gamma^2, $$
Furthermore, if the assumption of the second statement of the lemma holds, then  for every $A>1$, we also have:
$$\limsup_{n\to \infty}\int_{[A,\infty)} \Lambda_{0}(\dd y_0) \exp(\gamma y_0) \leq \exp(-\delta A) \limsup_{n\to \infty}\int_{[1,\infty)} \Lambda_{0}(\dd y_0) \exp((\gamma+\delta) y_0)=0 ,$$
and so in this case taking the limit $A\to\infty$, we get $\lim_{n\to \infty} I^{(0)}_n= \psi(\gamma)-a\gamma -\frac{1}{2}\sigma^{2}\gamma^2$. It remains to study the integral $I_1^{(n)}$. In this direction note that, for $A,B\geq 1$ with $\boldsymbol{\Lambda}_1(\{\exists i\geq 1:~|y_i|=A\})=0$, it follows from  Assumption \ref{A:BK} that:
\begin{equation}\label{eq:Lambda:n:B:A:1}
\lim \limits_{n\to \infty}   \int_{\mathcal{S}_1} \boldsymbol{\Lambda}^{(n)}_1 (\dd  \mathbf{y}) ~ \sum_{i=1}^{B}\exp(\gamma y_i) \mathbf{1}_{y_i\in[-A, A]} = \int \boldsymbol{\Lambda}_1 (\dd  \mathbf{y}) ~ \sum_{i=1}^{B}\exp(\gamma y_i)\mathbf{1}_{y_i\in[-A, A]}.
\end{equation}
It follows again by monotone convergence that:
\begin{equation}\label{eq:Lambda:n:B:A:2}
\liminf \limits_{n\to \infty} \int_{\mathcal{S}_1} \boldsymbol{\Lambda}_1^{(n)}(\dd  \mathbf{y}) \sum_{i=1}^{\infty}\exp(\gamma y_i) \geq  \int_{\mathcal{S}_1} \boldsymbol{\Lambda}_1(\dd  \mathbf{y}) \sum_{i=1}^{\infty}\exp(\gamma y_i).  
\end{equation}
In particular, putting all together, under Assumption \eqref{A:BK} we always have $\liminf_{n\to \infty}\kappan(\gamma)\geq \kappa(\gamma)$. For the rest of the proof we assume that the assumptions of the second part of the lemma hold,  and remark that to conclude we need to show that $\lim\limits_{n\to \infty} I_1^{(n)}= \int \mathbf{\Lambda}_1(\dd  \mathbf{y}) \sum \limits_{i\geq 1}\exp(\gamma y_i).$ To this end, by \eqref{eq:Lambda:n:B:A:1}, it suffices to show:
$$\limsup_{A\to\infty}\limsup_{B\to\infty}\limsup_{n\to \infty}\int_{\mathcal{S}_1} \boldsymbol{\Lambda}^{(n)}_1 (\dd  \mathbf{y}) ~\Big( \sum_{i=1}^{\infty}\exp(\gamma y_i)- \sum_{i=1}^{B}\exp(\gamma y_i) \mathbf{1}_{y_i\in[-A, A]} \Big)=0. $$
In words, we  aim to approximate the $ \boldsymbol{ \Lambda}^{(n)}_1$-integral of $\sum_{i} \mathrm{e}^{\gamma y_{i}}$ by that of $ \sum_{i \leq B} \mathrm{e}^{\gamma y_{i}} \mathbf{1}_{ y_{i} \in[-A,A]}$ uniformly in $n$. In this direction, notice first that, for every $(y_i)_{i\geq 1}\in \mathcal{S}_1$, we have the inequalities 
 \begin{eqnarray*} \mathrm{e}^{-\delta A} \sum_{i \geq 1} \mathrm{e}^{(\gamma + \delta) y_{i}} &\geq&  \sum_{i \geq 1} \mathrm{e}^{\gamma  y_{i}} \mathbf{1}_{ y_{i}> A}\\
 \mathrm{e}^{-\delta A} \sum_{i \geq 1} \mathrm{e}^{(\gamma - \delta) y_{i}} &\geq&  \sum_{i \geq 1} \mathrm{e}^{\gamma  y_{i}} \mathbf{1}_{ y_{i}<- A}.  \end{eqnarray*}
Next, for $\mathbf{y}=(y_i)_{i\geq 1}\in \mathcal{S}_1$, set $ X_{A,B}(\mathbf{y}) :=\sum_{i > B} \mathrm{e}^{\gamma  y_{i}} \mathbf{1}_{ y_{i} \in[-A, A]}$ and remark that since the $y_{i}$ are non-increasing, when $X_{A,B}(\mathbf{y})>0$ we must have $\sum_{i=1}^{B} \mathrm{e}^{\gamma y_{i}} \geq B \mathrm{e}^{-\gamma A}$ and we infer that
 $$ \left( \sum_{i \geq1} \mathrm{e}^{\gamma  y_{i}} \right)^{1+\delta} \geq \mathbf{1}_{ X_{A,B}(\mathbf{y}) >0} \big(X_{A,B}(\mathbf{y}) + B \mathrm{e}^{-\gamma A} \big)^{1+ \delta} \geq X_{A,B}(\mathbf{y}) \left(B\mathrm{e}^{-\gamma A}\right)^{\delta}$$ as a consequence of the inequality $(x+ y)^{1+\delta}\geq x y^{\delta}$ for $x,y\geq 0$. Combining these observations, we deduce that 
  \begin{eqnarray*} && \int \boldsymbol{ \Lambda}^{(n)} ( \mathrm{d} \mathbf{y}) \left( \sum_{i \geq 1} \mathrm{e}^{\gamma y_{i}} -\sum_{i=1}^{B} \mathrm{e}^{\gamma y_{i}} \mathbf{1}_{ y_{i} \in [-A,A]}\right) \\
   & \leq & \mathrm{e}^{-\delta A} \int \boldsymbol{ \Lambda}^{(n)} ( \mathrm{d} \mathbf{y}) \sum_{i \geq 1} \left(\mathrm{e}^{(\gamma + \delta) y_{i}} + \mathrm{e}^{(\gamma - \delta) y_{i}}\right) +  \frac{ \mathrm{e}^{\delta \gamma A}}{  B^{\delta} }  \int \boldsymbol{ \Lambda}^{(n)} ( \mathrm{d} \mathbf{y}) \left(\sum_{i \geq 1} \mathrm{e}^{\gamma  y_{i}}\right)^{1+\delta},  \end{eqnarray*} and thanks to our assumptions, 
  \eqref{eq:Lambda:n:B:A:2} follows. This completes the proof of the lemma.
\end{proof}

We start with a lower bound on the growth of the smallest super-harmonic functions considered in Lemma \ref{eq:bis:finite}. \begin{lemma}  \label{eq:bis:finite:converse} Suppose that $0<\gamma_{-}<\gamma^{+}$ are such that $$ \limsup \limits_{n\to \infty}\kappan(\gamma_{-})>0  \quad \mbox{ and } \quad \limsup \limits_{n\to \infty}\kappan(\gamma_{+})<0.$$ Then 
for every $\varpi:\mathbb{N}\to \mathbb{R}_+$ with finite support, the function $l_{\GW}^{\varpi}$ satisfies 
 $$l_{\GW}^{\varpi}(n)=O(n^{\gamma_{+}}) \quad \mbox{ and } n^{\gamma_{-}} = O ( l_{\GW}^{\varpi}(n)).$$
\end{lemma}
In particular, under the Assumption \ref{A:BK}, the first condition is granted as soon as $\kappa(\gamma_{-}) \in (0, \infty]$ by the first point of Lemma \ref{lem:convkappacond}. In particular, under the further Assumption \ref{A:omega-} and of the convergence of the discrete cumulants towards $\kappa$ in the right vicinity of $\omega_{-}$ we have 
$$ l_{\GW}^{\varpi}(n) = n^{\omega_{-} + o(1)}.$$

\begin{proof}
The upper bound $l_{\GW}^{\varpi}(n)=O(n^{\gamma_{+}})$ is established in Lemma~\ref{eq:bis:finite}. For the lower bound, consider the function $\varphi(n)=n^{\gamma_{-}}$. Since 
$$\limsup_{n\to\infty}\kappan(\gamma_{-})>0,$$
we can find $M\geq 1$ such that $\boldsymbol{G}_{\GW}\varphi$ is positive on $B:=\mathbb{N}\setminus \{1,\dots, M-1\}$. Moreover,  without loss of generality, we may assume that $\varpi$  is supported on $B^{c}$ and in particular $-\boldsymbol{G}_{\GW}\varphi\leq \varpi$ on $B$. Finally, since $B^c$ is finite,  there exists a constant $C>0$ such that $l_{\GW}^{\varpi}(n)\ge \frac{1}{C}\,\varphi(n)$ on $B^c$, and  an application of Lemma \ref{lem:minimal:reverse} entails $n^{\gamma_{-}}\leq C\cdot l_{\GW}^{\varpi}(n) $, for every $n\geq 1$. This completes the proof of the lemma.
\end{proof}

If we further assume convergence of the cumulant function  we even have a sharper result:

\begin{proposition}\label{prop:si:tout:va:bien:pour:kappa:n} 
Suppose that there exist $0<\gamma_1<\gamma_0$ with
$$\lim_{n\to \infty}\kappan(\gamma)=\kappa(\gamma)<0\qquad \text{ for every }\gamma\in [\gamma_1, \gamma_0].$$ Then fix some $\gamma \in (\gamma_1, \gamma_0)$
and  take $\varpi(n)=n^{\gamma-\alpha}$. We have
$$
\lim \limits_{n\to \infty}\frac{l^{ \varpi}_\GW(n)}{n^{\gamma}}= -\frac{1}{\kappa(\gamma)}.
$$
\end{proposition}
\begin{proof} 
Fix $\gamma\in (\gamma_1,\gamma_0)$. Next, set $\varpi(n)\coloneqq n^{\gamma-\alpha}$ and $\psi(n)\coloneqq n^{\gamma}$, for $n\geq 1$. By Proposition~\ref{prop:APHI}, the function $l_{\GW}^{\varpi}$ takes finite values; recall also from \eqref{eq:ellsuperh} that $-\boldsymbol{G}_{\GW} l_{\GW}^{\varpi}=\varpi$. Since by assumption $\lim_{n\to \infty}\kappan(\gamma)=\kappa(\gamma)<0$, 
for any positive constants $c,C$ with  $c<-1/\kappa(\gamma)<C$, there is some sufficiently large $N\geq 1$ such that 
$$-c\cdot  \boldsymbol{G}_{\GW} \psi(n)<\varpi(n)\leq  -C \cdot \boldsymbol{G}_{\GW} \psi(n),\quad \text{ for } n>N.$$
Introduce  $\varpi'\coloneqq \mathbf{1}_{\{1,\dots, N\}}$; by Assumption \ref{A:all:GW} the function $ l_{\GW}^{\varpi'}$ is strictly positive everywhere and  recall from Lemma \ref{eq:bis:finite}  that $l_{\GW}^{\varpi'}= \mathcal{O}(n^{\gamma_1})= o(n^{\gamma})$.
We can then find $c', C'>0$ such that 
$$c  \psi(n)- c'  l_{\GW}^{\varpi'}(n)<l_\GW^\varpi(n)\leq  C  \psi(n)+ C'  l_{\GW}^{\varpi'}(n),\quad \text{ for } n\leq N.$$
We now write $\varpi''=\varpi+c^\prime\cdot \varpi'$, so that both $-c\cdot \boldsymbol{G}_{\GW} \psi(n)<\varpi''(n)$ for $n>N$ and $c\psi(n)< l^{\varpi''}_{\GW}(n)$ for $n\leq N$. 
We can readily check from Proposition  \ref{prop:massepresqueomega-} that  all the assumptions of Lemma~\ref{lem:minimal:reverse} are fulfilled. We conclude that 
$l^{\varpi}_{\GW}+l^{\varpi'}_{\GW} = l^{\varpi''}_{\GW}\geq c\cdot  \psi$ everywhere and therefore 
$$\liminf_{n\to \infty} n^{-\gamma} l^{\varpi}_{\GW}(n) \geq c.$$
Similarly, we write  $\phi=C\psi +C' l_{\GW}^{\varpi'}$ and use Lemma \ref{lem:minimal}  to check that 
$$\limsup_{n\to \infty} n^{-\gamma} l^{\varpi}_{\GW}(n) \leq C.$$
This completes the proof of the proposition since $c$ and $C$ can be taken arbitrary close to $-1/\kappa(\gamma)$.
\end{proof}

\section*{Comments and bibliographical notes} As we already said, Theorem \ref{Th:cvdecoreprod} can be seen as an extension of a result of Bertoin and Kortchemski \cite{BK14}  on scaling limits for  Markov chain with  values in $\N$, which in turn generalized a work by Haas and Miermont \cite{haas2011self} in the case of decreasing Markov chains. Notice however that our proof is self-contained and some key techniques are different, e.g.~we use here the convergence of modified characteristics, whereas \cite{BK14} relies on the convergence of infinitesimal generators. Although we believe it is only technical, our strong assumption \eqref{E:killingalakon} on the killing rate seems annoying to remove and would require to sharpen our understanding of scaling limits of Markov chains back to its roots.  The only assumption which we did not manage to transform into an analytic tractable one is Assumption \ref{A:hGWrv}. In particular, we did not find easy criteria to ensure that for finitely supported $\varpi$, the function $l^{ \varpi}_\GW$  is regularly varying with exponent $\omega_-$: in standard situations we only managed to identify its growth exponent thanks to Lemma \ref{eq:bis:finite} and its counterpart Lemma \ref{eq:bis:finite:converse}. However, we will see that in most examples, this assumption will be delivered by analytic combinatorics techniques. Finally, we note than the analytic assumptions are usually stronger that the probabilistic assumption of the previous chapter, see the remark after Theorem VIII.2.17 in Jacod and Shiryaev \cite{JS03}.

\chapter{Applications} \label{chap:applications}
We illustrate here the results of Chapter \ref{chap:6} by discussing families of integer-type Galton--Watson trees that, after appropriate rescaling, converge to some of the self-similar Markov branching trees described in Chapter \ref{chap:example}. In particular, we recover the results of Haas and Miermont \cite{HM12} in the setting of fragmentation trees. Furthermore, we showcase the full scope of our invariance principles by considering Galton--Watson trees that arise in the Markovian exploration (also known as peeling explorations) of random planar maps and in critical parking processes on trees. Finally, we briefly outline additional potential applications.

\section{Non-increasing integer type Galton--Watson processes}

We shall first focus on the \textbf{non-increasing case}. This refers  to the discrete version of that considered in Section \ref{sec:5.2}, namely  where particles of types $ n \geq 1$ only give rise to particles of types at most $n$. Equivalently, the measures $\boldsymbol{\Lambda}^{(n)}$ defined in \eqref{eq:deflambdan} are  supported by 
$$\big\{ (y_0, ( y_1, \dots)) : y_i \leq 0\text{ for all }i \geq 0\big\}.$$ In this situation,  the mean-reproduction matrix \eqref{eq:meanrepoductionmatrix} is triangular, and   it is easy seen that the mean  total size  $ \mathbb{E}^{\GW}_{n}(\# T_{\GW})$ is finite for all $n \geq 1$ if and only if 
 \begin{equation} \label{eq:triangularcritical} \mathbb{E}^\GW_n\big(Z_n(1)\big) < 1 \quad\text{and} \quad \mathbb{E}^\GW_n\big(Z_k(1)\big) < \infty, \quad \text{ for all } n>k \geq 1. \end{equation}

 \subsection{Recovering Haas--Miermont}
 Let us start with the framework of Haas--Miermont \cite{HM12}, which, as already mentioned in the introduction, was our main source of inspiration. This framework consists of the non-increasing and  \textbf{conservative} case where we additionally have 
 $$\boldsymbol{\Lambda}^{(n)} \left( \left\{(y_{0},(y_{1}, \dots )) \in \mathcal S:  \sum_{j=0}^{\infty} \e^{y_j}  \ne 1 \right\} \right) =0.$$

More precisely, in \cite{HM12}, the authors consider discrete random trees constructed recursively from a \textbf{splitting algorithm}. The latter is induced by a splitting rule  $(q_{n})_{n\geq 1}$, where each $q_{n}$ is the law of a random partition of $n$, recall that a partition of  $n\geq 2$ 
is a non-increasing  sequence of positive integers with sum $n$ and we make the convention that the only partitions of $1$ are $\{1\}$ and the empty partition $\varnothing$. Following \cite{HM12}, one  constructs a random rooted tree with $n$ leaves according to the following procedure. Start from a collection of $n$ indistinguishable balls, and with probability $q_n(\{n_1,\dots,n_p\})$, split the collection into $p$ sub-collections with $n_1,\dots,n_p$ balls each.  We then iterate the splitting operation independently for each sub-collection using the probability distributions $q_{n_j}$ for $j=1, \ldots, p$. The process stops when a single ball is erased with probability $q_1(\{\varnothing\})$.

\begin{figure}[!h]
 \begin{center}
 \includegraphics[width=15cm]{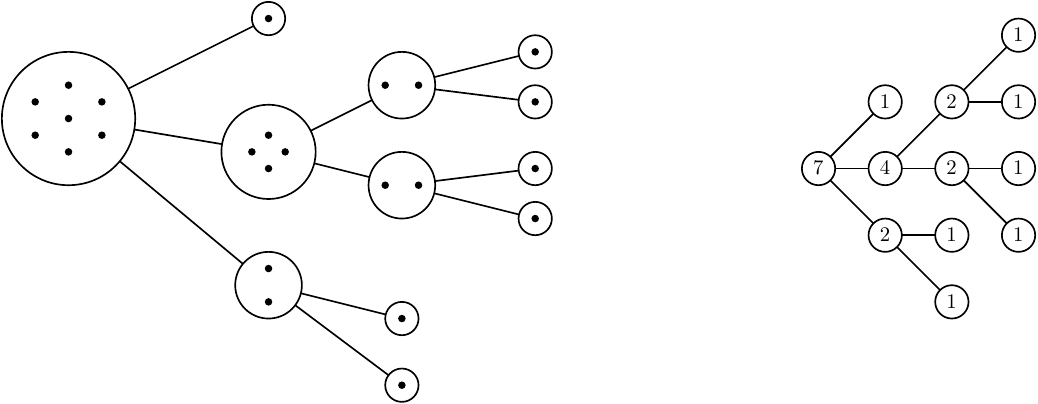}
 \caption{The multi-type Galton--Watson process induced by a split-rule $( q_n)_{n \geq 1}$.}
 \end{center}
 \end{figure}

If each vertex is labeled by its number of balls, this obviously yields a non-increasing and conservative multi-type Galton--Watson process, and the law of the underlying decorated tree $\texttt{T}_{\GW}=(T_\GW,  g_{\GW})$ is denoted by $ \Q_n^{\GW,  {\normalfont\text{\tiny{HM}}}}$. We shall always suppose that 
  \begin{eqnarray} \label{eq:stricttrivial} \forall n \geq 1, \qquad q_n( \{n\}) < 1,   \end{eqnarray} so that \eqref{eq:triangularcritical} is granted and in particular the underlying  trees $T_\GW$  are a.s.\ finite. For simplicity, we also assume  that $q_{1}(\varnothing) =1$, that is, vertices of label $1$ are leaves.  The discrete generalized L\'evy measures $ \boldsymbol{\Lambda}^{(n)}$ (see \eqref{eq:deflambdan}) are defined using the \textbf{locally largest selection rule}, i.e.~the legitimate child of a particle is the one with the largest label (breaking ties arbitrarily). We shall endow our decorated tree either with the measure 
  $ \uplambda_\GW^{ \mathbf{1}_{\{1\}}}$ which is the counting measure on the leaves  of $T_\GW$. In particular, the $\uplambda_\GW^{ \mathbf{1}_{\{1\}}}$-mass is constant equal to $n$ under $ \Q^{\GW,  {\normalfont\text{\tiny{HM}}}}_{n}$. We can now present the main result of \cite{HM12} as a corollary of our invariance principles:

  \begin{theorem}[{\cite[Theorem 5]{HM12}}]  \label{thm:haasmiermont}Suppose that the splitting rules $(q_n)_{n \geq 1}$ satisfy  \eqref{eq:stricttrivial} and that there exists a non-increasing conservative characteristic quadruplet   $(\sigma^{2}, \mathrm a,\boldsymbol{\Lambda}; \alpha)$, see Section \ref{sec:5.2}, with $\boldsymbol{\Lambda} \ne \delta_{(0,(-\infty, -\infty, \cdots))}$  such that  
  \begin{equation} \label{eq:cvdiscretefrag}  \int_{ \mathcal{S}} \boldsymbol{ \Lambda}^{(n)}( \mathrm{d} \mathbf{y})(1-  \mathrm{e}^{y_{0}}) \ f( \mathrm{e}^{{y_{0}}}, \mathrm{e}^{{y_{1}}}, \dots) \xrightarrow[n\to\infty]{} \int_{ \mathcal{S}}   \boldsymbol{ \Lambda}( \mathrm{d} \mathbf{y})(1-  \mathrm{e}^{y_{0}}) \ f( \mathrm{e}^{{y_{0}}}, \mathrm{e}^{{y_{1}}}, \dots),\end{equation}
   for any bounded continuous function $f:  \mathcal{S}\to \R$ (recall from \eqref{eq:deflambdan} that the renormalization is already included). In particular, \eqref{A:omega-} holds with $\omega_-=1$ so that $(T, g, \upmu)$ is well-defined under $\mathbb{P}_{1}$ and recall that the law of its equivalence class in $ \mathbb{T}_{m}$ is denoted by $ \mathbb{Q}_{1}$. Then the sequence of  the distributions  of the rescaled measured decorated trees 
 $$\left( T_\GW,\frac{d_{T_\GW}}{n^{\alpha}},\rho_{\GW}, \frac{g_\GW}{n}, \frac{ \uplambda^{ \mathbf{1}_{\{1\}}}_\GW}{n}\right) \mbox{ under } \Q^{\GW,  {\normalfont\text{\tiny{HM}}}}_{n},$$  
 converges as $n\to \infty$,  to the law  of 
the measured self-similar Markov tree $(\normalfont{\texttt{T}}, \upmu)$ under $\Q_1$.   \end{theorem}
Theorem 5 of  \cite{HM12} is actually slightly more general\footnote{The convergence in \cite{HM12} does not take into account the usc-decoration and takes place in the Gromov-Hausdorff-Prokhorov sense. However, in the conservative case, the difference is only cosmetic since the decoration $g$ is recovered from the $\upmu$-mass of fringe trees} than the result above since it allows for the presence of slowly varying function which we excluded in these pages for technical simplicity. This theorem had in turn  many applications,  to monotype Galton--Watson trees conditioned either by the number of leaves or by the number of vertices \cite{HM12},  to the cut-tree of large Galton--Watson trees  \cite{BM13}, or to growing $k$-ary trees in the spirit of R\'emy's algorithm \cite{haas2015scaling}.

\begin{proof} As we said above, we use the locally largest selection rule in order to define the various objects used in Chapters \ref{chap:6} and \ref{C:scaling}. Our goal is to apply Theorem \ref{T:mainunconds-mass}, and let us check the various requirements. Recall from Chapter \ref{sec:5.2}, that in the non-increasing conservative case, we have $\omega_-=1$ and $\kappa(\gamma) <1$ for all $\gamma >1$ as soon as the generalized L\'evy measure is non-trivial  $ \boldsymbol{\Lambda} \ne \delta_{(0, (-\infty, -\infty, \dots))}$. Assumption \ref{A:omega-} is plain for non-increasing conservative quadruplets and Assumption \ref{A:hGWrv} is equally trivial in the discrete conservative case since the $\uplambda_\GW^{ \mathbf{1}_{\{1\}}}$-mass is constant equal to $n$ under $ \Q^{\GW,  {\normalfont\text{\tiny{HM}}}}_{n}$. 

We shall now check Assumptions \ref{assum:main:gen:deco:repro} and \ref{A:all:existance:PHI} using the stronger  analytical criteria developed in Chapter \ref{C:scaling}.  In this direction, we remark that Assumption \ref{A:all:GW} (accessibility and sub-criticality of type $1$) is ensured by \eqref{eq:stricttrivial} and conservativeness. Let us next  show that Assumption \ref{A:BK} is also satisfied.   Indeed,  the first condition in Assumption \ref{A:BK} is given by  \eqref{eq:cvdiscretefrag} when one further requests the function $f$ to vanish on some neighborhood of~$\mathbf{0}$. 

The convergence \eqref{eq:cvdiscretefrag} is actually stronger, and since $ y \sim 1- \mathrm{e}^y$ in the neighborhood of $0$ it implies that 
$$  \int_{-1}^1 \Lambda_0^{(n)}( \mathrm{d}y) y \underset{ \mathrm{conservative}}{=} \int_{-1}^0 \Lambda_0^{(n)}( \mathrm{d}y) y   \xrightarrow[n\to\infty]{ \eqref{eq:cvdiscretefrag}}  \int_{-1}^0 \Lambda_0( \mathrm{d}y) y,$$ as well as 
$$  \int_{-1}^1 \Lambda_0^{(n)}( \mathrm{d}y) y^2 \underset{ \mathrm{conservative}}{=} \int_{-1}^0 \Lambda_0^{(n)}( \mathrm{d}y) y^2   \xrightarrow[n\to\infty]{ \eqref{eq:cvdiscretefrag}}  \int_{-1}^0 \Lambda_0( \mathrm{d}y) y^2.$$
This entails the convergence of the drift and variance term with $ \mathrm{a_{can}}=0$ and $\sigma^2=0$ and completes the verification of Assumption \ref{A:BK}. Finally, in the conservative case, the hypotheses of Lemma \ref{prop:si:tout:va:bien:pour:kappa:n} are easily checked for all $\gamma >1$, so that we have $\kappa^{(n)}(\gamma) \xrightarrow[n\to\infty]{} \kappa(\gamma)<1$ and in particular Assumption \ref{A:typex} holds. By Theorem \ref{Th:cvdecoreprod} we deduce that Assumption \ref{assum:main:gen:deco:repro} holds, whereas Lemma \ref{eq:bis:finite} grants Assumption \ref{A:all:existance:PHI}. We can therefore apply Theorem \ref{T:mainunconds-mass} which entails the desired the result. \end{proof}

\begin{remark} Several variants of this theorem can be obtained. For example, for every $\varpi:\mathbb{N}\to \mathbb{N}$ regularly varying with index $\gamma-\alpha$ for some $\gamma>1$,  Theorem \ref{T:mainunconds-length} instead of  \ref{T:mainunconds-mass} enables us to ensure that the sequence of the distributions  of the rescaled measured decorated trees 
 $$ \left( T_\GW,\frac{d_{T_\GW}}{n^{\alpha}}, \rho_\GW \frac{g_\GW}{n}, \frac{\uplambda^{\varpi}_{\GW}}{n^\alpha \varpi(n)}\right) \quad \mbox{under $\Q_n^{\GW,  {\normalfont\text{\tiny{HM}}}}$}$$
 converges as $n\to \infty$,  to the law  of 
the measured self-similar Markov tree $(\normalfont{\texttt{T}}, \uplambda^\gamma)$ under $\Q_1$. A transition appears at $\alpha=1$ and the length measure  $\uplambda_\GW^{\mathbf{1}}$ (do not confuse with $\uplambda^{ \mathbf{1}_{\{1\}}}_\GW$  used in the theorem above) converges either towards the harmonic measure if $\alpha <1$ or towards a length measure if $\alpha >1$, see \cite[Theorem 6]{HM12}.
\end{remark}

 \subsection{Generic-critical Galton--Watson trees conditioned by the height}\label{subsec:appli:GWTh}
 
In relation with Example \ref{ex:brownianheight}, we now treat the case of Galton--Watson trees conditioned by the height which is not covered by the Haas--Miermont framework since the splitting rule is not conservative anymore.

A famous result of Aldous \cite{Ald93} shows that after some proper rescaling,  a critical Galton--Watson tree with finite variance, conditioned in some sense to be large, converges in distribution to a Brownian CRT. When the conditioning is performed on the number of leaves or vertices, this can be shown using Theorem \ref{thm:haasmiermont}, see \cite{HM12} and Example \ref{ex:brownian}. We shall present below a precise argument tailored to our framework, when the conditioning  is performed on the height. 
Roughly speaking, we shall argue that  the Galton--Watson tree conditioned on having height $n$ and rescaled by a factor $1/n$, converges in distribution as $n\to \infty$ to the non-increasing self-similar Markov tree  of Example~\ref{ex:brownianheight}.

We consider a standard mono-type Galton--Watson process started from a single particle. We write $\mathbf{p}=(p_j)_{j\geq 0}$ for the reproduction law, where $p_j$ is the probability that a  particle has an offspring of size $j$; we  implicitly discard the degenerate case when $p_1=1$. We assume that the reproduction law is critical.

\begin{figure}[!h]
 \begin{center}
 \includegraphics[width=10cm]{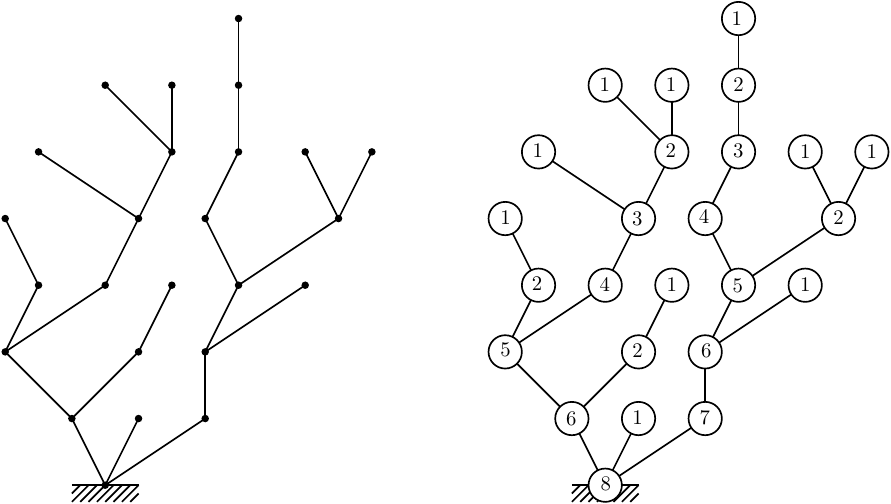}
 \caption{Assigning to each particle the type given by the height of its fringe subtree plus $1$ turns the monotype Galton--Watson conditioned on having height $h$ into a multi-type Galton--Watson under the law $ \Q_{h+1}^{\GW, \mathrm{height}}$. \label{fig:condheight}}
 \end{center}
 \end{figure}

We assign to each particle the type given by $1$ plus the height of the fringe-tree it generates (adding $1$ ensures that the type of a particle is always positive even when it has no offspring, in order to fit the framework of Chapter \ref{chap:6}). See Figure \ref{fig:condheight}. So the types of the progeny of a particle with type $n\geq 1$ are all at most $n-1$, and there is at least one particle in the progeny with type $n-1$ (recall that by convention the type $0$ is attributed to fictitious particles). {When conditioned on having height $h$, the mono-type tree can thus be seen as a labeled tree where the type of the root is equal to $h+1$.  It is readily seen from the branching property of the mono-type Galton--Watson process that, after forgetting the plane structure, this actually yields a non-increasing  integer-type Galton--Watson process}, and we use the notation $\Q^{\GW, \mathrm{height}}_n$ for the law of the underlying decorated tree $\texttt{T}_{\GW}=(T_\GW,  g_{\GW})$ when the initial particle has type $n$. Specifically, if we denote by $( \varphi(h))_{h\geq 1}$ the law of the height of the mono-type Galton--Watson tree plus one, then the probability that  the ancestor of a $ \mathbf{p}$-monotype Galton--Watson plane tree gives rise to $k \geq 0$ subtrees of height $h_{1}+1, \dots , h_{k}+1$ (from left to right) is plainly 
$$ p_{k} \prod_{i=1}^{k} \varphi(h_{i}).$$
We deduce that the law $\Q^{\GW, \mathrm{height}}_n$ of the multi-type Galton--Watson tree is  characterized by 
 \begin{equation} \label{eq:lawGWheight} \mathbb{Q}^{\GW, \mathrm{height}}_{n}\big(  \mathbf{Z}(1) = (i_{1},  \dots , i_{n-1}, 0,\dots)\big) = p_{i_{1} + \dots + i_{n-1}} \frac{\mathbf{1}_{ i_{n-1}\geq 1}}{\varphi(n)} \frac{(i_{1} + \dots + i_{n-1})!}{i_{1}! i_{2}! \dots i_{n-1}!}  \prod_{\ell=1}^{n-1} \big(\varphi(\ell)\big)^{i_{m}},  \end{equation}
 
Discarding the measure for simplicity we have:

\begin{proposition}[{Implied by \cite[Theorem 7.1]{LG10}}]  \label{prop:GWheight}The sequence of the distributions  of the rescaled measured decorated trees 
$$\left( T_\GW,\frac{d_{T_\GW}}{n}, \rho_\GW, \frac{g_\GW}{n}\right) \mbox{ under } \mathbb{Q}^{\GW, \mathrm{height}}_{n}$$
 converges as $n\to \infty$,  to the law  of 
the self-similar Markov tree $\normalfont{\texttt{T}}$  under $\Q_1$ with characteristic quadruplet $(0, -1, \boldsymbol{\Lambda}_{\mathrm{Height}}; 1)$ of Example \ref{ex:brownianheight}.
\end{proposition}
\begin{remark} \label{rek:finitevariancelegall}
The result above, in fact, follows from \cite[Theorem 7.1]{LG10}, which proves convergence using contour functions and also implies the convergence of the renormalized uniform measure towards  $\upmu$. This convergence can also be shown by checking Assumption \ref{A:hGWrv}, but we did not do so for simplicity. See \cite{kesten1996local} for a local limit theorem for the mass and height of generic critical Galton–Watson trees under a third moment condition. The latter can also be treated by considering the generating function of the total number of vertices and the height of the unconditioned Galton–Watson tree. A similar study when  $\mathbf{p}=(p_j)_{j\geq 0}$ is  is in the domain of attraction of a stable law, should be possible, and the arising  ssMt should correspond to those in Example \ref{ex:sampling}. \end{remark}

\begin{proof} The result will follow from an application of Theorem \ref{T:mainuncondw} for which we need to verify Assumptions  \ref{assum:main:gen:deco:repro} and  \ref{A:all:existance:PHI}. Through  Theorem \ref{Th:cvdecoreprod} and Lemma \ref{eq:bis:finite} it suffices to check Assumptions \ref{A:all:GW}, \ref{A:BK} and  \ref{A:typex}. Let us proceed. All that we will use below is the well-known fact that
 \begin{equation}\label{E:kestenneys}
 \varphi(n) \sim \frac{2}{ \mathrm{Var}(\mathbf{p}) n^2}\qquad\text{as }n\to \infty; 
 \end{equation}
where $\mathrm{Var}(\mathbf{p})=\sum_{j\geq 0} (j^2-j)p_j< \infty$, see for instance  \cite[Corollary 1]{kesten1966galton} or \cite[Corollary 1, p. 23]{AN72}. It is also convenient to introduce the notation  
 $$\Phi(n)=\sum_{1 \leq k \leq n}\varphi(k), \qquad n\geq 1,$$ 
 for the distribution function of the height (plus one) of the mono-type Galton--Watson tree; note that  \eqref{E:kestenneys} entails
\begin{equation}\label{E:kestenneys:2}1-\Phi(n)\sim \frac{2}{ \mathrm{Var}(\mathbf{p}) n}.
\end{equation}

\noindent $\bullet$ \textsc{Assumption \ref{A:all:GW}:}  Since the type $1$ is clearly admissible for every $n\geq 1$, it suffices to verify \eqref{eq:triangularcritical}. To this end, remark that  $\Q^{\GW, \mathrm{height}}_n(Z_n(1)=0)=1$, and  $\E_{n}^{\GW}(\sum_{k\geq 1} Z_{k}(1))\leq\varphi(n)^{-1} \sum_{j\geq 1}  j p_j <\infty$. 

\noindent $\bullet$ \textsc{Assumption~\ref{A:BK}:}  It is natural  in this setting to use the selection rule that declares  legitimate one child\footnote{Choose one at random in case of ties.} with type $n-1$ and  all the other illegitimate. Recall the general definition  \eqref{eq:deflambdan} for the measure  $\boldsymbol{\Lambda}^{(n)}$.  To start with,  observe that since the legitimate child of a particle with type $n$ has always the type $n-1$, there is 
 the identity
 \begin{equation} \label{Eq:LambHeight}
 \boldsymbol{\Lambda}^{(n)}(\dd y_0, \dd \mathbf{y})= \delta_{\log(1-1/n)}(y_0) \boldsymbol{\Lambda}_1^{(n)}(\dd \mathbf{y}),
 \end{equation}
where $\boldsymbol{\Lambda}_1^{(n)}$ stands for the push-forward of $\boldsymbol{\Lambda}^{(n)}$  by the second projection $(y,{\mathbf y})\mapsto {\mathbf y}$ from $\mathcal S$ to $\mathcal S_1$. The requirements
 $$\lim_{n\to \infty} n \int_{-1}^1 \Lambda_0^{(n)}( \mathrm{d}y) y= -1 \quad\text{and}
 \quad \lim_{n\to \infty} n  \int_{-1}^1  \Lambda_0^{(n)}( \mathrm{d}y) y^2= 0
 $$
in Assumption \ref{A:BK}  are then plain. 

We now turn our attention to the asymptotic behavior of $\boldsymbol{\Lambda}_1^{(n)}$, and in this direction, we first point out that the probability that a particle of type $n$ has at least two children with type $n-1$ (that is, at least one illegitimate child with type $n-1$) is 
 \begin{equation}\label{E:kestenneysp}
 \Q^{\GW, \mathrm{height}}_n\big(Z_{n-1}(1)\geq 2\big) \leq  \frac{1}{\varphi(n)}  \sum_{j\geq 2} p_j \left(^j_2\right) \varphi(n-1)^2
\underset{ \eqref{E:kestenneys}}{=} O(n^{-2}) \qquad \text{ as }n\to\infty. 
 \end{equation}
 More precisely, in the right hand side of \eqref{E:kestenneysp}, the denominator $\varphi(n)$ is the probability that the mono-type Galton--Watson process has height $n-1$, the sum $ \sum_{j\geq 2}p_j$ accounts for the unconditional offspring distribution of the initial particle,  and finally the term $ \left(^j_2\right) \varphi(n-1)^2$  stands for the mean number pairs of trees with height $n-2$ in a family of $j$ independent mono-type Galton--Watson trees. \par
 Similarly, for any $x>0$, the probability that under $ \Q_{n}^{\GW, \mathrm{height}}$ the ancestor particle produces at least $3$ particles with types strictly larger than $xn$ (supposed to be an integer to simplify notation) is bounded above by 

 \begin{align}\label{E:kestenneyspi}
  &\boldsymbol{\Lambda}_1^{(n)}\left(\left\{ \mathbf{y}\in \mathcal{S}_1: y_2> -\log(x)\right\} \right)  
 \nonumber \\
  &\leq  \frac{\varphi(n-2)}{\varphi(n-1)}  \sum_{j\geq 3} p_j j(j-1)  (1-\Phi(\lfloor  x n \rfloor -1) )\min((j-2)(1-\Phi(\lfloor  x n \rfloor -1)),1)\nonumber \\
&  \underset{  \eqref{E:kestenneys:2}}{=} o(n^{-1}), \qquad \text{ as }n\to\infty, 
 \end{align}
We infer that
 $$ \boldsymbol{\Lambda}_1^{(n)}\left(\left\{ \mathbf{y}\in \mathcal{S}_1: y_1> \log x\right\} \right)   = n \frac{\varphi(n-1)}{\varphi(n)} \sum_{j=2}^{\infty}p_j j(j-1) \Big(\Phi(n-2)-\Phi(xn)\Big) \Phi(xn)^{j-2}
+ o(1),  $$
where the quantity $o(1)$ accounts for the event that a particle of type $n$ has  two or more children with type $n-1$ or two illegitimate  with type greater than $x n$. It is now straightforward to deduce from \eqref{E:kestenneys:2} that
 \begin{equation} \label{Eq:equivheights}  
 \lim_{n\to \infty}\boldsymbol{\Lambda}_1^{(n)}\left(\left\{ \mathbf{y}\in \mathcal{S}_1: y_1> \log x \right\} \right) 
 =  \frac{2}{x}-2,
 \end{equation}
for example use Fatou's lemma to obtain a lower bound for the liminf and the trivial bound $\Phi(xn-1)^{j-2}\leq 1$ to obtain the same upper bound for the limsup.
Combining \eqref{Eq:LambHeight}, and \eqref{Eq:equivheights}, we see that  Assumption \ref{A:BK} holds for the characteristic quadruplet $(0, -1, \boldsymbol{\Lambda}_{\mathrm{Height}}; 1)$.
 
 \noindent $\bullet$ \textsc{Assumption \ref{A:typex}:} We claim that  the associated discrete cumulant function $\kappa^{(n)}(\gamma)$ converges to the continuous one $\kappa_{\mathrm{Height}}(\gamma)=-\gamma+2/(\gamma-1)$ for all $\gamma>1$. In particular, this implies  Assumption \ref{A:typex}. This can be directly checked. Fix $\gamma>1$ and note that since we have already verified Assumption \ref{A:BK},  an application of Lemma \ref{lem:convkappacond} gives  $\liminf_{n\to \infty}  \kappa^{(n)}(\gamma)\geq \kappa_{\mathrm{Height}}(\gamma)$. Moreover, calculations similar to the ones performed above that
  \begin{align*}
\kappa^{{(n)}}(\gamma) &\leq n(1-1/n)^\gamma  -n
+  \frac{n  \varphi(n-1)}{\varphi(n)}  \sum_{j\geq 2} p_j j (j-1)\sum_{k=1}^{n-1}(k/n)^\gamma  \varphi(k)\\
&\leq  n(1-1/n)^\gamma  -n
+  \frac{ 2\varphi(n-2)}{\varphi(n-1)}  \cdot \Big( n^{-1} \sum_{k=1}^{n-1}(k/n)^{\gamma-2}  \widetilde{\varphi}(k)\Big),
\end{align*}
where $ \widetilde{\varphi}(k)=k^{2} \mathrm{Var}(\mathbf{p})\varphi(k) /2$.  Writing:
$$ n^{-1} \sum_{k=1}^{n-1}(k/n)^{\gamma-2}  \widetilde{\varphi}(k)= \int_{0}^{1} \big(\lfloor nx\rfloor/n\big)^{\gamma-2} \widetilde{\varphi}(\lfloor nx\rfloor) \d x +o(1),  $$
an application of dominated convergence, combined with \eqref{E:kestenneys}, gives  $\limsup_{n\to \infty}  \kappa^{(n)}(\gamma)\geq \kappa_{\mathrm{Height}}(\gamma)$.  
\end{proof}

\section{Random planar maps} \label{sec:randomplanarmaps}
In this section we show that the $(\beta, \varrho)$-stable family of ssMt with no killing of Example \ref{ex:stablefamily} are the scaling limits of labeled trees appearing in Markovian explorations of certain random discrete surfaces (Boltzmann planar maps). This enables us to understand certain aspects of their large scale geometry.

The discret objects we consider in this section are planar maps. We recall that  a \textbf{planar map} is a combinatorial object, made by gluing of finitely many polygonal faces along their edges so that the resulting manifold has the topology of the sphere. Alternatively, they can be seen as equivalence classes of finite connected graphs properly drawn on the sphere viewed up to a continuous deformation of the sphere. They are key objects in combinatorics and are still the subject of many investigations. We refer to \cite[Chapters I,II and III]{CurStFlour} for background. In the applications below, our planar maps are always given with a distinguished oriented edge $\vec{e}$ called the \textbf{root}; the face lying on the right of $\vec{e}$ is called the root or external face and is seen as the  \textbf{boundary} of the map. The other faces are called the inner faces. For technical simplicity, we shall only consider \textbf{bipartite} planar maps, i.e. each face has an even degree. In particular, the degree of the external face, called the \textbf{perimeter} below, is always even.

\subsection{Scaling limits for peeling trees of $ \mathbf{q}$-Boltzmann maps}\label{sec:peelingnormal}

\begin{figure}[!h]
 \begin{center}
 \includegraphics[width=12cm]{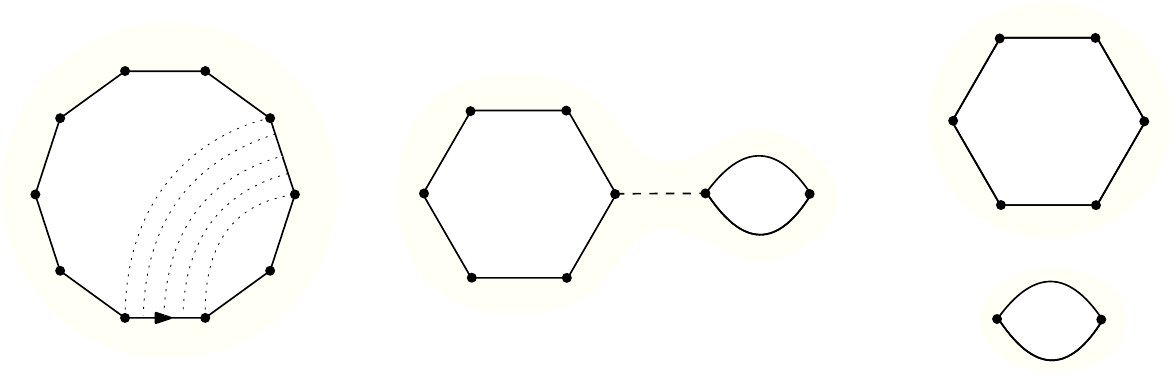}
 \caption{Illustration of case 2): if the root edge is adjacent on both sides to the external face (in which case the boundary of the external face is necessarily non-simple) then its erasure disconnects the map into two regions. \label{fig:peeling2}}
 \end{center}
 \end{figure}

Following the pioneer work of Tutte \cite{Tut63}, the recursive erasure of the root edge gives a mean to encode a planar map into a plane labeled tree. We refer to \cite[Section 14.3]{CurStFlour} for details and present here only the main idea. The algorithm is deterministic and proceeds as follows: imagine that a bipartite planar map $ \mathrm{m}$ is given so that the root face is of degree $2k$. We then erase the root edge; two topologically different cases may happen:
\begin{itemize} 
\item case 1) If the face of degree $2 \ell$ lying on the left of the root edge is not again the root face, then the erasure of the root edge leaves the map connected but increases the degree of the external face by $2 \ell -2$. See Figure \ref{fig:peeling1}.

\begin{figure}[!h]
 \begin{center}
 \includegraphics[width=12cm]{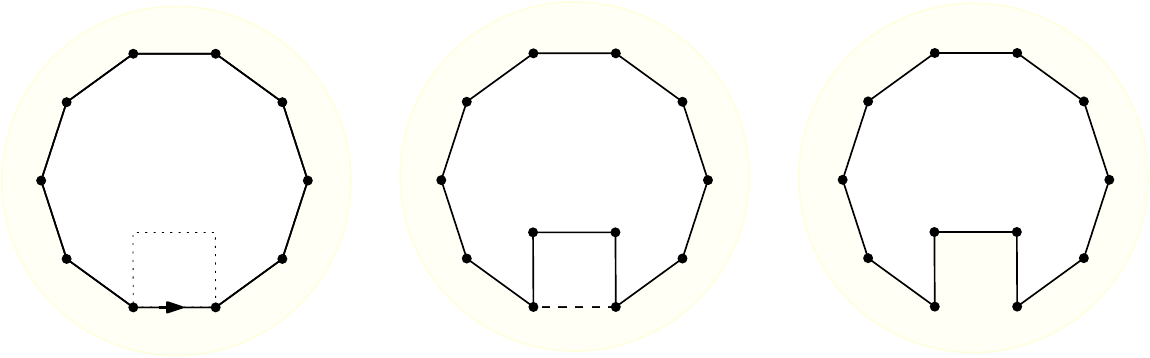}
 \caption{Illustration of case 1): the discovery of a new face of degree $2\ell$ increases the degree of the external face by $2\ell-2$. Notice that the boundary of the external face is not necessarily simple as in the drawing. \label{fig:peeling1}}
 \end{center}
 \end{figure}
 
\item case 2) If the root face is incident both on the left and right of the root edge, then the erasure of the latter splits the map into two connected components with external faces of degrees, say, $2k_{1}$ and $2k_{2}$ with $k_{1}+k_{2}=k-1$. Note that $k_{1}$ or $k_{2}$ (or both) may be equal to $0$. This corresponds to the vertex map $\dagger$ with a single vertex and no edge (hence boundary of perimeter $0$). See Figure \ref{fig:peeling2}.
\end{itemize}

\begin{figure}[!h]
 \begin{center}
 \includegraphics[width=15cm]{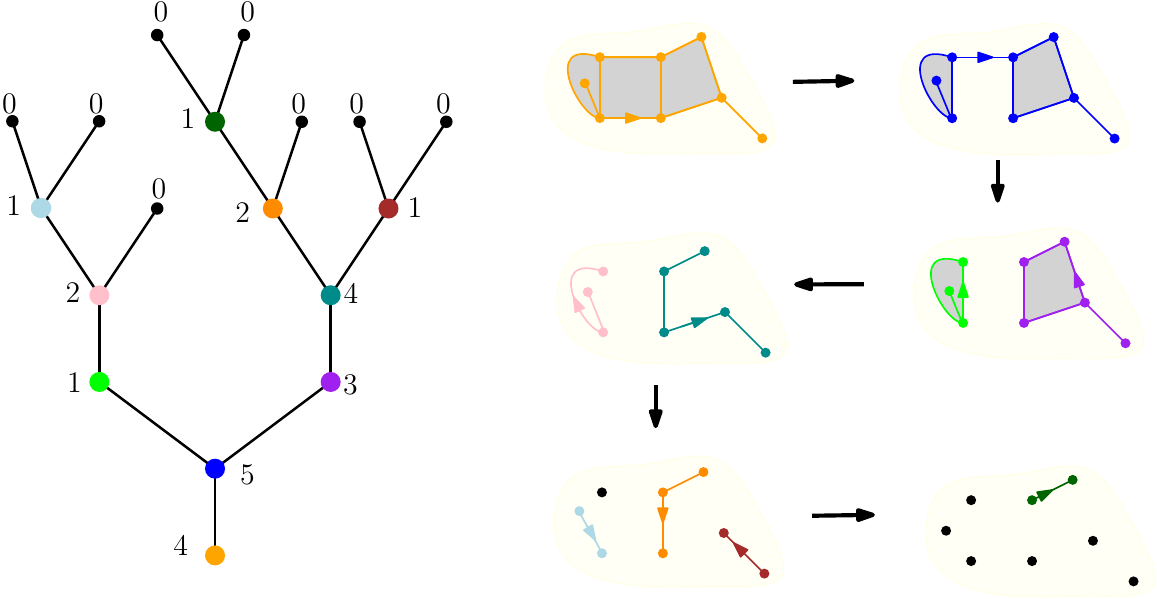}
 \caption{An example of a peeling exploration of a quadrangulation (all faces of degree 4) with  $3$ inner faces and a boundary of perimeter $10$. The resulting labeled plane tree is on the left. Its labels will later be shifted by one to fit the framework of Chapter \ref{chap:6}. \label{fig:branchingpeeling}}
 \end{center}
 \end{figure}

In both cases, we can select a new oriented edge on the boundary  of the at most two components and iterate the root erasure procedure. The choice of a new root edge depends on a \textbf{peeling algorithm}, denoted by $ \mathcal{A}$ in the following, which may depend on the past steps of exploration but not on the ``unexplored'' remaining maps (see  \cite[Section 4.1.2]{CurStFlour}), and different peeling algorithms may yield to different encodings. The iterative erasure of the root edge in each of the subsequent components can be recorded in a labeled plane tree that we denote by $$\mathrm{Peel}_{ \mathcal{A}}( \mathrm{m}),$$ see Figure \ref{fig:branchingpeeling}. Intuitively, each map with a boundary of length $2k$ appearing in the exploration is encoded by a vertex of label $k$ in  $\mathrm{Peel}_{ \mathcal{A}}( \mathrm{m})$. When case 2) happen, we represent  the component attached to the origin of the root edge on the left of the corresponding vertex in  $\mathrm{Peel}_{ \mathcal{A}}( \mathrm{m})$. In particular, $\mathrm{Peel}_{ \mathcal{A}}( \mathrm{m})$ is a finite plane binary tree whose vertices are labeled by $\{0,1,2, \dots \}$, the inner vertices of the tree correspond to the edges of the map, whereas the vertices labeled $0$ correspond to the vertices of the map (we shall shift the labels later to match the framework of the preceding chapters). See Figure \ref{fig:branchingpeeling} for an illustration.
 
 Such  explorations of maps are called \textbf{peeling explorations}, they have been systematically studied in \cite{CurStFlour}.  In particular, once $ \mathcal{A}$ is fixed, the encoding of the map into its labeled tree $ \mathrm{m} \mapsto \mathrm{Peel}_{ \mathcal{A}}( \mathrm{m})$ is indeed injective.  We now introduce Boltzmann measures on the set of planar maps. Given a  sequence of non-negative numbers  $\mathbf{q}= (q_{i} : i \geq 1)$, called the \textbf{weight sequence}, we can build the so-called  \textbf{$ \mathbf{q}$-Boltzmann measure} on the set of all finite bipartite planar maps defined by
\[\ w_{ \mathbf{q}}( \mathrm{m}):=\prod_{f \in \textup{Faces}( \mathrm{m})}q_{\frac{ \mathrm{deg}(f)}{2}}.\] The sequence $ \mathbf{q}$ is said to be \textbf{admissible} when $w_{ \mathbf{q}}$ gives a finite measure to the set of all rooted bipartite planar maps.  We suppose (aperiodicity condition) that for any $n \geq 1$ the set $ \mathcal{M}^{(n)}$ of all planar maps with perimeter (i.e. root face degree) $2n$ for $n \geq 1$ has a non-trivial $w_{ \mathbf{q}}$-weight which we denote by $W_{ \mathbf{q}}^{(n)}$ so that  
$$  \mathbb{W}_{ \mathbf{q}}^{(n)}:=\frac{w_{ \mathbf{q}}( \cdot \cap \mathcal{M}^{{(n)}})}{W_{ \mathbf{q}}^{{(n)}}},$$ is a probability distribution on $ \mathcal{M}^{(n)}$. We put $W_{ \mathbf{q}}^{(0)}=1$ for convention (this represent the weight of the vertex map). The following proposition is a consequence of the encoding presented above:

\begin{proposition} \label{Markov:peeling} Fix an admissible and aperiodic weight sequence $ \mathbf{q}$ and a peeling algorithm $ \mathcal{A}$. Then the law of $\mathrm{Peel}_{ \mathcal{A}}( \mathrm{m})$ under $\mathbb{W}_{ \mathbf{q}}^{(n)}$,  after forgetting the plane structure and  with all labels shifted by one, is an integer-type Galton--Watson process whose laws $ (\Q^{\GW, \mathrm{peel}}_{n})_{n \geq 1}$ are characterized by the only non trivial transitions:
$$ \begin{array}{ll} \mathbb{Q}^{\GW, \mathrm{peel}}_{n+1}( \mathbf{Z}(1)=\delta_{(n+k+1)}) = \displaystyle \frac{q_{k+1} W^{(n+k)}_{ \mathbf{q}}}{W^{(n)}_{ \mathbf{q}}}  & \mbox{ for } n \geq 1, k \geq 0\\\mathbb{Q}_{n+1}^{\GW, \mathrm{peel}}( \mathbf{Z}(1)=\delta_{n_{1}+1}+\delta_{n_2+1}) = (2- \mathbf{1}_{n_{1}=n_2}) \displaystyle \frac{ W^{(n_{1})}_{ \mathbf{q}} W_{ \mathbf{q}}^{{(n_2)}}}{W^{(n)}_{ \mathbf{q}}},  & \mbox{ for } \begin{array}{c}n \geq 1, \ \ n_{1},n_{2}\geq 0\\\mbox{ s.t. } n_{1}+n_{2} +1 =n.  \end{array}\\
 \mathbb{Q}^{\GW, \mathrm{peel}}_{1}( \mathbf{Z}(1)=  \mathbf{0})=1.& \end{array}$$
   \end{proposition}
\begin{proof} This is a direct consequence of the fact that the encoding is bijective and the multiplicative structure of the Boltzmann measure.  In particular, the factor $(2- \mathbf{1}_{n_{1}=n_{2}})$ is the symmetry factor arising when forgetting the plane structure of the branching. See Proposition 4.6 in \cite{CurStFlour} for more details.  \end{proof}

\begin{figure}[!h]
 \begin{center}
 \includegraphics[width=12cm]{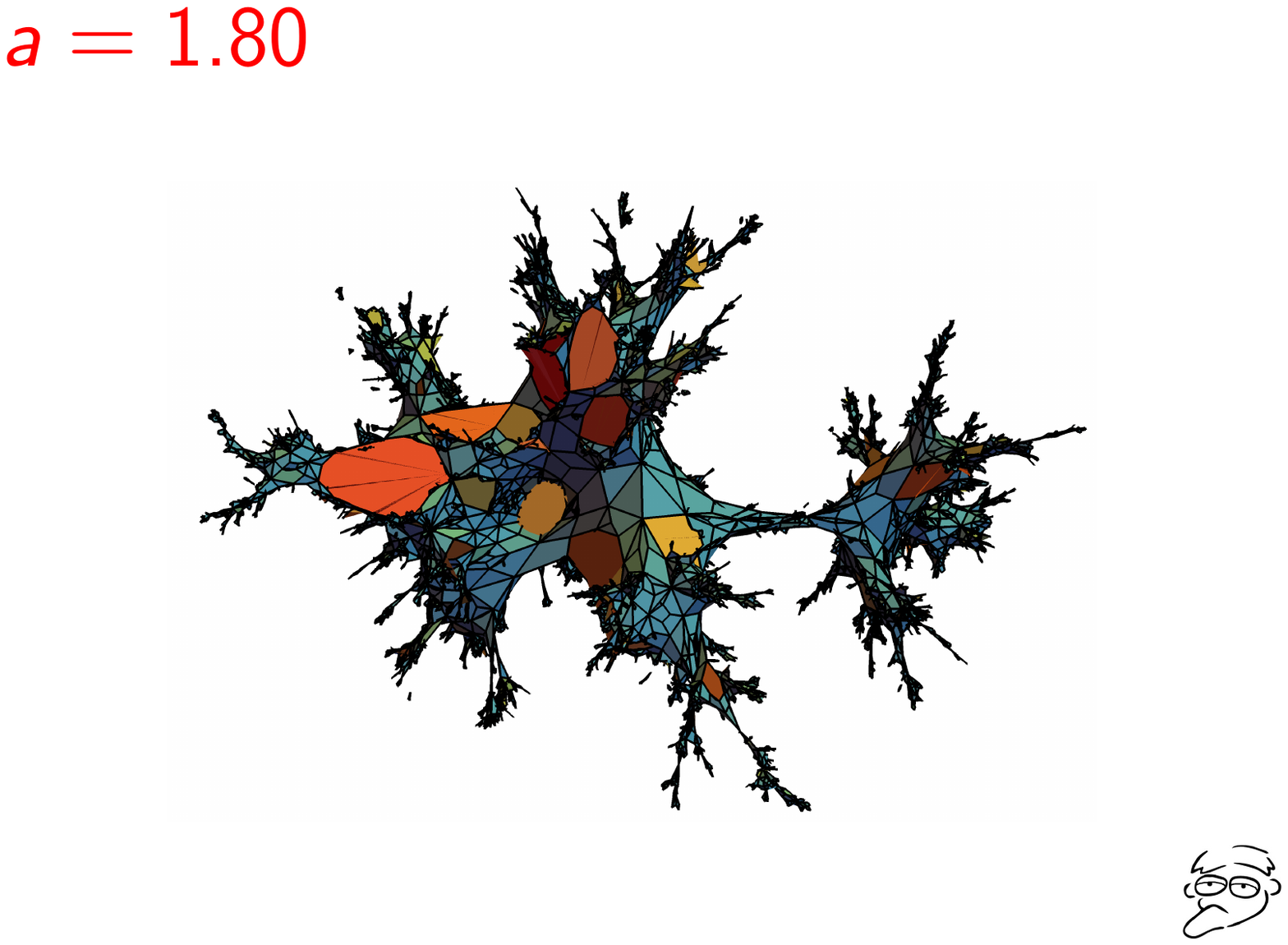}
 \caption{A simulation of large Boltzmann map with weight sequence of type  $\beta = 1.3$. We clearly see large faces persisting in the limit. \label{fig:stablemap}}
 \end{center}
 \end{figure}

The main result of this section (Theorem \ref{prop:GF}) is a scaling limit theorem for the decorated trees under $\Q^{\GW,\mathrm{peel}}_{n}$, but in order to obtain interesting limits, we shall require some properties on the weight sequence $ \mathbf{q}$. Specifically, we  introduce the function $$ f_{ \mathbf{q}}(z) = 1 + \sum_{k=1}^{\infty} q_{k} {2k-1 \choose k} z^{k},$$ and denote by $Z_{ \mathbf{q}}$ the smallest solution to the equation $f_{ \mathbf{q}}(x) = x$ which exists since $ \mathbf{q}$ is admissible (see Theorem 3.12 in \cite{CurStFlour}). Following \cite[Section 5.3]{CurStFlour}, we say that 
\begin{itemize}
\item $ \mathbf{q}$ is of \textbf{type $ \boldsymbol{\beta \in (1/2, 3/2]}$} if there exists $ \mathrm{c}>0$ such that
$$f_{ \mathbf{q}}(s)= Z_{\mathbf{q}}-\big(Z_{ \mathbf{q}}-s\big)+ \mathrm{c} \cdot \big(Z_{ \mathbf{q}}-s\big)^{\beta+1/2} +o\big(\big(Z_{ \mathbf{q}}-s\big)^{\beta+1/2}\big), \quad \text{as }s \uparrow Z_{\mathbf{q}}.$$
\item $ \mathbf{q}$ is  of \textbf{type} $ \boldsymbol{\beta = \frac{1}{2}}$ if there exists $0< \mathrm{c}<1$ such that
 $$f_{ \mathbf{q}}(s)= Z_{\mathbf{q}}-(1-\mathrm{c})\cdot \big(Z_{ \mathbf{q}}-s\big)+o(Z_q-s), \quad \text{as }s \uparrow Z_{\mathbf{q}}.$$
 \end{itemize}
 Notice the shift $\beta = a-1$ in the notation of \cite{CurStFlour}.  This terminology comes from the fact that in all those cases we have
  \begin{eqnarray} \label{eq:asymptWell} W^{{(\ell)}}_{ \mathbf{q}} \underset{\ell \to \infty}{\sim}   \frac{\mathrm{p}_{ \mathbf{q}}}{2} (c_{ \mathbf{q}})^{\ell+1} \ell^{-\beta-1},  \end{eqnarray}
for  $  c_{ \mathbf{q}} = 4 Z_ \mathbf{q}$ and some $ \mathrm{p}_{ \mathbf{q}}>0$ that will be used below. We refer to  Propositions 5.7, 5.9 and 5.10 in \cite{CurStFlour} for a proof. The construction of weight sequences of type $\beta \in (1/2, 3/2)$ needs some fine tuning, we refer to \cite[Section 7.2]{CRM22} for a discussion, and Section \ref{sec:O(n)} below for an example coming from the critical $O( \mathrm{n})$-model. The large scale geometries of random maps sampled according to $\mathbb{W}_{ \mathbf{q}}^{(n)}$ really differ depending on $\beta$:
\begin{itemize}
\item If $\beta=1/2$, then  informally, those random maps fold on their boundary and are tree-like, there are believed to converge after renormalization of the distances by $ n^{-1/2}$ towards the Brownian CRT of Example \ref{ex:brownian}. See \cite{Bet11} for the case of quadrangulations, and \cite[Section 4]{Richier17b} for a general discussion.
\item If $\beta=3/2$, it is expected that the random maps converge after renormalization of the distances by $n^{-1/4}$ towards the (free) Brownian disk, see Example \ref{ex:3/2stable}. This has been proved in  \cite{gwynne2019convergence} for quadrangulations with a simple boundary or \cite{BM15} for the case of critical regular weight sequences. We will show that their peeling trees converge in the scaling limit towards the Brownian Growth-Fragmentation tree of Example \ref{ex:3/2stable} but with self-similarity parameter $3/2$ instead of $1/2$ .
\item If $1/2<\beta<3/2$, the random maps should converge after renormalization of the distances by $n^{- \frac{1}{2\beta +1}}$ towards a variant of the stable gasket/carpet  of \cite{CRM22}. In particular, a transition appears at $\beta=1$: If $ \beta \in (1/2,1)$, the so-called dense phase, the stable gasket has cut-points; whereas if $\beta \in [1,3/2)$ it is almost surely homeomorphic to the Sierpinski carpet, see Figure \ref{fig:stablemap} and \cite[Section 8]{CRM22}. This topological phase transition is reminiscent of the same transition for the conformal loop ensembles of parameter $\kappa \in (8/3,8)$ \cite{SheCLE}, and a strong link is supposed to hold  when 
 $$ \beta = \frac{4}{ \kappa}.$$

\begin{figure}[!h]
 \begin{center}
 \includegraphics[height=4cm]{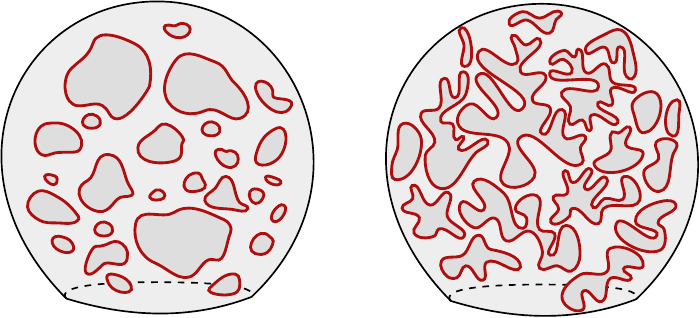}
 \caption{A schematic illustration of $ \mathbf{q}$-Boltzmann random planar maps of type $\beta \in (1/2;3/2)$ in the dilute (left) and dense phase (right) where their scaling limits are respectively the stable carpet and stable gasket of \cite{CRM22}.}
 \end{center}
 \end{figure}
 
  In terms of the peeling trees of the maps, no such transition appears and we will show that they converge in the scaling limit towards the ssMt of Example \ref{ex:stablefamily} with no killing. This will be used in Corollary \ref{cor:cactus} to deduce that the diameter of the \textbf{dual} of those maps is tight once renormalized by $ n^{-(\beta-1)}$ when $\beta > 1$. This is a first step towards understanding the scaling limits of the dual maps, see \cite{kammerer2024scaling} for a recent progress in this direction.
\end{itemize}

In order to state our scaling limits on the peeling trees, let us introduce a few scaling constants: Recalling \eqref{eq:asymptWell} we let $$  \mathrm{b}^{\mathrm{vertices}}_{ \mathbf{q}} = \frac{2}{c_{ \mathbf{q}}  \mathrm{p}_{ \mathbf{q}} \sqrt{\pi}} \quad \mbox{ and }\quad \mathrm{d}_{ \mathbf{q}} = \frac{\pi \cdot \mathrm{p}_{ \mathbf{q}}}{\Gamma(1+\beta)},$$ the first one being defined in  \cite[Eq. (10.7)]{CurStFlour} and the second one in \cite[Proposition 6.6]{BBCK18} (unfortunately called $ \mathrm{c}_{ \mathbf{q}}$). Those constants depend on the fine details of the weight sequence $ \mathbf{q}$. We then have:

\begin{theorem}[Scaling limits for peeling trees]  \label{prop:GF} Fix a peeling algorithm $ \mathcal{A}$ and let $ \mathbf{q}$ be a weight sequence of type $\beta \in [1/2;3/2]$. Considering the integer-type Galton--Watson tree of law $ \Q^{\GW,\mathrm{peel}}_{n}$ described in Proposition \ref{Markov:peeling}, the sequence of the distributions  of the rescaled measured decorated trees 

$$\left( T_\GW, \frac{d_{T_\GW}}{n^{\beta}} \cdot \mathrm{d}_{ \mathbf{q}} , \rho_\GW,  \frac{g_\GW}{n}, \frac{ \uplambda^{ \mathbf{1}_{\{1\}}}_\GW}{n^{\beta+1/2}} \cdot \frac{1}{\mathrm{b}^{ \mathrm{vertices}}_{ \mathbf{q}}}\right) \mbox{ under } \Q^{\GW,\mathrm{peel}}_{n}$$
 converges as $n\to \infty$,  to the law  of 
the measured self-similar Markov tree $(\normalfont{\texttt{T}},\upmu)$  under $\Q_1$ associated with the characteristics $(0,  \mathrm{a}_{ \beta-1/2,\beta}, \boldsymbol{\Lambda}_{ \beta-1/2,\beta}; \beta)$ of  Example \ref{ex:stablefamily} when $ \texttt{b}=1/2, \texttt{a} = \beta-1/2$ so that $\frac{1}{2} \leq \beta = \texttt{a} + \texttt{b} \leq \frac{3}{2}$ and $\varrho$ such that $\beta(1-\varrho)= 1/2$ (i.e. there is no killing term), see  \eqref{eq:nokiling} and \eqref{eq:kappanokilling}.

\end{theorem}
\begin{remark} A related but weaker convergence  has already been proved in \cite[Section 6]{BBCK18} in terms of growth-fragmentation processes. That convergence did not  guarantee the convergence in the Gromov--Hausdorff topology for the underlying tree. However, we will employ similar inputs in our proof.
\end{remark}

\begin{proof} The reader will have recognized an application of Theorem \ref{T:mainunconds-mass}, from which (by Theorem \ref{Th:cvdecoreprod} and Proposition \ref{prop:PHI:ex}) it suffices to check Assumptions \ref{A:all:GW}, \ref{A:BK} and  \ref{A:typex} as well as Assumption \ref{A:hGWrv}. Let us proceed:

\noindent $\bullet$ \textsc{Assumptions \ref{A:hGWrv} and \ref{A:all:GW}:} Recall from the peeling encoding (and the shift of label made in Proposition \ref{Markov:peeling}) that the particles of type $1$  in $T_{\GW}$ correspond to the vertices of the underlying map. By  \cite[Proposition 10.4]{CurStFlour} we have
$$\E_n^{\GW}\big(\uplambda_{\GW}^{ \mathbf{1}_{\{1\}}}(T_\GW)\big) \underset{n\to \infty}{\sim}  \mathrm{b}^{\mathrm{vertices}}_{ \mathbf{q}} \cdot  n^{\beta+1/2},$$
and $n^{-\beta-1/2}\cdot \uplambda_{\GW}^{ \mathbf{1}_{\{1\}}}(T_\GW)$, under $\Q_n^{\GW,\mathrm{peel}}$, converges in distribution to a random variable with mean $\mathrm{b}^{\mathrm{vertices}}_{ \mathbf{q}}$, as $n\to \infty$. In particular, this entails Assumption  \ref{A:hGWrv}. Since type $1$ is clearly accessible, the subcriticality is implied by the above display and so Assumption \ref{A:all:GW} clearly holds.

\noindent $\bullet$ \textsc{Assumption \ref{A:BK}:} We shall now establish the convergence of the characteristics. For this, we shall naturally use the locally largest exploration as selection rule in the discrete setting. We write $(\sigma^2=0, \mathrm{a}, \boldsymbol{\Lambda} ; \beta)$ for the characteristic quadruplet (locally largest bifurcator) associated with the ssMt $(T,g)$ of Example \ref{ex:stablefamily} with $\beta \in [1/2,3/2]$, no-killing and with the normalization \eqref{eq:nokiling} and \eqref{eq:kappanokilling}. The dependency of $  \mathrm{a}, \boldsymbol{ \Lambda}$ on $\beta$ is implicit here. If $F : \mathcal{S} \to \mathbb{R}$ is measurable then by Proposition \ref{Markov:peeling} we have 
  \begin{eqnarray}  \int \boldsymbol{\Lambda}^{(n)}(\d \mathbf{y})  F( \mathrm{e}^{y_0}, \mathrm{e}^{y_1},\dots) 
 &=&n^\beta \cdot \sum_{k \geq 0}\frac{q_{k+1} W^{(n-1+k)}_{ \mathbf{q}}}{W^{(n-1)}_{ \mathbf{q}}} F\Big(\frac{n+k}{n},0,0,\cdots\Big) \nonumber \\
& +& 2 n^{\beta}\cdot \sum_{k>\frac{n-1}{2}}^{n-2} \frac{ W^{(n-k-2)}_{ \mathbf{q}} W_{ \mathbf{q}}^{{(k)}}}{W^{(n-1)}_{ \mathbf{q}}} F\left( \frac{k+1}{n} , 1-\frac{k+1}{n}, 0, 0, \cdots\right)\nonumber \\
 & +& \mathbf{1}_{\frac{n-2}{2}\in \mathbb{N}} n^{\beta}\cdot \frac{ W^{(\frac{n-2}{2})}_{ \mathbf{q}} W_{ \mathbf{q}}^{{(\frac{n-2}{2})}}}{W^{(n-1)}_{ \mathbf{q}}}F\Big( \frac{n}{2}, \frac{n}{2}, 0, 0, \cdots \Big). \nonumber \end{eqnarray}
We first notice that the last term, corresponding to the event when a particle of type $n$ has  two  children with type $n/2$, is negligible using the asymptotic  \eqref{eq:asymptWell}. It is also convenient to introduce  $h(m):= W_{\mathbf{q}}^{(m-1)} c_{\mathbf{q}}^{-(m-1)}$, for every $m\geq 1$, and 
$$\nu(m)\coloneq \left\{ \begin{matrix}
{q}_{m+1} c_{ \mathbf{q}}^{m} &\text{if }m\geq 0,\\
2 W_{\mathbf{q}}^{(-m-1)}  c_{ \mathbf{q}}^{m} &\text{if }m<0.
\end{matrix} \right.
$$
In particular,  we have
 \begin{eqnarray}\int \boldsymbol{\Lambda}^{(n)}(\d \mathbf{y})  F( \mathrm{e}^{y_0}, \mathrm{e}^{y_1},\dots) &=&  
n^\beta\cdot  \sum_{k \geq 0}  \frac{h(n+k)}{h(n)} \nu(k)  F\Big(\frac{n+k}{n},0,\cdots\Big)\nonumber \\
&+& n^\beta\cdot
 \sum_{1 \leq k < \frac{n}{2}} \nu(-k) \frac{h(n-k)}{h(n)} F\left( \frac{n-k}{n} , \frac{k}{n}, 0,  \cdots\right)\nonumber \\
 &+& o(1). \label{eq:lambdanpeel} \end{eqnarray}
The first key fact is that by \cite[Lemma 5.2 and Proposition 10.1]{CurStFlour} the measure $\nu$ is a probability measure which is in the strict domain of attraction of the $\beta$-stable law. 
Namely, we have  \begin{equation} \label{eq:tailnu}\nu([k, \infty)) \sim \mathrm{p}_{\mathbf{q}}\cos\big((\beta+1)\pi\big) \beta^{-1} \cdot k^{-\beta} \:\: \mbox{  and }   \:\: \nu((-\infty, -k]) \sim \mathrm{p}_{\mathbf{q}} \beta^{-1} \cdot k^{-\beta}, \quad \mbox{as $k \to \infty$}.  \end{equation} When $F$ is a continuous function which vanishes in a neighborhood of $(1,0,\cdots)$, then in \eqref{eq:lambdanpeel} we can restrict to values of $k$ such that $ k > \varepsilon n$ and $n-k > \varepsilon n$, so that using \eqref{eq:asymptWell} and the above remarks, it is easily seen that 
  \begin{eqnarray*} \int \boldsymbol{\Lambda}^{(n)}(\d \mathbf{y})  F( \mathrm{e}^{y_0}, \mathrm{e}^{y_1},\dots)   &\xrightarrow[n\to\infty]{  }&  \mathrm{p}_{\mathbf{q}} \int_{1/2}^1 \frac{\mathrm{d}x}{(x(1-x))^{\beta+1}}F( x,1-x, 0,\cdots)\\  &+& \mathrm{p}_{\mathbf{q}} \cos\big((\beta+1)\pi\big) \int_{0}^\infty \frac{\mathrm{d}x}{(x(1+x))^{\beta+1}}F( 1+x,0, 0,\cdots)\\  &\underset{ \eqref{eq:nokiling}}{=}&  \mathrm{d}_{ \mathbf{q}}\cdot \int \boldsymbol{\Lambda}(\d \mathbf{y})  F( \mathrm{e}^{y_0}, \mathrm{e}^{y_1}, \dots).    \end{eqnarray*}
 The convergence of the generalized L\'evy measure in Assumption \ref{A:BK} is thus granted. The variance term is easily dealt with. Specifically, using \eqref{eq:asymptWell} it follows that for some constant $C>0$ (that may vary from line to line) we have

  \begin{align}\int_{(- \varepsilon, \varepsilon)} \Lambda_0^{(n)}(\d y) y^2 & \leq C n^\beta \cdot \Big(\sum_{0\leq k \leq 2 \varepsilon n } \nu(k) \log\big(\frac{n+k}{n}\big)^2 + \sum_{1\leq k \leq   2 \varepsilon n} \nu(-k) \log\big(\frac{n-k}{n}\big)^2\Big) +o(1) \nonumber \\
  & \leq  C n^\beta\cdot \Big(\sum_{0\leq k\leq  2 \varepsilon n} \nu(k) \big( \frac{k}{n}\big)^2 + \sum_{1\leq k\leq 2 \varepsilon n} \nu(-k) \big( \frac{k}{n}\big)^2 \Big) +o(1) \nonumber \\
  &\leq C \varepsilon^{ 2- \beta}, \label{eq:controlvariance}
   \end{align}
    for all $\varepsilon\in (0,1/2)$ and $n \geq 1$ large enough; where to obtain the last inequality we use  \eqref{eq:tailnu} and an Abel summation.
   This implies the convergence of the variance term in Assumption \ref{A:BK} with $\sigma^2=0$. It remains the most delicate convergence of the drift term. When $\beta \in [1/2,1)$ the same calculation as in the previous display yields 
   $$ \int_{( -\varepsilon, \varepsilon)} \Lambda_0^{(n)}(\d y) y \leq C \varepsilon^{1- \beta}.$$ In particular, using \eqref{eq:tailnu} again we deduce 
   $$\int_{(- 1,1)} \Lambda_0^{(n)}(\d y) y \xrightarrow[n\to\infty]{} \int_{(- 1, 1)} \Lambda_0(\d y) y,$$
   which implies the convergence of the drift term where $ \mathrm{a_{can}}=0$ in this case. However, when $\beta \in [1,3/2]$ delicate compensations are taking place and we need to be more careful. We shall use the key fact that for $\beta \geq 1$, the law $\nu$ is in the domain of attraction of the aforementioned $\beta$-stable law \textit{without centering}, in particular $\nu$ has mean zero when $\beta >1$. Furthermore we shall rely on the following lemma which sharpens \eqref{eq:asymptWell}:
   \begin{lemma}[Strong ratio limit theorem] \label{lem:strongratio} If $ \mathbf{q}$ is of type $\beta \in [1/2, 3/2]$ then, there exists $C>0$ such that  
   $$  \left| \frac{h(n+k)}{h(n)}-1\right| \leq C\frac{|k|}{n},$$
   for every $n\geq 1$ and all $k \in [-n/2, ( \mathrm{e}-1)n]$.
   \end{lemma}
   \begin{proof} The proof is based on the representation of $W^{(\ell)}_{ \mathbf{q}}$ as a Laplace integral. Specifically, from \cite[Eq. (5.7)]{CurStFlour} we read that 
  \begin{eqnarray} W^{(\ell)}_{ \mathbf{q}}  = \frac{\ell}{\ell+1} {2 \ell \choose \ell} \int_0^1 \mathrm{d}u \frac{(\phi(u))^{\ell+1}}{u^2}  = \frac{\ell}{\ell+1} {2 \ell \choose \ell}  (Z_ { \mathbf{q}})^{\ell+1} \int_0^1 \mathrm{d}u \frac{(\tilde{\phi}(u))^{\ell+1}}{u^2}, \end{eqnarray} where     $\phi$ the inverse function of $ z \mapsto z/ f_{ \mathbf{q}}(z) : [0,Z_{ \mathbf{q}}] \to [0,1]$ so that $\tilde{\phi} = \phi/Z_{ \mathbf{q}}$ is an increasing function $[0,1] \to [0,1]$, equal to $1$ at $1$ and which under the assumption on $ \mathbf{q}$ via $f_{ \mathbf{q}}$ satisfies $\tilde{\phi}(1-x)=1  - \mathrm{A} \cdot x^{\tilde{\beta}} + o(x^{\tilde{\beta}})$ for some $ \mathrm{A} >0$ and where $\tilde{\beta} = ( \beta-1/2)^{-1}$. Applying Lemma \ref{lem:laplaceratio} with $g(x) = \tilde{\phi}^{2}(x)/x^{2}$ (so that it is continuous at $0$) and $f(x) = \tilde{\phi}(x)$ we deduce that 
  $$ \frac{\int_0^1 \mathrm{d}u \frac{(\tilde{\phi}(u))^{n+1}}{u^2}}{\int_0^1 \mathrm{d}u \frac{(\tilde{\phi}(u))^{n}}{u^2}} \in\left[1 - \frac{C}{n}, 1+  \frac{C}{n} \right],$$ and the lemma follows easily for general $k$ using the penultimate display and the definition of $h$. \end{proof}

   Coming back to the drift term in the case $\beta \geq 1$, the above lemma and \eqref{eq:asymptWell} imply that there exists  yet another $C>0$ such that 
   $$\forall k \geq -n/2, \quad \left| \frac{h(n+k)}{h(n)}-  \left(\frac{n+k}{n}\right)^{-\beta-1} \right| \leq  C\cdot \Big(  \frac{|k|}{n} \wedge \delta_{n} \Big), \quad \mbox{ where } \delta_{n}  \xrightarrow[n\to\infty]{} 0.$$
   Plugging back inside \eqref{eq:lambdanpeel}, we deduce that 
    \begin{eqnarray*} \left|\int_{- 1}^{ 1} \Lambda_0^{(n)}(\d y) y - n^\beta \sum_{-n/2\leq k \leq (\mathrm{e}-1) n }  \nu(k) \left( \frac{n+k}{n} \right)^{-\beta-1} \log \left(\frac{n+k}{n}\right) \right| &&  \end{eqnarray*}  
    \begin{eqnarray*}
    &\leq & C\cdot  n^{\beta} \Big(\sum_{-n/2\leq k \leq ( \mathrm{e}-1)n } \nu(k) \Big(  \frac{k}{n} \wedge \delta_{n}\Big) \log\Big(\frac{n+k}{n}\Big) \Big) +o(1). \\
\end{eqnarray*}
The function $\left|\log(\frac{n+k}{n})\right|$, over $-n/2\leq k\leq ( \mathrm{e}-1)n $, can itself be bounded above by $ C \frac{k}{n}$ for possibly another constant $C>0$. Using \eqref{eq:tailnu} the above term is controlled  as in \eqref{eq:controlvariance} and shown to converge to $0$ as $n \to \infty$. On the other hand, one can write 
  \begin{equation}n^\beta\cdot \sum_{-n/2\leq k \leq(\mathrm{e}-1) n }  \nu(k) \left( \frac{n+k}{n} \right)^{-\beta-1} \log \left(\frac{n+k}{n}\right) =  n^\beta\cdot \sum_{k =-n/2}^{(\mathrm{e}-1) n }  \nu(k) \varphi(k/n) ,  \end{equation}
with $\varphi(x) = (1+x)^{-\beta-1} \log (1+ x) = x - (\beta+3/2) x^{2} + o(x^{2})$.  Now, recall from \cite[Lemma 5.2 and Proposition 10.1]{CurStFlour} that the measure $\nu$ is in the strict domain of attraction \textit{without centering} of the $\beta$-stable law with positivity parameter $\varrho$ satisfying $\beta(1- \varrho)=1/2$. In particular, when $\beta >1$ this implies that $\nu$ is centered and in the case $\beta=1$ this is equivalent to a more subtle truncation centering, see \cite[Proposition 10.1]{CurStFlour} or \cite[Proposition 2]{BCMcauchy}. By standard result e.g.\ \cite[Theorem 3.4, Chap VII]{JS03}, this implies the convergence of infinitesimal generators and in particular $n^\beta\cdot \sum_{-n/2\leq k\leq (\mathrm{e}-1) n }  \nu(k) \varphi(k/n)$ converges to
  \begin{eqnarray*}  \mathrm{p}_{  \mathbf{q}}\left(  \mathrm{A} \varphi^{\prime}(0)+   \int  \left( \frac{ \mathrm{d}y}{|y|^{\beta+1}} \mathbf{1}_{-1/2 < y <0} +  \cos((\beta+1) \pi) \frac{ \mathrm{d}y}{|y|^{\beta+1}} \mathbf{1}_{0 < y < \mathrm{e}-1}\right)  \big( \varphi(y)- y \varphi^{\prime}(0)\big)\right),  \end{eqnarray*} 
  as $n\to \infty$,  where
$$ \mathrm{A}= \frac{1- \cos ( \pi( \beta+1))}{\beta-1}$$ so that $\sigma^{2}=0, \mathrm{A},  \frac{ \mathrm{d}y}{|y|^{\beta+1}} \mathbf{1}_{y <0} + \cos ( \pi(\beta+1)) \frac{ \mathrm{d}y}{|y|^{\beta+1}} \mathbf{1}_{y >0}$ are the characteristics of a $\beta$-stable L\'evy process with the truncation used in the L\'evy Khintchine formula \eqref{E:LKfor}. It is then a straight-forward but tedious calculation to see that this matches \eqref{drift:a:b}.

   \noindent $\bullet$ \textsc{Assumption \ref{A:typex}:} Let us prove that the discrete cumulants $\kappa^{(n)}(\gamma)$ converge to the associated continuous cumulant $\kappa(\gamma)$ of \eqref{eq:kappanokilling}  for every $\gamma\in (\beta, 2\beta+1)$. Of course this entails Assumption \ref{A:typex}, since $\kappa(\gamma)<0$ for $\gamma\in(\beta+1/2, \beta +3/2)$. Since we have verified Assumption \ref{A:BK} above, one can use Lemma \ref{lem:convkappacond} whose hypotheses for $\gamma \in (\beta, 2\beta+1)$ are easily verified in our case using \eqref{eq:tailnu}.\end{proof}

\subsection{Applications to distances on the dual maps}
In the preceding section, we proved convergence of the peeling trees encoding  random Boltzmann planar map for \textit{any} peeling algorithm. By carefully choosing the peeling algorithm, this can imply scaling limits for distances in the \textbf{dual} map. This will  require to perform a Lamperti transformation (a.k.a.~subordination) in the discrete setting to change the self-similarity parameter  from $\beta$ to $\beta-1$. \medskip 

\begin{figure}[!h]
 \begin{center}
 \includegraphics[width=6cm]{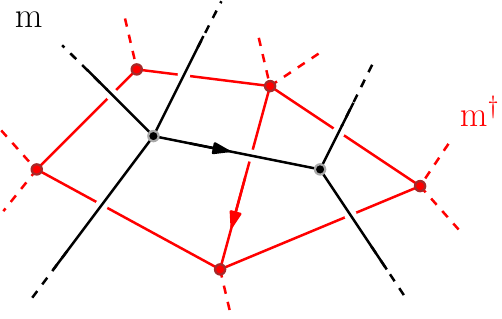}
 \caption{A (piece of a) map and (a piece of) its dual, notice that the roles of vertices and faces are exchanged. The root edge is transferred naturally from $ \mathrm{m}$ to $ \mathrm{m}^{\dagger}$.  \label{fig:dual}}
 \end{center}
 \end{figure}

\begin{figure}[!h]
 \begin{center}
 \includegraphics[width=16cm]{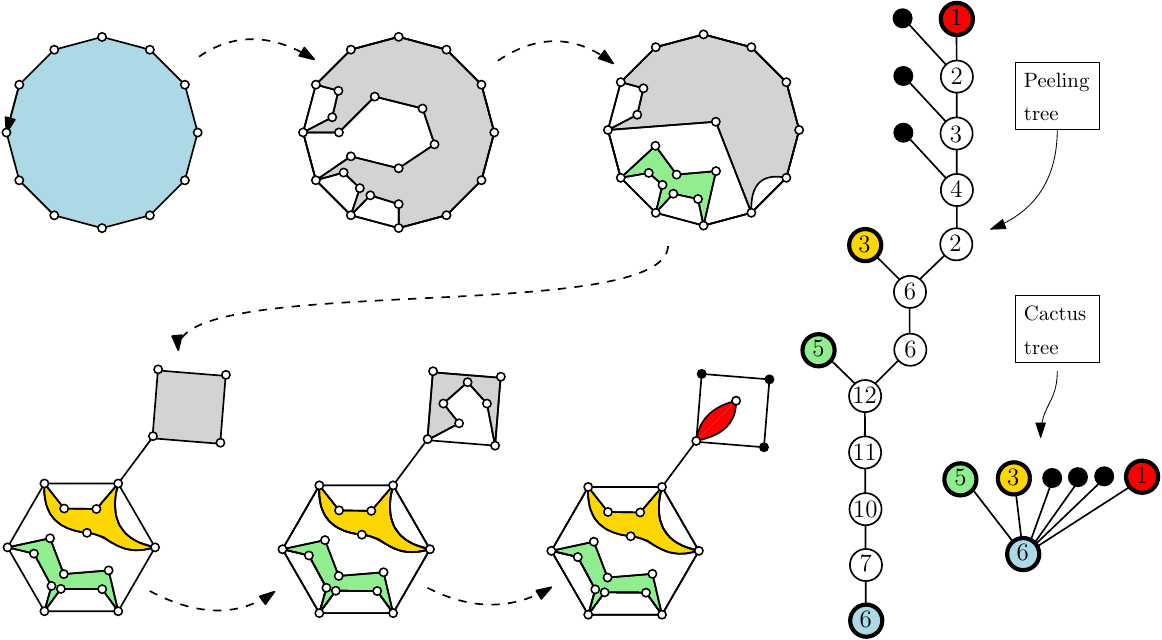}
 \caption{Illustration of the peeling by layers around one face stopped at steps $4,6,7,8$ and $11$. We peel all the incident edges of the face (here turning in counterclockwise order around the face) until all the faces adjacent to the initial boundary have been discovered. The vertices corresponding to vertices of  the $ \mathrm{Cactus}$ tree  are highlighted by thick circles on the  peeling tree. Right: By ignoring all labeled vertices that are not in the cactus tree yields to a single step for the multi-type Galton--Watson tree corresponding to the cactus tree. \label{fig:layersoneface}}
 \end{center}
 \end{figure}

\begin{figure}[!h]
 \begin{center}
 \includegraphics[width=16cm]{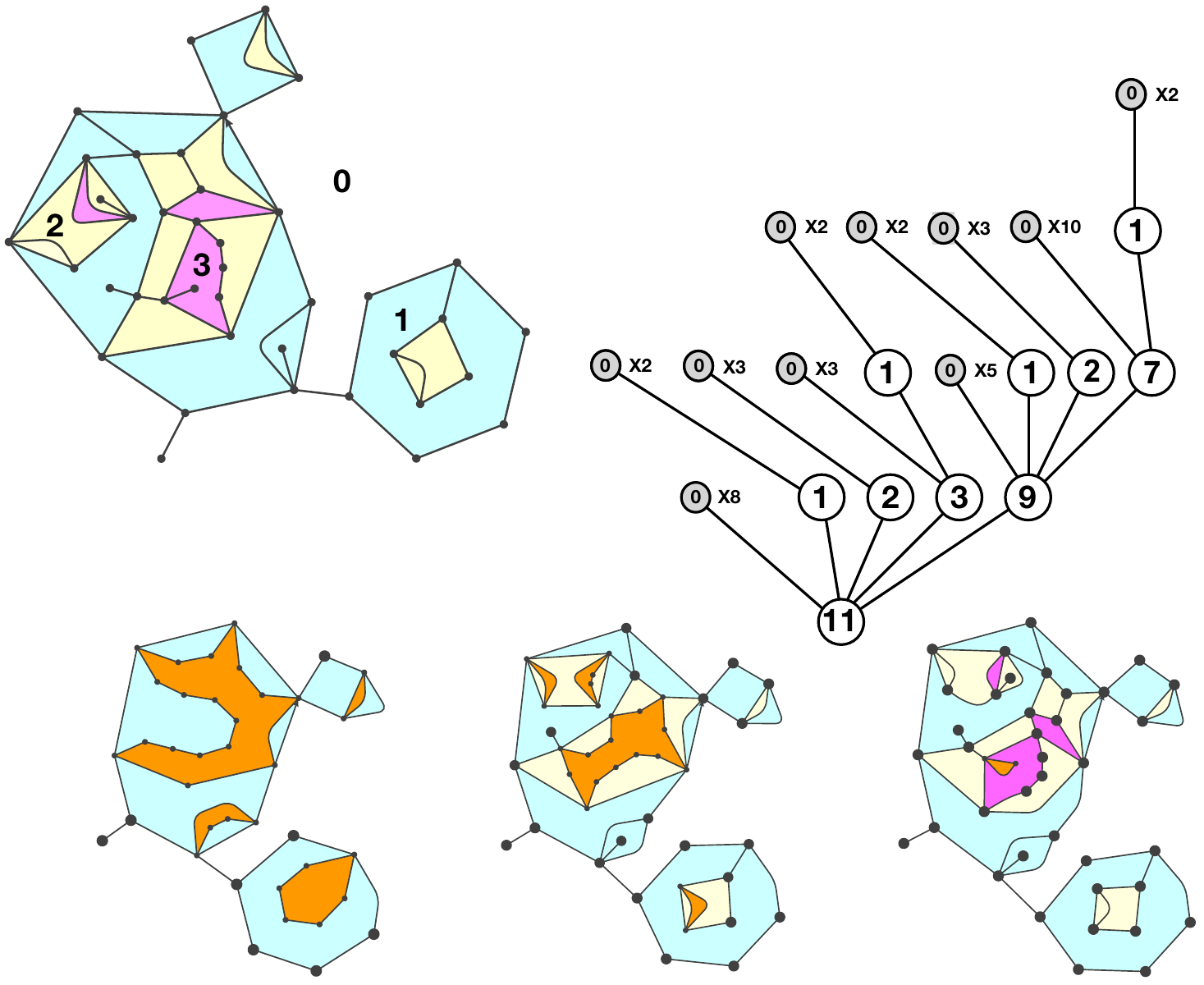}
 \caption{Example of the peeling by layers of a planar map (upper left corner) using the dual distance. We stopped the peeling when the (dual) ball of radius $1,2$ and $3$ have been discovered. The orange regions are the holes of these balls, i.e.~the regions which are still to be explored by the algorithm.  On the right is depicted the subordinated labeled tree representing the perimeters of the holes of those balls, as well as the vertices of label $0$ which correspond to the vertices of the map. \label{fig:labeleddual}}
 \end{center}
 \end{figure}

\begin{figure}[!h]
 \begin{center}
 \includegraphics[height=7.5cm,angle=-3]{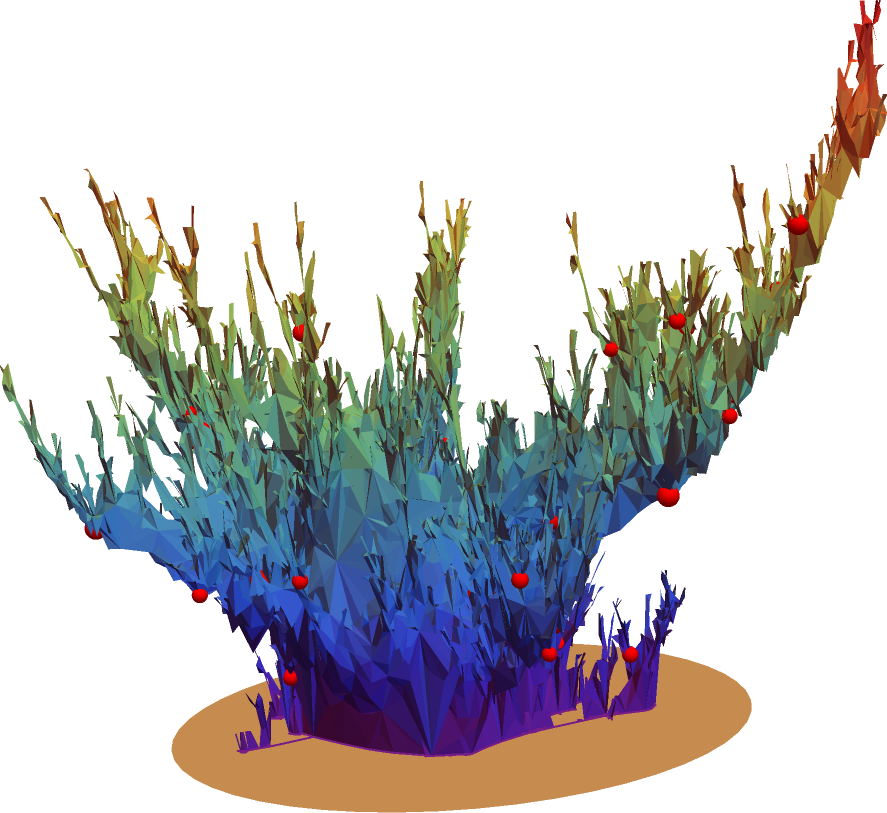}
  \includegraphics[height=7.5cm]{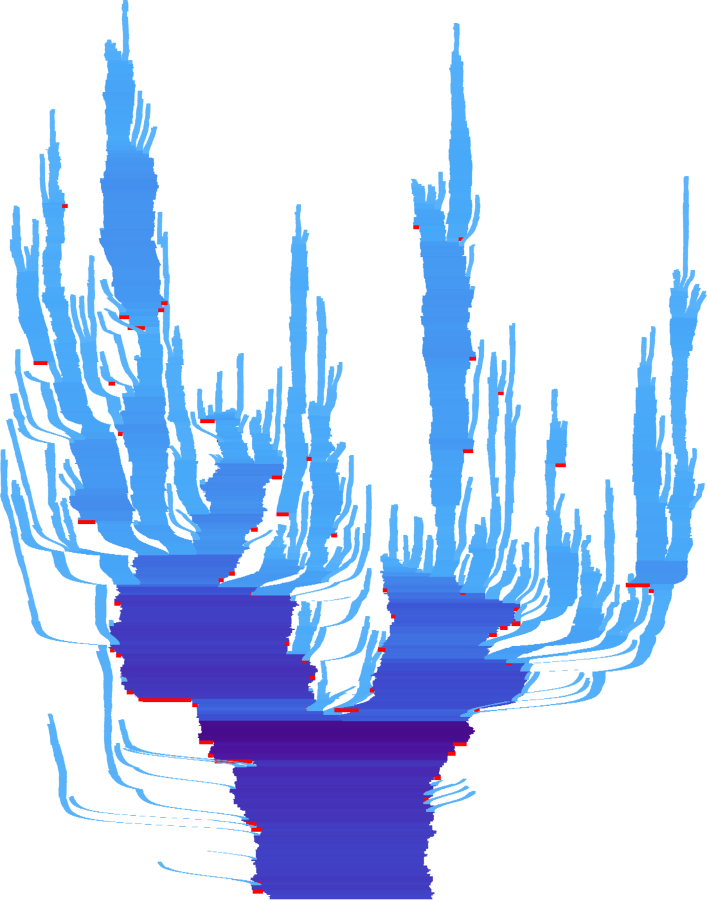}
 \caption{Left: A cactus representation of a random planar map with certain vertices of high degrees (which are the red dots), where the height of a vertex is its distance to the orange boundary. Right: A (different) simulation of the growth-fragmentation process describing the scaling limit of its perimeters at heights (the red part corresponds to positive jumps of the process). \label{fig:cactuses}}
 \end{center}
 \end{figure}

Recall that the \textbf{dual} map $ \mathrm{m}^{\dagger}$ of a planar map $ \mathrm{m}$ is the planar map obtained by exchanging the roles of vertices and faces but keeping the same incidence relations, see Figure \ref{fig:dual} for an illustration.
The \textbf{peeling by layer on the dual} is a peeling algorithm $ \mathcal{A}_{ \mathrm{dual} }$ that discovers a map $ \mathrm{m}$ by concentric balls around the external face for the distance on the dual map $ \mathrm{m}^\dagger$. See Figure \ref{fig:layersoneface} and \cite[Chapter 13]{CurStFlour} or \cite[Section 6.5]{BBCK18} for a precise definition. When applied to a fixed map $ \mathrm{m}$, it yields a labeled plane tree  $ \mathrm{Peel}_{ \mathcal{A}_{ \mathrm{dual}}}( \mathrm{m})$ as in the preceding section. As in \cite{BCK18} we shall take snapshots of the process when the balls of radius $r=0,1,2,\dots$ around the external face in $  \mathrm{m}^\dagger$ have been discovered. Specifically, as in \cite[Section 6.5]{BBCK18} we consider the balls $ \mathrm{Ball}_{r}^{\dagger}( \mathrm{m})$ seen as submaps of $ \mathrm{m}$. Those submaps may have several holes that are associated with some vertices of the peeling tree $ \mathrm{Peel}_{ \mathcal{A}_{ \mathrm{dual}}}( \mathrm{m})$ that correspond to completion times of a layer are that we call the \textbf{cactus} vertices. We shall also declare that all vertices labeled $0$ (corresponding to the discovery of a vertex in the map $ \mathrm{m}$) are cactus vertices. See Figures \ref{fig:layersoneface} and \ref{fig:labeleddual}.  By contracting all other vertices and forgetting the planar structure, we obtain a new labeled tree that we denote by 
$$\mathrm{Cactus}( \mathrm{m}^{\dagger}).$$
The terminology ``cactus'' comes from the fact that the height of a vertex in $\mathrm{Cactus}( \mathrm{m}^{\dagger})$ correspond to the distance, in the dual map, of the corresponding boundary to the external face minus $1$. See Figure \ref{fig:cactuses} and Example \ref{ex:3/2stable} for the continuum analog.

Since $\mathrm{Cactus}( \mathrm{m}^{\dagger})$ is a discrete ``subordination'' of the peeling tree $ \mathrm{Peel}_{ \mathcal{A}_{ \mathrm{dual}}}( \mathrm{m})$ for the algorithm $ \mathcal{A}_{ \mathrm{dual}}$, it is easy to check that under $\mathbb{W}_{ \mathbf{q}}^{(n)}$, it is again a multi-type Galton--Watson tree, and after a shift of all labels by one, its law is denoted by $( \Q^{\GW, \mathrm{cac}}_{n})_{n \geq 1}$. Describing precisely the kernels of $( \Q^{\GW, \mathrm{cac}}_{n})_{n \geq}$ is combinatorially  too demanding, and we shall not follow this path.

 Our goal is now to describe the scaling limits of the decorated tree $(T_{\GW},g, \uplambda^{ \mathbf{1}_{ \{1\}}})$ under $ \Q^{\GW, \mathrm{cac}}_{n}$.   It is known from \cite{BC16}  or \cite[Theorem 6.8]{BBCK18} that the ``typical height'' of these trees under $ \Q^{\GW, \mathrm{cac}}_{n}$ is of order $n^{\beta-1}$ when $\beta  \in (1, 3/2]$ and of logarithmic order when $\beta  \in (1/2, 1)$, see \cite{kammerer2024scaling}. In the following, we will then restrict to  the case $\beta >1$. To state the theorem, we need yet another scaling constant, see \cite[Eq. (11.1)]{CurStFlour} (this constant is  denoted by $ \mathrm{a}_{ \mathbf{q}}$ in the notation of  \cite[Section 6.5]{BBCK18})
 
  $$  1 + \mathrm{g}_{ \mathbf{q}} := \frac{1}{2}\left( 1 + \sum_{k\geq 0}^{\infty} (2k+1) q_{k+1} {c}_{ \mathbf{q}}^{-k}\right).$$

 We show convergence of the rescaled trees towards the same ssMt as in the previous theorem but with self-similarity exponent $\alpha = \beta-1$ instead of $\beta$.

\begin{corollary}[Scaling limits for dual cactus tree]  \label{cor:cactus} Fix the peeling algorithm $ \mathcal{A}_{ \mathrm{dual}}$ and let $ \mathbf{q}$ be a weight sequence of type $\beta \in (1;3/2]$. Considering the integer-type Galton--Watson tree of law $ \Q^{\GW, \mathrm{cac}}_{n}$ described above, the sequence of the distributions  of the rescaled measured decorated trees 

$$\left( T_\GW, \frac{d_{T_\GW}}{n^{\beta-1}} \cdot \frac{ {c}_{ \mathbf{q}}}{ 1 + \mathrm{g}_{ \mathbf{q}}}, \rho_\GW, \frac{g}{n}, \frac{ \uplambda^{ \mathbf{1}_{\{1\}}}_\GW}{n^{\beta+1/2}} \cdot \frac{1}{\mathrm{b}^{ \mathrm{vertices}}_{ \mathbf{q}}}\right) \mbox{ under } \Q^{\GW, \mathrm{cac}}_{n}$$
converges as $n\to \infty$,  to the law  of 
the measured self-similar Markov tree $(\normalfont{\texttt{T}},\upmu)$  under $\Q_1$ with the characteristic quadruplet $(0,  \mathrm{a}_{ \beta-1/2,\beta}, \boldsymbol{\Lambda}_{ \beta-1/2,\beta}; \beta-1)$  of  Example \ref{ex:stablefamily} when $ \texttt{b}=1/2, \texttt{a} = \beta-1/2$ so that $\frac{1}{2} \leq \beta = \texttt{a} + \texttt{b} \leq \frac{3}{2}$ and $\varrho$ such that $\beta(1-\varrho)= 1/2$ (i.e. there is no killing term), see  \eqref{eq:nokiling} and \eqref{eq:kappanokilling}.
\end{corollary}

\begin{remark} By using Theorem \ref{T:mainunconds-length} instead of Theorem \ref{T:mainunconds-mass} (and using the same proof minus the verification of Assumption \ref{A:hGWrv}) we obtain the convergence of the renormalized length measures $$\frac{ \uplambda^{ \varpi}_\GW}{ n^{\beta-1} \varpi(n)}$$ towards its continuous analog for any regular varying weight function with index $\gamma >1$. In particular, this implies the convergence of the sum of the $\gamma$-powers of the length of the components of the boundaries of dual balls of radius $r$ centered at the external face for $r \geq 0$, which is a new result.
\end{remark}

\begin{proof}  We will again apply our Theorem \ref{T:mainunconds-mass}, but contrary to the preceding cases, we will build upon Theorem \ref{prop:GF} to check,\textit{ using probabilist arguments} the Assumptions   \ref{assum:main:gen:deco:repro} and \ref{A:all:existance:PHI} (as well as \ref{A:hGWrv}) instead of the analytic ones Assumptions \ref{A:all:GW}, \ref{A:BK} and  \ref{A:typex}. 

\noindent $\bullet$ \textsc{Assumption \ref{A:hGWrv}}: Since the number of vertices of label $1$ is again the number of vertices of the map $ \mathrm{m}$, and since those vertices are kept in the transformation going from the peeling tree to the cactus tree, Assumption \ref{A:hGWrv} (as well as \ref{A:all:GW}) is again ensured by \cite[Proposition 10.4]{CurStFlour}.

\noindent $\bullet$ \textsc{Assumption \ref{A:all:existance:PHI}}: We will now construct superharmonic functions to check Assumption \ref{A:all:existance:PHI}. Fix $ \beta < \gamma < 2 \beta+1$. From the convergence of the discrete cumulants $ \kappa^{(n)}$ established in the proof of Theorem \ref{prop:GF}, by Proposition \ref{prop:APHI}, we can find a function $\phi:\mathbb{N}\to \mathbb{R}_+^*$ strictly superharmonic such that $\phi(n) \asymp n^{\gamma}$ and 
  \begin{eqnarray} \label{eq:diminish} \limsup \limits_{n\to \infty} n^{\alpha-\gamma}\boldsymbol{G}^{ \mathrm{peel}}_{\GW} \phi(n)<0,  \end{eqnarray}
these properties are obviously referring to the multi-type Galton--Watson tree of laws $( \Q_{n}^{\GW, \mathrm{peel}})_{n \geq 1}$ associated with the peeling trees. Let us prove that the same function $\phi$ has the desired properties but for the cactus tree of laws $( \Q_{n}^{\GW, \mathrm{cac}})_{n \geq 1}$.
Since the kernels of $(\Q_{n}^{\GW, \mathrm{cac}})$ are obtained by ``iterating '' a (random but)  finite number of times the kernels  of $( \Q_{n}^{\GW,\mathrm{peel}})_{n \geq 1}$ it is clear by Fatou's lemma that $\varphi$ remains a strictly super-harmonic function for $(\Q_{n}^{\GW, \mathrm{cac}})$. More precisely, suppose that we start from a face of perimeter $2n$, and peel around it using the algorithm $ \mathcal{A}_{ \mathrm{dual}}$ until we found the ball of dual distance $1$ from the external face after $K_{n}$ peeling steps. We  denote by $n = X^{{(n)}}_{1}, X^{{(n)}}_{2}, \dots , X^{{(n)}}_{K_{n}}$  the half-perimeters  of the holes peeled at each step in this process. If we write $\boldsymbol{G}_{\GW}^{ \mathrm{cac}}$ the generator (see \eqref{Eq:Laplacediscret}) for the multi-type tree of laws $( \Q_{n}^{\GW, \mathrm{cac}})_{n \geq 1}$ then we can write:
$$  \boldsymbol{G}_{\GW}^{ \mathrm{cac}}(\varphi)(n) =  \mathbb{E}\left[\sum_{i=1}^{K_n} \boldsymbol{G}_{\GW}(\varphi)(X^{(n)}_i) \right] \underset{ \eqref{eq:diminish}}{\leq} \mathbb{E}\left[\sum_{i=1}^{K_n} \mathbf{1}_{X^{(n)}_i > n/2}\right] (-  \mathrm{c}_{1} \cdot n^{\gamma-\beta }),$$ for some $ \mathrm{c}_{1}>0$ for all $n$ large enough. It follows from the discussion in \cite[Section 4.2]{BC16}, that as $n \to \infty$ peeling all the edges adjacent to a face of perimeter $2n$ necessitate roughly $ \approx \frac{1}{ 1+\mathrm{g}_{ \mathbf{q}}} n$ steps, and that most of those steps take place in an undiscovered hole of half-perimeter $ \approx 2n$. This implies that with our notation we have  $\mathbb{P}\left(\sum_{i=1}^{K_n} \mathbf{1}_{X_i^{(n)} > n/2} \geq \mathrm{c}_{2} \cdot n\right) \to 1$ as $n \to \infty$. Combining this with the preceding display yields that for $n$ large enough we have
$$ \boldsymbol{G}_{\GW}^{ \mathrm{cac}}(\varphi)(n) \leq -\mathrm{c}_{3} \cdot n^{\gamma-(\beta-1)},$$ for some $ \mathrm{c}_{3}>0$. Performing this operation for two values of $\gamma \in (\beta, 2 \beta+1)$ implies Assumption \ref{A:all:existance:PHI}.

\noindent \textsc{Assumption   \ref{assum:main:gen:deco:repro}:} The convergence of the evolution of the locally largest evolution under $( \Q^{\GW, \mathrm{cac}}_{n})_{n \geq 1}$ can be deduced, using a discrete Lamperti transformation, from that under $( \Q^{\GW,\mathrm{peel}}_{n})_{n \geq 1}$ which is controlled thanks to Theorem \ref{prop:GF}. This has been first performed in \cite[Proposition 6.6 item 2]{BBCK18}. More precisely, using the locally largest selection rule, let us consider the decoration-reproduction process $(f, \eta)$ under $ P_{n}^{\GW, \mathrm{cac}}$ and their rescaled versions obtained by \eqref{Eq:scalefeta} with $\alpha = \beta-1$. 
Since our non-trivial splittings under $ P_{n}^{\GW}$ are conservative, the reproduction process can be recovered from the decoration process. This is still asymptotically true in the case of $ P_{\GW, \mathrm{cac}}^{(n)}$ (see after \eqref{Eq:scalefeta} for the notation) and adapting \cite[Corollary 18]{BCK18} to our case (we leave the details to the reader here),  implies with the notation of Lemma \ref{lem:convreprodecolight}  that for any $ \varepsilon>0$   $$P_{\GW, \mathrm{cac}}^{(n), \varepsilon} \xrightarrow[n\to\infty]{(d)}P_{1}^{ \varepsilon},$$ where the convergence holds in the sense of the topology used in Assumption \ref{assum:main:gen:deco:repro}. One can then withdraw the cutoff at $ \varepsilon$ using Assumption \ref{A:all:existance:PHI} by Lemma \ref{lem:convreprodecolight}, and this proves Assumption \ref{assum:main:gen:deco:repro}. \end{proof}

 The previous result is a first step towards understanding the whole geometry of the dual maps $ \mathrm{m}^\dagger$ under $\mathbb{W}_{ \mathbf{q}}^{(n)}$. In particular it implies that  $$  \left\{n^{-(\beta-1)} \cdot \mathrm{Diam}( \mathrm{m}^\dagger) \mbox{  under }\mathbb{W}_{ \mathbf{q}}^{(n)}, n \geq 0 \right\} \mbox{is tight},$$ which is a new result. Compared to the primal distances see \cite{CRM22}, the dual distances are not encoding using Schaeffer-type constructions and are more difficult to control, see the recent work \cite{kammerer2024scaling} for scaling limits candidates whose construction implicitly relies on our ssMt.  As the Brownian growth-fragmentation of Example \ref{ex:3/2stable} has been directly identified within the Brownian disk in \cite{le2020growth}, it would be desirable to find the ssMt of Proposition \ref{prop:GF} (or their related analogs in the Corollary \ref{cor:cactus}) directly within the stable carpet/gasket.

\subsection{$O(n)$-decorated random  quadrangulations}
 \label{sec:O(n)}

In this section, we study the peeling process on random maps decorated with an $O( \mathrm{n})$ model. As we are going to see, these decorated maps can be seen as an ``overlay'' on the random maps with weight sequence of type $ \beta \in (1/2,3/2)$ in a similar fashion as the ssMt of Example \ref{ex:overlay} were an ``overlay'' of the ssMt of Example \ref{ex:stablefamily} without killing.

We shall focus on the case of quadrangulations decorated with loops, and we will employ a version of the peeling process introduced by Budd \cite{BudOn} to analyze this model. We refer to \cite{BBG11,BudOn,CCMOn} for details (and to \cite{korzhenkova2022exploration} for the triangular case). Let us recall the model: a rigid \textbf{loop-decorated} quadrangulation is a pair $(  \mathrm{m}_{\square}, \mathbf{l})$ where $ \mathrm{m}_{\square}$ is a quadrangulation, i.e.~a rooted planar map such that 
all  inners faces are  quadrangles, and  $\mathbf{l}:=\{l_1,\dots, l_k\}$ is a collection of disjoint (unoriented) loops on the dual map avoiding the root face. The term rigid comes from the fact that we further assume that loops can only cross  quadrangles through opposite sites, see Figure \ref{fig:takingthegasket} below. For simplicity, in the rest of the section, we will omit the mention to the term rigid.

The pair $ ( \mathrm{m}_{\square},  \mathbf{l})$  is called a loop-decorated quadrangulation, and we denote the set of all loop-decorated quadrangulation with perimeter $k$ by $\mathcal{Q}_{\mathrm{loop}}^{(k)}$. We can then define a natural measure on such configurations by putting 
$$ w_{x,y, \mathrm{n}}\big( (  \mathrm{m}_{\square}, \mathbf{l})\big) =  x^{|  \mathrm{m}_{\square}|} y^{| \mathbf{l}|}  \mathrm{n}^{\# \mathbf{l}},$$
where $x,y>0$ and $ \mathrm{n} \in (0,2)$ are parameters\footnote{It is customary to use the notation $ \mathrm{n}$ for the non-local weight for each loop, but the reader should not confuse it with the size parameter $n \to \infty$}, and further $|  \mathrm{m}_{\square}|$ is the number of faces of the quadrangulation that are not visited by a loop, $| \mathbf{l}|$  the total length of the loops and $\# \mathbf{l}$  the number of loops.  Provided that the measure $w_{x,y, \mathrm{n}}$ has finite total mass (the parameters are then called admissible) one can normalize it to define a probability measure $$\mathbb{W}^{(k)}_{ x,y, \mathrm{n}} = \frac{w_{x,y, \mathrm{n}}( \cdot \cap \mathcal{Q}^{{(k)}}_{\mathrm{loop}})}{W_{x,y, \mathrm{n}}^{(k)}} , \quad \mbox{ with } W_{x,y, \mathrm{n}}^{{(k)}} = w_{x,y, \mathrm{n}}( \mathcal{Q}^{{(k)}}_{\mathrm{loop}}), $$ on  loop-decorated random maps with half-perimeter $2k$, where we recall that the perimeter of the map corresponds to the degree of the external face. As in the preceding section, we put $W_{ x,y, \mathrm{n}}^{(0)}=1$ for convention as this represents the weight of the vertex map. These random decorated maps are very closely related to the Boltzmann maps studied in the preceding sections through their gaskets. The gasket $\mathrm{Gask}(\mathrm{m}_{\square}, \mathbf{l})$ is the bipartite rooted map obtained by pruning off the interiors of the outermost loops with respect to the root edge. Equivalently, the gasket $\mathrm{Gask}(\mathrm{m}_{\square}, \mathbf{l})$
consists of all the vertices that can be connected to the root edge by a path that does not cross a loop, as well as all the edges connecting such vertices; see Figure~\ref{fig:takingthegasket}). Recalling the setup of Section \ref{sec:peelingnormal}, it is easy to see that under $\mathbb{W}^{(k)}_{ x,y, \mathrm{n}}$, the map $\mathrm{Gask}(  \mathrm{m}_{\square}, \mathbf{l})$ is actually a $ \mathbf{q}$-Boltzmann map of law $  \mathbb{W}_{ \mathbf{q}}^{(k)}$ for the weight sequence $ \mathbf{q}\equiv \mathbf{q}(x,y, \mathrm{n})$ given by 
  \begin{eqnarray} \label{eq:qgasket}  {q}_{k} = x \mathbf{1}_{k=2} + y^{2k}\mathrm{n} W_{ x, y, \mathrm{n}}^{(k)},  \end{eqnarray}
 see \cite[Eq. (2)]{BudOn}, and with this choice of weight measure we even have the equalities $W_{ x,y, \mathrm{n}}^{(k)} = W_{ \mathbf{q}}^{(k)}$ for all $k \geq 0$, see \cite[Theorem 1]{BudOn}. This will allow us to apply results of Section \ref{sec:peelingnormal} to the study of
 loop-decorated quadrangulation.
  
\begin{figure}[!h]
 \begin{center}
 \includegraphics[height=5cm]{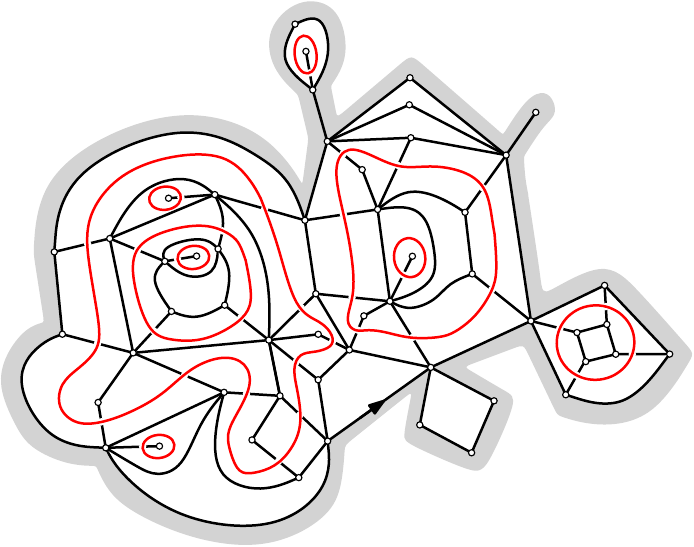} \hspace{0.5cm}
  \includegraphics[height=3cm]{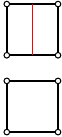}\hspace{0.5cm}
   \includegraphics[height=5cm]{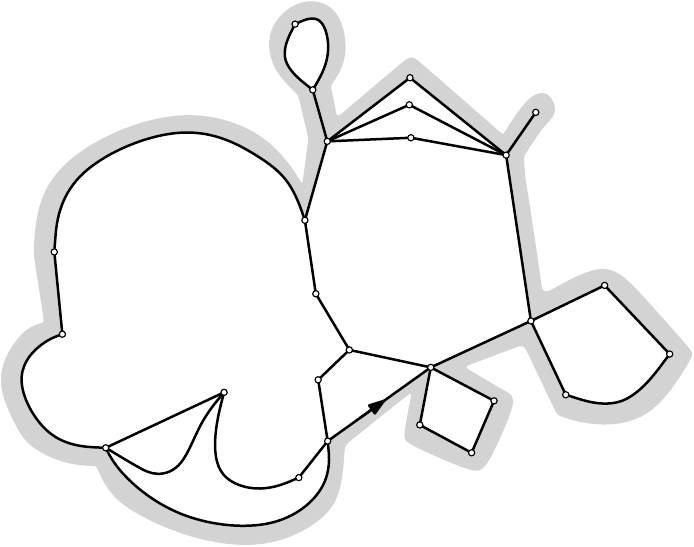}
 \caption{ A loop-decorated quadrangulation (with a boundary)  built out of the two types of faces depicted in the middle. On the right, the gasket of the map. \label{fig:takingthegasket}}
 \end{center}
 \end{figure}
 
Let us describe the phase diagram of the loop-decorated map $(  \mathrm{m}_{\square}, \mathbf{l})$ obtained in \cite{BBG11,BudOn} through the lens of its gasket. This will allow to identify the different regimes of interest for the parameters $(x,y, \mathrm{n})$ in order to obtain non-degenerated scaling limits. Fix $ \mathrm{n} \in (0,2)$\footnote{The case $n =2$ is a boundary case, see \cite{kammerer2024gaskets} for recents results}.  For most of the parameters $(x,y)$, the model is either non-admissible or the weight sequence $ \mathbf{q}(x,y, \mathrm{n})$ is of type $ \beta=\frac{1}{2}$, that is the gasket is subcritical.  It is only when $ x \leq \frac{1}{12}$ that there exists a unique value $y(x)$ so that $ \mathbf{q}(x,y, \mathrm{n})$ is of type $\ne 1/2$: The line $(x,y(x))$ is further split in two parts, the right one where the weight sequence $ \mathbf{q}(x,y, \mathrm{n})$ is  of type $ 3/2$, and the left part where it is of type 
\begin{equation}\label{beta:o:n:minus}
 \beta = 1 - \frac{1}{\pi} \arccos( \mathrm{n}/2).
 \end{equation}
At the intersection of these two curves sits a single point $(x^{*}, y^{*})$ where the gasket is a Boltzmann map with weight sequence $ \mathbf{q}(x,y, \mathrm{n})$  of type  
\begin{equation}\label{beta:o:n:plus} \beta = 1 + \frac{1}{\pi} \arccos( \mathrm{n}/2).
\end{equation}
We refer to \cite{BBG11,BudOn} for details. The case $\beta = 1 - \frac{1}{\pi} \arccos(\mathrm{n}/2) \in ( \frac{1}{2};1)$ lies in the so-called dense phase whereas $\beta = 1 + \frac{1}{\pi} \arccos(\mathrm{n}/2) \in (1; \frac{3}{2})$ belongs to the dilute phase. In particular, thanks to \eqref{eq:asymptWell}, we have the following asymptotics
\begin{equation}\label{eq:asym:W:q:h:n}
W_{x,y, \mathrm{n}}^{(k)}=W_{ \mathbf{q}}^{(k)} \underset{k\to \infty}{\sim}   \frac{\mathrm{p}_{ \mathbf{q}}}{2} (c_{ \mathbf{q}})^{k+1} k^{-\beta-1}, 
\end{equation}
 and due to \cite[Lemma 5.2 and Proposition 10.1]{CurStFlour} we in fact have $c_{ \mathbf{q}}= y^{-1}$. In the rest of this section, we shall focus on the non-generic case and suppose that $\beta \in (1/2,3/2) \backslash \{1\}$.

\begin{figure}[!h]
 \begin{center}
 \includegraphics[height=6cm]{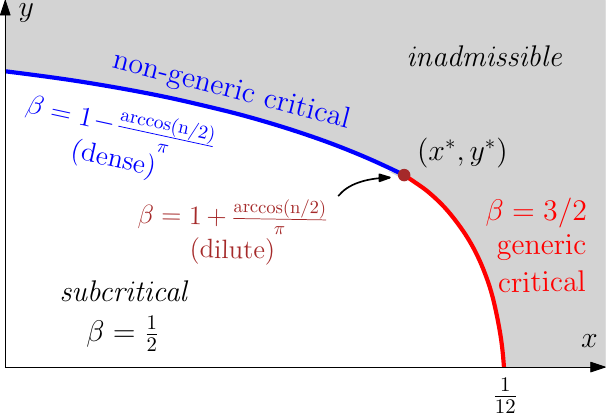}

 \caption{The phase diagram in $(x,y)$ for a fixed $ \mathrm{n}$.}
 \end{center}
 \end{figure}
\begin{figure}[!h]
  \begin{center}
  \includegraphics[width=15.5cm]{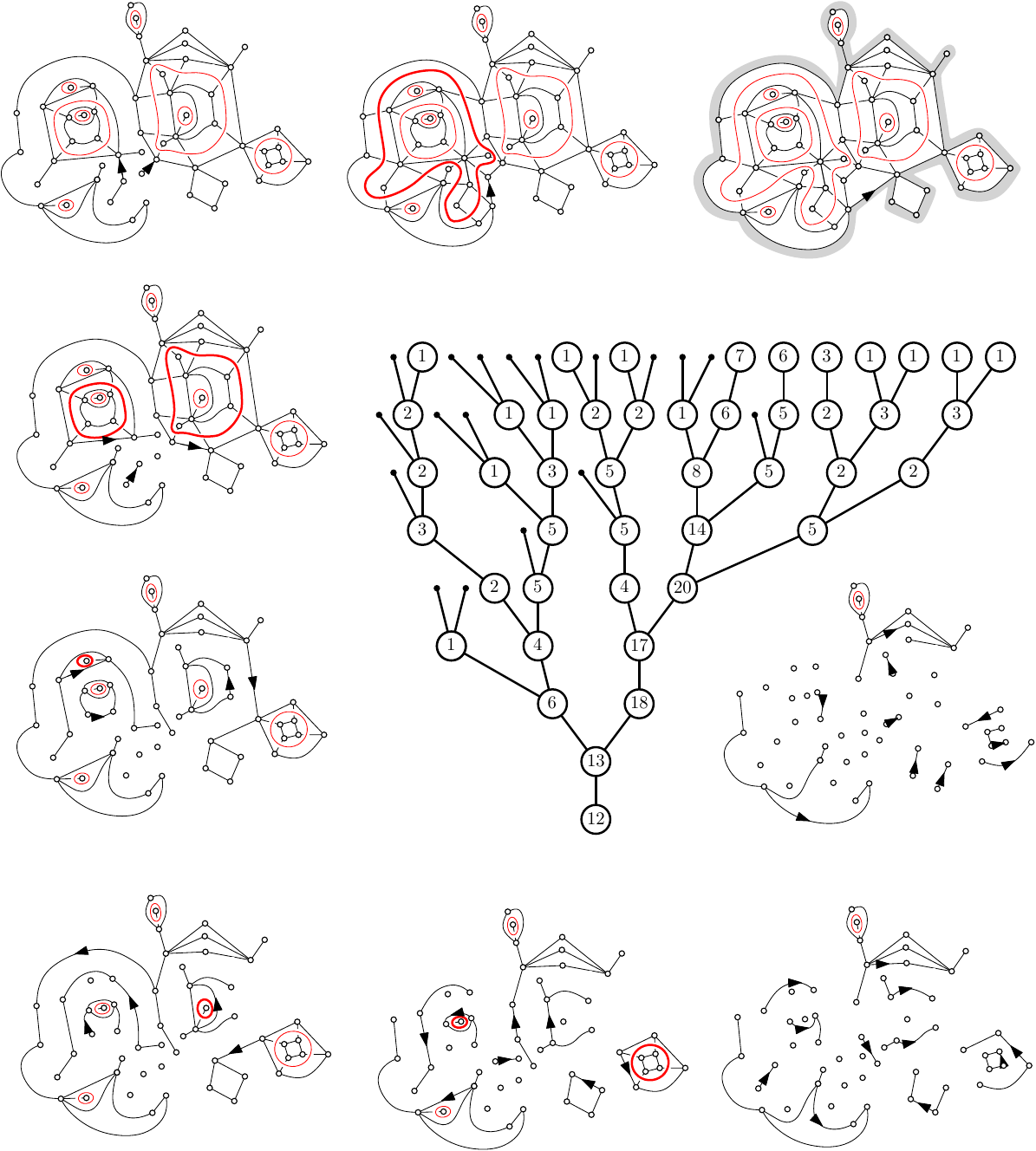}
  \caption{Illustration of the set-by-step (unfinished) peeling of a loop decorated quadrangulation. For better visibility, the leaves are represented as black dot instead of vertices carrying a label $0$.\label{fig:peelingOn}}
  \end{center}
  \end{figure}

 Following Budd \cite{BudOn}, one can adapt the peeling process to loop-decorated quadrangulations. Given an algorithm $ \mathcal{A}$ this yields a random plane tree $ \mathrm{Peel}_{ \mathcal{A}}(  \mathrm{m}_{\square}, \mathbf{l})$ whose vertices are labeled by integers $\{0,1,2, \dots \}$ corresponding to the half-perimeters of the boundaries encountered in the exploration. The important change  is that when we discover a new quadrangle containing a part of a loop, we immediately reveal the whole loop and the quadrangles it traverses, see Figure~\ref{fig:peelingOn}.
 
Namely,  suppose that we have a loop decorated quadrangulation with perimeter $2k$ and we observe that we must be in one of the mutually exclusive cases when erasing the root edge:
\begin{itemize} 
\item case 1) If the face  lying on the left of the rooted edge  is a quadrangle which is not visited by a loop, then the erasure of the root edge leaves the map connected but increases the degree of the external face by $2$.
\item case 2) If the face  lying on the left of the rooted edge  is a quadrangle which is visited by a loop $\ell$, then we erase the root edge and all the edges traversed by $\ell$. This operation splits the map into two connected components with perimeters $2k+2l-2$ and $2l$, where $2l$ is the length of the loop $\ell$.
\item case 3) If the root face is incident both on the left and right of the root edge, then the erasure of the root edge splits the map into two connected components with perimeters, say, $2k_{1}$ and $2k_{2}$ with $2k_{1}+2k_{2}=2k-2$.
\end{itemize}
Next, we root each connected component at an oriented edge of their boundary, following a specific algorithm $\mathcal{A}$, and continue the exploration recursively until there are no more edges. For a more detailed presentation, we refer to \cite{BudOn}.  As in Section \ref{sec:peelingnormal}, we can encode the previous exploration in a decorated  labeled plane tree   $$   \mathrm{Peel}_{ \mathcal{A}}(  \mathrm{m}_{\square} , \mathbf{l})$$
\clearpage
 where inner-vertices of the tree correspond to boundaries and the type corresponds to the half-perimeter. Again, the leaves, carrying the label $0$, correspond to the actual vertices of the underlying map. Since a picture is worth a thousand words we encourage the reader to spend some time  looking at  Figure \ref{fig:peelingOn}.

We then have the analog of Proposition \ref{Markov:peeling}:
\begin{proposition}  \label{Markov:O(n)}Fix  admissible parameters $( x,y, \mathrm{n})$ and denote by $ \mathbf{q}= \mathbf{q}(x,y, \mathrm{n})$ the weight sequence given in \eqref{eq:qgasket}. Fix a peeling algorithm $ \mathcal{A}$ for loop-decorated rigid quadrangulations. Then the law of $\mathrm{Peel}_{ \mathcal{A}}( \mathrm{m}_{\square}, \mathbf{l})$ under $ \mathbb{W}^{(k)}_{x,y, \mathrm{n}}$, with all labels shifted by one, is an integer-type Galton--Watson process whose laws $ (\Q^{\GW, \mathrm{loop}}_{n})_{n \geq 1}$ are characterized by the only non trivial transitions:
$$ \begin{array}{ll} \mathbb{Q}^{\GW, \mathrm{loop}}_{n+1}( \mathbf{Z}(1)=\delta_{(n+1+2)}) = \displaystyle \frac{x W^{(n+1)}_{ \mathbf{q}}}{W^{(n)}_{ \mathbf{q}}}  & \mbox{ for } n \geq 1,\\
\mathbb{Q}^{\GW, \mathrm{loop}}_{n+1}( \mathbf{Z}(1)=\delta_{(n+1+k)} + \delta_{k+2}) = \displaystyle \frac{ y^{2(k+1)}\mathrm{n} W_{ \mathbf{q}}^{(k+1)} W^{(n+k)}_{ \mathbf{q}}}{W^{(n)}_{ \mathbf{q}}}  & \mbox{ for } n \geq 1, k \geq 0\\
\mathbb{Q}^{\GW, \mathrm{loop}}_{n+1}( \mathbf{Z}(1)=\delta_{n_{1}+1}+\delta_{n_2+1}) = (2- \mathbf{1}_{n_{1}=n_2}) \displaystyle \frac{ W^{(n_{1})}_{ \mathbf{q}} W_{ \mathbf{q}}^{{(n_2)}}}{W^{(n)}_{ \mathbf{q}}},  & \mbox{ for } \begin{array}{c}n \geq 1, \ \ n_{1},n_{2}\geq 0\\\mbox{ s.t. } n_{1}+n_{2} +1 =n.  \end{array}\\
 \mathbb{Q}^{\GW, \mathrm{loop}}_{1}( \mathbf{Z}(1)=  \mathbf{0})=1.& \end{array}$$
   \end{proposition}
   \begin{proof} The proof is mutatis-mutandis the same as for Proposition \ref{Markov:peeling}: This follows easily considering the probability of the three cases enumerated above and the Markov property established in \cite[Lemma 1]{BudOn}. \end{proof}
   
   Equipping the random tree $T_\GW$ with the measure $\uplambda^{ \mathbf{1}_{\{1\}}}_\GW$ we can now state the analog of Theorem \ref{prop:GF} in this case:
   
   \begin{corollary}[Scaling limits for peeling trees of loop-decorated $ O( \mathrm{n})$ quadrangulations] 
   
   \label{prop:peelingOn} Fix  $(x,y, \mathrm{n})$  parameters with $ \mathrm{n}\in (0,2)$ so that the associated weight sequence $ \mathbf{q}$ by \eqref{eq:qgasket} is of type $\beta \in (1/2,3/2) \backslash \{1\}$ and fix a peeling algorithm $ \mathcal{A}$ for loop decorated quadrangulations.  Let us write  $\mathbb{Q}_{1}^{( \beta), \mathrm{over}}$ for the law of the ssMt of Example \ref{ex:overlay} with the corresponding $\beta$ and self-similarity parameter $\alpha= \beta$. 
  Considering the integer-type Galton--Watson tree of law $ \Q^{\GW, \mathrm{loop}}_{n}$ described in Proposition \ref{Markov:O(n)}, there exists $\mathrm{v}_{x,y, \mathrm{n}}>0$ so that the sequence of the distributions  of the rescaled measured decorated trees 

$$\left( T_\GW, \frac{d_{T_\GW}}{n^{\beta}} \cdot \mathrm{d}_{ \mathbf{q}}, \rho_\GW, \frac{g}{n}, \frac{ \uplambda^{ \mathbf{1}_{\{1\}}}_\GW}{n^{\min( 2\beta,2)}} \cdot \mathrm{v}_{x,y, \mathrm{n}}\right) \mbox{ under } \Q^{\GW, \mathrm{loop}}_{n}$$
converges as $n\to \infty$,  to the law  of 
the measured self-similar Markov tree $(\normalfont{\texttt{T}},\upmu)$  under $\Q_1^{(\beta), \mathrm{over}}$.
\end{corollary}
\begin{remark}[Critical case] The boundary case $n=2$ yielding to $\beta=1$ is excluded here, see \cite{kammerer2024gaskets} and \cite[Theorem A, item b]{aidekon2024scaling} for recent enumerative results. In fact, the ssMt of Example~\ref{ex:overlay} does not satisfy the Assumption \ref{A:gamma0} in this case since $\kappa \geq 0$, so that its existence is not a priori guarantee by Proposition \ref{P:constructionomega-}. See Example \ref{ex:ADS} for a discussion of this case.
\end{remark}
\begin{proof} The proof is again an application of our Theorem \ref{T:mainunconds-mass} but many of the verifications have been done in Theorem \ref{prop:GF}.  Indeed, if we use the locally largest rule as selection rule then the decoration  $f$ of the decoration process $(f, \eta)$ under $ \P^\GW_n$ is the same as for the labeled tree corresponding to the peeling restricted to the gasket. The reproduction is changed as follows: each time the decoration process makes a positive jump of $k \geq 0$, the reproduction process gains an atom of size $k+1$. With this observation at hands, as in the proof of Theorem \ref{prop:GF}, let us check Assumptions \ref{A:all:GW}, \ref{A:BK} and  \ref{A:typex} as well as Assumption \ref{A:hGWrv}. 

\noindent $\bullet$ \textsc{Assumptions \ref{A:hGWrv} and \ref{A:all:GW}:}  Perhaps the most challenging assumptions to check are luckily for us proved in the literature. By  \cite[Proposition 9]{BudOn}\footnote{Notice the typo there where $\max$ should be a $\min$.} we have that
$$\E_n^{\GW, \mathrm{loop}}\big(\uplambda_{\GW}^{ \mathbf{1}_{\{1\}}}(T_\GW)\big) \underset{n\to \infty}{\sim}  \mathrm{v}_{x,y, \mathrm{n}} \cdot  n^{\min( 2\beta,2)},$$
and the uniform integrability of  $n^{-\min( 2\beta,2)}\cdot \uplambda_{\GW}^{ \mathbf{1}_{\{1\}}}(T_\GW)$, under $\Q_n^{\GW, \mathrm{loop}}$ has very recently been proved in \cite[Theorem A]{aidekon2024scaling}. In particular, this entails Assumption  \ref{A:hGWrv}. Since type $1$ is clearly accessible, the subcriticality is implied by the above display and so Assumption \ref{A:all:GW} clearly holds.

\noindent $\bullet$ \textsc{Assumption \ref{A:BK}:} The convergence of the drift and variance term is implied by the calculations of the proof of Theorem \ref{prop:GF} since the decoration process under $ \Q^{\GW, \mathrm{loop}}_n$ is identical to the decoration under $\Q^{\GW, \mathrm{peel}}_n$ with the proper choice of $ \mathbf{q}$. It remains to check the weak convergence of the generalized L\'evy measure $ \boldsymbol{ \Lambda}^{(n)}$ taking into account the reproduction event in case of positive jumps of the decoration process. This is plain since  for every continuous $F$ with compact support and vanishing in the neighborhood of $(1, 0, 0,\dots)$ it suffices to modify the first line of \eqref{eq:lambdanpeel} into 

 \begin{eqnarray*} n^\beta \sum_{k \geq 0}  \frac{h(n+k)}{h(n)} \nu(k)  F\Big(\frac{n+k}{n}, \frac{k+2}{n},\cdots\Big)
  \end{eqnarray*}
  which by \eqref{eq:tailnu} tends as desired to 
 $$\mathrm{p}_{\mathbf{q}} \cos\big((\beta+1)\pi\big) \int_{0}^\infty \frac{\mathrm{d}x}{(x(1+x))^{\beta+1}}F( 1+x,x, 0,\cdots).$$
 
    \noindent $\bullet$ \textsc{Assumption \ref{A:typex}:} As in the proof of Theorem \ref{prop:GF} one proves that the discrete cumulants $\kappa^{(n)}(\gamma)$ converge to the associated continuous cumulant $\kappa(\gamma)$ of \eqref{ex:overlay}  for every $\gamma\in (\beta, 2\beta+1)$. Since we have verified Assumption \ref{A:BK} above, one can use Lemma \ref{lem:convkappacond} whose hypotheses for $\gamma \in (\beta, 2\beta+1)$ are easily verified in our case using \eqref{eq:tailnu} again. \end{proof}

\section{Parking on (random) trees} \label{sec:parking}

In this section, we show that the conservative ssMt from Example \ref{ex:stablefamily} also arise as scaling limits of fully parked components in the model of parking on (random) trees. Since this model is more recent than that of random planar maps, our results are somewhat less complete than those of the previous section—particularly concerning the scaling limit of the counting measure. We begin by recalling the model of parking on rooted trees, as recently popularized by Lackner and Panholzer \cite{LaP16}.

\newcommand{\car}{{\includegraphics[width=0.4cm]{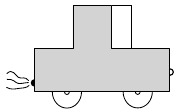}}}
\newcommand{\smallcar}{{\includegraphics[width=0.25cm]{images/car}}}

We consider a finite plane tree $ \mathrm{t}$ whose vertices will be interpreted as free parking spots, each spot accommodating at most 1 car, together with a configuration $\car_u \geq 0$ for $u \in \mathrm{t}$ representing the number of cars arriving on each vertex. Each car tries to park on its arrival vertex, and if the spot is occupied, it travels downward towards the root of the tree until it finds an empty vertex to park. If there is no such vertex on its way, the car exits the tree through the root $ \varnothing$. An Abelian  property of the model shows that the final configuration and the flux of cars going through any edge of the tree does not depend upon the order chosen to park the cars.

\begin{figure}[!h]
 \begin{center}
 \includegraphics[height=7cm]{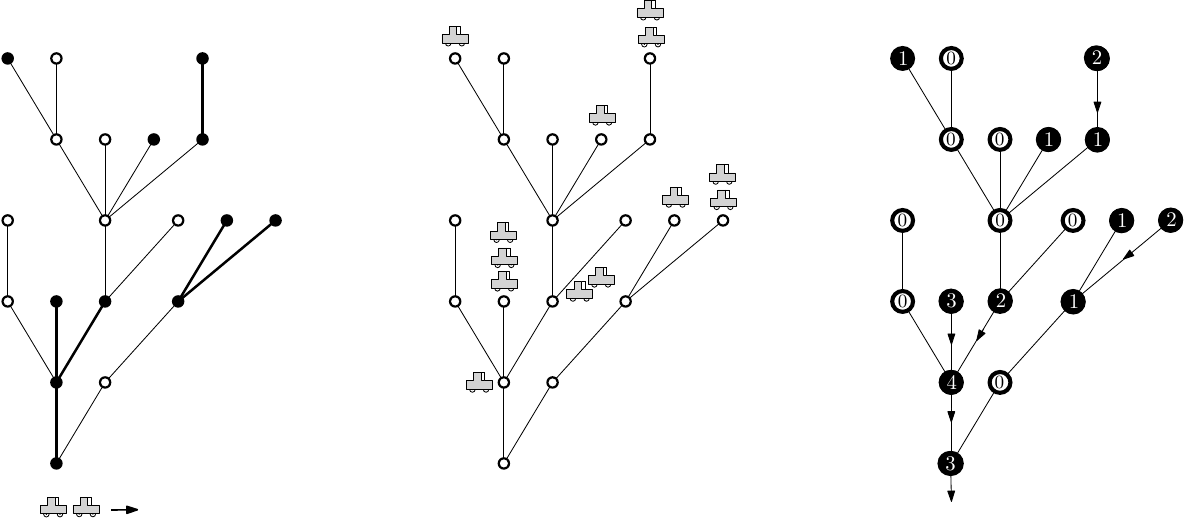}
 \caption{Illustration of the parking process  on a rooted tree: on the middle the initial configuration,  on the left,  the final configuration with the $2$ outgoing cars, and on the right, the same configuration where the vertices are labeled by the number of cars that visited each vertex. The fully parked components are in thick black lines on the left. \label{fig:parkinglabeled}}
 \end{center}
 \end{figure}

\begin{figure}[!h]
 \begin{center}
 \includegraphics[height=3.2cm]{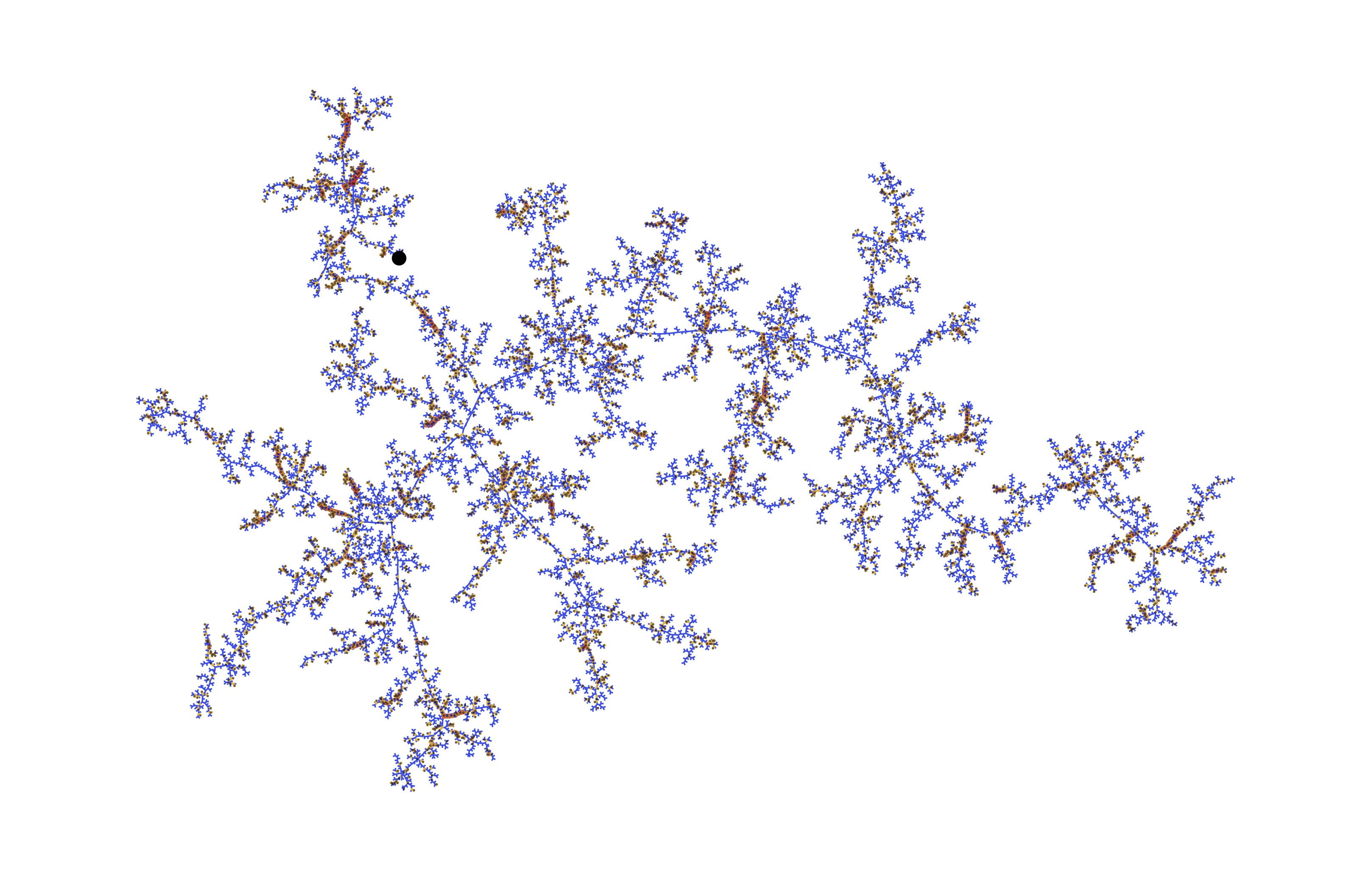}
  \includegraphics[height=3.2cm]{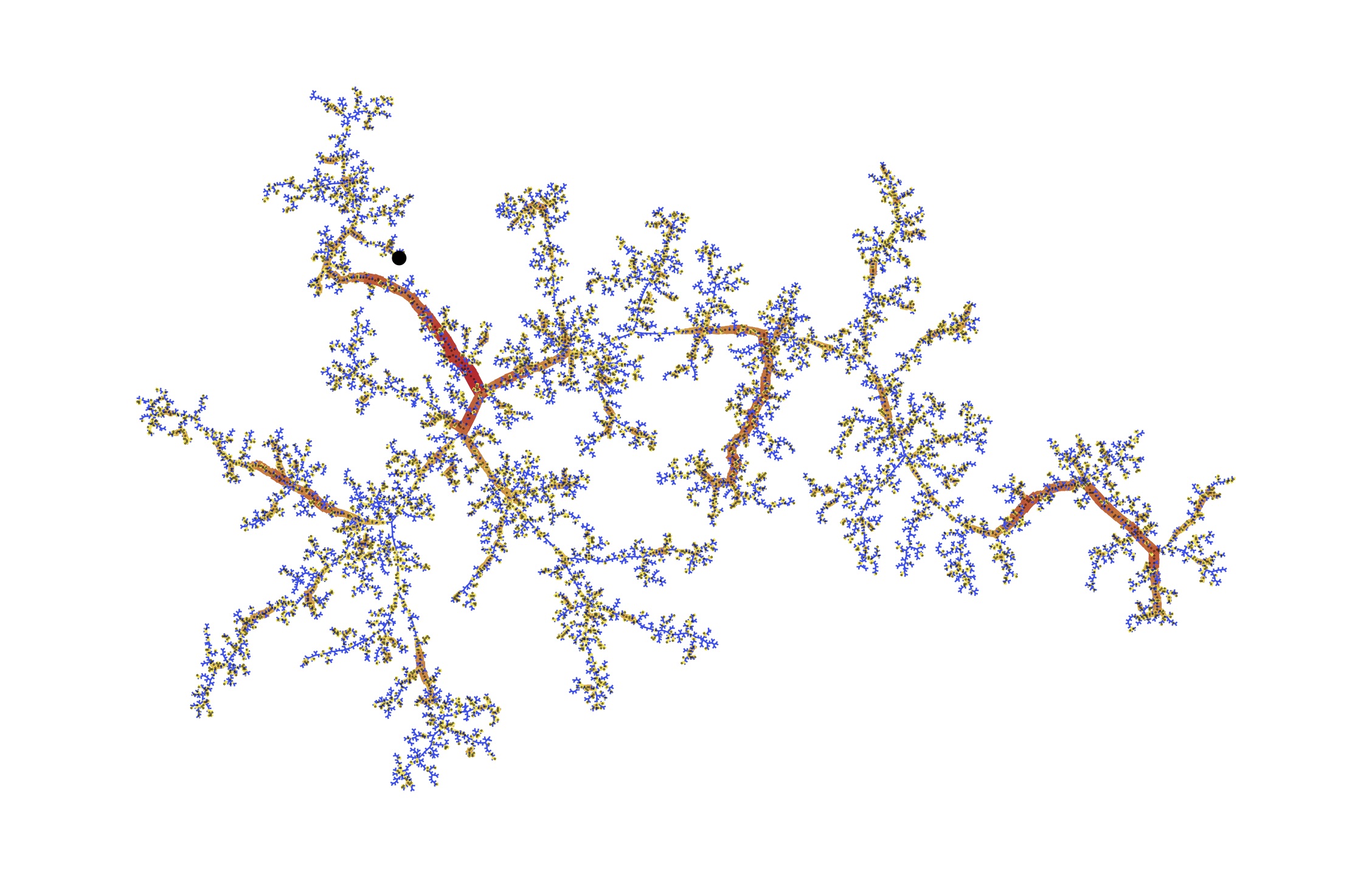}
   \includegraphics[height=3.2cm]{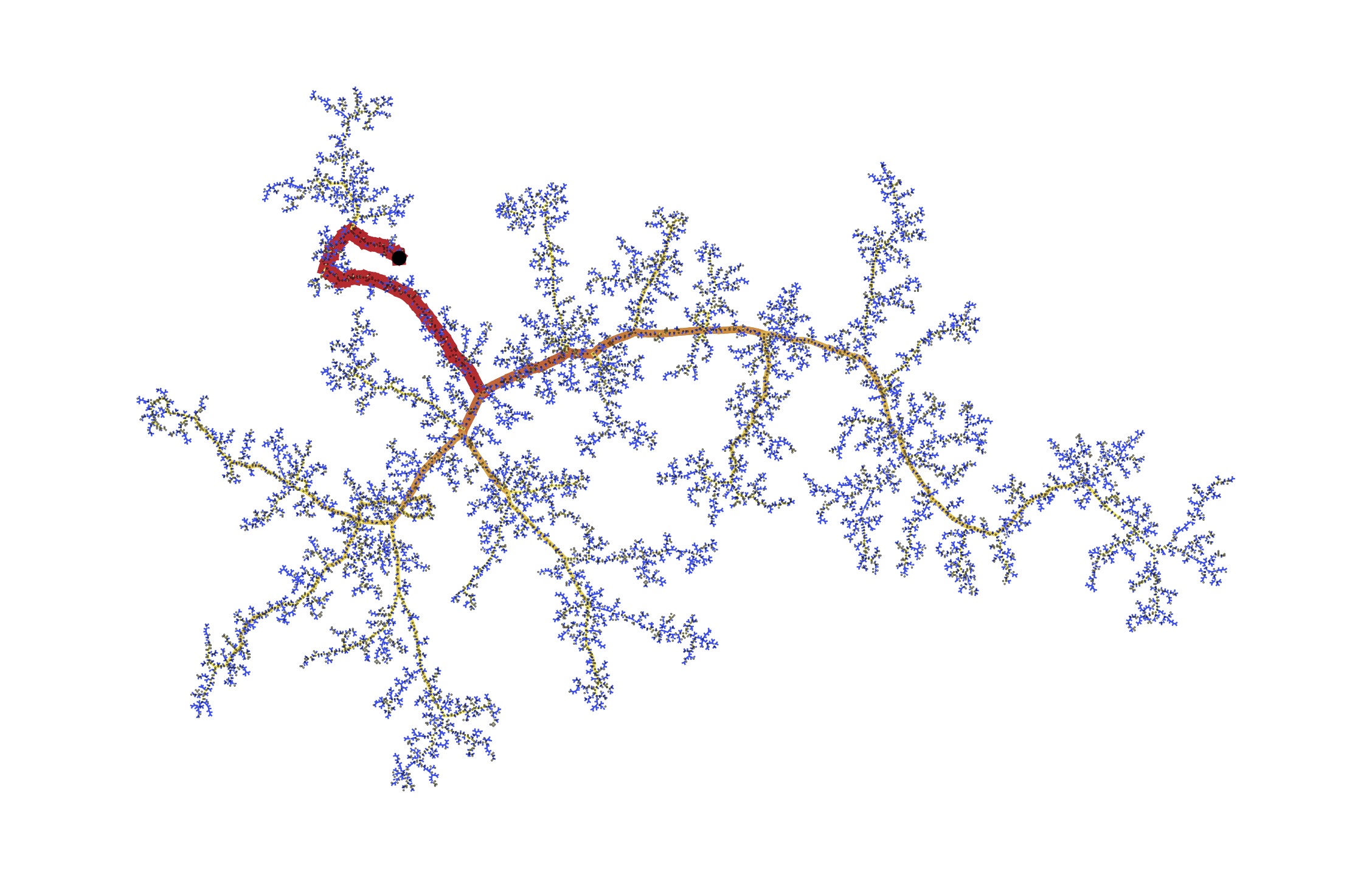}
 \caption{Simulation of  the parking  process on a random large uniform plane tree left:subcritical, middle:critical,  right:supercritical. The color (from blue to red) and the width (from small to large) of the edges represent the flux of cars. The root of the tree is indicated by a black dot.\label{fig:parkingsimu}}
 \end{center}
 \end{figure}

Let us now suppose that $ \mathrm{t}= T$ is a (monotype) Galton--Watson tree with offspring distribution $\xi$  and that conditionally on $ T$, the car arrivals $\car_u$ are i.i.d.~distributed according to some law $\mu$. After performing the parking of the cars, each vertex $u$ of the tree $T$ will be labeled by the number of cars $ \mathrm{Cars}(u)$ having visited this vertex in the parking (recall that by the Abelian property of the model,  this quantity does not depend upon the order in which we try to park the cars), see Figure \ref{fig:parkinglabeled}. We denote the obtained labeled tree by 
$$ \mathrm{Park}\big(T,(\car_u)_{u \in T}\big).$$ Note that a particle labeled $p\geq 0$ cannot give rise to a particle labeled $q \geq p+1$. More precisely, if a vertex visited by $p \geq 0$ cars on which $k\geq 0$ cars arrive gives rise to $i\geq 0$ vertices visited by $p_{1}, \dots , p_{i}$ cars each, then we must have 
  \begin{eqnarray} \label{eq:parking} p = k +\sum_{\ell =1}^{i} ( p_{\ell}-1)_{+},  \end{eqnarray}   where we recall that $(x)_{+} = \max(x,0)$.  We write  $\Omega^{(p)}$ for the probability (under the monotype Galton--Watson measure and the iid car arrivals) that the number of cars $ \mathrm{Cars}(\rho)$ visiting the root vertex $\rho$ of $T$ is equal to $p
\geq 0$. Then the branching property of the monotype Galton--Watson tree shows that conditionally on $ \mathrm{Cars}(\rho)=p$, if $p_{1}, ... , p_{i}$ are as above, the probability that the ancestor labeled $p$ gives rise to $i\geq 0$ subtrees with root labels $p_{1}, ... , p_{i} \geq 0$ enumerated from left-to-right is 

$$ \xi_{i} \cdot  \mu_{p- \sum_{\ell=1}^{i} (p_{\ell}-1)_{+}}  \cdot \frac{\prod_{\ell=1}^{i} \Omega^{(p_{\ell})}}{\Omega^{(p)}}.$$ We deduce from the last display that after forgetting the plane structure, we get a multi-type Galton--Watson tree, more precisely:

\begin{proposition}[Parked trees are integer-type Galton--Watson]  \label{prop:multiGWparking} For $p \geq 0$, denote by $\Omega^{(p)}$ the probability that the number of cars $ \mathrm{Cars}(\rho)$ visiting the root vertex $\rho$ of $T$ is equal to $p
\geq 0$. Then under $  \mathbb{P}( \cdot \mid  \mathrm{Cars}(\rho) = p)$ the random labeled tree $\mathrm{Park}(T,(\car_u)_{u \in T})$ is an integer-type Galton--Watson tree (where the type $0$ is allowed) whose laws $(\Q^{\GW, \mathrm{park}}_n)_{n \geq 0}$ are described by the following property:
 \begin{eqnarray}   && \mathbb{Q}^{\GW, \mathrm{park}}_{n}\big(  \mathbf{Z}(1) = (i_{0}, i_{1}, \dots , i_{n}, i_{n+1},0,0, \dots)\big)\nonumber \\ &=& \xi_{i_{0} + \dots + i_{n+1}} \cdot   \mu_{n- \sum_{\ell=0}^{n+1} i_\ell\
 (\ell-1)_{+}} \cdot  \frac{(i_{0} + \dots + i_{n+1})!}{i_{0}! i_{1}! \dots i_{n+1}!}  \frac{\prod_{\ell=0}^{n+1} \big(\Omega^{(\ell)}\big)^{i_{\ell}}}{{\Omega^{(n)}}}.  \label{eq:lawGWpark}\end{eqnarray}
\end{proposition}
Exactly as in the previous section, we now require a fine-tuning so that the above multi-type Galton--Watson tree converge to some interesting ssMt in the scaling limit. Specifically, it has been shown in recent works \cite{ConCurParking,CH19}  that if $\xi$ is critical (mean $1$), has finite variance $\Sigma_{\xi}^{2} \in (0, \infty)$ and if the car arrival distribution $\mu$ has mean $m_{\mu}>0$ and variance $\sigma_{\mu}^{2} \in (0,\infty)$ then the parameter $\Theta$ defined by 
  \begin{eqnarray} \label{eq:criticalparking} \Theta = (1-m_{\mu})^{2}- \Sigma_{\xi}^{2}( \sigma_{\mu}^{2} + m_{\mu}^{2}-m_{\mu}),  \end{eqnarray}
indicates whether the parking is subcritical or supercritical. More precisely, when $\Theta > 0$ (the subcritical case), we have $\mathbb{E}[\mathrm{Cars}(\rho)] < \infty$ under the unconditioned monotype Galton–Watson law, whereas if $\Theta < 0$ (the supercritical case), then $\mathbb{E}[\mathrm{Cars}(\rho)] = \infty$. Heuristically, in the subcritical regime, most cars manage to find a parking spot, resulting in a small outgoing flux of cars. In contrast, in the supercritical case, significant traffic jams occur, leading to an outgoing flux proportional to the size of the underlying tree. See Figure \ref{fig:parkingsimu}.

  \begin{center} \hrulefill \textit{~In the rest of this section we suppose critical parking, i.e. that $ \Theta =0$.} \hrulefill  \end{center}

\subsection{Universality of fully parked trees and the Brownian growth-fragmentation tree}\label{subsection:universality:fully}

Under the above assumptions, it is tempting to conjecture that under $ \mathbb{Q}_{n}^{\GW, \mathrm{park}}$, the rescaled decorated tree $\texttt{T}_{\GW}=(T_\GW,  g_{\GW})$ will converge to an interesting ssMt. But this is false! The reason is that although we do have convergence of the decoration-reproduction processes (Assumption \ref{assum:main:gen:deco:repro}), unfortunately Assumption \ref{A:all:existance:PHI} and \eqref{eq:non:explo:1} do not hold: Indeed, if we forget about the conditioning\footnote{In the supercritical case, the conditioning on $\{ \mathrm{Cars}(\rho) =p\}$ is strong enough to ensure that $ \mathbb{E}[|T| \mid  \mathrm{Cars}(\rho)=p ] < \infty$ but not in the critical case} of the label of the root vertex, the underlying tree is a critical monotype Galton--Watson tree for which $ \mathbb{E}[|T|] = \infty$ and which does not possess non-trivial finite super-harmonic function. The remedy is simple enough: we shall just truncate the tree at the vertices reaching label $0$ (no car visiting it) which we see as fictitious particles (that is why we did not shift the labels before). This boils down to consider the so-called \textbf{fully parked} component of the root in $\mathrm{Park}(T,(\car_u)_{u \in T})$, see \cite{chen2021enumeration,ConCurParking} and below. We shall denote the associated laws by $ (\mathbb{Q}_{n}^{\GW, \mathrm{fpark}})_{n \geq 1}$, obtained from the previous proposition by just cutting-off particles of label $0$.

Once we restrict to the fully parked component of the root, we have the following universality result in the bounded degree case:
\begin{theorem}[{\cite{ConCur25}}] \label{thm:ConCur25} In the critical parking case, when the distribution $\mu, \xi$ have bounded support, there exist  constants $c_{d}, c_{v}>0$ such that if we consider  the integer-type Galton--Watson tree of law $ \Q^{\GW,\mathrm{fpark}}_{n}$ described above, the sequence of the distributions  of the rescaled measured decorated trees  $$\left( T_\GW, \frac{d_{T_\GW}}{n^{3/2}} \cdot c_d, \rho_\GW, \frac{g}{n}, \frac{ \uplambda_\GW}{n^{2}} \cdot c_v\right) \mbox{ under } \Q^{\GW,\mathrm{fpark}}_{n}$$
converges as $n\to \infty$, to the law  of 
the  Brownian growth-fragmentation tree with index $3/2$, i.e.~the measured self-similar Markov tree $(\normalfont{\texttt{T}},\upmu)$  under $\Q_1$ with characteristic quadruplet 
$$\big( 0,   \mathrm{a}_{\mathrm{BroGF}}, \boldsymbol{ \Lambda}_{\mathrm{BroGF}} ;  \frac{3}{2}\big),$$ see Example \ref{ex:3/2stable} and the discussion after it.
\end{theorem}
The proof of the above theorem is heavily inspired by this work and may be seen as the discrete counter part to our Section \ref{sec:examplespinal}. In short, it consists in establishing spine decomposition in the discrete setting revealing that the labels along branches of those trees are intimately connected to random walk. Obviously, it is expected that the bounded degree assumption in the  above result could be slightly relaxed. We shall however see in the next section that if one allows car arrivals with heavy tails then we fall in another universality classes of ssMt. 

\subsection{Critical parking with stable car arrivals}

 In the following, we focus exclusively on the case when $\xi$ is the critical geometric distribution $\xi(k) = 2^{-k-1}$ for $k \geq0$ (so that $T$ is a uniform plane tree conditionally on its size) to be able to use the results of Chen \cite{chen2021enumeration}. We also restrict to critical case $\Theta=0$ i.e. according to \eqref{eq:criticalparking} when $2 \sigma_{\mu}^{2} + m_{\mu}^{2} =1$, and when furthermore the car arrival law has a heavy tail $ \mu(k) \sim \mathrm{cst}\cdot k^{-2 \beta-1},  \mbox{with }\beta \in (1, \frac{3}{2})$ for some constant  $ \mathrm{C}_{\mu}>0$. Let us record these assumptions for further use
  \begin{eqnarray} \label{eq:assumlinxiao} \underbrace{\left\{ \begin{array}{l} \xi(k) = 2^{-k-1}\\ \forall k \geq 0 \end{array}\right.}_{ \mbox{ geo. critical trees}}, \qquad  \underbrace{2 \sum_{k \geq 0} \mu(k) k^{2} - \left(\sum_{k \geq 0} \mu(k) k\right)^{2} =1,}_{ \mbox{critical parking}} \qquad  \underbrace{\left\{\begin{array}{l} \mu(k) \underset{k \to \infty}{\sim} \mathrm{C}_{\mu}\cdot k^{-2 \beta-1}\\ \beta \in (1,\frac{3}{2})\\
   \mathrm{C}_{\mu}>0. \end{array}\right.}_{ \mbox{cars heavy tail}}  \end{eqnarray}

\begin{theorem}[Scaling limit of critical non-generic fully parked trees] \label{thm:parkingstable}Fix $\beta \in (1, \tfrac 32)$ and assume \eqref{eq:assumlinxiao}.  Considering the integer-type Galton--Watson tree of law $ \Q^{\GW,\mathrm{fpark}}_{n}$ described above, the sequence of the distributions  of the rescaled decorated trees 
$$\left( T_\GW,\frac{d_{T_\GW}}{n^{\beta}} \cdot 2 \sqrt{\mathrm{C}_{\mu} |\Gamma(-2\beta)|} , \rho_\GW,\frac{g_{\GW}}{n}\right) \mbox{ under } \Q^{\GW,\mathrm{fpark}}_{n}$$
converges as $n\to \infty$, to the law  of the self-similar Markov tree $\normalfont{\texttt{T}}$  under $\Q_1$ associated with the characteristics $(0,  \mathrm{a}_{ 1,\beta-1}, \boldsymbol{\Lambda}_{ 1,\beta-1}; \beta)$  of  Example \ref{ex:stablefamily} when $ \texttt{a}=1, \texttt{b} = \beta-1$ so that $\frac{1}{2} < \beta = \texttt{a} + \texttt{b} < \frac{3}{2}$ and $\varrho = \frac{1}{\beta}$ (i.e.~spectrally negative and conservative case, but with presence of a killing term), see  \eqref{eq:GLMstablespecneg}.

\end{theorem}
\begin{remark}  We  conjecture that the above convergence holds jointly for the natural mass measure on the leaves, but this would require to check Assumption \ref{A:hGWrv} which we have not adressed in these pages. \end{remark}

The proof of Theorem \ref{thm:parkingstable} will occupy the remainder of Section \ref{subsection:universality:fully}. We begin by recalling results of Chen \cite{chen2021enumeration} on the exact and asymptotic enumeration of plane fully parked trees. These results are essential to our proof, as they provide access to the quantities $\Omega^{(p)}$ through a decomposition of parked trees into fully parked components, as developed by Chen and Contat \cite{chen2024parking}. Readers less familiar with analytic combinatorics may treat the following section as a black box for enumeration purposes and proceed directly to the proof of Theorem \ref{thm:parkingstable}.

\subsubsection{Chen's numeration of plane fully parked trees}

\begin{figure}[!h]
 \begin{center}
 \includegraphics[width=10cm]{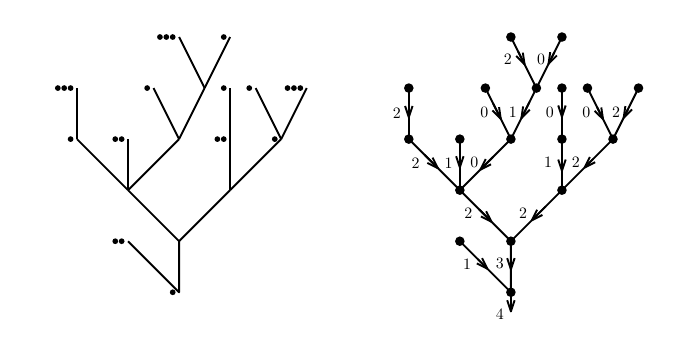}
 \caption{Left: A parking configuration. Right: the resulting fully-parked tree with flux of car on each edge. \label{fig:fpt}}
 \end{center}
 \end{figure}

A fully-parked tree is a plane tree $ \mathrm{t}$ together with cars arrivals $(\car_u : u \in \mathrm{t})$ such that after the parking process, every vertex contains a car, that is, with the notation introduced above that $ \mathrm{Cars}(u) \geq 1$ for every $u \in \mathrm{t}$, see Figure~\ref{fig:fpt} for an example. In particular $ \mathrm{Cars}(\rho)-1$ is the number of out-going cars that did not manage to park. Given a probability distribution $\mu$, we shall enumerate fully-parked trees by size and flux of out-going cars. 

In this direction, we form the generating function with a weight $x$ per vertex and $y$ per outgoing car at the root:
\begin{equation}
F(x,y) \equiv F(x,y;\mu) = \sum_{\begin{subarray}{c}( \mathrm{t}, (\smallcar_u)_{u \in \mathrm{t}})\\ \mathrm{fully-parked} \end{subarray}} x^{ \# \mathrm{t}} y^{ \mathrm{Cars}(\rho)-1} \prod_{u \in \mathrm{t}} \mu_{\smallcar_u},
\end{equation}
where the sum runs over all fully-parked configurations. The weight factor $y^{ \mathrm{Cars}(\rho)-1}$ is to match the definition of \cite{chen2021enumeration} and  notice that $F(0,y) = 0$. Introducing $B(z) = \sum_{k \geq 0} \mu_k z^k$ for the generating function of $\mu$, a  recursive decomposition ``\`a la Tutte'' gives the following equation of formal power series in $ \mathbb{C}[[x,y]]$:
\begin{equation}\label{eq:F(x,y)}
F(x,y) = \frac xy \Big( \frac{B(y)}{1-F(x,y)} - \frac{\mu_0}{1-F_0(x)} \Big)
\end{equation}
where $F_0(x):= F(x,0)$ is the generating function of fully parked trees without overflow. It is not hard to see that those equations uniquely determines the series $F(x,y)$ order by order in $x$ (with the initial condition $F(0,y)=0$), see \cite[Section 2]{chen2021enumeration} for details.

We now recall the parametric solution for $F(x,y)$ i.e.~express it in a new system of coordinates $ \hat{x}$. This is done in \cite[Section 2]{chen2021enumeration} inspired by the method  from Bousquet-M\'elou \& Jehanne \cite{BMJ06}, see also \cite[Section 2.2]{chen2024parking}.  
To do that, let $ \hat{x}(Y)$ defined by 
\begin{eqnarray}\label{eq:F_0(x) GOD} \hat{x}(Y) =x= \frac{Y \ B(Y)}{(B(Y)+Y B'(Y))^{2}},  \end{eqnarray}
then we have   \begin{eqnarray} \hat{F}_{0}(Y) = F_{0}(\hat{x}(Y))= 1- \frac{\mu_{0} B(Y)}{\big((B(Y))^{2}- (Y B'(Y))^{2}\big)}.  \end{eqnarray} 

The previous displays can be interpreted as a parametrization of $x$ by the auxiliary variable $Y$ (at least as long as this parametrization is legal), so that in this new coordinate system the function $F_{0}$ is very simple to express.  Once we have access to $F_{0}$, the equation \eqref{eq:F(x,y)} becomes a simple quadratic equation in $F(x,y)$ which has a unique solution that is a Taylor series in $y$:
\begin{equation}\label{eq:F(x,y) GOD}
F( \hat{x}(Y),y) = \frac12 + \frac{ \sqrt{ (\varphi(Y)+y)^2 - 4y B(y)\cdot \hat x(Y) } - \varphi(Y) }{2y}
\end{equation}
where $\varphi(Y) := Y \cdot \frac{1- \psi_B(Y)}{1+ \psi_B(Y)}$ and $\psi_B(Y) := \frac{Y B'(Y)}{B(Y)}$.

 \paragraph{Criticality.} The parametrization of $x$ via $Y$ in \eqref{eq:F_0(x) GOD} enables us to understand the critical role played by \eqref{eq:assumlinxiao}. Let us sketch the phenomenon for the reader's convenience and refer to \cite{chen2021enumeration} for details. Since $B$ has only non-negative coefficients, it is easy to see that $ Y \mapsto \hat{x}(Y)$ is bimodal (increasing then decreasing) over $ \mathbb{R}_{+}$, in particular the parametrization of $x$ by $Y$ is valid as long as $\hat{x}(Y)$ is strictly increasing. The criticality condition for the parking process is equivalent to asking that 
  \begin{eqnarray} \label{eq:criticparking}  \partial_{Y}\hat{x}(Y)|_{Y=1} =0 \iff (1- B'(1))^{2}- 2 B''(1) = 1 - m_\mu^{2} -2 \sigma_\mu^{2} =0.  \end{eqnarray}

In this case, since $B$ has radius of convergence exactly $1$, it is easy to see that $x_c := \hat{x}(1)$ is also the radius of convergence of $F_0(x)$. We then have the following universal formulas:
  \begin{eqnarray} \label{eq:universalformulas} x_{c} = \left(\frac{1}{1+m}\right)^{2}, \quad
  \frac{1}{ \sqrt{x_{c}}} F(x_{c},1) = m \quad \mbox{ and } \quad  \frac{1}{ \sqrt{x_{c}}} \partial_{y} F(x_{c},1) = \frac{1}{2}(1-m).  \end{eqnarray}
The last equation $\frac{1}{ \sqrt{x_{c}}} \partial_{y} F(x_{c},1) = \frac{1}{2}(1-m)$ is in fact true for all critical parking on critical Galton--Watson trees, see \cite[Theorem 1(ii)]{CH19}. Finally, when $ \mu$ is critical and has tails as in \eqref{eq:assumlinxiao}, we have the following asymptotics from \cite[Theorem 3]{chen2021enumeration} and a few straightforward simplifications:
  \begin{eqnarray}  \frac{1}{ \sqrt{x_{c}}}[y^{p}] F(x_{c},y) \underset{p \to \infty}{\sim} \frac{\sqrt{  \mathrm{C}_{\mu} |\Gamma(-2\beta)|}}{\Gamma(-b)} p^{-1-\beta}.  \label{eq:taillinxiao}\end{eqnarray}

\paragraph{Back to $\Omega^{(p)}$.} Now that we have the necessary background on enumeration of plane fully parked trees, let us translate this into knowledge on $\Omega^{(p)}$, which we recall is the probability (under the monotype critical geometric Galton--Watson tree and i.i.d.cars arrivals of law $\mu$) that the number of cars visiting the root is equal to $p \geq 0$. This relies on a decomposition of the parked tree $ \mathrm{Park}(T, (\car_u)_{u \in T})$ used in \cite{chen2024parking}, see also \cite{aldous2023parking}.

  \begin{proposition}[{\cite[Lemma 3 and Proposition 2]{chen2024parking} in the special case $q=1/2$}] \label{prop:alicelinxiao} Suppose that $\mu$ is critical i.e.~that $2 \sigma_{\mu}^{2} + m_{\mu}^{2}= 1$. Let $ {T}$ be a critical geometric $(1/2)$ Galton--Watson tree with i.i.d.~car arrivals of law $\mu$ on each vertex. Then the law of the number of cars visiting the root  satisfies
  $$ \Omega^{(p)}  = \mathbb{P}(\mathrm{Cars}(\rho)=p) = \left\{\begin{array}{ll} \displaystyle \frac{1}{\sqrt{x_{c}}} [y^{p+1}] F(x_{c},y), & \mbox{ if } p \geq 1, \\ 
\displaystyle  1-m & \mbox{ if } p=0. \end{array} \right. $$
In particular, \eqref{eq:taillinxiao} can be rewritten in the form:
\begin{equation}\label{asymp:Omega:p}
\Omega^{(p)} \underset{p \to \infty}{\sim} \frac{\sqrt{  \mathrm{C}_{\mu} |\Gamma(-2\beta)|}}{\Gamma(-b)} p^{-1-\beta}. 
\end{equation}
 \end{proposition}
Here also, the last item of the proposition,  $ \mathbb{P}( \mathrm{Cars}(\rho)=0) =1-m$, is a general fact for critical parking on critical Galton--Watson trees, see \cite[Theorem 1(iii)]{CH19}.

\subsubsection{Proof of Theorem \ref{thm:parkingstable}}
Now that we have gathered understanding on the function $\Omega^{(\cdot)}$ driving $\Q^{\GW,\mathrm{fpark}}_{n}$ via  \eqref{eq:lawGWpark}, we can dive into the proof of Theorem \ref{thm:parkingstable} for which the reader will have recognized yet another application of our Theorem \ref{T:mainuncondw} for which we need to check  Assumptions  \ref{assum:main:gen:deco:repro} and  \ref{A:all:existance:PHI}. We start by checking the convergence of the generalized L\'evy measure since it is the first time we encounter a non-trivial killing term in the limit.

\begin{proposition}[Convergence of the generalized L\'evy measure]  \label{prop:measurekilling}For the locally largest selection rule we have $ \boldsymbol{\Lambda}^{(n)} \to 2 \sqrt{ \mathrm{C}_\mu | \Gamma(-2\beta)|} \boldsymbol{\Lambda}_{1, \beta-1}$ as $n \to \infty$ in the sense of Assumption \ref{A:BK} with the notation of \eqref{eq:GLMstablespecneg}.
\end{proposition}
The reader may be puzzled, as we only prove the convergence of the generalized Lévy measure and not the convergence of the other characteristics required for Assumption \ref{A:BK}. The reason is that, although it should in principle be possible to also establish the convergence of the drift and variance, we will instead verify directly Assumption \ref{assum:main:gen:deco:repro}  combining  only the convergence of the generalized Lévy measure and connections with random walks.
\begin{proof}
Fix a continuous function $F$ vanishing in the neighborhood of $(1,0,0,...)$ and suppose that $F$ is also continuous at $(0,0,...)$. Start with a particle of label $n$ and denote by $c\geq 0$ the number of car arrivals on this vertex, $k\geq 0$ its number of children as well as $p_1, ..., p_k$ the labels of its offspring ordered from left to right. We write $p_i^*$ for the ordered decreasing sequence. 
We start with an easy lemma. Once $k\geq 1$ and the car arrivals $c\geq 0$ are fixed, one has to choose the flux $p_1, ..., p_k$ so that $\sum (p_i-1)_+ = n-c$. We first prove: 
\begin{lemma} \label{lem:measurekilling} For any $k \geq 0$ we have 
  \begin{eqnarray*}  \frac{n^\beta}{\Omega^{(n)}} \sum_{\begin{subarray}{c} 0 \leq p_1, ... , p_k\leq n \\ (p_1-1)_+ +  ... + (p_k-1)_+ =n\end{subarray}} \prod_{i=1}^k \Omega^{(p_i)} F\left( \frac{p^*_1}{n}, \frac{p_2^*}{n}, ... \right)&& 
  \end{eqnarray*}
  \begin{eqnarray} \label{eq:faciledabord} &\xrightarrow[n\to\infty]{}& \displaystyle k(k-1)   \int_{[1/2,1]} F( x, 1-x,0,0,...)   \frac{\sqrt{  \mathrm{C}_{\mu} |\Gamma(-2\beta)|}}{\Gamma(-b)} \cdot \frac{ \mathrm{d}x}{(x(1-x))^{-1-\beta}}. \end{eqnarray}
In particular, the splitting become binary and conservative in the limit. Furthermore, the left-hand side is uniformly bounded by some constant (depend on $F$) times $k^{4+\beta}$. 
\end{lemma} 
\begin{proof} Fix $A >0$ large, $m \in \{2,3, ... \}$ and consider the contribution when there are at least $m$ number $p_{i}$'s which are larger than $A$:
$$ S(n,k,A ; m) := \sum_{\begin{subarray}{c} 0 \leq p_1, ... , p_k\leq n \\ (p_1-1)_+ +  ... + (p_k-1)_+ =n\\ p_{i}\geq A \mathrm{\  for\ at \ least \  } m \mathrm{ \ indices} \end{subarray}} \prod_{i=1}^k \Omega^{(p_i)}.$$ Noticing that the largest $p_{i}$ must satisfy $ (p_{i}-1)_{+} \geq n/k$ we can upper bound the previous display by
  \begin{eqnarray*}S(n,k,A ; m) &\leq& m {k \choose m} \sum_{\begin{subarray}{c} 0 \leq p_1, ... , p_k\leq n \\ (p_1-1)_+ +  ... + (p_k-1)_+ =n\\ p_{i}\geq A, \forall i \leq m\\ 
p_{1} \geq n/k-1  \end{subarray}} \prod_{i=1}^k \Omega^{(p_i)}\\
&\leq  &m {k \choose m} \cdot \Big( \sup_{ \ell \geq n/k-1} \Omega^{{(\ell)}} \Big)\cdot \Big( \sum_{p_{2} \geq A} \Omega^{(p_{2})} \cdots \sum_{p_{m} \geq A} \Omega^{(p_{m})}\Big)\cdot \Big( \underbrace{\sum_{p_{m+1}, ... , p_{k}}  \prod_{i=m+1}^k \Omega^{(p_i)}}_{=1}\Big)\\
&\leq &  m {k \choose m}  \mathrm{Cst} \cdot (n/k)^{-\beta-1}   ( \mathrm{Cst}\cdot A^{-\beta})^{m-1} \leq \mathrm{Cst}^{m} k^{m+1+\beta+1} n^{-\beta-1} A^{-(m-1)\beta},  \end{eqnarray*}
where to obtain the last line we use the asymptotic \eqref{asymp:Omega:p}  and $ \mathrm{Cst}>0$ is a constant  only depending on the sequence $(\Omega^{{(p)}})$.  Once divided by $\Omega^{(n)}$ the term in $n^{-\beta-1}$  disappear and we see that the contribution of those splittings for which at least three $p_{i}$'s are of order $n$ is $O(n^{-2 \beta})$ which, even multiplied by $n^{\beta}$ tends to $0$. The convergence and the upper bound easily follows from those considerations. We leave the details to the reader.
\end{proof}

Given this result, one can condition on $c,k$ first and apply \eqref{eq:faciledabord}
  \begin{eqnarray*}
&&  n^\beta \int_{ \mathcal{S}} \boldsymbol{\Lambda}^{(n)}( \mathrm{d}  \mathbf{y}) F(  \mathrm{e}^{y_0},  \mathrm{e}^{y_1}, ...)\\
&=&  \frac{n^\beta}{\Omega^{(n)}} \sum_{c \geq 0} \mu_c \sum_{k \geq 0} 2^{-k-1} \sum_{\begin{subarray}{c}0\leq p_1, ... , p_k  \\
\sum (p_i-1)_+ = n-c \end{subarray}} \prod_{i=1}^k \Omega^{(p_i)} F\left( \frac{p^*_1}{n}, \frac{p_2^*}{n}, ... \right)\\
&=&  \frac{n^\beta}{\Omega^{(n)}} \sum_{c \geq 0} \mu_c \sum_{k \geq 0} 2^{-k-1}  \frac{\Omega^{(n-c)}}{(n-c)^\beta} \cdot  \underbrace{\frac{(n-c)^\beta}{\Omega^{(n-c)}}\sum_{\begin{subarray}{c}0\leq p_1, ... , p_k  \\
\sum (p_i-1)_+ = n-c \end{subarray}} \prod_{i=1}^k \Omega^{(p_i)} F\left( \frac{p^*_1}{n}, \frac{p_2^*}{n}, ... \right)}_{}  \end{eqnarray*}
For $c$ of order $o(n)$ the underbraced sum is converging towards $k(k-1) I_F$ where $I$ stands for the integral in \eqref{eq:faciledabord}. Using the domination by $  \mathrm{Cst}\cdot k(k-1)k^\beta I_F$, Proposition \ref{prop:alicelinxiao} and dominated convergence shows the convergence towards $ \mathrm{Var}(\xi) \cdot I_F$ when restricting to small (say $o(n)$) values of $c$. However, for those large values of $c=n-\ell$ where $\ell$ is of order $1$ the underbraced integral converges towards 
$$ \frac{\ell^\beta}{\Omega^{(\ell)}}F(0,0,0,...) \sum_{\begin{subarray}{c}0\leq p_1, ... , p_k  \\
\sum (p_i-1)_+ = \ell \end{subarray}} \prod_{i=1}^k \Omega^{(p_i)}.$$
 For those values of $c$, we have
$$ \frac{n^\beta}{\Omega^{(n)}}   \mu_{n-\ell}  \frac{\Omega^{(\ell)}}{(\ell)^\beta} \xrightarrow[n\to\infty]{ \mathrm{Prop. \ } \ref{prop:alicelinxiao}} \frac{ \mathrm{C}_{\mu}}{\frac{\sqrt{  \mathrm{C}_{\mu} |\Gamma(-2\beta)|}}{\Gamma(-b)}}  \frac{\Omega^{(\ell)}}{(\ell)^\beta}.$$
The middle terms $ 0 << c << n$ are easily neglected using crude bounds and we get finally that 
  \begin{eqnarray*} n^\beta \int_{ \mathcal{S}} \boldsymbol{\Lambda}^{(n)}( \mathrm{d}  \mathbf{y}) F( y_0, y_1, ...)  &\xrightarrow[n\to\infty]{}&  2 \cdot  \frac{\sqrt{  \mathrm{C}_{\mu} |\Gamma(-2\beta)|}}{\Gamma(-b)} \cdot \int_{[1/2,1]} F( x, 1-x,0,0,...)  \frac{ \mathrm{d}x}{(x(1-x))^{-1-\beta}} \\
  & +&  \frac{ \mathrm{C}_{\mu}}{\frac{\sqrt{  \mathrm{C}_{\mu} |\Gamma(-2\beta)|}}{\Gamma(-b)}} F(0,0,0,...) \underbrace{\sum_{k \geq 0} 2^{-k-1}\sum_{\ell=0}^\infty  \sum_{\begin{subarray}{c}0\leq p_1, ... , p_k  \\
\sum (p_i-1)_+ = \ell \end{subarray}} \prod_{i=1}^k \Omega^{(p_i)}}_{=1}\\
&\underset{ \eqref{eq:GLMstablespecneg}}{=}& 2 \sqrt{ \mathrm{C}_\mu | \Gamma(-2b)|} \int_{\mathcal S} F\big( \mathrm{e}^{y_0}, (\mathrm{e}^{y_{1}}, \mathrm{e}^{y_{2}}, \cdots)\big) \boldsymbol{\Lambda}_{1, \beta-1}( \dd y_0, \dd \mathbf{y}),  \end{eqnarray*}
where the last equality follows from classical formulas involving Gamma functions.
\end{proof}

We now introduce a random walk similar to the one used in the proof of Theorem \ref{prop:GF}, see also \cite{ConCur25} for a general setup. Specifically, let us define a measure on $ \mathbb{Z} \times \mathbb{Z}_{+}^{ \mathbb{Z}_{+}}$ by putting 
  \begin{eqnarray} \label{eq:defnubar} & &\bar{\nu}\big( q,( i_{0}, \dots , i_{n},0,\dots)\big)\\
  & :=&  2^{i_{0}+ \dots + i_{n}} \cdot \mu_{1-q-\sum_{\ell=0}^{n} i_{\ell}(\ell-1)_{+} } \cdot   \frac{(i_{0}+ \dots + i_{n}+1)!}{i_{0}!i_{1}!\dots i_{n}!} \prod_{\ell=0}^{n} \big(\Omega^{(\ell)}\big)^{i_{\ell}}.  \end{eqnarray} 
This measure is heuristically obtained by looking at the evolution of a particle of very large label $p$ for the multi-type Galton--Watson measure:  Using the polynomial asymptotic on $\Omega^{{(p)}}$ it is easy to see from \eqref{eq:lawGWpark} that when $p \to \infty$, with probability $\bar{\nu}\big( q,( i_{0}, \dots , i_{n},0,\dots)\big)$ this particle will give rise to a single largest particle of label $p+q$ and $i_{\ell}$ particles of label $\ell \geq 0$ (i.e. so that $\ell$ cars visit those vertices, and $(\ell-1)_{+}$ of them will visit the vertex of label $p$).  By  the compatibility condition \eqref{eq:parking}, for this to be possible, $k=1-q-\sum_{0\leq \ell \leq n} i_\ell(\ell-1)_+$ cars must necessary arrive on the vertex. We stress that vertices of label $0$ may appear, i.e. $i_{0}$ may be positive, but recall that those vertices are immediately frozen. Remark then that necessarily $q \leq 1$, i.e.~that the projection $\bar{\nu}_{0}$ of $\bar{\nu}$ on the first coordinate is a measure on $\{\dots,-3,-2,-1-,0,1\}$.  With the notation of the previous section, if $X\sim \bar{\nu}_{0}$ then an application of  \eqref{eq:lawGWpark},  and expressing $ \Omega^{(\cdot)}$ using Proposition \ref{prop:alicelinxiao}, show that $1-X$ has generating function given by 
$$ \sum_{k \geq 1} 2^{-k-1} \cdot k \cdot B(y) \left( \frac{1}{\sqrt{x_{c}}}F(x_{c},y) + (1-m)\right)^{k-1} = \frac{B(y)}{(\frac{1}{\sqrt{x_{c}}}F(x_{c},y) -m-1)^{2}}.$$ In particular it follows from \eqref{eq:universalformulas} that $\bar{\nu}_{1}$ hence $\bar{\nu}$ is indeed a probability measure and that $X$ is centered. Furthermore,  from \eqref{eq:taillinxiao} we see that $X$ is in the strict domain of attraction of the spectrally negative $\beta$-stable law. We can already use $\bar{\nu}$ to check:\\
 
\noindent $\bullet$ \textsc{Assumption \ref{A:all:existance:PHI}:} The function $ l_{\GW}^{ \mathbf{1}_{  \mathbb{N}}}(n) = \mathbb{Q}^{\GW, \mathrm{fpark}}_{n}( \# T_{\GW})$ is clearly a super-harmonic function for $ \Q^{\GW, \mathrm{fpark}}$ and according to \cite[Theorem 3]{chen2021enumeration} for $n \geq 1$ we have 
$$ \mathbb{Q}^{\GW, \mathrm{fpark}}_{n}( \# T_{\GW}) \underset{ \mbox{ \tiny \cite[Lemma 3]{chen2024parking}}}{=} \frac{ \frac{1}{ \sqrt{x_{c}}} [y^{n}]\partial_{x} F(x_{c},y)}{ \frac{1}{ \sqrt{x_{c}}} [y^{n}] F(x_{c},y)} \\
\underset{\mbox{\tiny\cite[Theorem 3]{chen2021enumeration}} }{\sim} \mathrm{Cst}\cdot n^{2 \beta-1},$$ for some constant  $ \mathrm{Cst}>0$ and where $2 \beta-1$ corresponds to the first root $\omega_{-}$ of the cumulant of the ssMt described in  Example \ref{ex:stablefamily} when $ \texttt{a}=1, \texttt{b} = \beta-1$. This provides the second super-harmonic function corresponding to the exponent $\gamma_{1}$ in Assumption \ref{A:all:existance:PHI}. The first one, corresponding to the exponent $\gamma_{0}$ is provided by $$\phi_{0}(n) :=  \frac{1}{\Omega^{(n)}} \overset{ \mathrm{Prop.\ } \ref{prop:alicelinxiao} \ \&\  \eqref{eq:taillinxiao}}{\underset{n \to \infty}{\sim}} \mathrm{Cst'} \ n^{\beta+1}$$ for some other constant $ \mathrm{Cst'}>0$ where indeed $\omega_{-} = 2\beta-1 < \beta+1$. To check that $\phi_{0}$ is super-harmonic, we use a calculation similar to the one performed in \cite[Proof of Theorem 4]{ConCur25}: notice that for $n \geq 1$ using \eqref{eq:lawGWpark} we have 
 \begin{eqnarray*} && \Q_{n}^{\GW, \mathrm{fpark}}( \phi_{0}( \mathbf{Z}(1)))\\
 & \underset{ \eqref{eq:lawGWpark}}{=}&
  \sum_{i_{0}, \dots , i_{n+1}}2^{i_{0} + \dots + i_{n+1}-1} \cdot  \displaystyle \mu_{n- \sum_{\ell=0}^{n+1} i_\ell\
 (\ell-1)_{+}} \cdot  \frac{(i_{0} + \dots + i_{n+1})!}{i_{0}! i_{1}! \dots i_{n+1}!}  \frac{\prod_{\ell=0}^{n+1} \big(\Omega^{(\ell)}\big)^{i_{\ell}}}{{\Omega^{(n)}}} \sum_{\ell=0}^{n+1} i_{\ell} \phi_{0}(\ell)\\
  & \underset{\phi_{0}(\cdot) = 1/\Omega^{(\cdot)}}{=}&
 \frac{1}{\Omega^{(n)}}\sum_{i_{0}, \dots , i_{n+1}}2^{i_{0} + \dots + i_{n+1}-1}   \sum_{\ell=0}^{n+1}   \mathbf{1}_{i_{\ell} \geq 1}  \displaystyle \mu_{n- \sum_{\ell=0}^{n+1} i_\ell\
 (\ell-1)_{+}} \\ && \qquad \times   \frac{(i_{0} + \dots + (i_{\ell}-1) + i_{n+1})!}{i_{0}! \dots (i_{\ell}-1)! i_{1}! \dots i_{n+1}!}  \big(\Omega^{(0)}\big)^{i_{0}} \dots \big(\Omega^{(\ell)}\big)^{i_{\ell}-1} \dots \big(\Omega^{(n)}\big)^{i_{n+1}}.  \end{eqnarray*}
We now perform the change of variable $j_{0} = i_{0}, j_{1}= i_{1}, \dots , j_{\ell} = i_{\ell}-1, \dots , j_{n+1} = i_{n+1}$ and $\ell = n+q$ so that the previous display becomes
 \begin{eqnarray*} 
\frac{1}{\Omega^{(n)}}\sum_{q \geq -n} \sum_{j_{0}, \dots , j_{n+1}}2^{j_{0} + \dots + j_{n+1}}    \displaystyle \mu_{n- (n+q-1)_{+} - \sum_{\ell=0}^{n+1} j_\ell\
 (\ell-1)_{+}} \cdot  \frac{(j_{0} + \dots + j_{n+1}+1)!}{j_{0}!  \dots j_{n+1}!} \prod_{k \geq 0} \big(\Omega^{(k)}\big)^{j_{k}}.  \end{eqnarray*}
 In the above display, when $n+q \geq 1$ we have $n-(n+q-1)_{+} = 1-q$ and we recover from \eqref{eq:defnubar} the definition of $\bar{\nu}\big( q,( j_{0}, \dots , j_{n+1},0,\dots)\big)$ in particular we must have $1-q - \sum_{\ell=0}^{n+1} j_\ell (\ell-1)_{+} \geq 0$. When $q=-n$ we have a slight shift and get instead $\bar{\nu}\big( 1-n,( j_{0}, \dots , j_{n+1},0,\dots)\big)$. In particular, recalling that $\overline{\nu}_{0}$ is the push-forward of $\bar{\nu}$ on its first coordinate, we finally get that 
   \begin{eqnarray*}
 \Q_{n}^{\GW, \mathrm{fpark}}( \phi_{0}( \mathbf{Z}(1))) & = & \frac{1}{\Omega^{(n)}} \left(\sum_{q \geq -n+1} \bar{\nu}_{0}(q)  + \bar{\nu}_{0} (1-n)\right).
 \end{eqnarray*}
 Since we already remarked that $\bar{\nu}_{0}$ is in the strict domain of attraction of the $\beta$-stable law, we deduce that the parenthesis in the previous display falls short of being $1$ up to a factor of order $n^{-\beta}$: It follows that for some $ \mathrm{C}>0$ we have 
$$\boldsymbol{G}_{\GW} \phi_0(n) \leq \phi_{0}(n) \left( 1- \frac{ \mathrm{C}}{n^{\beta}}\right)$$
so that $\phi_{0}(n)$ is indeed super-harmonic, $\phi_{0}(n) \asymp n^{\gamma_{0}}$ with $\gamma_{0}= \beta+1$ and (recalling that the self-similarity parameter $\alpha$ is equal to $\beta$) 
$$ \limsup_{n \to \infty} n^{\alpha-\gamma_0}\boldsymbol{G}_{\GW} \phi_0(n) < 0,$$ as required to check Assumption \ref{A:all:existance:PHI}.

\noindent $\bullet$ \textsc{Assumption \ref{assum:main:gen:deco:repro}:} We finally check the convergence of the decoration-reproduction processes to complete the proof, the proof is inspired from \cite{ConCur25}. For this, we again use the probability measure $ \bar{\nu}$. Specifically, we use it to create a decoration-reproduction process $(f, \eta)$ which under the probability measure 
\newcommand{\SRW}{ {\normalfont\text{\tiny{RW}}}}
${P}_{n}^{ \SRW}$ is defined as follows:
Let $(X_{i} ; (Y_{j}^{(i)})_{j \geq 1})$ be i.i.d.~copies of law $\bar{\nu}$ and then define 
$$ f(0) = n, \quad f(k) = f(0) + \sum_{\ell=1}^{k} X_{i}$$
and kill the process (and set it to $0$) at $z :=  \inf\{ \ell \geq 0 : f(k) \leq 0\}$ the first time when it enters $ \mathbb{Z}_{\leq 0}$. The decoration process is then 
$$ \eta = \sum_{i=1}^{z} \sum_{j \geq 1}\delta_{(i,Y_{j}^{(i)})}.$$
This decoration-reproduction process modeled on iid increments can then be used to recover information about the decoration-reproduction process under $ P_{n}^{\GW}$ for the locally largest selection rule. 
If $(f, \eta)$ is the underlying decoration-reproduction process, for $t \geq 0$ we shall denote by $ \mathcal{L}_{t}$  the event where $(f, \eta)$ corresponds to a strict locally largest exploration up to time $t$, i.e. if $z>t$ and 
 $$f(j) > x \quad \mbox{for any atom\ } (j,x) \mbox{	 \ of } \eta, \mbox{  with } 0 \leq j \leq t.$$ 
 \begin{lemma}[Similar to {\cite[Proposition 10]{ConCur25}}] \label{lem:h-transformnu} For any positive function $F$ we have
 \begin{eqnarray} \label{eq:h-transformnu} E_{n}^{\GW}\left( F(f|_{[0,t]}, \eta|_{[0,t]})\mathbf{1}_{  \mathcal{L}_{t}} \right) =  E_{n}^{\SRW}\left(F(f|_{[0,t]}, \eta|_{[0,t]}) \mathbf{1}_{  \mathcal{L}_{t}} \frac{\Omega^{(f(t))}}{\Omega^{(f(0))}}\right). \end{eqnarray}
 A similar equation holds if we restrict to (non-strict) locally largest exploration but involves a combinatorial factor accounting for the possible choices of exploration in case of ties.
 \end{lemma}
 
 \begin{proof} By the product form of the above displays, it is sufficient to check the equality in the case of $t=1$. Restricting to locally largest exploration with no ties, 
we can use   \eqref{eq:lawGWpark}, to rewrite  $P_{n}^{\GW}( f(1) = n+q)$, for $n+q \geq 1$, in the form:
\begin{align*}
  \sum_{i_{0}, ... , i_{n+q-1}} 2^{i_{0} + \dots + i_{n+q-1}} \cdot   &\mu_{1-q- \sum_{\ell=0}^{n+q-1} i_\ell\
 (\ell-1)_{+}}   \frac{(i_{0} + \dots + i_{n+q-1} +1)!}{i_{0}! i_{1}! \dots i_{n+q-1}!}  \frac{\prod_{\ell=0}^{n+q-1} \big(\Omega^{(\ell)}\big)^{i_{\ell}}}{{\Omega^{(n)}}} \Omega^{(n+q)}\\
 &=  \frac{\Omega^{(n+q)}}{\Omega^{(n)}} \bar{\nu}(q; (i_{0}, ... , )) \mathbf{1}_{ i_{k} =0, \forall k \geq n+q} \mathbf{1}_{q > n} = P_{n}^{\SRW}(f(1)=n+q). \end{align*}
 \end{proof}
 
The scaling limit of the decoration-reproduction process under $P_{n}^{\SRW}$ is easy to get since under $ P_{n}^{\SRW}$ the process $f$ is a $\nu$-random walk started from $n$. Recall that $\nu$ is in the strict domain of  attraction of the $\beta$-stable spectrally negative law, and under $ \mathbf{P}_{1}$ we let $(\xi_{t} : t \geq 0)$ be a $\beta$-stable spectrally negative L\'evy process (see \eqref{eq:levymeasurestable} in Section \ref{sec:spinalex} for the normalization), then by \cite[Theorem 15.17]{Kal07} or \cite[Chapter VII]{JS03} we have convergence in the Skorokhod sense:
$$ \left(\frac{f(\lfloor n^{\beta}t \rfloor)}{ n}\right)_{t \geq 0}  \mbox{\  under } {P}^{\SRW}_{n} \quad \xrightarrow[n\to\infty]{(d)} \quad  \left(  {c} \cdot \xi_{t}\right)_{t \geq 0}  \mbox{ \ under } \mathbf{P}_{1}, \qquad \mbox{ where } {c} >0,$$
{where the constant $c = \big( \sqrt{ \mathrm{C}_{\mu }|\Gamma(-2\beta)|}\big)$ is computed from \eqref{eq:taillinxiao} and \eqref{eq:levymeasurestable}.} The preceding convergence is easily extended to get a scaling limit of the decoration-reproduction processes under $ {P}_n^{ \SRW}$:

\begin{lemma}[Scaling limit of the decoration-reproduction processes]  \label{lem:scalinglimit} Recall from \eqref{Eq:scalefeta} the notation $(f^{(n)}, \eta^{(n)})$ for the rescaled decoration-reproduction process with $\alpha= \beta \in (1, \frac{3}{2})$. Then we have 
  \begin{eqnarray}  \left( f^{(n)} , \eta^{(n)}\right) \mbox{ \ under } {P}_{n}^{ \SRW} \quad \xrightarrow[p\to\infty]{(d)} \quad  \Big( ({c}\cdot  \xi^{}_{t})_{t \geq 0} , \sum_{\begin{subarray}{c}t \geq 0 \\    \Delta \xi^{}_{t}>0 \end{subarray}} \delta_{t, {c}\cdot \Delta \xi^{}_{t}} \Big) \mbox{ \ under } \mathbf{P}_{1},  \label{lem:scalinglimitneutral} \end{eqnarray}   where the convergence holds for the product topology with the Skorokhod topology for the first coordinate and  the vague convergence of measures on $ \mathbb{R}_{+} \times \mathbb{R}_{+}^{*}$ for the second one.
\end{lemma}
\begin{proof} The convergence of the first coordinate is granted by the above discussion. We thus need to show that the atoms of the measure are prescribed by the jumps of the process. 
Since $f$ cannot increase by more than $1$,  the proof reduces to showing that a large negative jump $\Delta f(i) = X_{i}$ corresponds to a single $Y^{(i)}_j$ for some $j \geq 1$ (as opposed to  several of them whose total sum would roughly be equal to $\Delta f(i)$). However, the calculations performed in Proposition \ref{prop:measurekilling} and Lemma \ref{lem:measurekilling} ensures that for any $ \varepsilon>0$ we have
$$ \sup_{r \geq 0 }{P}^{\SRW}_{r}\left( \sum_{j\geq 1 }Y^{1}_{j} - \max_{j \geq 1} Y^{1}_{j} \geq \varepsilon p\right) = o(p^{-\beta}),$$ so that a union bound completes the proof.\end{proof}

Armed with those lemmas, we can now establish Assumption \ref{assum:main:gen:deco:repro}, i.e.~the convergence of the laws  $P_{\GW}^{(n)}$ of the rescaled  decoration-reproduction process (for the locally largest selection rule).  Since we already established Assumption \ref{A:all:existance:PHI}, by Lemma \ref{lem:convreprodecolight}, it suffices to establish the convergence of the law $P_{\GW}^{{(n), \varepsilon}}$ for the exploration freezed when dropping below level $  \varepsilon>0$. We just focus on the decoration process for simplicity, since the reproduction component is, in the limit, only consisting of the (negative of) the jumps of the process by Lemma \ref{lem:scalinglimit}. Fix $ \varepsilon>0$ and $F$ a bounded continuous function of the Skorokhod topology and $G$ a bounded continuous function on $[0, \varepsilon]$. We suppose that $F(\omega)=0$ vanishes as soon as the process $\omega$ drops below level $  \varepsilon>0$ or if the process makes a negative jump such that $ \omega(t-)\leq \omega(t)/2$. Recall from Lemma \ref{lem:convreprodecolight} the notation $\theta^{ \varepsilon} = \inf\{ u \geq 0 : f(u) \leq \varepsilon\}$. Recalling \eqref{Eq:scalefeta}, we can write
  \begin{eqnarray*} &&E_{\GW}^{(n)}\left[ F\left( f^{}|_{ [0, \theta^{ \varepsilon})}\right) G(f(\theta^{ \varepsilon}))\right] \\
  &=& \sum_{k \geq 1} E_{n}^{\GW}\left[ F\left( f^{(n)}|_{[0, k n^{-\beta})} \right)  \cdot E_{f(k-1)}^{\GW}\left[G( n^{-1} f(1))  \mathbf{1}_{f(1) \leq \varepsilon n}\right]\right]\\
    & \underset{ \mathrm{Lem.}\  \ref{lem:h-transformnu}}{=}& \sum_{k \geq 1} E_{n}^{\SRW}\left[ F\left( f^{(n)}|_{[0, k n^{-\beta})} \right) \frac{\Omega^{(f( k-1))}}{\Omega^{(n)}} \cdot E_{f(k-1)}^{\GW}\left[G( n^{-1} f(1))  \mathbf{1}_{f(1) \leq \varepsilon n}\right]\right].
  \end{eqnarray*}
   One can now apply the asymptotics on $\Omega^{(\cdot)}$, Lemma \ref{lem:scalinglimit} and Proposition \ref{prop:measurekilling} to see that the previous display converges as $n \to \infty$ to 
\begin{eqnarray*}    &&   2 \sqrt{ \mathrm{C}_\mu | \Gamma(-2\beta)|}  \int_{0}^{\infty} \mathrm{d}t  ~ \mathbf{E}_{1}\left[F( \xi|_{[0,t)})\mathcal{}  \cdot (\xi_{t-})^{-\beta-1} \cdot \int_{-\infty}^{\log ( \varepsilon \xi_{t-}^{-1})} \Lambda_{0} ( \mathrm{d}x) G( \xi_{t-} \mathrm{e}^{x}) \right],
    \end{eqnarray*}
     where $\Lambda_{0}$ is the projection on the first coordinate of the generalized L\'evy measure $\boldsymbol{\Lambda}_{1, \beta-1}$ appearing in Proposition \ref{prop:measurekilling}. One then checks that the law of the process defined by the above formula is indeed the law of the pssMp used for the decoration in Example \ref{ex:stablefamily} when $ \texttt{a}=1, \texttt{b} = \beta-1$, which is stopped when dropping below level $ \varepsilon$. See \cite[Proof of Lemma 17]{ConCur25} for similar calculations. We leave the details to the reader.

 \section{Miscellaneous}
 We conclude this chapter devoted to discrete applications with a few models of random labeled trees which we believe should converge towards ssMt. We do not claim any mathematical result and hope that these questions will trigger new works both in combinatorics and in probability. We start with yet another application to the peeling process on random planar maps.
  
  \subsection*{Simple peeling}
The peeling we used in Section \ref{sec:randomplanarmaps}  is called the ``lazy'' peeling exploration; it has been introduced by Budd \cite{Bud15}, see \cite{CurStFlour} for a comprehensive survey about its applications. It is more flexible that the  \textbf{simple} peeling process of Angel \cite{Ang03} which was actually used first in the mathematic literature. The main difference is that the ``lazy'' peeling process is based on maps with general boundaries and the transitions can split the map in at most two parts, whereas the simple peeling process is based on maps with a simple boundary where the transitions can be very complicated. We refer to \cite[Section 4.3]{CurStFlour} for a detailed comparison. When the random maps have small face degrees, we expect that the labeled trees associated to the simple peeling process  are in the same universality classes as in that of Theorem \ref{prop:GF} in the case $\beta= \frac{3}{2}$. However, in the case of Boltzmann measures (with a simple boundary) with a critical weight sequence of type $\beta \in (1/2, 3/2)$, the simple peeling transitions are considerable more involved and may yield to new ssMt with non binary generalized L\'evy measures. We leave those questions open and believe that  new important objects are to be found there.

 \subsection*{Fighting fish} 
\begin{figure}[!h]
  \begin{center}
  \includegraphics[width=15cm]{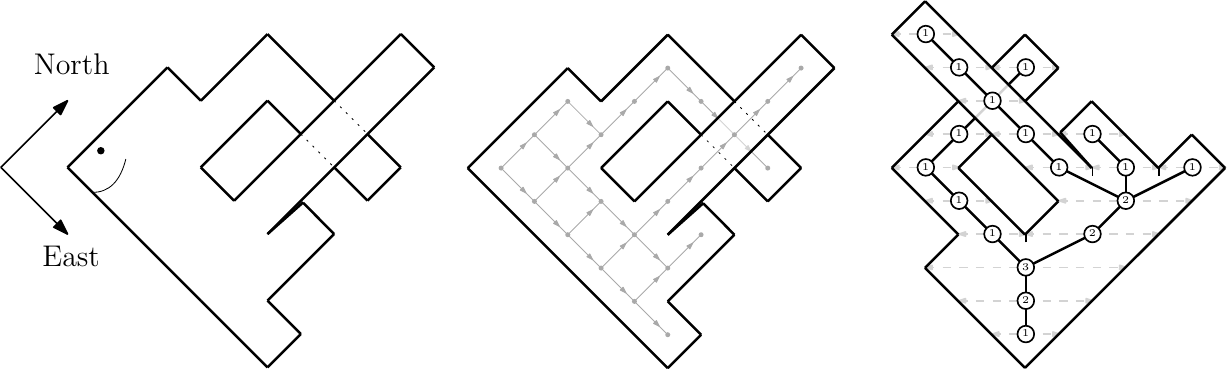}
  \caption{A fighting fish obtained by the oriented gluing of (possibly branching) squares and its representation (not bijective) as a labeled plane tree. \label{fig:fighting}}
  \end{center}
  \end{figure}

A fighting fish is a branched surface over $ \mathbb{Z}^{2}$ that is obtained by gluing together flexible unit squares along their edges in an oriented way (only North or East gluings), resulting in independent branches that may overlap. The model has been introduced and studied in \cite{duchi2017fighting,duchi2017fightingbis} mainly from the enumerative point of view, and many connections with the theory of planar map counting have been drawn \cite{cioni2023direct,duchi2022bijection,duchi2023bijection}. A labeled tree can naturally be associated (non bijectively) to any fighting fish, see Figure \ref{fig:fighting}, and they are multi-type Galton--Watson trees when the underlying fighting fish is taken under the (critical) Boltzmann measure i.e. when each square comes with a weight $ \frac{4}{27}$. The enumerative results cited above suggested that those random trees (when started from a large decoration) converge after scaling towards the same ssMt as in Theorem \ref{prop:GF} in the case $\beta= \frac{3}{2}$.

\subsection*{Weil-Petersson random hyperbolic surfaces}
\label{sec:WP}  
To conclude, let us show that multi-type Galton--Watson trees are also underlying the geometry of random hyperbolic surfaces. The following discussion is inspired from \cite{BBCPpeeling}. For $n \geq 2$ we denote by $ \mathcal{M}(n,p)$ the moduli space of hyperbolic surfaces of genus $0$, with $n$ punctures and a boundary of length $p$. This is a natural finite measure $\mu_{n,p}$  on each space $\mathcal{M}(n,p)$ called the Weil-Petersson measure and whose total mass is denoted by $ \mathrm{WP}_{n,p}$. We can thus consider the Boltzmann measure on $\cup_{n \geq 2}\mathcal{M}(n,p)$ defined by 
$$ \sum_{n \geq 2} x^{n} \mu_{n,p},$$ which is well-defined as long as $ \sum_{n \geq 2}  \mathrm{WP}_{n,p}\cdot x^n < \infty$. The radius of convergence of this series is known and it is finite at this point, enabling us to consider the critical Boltzmann measure, cf. \cite{MR1234274}. A re-interpretation in \cite{BBCPpeeling} of the work of Mirzakhani \cite{mirzakhani2007weil} shows that similarly as what we can do in planar maps, given a peeling algorithm $ \mathcal{A}$ we can encode almost  all (i.e.~up to measure zero) surfaces of $\cup_{n \geq 2}\mathcal{M}(n,p)$ by a labeled plane tree. This labeled tree encodes a decomposition of the underlying hyperbolic surface into a pant decomposition and the labels of the vertices correspond to the hyperbolic lengths of the pants in the decomposition.
\begin{figure}[!h]
 \begin{center}
 \includegraphics[width=12cm]{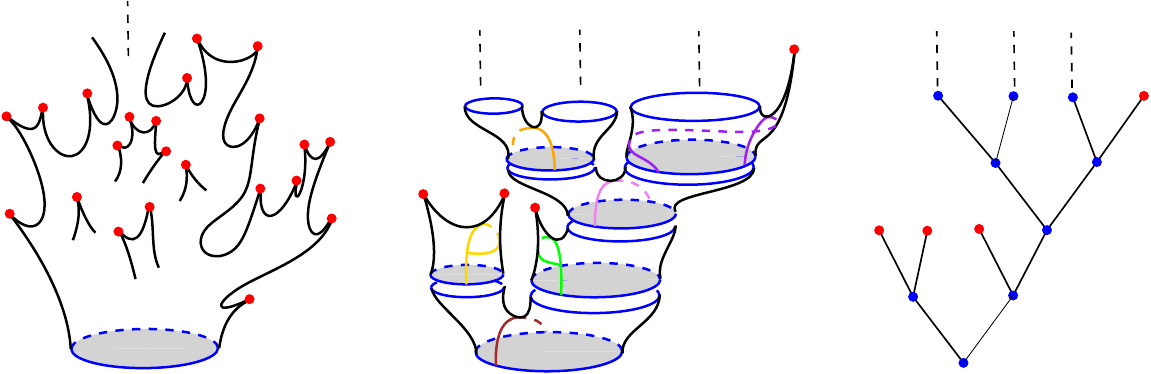}
 \caption{From hyperbolic surfaces sampled according to Weil-Peterson measure to a decomposition in pair of pants following Mirzakhani \cite{mirzakhani2007weil}.}
 \end{center}
 \end{figure}
 Although this does not fall directly in our setup since the labels of the tree are positive reals instead of integer, our techniques also apply in this setting and enable to prove a scaling limit result for those tree similar to that of Theorem \ref{prop:GF} in the case $\beta= \frac{3}{2}$. We hope to have sparked the reader's curiosity.

\chapter{Appendix} \label{chap:appendix}
This appendix contains a few technical results which may be of broader interest. We decided to single them out for clarity and accessibility; the sections below are independent one from the others. 
\section{Moments of exponential functionals of a L\'evy process}
 Let $\xi$ be a Lévy process (possibly with killing). We write $P$ for its law and $\psi$ for  its Laplace exponent as defined in Section \ref{sec:2.2}; see   \eqref{eq:levykhintchine}. 
 We consider the so-called exponential functionals
 $$I(\beta)\coloneqq \int_{0}^{\zeta} \exp(\beta \xi(t)) \d t,$$
 where $\beta >0$ is a parameter and  $\zeta$ stands for the lifetime of $\xi$.
These  random variables have been much studied in the literature on L\'evy processes, see e.g. \cite{bertoin2005exponential} for background.
 
 \begin{lemma}\label{gene:levy:rivero} Assume that 
 $\psi(\gamma)<0$ for some $\gamma>0$. Then, we have
\begin{equation}\label{eq:sup:rivero:k}
E\Big(\sup_{t\geq 0} \exp\big(\gamma\xi(t)\big)\Big)<\infty \quad \text{ and }\quad  E\Big(I(\beta)^{\gamma/\beta}\Big)<\infty,
\end{equation}
for every $\beta>0$. 
\end{lemma}
\begin{proof}
We work under the assumption of the lemma and remark that, since $\psi(\gamma)<0$, 
$$E\Big(I(\gamma)\Big)=\int_{0}^{\infty} E\big(\exp(\gamma \xi(t)), t<\zeta \big)\d t =\int_{0}^{\infty} \exp(t\psi(\gamma))\d t=-\frac{1}{\psi(\gamma)}<\infty.$$
The idea is to use the variable $I(\gamma)$ to bound the expectation in the left-hand side of  \eqref{eq:sup:rivero:k}. In this direction, we fix $b>0$ such that $P_1(I(\gamma)\geq b)>0$, and, for every $r>0$, we introduce the stopping time $T_r:=\inf\{s\geq 0:~\xi(s)\geq \log(r)/\gamma\}$. An application of the strong  Markov property then yields
\begin{align*}
P\big(I(\gamma)\geq b\cdot r\big)&\geq P\Big(\big\{T_r<\infty\big\}\cap \Big\{\int_{T_r}^{\zeta} \exp\big(\gamma (\xi(t)-\xi(T_r))\big)\d t\geq b\cdot r\exp\big(-\gamma \xi(T_r)\big)\Big\}\Big)\\
&\geq P\big(T_r<\infty\big)\cdot P\big(I(\gamma)\geq b\big),
\end{align*}
where in the second line we used that $\xi(T_r)\geq  \log(r)/\gamma$ on the event $\{T_r<\infty\}$. Henceforth, we have
 $$P\Big(\sup_{t\geq 0}  \exp(\gamma \xi(t))\geq  r\Big)\leq \frac{P\big(I(\gamma)\geq b r\big)}{P(I(\gamma)\geq  b)},\quad \text{ for }r>0.$$
Integrating over $r>0$ yields 
 \begin{align}\label{eq:sup:xi:gamma}
E\Big(\sup_{t\geq 0} \exp\big(\gamma\xi(t)\big)\Big)\leq \frac{E(I(\gamma))}{b \cdot P(I(\gamma)\geq b)}<\infty.
 \end{align}

We now turn our attention to the second inequality in \eqref{eq:sup:rivero:k},  which has already been  established by Rivero \cite[Lemma 2]{rivero2007recurrent} when  $\beta\geq \gamma$. Suppose now  $\beta< \gamma$ and consider the integer part of $\gamma/\beta$,  $k=\lfloor \gamma/\beta\rfloor \geq 1$.  Since $\psi(\gamma)<0$ and $\psi(0)\leq 0$, the convexity of $\psi$ entails that $\psi(\gamma-\ell\beta)\leq 0$ for every integer $1\leq \ell\leq k$. Thanks to \cite[Lemma 2.1]{Maulik:Zwart},  we have that
$$E\Big(I(\beta)^{\gamma/\beta}\Big)=\Big(\prod_{0\leq\ell\leq k} \frac{\gamma-\ell \beta}{-\psi(\gamma-\ell \beta)}\Big)\cdot E\Big(I(\beta)^{(\gamma-k\beta)/\beta}\Big).$$
Since $\gamma-k\beta\leq \beta$, we can now apply \cite[Lemma 2]{rivero2007recurrent} to infer that the previous display is finite. This completes the proof of the lemma.
\end{proof}

\section{A law of large numbers for uniformly integrable variables}
 Let $(a_{i,j})_{i,j\in \N}$ be null array of nonnegative random variables on some probability space $(\Omega, \mathcal{F}, \mathbb{P})$, that is, for any $\varepsilon >0$,
$$\#\big\{(i,j)\in \N^2: a_{i,j}>\varepsilon\big\} < \infty \qquad \text{a.s.}$$ 
 Consider   a uniformly integrable family $\{\upbeta_{i,a}: i\in \N \, \text{ and }\, a\geq 0\}$ of nonnegative random variables with $\E(\upbeta_{i,a})=1$.
Let $(b_{i,j})_{i,j\in \N}$ be a second array of random variables, such that for every $i,j\in \N$,
conditionally  given the sequence $(a_{i,k})_{k\in \N}$, the r.v.~$b_{i,j}$ has the same law as $\upbeta_{i,a}$ for $a=a_{i,j}$, independently of the other variables
$(b_{i,k})_{k\neq j}$.

\begin{lemma} \label{L:tailored} 
If the sequence $A_i\coloneqq \sum_{j\geq 1} a_{i,j}$ converges in $L^1(\P)$ as $i\to \infty$ to some random variable $A$, then we have also
$$\lim_{i\to \infty} \sum_{j\geq 1} a_{i,j}b_{i,j}= A, \qquad \text{in }L^1(\P).$$
\end{lemma}

\begin{proof} It suffices to check the conclusion of claim with the weaker convergence in probability, since then the stronger convergence in $L^1$ follows from the comparison of mathematical expectations and  an application of the Scheff\'e lemma.
Fix some arbitrarily small $\varepsilon >0$. 
To start with, the assumption of uniform integrability for the family $(\upbeta_{i,a})$ enables us to choose some $n\geq 1$ sufficiently large such that 
$$m_i(a)\coloneqq \E\big(\upbeta_{i,a}\indset{\upbeta_{i,a}\leq n}\big)\in(1- \varepsilon, 1],\qquad \text{for all } i\geq 1 \text{ and }a\geq 0.$$
We then set $b'_{i,j}\coloneqq b_{i,j}\indset{b_{i,j}\leq n},$
and note that 
$$\E\Big( \sum_{j\geq 1} a_{i,j} \big( b_{i,j}-b'_{i,j}\big) \Big)\leq \varepsilon \E(A_i). $$
We next write
$$\sum_{j\geq 1} a_{i,j} (b'_{i,j} -1)= \sum_{j\geq 1} a_{i,j} (m_i(a_{i,j})-1)+ \sum_{j\geq 1}a_{i,j}(b'_{i,j}-m_i(a_{i,j}))  .$$
On the one-hand, we have 
$$0\leq \sum_{j\geq 1} a_{i,j}(1- m_i(a_{i,j})) \leq \varepsilon A_i.$$
On the other hand, introduce the event
$$\Lambda_i(\varepsilon)\coloneqq \big\{a_{i,j}\leq \varepsilon n^{-2} \text{ for all }j\geq 1\big\}.$$
Since conditionally on $(a_{i,j})_{j\geq 1}$, the variables $b'_{i,j}-m_i(a_{i,j})$ are independent, centered, and bounded in absolute value by $n$, we have
$$\E\Big( \Big| \sum_{j\geq 1}a_{i,j}(b'_{i,j}-m_i(a_{i,j})) \Big|^2 \indset{\Lambda_i(\varepsilon)}\Big) \leq \E\Big(  \sum_{j\geq 1}a_{i,j}^{2} n^{2} \indset{\Lambda_i(\varepsilon)}\Big)\leq \E\Big(  \sum_{j\geq 1}a_{i,j} \cdot \varepsilon n^{-2} \cdot n^{2}\Big) =  \varepsilon \E(A_i).$$

Putting the pieces together, we have shown that 
$$\E\Big( \Big| \sum_{j\geq 1}a_{i,j}(b_{i,j}-1) \Big| \indset{\Lambda_i(\varepsilon)}\Big) \leq 2\varepsilon \E(A_i) +\sqrt{\varepsilon \E(A_i)}.$$
To conclude the proof, it suffices to observe that $\lim_{i\to \infty} \P(\Lambda_i(\varepsilon))=1$, since $(a_{i,j})$ is a null array.
\end{proof}

\section{Tail of first passage time of a supermartingale with polynomial drift}
Let  $(X_k)_{k \geq 0}$ be a nonnegative supermartingale in some filtered probability space $(\Omega, \mathcal{F}, (\mathcal{F}_k)_{k\geq 0}, \mathbb{P})$. 
We assume that $ \mathbb{E}(X_0) < \infty$ and then  $ \mathbb{E}(X_k) < \infty$ for all $k \geq 0$.
We are interested in the tail distribution of the first passage time under $1/2$, $\tau = \inf \{ k \geq 0 : X_k\leq1/2\}$.
\begin{lemma}[]  \label{lem:supermartingaleextinct} 
Suppose that $\tau<\infty$ a.s., and  that there exist two constants $ \mathrm{c}>0$ and $\delta  \in (-1, \infty)$ such that 
\begin{equation}\label{eq:lem:unif:mart:delta}
\mathbb{E}\big(X_{k+1}-X_k \mid \mathcal{F}_{k}\big) \leq -  \mathrm{c} \cdot X_{k}^{-\delta},\quad \text{ for }\, k\leq \tau.
\end{equation}
Then we have
 $$ \sup \limits_{x>1}x^{\frac{1}{1+\delta}}\cdot \mathbb{P}( \tau \geq x X_0^{1+ \delta})\leq C(\delta),$$
where $C(\delta)=(2^{-1}\wedge \mathrm{c})^{-\frac{1}{1+\delta}}$, if $\delta \in(-1,0]$, and $C(\delta)=1+ \frac{2^{1+2\delta}}{(2^{\delta}-1)\mathrm{c}} $, if $\delta>0$.
\end{lemma}

\begin{proof}
By working under the conditional law $\P(\cdot \mid X_0)$, there is no loss of generality supposing  that  $X_0$ is deterministic, and we write for convenience $X_0=n$.
 In the case when $-1<\delta \leq 0$, the result follows from work of Aspandiiarov and Iasnogorodski \cite[Theorem 2]{aspandiiarov1999general}. More precisely, writing $X_k = f(X_k^{1+\delta})$ with $f(x) = x^{ \frac{1}{1+\delta}}$ which is convex in the neighborhood of $+\infty$ because $\delta < 0$, and since $f'(X_k^{1+\delta}) = (1+ \delta) X_k^{-\delta}$, we can apply  the  last display of the proof of \cite[Theorem 2]{aspandiiarov1999general} to obtain
$$ \mathbb{E}\big(f(\lambda \tau)\big) \leq  \ f(X_0^{1+\delta}) =  X_0,$$
for every $\lambda< 2^{-1}\wedge \mathrm{c}$, where $\mathrm{c}$ is the constant appearing in the statement. The desired result now follows by an application of Markov inequality,  since $f$ is increasing.
 
Let us move to  the case $\delta >0$,  and  introduce the sequence of stopping times
$$S_m:=\inf\big\{k\geq 0:~X_k\geq 2^m n\big\}\wedge \tau, \quad m\geq 0. $$
 The supermartingale property of $X$ ensures that $2^{m}n\mathbb{P}(S_m<\tau)\leq \mathbb{E}(X_{S_m})\leq n$ and thus:
\begin{eqnarray} \mathbb{P}( S_m< \tau) \leq 2^{-m}, \quad \text{for } m\geq 0. \label{eq:superma}  \end{eqnarray} 
Fix now $x>1$ and  let  $m_x$  be the smallest integer so that $(2^{m_x})^{1+\delta} \geq x$. Combining  \eqref{eq:superma} with the definition of $m_x$ and Markov inequality, we infer that
$$\mathbb{P}(\tau\geq x n^{1+\delta})\leq \mathbb{P}(S_{m_x}<\tau)+ \mathbb{P}(S_{m_x}\geq  x n^{1+\delta})\leq x^{-\frac{1}{1+\delta}}+ \frac{\mathbb{E}\big(S_{m_x}\big)}{x n^{1+\delta}}. $$
 It remains to control the quantity $\mathbb{E}\big(S_{m_x}\big)$. To this end, we write:
\begin{equation}\label{eq:teles:S:mx}
\mathbb{E}\big(S_{m_x}\big)=\sum \limits_{m=0}^{m_x-1} \mathbb{E}\big(S_{m+1}-S_m\big).
\end{equation}
We now remark that, for every $m\geq 0$, we have $X_k\leq 2^{m+1} n$ for $S_m\leq k< S_{m+1}$, and then by \eqref{eq:lem:unif:mart:delta} the process
$$\Big(X_{(S_m+k)\wedge S_{m+1}}+ \mathrm{c} n^{-\delta} 2^{-\delta (m+1)} k\Big)_{k\geq 0} $$
is a supermartingale, where $ \mathrm{c}$ is the constant appearing in the statement of the lemma. Since $S_{m+1}$ is almost surely finite, we get
$$\mathbb{E}\big(X_{S_{m+1}}\big)+ \mathrm{c} n^{-\delta} 2^{-\delta(m+1)}\mathbb{E}\big(S_{m+1}-S_m\big)\leq \mathbb{E}(X_{S_m}). $$ As a result, we obtain
$$\mathbb{E}\big(S_{m+1}-S_m\big)\leq \mathrm{c}^{-1} n^{\delta} 2^{\delta(m+1)}\cdot \mathbb{E}\big(X_{S_m}-X_{S_{m+1}}\big),\quad m\geq 0.$$ Moreover, since
on the event $\{\tau\leq S_m\}$, we have $X_{S_m}=X_{S_{m+1}}=X_{\tau} \in [0,1/2)$, and on the complementary event   $\{\tau> S_m\}$, we have $X_{S_m}-X_{S_{m+1}}\leq 2^{m+1} n$, an application of \eqref{eq:superma} gives:
\begin{align*}
\mathbb{E}\big(S_{m+1}-S_m\big)\leq \mathrm{c}^{-1} n^{1+\delta} 2^{(1+\delta)(m+1)}\cdot \mathbb{P}\big(S_m<\tau\big)\leq \mathrm{c}^{-1} n^{1+\delta} 2^{1+\delta(m+1)},\quad \text{for } m\geq 0.
\end{align*}
Coming back to \eqref{eq:teles:S:mx}, by the previous display we derive that:
$$\mathbb{E}\big(S_{m_x}\big)\leq  2^{1+\delta} \mathrm{c}^{-1} n^{1+\delta} \sum \limits_{m=0}^{m_x-1}2^{\delta m} \leq \frac{2^{1+\delta}}{(2^{\delta}-1)\mathrm{c}} n^{1+\delta} 2^{\delta m_x }\leq \frac{2^{1+2\delta}}{(2^{\delta}-1)\mathrm{c}} n^{1+\delta} x^{\frac{\delta}{1+\delta}}$$
where to obtain the last inequality we used $2^{m_x}\leq 2x^{\frac{1}{1+\delta}}$. A direct computation now gives:
$$x^{ \frac{1}{1+\delta}}\cdot \mathbb{P}(\tau\geq x n^{1+\delta})\leq  1+ \frac{2^{1+2\delta}}{(2^{\delta}-1)\mathrm{c}}, $$
since $x>1$ is arbitrary we obtain the desired result.
\end{proof}

\section{A ratio limit and Laplace method}

 \begin{lemma}[] \label{lem:laplaceratio} Let $f,g : [0,1] \to \mathbb{R}$ be two continuous functions so that $f(1)=g(1)=1$. We suppose furthermore that $|f(x)|<1$ for $x\in[0, 1)$,  and that there exist $C>0$  and $\beta >0$ such that 
 $$1-f(1-x) \sim C x^{\beta}  \quad \mbox{ as } x \to 0+.$$
If we set 
 $$ u_{n} = \int_{0}^{1} \ g(x) f(x)^{n}  \mathrm{d}x, \quad \mbox{ for } n \geq 1,$$
 then 
 $$ u_{n} \sim C^{-1/\beta} \Gamma(1 + 1/a) \cdot n^{-1/\beta} \quad \mbox{ as $n \to \infty$},$$ and there exists $C'>0$ such that
  $$\left| \frac{u_{n+1}}{u_{n}}-1 \right| \leq \frac{C'}{n},\qquad \text{ for all $n \geq 1$.}$$ 
   \end{lemma}
   \begin{proof} The first assertion stems from the classical Laplace's method. Regarding the second,  note that
   $$ \left| u_{n}-u_{n+1} \right| \leq  \|g \|_{\infty}  \int_0^{1} \big(1-f(1-x)\big) f(1-x)^{n} \d x.$$
 Our assumptions ensure the existence of  $c_1,c_2,\delta>0$ such that 
 $$\big(1-f(1-x)\big) f(1-x)^{n}\leq c_1 x^{\beta}\exp(-c_2 n x^\beta)\qquad \text{ for every $n\geq 1$ and $x\in [0,\delta]$.}$$
 As a consequence,
 $$\int_0^{\delta}\big(1-f(1-x)\big) f(1-x)^{n} \d x~ \leq c_1 \int_0^\infty x^{\beta} \exp(-c_2 n x^{\beta}) \d x =c_1 \beta^{-2} \Gamma(\beta^{-1}) \cdot \big(c_2 n\big)^{-\frac{1+\beta}{\beta}}. $$
On the other hand,  $\int_\delta^1 \big(1-f(1-x)\big) f(1-x)^{n}\d x$ decays exponentially fast as $n\to \infty$, since $|f(1-x)|$ remains bounded away from $1$ for $x\in [\delta, 1]$. This implies that
$$\left| u_{n}-u_{n+1} \right|=O(n^{-\frac{1+\beta}{\beta}}), \quad \text{ as } n\to \infty,$$
and the desired result now follows by the first item. 
  \end{proof}
\renewcommand{\bibname}{References}
\addcontentsline{toc}{chapter}{References}

 \bibliographystyle{siam}
\bibliography{bibli}

\end{document}